\documentclass[12pt]{tac}   

\usepackage{latexsym, amsmath, amsfonts}
\usepackage{bbm}
\usepackage{enumerate, amssymb, xspace}
\usepackage[mathscr]{eucal}
\usepackage{graphicx}
\usepackage{rotating}
\usepackage{hyperref}
\usepackage{color,colordvi}
\usepackage{mathtools}   

\usepackage{tikz} 
\usepackage{tikz-cd}
\usetikzlibrary{matrix,arrows,tqft,calc,decorations,arrows,shapes}
\usetikzlibrary{decorations.markings}
\usetikzlibrary{decorations.pathmorphing}
\usetikzlibrary{arrows.meta}  
\usetikzlibrary{calc} 
 
\usepackage{makeidx}
\usepackage{pifont}
\usepackage{psfrag}
\usepackage{stackrel}

\newcommand\arrowDefect   {stealth'}
\newcommand\arrowDefectB  {stealth' reversed}
\newcommand\colorCircle   {brown!95}
\newcommand\colorCircleL  {brown!65}
\newcommand\colorDefect   {purple!93}
\newcommand\colorFrame    {blue!50!black}
\newcommand\colorFreebdy  {purple!66!black}
\newcommand\colorFTC      {red!70!black}
\newcommand\colorGraph    {blue!69!black}
\newcommand\colorHol      {blue!66!black}
\newcommand\colorSigma    {yellow!16}
\newcommand\colorSigmaDark{yellow!26}
\newcommand\colorSigmaDARK{yellow!32}
\newcommand\colorSigmaLight{yellow!12}
\newcommand\colorTransp   {purple!50}
\newcommand\fillBdyCircle {green!8}
\newcommand\fillCobord    {yellow!21}
\newcommand\transvL  {1.4pt}    
\newcommand\transvR  {0.10}     
\newcommand\transv[1]{\filldraw[\colorCircle,fill=white,line width=\transvL] (#1) circle (\transvR) }
\newcommand\transV[2]{\filldraw[\colorCircle,fill=white,line width=\transvL] (#1) circle (\transvR) node #2}
\newcommand\widthDefect   {1.3pt}
\newcommand\scopeArrow[2] {\begin{scope}[decoration={markings,mark=at position #1
                   with \arrow{#2}}]}  
\newcommand\drawOrientation {\draw[-{>[scale=1.5, length=1.6, width=2.9]},semithick]
                   (0,0) arc[radius=0.13,start angle=75,end angle=410]}
\newcommand\drawOrientationLarge {\draw[-{>[scale=2.5, length=1.6, width=2.9]},semithick]
                   (0,0) arc[radius=0.27,start angle=75,end angle=410]}

\newcommand\Cite[2] {\cite[#1]{#2}}

\theoremstyle {definition}
\newtheorem{thm}{Theorem}[section]

\newtheorem{Lemma}{Lemma}
\newtheorem{Proposition}{Proposition}
\newtheorem{cor}{Corollary}

\newtheorem{Definition}{Definition}

\newtheoremrm{Proofmain}{Proof of Theorem \ref{thm:main}}
\newtheoremrm{Remark}{Remark}
\newtheoremrm{Remarks}{Remarks}
\newtheoremrm{Example}{Example} 

\makeatletter
\newcommand*{\relrelbarsep}{.386ex}
\newcommand*{\relrelbar}{%
  \mathrel{%
      \mathpalette\@relrelbar\relrelbarsep
        }%
        }
\newcommand*{\@relrelbar}[2]{%
  \raise#2\hbox to 0pt{$\m@th#1\relbar$\hss}%
    \lower#2\hbox{$\m@th#1\relbar$}%
    }
\providecommand*{\rightrightarrowsfill@}{%
      \arrowfill@\relrelbar\relrelbar\rightrightarrows
      }
\providecommand*{\leftleftarrowsfill@}{%
       \arrowfill@\leftleftarrows\relrelbar\relrelbar
        }
\providecommand*{\xrightrightarrows}[2][]{%
          \ext@arrow 0359\rightrightarrowsfill@{#1}{#2}%
          }
\providecommand*{\xleftleftarrows}[2][]{%
    \ext@arrow 3095\leftleftarrowsfill@{#1}{#2}%
    }
\makeatother
\renewcommand\AA             {\mathbb{A}}
\newcommand\acalm            {{}_\cala \calm}
\newcommand\act              {\triangleright}
\newcommand\basicdisk        {straight disk}
\newcommand\BB               {\mathbb{B}}
\newcommand\Bfrdef           {\mathrm{Bord}_2^{\mathrm{def},0}}
\newcommand\Bfrdefgfine      {\mathrm{Bord}_2^{\mathrm{def},0,\mathrm{fine}}}
\newcommand\Bfrdefcl         {\mathrm{Bord}_{2,\mathrm{cl}}^{\mathrm{def},0}}
\newcommand\Bfrdefdec        {\mathrm{Bord}_2^{\mathrm{def}}}
\newcommand\Bfrdefdecgfine   {\mathrm{Bord}_2^{\mathrm{def},\mathrm{fine}}}
\newcommand\be               {\begin{equation}}
\newcommand\bearl            {\begin{array}{l}}
\newcommand\bearll           {\begin{array}{ll}}
\newcommand\bigboxtimes      {\mbox{\large$\boxtimes$}}
\newcommand\boti             {\,{\boxtimes}\,}
\newcommand\Boti[1]          {{\stackrel{#1}\boxtimes}}
\newcommand\bs[1]            {^{{\sss(}#1{\sss)}}}
\newcommand\cala             {{\mathcal A}}
\newcommand\Cala             {{\!\mathcal A}}
\newcommand\calb             {{\mathcal B}}
\newcommand\Calb             {{\!\mathcal B}}
\newcommand\calc             {{\mathcal C}}
\newcommand\cald             {{\mathcal D}}
\newcommand\calk             {{\mathcal K}}
\newcommand\calm             {{\mathcal M}}
\newcommand\Calm             {{\!\mathcal M}}
\newcommand\calmopp          {{\CC{\mathcal M}}}
\newcommand\caln             {{\mathcal N}}
\newcommand\calnopp          {{\CC{\mathcal N}}}
\newcommand\calx             {{\mathcal X}}
\newcommand\caly             {{\mathcal Y}}
\newcommand\calz             {{\mathcal Z}}
\newcommand\cc[1]            {\overline{{#1}}}  
\newcommand\CC[1]            {\overline{{#1}}} 
\newcommand\cenf[1]          {#1} 
\newcommand\cent             {\mathcal{Z}}
\newcommand\cir              {\,{\circ}\,}
\newcommand\coevl            {\mathrm{coev}^{\rm l}}
\newcommand\coevr            {\mathrm{coev}^{\rm r}}
\newcommand\Comod            {\text{\rm Comod}}
\newcommand\complex          {{\ensuremath{\mathbbm C}}}
\newcommand\csfun            {relative block functor}
\newcommand\cspfun           {relative pre-block functor}
\newcommand\DA               {D_{{\!\!\mathcal A}}}
\newcommand\DD               {\mathbb{D}}
\newcommand\DDsAMN           {\DD^{\uparrow}_{\Cala,\calm,\caln}}
\newcommand\DDt              {\DD_\text{\rm tr}}
\newcommand\DDtad            {\DD^\text{\rm tad}}
\newcommand\defectonemanifold{defect one-ma\-ni\-fold}
\newcommand\defectsurface    {defect surface}
\newcommand\DR               {\Phi}
\newcommand\dsty             {\displaystyle }
\newcommand\Edg              {E_\DD}
\newcommand\ee               {\end{equation}}
\newcommand\eear             {\end{array}}
\newcommand\eps              {\varepsilon}
\newcommand\eq               {\,{=}\,}
\newcommand\End              {{\rm End}}
\newcommand\Enumerate        {\renewcommand\leftmargini{1.34em} \begin{enumerate}} 
\newcommand\evl              {\mathrm{ev}^{\rm l}}
\newcommand\evr              {\mathrm{ev}^{\rm r}}
\newcommand\filldisk         {fillable disk}
\newcommand\filldiskrep      {fillable-disk replacement}
\newcommand\findim           {fini\-te-di\-men\-si\-o\-nal}
\newcommand\FinTensL         {\mathcal F\!\imath n\mathcal C\!at_\otimes^{\rm l.e.}}
\newcommand\Fun              {{\mathcal Fun}}
\newcommand\Funle            {{\mathcal Lex}}
\newcommand\FUNle            {{\mathcal{LE}\hspace*{-1.5pt}\mathcal{X}}}
\newcommand\FUNre            {{\mathcal{RE}\hspace*{-1.5pt}\mathcal{X}}}
\newcommand\Funlebal         {{\mathcal Lex}^{\rm bal}}
\newcommand\FunleSigmavect   {\Funle(\ZZ(\partial\Sigma),\vect)}
\newcommand\Funre            {{\mathcal Rex}}
\newcommand\GammaD           {\Gamma_{\!\rm \text{tot}}}
\newcommand\gfine            {{gluable fine}}
\newcommand\ggfine           {{\emph{gluable} fine}}
\newcommand\hatPhile         {\widehat\Phi^{\rm l}}
\newcommand\hatPsile         {\widehat\Psi^{\rm l}}
\newcommand\hoco             {\circ}   
\newcommand\hol              {\mathrm{hol}}
\newcommand\holc             {\gamma^\text{\rm c}}  
\newcommand\holcc            {\gamma^\text{\rm cc}} 
\newcommand\Hom              {{\rm Hom}}
\newcommand\id               {\mathrm{id}}
\newcommand\Id               {\mathrm{Id}}
\newcommand\II               {{\mathcal I}}
\newcommand\iN               {\,{\in}\,}
\newcommand\ind              {\mathrm{ind}}
\newcommand\Itemize          {\renewcommand\leftmargini{1.34em} \begin{itemize}} 
\newcommand\Itemizeiii       {\renewcommand\leftmargini{2.14em} \begin{itemize}} 
\newcommand\IZ               {I}   
\newcommand\JJ               {{\mathbb J}}
\newcommand\ko               {{\ensuremath{\Bbbk}}}
\renewcommand\L              {\mathrm L}
\newcommand\la               {{\rm l.a.}}
\newcommand\LII              {{\mathbb I}}
\newcommand\LL               {\mathbb L}
\newcommand\mute             {silent}
\newcommand\N                {\mathbb{N}}
\newcommand\Nat              {\mathrm{Nat}}
\newcommand\nn[1]            {N_{{#1}}}  
\newcommand\nxl[1]           {\\[.#1em]}
\newcommand\Nxl[1]           {\\[-1.2em]\\[.#1em]}
\newcommand\om               {\mho}
\newcommand\one              {{\bf1}}

\newcommand\opp              {^\text{\rm opp}}
\newcommand\Opp              {^{\mathrm{opp}_{}}}
\newcommand\oti              {\,{\otimes}\,}
\newcommand\otik             {\,{\otimes_\ko}\,}
\newcommand\partialf         {\partial_\text{\rm free}}
\newcommand\partialg         {\partial_\text{\rm glue}}
\newcommand\partialgg        {\partial_\text{\rm g}}
\newcommand\partiali         {\partial_-}
\newcommand\partialo         {\partial_+}
\newcommand\partialr         {\partial_\text{\rm fill}}
\newcommand\pift             {{\rm N}^{\rm r}}
\newcommand\pifu             {{\rm N}^{\rm l}}
\newcommand\Phile            {\Phi^{\rm l}}
\newcommand\Phire            {\Phi^{\rm r}}
\newcommand\PP               {{\mathbb P}}  
\newcommand\Ppop             {{\mathbb P}_\text{\!\rm p.o.p.}}  
\newcommand\pre              {\text{\rm pre}_{}}
\newcommand\Psile            {\Psi^{\rm l}}
\newcommand\Psire            {\Psi^{\rm r}}
\newcommand\QQ               {\mathbb Q}
\newcommand\R                {\mathbb{R}}
\newcommand\ra               {{\rm r.a.}}
\newcommand\refi             {_\text{\rm ref}}
\renewcommand\SS             {\mathbb{S}}
\newcommand\sss              {\scriptscriptstyle}
\newcommand\TT               {\mathbb T}
\newcommand\TX               {T^{\calx}}
\newcommand\TY               {T^{\caly}}
\newcommand\TZ               {T}
\newcommand\Times            {\,{\times}\,}
\newcommand\To               {\,{\to}\,}
\newcommand\tocong           {\xrightarrow{~\cong~}}
\newcommand\tosim            {\xrightarrow{~\simeq~}}

\newcommand\Tr               {\cup}
\newcommand\tadpolecircle    {tadpole circle}
\newcommand\tadpoledisk      {tadpole disk}
\newcommand\transdisk        {transparent disk}
\newcommand\twocell          {$2$-patch}
\newcommand\twocells         {$2$-patches}
\newcommand\undefectonemanifold {unlabeled \defectonemanifold}
\newcommand\undefectsurface  {unlabeled \defectsurface}
\newcommand\UU               {\mathrm U}
\newcommand\veco             {\ast}   
\newcommand\vect             {\ensuremath{\mathrm{vect}}}
\newcommand\vSigma           {\varSigma}
\newcommand\Vee              {{}^{\vee\!}}
\newcommand\widehathol       {\Gamma}
\newcommand\XX               {{\mathbb X}} 
\newcommand\Z                {\mathbb{Z}}
\newcommand\ZE               {{\mathrm E}} 
\newcommand\Zfine            {\ZZ_{\!\text{\rm fine}}}
\newcommand\Zpre             {\ZZ^\text{\rm pre}_{\phantom|}}
\newcommand\ZpreRE           {{\widetilde\ZZ^\text{\rm pre}}}
\newcommand\Zprev            {\ZZ^\text{\rm pre}}
\newcommand\Zpreprime        {\widehat{\ZZ^\text{\rm pre}_{\phantom|}}}
\newcommand\Zprime           {\widehat{\ZZ}}
\newcommand\ZZ               {{\mathrm T}} 

\newcommand\Arep[2]  {{\em #2}, available at {\tt #1}}
\newcommand\Bi[2]    {\bibitem[#2]{#1}}
\newcommand\inBO[9]  {{\em #9}, in:\ {\em #1}, {#2}\ ({#3}, {#4} {#5}), p.\ {#6--#7} {\tt [#8]}}
\newcommand\J[7]     {{\em #7}, {#1} {#2} ({#3}) {#4--#5} {{\tt [#6]}}}
\newcommand\JN[7]    {{\em #7}, {#1} ({#3}) {#4--#5} {{\tt [#6]}}}
\newcommand\JO[6]    {{\em #6}, {#1} {#2} ({#3}) {#4--#5} }
\newcommand\BOOK[4]  {{\em #1\/} ({#2}, {#3} {#4})}
\newcommand\PhD[2]   {{\em #2}, Ph.D.\ thesis #1}
\newcommand\Prep[2]  {{\em #2}, preprint {\tt #1}}
\newcommand\uPrep[2] {{\em #2}, unpublished preprint {\tt #1}}

\begin{document}

\numberwithin{equation}{section}
\author{J\"urgen Fuchs, Gregor Schaumann, and Christoph Schweigert}
\thanks{We thank Alain Brugui\`eres and Tobias Dyckerhoff for discussions and Nils Carqueville,
Julian Farnsteiner, C\'esar Galindo,
Eilind Karlsson and Vincent Koppen for helpful comments on the manuscript.
JF is supported by VR under project no.\ 2017-03836. CS is partially supported by the
RTG 1670 ``Mathematics inspired by String theory and Quantum Field Theory''
and by the Deutsche Forschungsgemeinschaft (DFG, German Research Foundation) under
Germany's Excellence Strategy - EXC 2121 ``Quantum Universe''- QT.2.
Finally we thank the editor (Joachim Kock) and the five anonymous referees
for their extremely careful handling of the submission.}
\address{Teoretisk fysik, Karlstads Universitet\\ 
Universitetsgatan 21, 651\,88 Karlstad, Sweden
\\[5pt] 
Mathematische Physik, Institut f\"ur Mathematik, Universit\"at W\"urzburg\\ 
Emil-Fischer-Stra\ss e 31, 97\,074 W\"urzburg, Germany
\\[5pt]
Bereich Algebra und Zahlentheorie, Fachbereich Mathematik, Universit\"at Hamburg\\
Bundesstra\ss e 55, 20\,146\, Hamburg, Germany
\\[5pt]
~}
\title{A modular functor from state sums for finite tensor categories and their bimodules}
\copyrightyear{2022}
\eaddress{juerfuch@kau.se\CR
gregor.schaumann@uni-wuerzburg.de\CR christoph.schweigert@uni-hamburg.de}
\maketitle
\begin{abstract}
We construct a modular functor which takes its values in the monoidal
bicategory of finite categories,
left exact functors and natural transformations. The modular functor is defined on bordisms
that are $2$-framed. Accordingly we do not need to require that the finite categories 
appearing in our construction are semisimple, nor that the finite tensor categories that 
are assigned to two-dimensional strata are endowed with a pivotal 
structure. Our prescription can be understood as a state-sum construction. The state-sum
variables are assigned to one-dimensional strata and take values in bimodule categories
over finite tensor categories, whereby we also account for the presence of boundaries and 
defects. Our construction allows us to explicitly compute functors associated to surfaces and
representations of mapping class groups acting on them.
\end{abstract}

 \tableofcontents

\section{Introduction}

Finite tensor categories are linear rigid monoidal categories obeying certain finiteness 
conditions. They arise in various contexts in representation theory, e.g.\ as categories
of finite-di\-men\-sional representations of finite-dimensional Hopf algebras, or
as representation categories of suitable vertex algebras or nets of observable algebras.
They also appear naturally in 
rigorous approaches to low-dimensional conformal and topological quantum field theories.
Fusion categories, i.e.\ finite tensor categories that are semisimple, form an important
subclass. They occur e.g.\ in the classification of subfactors and in topological field 
theory and its applications, such as the description of local invariants of 
knots and links, the study of topological phases of matter, and quantum gravity,
and as renormalization group fixed points of string net models. 

By now a comprehensive body of mathematical results -- which one may collectively refer
to as ``categorified representation theory'' -- has been built around finite tensor
categories. Remarkably, many results in categorified representation theory do not
require semisimplicity, but rather rely on the finiteness properties. Let us illustrate 
what we mean by categorified representation theory: Thinking of finite tensor categories
as a categorification of rings, it is natural to study their module and bimodule categories. 
Direct applications of these are in the description of defects and boundaries in topological
field theories. On the mathematical side, module and bimodule categories lead to a rich algebraic
structure. Specifically, invertible (not necessarily semisimple) bimodule categories give
rise to a (higher categorical
variant) of a group, the Brauer--Picard group, which plays a central role in the construction
of equivariant modular functors and thus of orbifold theories.
The bicategory of module categories leads to a bicategorical variant of Morita theory. 

An important structure in categorified representation theory is the Drinfeld center. 
For instance, Morita equivalent finite tensor categories have equivalent Drinfeld centers, 
and the Brauer--Picard group can be computed in terms of braided autoequivalences of the
Drinfeld center \Cite{Thm.\,4.1}{daNik}. 
While finite tensor categories can be endowed with interesting
additional features, e.g.\ with a pivotal structure, a braiding, or a ribbon structure, it
is worth pointing out that there is a rich theory already without assuming any such extra
features. Accordingly we take finite tensor categories without additional structure
as the starting point of the present paper. 

A comprehensive algebraic theory calls for an organizing principle. Indeed, modular functors
should provide such a principle. 
(For the precise notion of modular functor we are using in this paper, see Definition
\ref{def:modularfunctor}.) This is not surprising; many algebraic theories turn out
to have organizing principles that can be expressed in terms of geometric structure. For
instance, when working with associative algebras it is helpful to be aware of aspects of
rooted trees, while the theory of Frobenius algebras becomes very transparent in the light
of two-dimensional oriented topological field theory. In a similar vein, for categorified
representation theory it is commonly agreed that variants of extended three-dimensional 
topological field theories based on state-sum constructions should play a role. For our
purposes, we have an (extended) topological field theory or -- for more general input data
that are not required to be semisimple -- a modular functor in mind that is set up at the
level of bicategories: it assigns to one-dimensional structures categories, to two-dimensional
structures functors, and to elements of mapping class groups natural transformations. Such
topological field theories or modular functors are of intrinsic mathematical interest and 
have at the same time a lot of applications, ranging from physics to computer science. 

\medskip

On the other hand, the standard approach to state-sum constructions, as pioneered by
Turaev--Viro \cite{tuvi} and Barrett--Westbury \cite{bawe2},
is unsatisfactory when it comes to ``explaining'' categorified algebra:
\Itemize

\item
The Turaev--Viro--Barrett--Westbury construction is based on
fusion categories that are pivotal (and even spherical). In contrast, various non-trivial
aspects of categorified representation theory do not require a pivotal structure, which should
therefore better be treated as an additional feature. (This point of view is also 
advocated in \cite{doSs3}.)

\item
Turaev--Viro theory based on a fusion category $\cala$ assigns to a circle the Drinfeld center 
$\calz(\cala)$. The fusion category $\cala$ itself, on the other hand, is
effectively invisible, in the sense that Morita equivalent spherical fusion categories
give the same extended topological field theory at the bicategorical level.

\item
As already pointed out, much of categorified representation theory works beyond the realm
of fusion categories, for the larger class of finite tensor categories which enjoy 
analogous finiteness properties as fusion categories, but are not necessarily semisimple.
\end{itemize}

The central goal we achieve in this paper is a geometric framework that governs
the categorified representation theory of finite tensor categories
and their finite (bi)module categories. To overcome the shortcomings of
conventional state-sum constructions, we work in the following setting:
\Itemize

\item
To allow for finite tensor categories that are not semisimple, we construct 
a modular functor rather than a 3-2-1-extended topological field theory. Specifically,
we do not formulate a theory for arbitrary three-manifolds with corners, but
restrict ourselves to surfaces and to actions of their mapping class groups.
The idea that non-semisimple categories only allow one to deal with a restricted class of
three-manifolds (or with three-manifolds having additional structure, see e.g.\
\Cite{Def.\,3.3}{bcgp}) -- in our case, at least cylinders twisted by
the action of an element of the mapping class group -- is not new, see e.g.\ \cite{doSs3}.

\item
One way to make the finite tensor category itself, rather than merely its Drinfeld center,
visible in the construction, is to consider the extension of the theory to the point
\cite{doSs3}. In our construction we
expose the finite tensor category, as well as its module and bimodule categories, 
by instead extending the category of cobordisms to include boundaries and defects. This 
modification is not new either. Indeed, it is known \cite{fusV} that boundary conditions 
for a topological field theory of Reshetik\-hin--Tu\-ra\-ev type based on a modular tensor 
category $\calc$ correspond to Witt-trivializations, i.e.\ to braided equivalences 
$\calc \,{\simeq}\, \calz(\cala)$ to the Drinfeld center of some fusion category $\cala$, 
which, in turn, has a direct interpretation as a category of
Wilson lines associated with a specific boundary condition.
\\
That we include defects has a further benefit: There are defects between any two topological
field theories of Turaev--Viro type, and hence our construction encompasses in a single theory
the Turaev--Viro theories for all choices of fusion categories.

\item
Finally, in order to do without a pivotal structure on the algebraic side, we supplement
structure on the geometric side and work with $2$-framed manifolds rather than with oriented
ones.\,%
  \footnote{~In Remark \ref{rem:TV} we sketch how our construction should be adapted
  so as to apply to a monoidal bicategory of oriented cobordisms.}
Again, this approach has been advocated before, see \cite{doSs3} as well as \cite{kupe3}.
In \cite{doSs3} the existence of a framed theory is shown with the help of the 
cobordism hypothesis. This approach is non-constructive and thus difficult to compare 
with other approaches. Our approach is in the framed setting, too, but it is 
constructive. We expect it to be the 2-1 framed defect theory with mapping class group 
actions (sometimes called a $(2{+}\epsilon$)-theory) that extends to a framed fully
extended  defect 2-d theory.
\end{itemize}

Several frameworks for addressing our goal may come to mind, such as factorization algebras, 
Kita\-ev-type state-sum models, or constructions based on fully extended topological field 
theories that invoke the cobordism hypothesis. The approach taken in the present paper
provides an explicit state-sum construction in a purely categorical setting. This avoids
the introduction of extra structure, and it does not invoke the cobordism hypothesis,
but has structural similarities with constructions familiar from factorization algebras.
Indeed one might even hope that our construction, for instance
Proposition \ref{Proposition:main-disk-repl}, will allow one to compare state-sum
models, factorization algebras and the fully local theory in our concrete setting. 

Our construction
is tailored to the specific target bicategory of finite tensor categories with
the Deligne product as the monoidal structure and with left exact (or, alternatively, 
right exact) functors as 1-morphisms, and uses the full power of that structure.
Defects are built in from the start, as the carriers of the state-sum variables.
Being very concrete, our approach leads directly to fully explicit
computational prescriptions for specific situations of interest. Our findings are in line with
the results, conjectures and expectations in other approaches,
albeit the direct comparison between different frameworks is far from straightforward.

\medskip

Even given the clear program outlined above, the right definitions and a full construction 
still turn out to be subtle. Accordingly, in a sense, our first important insight is
Definition \ref{def:Bordcats}: it specifies a monoidal bicategory $\Bfrdefdec$ of $2$-framed 
defect cobordisms that suits our purposes. The description of a modular functor 
in Definition \ref{def:modularfunctor} is then standard, and it follows 
from our general goal that we aim for a modular functor with values in the monoidal bicategory
$\mathcal S \,{=}\, \FUNle$, having finite categories as objects, left exact functors as 
1-morphisms, and the Deligne product as the monoidal structure.
(There is also a variant that instead uses as target the bicategory of
finite categories with right exact functors.)
Theorem \ref{thm:main}, the main result of this paper, then asserts 
the existence of such a modular functor, i.e.\ of a symmetric monoidal bifunctor
  \be
  \ZZ:\quad \Bfrdefdec \xrightarrow{~~} \,\mathcal S \,.
  \ee

The unsuspicious words ``symmetric monoidal bifunctor'' specify quite a lot of structure and 
properties. (This partly explains the length of the present paper.) It includes, 
for instance -- in the form of the horizontal composition of 1-morphisms -- factorization of 
the modular functor under the gluing of surfaces, see Theorem \ref{thm:Facto}.

\medskip

Let us now summarize the main line of our arguments:
\Itemize

\item
We have to assign a finite linear category to each object of $\Bfrdefdec$, i.e.\ to
certain one-manifolds with additional structure. These finite categories are
suitable generalizations of Drinfeld centers. This is the topic of Section 
\ref{sec:bdrycat}.

\item
Next we must assign a left exact functor to each 1-morphism of $\Bfrdefdec$, i.e.\ to
bordisms with extra structure, which we call \emph{\defectsurface s}. This
is achieved in the form of a state-sum construction and follows the standard
three-step pattern of such constructions: For a surface with boundary, one first 
constructs a ``big'' vector space -- actually, a linear functor.
We call this functor the \emph{pre-block functor}.
 \\[.2em]
That we work with categories that are not necessarily semisimple forces us to work 
systematically with natural notions from category theory. Specifically, to construct left 
exact functors, we use Hom functors and implement the sum over states by taking coends.
Indeed we would raise the claim that the systematic use of category-theoretic concepts allows
for a substantial conceptual clarification, even when dealing with semisimple categories. 
 \\[.2em]
The construction of pre-block functors occupies the first part of Section \ref{sec:preblock}.

\item
The second step in a state-sum construction consists in imposing an appropriate flatness
condition. It is one of the novel insights of this paper that when making use of the 
$2$-framing on the surfaces, one can enforce flatness of holonomies without assuming the 
existence of a pivotal structure on the finite tensor categories. To be able to impose 
flat holonomy, the defect network of the surface must be such that each of its \twocells\
has the topology of a disk; we call a surface of this type a \emph{fine} \defectsurface.
The solution of the flatness condition on the pre-block functor for a fine surface $\Sigma$
gives another left exact functor, which assigns to $\Sigma$ in a functorial way subspaces
of the big vector spaces (which one might think of as `spaces of ground states').
Constructing these functors for all fine surfaces is the second main subject of Section
\ref{sec:preblock}.

\item 
A modular functor must, of course, assign a functor to \emph{any} defect surface, not just
to fine ones. In Section \ref{sec:5} we explain how to define such functors, which we call
\emph{block functors}, to surfaces with a defect network that is not necessarily fine. To
this end we introduce the notion of a \emph{refinement} of a \defectsurface. We show that
refinements to fine surfaces exist, and then proceed to set up, first for disks, a system
of isomorphisms for functors associated to fine refinements in such a manner that we can
define the block functors for disks as limits. The block functors for general
\defectsurface s are then constructed from the block functors for disks. We also show that
the so defined block functors obey factorization and study actions of mapping class groups.

\item
When combined, these results imply our main Theorem \ref{thm:main}. A concise summary
of the proof is given in Section \ref{prf:Proofmain}.

\end{itemize} 

Instead of taking finite tensor categories with left exact functors as the target bicategory
of the modular functor, we could have chosen the bicategory of finite tensor categories with
\emph{right} exact functors. Each of these two bicategories is monoidal, with the Deligne
product as a symmetric monoidal structure. A duality between the left and right exact functors
is provided by the Eilenberg--Watts functors that were studied in \cite{shimi7,fScS2}. According
to this duality, the left exact Hom functor gets replaced by the vector space dual of the Hom
functor, which is right exact, and coends in the state-sum construction must be replaced by
ends. Beyond this aspect, the Eilenberg--Watts calculus also plays a significant role in our
approach and in interpreting our results. For instance, it makes it easy to describe how the 
modular functor provides, via the fusion of boundary insertions, a composition on Deligne 
products $\cc\calm\boti\caln$ (see Proposition \ref{prop:preblock-compose}); for other 
applications, see e.g.\ Example \ref{ex:4.4-part1} or Corollary \ref{cor:modulenattrafo}.

Several complementary results help to make block functors computable. For instance,
in Theorem \ref{thm:fusion-rel-Del} we show that the fusion of defect lines corresponds
to a relative Deligne product of bimodule categories (which depends on the framings involved),
while Proposition \ref{prop:block-compose} tells us that a pair of gluing boundaries can be 
combined to a single one, in a way that is described by the composition of functors (see the
picture \eqref{eq:PsiG.PsiF=PsiGF} below), without changing the block functor. 
As a consequence, a specific simple \defectsurface, which we call the
`\basicdisk' (displayed in picture \eqref{pic:ex:4.4}), is of particular importance;
we show that its pre-block spaces consist of natural transformations
(see formula \eqref{eq:preblock442}) and that its block spaces are the corresponding
module natural transformations (Corollary \ref{cor:modulenattrafo}).
The occurrence of module natural transformations is the simplest instance in which our
construction contributes to the program of geometrically realizing categorified
representation theory. Also, in view of the construction of block functors for general
\defectsurface s from those for disks, one may think of block functors as giving spaces
that constitute a huge generalization of spaces of module natural transformations. 

We also provide details for a few situations of specific interest, a sample being the following:
The functors of braided induction (or $\alpha$-induction) appear naturally in the block
functor for disks with a free boundary (Example \ref{exa:braidedinduction});
the `transmission functor' which was considered in \Cite{Sect.\,5.1}{enoM} is obtained as the
block functor for a cylinder with a circumferential defect line (Example \ref{exa:ENOM});
a variant of the twist (involving the double-dual functor) on objects of a braided
monoidal category appears in the natural transformation that is obtained from a Dehn twist 
on the cylinder over a circle (Proposition \ref{Proposition:dehntwist});
the block functor for a two-sphere with one circular defect line and any number of 
insertions (Example \ref{ex:3holedsphere}) has bimodule natural transformations as its values;
and a pair of pants (Example \ref{ex:S1S2S3}) realizes, via the Eilenberg--Watts equivalences,
the composition of bimodule functors as well as, in the situation that all defects involved
are transparent, the tensor product in the Drinfeld center (which is the functor that also
the standard Turaev--Viro construction assigns to a pair of pants).
As examples for the actions of mapping class groups, we consider a braiding move on
a three-punctured sphere (Proposition \ref{prop:braid}) and a Dehn twist
(Proposition \ref{Proposition:dehntwist}).


\section{Framed defect manifolds} \label{sec:frdbord}

In this section we define a bicategory $\Bfrdefdec$ of two-dimensional
$2$-framed bordisms with labeled defects. To this end we introduce in a first step
a geometric bicategory $\Bfrdef$ having unlabeled defects.
All manifolds considered below are assumed to be smooth and oriented.

\subsection{Framed defect bordisms} \label{sec:Bfrdef}

We consider manifolds which can have boundaries and corners and can contain defect lines,
and which in addition are endowed with a framing. Before describing the bicategory $\Bfrdef$
of all such manifolds, we first consider a sub-bicategory $\Bfrdefcl$ of manifolds without
corners. We start by giving representatives for the morphisms of this sub-bicategory.

Denote by $I \,{=}\, [0,1]$ the standard interval and by $\SS^1$ the
standard circle; both endowed with their standard orientation. We write
  \be
  \bearl
  I^{\sqcup n} := I\,{\sqcup}\, I\,{\sqcup}\,\ldots\,{\sqcup}\, I \qquad {\rm  and} 
  \Nxl7
  (\SS^1)^ {\sqcup n} := \SS^1\,{\sqcup}\, \SS^1\,{\sqcup}\,\ldots\,{\sqcup}\, \SS^1
  \eear
  \ee
for the corresponding finite disjoint unions, with $n \,{\in}\, \mathbb Z_{\ge 0}$ and 
with $I^{\sqcup 0}$ and $(\SS^1)^{\sqcup 0}$ being the empty set.
Let $\vSigma$ be a compact oriented surface, possibly with boundary. We endow $\vSigma$ with
further structure that accounts for defect lines and a compatible framing. To this end we 
first introduce the additional datum of an embedding 
$\delta\colon I^{\sqcup n} \sqcup (\SS^1)^{\sqcup m}\,{\to}\, \vSigma$ 
for some $m,n \,{\in}\, \mathbb Z_{\ge 0}$ that is subject to the following restrictions:
We require that the end points of each interval are mapped to
the boundary, i.e.\ $\delta(\{0,1\}^{\sqcup n}) \,{\subset}\, \partial\vSigma$; that all other
points of the image of $\delta$ lie in the interior of $\vSigma$; and that each connected
component of $\partial\vSigma$ must contain at least one end point of one of the intervals.\,%
 \footnote{~The reason for imposing this requirement will become clear once we decorate
 the intervals with algebraic data.}

We call the image of $\delta$ in $\vSigma$
and, by abuse of language, also the map $\delta$ itself, the (set of)
\emph{unlabeled defect lines} of $\vSigma$. Note that each of the defect lines inherits an
orientation from the standard orientation of the interval $I$ or the circle $\SS^1$,
respectively. As an illustration, the following picture shows a situation in which the
underlying surface $\vSigma$ is a sphere with three holes and in which $\vSigma$
contains three unlabeled defect lines:
  \def \locpa  {1.8}  
  \be
  \raisebox{-5.2em}
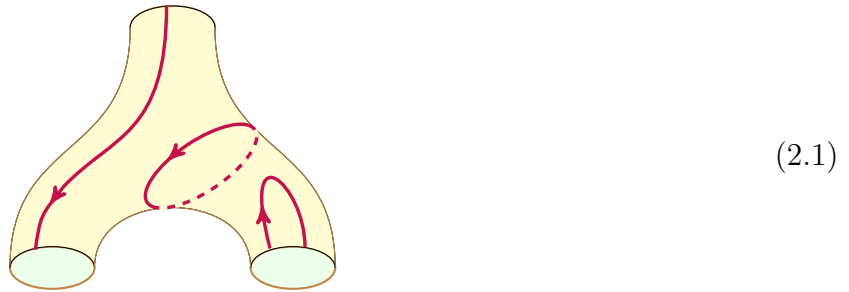
{\begin{tikzpicture}[scale=\locpa,tqft/cobordism/.style={draw}] 
  \begin{scope} [tqft/every boundary component/.style={fill=\fillBdyCircle,thick,
                draw=\colorCircle,scale=\locpa},
                tqft/cobordism edge/.style={draw,\colorCircle}]
  \pic[tqft/pair of pants,fill=\fillCobord,name=Pants,at={(0,0)},scale=\locpa];
  \scopeArrow{0.83}{\arrowDefect}
  \draw[line width=\widthDefect,color=\colorDefect,postaction={decorate}] (-0.05,0.18)
       .. controls (-0.05,-0.9) and (-0.5,-0.9) .. (-0.9,-1.32) 
       .. controls (-1.1,-1.57) .. (-1.14,-1.85);
  \end{scope}
  \scopeArrow{0.27}{\arrowDefect}
  \draw[line width=\widthDefect,color=\colorDefect,postaction={decorate}] (.81,-1.84)
       .. controls (.55,-0.84) and (1.1,-1.3) .. (1.1,-1.83);
  \end{scope}
  \scopeArrow{0.33}{\arrowDefect}
  \draw[line width=\widthDefect,color=\colorDefect,postaction={decorate},rotate=34]
       (-.45,-1.09) ellipse (.54 and .21);
  \draw[dashed,line width=1.2*\widthDefect,color=\fillCobord,rotate=34]
       ( .09,-1.09) arc (360:200: .54 and .21);
  \end{scope}
  \end{scope}
  \end{tikzpicture}}
  \ee
~

We allow only for pairs $(\vSigma,\delta)$ of surfaces with defect lines that 
can be endowed with the additional structure of a 2-\emph{framing}, that is, with a 
non-vanishing vector field $\chi$ on $\vSigma$ that along each defect line is parallel to it 
and whose direction matches the orientation of the defect line. Together with the 
orientation of $\vSigma$, the vector field $\chi$ determines a trivialization of the
tangent bundle of $\vSigma$, unique up to homotopy, hence the term $2$-framing for $\chi$. 

\begin{Definition}\label{def:undefsurf0}\label{def:fine}
\Itemizeiii
\item[{\rm (i)}]
An \emph{\undefectsurface\ without corners} is a triple $(\vSigma,\delta,\chi)$ consisting of a 
compact oriented surface $\vSigma$ without corners, and
with unlabeled defect lines $\delta$ and a $2$-framing $\chi$ on $(\vSigma,\delta)$.

\item[{\rm (ii)}]
An \undefectsurface\ $(\vSigma,\delta,\chi)$ is called \emph{fine} iff each connected component 
of $\vSigma \,{\setminus} \{ \delta \}$ is topologically a disk.

\item[{\rm (iii)}]
A fine \undefectsurface\ $(\vSigma,\delta,\chi)$ is called \emph{\gfine} iff
the boundary of every connected component of $\vSigma \,{\setminus} \{ \delta \}$
contains at most one connected component of the boundary of $\vSigma$.

\end{itemize}
\end{Definition}

\begin{Remarks}\label{rem:...iii:gfine}
\Itemizeiii
\item[(i)]
That we here augment the terminology by the qualification \emph{unlabeled} is due to the fact
that we will want to be able to distinguish between different types of defect lines
and accordingly will decorate them, in Section \ref{sec:labeling}, with suitable labels. 

\item[(ii)]
We do not allow for junctions of defect lines. Instead, any putative point in $\vSigma$ at which
defect lines would meet is realized as a boundary circle on which those defect lines end. 
Hereby we avoid the use of stratified manifolds which e.g.\ appear in the approach of
\cite{camSc,carS}.
(The size of such a boundary circle is immaterial, though. If we imagine to shrink it to a 
point, we could assign the category that any given modular functor associates to the
boundary circle instead also to a vertex at which the defect lines meet,
compare Section 2.4 of \cite{carS}.)

\item[(iii)]
The construction presented in this paper is a state-sum construction. We exhibit
(in Sections \ref{sec:bdrycat} and \ref{sec:preblock}) what it assigns to decorated 
one- and two-manifolds. State-sum constructions are generally believed to extend down
to the point, to which in our case they should assign an object in a symmetric 
monoidal tricategory. Now obviously, a point will have a decoration as well, namely
a finite tensor category $\cala$. It is thus natural to expect that the relevant 
tricategory is the tricategory of finite tensor categories and their bimodule
categories. This tricategory is also the one used in \cite{doSs3}, see Remark 
\ref{rmk:mainthm}(i). In the present paper we refrain from introducing tricategorical 
concepts and therefore do not investigate this issue any further.

\item[(iv)]
The \undefectsurface\ resulting from the gluing of two fine \undefectsurface s is not
necessarily fine. In contrast, the \undefectsurface\ resulting from the gluing of two 
\ggfine\ \undefectsurface s is again \gfine.
\end{itemize}
\end{Remarks}

As an illustration, the following picture indicates a framing for a closed defect surface
whose underlying surface is an annulus and which has two defect lines (here and below, the
surface is embedded in the plane and is taken to inherit the standard orientation of the plane):
    \def\locpa  {0.81}   
    \def\locpA  {2.5}    
    \def\locpB  {1.65}   
    \def\locpb  {0.62}   
    \def\locpc  {2.2pt}  
    \def\locpd  {0.12}   
  \be
  \raisebox{-5.9em}
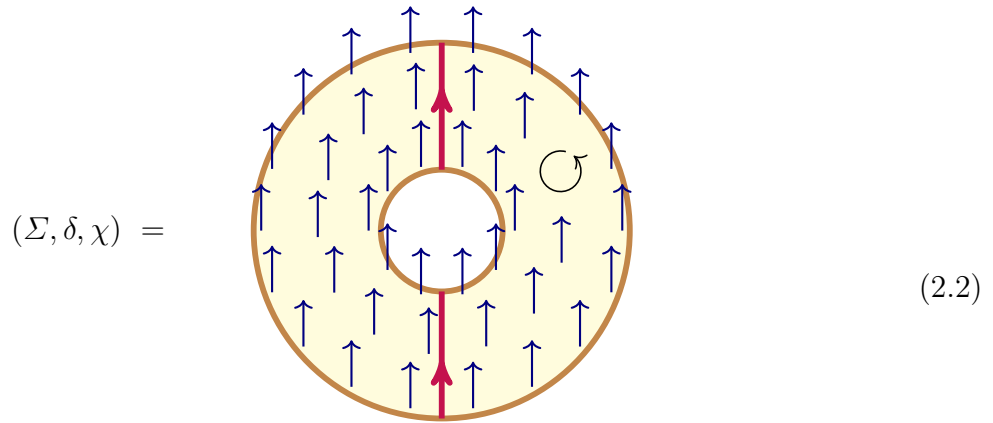
{\begin{tikzpicture}
  \scopeArrow{0.158}{} \scopeArrow{0.65}{}
  \filldraw[fill=\colorSigma,line width=\locpc,draw=\colorCircle,postaction={decorate}]
       node[left=3.5cm,color=black] {$(\vSigma,\delta,\chi) ~=$} (0,0) circle (\locpA);
  \end{scope} \end{scope}
  \scopeArrow{0.18}{} \scopeArrow{0.65}{}
  \filldraw[fill=white,line width=\locpc,draw=\colorCircle,postaction={decorate}]
       (0,0) circle (\locpa);
  \end{scope} \end{scope}
  \begin{scope}[shift={(\locpB,\locpa+0.24)}] \drawOrientationLarge; \end{scope}
  \scopeArrow{0.46}{\arrowDefect}
  \draw[line width=\locpc,color=\colorDefect,postaction={decorate}] (0,-\locpA) -- (0,-\locpa);
  \end{scope}
  \scopeArrow{0.65}{\arrowDefect}
  \draw[line width=\locpc,color=\colorDefect,postaction={decorate}] (0,\locpa) -- (0,\locpA);
  \end{scope}
  \draw[thick,color=\colorFrame,->] (0:1.2*\locpa) -- +(0,\locpb);
  \draw[thick,color=\colorFrame,->] (180:1.2*\locpa) -- +(0,\locpb);
  \foreach \xx in {36,72,...,144} \draw[thick,color=\colorFrame,->]
       (\xx:1.1*\locpa) -- +(0,\locpb);
  \foreach \xx in {216,252,...,324} \draw[thick,color=\colorFrame,->]
       (\xx:1.1*\locpa) -- +(0,\locpb);
  \foreach \xx in {48,75,...,345} \draw[thick,color=\colorFrame,->]
       (\xx:\locpB) -- +(0,\locpb);
  \foreach \xx in {0,20,...,340} \draw[thick,color=\colorFrame,->]
       (\xx:0.96*\locpA) -- +(0,\locpb);
  \end{tikzpicture}}~~
  \label{eq:picture-closed-annulus}
  \ee ~

To introduce the \emph{objects} of the category $\Bfrdefcl$, we consider the restriction of
the so defined structure to the boundary circles (together with little collars around them). 
A point on the boundary $\partial\vSigma$ is said to be
\emph{marked} iff it is in the image of $\delta$, i.e.\ is the end point of a defect line. 
We label a marked point by $+1$ and call it a \emph{positive point} if it is the image of 
an initial point $0\,{\in}\, I$  of a defect line, and label it by $-1$ and call it a
\emph{negative point} if it is the image of an end point $1\,{\in}\, I$. We also call a 
marked point on $\partial\vSigma$ -- or, more generally, on a one-manifold -- together
with a sign $\pm1$ an (unlabeled) \emph{defect point}, and the closure of the 
interval along a boundary circle between two neighboring defect points a \emph{segment} $s$.
To determine the structure induced on the boundary $\partial \vSigma$ by the $2$-framing of an
\undefectsurface, we make use of the fact that in order to glue bordisms of smooth manifolds
along boundary circles, the circles need to be endowed with collars. 
Concretely, a connected component of $\partial\vSigma$ is to be
considered with (the germ of) an embedding $\SS^1 \,{\times}\, I \,{\to}\, \vSigma$. Thus
we obtain a non-vanishing vector field on $\SS^1 \,{\times}\, I$ as the pullback of 
a $2$-framing $\chi$ on $(\vSigma,\delta)$. Note that at each defect point on
$\partial\vSigma$ the framing vector field $\chi$ on $\vSigma$
provides a non-vanishing vector; by our requirement that near any defect line $\delta$ the 
vector field $\chi$ is parallel to $\delta$ with matching direction, the so obtained vector
at a defect point $p \,{\in}\, \delta \,{\cap}\,\partial\vSigma$ is outward-pointing iff
$\delta$ is oriented towards $\partial\vSigma$, i.e.\ iff $p$ is a negative point.

Further, it is natural to allow also for one-manifolds with boundary, and as a consequence
admit corresponding \defectsurface s (see below) as well, which then have (one-dimensional)
boundaries and corners. Denoting by $\underline{\R}$ the trivial 1-dimensional vector
bundle with oriented fiber over a given base, we then arrive at

\begin{Definition}
\Itemizeiii

\item[{\rm (i)}]
An \emph{\undefectonemanifold} is a triple $\LL \eq (\L,\epsilon,\chi)$ consisting of
an oriented compact one-manifold $\L$, possibly with boundary, a finite set
$\epsilon \,{\supseteq}\, \partial\L$ of defect points, and a non-vanishing 
vector field $\chi \,{\in}\, \Gamma(T \L{\oplus}\underline{\R})$.
\\
Further, at each defect point $p \,{\in}\, \epsilon$ the component of $\chi(p)$ in $T \L$ is
required to vanish, and the component of $\chi(p)$ in $\underline\R$ must be positive iff 
$p$ is a positive point.

\item[{\rm (ii)}]
A \emph{closed \undefectonemanifold} is an \undefectonemanifold\ with empty boundary.

\item[{\rm (iii)}]
A morphism of \undefectonemanifold s is a diffeomorphism of manifolds that preserves the
non-vanishing vector field. 
\end{itemize}
\end{Definition}

Let us illustrate a few typical situations of \undefectonemanifold s by pictures.
In such pictures we use the following conventions. 
We draw a one-manifold $\L$ as embedded in the paper plane $\R^2$. For the tangent 
space $T_p\L$ at any point $p\,{\in}\,\L$ we adopt the standard convention to depict it as
the tangential affine line in $\R^2$.

To graphically represent a section in the bundle $T\L{\oplus}\underline\R$ we must in
addition specify the direction of $\underline\R$; it suffices to do this for an interval
and for a circle. For an \emph{interval} embedded horizontally in the plane and oriented
from left to right, we take the positive direction of $\underline\R$
to point upwards. Thus at any positive marked point on the interval the vector fields 
we consider point upwards, and at any negative marked point they point downwards;
in particular, between a positive and a negative point the vector field is tangential in 
at least one point. The following picture shows examples of unlabeled defect intervals:
    \def\locpa  {3.0}    
    \def\locpb  {0.8}    
    \def\locpc  {2.2pt}  
    \def\locpd  {0.02}   
    \def\locpe  {4.9}    
    \def\locpf  {1.6}    
    \def\locpF  {0.75}   
    \def\locpx  {0.08}   
  \be
  \raisebox{-4.5em}{\begin{tikzpicture}
  \draw[line width=\locpc,color=\colorCircle] (-\locpd,0) -- (\locpa+\locpd,0);
  \foreach \xx in {0,1,...,10}
     \draw[thick,color=\colorFrame,->] (0.1*\locpa*\xx cm,0) -- (0.1*\locpa*\xx cm,\locpb);
  \draw[line width=\locpc,color=\colorCircle,->] (0.27*\locpa,0) -- (0.29*\locpa,0);
  \filldraw[color=\colorDefect] (0,0) circle (\locpx) node[below] {$\scriptstyle+$}
                                (\locpa,0) circle (\locpx) node[below] {$\scriptstyle+$};
  \begin{scope}[shift={(\locpe,\locpF)}]
  \draw[line width=\locpc,color=\colorCircle] (-\locpd,0) -- (\locpa+\locpd,0);
  \foreach \xx in {0,1,...,10}
     \draw[thick,color=\colorFrame,->] (0.1*\locpa*\xx cm,0) -- (0.1*\locpa*\xx cm,-\locpb);
  \draw[line width=\locpc,color=\colorCircle,->] (0.27*\locpa,0) -- (0.29*\locpa,0);
  \filldraw[color=\colorDefect] (0,0) circle (\locpx) node[above=-2pt] {$\scriptstyle-$}
                                (\locpa,0) circle (\locpx) node[above=-2pt] {$\scriptstyle-$};
  \end{scope}
  \begin{scope}[shift={(0,-\locpf)}]
  \draw[line width=\locpc,color=\colorCircle,->] (0.83*\locpa,0) -- (0.85*\locpa,0);
  \draw[line width=\locpc,color=\colorCircle] (-\locpd,0) -- (\locpa+\locpd,0);
  \foreach \xx in {0,1,...,4} \draw[thick,color=\colorFrame,->]
             (0.06*\locpa*\xx cm,0) -- +(90-18*\xx:\locpb);
  \draw[thick,color=\colorFrame,->] (0.3*\locpa,0) -- +(7:\locpb);
  \draw[thick,color=\colorFrame,->] (0.42*\locpa,0) -- +(-24:\locpb);
  \foreach \xx in {7,8,9,10} \draw[thick,color=\colorFrame,->]
             (-0.4*\locpa+0.14*\locpa*\xx,0) -- +(90-18*\xx:\locpb);
  \filldraw[color=\colorDefect] (0,0) circle (\locpx) node[below] {$\scriptstyle+$}
                                (\locpa,0) circle (\locpx) node[above=-2pt] {$\scriptstyle-$};
  \end{scope}
  \begin{scope}[shift={(\locpe,-\locpf)}]
  \draw[line width=\locpc,color=\colorCircle,->] (0.83*\locpa,0) -- (0.85*\locpa,0);
  \draw[line width=\locpc,color=\colorCircle] (-\locpd,0) -- (\locpa+\locpd,0);
  \foreach \xx in {0,1,...,4} \draw[thick,color=\colorFrame,->]
             (0.06*\locpa*\xx cm,0) -- +(-90+18*\xx:\locpb);
  \draw[thick,color=\colorFrame,->] (0.3*\locpa,0) -- +(-7:\locpb);
  \draw[thick,color=\colorFrame,->] (0.42*\locpa,0) -- +(24:\locpb);
  \foreach \xx in {7,8,9,10} \draw[thick,color=\colorFrame,->]
             (-0.4*\locpa+0.14*\locpa*\xx,0) -- +(-90+18*\xx:\locpb);
  \filldraw[color=\colorDefect] (0,0) circle (\locpx) node[above=-2pt]{$\scriptstyle-$}
                                (\locpa,0) circle (\locpx) node[below]{$\scriptstyle+$};
  \end{scope}
  \end{tikzpicture}}
  \label{eq:exa-intervals}
  \ee
~\\
For a \emph{circle} embedded in $\R^2$ we fix conventions by requiring that
the trivial bundle $\underline\R$ is outward-pointing. Then the vector fields of our interest
point outwards at any positive marked point and inwards at any negative one. Again, between a
positive and a negative point the vector field has to be tangential in at least one point. 
Here are simple examples of $2$-framed circles:
    \def\locpa  {1.45}   
    \def\locpb  {0.71}   
    \def\locpc  {2.2pt}  
    \def\locpd  {0.12}   
    \def\locpe  {3.9}    
  \begin{eqnarray}
  \begin{tikzpicture}
  \begin{scope}[shift={(-\locpb,0)}]
  \scopeArrow{0.15}{<} \scopeArrow{0.65}{<}
  \draw[line width=\locpc,color=\colorCircle,postaction={decorate}] (0,0) circle (\locpa);
  \end{scope} \end{scope}
  \filldraw[color=\colorDefect] (0,\locpa)  circle (\locpd) node[below=1pt] {$+$}
                                (0,-\locpa) circle (\locpd) node[above=1pt] {$+$};
  \foreach \xx in {0,20,...,360}
           \draw[thick,color=\colorFrame,->] (\xx:\locpa) -- (\xx:\locpa+\locpb);
  \end{scope}
  \begin{scope}[shift={(\locpe,0)}]
  \scopeArrow{0.11}{<} \scopeArrow{0.64}{<}
  \draw[line width=\locpc,color=\colorCircle,postaction={decorate}] (0,0) circle (\locpa);
  \end{scope} \end{scope}
  \filldraw[color=\colorDefect] (0,\locpa)  circle (\locpd) node[above=-1pt] {$-$}
                                (0,-\locpa) circle (\locpd) node[below=-1pt] {$-$};
  \foreach \xx in {0,24,...,360}
           \draw[thick,color=\colorFrame,->] (\xx:\locpa) -- (\xx:\locpa-\locpb);
  \end{scope}
  \begin{scope}[shift={(2*\locpe,0)}]
  \scopeArrow{0.16}{<} \scopeArrow{0.64}{<}
  \draw[line width=\locpc,color=\colorCircle,postaction={decorate}] (0,0) circle (\locpa);
  \end{scope} \end{scope}
  \filldraw[color=\colorDefect] (0,\locpa)  circle (\locpd) node[below=1pt] {$+$}
                                (0,-\locpa) circle (\locpd) node[below=-1pt] {$-$};
  \foreach \xx in {0,20,...,360}
           \draw[thick,color=\colorFrame,->] (\xx:\locpa) -- +(0,\locpb);
  \end{scope}
  \begin{scope}[shift={(3*\locpe,0)}]
  \scopeArrow{0.143}{<} \scopeArrow{0.656}{<}
  \draw[line width=\locpc,color=\colorCircle,postaction={decorate}] (0,0) circle (\locpa);
  \end{scope} \end{scope}
  \filldraw[color=\colorDefect] (0,\locpa)  circle (\locpd) node[above=-1pt] {$-$}
                                (0,-\locpa) circle (\locpd) node[above=1pt] {$+$};
  \foreach \xx in {0,20,...,360}
           \draw[thick,color=\colorFrame,->] (\xx:\locpa) -- +(0,-\locpb);
  \end{scope}
  \end{tikzpicture}
  \nonumber
  \end{eqnarray}
  ~\\[-3.2em]
  \be
  \phantom{\int}
  \label{eq:exa-circles}
  \ee

Up to homotopy keeping the vector field $\chi$ transversal at the defect points, $\chi$ is
determined by its winding between neighboring defect points. When counted in units of $\pi$
in the direction of the orientation, the winding is an integer. We call this integer the
\emph{framing index} of the segment $s$ and denote it by
  \be
  \ind_\chi(s) \in \Z \,.
  \ee
The index is an even integer if the neighboring defect points have the same sign, and odd 
otherwise. 
As an illustration, the values of the index for the intervals in \eqref{eq:exa-intervals} are
    \def\locpa  {3.0}    
    \def\locpA  {1.45}   
    \def\locpb  {0.8}    
    \def\locpc  {2.2pt}  
    \def\locpd  {0.02}   
    \def\locpe  {7.2}    
    \def\locpE  {6.38}   
    \def\locpf  {2.1}    
    \def\locpg  {1.51}   
    \def\locpx  {0.09}   
    \def\locpX  {0.11}   
  \be
  \raisebox{-4.9em}{\begin{tikzpicture}[scale=1.1]
  \foreach \xx in {0,1,...,10}
     \draw[thick,color=\colorFrame,->] (0.1*\locpa*\xx cm,0) -- (0.1*\locpa*\xx cm,\locpb);
  \draw[line width=\locpc,color=\colorCircle] 
     node[anchor=east,color=black] {$\ind_\chi\Big($} (-\locpd,0) -- (\locpa+\locpd,0)
     node[right=2pt,color=black] {$\Big) \,=\, 0 \,= $};
  \draw[line width=\locpc,color=\colorCircle,->] (0.27*\locpa,0) -- (0.29*\locpa,0);
  \filldraw[color=\colorDefect] (0,0) circle (\locpx) node[below] {$\scriptstyle+$}
                                (\locpa,0) circle (\locpx) node[below] {$\scriptstyle+$};
  \begin{scope}[shift={(\locpE,0)}]
  \foreach \xx in {0,1,...,10}
     \draw[thick,color=\colorFrame,->] (0.1*\locpa*\xx cm,0) -- (0.1*\locpa*\xx cm,-\locpb);
  \draw[line width=\locpc,color=\colorCircle] 
     node[anchor=east,color=black] {$\ind_\chi\Big($} (-\locpd,0) -- (\locpa+\locpd,0)
     node[anchor=west,color=black] {$\Big) \,,$};
  \draw[line width=\locpc,color=\colorCircle,->] (0.27*\locpa,0) -- (0.29*\locpa,0);
  \filldraw[color=\colorDefect] (0,0) circle (\locpx) node[above=-2pt] {$\scriptstyle-$}
                                (\locpa,0) circle (\locpx) node[above=-2pt] {$\scriptstyle-$};
  \end{scope}
  \begin{scope}[shift={(0,-\locpf)}]
  \draw[line width=\locpc,color=\colorCircle]
             node[left,color=black] {$\ind_\chi\Big($} (-\locpd,0) -- (\locpa+\locpd,0)
             node[right,color=black] {$\Big) \,=\, {-}1 \,,$};
  \foreach \xx in {0,1,...,4} \draw[thick,color=\colorFrame,->]
             (0.06*\locpa*\xx cm,0) -- +(90-18*\xx:\locpb);
  \draw[thick,color=\colorFrame,->] (0.3*\locpa,0) -- +(7:\locpb);
  \draw[thick,color=\colorFrame,->] (0.42*\locpa,0) -- +(-24:\locpb);
  \foreach \xx in {7,8,9,10} \draw[thick,color=\colorFrame,->]
             (-0.4*\locpa+0.14*\locpa*\xx,0) -- +(90-18*\xx:\locpb);
  \draw[line width=\locpc,color=\colorCircle,->] (0.83*\locpa,0) -- (0.85*\locpa,0);
  \filldraw[color=\colorDefect] (0,0) circle (\locpx) node[below] {$\scriptstyle+$}
                                (\locpa,0) circle (\locpx) node[above=-2pt] {$\scriptstyle-$};
  \end{scope}
  \begin{scope}[shift={(\locpe,-\locpf)}]
  \draw[line width=\locpc,color=\colorCircle] 
             node[left,color=black] {$\ind_\chi\Big($} (-\locpd,0) -- (\locpa+\locpd,0)
             node[right,color=black] {$\Big) \,=\, 1 \,,$};
  \foreach \xx in {0,1,...,4} \draw[thick,color=\colorFrame,->]
             (0.06*\locpa*\xx cm,0) -- +(-90+18*\xx:\locpb);
  \draw[thick,color=\colorFrame,->] (0.3*\locpa,0) -- +(-7:\locpb);
  \draw[thick,color=\colorFrame,->] (0.42*\locpa,0) -- +(24:\locpb);
  \foreach \xx in {7,8,9,10} \draw[thick,color=\colorFrame,->]
             (-0.4*\locpa+0.14*\locpa*\xx,0) -- +(-90+18*\xx:\locpb);
  \draw[line width=\locpc,color=\colorCircle,->] (0.83*\locpa,0) -- (0.85*\locpa,0);
  \filldraw[color=\colorDefect] (0,0) circle (\locpx) node[above=-2pt] {$\scriptstyle-$}
                                (\locpa,0) circle (\locpx) node[below] {$\scriptstyle+$};
  \end{scope}
  \end{tikzpicture}}
  \ee
  ~\\
Similarly, the index of both segments of the first two circles in \eqref{eq:exa-circles} 
is zero and the index of both segments of  the other two circles in \eqref{eq:exa-circles}
is equal to 1, while the index of the single segment of each of the two circles
  \be
  \raisebox{-5.3em}
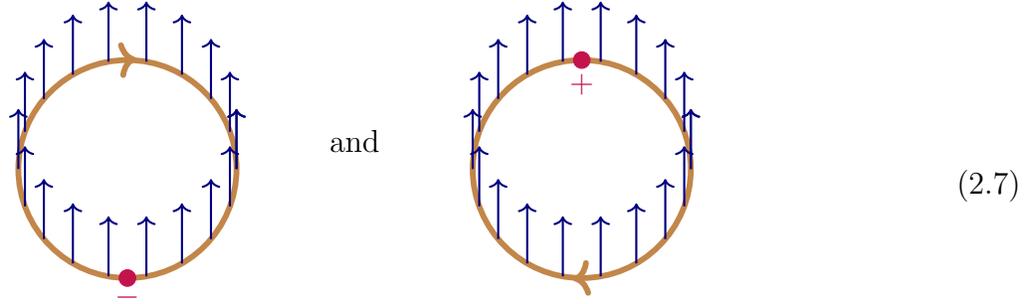
{\begin{tikzpicture}
  \begin{scope}[shift={(-2*\locpg,-\locpF)}]
  \scopeArrow{0.262}{<}
  \draw[line width=\locpc,color=\colorCircle,postaction={decorate}] (0,0) circle (\locpA);
  \end{scope}
  \filldraw[color=\colorDefect] (0,-\locpA) circle (\locpX) node[below=-1pt] {$-$};
  \foreach \xx in {0,20,...,360} \draw[thick,color=\colorFrame,->] (\xx:\locpA) -- +(0,\locpb);
  \end{scope}
  \begin{scope}[shift={(2*\locpg,-\locpF)}]
  \scopeArrow{0.762}{<}
  \draw[line width=\locpc,color=\colorCircle,postaction={decorate}] (0,0) circle (\locpA);
  \end{scope}
  \filldraw[color=\colorDefect] (0,\locpA)  circle (\locpX) node[below=1pt] {$+$};
  \foreach \xx in {0,20,...,360} \draw[thick,color=\colorFrame,->] (\xx:\locpA) -- +(0,\locpb);
  \end{scope}
  \node at (0,-0.5*\locpF) {and};
  \end{tikzpicture}}
  \label{eq:2-2-circle}
  \ee
is equal to 2.

To capture the information contained in the indices of the segments of a defect
one-manifold we introduce the following concept:

\begin{Definition} \label{def:integerref}
Let $\epsilon \,{=}\, (\epsilon_i)_{i=1,2,...,n}$ be an $n$-tuple of signs, which we consider
as either linearly or cyclically ordered. An $n$-tuple 
$\kappa \,{=}\, (\kappa_{i})_{i=1,2,...,n} \,{\in}\, \Z^n_{\phantom I}$ of integers 
is said to be a tuple of \emph{framing indices} associated with a linearly ordered $n$-tuple 
$\epsilon$ of signs iff, for every $i\,{\in}\,\{1,2,...\,,n-1\}$,
$\kappa_{i}$ is even if the product $\epsilon_{i+1}\,\epsilon_i$ is positive,
while $\kappa_{i}$ is odd if $\epsilon_{i+1}\,\epsilon_i$ is negative.
If $\epsilon$ is considered as cyclically ordered, then we impose in addition
the same rule on $\kappa_{n}$ as a function of the product $\epsilon_n\,\epsilon_1$.
\end{Definition}

Later on, only the homotopy class of the vector field will matter; accordingly,
the datum $\chi$ of an
\undefectonemanifold\ is equivalent to the datum of a tuple $\kappa$ of framing indices for
the signs $\epsilon$ of the defect points. (This motivates the terminology `framing index';
compare also Remark \ref{rem:framing} below.) Accordingly, we will also use the notation
$(\L,\epsilon,\kappa)$ in place of $(\L,\epsilon,\chi)$ for $2$-framed defect one-manifolds.

\medskip

We still have to introduce the general \undefectsurface s that can have general
\undefectonemanifold s as their boundary components. This is done as follows. Again we start
with a compact oriented surface $\vSigma$, now possibly with corners. Again there is
an embedding $\delta\colon I^{\sqcup n} \,{\sqcup}\, (\SS^1)^{\sqcup m}\,{\to}\, \vSigma$
as an additional structure. But now we allow for more general embeddings than before: the
image of an interval or a circle is also allowed to be contained in a connected component 
of the boundary $\partial\vSigma$. 
(Thus the intersection of the image of an interval with $\partial\vSigma$ is either empty, 
equal to the image of the end points, or equal to the image of the interval.)
If this is the case, then we call the image of the interval or circle an unlabeled 
\emph{free boundary}. The end points of a free boundary interval are corners of $\vSigma$;
they constitute additional marked points of $\Sigma$, to which we still refer as defect 
points. Still at each defect point on $\partial\vSigma$ the vector field $\chi$ 
on $\vSigma$ provides a non-vanishing vector; for $\delta$ a free boundary the so obtained 
vector at a defect point $p \,{\in}\, \delta \,{\cap}\,\partialg\vSigma$ is outward-poin\-ting
iff $\delta$ is oriented towards $\partialg\vSigma$, i.e.\ iff $p$ is a negative point.
If a circle is not mapped by $\delta$ to a boundary
component of $\vSigma$, then its image has again to be contained in the interior of $\vSigma$. 

Further, suppose that a connected component 
of $\partial\vSigma$ contains at least one free boundary segment. A connected component of the
complement of the union of the free boundaries of that connected component is then called 
a \emph{gluing interval} (see the picture \eqref{eq:example-gluing-int} below).
A connected component of $\partial\vSigma$ that does not contain
any free boundary is called a \emph{gluing circle}. If an interval is not mapped
by $\delta$ to a boundary component of $\vSigma$, then its end points must be mapped to a
gluing circle or gluing interval, and its interior to points in the interior of $\vSigma$. 
The images of the latter types of circles and intervals, which have non-empty
intersection with the interior of $\vSigma$, are called \emph{unlabeled defect lines}.
As an illustration, the following picture shows a defect surface whose underlying surface 
is a sphere with four holes and which has one unlabeled defect line, one unlabeled free boundary 
interval and one unlabeled free boundary circle, and one gluing interval and two gluing circles:
 \def \locpa  {1.6}  
  \be
  \raisebox{-5.3em}
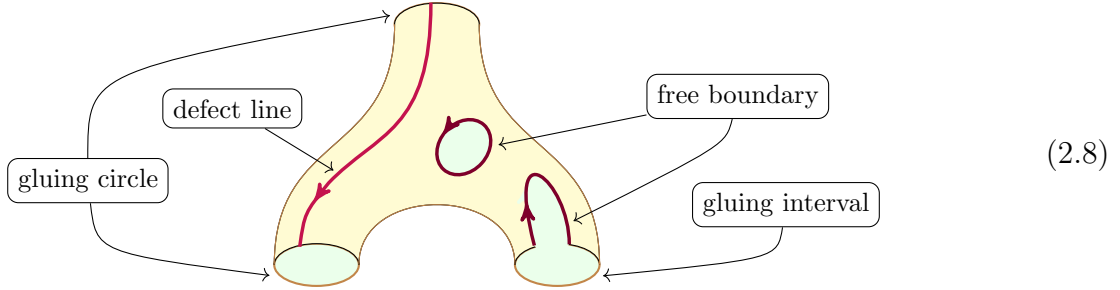
{\begin{tikzpicture}[scale=\locpa,tqft/cobordism/.style={draw}] 
  \begin{scope} [tqft/every boundary component/.style={fill=\fillBdyCircle,thick,
                draw=\colorCircle,scale=\locpa},
                tqft/cobordism edge/.style={draw,\colorCircle}]
  \pic[tqft/pair of pants,fill=\fillCobord,name=Pants,at={(0,0)},scale=\locpa];
  \scopeArrow{0.83}{\arrowDefect}
  \draw[line width=\widthDefect,color=\colorDefect,postaction={decorate}] (-0.05,0.18)
       .. controls (-0.05,-0.9) and (-0.5,-0.9) .. (-0.9,-1.32) 
       node(NodeD1)[left=-11pt,yshift=2pt] {} .. controls (-1.1,-1.57) .. (-1.14,-1.85);
  \end{scope}
  \scopeArrow{0.27}{\arrowDefect}
  \filldraw[line width=\widthDefect,fill=\fillBdyCircle,draw=\colorFreebdy,postaction={decorate},rotate=48]
        (-.61,-.85) ellipse (.25 and .20) node(NodeF2)[right=4pt] {};
  \end{scope}
  \scopeArrow{0.27}{\arrowDefect}
  \fill[line width=\widthDefect,\fillBdyCircle,postaction={decorate}]
        (.81,-1.84) .. controls (.55,-0.84) and (1.1,-1.3) .. (1.1,-1.82) 
         .. controls (1,-2.1) .. cycle;
  \draw[line width=\widthDefect,draw=\colorFreebdy,postaction={decorate}]
        (.81,-1.84) .. controls (.55,-0.84) and (1.1,-1.3) .. (1.1,-1.83)
	node(NodeF1)[above=2pt,xshift=-3pt] {};
  \node(NodeDx) [draw=black,rectangle, rounded corners]
        at (-1.7,-0.7) {\footnotesize defect line};
  \node(NodeFx) [draw=black,rectangle, rounded corners]
        at (2.5,-0.6) {\footnotesize free boundary};
  \node(NodeGx) [draw=black,rectangle, rounded corners]
        at (-2.9,-1.3) {\footnotesize gluing circle};
  \node(NodeGy) [draw=black,rectangle, rounded corners]
        at (2.9,-1.5) {\footnotesize gluing interval};
  \node(NodeG1) at (-0.3,0.09) {};
  \node(NodeG2) at (-1.3,-2.1) {};
  \node(NodeG3) at (1.29,-2.1) {};
  \draw[->] (NodeFx) .. controls +(-0.1,-0.4) .. (NodeF1);
  \draw[->] (NodeFx) -- (NodeF2);
  \draw[->] (NodeDx) -- (NodeD1);
  \draw[->] (NodeGx) .. controls +(0,0.9) .. (NodeG1);
  \draw[->] (NodeGx) .. controls +(0.07,-0.6) .. (NodeG2);
  \draw[->] (NodeGy) .. controls +(-0.07,-0.45) .. (NodeG3);
  \end{scope}
  \end{scope}
  \end{tikzpicture}}
  \label{eq:example-gluing-int}
  \ee
  ~\\
A 2-\emph{framing} on a general surface $\vSigma$ containing defect lines $\delta$ is a
non-vanishing vector field $\chi$ on $\vSigma$ that is parallel to, and has the same
direction as, each defect line and each free boundary. (On the other hand, there is 
no such condition restricting the vector field near gluing segments.)
    
We can now generalize Definition \ref{def:undefsurf0}(i) to

\begin{Definition}
An \emph{\undefectsurface} is a triple $(\vSigma,\delta,\chi)$, where $\vSigma$ is a compact 
oriented surface, possibly with boundary and possibly with corners, $\delta$ is the union of 
unlabeled defect lines in $\vSigma$ and of unlabeled free boundary intervals 
on the boundary $\partial\vSigma$, and $\chi$ is a $2$-framing on $(\vSigma,\delta)$.
 \\
An \emph{isomorphism $\varphi\colon (\vSigma,\delta,\chi) \,{\to}\, (\vSigma',\delta',\chi') $ 
of \undefectsurface s} is a diffeomorphism of the underlying manifolds that respects the 
orientations and the vector fields and that maps defect lines bijectively to defect lines.
\end{Definition}

A corner of $\vSigma$ is necessarily one of the end points of a free boundary interval;
as a consequence, the vector field at a corner is parallel to that free boundary.

Given an \undefectsurface\ $(\vSigma,\delta,\chi)$, we split its boundary as
  \be
  \partial\vSigma = \partialg\vSigma \cup \partialf\vSigma
  \ee
into the two parts
that consist of gluing segments and of free boundary segments, respectively. (Each of the
two parts can be empty; their intersection $\partialg\vSigma \,{\cap}\, \partialf\vSigma\,{=}\,
\partial(\partialf\vSigma)$ is the set of corners of $\vSigma$.)
We refer to $\partialg\vSigma$ as the \emph{gluing boundary} of $\vSigma$. The gluing boundary
$\partialg\vSigma$ becomes in the following manner an \undefectonemanifold.
The embedding $\iota\colon \partial\vSigma \,{\hookrightarrow}\, \vSigma$ gives rise
to an embedding $T_p(\partial\vSigma) \,{\hookrightarrow}\, \iota^*(T_p\vSigma)$
of the tangent space at every point $p \,{\in}\, \partial\vSigma$. Further, by using the
inward-pointing normal $n_p$ (with respect to some chosen auxiliary metric on $\vSigma$)
at $p$ one can then identify $(0,\xi) \,{\in}\, T_p(\partial\vSigma){\oplus}\underline{\R}$
with $\xi n_p \,{\in}\, \iota^*(T_p\vSigma)$. This provides an isomorphism of the tangent
bundle of $\vSigma$, restricted to the boundary, with $T(\partial\vSigma){\oplus}\underline{\R}$.
This way the $2$-framing on $\vSigma$ induces a $2$-framing on the boundary, whereby in particular
the gluing boundary $\partialg\vSigma$ is endowed with the structure of a (not necessarily
connected) $2$-framed defect one-manifold.
We denote the so obtained \undefectonemanifold\ by $\partialg(\vSigma,\delta,\chi)$.

Instead of using the inward-pointing normal $n_p$, yielding
$\partial_\text{inward}(\vSigma,\delta,\chi) \,{\equiv}\, \partialg(\vSigma,\delta,\chi)$
we could as well use the outward-pointing normal $-n_p$. This would yield
another \undefectonemanifold\ $\partial_\text{out}(\vSigma,\delta,\chi)$
that differs from $\partial_\text{inward}(\vSigma,\delta,\chi)$ by replacing the
vector field $\psi$ on $\partialg(\vSigma,\delta,\chi)$ by
$\psi_\text{outward} \,{:=}\, \overline{\psi}$, where the overbar denotes
a flip of the $\underline{\R}$-coordinate, i.e.
  \be
  \overline\psi(p) = (\alpha,{-}\xi) ~~:\Longleftrightarrow~~
  \psi(p) = (\alpha,\xi) ~\in T(\partial\vSigma)\,{\oplus}\,\underline{\R} \,.
  \label{eq:flipvectorfield}
  \ee
This motivates the following

\begin{Definition} \label{def:opposite-unlabeled}
The \emph{opposite} $\overline{(\L,\epsilon,\chi)}$ of an \undefectonemanifold\
$(\L,\epsilon,\chi)$ consists of the manifold $\L$ taken with opposite orientation, 
which we denote by $\cc{\L}$, the same marked points but with flipped signs,
and of the flipped vector field in the sense of \eqref{eq:flipvectorfield}, i.e.
  \be
  \overline{(\L,\epsilon,\chi)} \,:=\, (\overline\L,-\epsilon,\overline\chi) \,. 
  \label{eq:overlineL}
  \ee
\end{Definition}

The following picture shows examples of a defect circle and a defect interval together 
with their opposites (recall that only the homotopy class of the vector field matters):
    \def\locpa  {1.45}   
    \def\locpA  {2.16}   
    \def\locpb  {0.71}   
    \def\locpB  {0.74}   
    \def\locpc  {2.2pt}  
    \def\locpd  {0.12}   
    \def\locpD  {0.02}   
    \def\locpe  {2.5cm}  
    \def\locpE  {15pt}   
    \def\locpf  {3.1}    
    \def\locpF  {8.6}    
    \def\locpg  {3.0}    
    \def\locpx  {0.10}   
  \be
  \raisebox{-6.1em}{\begin{tikzpicture}
  \scopeArrow{0.144}{<} \scopeArrow{0.422}{<} \scopeArrow{0.618}{<}
  \draw[line width=\locpc,color=\colorCircle,postaction={decorate}]
       node[left=\locpe,color=black] {$(\L,\epsilon,\chi) ~=$} (0,0) circle (\locpa);
  \end{scope} \end{scope} \end{scope}
  \foreach \xx in {80,60,...,-90} \draw[thick,color=\colorFrame,->] (\xx:\locpa) -- +(0,-\locpb);
  \foreach \xx in {100,119,...,180} \draw[thick,color=\colorFrame,->] (\xx:\locpa) -- (\xx:\locpB);
  \draw[thick,color=\colorFrame,->] (190:\locpa) -- +(0.983*\locpb,0.18*\locpb);
  \draw[thick,color=\colorFrame,->] (212:\locpa) -- +(0.707*\locpb,0.707*\locpb);
  \draw[thick,color=\colorFrame,->] (233:\locpa) -- +(0.18*\locpb,0.983*\locpb);
  \draw[thick,color=\colorFrame,->] (248:\locpa) -- +(-0.983*\locpb,0.18*\locpb);
  \draw[thick,color=\colorFrame,->] (258:\locpa) -- +(-0.955*\locpb,-0.29*\locpb);
  \draw[thick,color=\colorFrame,->] (265:\locpa) -- +(-0.525*\locpb,-0.85*\locpb);
  \filldraw[color=\colorDefect] (0,\locpa)  circle (\locpd) node[above] {$-$} 
                                (-\locpa,0) circle (\locpd) node[left=1pt] {$-$} 
                                (0,-\locpa) circle (\locpd) node[above=1pt] {$+$};
  \begin{scope}[shift={(\locpF,0)}]
  \scopeArrow{0.175}{<} \scopeArrow{0.403}{<} \scopeArrow{0.694}{<}
  \draw[line width=\locpc,color=\colorCircle,postaction={decorate}]
       node[left=\locpe,color=black]{$\Longrightarrow\quad~\overline{(\L,\epsilon,\chi)}~=$}
       (0,0) circle (\locpa);
  \end{scope} \end{scope} \end{scope}
  \foreach \xx in {87,67,...,-93} \draw[thick,color=\colorFrame,->] (\xx:\locpa) -- +(0,\locpb);
  \foreach \xx in {104,119,...,179} \draw[thick,color=\colorFrame,->] (\xx:\locpa) -- (\xx:\locpA);
  \draw[thick,color=\colorFrame,->] (193:\locpa) -- +(-0.983*\locpb,0.18*\locpb);
  \draw[thick,color=\colorFrame,->] (207:\locpa) -- +(-0.955*\locpb,0.29*\locpb);
  \draw[thick,color=\colorFrame,->] (227:\locpa) -- +(-0.866*\locpb,0.5*\locpb);
  \draw[thick,color=\colorFrame,->] (245:\locpa) -- +(-0.525*\locpb,0.85*\locpb);
  \draw[thick,color=\colorFrame,->] (259:\locpa) -- +(-0.394*\locpb,0.92*\locpb);
  \filldraw[color=\colorDefect] (0,\locpa)  circle (\locpd) node[below=1pt] {$+$} 
                                (-\locpa,0) circle (\locpd) node[right=1pt] {$+$} 
                                (0,-\locpa) circle (\locpd) node[below] {$-$};
  \end{scope}
  \begin{scope}[shift={(\locpE-\locpe,-\locpf)}]
  \draw[line width=\locpc,color=\colorCircle] 
       node[left=\locpE,color=black]{$(\L,\epsilon,\chi)~=$} (-\locpD,0) -- (\locpg+\locpD,0);
  \foreach \xx in {0,1,...,4} \draw[thick,color=\colorFrame,->]
             (0.06*\locpg*\xx,0) -- +(90-18*\xx:\locpb);
  \draw[thick,color=\colorFrame,->] (0.3*\locpg,0) -- +(7:\locpb);
  \draw[thick,color=\colorFrame,->] (0.42*\locpg,0) -- +(-24:\locpb);
  \foreach \xx in {7,8,9,10} \draw[thick,color=\colorFrame,->]
             (-0.4*\locpg+0.14*\locpg*\xx,0) -- +(90-18*\xx:\locpb);
  \draw[line width=\locpc,color=\colorCircle,->] (0.83*\locpg,0) -- (0.85*\locpg,0);
  \filldraw[color=\colorDefect] (0,0) circle (\locpx) node[below] {$\scriptstyle+$}
                                (\locpg,0) circle (\locpx) node[above=-1pt] {$\scriptstyle-$};
  \end{scope}
  \begin{scope}[shift={(\locpF cm+\locpE-\locpe,-\locpf)}]  %
  \draw[line width=\locpc,color=\colorCircle]
       node[left=\locpE,color=black]{$\Longrightarrow\quad~\overline{(\L,\epsilon,\chi)}~=$}
       (-\locpD,0) -- (\locpg+\locpD,0);
  \foreach \xx in {0,1,...,4} \draw[thick,color=\colorFrame,->]
             (0.06*\locpg*\xx,0) -- +(-90+18*\xx:\locpb);
  \draw[thick,color=\colorFrame,->] (0.3*\locpg,0) -- +(-7:\locpb);
  \draw[thick,color=\colorFrame,->] (0.42*\locpg,0) -- +(24:\locpb);
  \foreach \xx in {7,8,9,10} \draw[thick,color=\colorFrame,->]
             (-0.4*\locpg+0.14*\locpg*\xx,0) -- +(-90+18*\xx:\locpb);
  \draw[line width=\locpc,color=\colorCircle,->] (0.83*\locpg,0) -- (0.85*\locpg,0);
  \filldraw[color=\colorDefect] (0,0) circle (\locpx) node[above=-1pt] {$\scriptstyle-$}
                                (\locpg,0) circle (\locpx) node[below] {$\scriptstyle+$};
  \end{scope}
  \end{tikzpicture}}
  \label{pic_oppositeL}
  \ee
  ~

To summarize, we have obtained two bicategories $\Bfrdef$ and $\Bfrdefcl$:
\Itemize
\item
Objects of $\Bfrdef$ are \undefectonemanifold s,
objects of $\Bfrdefcl$ are closed \undefectonemanifold s.
(Recall that an \undefectonemanifold\ is oriented and endowed with a $2$-framing.)

\item
A 1-morphism $\LL \,{\to}\, \LL'$ in $\Bfrdef$ is an \undefectsurface\ $\vSigma$ together
with a \emph{boundary parametrization}, i.e.\ an isomorphism 
$\partialg\vSigma \,{\xrightarrow{~\cong~}}\, \overline{\L'_{\phantom|}} \,{\sqcup}\, \L$ of
\undefectonemanifold s that extends to a small collar over $\partialg\vSigma$.
The 1-morphisms of $\Bfrdefcl$ are those of $\Bfrdef$ for which the
\undefectsurface\ $\vSigma$ is closed.

\item
Composition of 1-morphisms is given by gluing an incoming boundary and an outgoing boundary
which are each others' opposites. These boundaries can consist of gluing circles or gluing
intervals; the gluing has to account for the parametrizations of the boundaries.    

\item
A 2-morphism $\varphi\colon \vSigma \,{\to}\, \vSigma' $ in $\Bfrdef$ between two 1-morphisms
is represented by an isomorphism  $\varphi$ of \undefectsurface s that respects the
boundary parametrizations.
Two isomorphisms $\varphi_0,\varphi_1\colon \Sigma \,{\to}\, \Sigma'$ represent the same
2-morphism if and only if there is an isotopy $h\colon \Sigma\Times [0,1] \,{\to}\, \Sigma'$ 
with $h(-,0) \,{=}\, \varphi_0$ and $h(-,1) \,{=}\, \varphi_1$ that satisfies 
$\delta' \,{=}\, h(-,t)(\delta)$ and $\chi' \,{=}\, (h_t)_*(\chi)$ for all
$t \,{\in}\, [0,1]$. 2-morphisms in $\Bfrdefcl$ are defined accordingly.

\item
The vertical composition of 2-morphisms is induced by composition of isomorphisms.
The horizontal composition of morphisms is given by gluing of surfaces
along gluing circles or gluing intervals. 

\end{itemize}

\noindent
These bicategories are symmetric monoidal, with the monoidal structure given by disjoint union.
(See \cite{schom1} for the definition of a symmetric monoidal bicategory, and Section 
3.1.4 therein for an outline of how to construct the symmetric monoidal bicategory $\Bfrdef$.) 
The full subcategories of $\Bfrdef$ and $\Bfrdefcl$ consisting of \emph{fine}
\undefectsurface s (with and without corners, respectively) are symmetric monoidal as well.

The values of the modular functor on 2-morphisms provide us with representations of the 
relevant mapping class groups of \defectsurface s. This will be analyzed in detail in 
Section \ref{sec:mcg}.

\begin{Remarks}\label{Remark:structure-def-bord}
\Itemizeiii
\item[(i)]
To get a well-defined horizontal composition of 2-morphisms, we should work with 
collars. This is standard \Cite{Thm.\,1.3.12}{KOck}, and nothing new happens in our 
context. Accordingly we suppress this issue.

\item[(ii)]
By definition of the 2-morphisms, in case the \undefectsurface\ $\vSigma$ does not have a
free boundary and does not have any defects, the 2-morphisms from $\vSigma$ to itself 
form the framed mapping class group of $\vSigma$.

\item[(iii)]
Consider a \defectsurface\ $(\Sigma,\delta,\chi)$ with defects $\delta$ and $2$-framing $\chi$,
and the same underlying surface with the same defect lines but with another $2$-framing $\chi'$, 
together with a homotopy $\chi_{t}$ from 
$\chi$ to $\chi'$, that is, $\chi_{t}$ is a smooth family $\chi_{t}\colon \Sigma \To T\Sigma$ of 
framing vector fields  for $t \iN [0,1]$, such that $\chi_{0} \,{=}\, \chi$ and $\chi_{1} \,{=}\, \chi'$.
Then $(\Sigma, \delta, \chi)$ and $(\Sigma, \delta, \chi')$ are isomorphic in $\Bfrdef$,
i.e.\ there exists an automorphism $\varphi$ of $\Sigma$ which preserves $\delta$ and satisfies
$T\varphi(\chi) \,{=}\, \chi'$. This can be seen by considering the cylinder $\Sigma
\,{\times}\, [0,1]$ over $\Sigma$ with the vector field $\widetilde{\chi} \,{=}\, (\chi_{t},t)$:
Since $\Sigma \,{\times}\, [0,1]$ is compact, the vector field $\widetilde{\chi}$ has a complete 
flow $\varphi_{t}\colon \Sigma \,{\cong}\, (\Sigma,0) \,{\to}\, (\Sigma,t) \,{\cong}\,\Sigma$,
which is a 1-parameter family of automorphisms of $\Sigma$ that preserve $\delta$ 
(since each $\chi_{t}$ is tangential to $\delta$) with $\varphi_{0} \,{=}\, \id_{\Sigma}$. It
follows from the flow equation 
$\frac{\mathrm d}{\mathrm dt} \varphi_{t}\big|_{t=0}^{} \,{=}\, \chi_{t}$ that
$T\varphi_{1}(\chi_{0}(p)) \,{=}\, \frac{\mathrm d}{\mathrm dt}\varphi_{1}(\varphi_{t}(p))
\big|_{t=0} \,{=}\, \frac{\mathrm d}{\mathrm dt}\varphi_{1+t}(p)\big|_{t=0} \,{=}\, \chi_{1}(p)$.

\item[(iv)]
Since, as noted in Remark \ref{rem:...iii:gfine}(iii), \gfine\ \defectsurface s compose
to a \gfine\ \defectsurface, there is a symmetric monoidal
sub-bicategory $\Bfrdefgfine$ whose 1-morphisms are \gfine\ \defectsurface s.
\end{itemize}
\end{Remarks}


\subsection{Labels for defect bordisms}\label{sec:labeling}

We are now going to assign an additional algebraic datum to each connected component of the 
complement of the defect lines and boundaries in a defect surface. Afterwards we also
assign a corresponding datum to each defect line and to each free boundary segment.
Before formulating this prescription, several further concepts need to be
recalled. All algebraic categories of our interest are assumed to be finite, 
abelian and linear over a fixed algebraically closed field $\ko$. Similarly we 
require functors and natural transformations to be linear, unless specified otherwise.
For the notion of a \emph{finite tensor category} see e.g.\ \cite{EGno}. We will heavily
use that every object of such a category has a left and a right dual; we do not assume any
relation between the two duals.  Our conventions concerning dualities of a rigid category
$\calc$ are as follows.  The right dual of an object $c$ is denoted by $c^\vee$, and the
right evaluation and coevaluation are morphisms
  \be
  \evr_c \in \Hom_\calc (c^\vee \otimes c,\one) \qquad{\rm and}\qquad
  \coevr_c \in \Hom_\calc (\one ,c \otimes c^\vee) \,,
  \ee
while the left evaluation and coevaluation are
  \be
  \evl_c \in \Hom_\calc (c \otimes \Vee c,\one) \qquad{\rm and}\qquad
  \coevl_c \in \Hom_\calc( \one ,\Vee c \otimes c)
  \ee
with $\Vee c$ the left dual of $c$.

Further recall that a (left) module category over a finite tensor category $\cala$ (or, for short,
an $\cala$-mo\-dule), is a finite linear category $\calm \eq \acalm$ together with a bilinear
functor, exact in the first variable, from $\cala\Times\calm$ to $\calm$, which
we call the \emph{action} of $\cala$ and just denote by a dot `$.$', as well as with
natural isomorphisms $\mu$ and $\lambda$ with components $\mu_{a,b,m}\iN 
\Hom_\calm ((a\oti b).m,a.(b.m))$ and $\lambda_m \iN \Hom_\calm (\one_\Cala .m, m)$ 
that satisfy pentagon and triangle relations
analogous to the associator and unit constraint of a monoidal category.
Right $\cala$-modules and $\cala$-$\calb$-bimodules are defined analogously. For ease of
notation, we will use the symbol $\calm$ both for module and for bimodule categories.
It is natural to consider a bimodule category $_{\cala}\calm_{\calb}$ as a 1-morphism 
$\cala \,{\to}\, \calb$ in the tricategory $\FinTensL$ that has finite tensor categories 
as objects, finite bimodule categories as 1-morphisms, and categories
$\Funle_{\Cala,\calb}({}_{\cala}\calm_{\calb},{}_{\cala}\caln_{\calb})$ of left exact 
bimodule functors and bimodule natural transformations as 2- and 3-morphisms, respectively. 
(Alternatively, one could consider a tricategory with categories
$\Funre_{\Cala,\calb}({}_{\cala}\calm_{\calb},{}_{\cala}\caln_{\calb})$
of right exact bimodule functors as 2- and 3-mor\-phisms. In this paper we focus on 
the formulation with left exact functors.)

There is then an obvious notion of a \emph{composable string} of bimodule categories, and 
likewise it is clear what a \emph{cyclically composable string} of bimodule categories is.
We can also allow for left and right modules, respectively, as the ends of
a string of composable bimodule categories, by considering a 
left $\cala$-module as an $\cala$-\vect-bimodule and a right $\calb$-module as a
\vect-$\calb$-bimodule; thus e.g.\ a right module category $\calm_{\Cala}$,
a bimodule category ${}_{\cala}\caln_{\calb}$ and a left module category ${}_{\calb}\calk$
form a composable string that constitutes a 1-morphism $\vect \,{\to}\, \vect$ in $\FinTensL$.

\medskip

Let now $\vSigma$ be a surface with defect lines $\delta$. We denote by 
  \be
  \vSigma \bs1 := \delta\,{\cup}\,\partial\vSigma
  \ee
the union of the defect lines and the boundary of $\vSigma$. 

\begin{Definition} \label{def:twocell}
By a \emph{\twocell} of $\vSigma$ we mean a connected component of the complement of 
$\vSigma \bs1$ in $\vSigma$ together with the adjacent subset of $\vSigma \bs1$.
\end{Definition}
	 
A \defectsurface\ is called fine iff the underlying \undefectsurface\ is fine, i.e.\
iff every \twocell\ is topologically a disk.
A \defectsurface\ is called \gfine\ iff the underlying \undefectsurface\ is \gfine.

We make the following assignments,
which are in line with existing literature (see e.g.\ Table 1 in \cite{kiKon}):
\Itemize
\item
To a \twocell\ we assign finite tensor category.\,%
  \footnote{~It is worth noting that we do not require the existence of a pivotal structure on
  these categories. The additional geometric structure of a framing 
  allows us to dispense with pivotal structures; compare Remark \ref{rem:pivot}.
  Also note that we do not make any assumption about the topology of the \twocells; 
  see however the notion of a fine defect surface in Definition \ref{def:undefsurf0}.}

\item
To a defect line that separates \twocells\ labeled by finite tensor categories
$\cala$ and $\calb$ we assign an $\cala$-$\calb$- or $\calb$-$\cala$-bimodule category,
depending on the relative orientations of the defect line and the adjacent \twocells\
(see the picture \eqref{eq:fixing-convent} below).

\item
Similarly, to a \emph{free boundary} we assign a left or right module category
over the monoidal category associated to the adjacent \twocell.
\end{itemize}

\noindent
The following picture fixes uniquely our convention for the bimodule
categories assigned to defect lines:
 \def\locpc  {2.2pt}  
  \be
  \raisebox{-4.0em}
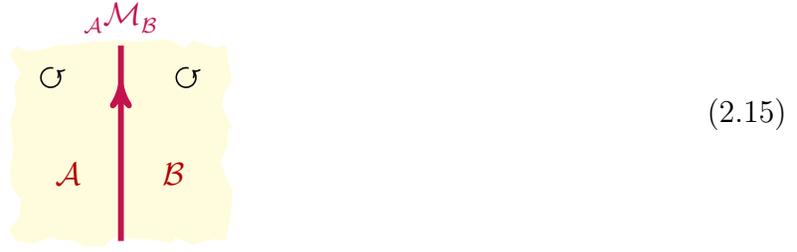
{\begin{tikzpicture}[scale=1.2]
  [decoration={random steps,segment length=2mm}]
  \fill[\colorSigma,decorate] (-1.4,0) rectangle (1.4,-2.6);
  \scopeArrow{0.81}{\arrowDefect}
  \draw[line width=\locpc,color=\colorDefect,postaction={decorate}] (0,-2.6) -- (0,0)
     node[above,color=\colorDefect] {${}_\cala^{}\calm_\calb^{}$};
  \end{scope}
  \begin{scope}[shift={(-0.9,-0.3)}] \drawOrientation; \end{scope}
  \begin{scope}[shift={(0.9,-0.3)}] \drawOrientation; \end{scope}
  \draw[thick,color=\colorFTC]  (0.7,-1.7) node {$\calb$};
  \draw[thick,color=\colorFTC] (-0.7,-1.7) node {$\cala$};
  \end{tikzpicture}}
  \label{eq:fixing-convent}
  \ee
  ~\\
Similarly, for free boundaries our convention is fixed by the following pictures:
  \be
  \raisebox{-3.5em}
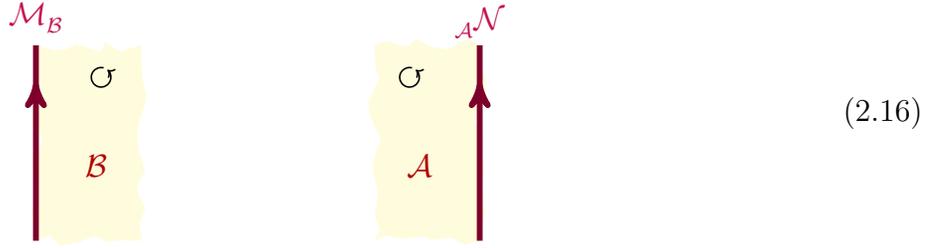
{\begin{tikzpicture} 
  [decoration={random steps,segment length=2mm}]
  \fill[\colorSigma,decorate] (0,0) -- (0,-2.6) -- (1.4,-2.6) -- (1.4,0) -- cycle;
  \scopeArrow{0.81}{\arrowDefect}
  \draw[line width=\locpc,color=\colorFreebdy,postaction={decorate}]
     node[above,color=\colorDefect] {$\calm_\calb^{}$} (0,-2.6) -- (0,0);
  \end{scope}
  \begin{scope}[shift={(0.9,-0.3)}] \drawOrientation; \end{scope}
  \draw[thick,color=\colorFTC]  (0.8,-1.6) node[fill=\colorSigma] {$\calb$};
  \begin{scope}[shift={(5.9,0)}]
  \fill[\colorSigma,decorate] (-1.4,0) rectangle (0,-2.6);
  \scopeArrow{0.81}{\arrowDefect}
  \draw[line width=\locpc,color=\colorFreebdy,postaction={decorate}]
     node[above=-2pt,color=\colorDefect] {${}_\cala^{}\caln$} (0,-2.6) -- (0,0);
  \end{scope}
  \begin{scope}[shift={(-0.9,-0.3)}] \drawOrientation; \end{scope}
  \draw[thick,color=\colorFTC] (-0.8,-1.6) node[fill=\colorSigma] {$\cala$};
  \end{scope}
  \end{tikzpicture}}
  \ee

The labeling of the building blocks of $\vSigma$ induces a labeling of the segments of 
$\partialg\vSigma$ and thereby determines an assignment of labels for 
\defectonemanifold s: Defect points are labeled by bimodule categories, and 
free boundary segments as well as their end points by module categories, in such a way that,
together with the orientation of the defect points, the labels form a composable string of 
bimodule categories. 

We summarize our prescriptions in

\begin{Definition} \label{def:Bordcats}
$\Bfrdefdec$ is the following symmetric monoidal bicategory:
\Itemizeiii

\item[{\rm (i)}]
Objects of $\Bfrdefdec$, called \emph{\defectonemanifold s}, are tuples
$\LL \,{=}\, \big( \L,\epsilon,\chi,\{\calm_j\} \big)$ given by an \undefectonemanifold\ 
together with an assignment $\{\calm_j\}$ of labels to its marked points, consisting of a
bimodule category for each defect point forming the end of a defect line, and a module category
for each defect point at the end of a free boundary segment, in such a way that the finite tensor
categories involved in consecutive marked points match $($when taking orientations into account$)$.

\item[{\rm (ii)}]
$1$-morphisms $\LL\,{\to}\,\LL'$ of $\,\Bfrdefdec$, called \emph{\defectsurface s}, are
tuples 
  \be
  \Sigma = \big( \vSigma,\delta,\chi,\{\cala_k,\calm_l\} \big)
  \ee
consisting of an
\undefectsurface\ and an assignment of labels, together with an isomorphism
$\partialg\vSigma \,{\xrightarrow{~\cong\,}}\, \L \,{\sqcup}\, \overline{\L'}$ of
\defectonemanifold s, such that the labels in the interior and on the boundary of $\vSigma$ match.

\item[{\rm (iii)}]
A $2$-morphism from a $1$-morphism $\LL\,{\to}\,\LL'$ given by 
$\Sigma \,{=}\, \big( \vSigma,\delta,\chi,\{\cala_k,\calm_l\} \big)$ to a $1$-mor\-phism
with the same source and target given by $\Sigma' \,{=}\, \big(
\vSigma',\delta',\chi',\{\cala'_k,\calm'_l\} \big)$ is represented by an isomorphism
$\varphi\colon \vSigma \,{\to}\, \vSigma'$ of \undefectsurface s that preserves the
labels $\{\cala_k,\calm_l\}$ of the various strata. We call such an isomorphism a
\emph{morphism of defect surfaces}. Two morphisms  $\varphi,\varphi'\colon \Sigma \,{\to}\,
\Sigma'$ are equivalent iff they are equivalent as morphisms of \undefectsurface s. 

\item[{\rm (iv)}]
The vertical composition of $\,2$-morphisms is induced by composition of isomorphisms.
The horizontal composition of morphisms is given by gluing of surfaces
along gluing circles or gluing intervals. 

\item[{\rm (v)}]
The monoidal structure is given by disjoint union.
\end{itemize}
\end{Definition}

\begin{Remarks}
\Itemizeiii
\item[(i)]
We write generically $\{\cala_i,\calm_j\}$ etc.\ for the relevant 
sets of finite tensor and (bi)module categories.
For brevity, below we will usually suppress these additional data (which are objects
and 1-morphisms, respectively, of the tricategory $\FinTensL$) in our notation.

\item[(ii)]
Analogously as for unlabeled surfaces (see Remark \ref{rem:...iii:gfine}(iii)), gluing two
\emph{\gfine} \defectsurface s gives again a \gfine\ \defectsurface. $\Bfrdefdec$ has 
thus a symmetric monoidal
full sub-bicategory $\Bfrdefdecgfine$ whose 1-morphisms are \gfine\ \defectsurface s.
\end{itemize}
\end{Remarks}

The notion of opposite one-manifold extends as follows from the unlabeled to the labeled case:

\begin{Definition} \label{def:opposite}
The \emph{opposite} $\overline{\LL}$ of a (labeled) \defectonemanifold\ $\LL$ is the
opposite $\overline{(\L,\epsilon,\chi)}$ of the underlying unlabeled manifold
$(\L,\epsilon,\chi)$ together with the same assignment of labels.
\end{Definition} 

Two \defectsurface s $\Sigma_1$ and $\Sigma_2$ can be glued along a \defectonemanifold\ $\LL$
if and only if the corresponding \defectonemanifold\ $\LL_1 \,{\subset}\, \Sigma_1$ and
$\LL_2 \,{\subset}\, \Sigma_2$ are opposite to each other. To be able to work with smooth
manifolds, gluing is actually along collars. This is demonstrated in the following picture
which shows the situation locally around a defect point on a gluing segment:
    \def\locpa  {0.95}   
    \def\locpA  {1.65}   
    \def\locpb  {0.8}    
    \def\locpB  {2.35}   
    \def\locpc  {2.2pt}  
    \def\locpX  {0.11}   
    \def\locpy  {0.32}   
    \def\locpY  {0.08}   
  \be
  \raisebox{-4.6em}
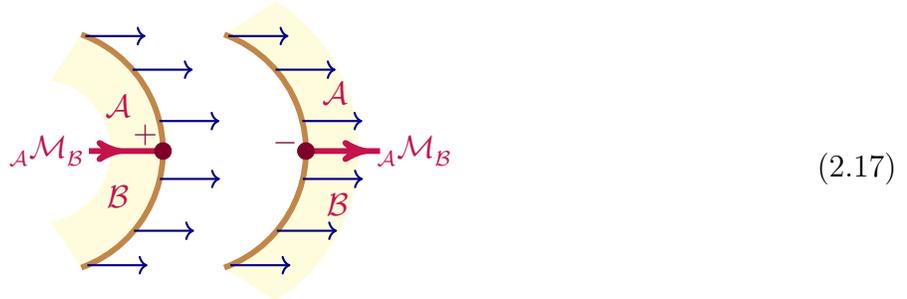
{\begin{tikzpicture}
  \fill[\colorSigma] (-70:\locpA) arc (-70:70:\locpA)
       -- (81:\locpa) arc (81:-81:\locpa) -- cycle;
  \draw[line width=\locpc,color=\colorCircle] (-70:\locpA) arc (-70:70:\locpA);
  \scopeArrow{0.45}{\arrowDefect}
  \draw[line width=\locpc,color=\colorDefect,postaction={decorate}] (0.7*\locpa,0)
       node[left=-3pt,color=\colorDefect] {${}_\cala^{}\calm_\calb^{}$} -- 
       node[above=8pt,xshift=-3pt,color=\colorDefect] {$\cala$}
       node[below=8pt,xshift=-3pt,color=\colorDefect] {$\calb$} (\locpA,0);
  \end{scope}
  \filldraw[color=\colorFreebdy] (\locpA,0) circle (\locpX) node[left=-2pt,yshift=6pt] {$-$};
  \foreach \xx in {-67,-40,...,70}
           \draw[thick,color=\colorFrame,->] (\xx:\locpA) -- +(\locpb,0);
  \begin{scope}[shift={(1.9,0)}]
  \fill[\colorSigma] (-70:\locpA) arc (-70:70:\locpA)
       -- (58:\locpB) arc (58:-58:\locpB) -- cycle;
  \draw[line width=\locpc,color=\colorCircle] (-70:\locpA) arc (-70:70:\locpA);
  \scopeArrow{0.84}{\arrowDefect}
  \draw[line width=\locpc,color=\colorDefect,postaction={decorate}] (\locpA,0) -- 
       node[above=13pt,xshift=-3pt,color=\colorDefect] {$\cala$} 
       node[below=11pt,xshift=-2pt,color=\colorDefect] {$\calb$}
       (\locpB+0.3*\locpa,0) node[right=-6pt,color=\colorDefect] {${}_\cala^{}\calm_\calb^{}$};
  \end{scope}
  \filldraw[color=\colorFreebdy] (\locpA,0) circle (\locpX) node[left=-1pt,yshift=3pt] {$+$};
  \foreach \xx in {-67,-40,...,70}
           \draw[thick,color=\colorFrame,->] (\xx:\locpA) -- +(\locpb,0);
  \end{scope}
  \end{tikzpicture}}
  \ee

\medskip 

The main result of this article is the construction of two specific modular functors.

\begin{Definition}\label{def:modularfunctor}
Let $\mathcal S$ be a symmetric monoidal bicategory. An $\mathcal S$-valued
\emph{framed modular functor} 
(with decoration data in finite categories) is a symmetric monoidal $2$-functor
  \be
  \ZZ:\quad \Bfrdefdec \xrightarrow{~} \mathcal S \,.
  \ee
\end{Definition}

Given the prescriptions for labels present in our setting, a natural choice of a target 
bicategory $\mathcal S$ is $\FUNle$, i.e.\ the bicategory
that has as objects finite \ko-linear categories, as 1-morphisms left exact functors
and as 2-morphisms natural transformations, with monoidal structure given by the
Deligne product. The restriction of functors to left exact ones is due to the fact that
the Deligne product of left exact functors 
is defined and provides the symmetric monoidal structure at the level of 1- and 
2-morphisms. The same also applies to right exact functors, and indeed we could have chosen
instead the symmetric monoidal bicategory $\FUNre$ whose morphisms are right exact functors.

To summarize, we will show:

\begin{thm} \label{thm:main}
There exists a state-sum construction that provides an explicit framed modular functor 
with values in $\FUNle$, as well as a framed modular functor with values in $\FUNre$.
\end{thm}

Several comments are in order:

\begin{Remarks} \label{rmk:mainthm}
\Itemizeiii
 
\item[(i)]
The existence of a similar functor for manifolds without defects has been
shown in \Cite{Cor.\,5}{doSs3}, invoking the cobordism hypothesis. 
Our approach is more direct and and can be directly compared with state-sum constructions.  

\item[(ii)]
The definition of a modular functor implies in particular that our construction is compatible
with factorization or, more specifically, with the gluing of \defectsurface s along 
gluing intervals or gluing circles. Indeed, the composition of two 1-morphisms in 
$\Bfrdefdec$, i.e.\ defect surfaces $\Sigma$ and $\Sigma'$, is by gluing the `outgoing'
part $\partialo\Sigma$ of the gluing boundary of $\Sigma$ with the `incoming' part
$\partiali\Sigma'$ of the gluing boundary of $\Sigma'$. For details, see Section \ref{sec:facto}.

\item[(iii)]
Our prescription may be seen as a kind of Turaev--Viro construction. Indeed, a standard 
semisimple Turaev--Viro construction corresponds to specializing our prescription by using 
only ``transparent'' labelings and assuming in addition the presence of a pivotal structure
(which allows one to eliminate
the framing), see Remark \ref{rem:TV}. Now in the standard Turaev--Viro situation a crucial
property is the independence of the choice of a triangulation. In contrast, the framed 
modular functor considered here is defined on arbitrary \defectsurface s
for which all \twocells\ are contractible, without picking a triangulation.

However, as explained in Section \ref{sec:refin-defect-surf}, the construction
of the modular functor does make use of such ``refining'' triangulations, all parts of 
which are labeled by transparent labels. As we show in Section \ref{sec:nonfine}, our 
construction is compatible with transparency, in the sense that the structure we have 
at our disposal is sufficient to define a block functor that is independent of the 
triangulation as a (co)limit. (Thus not only the precise position of transparently 
labeled defects is irrelevant -- this invariance up to homotopy is valid for any 
topological defect -- but not even their combinatorial configuration, i.e.\ the 
particular choice of refining triangulation, matters. These properties justify 
the qualification ``transparent''.)

\item[(iv)]
Mapping class group elements are specific 2-morphisms in $\Bfrdefdec$. Our construction thus
provides representations of mapping class groups. This is studied in Section \ref{sec:mcg}.
Note that since we are dealing with a state-sum construction, all mapping class
group elements are represented by genuine linear actions, so that there is no
need for a central extension of the mapping class groups.

\item[(v)]
The decoration data of the modular functor are themselves categories, e.g.\ a bimodule
category $\calm$ for a defect line. If we assign different bimodule categories $\calm$ 
and $\calm'$ to a defect line, the modular functor provides us with two functors. Given 
moreover a bimodule functor $\calm \,{\to}\, \calm'$, it is
reasonable to ask how these two functors are related. This relationship may be 
seen as a functoriality with respect to the decoration data; it is indeed realized
in our construction, but we will not discuss it in the present paper.

\item[(vi)]
The left exact version of the functor is compatible with operations on defect
labels in the following sense. As shown in Proposition \ref{prop:block-compose},
the contraction of a defect line is implemented by a composition of functors, while 
Theorem \ref{thm:fusion-rel-Del} implies that the fusion of two parallel
defect lines corresponds to a variant of the relative Deligne product of bimodule
categories, which is the composition of 1-morphisms in $\FinTensL$.
Thus in particular our construction is compatible with the identities in $\FinTensL$.
This, in turn, is implicit in the construction of the functor via refining triangulations.
\end{itemize}
\end{Remarks}


\section{Assigning categories to \defectonemanifold s} \label{sec:bdrycat}

The goal of this section is to define our modular functor on objects, that is,
to associate to any \defectonemanifold\ $\LL$ a finite \ko-li\-ne\-ar category $\ZZ(\LL)$.
We call these categories \emph{gluing categories}, because they are assigned to boundary
segments of \defectsurface s along which these can be glued together to form more complicated
\defectsurface s. The gluing category $\ZZ(\LL)$ will be defined as the category of objects in 
a Deligne product, endowed with the additional structure of \emph{balancings}. More concretely,
we first take the Deligne product of all categories that are assigned to the marked points of
$\LL$. We then take objects in this Deligne product together with a balancing for each 
monoidal category that is assigned to a segment of $\LL$. Such a balancing allows
one to swap the action of objects in a monoidal category from one (bi)module category in
a composable string to a neighboring one.

\subsection{Twisted bimodule categories}\label{sec:duals}

Let $\cala_1$ and $\cala_2$ be finite tensor categories and $\calm$ an
$\cala_1$-$\cala_2$-bimodule category. The double left and right dual functors of a finite tensor
category have a natural monoidal structure, so that we can twist the left and right actions on
$\calm$ by powers of the double left or right dual of $\cala_1$ and $\cala_2$, respectively. 
To describe these twisted actions comfortably, we introduce the following notation.
Given an object $a$ of a finite tensor category $\cala$, we use the shorthand ${}^{[\kappa]}a$,
for $\kappa \,{\in}\,\N$, for the $\kappa$-fold left dual $^{\vee\vee...\vee\!}a$ of $a$,
and analogously $a^{[\kappa]}$ for the $\kappa$-fold right dual; we also write
$^{[0]\!}a \,{=}\, a \,{=}\, a^{[0]}$. Thus the double left dual functor maps objects as
$a \,{\xmapsto{~~}}\, {}^{\vee\vee\!}a \,{=}\, {}^{[2]}a $. Further, in view of the canonical 
isomorphisms $ {}^{\vee\!}(a^{\vee}) \,{\cong}\, a \,{\cong}\, (^{\vee\!}a)^{\vee}$ it is natural
to extend these definitions by taking $^{[-\kappa]}a$ for $\kappa\,{\in}\, \N$ to be the
$\kappa$-fold right dual, $^{[-\kappa]}a \,{=}\, a^{[\kappa]}$, and vice versa,

Now for any pair $(\kappa_{1},\kappa_{2}) \,{\in}\, 2\Z \,{\times}\, 2\Z$ we denote by
${}_{}^{\kappa_{1}\!\!}\calm{}^{\kappa_{2}}_{}$ the bimodule category for which the left 
and right actions on $\calm$ are twisted by the $\kappa_{1}$- and $\kappa_{2}$-fold left 
and right dual, respectively, i.e.\ $\cala_1$ acts as $m \,{\mapsto}\, {}^{[\kappa_1]}a_1.m$
and $\cala_2$ as $m \,{\mapsto}\, m.a_{2}^{[\kappa_2]}$. We also abbreviate
${}_{}^{\kappa\!\!}\calm{}^{0} \,{\equiv}\, {}_{}^{\kappa\!\!}\calm$
and ${}_{}^{0\!\!}\calm{}^{\kappa} \,{\equiv}\, \calm{}^{\kappa}$.
Similarly, for every pair of \emph{odd} integers $\kappa_1$ and $\kappa_2$, the opposite category
$\calm \Opp$ can be endowed with the structure of an $\cala_2$-$\cala_1$-bimodule by setting
  \be
  a_2\,.\,\cc m := \cc{m\,.\,{}^{[\kappa_1]}a_2} \qquad\text{and}\qquad
  \cc m\,.\,a_1 := \cc{ a_{1}^{[\kappa_2]}.m} \,.
  \label{eq:opp-module}
  \ee
Here we write $\cc x$\ for the object $x\,{\in}\,\calm$ seen as an object in $\calm\Opp$. 
Whenever convenient we will from now on also use the
notation $\cc{\calx}$ for the opposite category of any category $\calx$, as well as
  \be
  \calx^\epsilon := \left\{\begin{array}{ll}
  \calx & \text{ for } \epsilon \,{=}\, {+}1 \,,
  \nxl2
  \CC\calx ~& \text{ for } \epsilon \,{=}\, {-}1 \,. \eear \right.
  \label{eq:M-epsilon}
  \ee
Further, we denote the bimodule categories
with actions \eqref{eq:opp-module} by ${}_{}^{\kappa_1\!}\CC\calm{}^{\kappa_2}_{}$, and
analogously for left and right modules. Note that the so obtained $\Z{\times}\Z$-torsor 
of bimodule categories does not have a natural section. Allowing also for the modules
${}_{}^{\kappa_1\!}\CC\calm{}^{\kappa_2}_{}$, we can consider cyclically or linearly
composable strings of the more general form $(\calm_i^{\epsilon_i})_{1\leq i\leq n}$.

\begin{Remark}\label{rem:pivot}
A pivotal structure on a finite tensor category $\cala$, if it exists, furnishes a monoidal
equivalence ${?}^{\vee\vee} {\Rightarrow}\, \Id_\Cala$ and thus allows one to
identify all the module categories that result from twists by powers of the double dual. But
we do not assume the existence of a pivotal structure and accordingly use all 
twisted bimodule structures ${}_{}^{\kappa_{1}\!\!}\calm{}^{\kappa_{2}}_{}$ and
${}_{}^{\kappa_{2}\!}\CC\calm{}^{\kappa_{1}}_{}$ in our construction. However, as will be
explained in Section \ref{sec:balanc-pair-bimod}, there is a canonical 4-periodicity in 
the so obtained family of bimodule categories. 
\end{Remark}


\subsection{Balancings for bimodule categories}\label{sec:bal}

A connected \defectonemanifold\ $\LL$ with $n\,{>}\,0$ marked points comes by definition 
with an $n$-tuple of bimodule categories $(\calm_i)_{1\leq i\leq n}$ such that the string
$(\calm_i^{\epsilon_i})_{1\leq i\leq n}$ is either cyclically or linearly composable. 
To be 
able to introduce the \ko-linear category associated to a \defectonemanifold\ $\LL$, one further
ingredient is needed: \emph{twisted balancings} for strings of composable bimodules categories.

\begin{Definition}
\Itemizeiii

\item[{\rm (i)}]
Let $\cala$ be a monoidal category and $\calm$ an $\cala$-bimodule. A \emph{balancing} for an 
object $m\iN\calm$ is a natural family $(\sigma \eq (\sigma_a\colon a.m\To m.a)_{a\in \cala})$ of 
isomorphisms in $\calm$ such that $\sigma_1 \eq \id_m$ 
$($up to structure isomorphisms$)$ and such that the diagram
  \be
  \begin{tikzcd}
  (a{\otimes}a')\,.\,m \ar{rr}{\sigma_{aa'}} \ar{dr}[left,yshift=-4pt]{a.\sigma_{a'}}
  &~&  m\,.\,(a{\otimes}a')
  \\
  ~& a\,.\,m\,.\,a' \ar{ur}[right,yshift=-3pt]{\sigma_a.a'} &~
  \end{tikzcd}
  \label{eq:tensor-balancing}
  \ee
commutes for all $a,a' \iN \cala$.
$($For brevity we omit the constraint morphisms of the bimodule category $\calm.)$

\item[{\rm (ii)}]
The \emph{category} $\calz_\Cala(\calm)$ \emph{of objects with balancings} for an
$\cala$-bimodule category $\calm$ has as objects 
pairs $(m,\sigma)$ consisting of an object of $\calm$ and a balancing. 
\\
The morphisms $\Hom_{\calz_\Cala(\calm)}((m,\sigma),(m',\sigma'))$ are those morphisms 
$m\,{\xrightarrow{~f\,}}\, m'$ in $\calm$ for which the diagram
  \be
  \begin{tikzcd}
  a\,.\,m \ar{r}{\sigma} \ar{d}[left]{a.f} & m\,.\,a\ar{d}{f.a}
  \\
  a\,.\,m'\ar{r}{\sigma'} & m'.\,a
  \end{tikzcd}
  \ee
commutes for all $a\iN\cala$.
\end{itemize}
\end{Definition}

\begin{Remark} 
In case $\cala$ has a duality we actually do not have to require that balancings
are \emph{iso}morphisms. Indeed, if $\calm$ is a bimodule category over 
a monoidal category $\cala$ with right duals, then any natural morphism
$\sigma$ satisfying the tensoriality \eqref{eq:tensor-balancing} for $m\iN\calm$
has a two-sided inverse given by $ \sigma_a' \,{:=}\, (a.m.\evr_a) \,{\circ}\,
(a.\sigma_{a^\vee}.a) \,{\circ}\, (\coevr_a.m.a) $.
This can be verified directly by combining \eqref{eq:tensor-balancing} 
and the naturality of $\sigma$.
\end{Remark}

\begin{Definition}\label{def:framedcenter}
\Itemizeiii

\item[{\rm (i)}]
Let $(\calm_i^{\epsilon_i})_{i=1,\ldots,n}$ be a string of $n$ cyclically com\-posa\-ble 
bimodule categories, and let $\kappa$ be an $n$-tuple of framing indices associated 
with $\epsilon$. Then the $\kappa$-\emph{framed center} $($or \emph{framed center}, 
for short$)$, denoted by
  \be 
  \cenf{  \calm_1^{\epsilon_1} \!\stackrel{\kappa^{}_{1}}\boxtimes\!
  \calm_2^{\epsilon_2} \!\stackrel{\kappa^{}_{2}}\boxtimes\!
  \calm_3^{\epsilon_3} \!\stackrel{\kappa_{3}^{}} \boxtimes \cdots
  \!\!{\stackrel{\kappa_{n-1}}\boxtimes} \!\!  
  \calm_n^{\epsilon_n} \,{\stackrel{\kappa_{n}}\boxtimes} } \,,
  \ee
is the \emph{category of twisted balancings}: Objects are objects of
$\,\bigboxtimes_i \calm_i^{\epsilon_i}$ $($cyclically com\-po\-sa\-ble$)$ 
together with a balancing for each action of a finite tensor category involved.
 \\[.2em]
For an object of $\,\bigboxtimes_i \calm_i^{\epsilon_i}$ of the form $m_1^{\epsilon_1}\boti
m_2^{\epsilon_2}\boti\cdots\boti m_{n-1}^{\epsilon_{n-1}}\boti m_n^{\epsilon_{n}}$ and for
the case of even values $\kappa_i$ and $\epsilon_{1}\,{=}\,\epsilon_{i}\,{=}\,\epsilon_{i+1}
\,{=}\,\epsilon_{n}\,{=}\,1$, the balancing consists of coherent 
bimodule isomorphisms $($which also take care of hexagon-type constraints$)$
  \be
  \bearll
  m_{i}.a \boxtimes m_{i+1} \xrightarrow{~\cong~} m_{i}\boxtimes {}^{[\kappa_i -2]}a. m_{i+1}
  \qquad \text{for}~~ i \,{\in}\, \{1,2,...\,,n{-1}\} \qquad \text{and}
  \Nxl2
  m_{1} \boxtimes \cdots \boxtimes m_{n}.a \xrightarrow{~\cong~}
  {}^{[\kappa_{n}-2]} a.m_{1} \boxtimes \cdots \boxtimes m_n \,.
  \eear 
  \label{eq:isomo-tw-center}
  \ee
In case $\epsilon_i\,{=}\,\epsilon_{i+1}\,{=}\,{-}1$, the balancing is $\cc{a.m_i} \boti
\cc{m_{i+1}} \,{\xrightarrow{~\cong~}}\, \cc{m_i}\boti \cc{m_{i+1}.{}^{[\kappa_i+2]} a}$,
and analogously if $\epsilon_{1}\,{=}\,\epsilon_{n}\,{=}\,{-}1$. If $\kappa_i$ is odd,
then we either have $\epsilon_i\,{=}\,{-}1$ and $\epsilon_{i+1}\,{=}\,1$ and
deal with an object of the form $\cc{m_{i}} \boti m_{i+1}$, and it is then
understood that the isomorphism is 
$\cc{a.m_i} \boti m_{i+1} \,{\cong}\, \cc{m_i} \boti {}^{[\kappa_{i}]}a. m_{i+1}$,
or else $\epsilon_i\,{=}\,1$ and $\epsilon_{i+1}\,{=}\,{-}1$ and the isomorphism is
$m_i.a \boti \cc{m_{i+1}} \,{\cong}\, m_i \boti \cc{ m_{i+1}. {}^{[\kappa_{i}]}a}$.
 \\[.2em]
Morphisms in the category are morphisms of
$\bigboxtimes_i \calm_i^{\epsilon_i}$ that are compatible with the balancings.
\item[{\rm (ii)}]
Similarly, for a collection $(\calm_i^{\epsilon_i})_{i=1,\ldots,n}$ of linearly composable
bimodule categories, and for $\kappa$ an $n$-tuple of integers that refines the signs
$\epsilon$, the $\kappa$-\emph{framed center} is again defined as the corresponding
category of twisted balancings, with only $n{-}1$ balancings involved.
\end{itemize}
\end{Definition}

\begin{Remarks}
\Itemizeiii

\item[(i)]
We slightly abuse notation by omitting the bracketing for objects in multiple ordinary Deligne
products. This is unproblematic because different bracketings are related by canonical 
coherent isomorphisms.

\item[(ii)]
It is sufficient to specify, as done in \eqref{eq:isomo-tw-center}, the balancings only
for objects that are of $\boxtimes$-fac\-torized form. Below we will often analogously use
$\boxtimes$-factorized objects as placeholders for generic objects.

\item[(iii)]
In the special case that $\calm \,{=}\, \cala$ is a finite tensor category,
regarded as a bimodule category over itself, a balancing is a half-braiding and 
$\cala\,{\Boti{2}} \,{=}\, \cent(\cala)$ 
is the Drinfeld center of $\cala$. This justifies the terminology ``framed center''. 
Applying the tensor product of $\cala$ to two objects $b_i \,{\in}\, \cala\,{\Boti{\kappa_i}}$
gives an object in $\cala\Boti{\kappa_1{+}\kappa_2{-}2}$.
In particular, only the category $\cent(\cala) \,{=}\, \cala\,\Boti{2}$ is monoidal,
while for general $\kappa$ there are mixed tensor products $\big(\cala\,\Boti{\kappa_1}\big)\boti
\big(\cala\,\Boti{\kappa_2}\big) \,{\to}\, \cala\Boti{\kappa_1{+}\kappa_2{-}2}$.

\item[(iv)]
Instead of ordering the bimodule
categories as $\calm_1^{\epsilon_1} \cdots \calm_n^{\epsilon_n}$ we could as well order
them as $\calm_n^{\epsilon_n} \cdots \calm_1^{\epsilon_1}$. Accordingly for each pair
of consecutive bimodules we can interpret the balancing in two ways, namely as swapping
the action of the relevant finite tensor category from $\calm_i^{\epsilon_i}$ to
$\calm_{i+1}^{\epsilon_{i+1}}$ or from $\calm_{i+1}^{\epsilon_{i+1}}$ to $\calm_i^{\epsilon_i}$
corresponding, respectively, to the two pictures
    \def\locpa  {0.8}    
    \def\locpA  {2.3}    
    \def\locpb  {0.95}   
    \def\locpB  {1.67}   
    \def\locpc  {2.2pt}  
    \def\locpd  {0.12}   
    \def\locpx  {41}     
    \def\locpX  {34}     
    \def\locpy  {325}    
    \def\locpY  {332}    
  \be
  \raisebox{-3.9em}{\begin{tikzpicture}
  \scopeArrow{0.59}{>}
  \draw[line width=\locpc,\colorCircle] (0,\locpa) arc (90:-72:\locpa);
  \end{scope}
  \filldraw[color=\colorDefect] (\locpx:\locpa) circle (\locpd);
  \filldraw[color=\colorDefect] (\locpy:\locpa) circle (\locpd);
  \scopeArrow{0.88}{\arrowDefect}
  \draw[line width=\locpc,color=\colorDefect,postaction={decorate}]
           (\locpx:\locpa) -- (\locpx:\locpA)
           node[above=-3pt,color=\colorDefect] {$~\calm_\Cala$};
  \draw[line width=\locpc,color=\colorDefect,postaction={decorate}]
           (\locpy:\locpa) -- (\locpy:\locpA)
           node[below=-6pt,color=\colorDefect] {$~~{}_\cala\caln$};
  \end{scope}
  \draw[\colorHol,thick,->] plot[smooth,tension=0.55] coordinates{
          (\locpX:\locpB) (\locpX:1.1*\locpb) (\locpX-16:\locpb)
          (0:\locpb) (\locpY+16:\locpb) (\locpY:1.1*\locpb) (\locpY:\locpB) };
  \begin{scope}[shift={(2*\locpA-0.2,0)}] \node at (0,0) {and}; \end{scope}
  \begin{scope}[shift={(2*\locpA+1.7,0)}]
  \scopeArrow{0.59}{>}
  \draw[line width=\locpc,\colorCircle] (0,\locpa) arc (90:-72:\locpa);
  \end{scope}
  \filldraw[color=\colorDefect] (\locpx:\locpa) circle (\locpd);
  \filldraw[color=\colorDefect] (\locpy:\locpa) circle (\locpd);
  \scopeArrow{0.88}{\arrowDefect}
  \draw[line width=\locpc,color=\colorDefect,postaction={decorate}]
           (\locpx:\locpa) -- (\locpx:\locpA)
           node[above=-3pt,color=\colorDefect] {$~\calm_\Cala$};
  \draw[line width=\locpc,color=\colorDefect,postaction={decorate}]
           (\locpy:\locpa) -- (\locpy:\locpA)
           node[below=-6pt,color=\colorDefect] {$~~{}_\cala\caln$};
  \end{scope}
  \draw[\colorHol,thick,<-] plot[smooth,tension=0.55] coordinates{
          (\locpX:\locpB) (\locpX:1.1*\locpb) (\locpX-16:\locpb)
          (0:\locpb) (\locpY+16:\locpb) (\locpY:1.1*\locpb) (\locpY:\locpB) };
  \end{scope}
  \end{tikzpicture}}
  \label{pic:twowaysofbalancing}
  \ee
(Here we also indicate defect lines attached at the defect points, in order to indicate how
the gluing segment may appear as part of the boundary of a \defectsurface. This will be done
analogously also in other pictures.)
The orientation of the \defectonemanifold s provides us with one particular ordering, e.g.\
we swap from $\calm_i^{\epsilon_i}$ to $\calm_{i+1}^{\epsilon_{i+1}}$
in the situation shown in the picture \eqref{pic:L-eps-kappa-123} below.

\item[(v)]
By direct calculation one checks that there is an equivalence
  \be
  \cenf{(\calm \,{\stackrel\kappa\boxtimes}\, \caln)}\Opp \,\simeq\,
  \cenf{\CC\caln \,{\stackrel{-\kappa}\boxtimes}\, \CC\calm}
  \label{eq:MkNopp=N-kM}
  \ee
for any pair of a right $\cala$-module $\calm$ and left $\cala$-module $\caln$ and 
any index $\kappa$.  
\end{itemize}
\end{Remarks}

The framed center for the case of a single bimodule will be particularly relevant,
so we give it a separate name:

\begin{Definition}\label{def:Znkappa}
Let $\cala$ be a finite tensor category and $\calm$ a finite $\cala$-bimodule. For
$\kappa \,{\in}$
               $2\Z$, the $\kappa$-\emph{twisted center}
$\cent_{}^\kappa(\calm) \,{\equiv}\, \cent_\Cala^\kappa(\calm)$ is the category that has as
objects pairs $(m,\sigma)$ consisting of an object $m \,{\in}\, \calm$ and a twisted balancing
$\sigma\,{=}\,(\sigma_a)$ with $\sigma_a\colon a.m\,{\xrightarrow{\,\cong\,}}\,m.a^{[\kappa-2]}$, 
i.e.\ $\cent_{}^\kappa(\calm) \,{=}\, \calm \,{\stackrel\kappa\boxtimes}\,$.
\end{Definition}
   
The forgetful functor $\cent^\kappa(\calm)\,{\to}\,\calm$ is an exact functor of finite 
categories, hence it has a left and a right adjoint. It is convenient to express the adjoints
using the language of (co)ends, which we review in Appendix \ref{app:1}.
A right adjoint is given by the co-induction functor 
  \be
  \bearll
  \IZ_{[\kappa]}: & \calm \,\xrightarrow{~~}\, \cent^\kappa(\calm) \,,
  \Nxl2 &
  ~ m \,\xmapsto{~~} \int_{\!a\in\cala} a\,.\,m\,.\,a^{[\kappa-1]} \,,
  \label{eq:Fkappa}
  \eear
  \ee
and a left adjoint by the corresponding induction functor $\IZ^{[\kappa]}$,
see Corollary \ref{cor:ind-co-left}.

\begin{Remarks} \label{Remark:twisted-center}
The following statements about these categories follow directly from the definitions:
\Itemizeiii

\item[(i)]
For even $\kappa$ we have $\cent_\Cala^{\kappa+2}(\calm) \,{=}\, \cent(\calm^{\kappa})$,
i.e.\ the twisted center is the ordinary Drinfeld center for the $\cala$-module for
which the right $\cala$-action is twisted by the $\kappa{+}2$-fold dual. 

\item[(ii)]
The notation fits with the conventions for the twisted actions: for $\kappa \,{\in}\, 2 \Z$
we have
  \be
  \cenf{\caln_{\Cala}\, {\stackrel {\kappa+2} \boxtimes}\! {}_{\,\cala}\calm}
  = \cent(\caln_{\Cala}^{\kappa_{\phantom:}} {\boxtimes} {}_{\,\cala}\calm)
  = \cent(\caln_{\Cala}\, {\boxtimes}\,
  {}_{\cala\,\,}^{}\!\!\!{}^{\kappa_{\phantom:}\!\!\!\!}\calm) \,.
  \ee
Similarly, for odd $\kappa$ one finds (using the notation \eqref{eq:opp-module}
  \be
  \cenf{\caln \,{\stackrel \kappa \boxtimes}\, \cc{\calm}}
  = \cent(\caln \,{\boxtimes}\, {}^{\kappa\!}\cc{\calm}) 
  \ee
for two right modules $\caln_{\Cala}$ and $\calm_{\Cala}$, and analogously
  \be
  \cenf{\cc{\caln} \,{\stackrel \kappa \boxtimes}\, \calm}
  = \cent(\cc{\caln}{}^{\kappa}_{} \,{\boxtimes}\, \calm)
  \ee
for left modules.

\item[(iii)]
Let $\cala$ be a monoidal category with a right duality. Let $\calm$ be a left
$\cala$-module and $\caln$ a right $\cala$-module, whereby their Deligne product
$\caln \boti \calm$ is an $\cala$-bimodule. Then the category
  \be
  \cenf{\caln {\stackrel0\boxtimes} \calm} \,\simeq\, \caln \boxtimes_{\cala} \calm
  \ee
is the relative Deligne product (see e.g.\ \Cite{Sect.\,2.5}{fScS} for the definition). This
category can also be realized as the category of modules over $\int^{a\in\cala}\! 
a \boti {}^{\vee\!}a$, which has a natural structure of a Frobenius algebra $A$ in
$\cala \boti \CC\cala$. In this description the universal induction functor is the induction
functor for the algebra $A$. Recall \cite{fScS} that (contrary to statements in the literature) 
$\caln \,{\boxtimes_{\cala}}\, \calm$ is \emph{not} the center of $\caln$ and $\calm$, but rather
the twisted center $\caln {\stackrel0\boxtimes} \calm $, which with our conventions has objects
$n \boti m$ equipped with balancing $n.a \boti m \,{\cong}\, n \boti a^{\vee\vee}{.}\,m\,$.

\item[(iv)]
More generally, for an $\cala$-bimodule category $\calb$, the category
$\cent^2_\Cala(\calb) \,{\simeq}\, \calb_{Z_\Cala}$ 
is the category-valued trace of the bimodule category $\calb$, see \Cite{Sect.\,3}{fScS}.

\item[(v)] 
Just like the category-valued trace of a bimodule is defined via a universal property with
respect to balanced functors \Cite{Defs.\,3.2\,\&\,2.7}{fScS}, the twisted center 
$\cent^\kappa(\calm)$ of a bimodule category $\calm$ can be characterized by the universal 
property that for any finite category $\calx$,
pre-composition with the co-induction $\IZ_{[\kappa]}$ gives a distinguished equivalence
  \be
  \begin{array}{c}
  \Funle(\cent^\kappa(\calm),\calx) \,\xrightarrow{~\simeq\,}\, \Funle^{\kappa}(\calm,\calx)
  \Nxl2
  \hspace*{2.72em} F \,\longmapsto\, F \circ \IZ_{[\kappa]}
  \eear
  \label{eq:Lexk+2=LexZk}
  \ee
with $\Funle^{\kappa}(\calm,\calx)$ the category of $\kappa$-balanced functors, i.e.\ 
(compare Definition \ref{definition:bal-functor}) functors $F\colon \calm \To \calx$ with 
coherent isomorphisms $F(c.m) \,{\cong}\, F(m.c^{[\kappa]})$ (see Proposition
\ref{Cor2prop:EW:coendzFz} for the precise statement). A similar statement holds for the 
induction functor $\IZ^{[\kappa]}$, with right exact functors that are $\kappa{-}2$-balanced.
 
\item[(vi)]
For $\calm$ a right and $\caln$ a left
$\cala$-module, it follows from the definitions that there are distinguished equivalences
  \be
  \calm^{\kappa}_{} \,\Boti{\kappa'}\, \caln \simeq\, \calm \Boti{\kappa+\kappa'} \caln
  \simeq\, \calm \,\Boti{\kappa}\, {{}^{\kappa'\!\!}}_{}\caln 
  \label{eq:[k]l=k+l}
  \ee
for any pair $\kappa$, $\kappa'$ of even integers, and similar equivalences if $\kappa$ 
or $\kappa'$ are odd, for instance $\calm^{\kappa}_{} \,\Boti{\kappa'}\, \caln
\,{\simeq}\, \CC\calm \Boti{\,\kappa+\kappa'}\!\caln $
in case $\calm$ and $\caln$ are left modules, $\kappa$ is odd and $\kappa'$ even,
\end{itemize}
\end{Remarks}


\subsection{Balancing and (co)monads}\label{sec:balancing+monads}

Framed centers (which will play the role of gluing categories)
were introduced in Definition \ref{def:framedcenter}
as categories of balancings. It turns out that in order to check that these categories
have desirable features, such as cyclic invariance, it is convenient to express them
with the help of suitable (co)monads and their (co)modules.

Recall that for any monad one can consider a category of modules. If the monad is defined over 
a monoidal category and is given by tensoring with some algebra $A$, then a module 
over the monad is the same as a module over the algebra $A$.
An analogous statement is valid for comodules over a comonad 
(for some details see Appendix \ref{app:1}). Of 
particular interest to us is the central comonad, i.e.\ the endofunctor
  \be
  Z: \quad b \,\xmapsto{~~} \int_{a\in \cala}\! a\, {\otimes}\, b \,{\otimes}\, a^\vee 
  \ee
of a finite tensor category $\cala$, as well as the central monad
$ b \,{\xmapsto{~}} \int^{a\in \cala}\! a\, {\otimes}\, b \,{\otimes}\, \Vee a$
(these (co)ends exist, see \Cite{Thm.\,3.4}{shimi7}).
The Drinfeld center $\calz(\cala)$ is canonically equivalent \cite{brVi5} to the category of
modules over the central monad and to the category of comodules over the central comonad.
This description is based on the observation that if $\cala$ is a monoidal
category with right dualities, then a natural family of isomorphisms
$x.a\,{\to}\, a.x$ amounts to a dinatural family of morphisms $a^\vee\!.x.a\,{\to}\, x$.

The construction generalizes to $\cala$-bimodule categories $\calm$ as follows. Since all 
categories involved are finite, the end $Z_\Cala(m) \,{:=} \int_{a\in \cala} a.\,m\,.\,a^\vee$
and coend $\int^{a\in \cala}\! a .\, m \,. \Vee a $ exist for every $m\iN\calm$ and provide
endofunctors of $\calm$ that have a natural structure of a comonad and a monad, respectively. 
It is natural to generalize this construction further by allowing for \emph{twisted} balancings:
Thus let $\cala$ be a finite tensor category and $\calm$ a finite $\cala$-bimodule category. 
Then for any $\kappa \,{\in}\, 2\Z$ we have a comonad on $\calm$ given by the endofunctor
  \be
  Z_{[\kappa]}:\quad  m \,\xmapsto{~~} \int_{a\in\cala} a\,.\,m\,.\,a^{[\kappa-1]} \,,
  \label{eq:def:comonadZn}
  \ee
such that $Z_{[2]} \,{=}\, Z_{\Cala}$, and a monad given by the endofunctor
  \be
  Z^{[\kappa]}:\quad  m \,\xmapsto{~~} \int^{a\in\cala}\! a\,.\,m\,.\,a^{[\kappa-3]} \,.
  \label{eq:def:monadZn}
  \ee

As is shown in Corollary \ref{cor:ind-co-left}(ii) in the Appendix,
the categories of $Z_{[\kappa]}$-comodules and of $Z^{[\kappa]}$-modules are both equivalent
to the $\kappa$-twisted center $\cent_\Cala^\kappa(\calm)$ introduced in Definition 
\ref{def:Znkappa}.  In particular, the category $\calm^{Z_\Cala}$ of $Z_\Cala$-comodules 
is equivalent to the category $\cent_\Cala(\calm)$ of objects with balancings.

Now on the Deligne product that arises for the gluing category assigned to a connected 
$2$-fra\-med one-dimensional defect manifold $\LL$ we have to deal with \emph{several}
balancings, each of which is of the form \eqref{eq:isomo-tw-center}, namely one for each
consecutive pair of bimodule categories in the relevant string of composable bimodule categories.
It is not hard to see that the associated comonads (as well as the corresponding monads)
\emph{commute}, meaning that there are distributive laws, i.e.\ natural transformations
$\ell\colon Z'\,{\circ}\, Z \,{\Rightarrow}\, Z \,{\circ}\, Z'$ compatible with the comonad
(respectively monad) structure. (This amounts to two triangle and two pentagon identities. The
corresponding structural data are either trivial or induced from the structure morphisms of the
bimodule category.) Now for $Z$ and $Z'$ comonads, a distributive law endows the endofunctor 
$Z\,{\circ}\, Z'$ again with the structure of a comonad, with comultiplication
  \be
  Z\circ Z' \xrightarrow{~\,\Delta\circ\Delta'~}
  Z\circ Z\circ Z'\circ Z' \xrightarrow{~\,Z\circ \ell\circ Z'~} Z\circ Z'\circ Z\circ Z'
  \ee
and with counit $Z \,{\circ}\, Z' \,{\xrightarrow{~\eps\circ\eps'\,}}\, \id $, where
$\Delta$, $\Delta'$ and $\eps$, $\eps'$ are the coproduct and counit of $Z$ and $Z'$, 
respectively. Since the Deligne product of finite categories is symmetric, the independence
of the category $\ZZ(\LL,\epsilon,\sigma)$ on the linear order chosen is thus obvious.

\medskip

Let us describe more explicitly how we can write the framed center as a category of comodules
over commuting comonads, restricting for simplicity to the case of the framed center 
$\cenf{ \calm_1^{\epsilon_1} \!\stackrel{\kappa^{}_{1}}\boxtimes\!\calm_2^{\epsilon_2}}$ 
of just two categories. Denoting the (left or right)
actions of $\cala$ on these categories by $\act_{1}\colon \cala \,{\to}\, \End(\calm_{1})$
and $\act_{2}\colon \cala \,{\to}\, \End(\calm_{2})$, we have a functor 
  \be
  \act_{1}^{\epsilon_1} \boti \big(\act_2 \,{\circ}\, {}^{[f(\kappa_1)]}(-) \big)^{\epsilon_2}_{}:
  \quad \cala \boti \cc{\cala} \to \End(\calm_1^{\epsilon_1} \boti \calm_2^{\epsilon_2})
  \ee
with the function $f$ given by $f(\kappa_1) \,{=}\, \kappa_1\,{+}\,1\,{-}2\epsilon_1$ for even
$\kappa_1$ (i.e.\ for $\epsilon_1\,{=}\,\epsilon_2$) and by $f(\kappa_1) \,{=}\, \kappa_1\,{+}\,1$
for odd $\kappa_1$.
Taking the end of this functor defines the comonad 
$Z_{[\kappa_1]}$ on $\calm_1^{\epsilon_1} \boti \calm_2^{\epsilon_2}$; for even
$\kappa_1$ this is a special case of the comonad $Z_{[\kappa]}$ in \eqref{eq:def:comonadZn}.
Applying Proposition \ref{proposition:framed-center-univ}(ii) now gives

\begin{Lemma} \label{lemma:tw-center-cat-comod}
The category of comodules over the comonad $Z_{[\kappa]}$ on 
$\calm_1^{\epsilon_1} \boti \calm_2^{\epsilon_2}$ is equivalent to the framed center 
$\cenf{ \calm_1^{\epsilon_1} \!\stackrel{\kappa^{}}\boxtimes\! \calm_2^{\epsilon_2}}$.
\end{Lemma}

In this context it is worth recalling that the $\kappa$-twisted center $\cent^\kappa_{}(\calm)$
of an $\cala$-bimodule introduced in Definition \ref{def:Znkappa} comes with a universal functor
$\calm \,{\to}\, \cent^\kappa(\calm)$, namely the co-induction functor \eqref{eq:Fkappa}
which is right adjoint to the forgetful functor $\cent^\kappa(\calm)\,{\to}\,\calm$.
The comonad associated with this adjunction is precisely $Z_{[\kappa]}$, and analogously
the left adjoint corresponds to the monad $Z^{[\kappa]}$.
Moreover (see part (iii) of Proposition \ref{Cor2prop:EW:coendzFz}), for $\cala$ a
finite tensor category and $\calm$ an $\cala$-bimodule category we have an isomorphism
  \be
  \int^{z\in \cent^\kappa(\calm)}\! \cc z\boxtimes z
  \,\cong \int^{m\in\calm}\! \cc m\boxtimes Z_{[\kappa]}(m) 
  \label{eq:zkappa-zkappa}
  \ee
of objects in the category $\CC{\cent^\kappa(\calm)}\boti \cent^\kappa(\calm)$.
In particular, the object on the right hand side of \eqref{eq:zkappa-zkappa} is canonically 
an object in this category.


\subsection{Gluing categories}

We have now provided all algebraic ingredients needed for introducing the gluing categories.

\begin{Definition} \label{def:T(L,e,k)}
The \emph{gluing category} $\ZZ(\LL)$ that is assigned to a \defectonemanifold\
$\LL \,{=}\, (\L,\epsilon,\chi,\{\calm_i\})$ is defined as follows.
\Itemizeiii

\item[{\rm (i)}]
If the one-manifold $\L$ underlying $\LL$ is connected, as in 
    \def\locpa  {1.45}   
    \def\locpA  {2.16}   
    \def\locpb  {0.71}   
    \def\locpc  {2.2pt}  
    \def\locpd  {0.12}   
  \be
  \raisebox{-4.7em}{\begin{tikzpicture}
  \scopeArrow{0.34}{<} \scopeArrow{0.58}{<} \scopeArrow{0.85}{<}
  \draw[line width=\locpc,color=\colorCircle,postaction={decorate}]
       node[left=2.3cm,color=black] { $\LL ~~=$} (0,0) circle (\locpa);
  \end{scope} \end{scope} \end{scope}
  \draw[line width=\locpc,color=white,dashed] (0.866*\locpa,-0.5*\locpa) arc (-30:35:\locpa);
  \filldraw[color=\colorDefect]
       (0,\locpa) circle (\locpd) node[above]{$~ \calm_3^{\epsilon_3} $} 
       (-\locpa,0) circle (\locpd) node[right,yshift=-1pt]{$\calm_2^{\epsilon_2}$} 
       (0,-\locpa) circle (\locpd) node[below=1pt]{$~\calm_1^{\epsilon_1}$}
       (0.7*\locpa,0.79*\locpa) node[above]{$\kappa_3$}
       (-\locpa,-\locpa) node[above,xshift=4pt]{$\kappa_1$}
       (-\locpa,\locpa) node[below=4pt]{$\kappa_2$};
  \end{tikzpicture}~~
  \label{pic:L-eps-kappa-123}}
  \ee
then the gluing category is the corresponding $\kappa$-framed center introduced in
Definition $\ref{def:framedcenter}:$
  \be
  \ZZ(\LL) := \cenf{ \calm_1^{\epsilon_1} \,{\Boti{\kappa^{}_1}}\,
  \calm_2^{\epsilon_2} \,{\Boti{\kappa^{}_2}}\,
  \calm_3^{\epsilon_3} \,{\Boti{\kappa^{}_3}}\, \cdots } \,.
  \label{eq:TL=Z}
  \ee

\item[{\rm (ii)}]
If $\,\LL$ is a disjoint union of connected \defectonemanifold s $\LL_i$, then $\ZZ(\LL)$
is the Deligne product $\bigboxtimes_i \ZZ(\LL_i)$ of the $\kappa$-framed centers $\ZZ(\LL_i)$.
\end{itemize}
\end{Definition}

Implicitly, the prescription for the gluing categories assigned to
gluing intervals is completely determined by the prescription for gluing circles: 
just regard left and right $\cala$-modules as $\cala$-$\vect$- and as
$\vect$-$\cala$-bimodules, respectively. In the sequel we will tacitly make this
identification whenever convenient.
 \\
As an illustration, consider the following situations.

\begin{Example}\label{exa39}
Consider the following \defectonemanifold s $\LL$ and $\LL'$:
    \def\locpa  {1.45}   
    \def\locpb  {0.71}   
    \def\locpB  {0.78}   
    \def\locpc  {2.2pt}  
    \def\locpd  {0.12}   
    \def\locpg  {3.0}    
  \begin{eqnarray}
  \begin{tikzpicture}
  \scopeArrow{0.16}{<} \scopeArrow{0.64}{<}
  \draw[line width=\locpc,color=\colorCircle,postaction={decorate}]
       node[left=1.9cm,color=black] {$\LL ~=$} (0,0) circle (\locpa);
  \end{scope} \end{scope}
  \filldraw[color=\colorDefect]
       (0,\locpa)  circle (\locpd) node[above=1pt]{$+$}
                                   node[below=1pt,xshift=-2pt]{${}^{}_\cala\caln_\calb^{}$}
       (0,-\locpa) circle (\locpd) node[above=-1pt]{$-$}
                                   node[below,xshift=-2pt]{${}^{}_\cala\calm_\calb^{}$};
  \foreach \xx in {0,20,...,360}
           \draw[thick,color=\colorFrame,->] (\xx:\locpa) -- +(0,\locpb);
  \begin{scope}[shift={(2*\locpa+2.3,0)}]
  \scopeArrow{0.427}{>} \scopeArrow{0.747}{>}
  \draw[line width=\locpc,color=\colorCircle,postaction={decorate}]
       node[left=\locpE,color=black]{$\LL'~=$} (0,0) -- (2*\locpg,0);
  \end{scope} \end{scope}
  \draw[line width=\locpc,color=\colorCircle] (-\locpd,0) -- (\locpg+\locpd,0);
  \foreach \xx in {0,1,...,4} \draw[thick,color=\colorFrame,->]
             (0.06*\locpg*\xx cm,0) -- +(90-18*\xx:\locpB cm);
  \draw[thick,color=\colorFrame,->] (0.3*\locpg cm,0) -- +(7:\locpB cm);
  \draw[thick,color=\colorFrame,->] (0.42*\locpg cm,0) -- +(-24:\locpB cm);
  \foreach \xx in {7,8,9,10} \draw[thick,color=\colorFrame,->]
             (-0.4*\locpg cm+0.14*\locpg*\xx cm,0) -- +(90-18*\xx:\locpB cm);
  \begin{scope}[shift={(\locpg,0)}]
  \foreach \xx in {1,2,...,10} \draw[thick,color=\colorFrame,->]
                  (0.1*\locpg*\xx cm,0) -- (0.1*\locpg*\xx cm,-\locpB);
  \end{scope}
  \filldraw[color=\colorDefect]
       (0,0) circle (\locpd) node[below=1pt]{$+$} node[below=13pt]{$\calk_\Cala^{}$} 
       (\locpg,0) circle (\locpd) node[above=-1pt]{$-$}
                                  node[above=8pt,xshift=-1pt]{${}^{}_\calb\calm_\Cala^{}$}
       (2*\locpg,0) circle (\locpd) node[above=-1pt]{$-$}
                                  node[above=8pt,xshift=1pt]{$\caln^{}_\calb$};
  \end{scope}
  \end{tikzpicture}
  \nonumber \\[-3.1em]~
  \end{eqnarray}
  ~\\[-0.2em]
Here the two intervals of $\LL$ both have index $+1$.
Accordingly we associate to $\LL$ the category
  \be
  \ZZ(\LL) = \cenf{ \CC{{}^{}_\cala\calm_\calb^{}}
  \,\,{\stackrel{1}\boxtimes} {}^{}_{\,\cala}\caln_\Calb^{} \stackrel1\boxtimes } \,.
  \ee
Thus we deal with the two balancings $ \cc {a.m}\boti n\cong \cc m\boti {}^{\vee\!}a.n $ and
$ \cc m \boti n.b\simeq \cc{m.{}^{\vee}b}\boti n $. In the case of $\LL'$, 
the first segment has index $-1$ and we have
  \be
  \ZZ(\LL') = \cenf{ \calk_\Cala^{} \,{\stackrel{-1}\boxtimes}\,
  \CC{{}^{}_\calb\calm_\Cala^{}}  \,{\stackrel{0}\boxtimes}\, \CC{\caln^{}_\calb} } \,.
  \ee
\end{Example}

In view of Lemma \ref{lemma:tw-center-cat-comod} our prescription for associating categories
to \defectonemanifold s can be concisely summarized as follows: 
\Itemize
\item
Consider the Deligne product of the (bi)module categories involved in a composable string. 
\item
On the so obtained category there is a comonad, specified by the values of the indices,
whose comodules describe the relevant balancings.
Take the category of comodules over this comonad.
\end{itemize}

\medskip

The \defectonemanifold\ $\LL$ in Example \ref{exa39} and its variant with reversed orientation
actually play a fundamental role: they arise in particular when one modifies a 
\defectsurface\ locally by replacing a piece of line defect by a new gluing circle that is 
regarded as an incoming or outgoing boundary circle, respectively; 
pictorially, in the case of $\LL$ we have
    \def\locpa  {1.45}   
    \def\locpb  {0.5}    
    \def\locpB  {0.2}    
    \def\locpc  {2.2pt}  
    \def\locpd  {0.12}   
    \def\locpe  {70}     
    \def\locpf  {0.7}    
    \def\locpg  {3.7}    
    \def\locph  {0.61}   
  \be
  \raisebox{-4.1em}{\begin{tikzpicture}[scale=1.1]
  [decoration={random steps,segment length=2mm}]
  \fill[\colorSigma,decorate]
       (-3*\locpb,0) -- (3*\locpb,0) -- ++(\locpe:\locpg) -- ++(-6*\locpb,0) -- cycle;
  \scopeArrow{0.86}{\arrowDefect}
  \draw[line width=\locpc,color=\colorDefect,postaction={decorate}] (0,0) -- (\locpe:\locpg);
  \end{scope}
  \foreach \xx in {-2,-1,1,2} \foreach \yy in {0,1,2,3} \draw[thick,color=\colorFrame,->]
       (\xx*\locpb+1.8*\yy*\locpB,\yy*\locpf+\yy*\locpB) -- +(\locpe:\locpf);
  \begin{scope}[shift={(6*\locpb+2.7,0)}]
  [decoration={random steps,segment length=2mm}]
  \fill[\colorSigma,decorate]
       (-3*\locpb,0) -- (3*\locpb,0) -- ++(\locpe:\locpg) -- ++(-6*\locpb,0) -- cycle;
  \scopeArrow{0.17}{\arrowDefect} \scopeArrow{0.91}{\arrowDefect}
  \draw[line width=\locpc,color=\colorDefect,postaction={decorate}] (0,0) -- (\locpe:\locpg);
  \end{scope} \end{scope}
  \foreach \xx in {-2,-1,1,2} \foreach \yy in {0,1,2,3} \draw[thick,color=\colorFrame,->]
       (\xx*\locpb+1.8*\yy*\locpB,\yy*\locpf+\yy*\locpB) -- +(\locpe:\locpf);
  \scopeArrow{0.41}{<} \scopeArrow{0.95}{<}
  \filldraw[line width=\locpc,\colorCircle,fill=white,postaction={decorate}]
       (\locpe:2.59*\locph) node[left=2.3cm,yshift=7pt,color=black] {{\Large$\rightsquigarrow$}}
       circle (\locph);
  \end{scope} \end{scope}
  \filldraw[color=\colorDefect] (\locpe:1.59*\locph) circle (\locpd)
                                (\locpe:3.59*\locph) circle (\locpd);
  \end{scope}
  \end{tikzpicture}}
  \label{eq:newcircle-in-line}
  \ee
In this situation the new gluing circle inherits a $2$-framing, and this is precisely
the one of $\LL$ in Example \ref{exa39}. Similarly, when a defect line ends at a gluing 
circle with a single defect point and the vector field is analogous to the one in 
\eqref{eq:newcircle-in-line}, the $2$-framing of that gluing circle (when suitably oriented) 
has index $2$.  This observation justifies the following convention, to be used in the 
sequel throughout: \emph{If to a gluing circle with a single defect point no index label 
is attached, then it is meant to constitute a circle of index $+2$; any other gluing 
segment to which no index label is attached is meant to constitute a segment of index $+1$.} 
In pictures,
    \def\locpa  {0.95}   
    \def\locpc  {2.2pt}  
    \def\locpd  {0.12}   
  \be
  \raisebox{-2.6em}{\begin{tikzpicture}
  \scopeArrow{0.59}{<}
  \draw[line width=\locpc,\colorCircle,postaction={decorate}] (0,0) circle (\locpa);
  \end{scope}
  \filldraw[color=\colorDefect] (70:\locpa) circle (\locpd);
  \begin{scope}[shift={(3.5*\locpa,0)}]
  \scopeArrow{0.67}{<}
  \draw[line width=\locpc,\colorCircle,postaction={decorate}] (0,0) circle (\locpa);
  \end{scope}
  \filldraw[color=\colorDefect] (70:\locpa) circle (\locpd);
  \node at (220:1.25*\locpa) {$2$};
  \end{scope}
  \begin{scope}[shift={(8.5*\locpa,0)}]
  \scopeArrow{0.41}{<} \scopeArrow{0.95}{<}
  \draw[line width=\locpc,\colorCircle,postaction={decorate}] (0,0) circle (\locpa);
  \end{scope} \end{scope}
  \filldraw[color=\colorDefect] (70:\locpa) circle (\locpd) (250:\locpa) circle (\locpd);
  \end{scope}
  \begin{scope}[shift={(12*\locpa,0)}]
  \scopeArrow{0.41}{<} \scopeArrow{0.95}{<}
  \draw[line width=\locpc,\colorCircle,postaction={decorate}] (0,0) circle (\locpa);
  \end{scope} \end{scope}
  \filldraw[color=\colorDefect] (70:\locpa) circle (\locpd) (250:\locpa) circle (\locpd);
  \node at (-38:1.20*\locpa) {$1$};
  \node at (129:1.22*\locpa) {$1$};
  \end{scope}
  \node at (1.75*\locpa,0) {$\equiv$};
  \node at (6*\locpa,0) {and};
  \node at (10.25*\locpa,0) {$\equiv$};
  \end{tikzpicture}}~~
  \ee
	   
It remains to describe the category for the opposite, in the sense of 
{Definition} \ref{def:opposite}, of a \defectonemanifold:

\begin{Lemma}
\label{lemma:Op-circle}
The gluing category assigned to the opposite of a \defectonemanifold\ is the
opposite category, i.e.\ we have
  \be
  \ZZ(\overline{\LL}) = {\big(\ZZ(\LL)\big)}\Opp
  \ee
for any \defectonemanifold\ $\LL$.
\end{Lemma}

\begin{proof}
In view of the definition of the opposite manifold, it 
is sufficient to understand what happens for a single gluing segment (i.e.\ the situation 
displayed e.g.\ in picture \eqref{pic:twowaysofbalancing}). In that situation one deals with
the balanced product $\calm_\Cala\,{\stackrel\kappa\boxtimes}\, {}_\cala\caln$ of just
two factors, and for such a product the statement boils down to the isomorphism that was
already observed in formula \eqref{eq:MkNopp=N-kM}.
\end{proof}

\begin{Example} 
The $\kappa$-twisted center appears as a specific gluing category -- it is the category
assigned to a circle with one defect point and framing index $\kappa$, for $\kappa$ even, i.e.
    \def\locpa  {0.75}   
    \def\locpb  {0.81}   
    \def\locpc  {2.2pt}  
    \def\locpd  {0.12}   
    \def\locpe  {0.68}   
    \def\locpf  {0.32}   
  \be
  \raisebox{-2.0em}{\begin{tikzpicture}
  \scopeArrow{0.11}{<}
  \filldraw[line width=\locpc,\colorCircle,fill=white,postaction={decorate}]
           node[left=1.1 cm,yshift=0.13cm,color=black] {$\ZZ\Big($}
           (0,0) circle (\locpa) node[right=0.65cm,color=black]{$\kappa$}
           node[right=1.2cm,yshift=0.13cm,color=black]
	   {$\Big) ~=\, {{}_{\phantom|}\cent}^{\kappa}_{}({}^{}_\cala\calm_\Cala^{}) \,.$};
  \end{scope}
  \filldraw[color=\colorDefect] (0,\locpa) circle (\locpd);
  \scopeArrow{0.81}{\arrowDefect}
  \draw[line width=\locpc,color=\colorDefect,postaction={decorate}] (0,\locpa) -- +(0,\locpb)
           node[above=-4pt,color=\colorDefect] {${}^{}_\cala\calm_\Cala^{}$};
  \end{scope}
  \end{tikzpicture}}
  \label{eq:kappacircle}
  \ee
  ~
\end{Example}


\subsection{Canonical equivalences between gluing categories}

The following considerations can be used to considerably reduce the number of
\defectonemanifold s that we have to examine in detail: We analyze what happens when the 
orientation of a free boundary segment is flipped, and when two neighboring defect points 
are fused. Concerning the former issue we have

\begin{Proposition} \label{prop:lineflip}
Let $\calm$ and $\calk$ be left modules over a finite tensor category $\cala$. Let $\caln$ be
a right $\cala$-module, and denote by ${}^{\#\!}\caln \,{\equiv}\, {}_{}^{1}\cc{\caln}$ the
left $\cala$-module with 1-twisted left action $($in the convention of \eqref{eq:opp-module}$)$.
Up to canonical equivalence, the gluing category associated to a \defectonemanifold\
does not change if the orientation of a free boundary segment labeled by $\caln$ is flipped, 
in the sense that $\caln$ is replaced by ${}^{\#\!}\caln$ and simultaneously the orientation
of the boundary segment is inverted, corresponding to locally replacing the situation
    \def\locpa  {0.52}   
    \def\locpb  {1.4}    
    \def\locpc  {2.2pt}  
    \def\locpd  {0.12}   
    \def\locpe  {1.4}    
  \be
  \raisebox{-8.5em}{\begin{tikzpicture}
  [scale=1.2,decoration={random steps,segment length=2mm}]
  \fill[\colorSigma,decorate] (0,-\locpa-\locpb) rectangle (-\locpe,3*\locpa+2*\locpb);
  \fill[white] (0,0) circle (\locpa) (0,2*\locpa+\locpb) circle (\locpa);
  \scopeArrow{0.48}{>}
  \draw[line width=\locpc,\colorCircle,postaction={decorate}]
           (0,-\locpa) arc (270:90:\locpa)
           node[left=11pt,yshift=-5pt,color=black]{$\kappa$};
  \draw[line width=\locpc,\colorCircle,postaction={decorate}]
           (0,\locpa+\locpb) arc (270:90:\locpa)
           node[left=9pt,yshift=-5pt,color=black]{$\kappa'$};
  \end{scope}
  \scopeArrow{0.77}{\arrowDefect}
  \draw[line width=\locpc,color=\colorFreebdy,postaction={decorate}]
           (0,3*\locpa+\locpb) -- node[right=-2pt,yshift=-6pt]{$\calk$} +(0,\locpb);
  \draw[line width=\locpc,color=\colorFreebdy,postaction={decorate}]
           (0,-\locpa-\locpb) -- node[right=-2pt,yshift=-6pt]{$\calm$} +(0,\locpb);
  \end{scope}
  \scopeArrow{0.52}{\arrowDefectB}
  \draw[line width=\locpc,color=\colorFreebdy,postaction={decorate}]
           (0,\locpa) -- node[right=-2pt,yshift=9pt]{$\caln$} +(0,\locpb);
  \end{scope}
  \filldraw[color=\colorDefect] (0,\locpa) circle (\locpd) (0,-\locpa) circle (\locpd)
                                (0,\locpa+\locpb) circle (\locpd)
                                (0,3*\locpa+\locpb) circle (\locpd);
  \begin{scope}[shift={(1.9,\locpa+0.47*\locpb)}] \node at (0,0) {by}; \end{scope}
  \begin{scope}[shift={(4.8,0)}]
  [decoration={random steps,segment length=2mm}]
  \fill[\colorSigma,decorate] (0,-\locpa-\locpb) rectangle (-\locpe,3*\locpa+2*\locpb);
  \fill[white] (0,0) circle (\locpa) (0,2*\locpa+\locpb) circle (\locpa);
  \scopeArrow{0.48}{>}
  \draw[line width=\locpc,\colorCircle,postaction={decorate}]
           (0,-\locpa) arc (270:90:\locpa)
           node[left=10pt,yshift=-5pt,color=black]{$\kappa{+}1$};
  \draw[line width=\locpc,\colorCircle,postaction={decorate}]
           (0,\locpa+\locpb) arc (270:90:\locpa)
           node[left=10pt,yshift=-5pt,color=black]{$\kappa'{-}1$};
  \end{scope}
  \scopeArrow{0.77}{\arrowDefect}
  \draw[line width=\locpc,color=\colorFreebdy,postaction={decorate}]
           (0,3*\locpa+\locpb) -- node[right=-2pt,yshift=-6pt]{$\calk$} +(0,\locpb);
  \draw[line width=\locpc,color=\colorFreebdy,postaction={decorate}]
           (0,\locpa) -- node[right=-2pt,yshift=-6pt]{${}^{\#\!}\caln$} +(0,\locpb);
  \draw[line width=\locpc,color=\colorFreebdy,postaction={decorate}]
           (0,-\locpa-\locpb) -- node[right=-2pt,yshift=-6pt]{$\calm$} +(0,\locpb);
  \end{scope}
  \filldraw[color=\colorDefect] (0,\locpa) circle (\locpd) (0,-\locpa) circle (\locpd)
                  (0,\locpa+\locpb) circle (\locpd) (0,3*\locpa+\locpb) circle (\locpd);
  \end{scope}
  \end{tikzpicture}}
  \label{eq:lineflip}
  \ee
\end{Proposition}

\begin{proof}
By definition, the action of $\cala$ on the left $\cala$-module ${}^{\#\!}\caln$ is
determined through its action on $\caln$ by $a \,.\, \cc{n} \,{:=}\, \cc{n.{}^{\vee\!}a}$.
For the lower gluing segment in \eqref{eq:lineflip} we thus have an equivalence
  \be
  \cenf{\CC{{}_\cala \calm} \,{\stackrel{\kappa}\boxtimes}\, \CC{\caln_\cala}} \simeq
  \cenf{\CC{{}_\cala \calm}\,{\stackrel{\kappa+1}\boxtimes} {}_{\cala}^{~1}\cc{\caln} }\,{:}
  \ee
in both categories the balancing is
  \be
  \cc {a.x_\calm}\boti \cc{x_\caln}\cong \cc {x_\calm}\boti {}^{[\kappa+1]}a. \cc{x_\caln}
  = \cc{ x_\calm}\boti \cc{x_\caln. {}^{[\kappa+2]}a} \,.
  \ee
Similarly there is an equivalence
  \be
  \caln_\cala \,{\stackrel{\kappa'}\boxtimes}\, {}_\cala \calk
  \,{\simeq}\, \cc{ {}_{\cala}^{~1}\cc{\caln} \,{\Boti{\kappa'-1}}\! {}_\cala \calk}
  \ee
for the upper gluing segment.
\end{proof}
  
As an illustration, the following picture shows the framing with indices $\kappa\,{=}\,0$
and $\kappa'\,{=}\,2$ on the left hand side of \eqref{eq:lineflip} that results
in the straight framing on the right hand side:
    \def\locpa  {1.02}   
    \def\locpb  {1.6}    
    \def\locpc  {2.2pt}  
    \def\locpC  {1.2pt}  
    \def\locpd  {0.12}   
  \be
  \raisebox{-10.3em}{\begin{tikzpicture}
  \fill[\colorSigma] (0,3*\locpa+2*\locpb) -- (0,-\locpa-\locpb) -- (-3.5*\locpb,-\locpa-\locpb)
          -- (-2.8*\locpb,3*\locpa+2*\locpb) -- cycle;
  \begin{scope}[decoration={random steps,segment length=5mm}]
  \fill[\colorSigma,decorate] (-2*\locpb,3*\locpa+2*\locpb) -- (-2*\locpb,-\locpa-\locpb)
	  -- (-3.8*\locpb,-\locpa-\locpb) -- (-3.8*\locpb,1.6*\locpb)
          .. controls (-3.7*\locpb,2*\locpb)
          and (-3.2*\locpb,2.9*\locpb) .. (-3.2*\locpb,3*\locpa+2*\locpb) -- cycle;
  \end{scope}
  \scopeArrow{0.77}{\arrowDefect}
  \draw[line width=\locpc,color=\colorFreebdy,postaction={decorate}] (0,3*\locpa+\locpb) -- +(0,\locpb);
  \draw[line width=\locpc,color=\colorFreebdy,postaction={decorate}] (0,-\locpa-\locpb) -- +(0,\locpb);
  \end{scope}
  \scopeArrow{0.77}{\arrowDefectB}
  \draw[line width=\locpc,color=\colorFreebdy,postaction={decorate}] (0,\locpa) -- +(0,\locpb);
  \end{scope}
  \begin{scope}[shift={(0,2*\locpa+\locpb)}]
  \scopeArrow{0.83}{>}
  \draw[line width=\locpC,\colorFrame,postaction={decorate}]
           (100:\locpa) [out=100,in=270] to (-0.35,\locpa+\locpb);
  \end{scope}
  \scopeArrow{0.73}{>}
  \draw[line width=\locpC,\colorFrame,postaction={decorate}]
           (115:\locpa) [out=115,in=270] to (-0.75,\locpa+\locpb);
  \end{scope}
  \scopeArrow{0.63}{>}
  \draw[line width=\locpC,\colorFrame,postaction={decorate}]
           (130:\locpa) [out=130,in=270] to (-1.15,\locpa+\locpb);
  \end{scope}
  \scopeArrow{0.57}{<}
  \draw[line width=\locpC,\colorFrame,postaction={decorate}]
           (158:\locpa) [out=208,in=45] to (-1.73*\locpb,-0.5*\locpb)
                        [out=225,in=90] to (-2.3*\locpb,-2.2*\locpb) -- (-2.3*\locpb,-3.9*\locpb);
  \end{scope}
  \scopeArrow{0.57}{<}
  \draw[line width=\locpC,\colorFrame,postaction={decorate}]
           (180:\locpa) [out=213,in=50] to (-1.51*\locpb,-0.6*\locpb) 
                        [out=230,in=90] to (-2.1*\locpb,-2.1*\locpb) -- (-2.1*\locpb,-3.9*\locpb);
  \end{scope}
  \scopeArrow{0.59}{<}
  \draw[line width=\locpC,\colorFrame,postaction={decorate}]
           (224:\locpa) [out=166,in=90] to (-1.9*\locpb,-2.3*\locpb) -- (-1.9*\locpb,-3.9*\locpb);
  \end{scope}
  \draw[line width=\locpC,\colorFrame]  
           plot [smooth] coordinates{ (-2.5*\locpb,-3.9*\locpb) (-2.5*\locpb,-1.7*\locpb)
           (-2.33*\locpb,-1.0*\locpb) (-1.88*\locpb,-0.37*\locpb)
           (-1.32*\locpb,0.03*\locpb) (158:1.24*\locpa) (141:1.28*\locpa)
           (-1.56,\locpa+0.4*\locpb) (-1.6,\locpa+\locpb) };
  \end{scope}
  \draw[line width=\locpC,\colorFrame] plot [smooth] coordinates{
           (-1.7*\locpb,-\locpb-\locpa) (-1.7*\locpb,-0.4*\locpb)
           (-1.59*\locpb,0.6*\locpb) (131:2.92*\locpa) (111:2.64*\locpa) (125:\locpa)};
  \draw[line width=\locpC,\colorFrame,->] (-1.67*\locpb,0.08)+(260:0.1) -- (-1.67*\locpb,0.08);
  \scopeArrow{0.52}{>}
  \draw[line width=\locpC,\colorFrame,postaction={decorate}]
           (0,2*\locpa+\locpb)+(254:\locpa) [out=268,in=96] to (101:\locpa);
  \end{scope}
  \scopeArrow{0.66}{>}
  \draw[line width=\locpC,\colorFrame,postaction={decorate}]
           (0,2*\locpa+\locpb)+(237:\locpa) [out=279,in=100] to (112:\locpa); 
  \end{scope}
  \scopeArrow{0.66}{<}
  \draw[line width=\locpC,\colorFrame,postaction={decorate}]
           (140:\locpa) [out=110,in=0]  to (-0.9*\locpb,1.3*\locpb)
                        [out=180,in=70] to (-1.2*\locpb,\locpb)
                        [out=250,in=90] to (-1.5*\locpb,-\locpa-\locpb);
  \end{scope}
  \scopeArrow{0.64}{<}
  \draw[line width=\locpC,\colorFrame,postaction={decorate}]
           (158:\locpa) [out=125,in=0]  to (-0.91*\locpb,0.9*\locpb)
                        [out=180,in=72] to (-1.15*\locpb,0.5*\locpb)
                        [out=255,in=90] to (-1.32*\locpb,-\locpa-\locpb);
  \end{scope}
  \scopeArrow{0.60}{<}
  \draw[line width=\locpC,\colorFrame,postaction={decorate}]
           (177:\locpa) [out=140,in=0]  to (-0.91*\locpb,0.4*\locpb)
                        [out=180,in=75] to (-1.11*\locpb,0.1*\locpb)
                        [out=260,in=90] to (-1.14*\locpb,-\locpa-\locpb);
  \end{scope}
  \scopeArrow{0.59}{<}
  \draw[line width=\locpC,\colorFrame,postaction={decorate}]
           (196:\locpa) [out=184,in=0]  to (-0.88*\locpb,0.05*\locpb)
                        [out=180,in=75] to (-0.99*\locpb,-0.25*\locpb)
                        [out=270,in=90] to (-0.93*\locpb,-\locpa-\locpb);
  \end{scope}
  \scopeArrow{0.51}{<}
  \draw[line width=\locpC,\colorFrame,postaction={decorate}]
           (214:\locpa) [out=201,in=90] to (-0.73*\locpb,-\locpa-\locpb);
  \end{scope}
  \scopeArrow{0.58}{<}
  \draw[line width=\locpC,\colorFrame,postaction={decorate}]
           (229:\locpa) [out=221,in=90] to (-0.54*\locpb,-\locpa-\locpb);
  \end{scope}
  \scopeArrow{0.66}{<}
  \draw[line width=\locpC,\colorFrame,postaction={decorate}]
           (245:\locpa) [out=237,in=90] to (-0.36*\locpb,-\locpa-\locpb);
  \end{scope}
  \scopeArrow{0.75}{<}
  \draw[line width=\locpC,\colorFrame,postaction={decorate}]
           (260:\locpa) [out=250,in=90] to (-0.18*\locpb,-\locpa-\locpb);
  \end{scope}
  \draw[line width=\locpC,\colorFrame,->]
           (-2.48*\locpb,0.77*\locpb)+(262:0.1) -- (-2.48*\locpb,0.77*\locpb);
  \draw[line width=\locpC,\colorFrame] plot [smooth] coordinates{
           (-2.7*\locpb,-\locpa-\locpb) (-2.7*\locpb,0.6*\locpb)
           (-2.49*\locpb,0.7*\locpa+\locpb) (-1.88*\locpb,1.8*\locpa+\locpb)
           (-1.68,2.55*\locpa+\locpb) (-1.48,2.85*\locpa+\locpb)
           (-2.03,3*\locpa+1.4*\locpb) (-2.1,3*\locpa+2*\locpb) };
  \draw[line width=\locpC,\colorFrame,->]
           (-2.67*\locpb,0.93*\locpb)+(265:0.1) -- (-2.67*\locpb,0.93*\locpb);
  \draw[line width=\locpC,\colorFrame] plot [smooth] coordinates{
           (-2.9*\locpb,-\locpa-\locpb) (-2.9*\locpb,0.4*\locpb)
           (-2.79*\locpb,0.6*\locpa+\locpb) (-2.22*\locpb,1.75*\locpa+\locpb)
           (-2.43,2.55*\locpa+\locpb) (-2.08,3*\locpa+\locpb)
           (-2.58,3*\locpa+1.5*\locpb) (-2.65,3*\locpa+2*\locpb) };
  \draw[line width=\locpC,\colorFrame,->]
           (-2.864*\locpb,1.1*\locpb)+(265:0.1) -- (-2.864*\locpb,1.1*\locpb);
  \draw[line width=\locpC,\colorFrame] plot [smooth] coordinates{
           (-3.1*\locpb,-\locpa-\locpb) (-3.1*\locpb,0.6*\locpb)
           (-3.0*\locpb,0.8*\locpa+\locpb) (-2.52*\locpb,1.8*\locpa+\locpb)
           (-3.13,2.55*\locpa+\locpb) (-2.78,3*\locpa+\locpb)
           (-3.15,3*\locpa+1.5*\locpb) (-3.2,3*\locpa+2*\locpb) };
  \draw[line width=\locpC,\colorFrame,->]
           (-3.06*\locpb,1.27*\locpb)+(264:0.1) -- (-3.06*\locpb,1.27*\locpb);
  \draw[line width=\locpC,\colorFrame] plot [smooth] coordinates{
           (-3.3*\locpb,-\locpa-\locpb) (-3.3*\locpb,0.6*\locpb)
           (-3.23*\locpb,0.8*\locpa+\locpb) (-2.84*\locpb,1.8*\locpa+\locpb)
           (-3.93,2.55*\locpa+\locpb) (-3.68,3*\locpa+\locpb)
           (-3.77,3*\locpa+1.6*\locpb) (-3.8,3*\locpa+2*\locpb) };
  \draw[line width=\locpC,\colorFrame,->]
           (-3.25*\locpb,1.44*\locpb)+(261:0.1) -- (-3.25*\locpb,1.44*\locpb);
  \draw[line width=\locpC,\colorFrame] plot [smooth] coordinates{
           (-3.5*\locpb,-\locpa-\locpb) (-3.5*\locpb,0.8*\locpb)
           (-3.45*\locpb,\locpa+\locpb) (-5.1,1.8*\locpa+\locpb)
           (-4.7,2.7*\locpa+\locpb) (-4.5,3*\locpa+1.4*\locpb) (-4.5,3*\locpa+2*\locpb) };
  \draw[line width=\locpC,\colorFrame,->]
           (-3.46*\locpb,1.61*\locpb)+(261:0.1) -- (-3.46*\locpb,1.61*\locpb);
  \scopeArrow{0.40}{>}
  \filldraw[line width=\locpc,\colorCircle,fill=white,postaction={decorate}]
           (0,-\locpa) arc (270:90:\locpa);
  \end{scope}
  \scopeArrow{0.45}{>}
  \filldraw[line width=\locpc,\colorCircle,fill=white,postaction={decorate}]
           (0,\locpa+\locpb) arc (270:90:\locpa);
  \end{scope}
  \filldraw[color=\colorDefect] (0,-\locpa) circle (\locpd) (0,\locpa) circle (\locpd)
   (0,\locpa+\locpb) circle (\locpd) (0,3*\locpa+\locpb) circle (\locpd);
  \end{tikzpicture}}
  \ee ~\\[-.4em]

\begin{cor} \label{cor:bimodulelineflip}
Let $\cala$ and $\calb$ be finite tensor categories. Let $\calm$ and $\calk$ be
$\cala$-$\calb$-bimodules and $\caln$ a $\calb$-$\cala$-bimodule.
If we locally replace the situation
    \def\locpa  {0.52}   
    \def\locpb  {1.4}    
    \def\locpB  {1.2}    
    \def\locpc  {2.2pt}  
    \def\locpd  {0.12}   
    \def\locpe  {1.45}   
  \be
	\raisebox{-7.7em}{\begin{tikzpicture}[scale=1.1]
  [decoration={random steps,segment length=2mm}]
  \fill[\colorSigma,decorate] (\locpe,-\locpa-\locpB) rectangle (-\locpe,3*\locpa+\locpb+\locpB);
  \fill[white] (0,0) circle (\locpa cm);
  \fill[white] (0,2*\locpa+\locpb) circle (\locpa cm);
  \scopeArrow{0.48}{>}
  \draw[line width=\locpc,\colorCircle,postaction={decorate}]
           (0,-\locpa) arc (270:90:\locpa cm)
           node[left=11pt,yshift=-5pt,color=black]{$\kappa$};
  \draw[line width=\locpc,\colorCircle,postaction={decorate}]
           (0,\locpa) arc (90:-90:\locpa cm)
           node[right=9pt,yshift=6pt,color=black]{$\lambda$};
  \draw[line width=\locpc,\colorCircle,postaction={decorate}]
           (0,\locpa+\locpb) arc (270:90:\locpa cm)
           node[left=9pt,yshift=-5pt,color=black]{$\kappa'$};
  \draw[line width=\locpc,\colorCircle,postaction={decorate}]
           (0,3*\locpa+\locpb) arc (90:-90:\locpa cm)
           node[right=9pt,yshift=6pt,color=black]{$\lambda'$};
  \end{scope}
  \scopeArrow{0.67}{\arrowDefect}
  \draw[line width=\locpc,color=\colorFreebdy,postaction={decorate}]
           (0,3*\locpa+\locpb) -- node[right=-3pt,yshift=9pt]{$\calk$}
           node [left=16pt,yshift=-2pt,fill=\colorSigmaLight]{$\cala$}
           node[right=16pt,yshift=-2pt,fill=\colorSigmaLight]{$\calb$} +(0,\locpB);
  \end{scope}
  \scopeArrow{0.77}{\arrowDefect}
  \draw[line width=\locpc,color=\colorFreebdy,postaction={decorate}]
           (0,-\locpa-\locpB) -- node[right=-2pt,yshift=-6pt]{$\calm$} +(0,\locpB);
  \end{scope}
  \scopeArrow{0.52}{\arrowDefectB}
  \draw[line width=\locpc,color=\colorFreebdy,postaction={decorate}]
           (0,\locpa) -- node[right=-2pt,yshift=8pt]{$\caln$} +(0,\locpb);
  \end{scope}
  \filldraw[color=\colorDefect] (0,\locpa) circle (\locpd) (0,-\locpa) circle (\locpd)
                                (0,\locpa+\locpb) circle (\locpd)
                                (0,3*\locpa+\locpb) circle (\locpd);
  \begin{scope}[shift={(2.4,\locpa+0.5*\locpb)}] \node at (0,0) {by}; \end{scope}
  \begin{scope}[shift={(4.9,0)}]
  [decoration={random steps,segment length=2mm}]
  \fill[\colorSigma,decorate] (\locpe,-\locpa-\locpB) rectangle (-\locpe,3*\locpa+\locpb+\locpB);
  \fill[white] (0,0) circle (\locpa cm);
  \fill[white] (0,2*\locpa+\locpb) circle (\locpa cm);
  \scopeArrow{0.48}{>}
  \draw[line width=\locpc,\colorCircle,postaction={decorate}]
           (0,-\locpa) arc (270:90:\locpa cm)
           node[left=11pt,yshift=-5pt,color=black]{$\kappa{+}1$};
  \draw[line width=\locpc,\colorCircle,postaction={decorate}]
           (0,\locpa) arc (90:-90:\locpa cm)
	   node[right=9pt,yshift=6pt,color=black]{$\lambda{+}1$};
  \draw[line width=\locpc,\colorCircle,postaction={decorate}]
           (0,\locpa+\locpb) arc (270:90:\locpa cm)
           node[left=10pt,yshift=-7pt,color=black]{$\kappa'{-}1$};
  \draw[line width=\locpc,\colorCircle,postaction={decorate}]
           (0,3*\locpa+\locpb) arc (90:-90:\locpa cm)
	   node[right=9pt,yshift=6pt,color=black]{$\lambda'{-}1$};
  \end{scope}
  \scopeArrow{0.67}{\arrowDefect}
  \draw[line width=\locpc,color=\colorFreebdy,postaction={decorate}]
           (0,3*\locpa+\locpb) -- node[right=-3pt,yshift=9pt]{$\calk$}
           node [left=16pt,yshift=-2pt,fill=\colorSigmaLight]{$\cala$}
           node[right=16pt,yshift=-2pt,fill=\colorSigmaLight]{$\calb$} +(0,\locpB);
  \end{scope}
  \scopeArrow{0.77}{\arrowDefect}
  \draw[line width=\locpc,color=\colorFreebdy,postaction={decorate}]
           (0,-\locpa-\locpB) -- node[right=-2pt,yshift=-6pt]{$\calm$} +(0,\locpB);
  \draw[line width=\locpc,color=\colorFreebdy,postaction={decorate}]
           (0,\locpa) --
           node[right=-3pt,yshift=-8pt]{${}_{}^1\cc{\caln}{}_{}^{\,1}$} +(0,\locpb);
  \end{scope}
  \filldraw[color=\colorDefect] (0,\locpa) circle (\locpd) (0,-\locpa) circle (\locpd)
                  (0,\locpa+\locpb) circle (\locpd) (0,3*\locpa+\locpb) circle (\locpd);
  \end{scope}
  \end{tikzpicture}}
  \label{eq:lineflipBi}
  \ee
then up to canonical equivalence, the gluing category associated with the disjoint union of the
two \defectonemanifold s remains unchanged. 
\end{cor}

\begin{proof}
A statement analogous to Proposition \ref{prop:lineflip} holds for right modules.
Combining the two results for free boundary segments immediately gives the stated
result for the flip of the orientation of a defect is line.
\end{proof}

Also note that in \eqref{eq:lineflipBi} it is inessential that the gluing circles have
only two defect points. Indeed we can likewise replace the situation
    \def\locpa  {0.52}   
    \def\locpb  {1.2}    
    \def\locpB  {1.1}    
    \def\locpc  {2.2pt}  
    \def\locpd  {0.12}   
    \def\locpe  {1.45}   
    \def\locpf  {0.04}   
    \def\locpF  {1.31}   
  \be
  \raisebox{-7.1em}{\begin{tikzpicture}
  [scale=1.1,decoration={random steps,segment length=2mm}]
  \fill[\colorSigma,decorate] (\locpe,-\locpa-\locpB) rectangle (-\locpe,3*\locpa+\locpb+\locpB);
  \fill[white] (0,0) circle (\locpa cm);
  \fill[white] (0,2*\locpa+\locpb) circle (\locpa cm);
  \scopeArrow{0.48}{>}
  \draw[line width=\locpc,\colorCircle,postaction={decorate}]
           (0,-\locpa) arc (270:90:\locpa cm)
           node[left=11pt,yshift=-5pt,color=black]{$\kappa$};
  \draw[line width=\locpc,\colorCircle,postaction={decorate}]
           (0,\locpa) arc (90:-90:\locpa cm)
           node[right=9pt,yshift=6pt,color=black]{$\lambda$};
  \draw[line width=\locpc,\colorCircle,postaction={decorate}]
           (0,\locpa+\locpb) arc (270:90:\locpa cm)
           node[left=9pt,yshift=-5pt,color=black]{$\kappa'$};
  \draw[line width=\locpc,\colorCircle,postaction={decorate}]
           (0,3*\locpa+\locpb) arc (90:-90:\locpa cm)
           node[right=9pt,yshift=6pt,color=black]{$\lambda'$};
  \end{scope}
  \scopeArrow{0.67}{\arrowDefect}
  \begin{scope}[shift={(0,2*\locpa+\locpb)}]
  \draw[line width=\locpc,color=\colorFreebdy,postaction={decorate}]
           (115:\locpa) -- node[left=1pt,yshift=4pt]{$\calk_1$} (115:\locpa+\locpB);
  \draw[line width=\locpc,color=\colorFreebdy,postaction={decorate}]
           (65:\locpa) -- node[right=1pt,yshift=4pt]{$\calk_k$} (65:\locpa+\locpB);
  \filldraw[color=\colorDefect] (65:\locpa) circle (\locpd) (115:\locpa) circle (\locpd);
  \filldraw[color=black] (78:\locpF) circle (\locpf) (90:\locpF) circle (\locpf)
                         (102:\locpF) circle (\locpf);
  \end{scope} \end{scope}
  \scopeArrow{0.78}{\arrowDefect}
  \draw[line width=\locpc,color=\colorFreebdy,postaction={decorate}]
           (245:\locpa+\locpB) -- node[left=-1pt,yshift=-1pt]{$\calm_m$} (245:\locpa);
  \draw[line width=\locpc,color=\colorFreebdy,postaction={decorate}]
           (295:\locpa+\locpB) -- node[right=1pt,yshift=-3pt]{$\calm_1$} (295:\locpa);
  \filldraw[color=\colorDefect] (245:\locpa) circle (\locpd) (295:\locpa) circle (\locpd);
  \filldraw[color=black] (258:\locpF) circle (\locpf) (270:\locpF) circle (\locpf)
                         (282:\locpF) circle (\locpf);
  \end{scope}
  \scopeArrow{0.49}{\arrowDefectB}
  \draw[line width=\locpc,color=\colorFreebdy,postaction={decorate}]
           (0,\locpa) -- node[right=-2pt,yshift=7pt]{$\caln$} +(0,\locpb);
  \filldraw[color=\colorDefect] (0,\locpa) circle (\locpd) (0,\locpa+\locpb) circle (\locpd);
  \end{scope}
  \begin{scope}[shift={(2.4,\locpa+0.5*\locpb)}] \node at (0,0) {by}; \end{scope}
  \begin{scope}[shift={(4.9,0)}]
  [decoration={random steps,segment length=2mm}]
  \fill[\colorSigma,decorate] (\locpe,-\locpa-\locpB) rectangle (-\locpe,3*\locpa+\locpb+\locpB);
  \fill[white] (0,0) circle (\locpa cm);
  \fill[white] (0,2*\locpa+\locpb) circle (\locpa cm);
  \scopeArrow{0.48}{>}
  \draw[line width=\locpc,\colorCircle,postaction={decorate}]
           (0,-\locpa) arc (270:90:\locpa cm)
           node[left=10pt,yshift=-5pt,color=black]{$\kappa{+}1$};
  \draw[line width=\locpc,\colorCircle,postaction={decorate}]
           (0,\locpa) arc (90:-90:\locpa cm)
	   node[right=9pt,yshift=6pt,color=black]{$\lambda{+}1$};
  \draw[line width=\locpc,\colorCircle,postaction={decorate}]
           (0,\locpa+\locpb) arc (270:90:\locpa cm)
           node[left=10pt,yshift=-5pt,color=black]{$\kappa'{-}1$};
  \draw[line width=\locpc,\colorCircle,postaction={decorate}]
           (0,3*\locpa+\locpb) arc (90:-90:\locpa cm)
	   node[right=9pt,yshift=6pt,color=black]{$\lambda'{-}1$};
  \end{scope}
  \scopeArrow{0.67}{\arrowDefect}
  \begin{scope}[shift={(0,2*\locpa+\locpb)}]
  \draw[line width=\locpc,color=\colorFreebdy,postaction={decorate}]
           (115:\locpa) -- node[left=1pt,yshift=4pt]{$\calk_1$} (115:\locpa+\locpB);
  \draw[line width=\locpc,color=\colorFreebdy,postaction={decorate}]
           (65:\locpa) -- node[right=1pt,yshift=4pt]{$\calk_k$} (65:\locpa+\locpB);
  \filldraw[color=\colorDefect] (65:\locpa) circle (\locpd) (115:\locpa) circle (\locpd);
  \filldraw[color=black] (78:\locpF) circle (\locpf) (90:\locpF) circle (\locpf)
                         (102:\locpF) circle (\locpf);
  \end{scope} \end{scope}
  \scopeArrow{0.78}{\arrowDefect}
  \draw[line width=\locpc,color=\colorFreebdy,postaction={decorate}]
           (245:\locpa+\locpB) -- node[left=-1pt,yshift=-1pt]{$\calm_m$} (245:\locpa);
  \draw[line width=\locpc,color=\colorFreebdy,postaction={decorate}]
           (295:\locpa+\locpB) -- node[right=1pt,yshift=-3pt]{$\calm_1$} (295:\locpa);
  \filldraw[color=\colorDefect] (245:\locpa) circle (\locpd) (295:\locpa) circle (\locpd);
  \filldraw[color=black] (258:\locpF) circle (\locpf) (270:\locpF) circle (\locpf)
                         (282:\locpF) circle (\locpf);
  \end{scope}
  \scopeArrow{0.79}{\arrowDefect}
  \draw[line width=\locpc,color=\colorFreebdy,postaction={decorate}] (0,\locpa) --
           node[right=-3pt,yshift=-7pt]{${}_{}^1\cc{\caln}{}_{}^{\,1}$} +(0,\locpb);
  \end{scope}
  \filldraw[color=\colorDefect] (0,\locpa) circle (\locpd) (0,\locpa+\locpb) circle (\locpd);
  \end{scope}
  \end{tikzpicture}}
  \label{eq:lineflipBi2}
  \ee
for any numbers $m$ of incoming and $k$ of outgoing defect points.

Furthermore, by applying Proposition \ref{prop:lineflip} twice we see that 
the same gluing category is obtained when changing the framing in such a way that indices
$\kappa,\kappa'$ on consecutive gluing segments along a \defectonemanifold\ (as in
\eqref{eq:lineflip}) get replaced by the pair $\kappa{+}2,\kappa'{-}2$
and the category labeling the defect line gets twisted with the corresponding double dual,
thus e.g.\ replacing $\caln$ in \eqref{eq:lineflip} by $\caln^{-2}_{})$
(see also Appendix \ref{app:shift}).

\medskip 

Next we describe the effect of `fusing' two neighboring defect points on a \defectonemanifold.
In the pictures below we display -- as we already did in the picture \eqref{eq:kappacircle} --
defect lines that are attached to defect points, even though for now we only deal with 
\defectonemanifold s (this will be convenient when using our results later on to discuss
the fusion of defect lines).

\begin{Proposition} \label{Prop:fusion-at-gluing-circles}
For any triple $(i,j,k)$ of framing indices there are canonical equivalences 
    \def\locpa  {0.73}   
    \def\locpA  {1.9}    
    \def\locpb  {1.01}   
    \def\locpc  {2.2pt}  
    \def\locpd  {0.12}   
    \def\locpE  {0.88}   
    \def\locpf  {0.22}   
    \def\locpu  {90}     
    \def\locpx  {-45}    
    \def\locpy  {-135}   
    \def\locpz  {0.37}   
  \be
  \raisebox{-4.7em}{\begin{tikzpicture}
  \scopeArrow{0.01}{<} \scopeArrow{0.43}{<} \scopeArrow{0.74}{<}
  \filldraw[line width=\locpc,\colorCircle,fill=white,postaction={decorate}]
           node[left=1.7cm,color=black]{$\ZZ\Big($} (0,0) circle (\locpa)
           node[right=1.6cm,color=black]{$\Big)~~~~\simeq$}
           node[left=0.67cm,color=black]{$\scriptstyle j$}
           node[right=0.59cm,yshift=12pt,color=black]{$\scriptstyle k$}
           node[below=0.66cm,xshift=5pt,color=black]{$\scriptstyle i$};
  \end{scope} \end{scope} \end{scope}
  \filldraw[color=\colorDefect] (\locpu:\locpa) circle (\locpd)
           (\locpx:\locpa) circle (\locpd) (\locpy:\locpa) circle (\locpd);
  \scopeArrow{0.62}{\arrowDefectB}
  \draw[line width=\locpc,color=\colorDefect,postaction={decorate}]
           (\locpx:\locpa) -- (\locpx:\locpA)
           node[below=-2pt,xshift=-5pt,color=\colorDefect] {$~~~{}^{}_\calb\caln_\calc^{}$};
  \draw[line width=\locpc,color=\colorDefect,postaction={decorate}]
           (\locpy:\locpa) -- (\locpy:\locpA)
           node[below=-2pt,color=\colorDefect] {${}^{}_\cala\calm_\calb^{}~~~$};
  \end{scope}
  \scopeArrow{0.81}{\arrowDefect}
  \draw[line width=\locpc,color=\colorDefect,postaction={decorate}]
           (\locpu:\locpa cm) -- (\locpu:\locpA)
           node[above=-4pt,color=\colorDefect] {${}^{}_\cala\calk_\calc^{}$};
  \end{scope}
  \begin{scope}[shift={(5.8,\locpz)}] 
  \scopeArrow{0.47}{<} \scopeArrow{0.97}{<}
  \filldraw[line width=\locpc,\colorCircle,fill=white,postaction={decorate}]
       node[left=1.5cm,yshift=-\locpz cm,color=black]{$\ZZ\Big($}
       (0,0) circle (\locpa cm) node[right=1.5cm,yshift=-\locpz cm,color=black]{$\Big)$}
       node[left=0.64cm,yshift=-6pt,color=black]{$\scriptstyle j$}
       node[right=0.64cm,yshift=6pt,color=black]{$\scriptstyle k$};
  \end{scope} \end{scope}
  \filldraw[color=\colorDefect] (0,\locpa) circle (\locpd);
  \filldraw[color=\colorDefect] (0,-\locpa) circle (\locpd);
  \scopeArrow{0.81}{\arrowDefect}
  \draw[line width=\locpc,color=\colorDefect,postaction={decorate}] (0,\locpa) -- +(0,\locpb)
       node[above=-4pt,color=\colorDefect] {${}^{}_\cala\calk_\calc^{}$};
  \end{scope}
  \scopeArrow{0.66}{\arrowDefectB}
  \draw[line width=\locpc,color=\colorDefect,postaction={decorate}] (0,-\locpa) -- +(0,-\locpb)
       node[below=-5pt,color=\colorDefect] {$\cenf{\calm{\stackrel{-i}\boxtimes}\caln}~$};
  \end{scope}
  \end{scope}
  \end{tikzpicture}}
  \label{eq:Z(MNK)=Z(Z(MN)K)}
  \ee
and 
    \def\locpa  {0.73}   
    \def\locpA  {1.9}    
    \def\locpb  {1.02}   
    \def\locpc  {2.2pt}  
    \def\locpd  {0.12}   
    \def\locpE  {0.88}   
    \def\locpf  {0.22}   
    \def\locpu  {270}    
    \def\locpx  {135}    
    \def\locpy  {45}     
    \def\locpz  {0.37}   
  \be
  \raisebox{-6.1em}{\begin{tikzpicture}
  \scopeArrow{0.28}{<} \scopeArrow{0.51}{<} \scopeArrow{0.91}{<}
  \filldraw[line width=\locpc,\colorCircle,fill=white,postaction={decorate}]
           node[left=1.7cm,color=black]{$\ZZ\Big($} (0,0) circle (\locpa cm)
           node[right=1.6cm,color=black]{$\Big)~~~~\simeq$}
           node[left=0.61cm,yshift=-11pt,color=black]{$\scriptstyle k$}
           node[right=0.67cm,color=black]{$\scriptstyle j$}
           node[above=0.64cm,xshift=9pt,color=black]{$\scriptstyle i$};
  \end{scope} \end{scope} \end{scope}
  \filldraw[color=\colorDefect] (\locpu:\locpa) circle (\locpd);
  \filldraw[color=\colorDefect] (\locpx:\locpa) circle (\locpd);
  \filldraw[color=\colorDefect] (\locpy:\locpa) circle (\locpd);
  \scopeArrow{0.66}{\arrowDefectB}
  \draw[line width=\locpc,color=\colorDefect,postaction={decorate}]
           (\locpu:\locpa) -- (\locpu:\locpA)
           node[below=-2pt,color=\colorDefect] {${}^{}_\cala\calk_\calc^{}\,$};
  \end{scope}
  \scopeArrow{0.81}{\arrowDefect}
  \draw[line width=\locpc,color=\colorDefect,postaction={decorate}]
           (\locpx:\locpa) -- (\locpx:\locpA)
           node[above=-2pt,color=\colorDefect] {${}^{}_\cala\calm_\calb^{}~~$};
  \draw[line width=\locpc,color=\colorDefect,postaction={decorate}]
           (\locpy:\locpa) -- (\locpy:\locpA)
           node[above=-2pt,color=\colorDefect] {$~~{}^{}_\calb\caln_\calc^{}$};
  \end{scope}
  \begin{scope}[shift={(5.8,-\locpz)}] 
  \scopeArrow{0.47}{<} \scopeArrow{0.97}{<}
  \filldraw[line width=\locpc,\colorCircle,fill=white,postaction={decorate}]
       node[left=1.5cm,yshift=\locpz cm,color=black]{$\ZZ\Big($}
       (0,0) circle (\locpa cm) node[right=1.5cm,yshift=\locpz cm,color=black]{$\Big)$}
       node[left=0.64cm,yshift=-6pt,color=black]{$\scriptstyle k$}
       node[right=0.64cm,yshift=6pt,color=black]{$\scriptstyle j$};
  \end{scope} \end{scope}
  \filldraw[color=\colorDefect] (0,\locpa) circle (\locpd) (0,-\locpa) circle (\locpd);
  \scopeArrow{0.81}{\arrowDefect}
  \draw[line width=\locpc,color=\colorDefect,postaction={decorate}] (0,\locpa) -- +(0,\locpb)
           node[above=-2pt,color=\colorDefect] {$\cenf{\calm{\stackrel{i}\boxtimes}\caln}\,$};
  \end{scope}
  \scopeArrow{0.66}{\arrowDefectB}
  \draw[line width=\locpc,color=\colorDefect,postaction={decorate}] (0,-\locpa) -- +(0,-\locpb)
           node[below=-2pt,color=\colorDefect] {${}^{}_\cala\calk_\calc^{}\,$};
  \end{scope}
  \end{scope}
  \end{tikzpicture}}
  \label{eq:Z(KMN)=Z(K.Z(MN))}
  \ee
of gluing categories.
\end{Proposition}

\begin{proof}
The category on the left hand side of \eqref{eq:Z(MNK)=Z(Z(MN)K)} is
$\cenf{ \calk \,{\stackrel k\boxtimes}\, \CC\caln
\,{\stackrel i\boxtimes}\, \CC\calm \, {\stackrel j\boxtimes}}$, 
while the one on the right hand side is $\cenf{  \calk \,{\stackrel k\boxtimes}\,
\cenf{[\calm{\stackrel{-i}\boxtimes}\caln{]}\Opp} \, {\stackrel j\boxtimes}}$. We have
  \be
  \cenf{ \calk \,{\stackrel k\boxtimes}\, \CC\caln
  \,{\stackrel i\boxtimes}\, \CC\calm \, {\stackrel j\boxtimes}}
  \,\simeq\, \cenf{ \calk \,{\stackrel k\boxtimes}\, (\CC\caln \,{\stackrel
  i\boxtimes}\, \CC\calm) \, {\stackrel j\boxtimes}}
  \,\simeq\, \cenf{  \calk \,{\stackrel k\boxtimes}\,
  \cenf{[\calm{\stackrel{-i}\boxtimes}\caln{]}\Opp} \, {\stackrel j\boxtimes}} \,,
  \ee
where the first step uses the distributive law and the second step Lemma 
\ref{lemma:Op-circle}. Canonical functors realizing these equivalences are provided by
the functors which come with the universal property \eqref{eq:Lexk+2=LexZk} of the respective
twisted centers. The proof of the equivalence \eqref{eq:Z(KMN)=Z(K.Z(MN))} is analogous.
\end{proof}

These results about gluing categories tell us that the fusion of defect points provides a well
defined operation for the associated categories. They also suggest that the same holds for defect
\emph{lines}, with the fusion of two defects labeled by bimodule categories $_\calc\calm_\Cala$
and $_\cala\caln_\calb$ and separated by a strip with winding $i$ realizing a single
defect labeled by the $\calc$-$\calb$-bi\-mo\-dule $\cenf{\calm{\stackrel{-i}\boxtimes}\caln}$.
To see that this is indeed the case we will have to define the functors
assigned to \defectsurface s, not just the categories assigned to the gluing segments
on their boundary; this will be done in Section \ref{sec:fusion}.
Analogous equivalences as in Proposition \ref{Prop:fusion-at-gluing-circles} hold for
the gluing categories of circles with more than three defect points. We have for instance

\begin{Example}\label{ex:MNKfusion}
There is a canonical equivalence
    \def\locpa  {0.73}   
    \def\locpA  {1.9} 
    \def\locpb  {0.91} 
    \def\locpc  {2.2pt}  
    \def\locpd  {0.12}   
    \def\locpE  {0.88}   
    \def\locpf  {0.22}   
    \def\locpx  {54}     
    \def\locpX  {234}    
    \def\locpy  {126}    
    \def\locpY  {306}    
    \def\locpz  {0.01}   
  \be
  \raisebox{-5.5em}{\begin{tikzpicture}
  \scopeArrow{0.28}{<} \scopeArrow{0.50}{<} \scopeArrow{0.73}{<} \scopeArrow{0.98}{<}
  \filldraw[line width=\locpc,\colorCircle,fill=white,postaction={decorate}]
           node[left=1.44cm,color=black]{$\ZZ\Big($} (0,0) circle (\locpa)
           node[right=1.47cm,color=black] {$\Big)~~~~\simeq$}
           node[above=0.64cm,xshift=6pt,color=black]{$i$}
           node[below=0.61cm,xshift=5pt,color=black]{$-j$};
  \end{scope} \end{scope} \end{scope} \end{scope}
  \filldraw[color=\colorDefect] (\locpx:\locpa) circle (\locpd) (\locpy:\locpa) circle (\locpd)
           (\locpX:\locpa) circle (\locpd) (\locpY:\locpa) circle (\locpd);
  \scopeArrow{0.84}{\arrowDefect}
  \draw[line width=\locpc,color=\colorDefect,postaction={decorate}]
           (\locpx:\locpa) -- (\locpx:\locpA)
	   node[above=-2pt,color=\colorDefect] {$~~~\cenf{\caln{\stackrel{j}\boxtimes}\calk}~~$};
  \draw[line width=\locpc,color=\colorDefect,postaction={decorate}]
           (\locpy:\locpa) -- (\locpy:\locpA)
           node[above=-2pt,color=\colorDefect] {$\calm~$};
  \end{scope}
  \scopeArrow{0.62}{\arrowDefectB}
  \draw[line width=\locpc,color=\colorDefect,postaction={decorate}]
           (\locpX:\locpa) -- (\locpX:\locpA)
	   node[below=-7pt,xshift=-13pt,color=\colorDefect]
           {$~~~\cenf{\calm{\stackrel{i}\boxtimes}\caln}$};
  \draw[line width=\locpc,color=\colorDefect,postaction={decorate}]
           (\locpY:\locpa cm) -- (\locpY:\locpA cm)
           node[below=-2pt,color=\colorDefect] {$~\calk$};
  \end{scope}
  \begin{scope}[shift={(5.9,-\locpz)}] 
  \scopeArrow{0.44}{<} \scopeArrow{0.94}{<}
  \filldraw[line width=\locpc,\colorCircle,fill=white,postaction={decorate}]
           node[left=1.7cm,yshift=\locpz cm,color=black]{$\ZZ\Big($}
           (0,0) circle (\locpa) node[right=1.6cm,yshift=\locpz cm,color=black]{$\Big)$};
  \end{scope} \end{scope}
  \filldraw[color=\colorDefect] (0,\locpa) circle (\locpd) (0,-\locpa) circle (\locpd);
  \scopeArrow{0.84}{\arrowDefect}
  \draw[line width=\locpc,color=\colorDefect,postaction={decorate}] (0,\locpa) -- +(0,\locpb)
           node[above=-2pt,color=\colorDefect]
	   {$\cenf{\calm{\stackrel{i}\boxtimes}\caln{\stackrel{j}\boxtimes}\calk}\,$};
  \end{scope}
  \scopeArrow{0.68}{\arrowDefectB}
  \draw[line width=\locpc,color=\colorDefect,postaction={decorate}] (0,-\locpa) -- +(0,-\locpb)
           node[below=-7pt,color=\colorDefect]
	   {$\cenf{\calm{\stackrel{i}\boxtimes}\caln{\stackrel{j}\boxtimes}\calk}~$};
  \end{scope}
  \end{scope}
  \end{tikzpicture}}
  \ee
Both gluing categories are equivalent to the category
$\Funle_{\cala,\calb} (\cenf{\calm \,{\Boti i}\, \caln \,{\Boti j}\, \calk} \,,
\cenf{\calm \,{\Boti i}\, \caln \,{\Boti j}\, \calk})$. 
\end{Example}

As described in Appendix \ref{sec:moduleEW}, the Eilenberg--Watts equivalences that will
be formulated in the next subsection can be lifted to twisted centers.
A common feature of the gluing categories in Example \ref{ex:MNKfusion} is that
they are canonically equivalent to categories of endofunctors. Thereby each of these
twisted centers contains a distinguished object, namely the one that corresponds to the
respective identity functor. We refer to those objects as \emph{distinguished fusion objects}.


\subsection{Balanced pairings for bimodule categories} \label{sec:balanc-pair-bimod}

We denote, for finite categories $\caln$ and $\calm$, by $\Funle(\caln,\calm)$ and
$\Funre(\caln,\calm)$ the finite category of left and right exact functors, respectively,
from $\caln$ to $\calm$. It will be convenient to re-ex\-press the gluing categories
associated with \defectonemanifold s, which we have introduced via Deligne products, in
terms of such functor categories. This is achieved with the help of a Morita invariant
formulation of Eilenberg--Watts type results, which exhibits explicit equivalences from
$\Funle(\caln,\calm)$ and $\Funre(\caln,\calm)$ to $\CC\caln \boti \calm$. More
specifically, there are (two-sided) adjoint equivalences \cite{shimi7,fScS2}
  \be
  \begin{array}{lll}
  \Phile_{} \equiv \Phile_{\!\caln,\calm}:
  & \calnopp \boti \calm \xrightarrow{~\simeq~} \Funle(\caln,\calm) \,,~
  & \cc n \boti m \longmapsto \Hom_\caln(n,-) \oti m \,,
  \Nxl5
  \Psile_{} \equiv \Psile_{\!\caln,\calm}:
  & \Funle(\caln,\calm) \xrightarrow{~\simeq~} \calnopp \boti \calm \,,~
  & F \longmapsto \int^{n\in\caln} \cc n \boti F(n) \,,
  \eear
  \label{Phile.Psile}
  \ee
and
  \be
  \begin{array}{lll}
  \Phire_{} \equiv \Phire_{\!\caln,\calm}:
  & \calnopp \boti \calm \xrightarrow{~\simeq~} \Funre(\caln,\calm) \,,~
  & \cc n \boti m \longmapsto {\Hom_\caln(-,n)}^{\!*}_{} \oti m \,,
  \Nxl5
  \Psire_{} \equiv \Psire_{\!\caln,\calm}:
  & \Funre(\caln,\calm) \xrightarrow{~\simeq~} \calnopp \boti \calm \,,~
  & G \longmapsto \int_{n\in\caln} \cc n \boti G(n) \,,
  \eear
  \label{Phire...Psire}
  \ee
where $(-)^*_{}$ denotes the dual vector space. We refer to these pairs of adjoint functors
as \emph{Eilenberg--Watts functors}. These equivalences of categories 
are compatible with the structure of bimodule categories over finite tensor categories 
\Cite{Sect.\,4}{fScS2}. For instance, the (right exact) \emph{Nakayama functor} 
  \be
  \pift := \int^{m\in\calm}\! {\Hom_\calm(-,m)}^* \oti m \,,
  \label{eq:Nr}
  \ee
i.e.\ the image $\Phire{\circ}\Psile(\id_\calm)$ of the identity functor, regarded as a left 
exact functor, in $\Funre(\calm,\calm)$, is a bimodule functor from $\calm$ to 
${}_{}^{2\!}\calm_{}^2$ \Cite{Thm.\,4.5}{fScS2}, and similarly for the left exact 
Na\-ka\-ya\-ma functor $\pifu \,{=}\,\Phile{\circ}\Psire(\id_\calm)$.

For later use we collect some properties of the left exact Nakayama functor of a finite
category $\calm$. We have
  \be
  \int_{m\in\calm}\! \cc{m} \boti m = \Psire(\id_\calm)
  \cong \Psile {\circ} \Phile {\circ} \Psire(\id_\calm) \cong \Psile(\pifu) 
  = \int^{m\in\calm}\!\! \cc{m} \boti \pifu(m) \,.
  \ee
Further, in case $\calm \,{=}\, \cala$ is a finite tensor category, we can express 
$\pifu$ as \Cite{Lemma\,4.10}{fScS2} $\pifu(a) \,{=}\, \DA \oti {}^{\vee\vee\!}a$
with $\DA$ the distinguished invertible object \Cite{Def.\,3.1}{etno2} of $\cala$, and thus
  \be
  \int_{\!a\in\cala} a^{\vee\vee} \boti \cc{D_\Cala \oti a}
  \,\cong \int^{a\in\cala}\! a \boti \cc{a} 
  \label{eq:disting-iso}
  \ee
and
  \be
  \int_{a\in\cala}\! \cc{a} \boti a
  \,\cong \int^{a\in\cala}\!\! \cc{a} \boxtimes \DA \oti {}^{\vee\vee\!}a \,.
  \label{eq:end-coend-D}
  \ee
By the category-theoretic version of Radford's $S^4$-theorem,
the invertible object $\DA$ comes \Cite{Thm.\,3.3}{etno2} with
coherent isomorphisms $\DA \oti a \,{\cong}\, a^{\vee\vee\vee\vee} {\otimes}\, \DA$
and can thus be regarded canonically as an object of the $-2$-twisted center,
  \be
  \DA \in \cent_{}^{-2}(\cala) \,.
  \label{eq:D-in-Z4}
  \ee
It follows that acting with $\DA$ is an equivalence
  \be
  \DA\,.\,{-} :\quad {}^{\kappa + 4\!}\calm \xrightarrow{~\simeq~} {}_{}^{\kappa\!\!}\calm
  \label{eq:DA:kappa+4}
  \ee
of $\cala$-modules for $\kappa \iN 2 \Z$. Hence there is a distinguished 4-periodicity in the
family of bimodule categories ${}_{}^{\kappa_{1}\!\!}\calm{}^{\kappa_{2}}_{}$ indexed by
$(\kappa_{1},\kappa_{2}) \,{\in}\, 2\Z \,{\times}\, 2\Z$, and similarly for odd $\kappa$. Moreover,
  \be
  a^{[\kappa]} \otimes \DA^{\otimes n} \cong \DA^{\otimes n} \otimes a^{[\kappa-4n]}
  \label{eq:D-and-powers}
  \ee
for all $\kappa, n \,{\in}\, \Z$, and owing to invertibility of $\DA^{}$,
there is a canonical isomorphism $\DA^{\vee \vee} \,{\cong}\, \DA^{}$.

\medskip

The Eilenberg--Watts equivalences allow us to switch back and forth between Deligne products 
and categories of half-exact functors and thereby to understand features of the former type 
of categories in terms of the latter, and vice versa. One application that will turn out
to be crucial below (see e.g.\ the calculation needed in Example \ref{ex:4.4})
is the following. Every left exact module functor 
$F\colon \caln_\Cala \,{\to}\, \calm_\Cala$ yields a balanced pairing 
  \be
  \Hom_\calm(-,F(-)): \quad \cc{\calm_\Cala} \boxtimes \caln_\Cala \to \vect \,,
  \label{eq:balanced-F-pairing}
  \ee
with the balancing obtained by combining the module structure of $F$ and the duality of 
$\cala$: $\Hom(n,F(m.a)) \,{\cong}\, \Hom(n. \Vee a, F(m))$.
This balancing is coherent with respect to the monoidal structure of $\cala$. Moreover, 
since $F$ ist left exact,
the Eilenberg--Watts equivalence \eqref{Phile.Psile} yields a natural isomorphism 
  \be
  \Hom_{\cc{\caln}\boxtimes\calm} (\cc{n} \boti m, \Psile(F)) \,{\cong}\, \Hom_\caln(n,F(m)) \,,
  \ee
i.e.\ the pairing \eqref{eq:balanced-F-pairing} is representable by the object $\Psile(F)$.
The balancing of the pairing transports to a balancing of the representing object,
thus yielding the following result, which one may think of as specific kinds of
substitution of variables rules for coends and ends:

\begin{Lemma} \label{lemma:balancing-psil(F)}
For a left exact bimodule functor $F\colon \calm \,{\to}\, \caln$ between 
$\calb$-$\cala$-bimodules $\calm$ and $\caln$, the coend
$\Psile(F) \,{=}\, \int^{m\in\calm} \cc{m} \boti  F(m) \,{\in}\, \calmopp\boti\caln$
comes with coherent isomorphisms
  \be
  \bearl\dsty
  \int^{m\in\calm}\! \cc{m} \boti F(m). a 
  \,\cong \int^{m\in\calm}\! \cc{m \,. \Vee a} \boti F(m)  \qquad\text{and}
  \Nxl3\dsty
  \int^{m\in\calm} \cc{m} \boti b.F(m)
  \,\cong \int^{m\in\calm} \cc{b^\vee .\,m} \boti F(m) \,.
  \label{eq:Bala-le}
  \eear
  \ee
Analogously, for a right exact bimodule functor $\,G\colon \calm \,{\to}\, \caln$ the end
$\,\Psire(G) \eq \int_{m\in\calm} \cc{m} \,{\boxtimes}$
                       $G(m)$ is equipped with coherent isomorphisms
  \be
  \bearl\dsty
  \int_{m\in\calm}\! \cc{m} \boti G(m) .a
  \,\cong \int_{m\in\calm} \cc{m\,.\, a^\vee} \boti G(m)  \qquad\text{and}
  \Nxl4\dsty
  \int_{m\in\calm} \cc{m} \boti b. G(m) 
  \,\cong \int_{m\in\calm} \cc{\Vee b . \,m} \boti G(m) \,.
  \label{eq:Bala-re}
  \eear
  \ee
\end{Lemma}

\begin{proof}
With the help of the representing object $\Psile(F)$ we define the balancing by
the requirement that the diagram 
  \be
  \begin{tikzcd}
  \Hom_{\cc{\calm}\boxtimes\caln}(\cc{m . a} \boti n, \Psile(F)) \ar{r} \ar{d}
  & \Hom_\caln(n,F(m.a)) \ar{d}
  \\
  \Hom_{\cc{\calm}\boxtimes\caln}(\cc{m} \boti n . \Vee a, \Psile(F)) \ar{r}
  & \Hom_\caln(n. \Vee a,F(m))
  \end{tikzcd}
\label{eq:def-bal-le}
  \ee
as well as the corresponding diagram for the left action commute.
This produces directly the isomorphisms in \eqref{eq:Bala-le}. 
The case of a right exact bimodule functor is shown analogously via its representing 
property of the dual Hom functor. 
\end{proof}

{}From this statement we obtain two types of balancings for the 
identity bimodule functor. We record them for later use:
  \be
  \int_{m \in \Calm}\!\! m.a\boti \cc m \,\cong \int_{m \in \Calm}\!\! m\boti \cc{m.a^\vee}
  \quad\text{~and~}\quad 
  \int_{m \in \Calm}\!\! b.m\boti \cc m \,\cong \int_{m \in \Calm}\!\! m\boti \cc{\Vee b.m} \,.
  \label{eq:mmbalancing}
  \ee

By applying the Eilenberg--Watts equivalence $\Phile$ from \eqref{Phile.Psile} to the 
result in Lemma \ref{lemma:balancing-psil(F)} we get

\begin{cor} \label{cor:lex-framed}
For bimodule categories ${}_{\cala}\calm_{\calb}$ and ${}_{\cala}\caln_{\calb}$,
the Eilenberg--Watts equivalences induce an equivalence 
  \be
  \cenf{ \CC\calm \,{\stackrel1\boxtimes}\, \caln \,{\stackrel1\boxtimes}}\,
  \,\simeq\, \Funle_{\Cala,\calb}(\calm,\caln) \,.
  \label{eq:equ-lex-framed-cent}
  \ee
of categories.
\end{cor}

In particular, for the regular bimodule category ${}_{\cala}\cala_{\cala}$ we obtain 
  \be
  \cenf{ \CC\cala \,{\stackrel1\boxtimes}\, \cala \,{\stackrel1\boxtimes}}\,
  \,\simeq\, \Funle_{\Cala,\cala}(\cala,\cala)  \,\simeq\, \calz(\cala) \,.
  \label{eq:Drinf-framed-cent-A}
  \ee

\begin{Remark}
Using the equivalence \eqref{Phire...Psire} instead of \eqref{Phile.Psile}, we could as well 
have expressed framed centers through categories of right exact instead of left exact 
module functors. The resulting expressions would, however, be somewhat more complicated 
than the one obtained in
\eqref{eq:equ-lex-framed-cent}, owing to the additional occurrence of double duals.
For instance, $\cenf{ \CC\calm \,{\stackrel1\boxtimes}\, \caln \,{\stackrel1\boxtimes}}\,
\,{\simeq}\, \Funre_{\Cala,\calb}(\calm,{}^{2\!}_{}\caln^{\,2}_{}) \,.$
\end{Remark}

\begin{Remark} \label{Remark:Kappa-twist}
The equivalences among framed centers obtained in formula \eqref{eq:[k]l=k+l} can be
combined with other results so as to yield various further distinguished equivalences.
Specifically,
together with Corollary \ref{cor:lex-framed} we arrive at distinguished equivalences
  \be
  \CC{\calm} \,{\Boti{\kappa}}\, \caln \,{\Boti{\kappa'}}
  \,\simeq\, \Funle_{\cala,\calb}({}_{}^{1-\kappa\!\!}\calm{}^{1-\kappa'},\caln)
  \,\simeq\, \Funle_{\cala,\calb}(\calm,{}_{}^{\kappa-1\!\!}\caln{}^{\kappa'-1}) \,.
  \label{eq:MbotikappaN=Lex}
  \ee
\end{Remark}


\section{Assigning functors to \defectsurface s} \label{sec:preblock}
 
Our construction of the framed modular functor on the level of 1-morphisms involves several
different collections of left exact functors and proceeds in three steps: We first introduce, 
in Section \ref{sec:pre-blocks}, auxiliary functors $\Zpre(\Sigma)$, to which we refer as
\emph{pre-block} functors. In Section \ref{sec:blocks}, in a second step we impose constraint
to construct \emph{fine block} functors $\Zfine(\Sigma)$ from the pre-block functors. As
indicated by the terminology this procedure only makes sense if the surface $\Sigma$ is fine in
the sense of Definition \ref{def:fine}. To get the actual \emph{block} functors $\ZZ(\Sigma)$
-- which then furnish a modular functor in the sense of Definition \ref{def:modularfunctor}
-- for arbitrary, not necessarily fine, surfaces we introduce a notion of \emph{refinement}
of a surface (Section \ref{sec:refin-defect-surf}) and obtain $\ZZ(\Sigma)$ as a limit over
all refinements of $\Sigma$ (Section \ref{sec:nonfine}).

The functors $\Zpre(\Sigma)$, $\Zfine(\Sigma)$ and $\ZZ(\Sigma)$ are, a priori, functors 
from $\ZZ(\partiali\Sigma)$ to $\ZZ(\partialo\Sigma)$, with
$\partiali\Sigma\,{\sqcup}\,{-}\partialo\Sigma \,{=}\, \partialg\Sigma$. 
However, as explained in Appendix \ref{app:1} (see Equation \eqref{eq:MN->MccNvect}),
for finite categories we have equivalences 
  \be
  \Funle(\calm,\caln) \,\xrightarrow{~~\simeq~~}\, \Funle(\calm\boti\CC\caln,\vect) \,.
  \label{eq:LexMN-LexMNbvect}
  \ee
Accordingly, we will focus our attention to the case that $\partialo\Sigma \,{=}\, \emptyset$,
in which we deal with functors from $\ZZ(\partialg\Sigma)$ to $\vect$. The general case
is then obtained directly by invoking the equivalences \eqref{eq:LexMN-LexMNbvect}.


\subsection{Pre-block functors}\label{sec:pre-blocks}

To obtain the pre-block functor for a \defectsurface\ $\Sigma$ with
$\partialo\Sigma \,{=}\, \emptyset$ we start from a Hom functor 
whose covariant arguments come from the gluing categories for the gluing
boundaries of $\Sigma$. The Hom functor pairs every such covariant variable with a 
contravariant one; we take a coend over each of the latter variables.

In more detail, the pre-block functors are constructed as follows. The boundary
$\partialg\Sigma$ is a \defectonemanifold. According to \eqref{eq:TL=Z} the gluing
category $\ZZ(\LL)$ assigned to a connected \defectonemanifold\ $\LL$ is a
$\kappa$-framed center of the form $\cenf{ \calm_1^{\epsilon_1} \,{\Boti{\kappa^{}_1}}\,
\calm_2^{\epsilon_2} \,{\Boti{\kappa^{}_2}}\, \calm_3^{\epsilon_3}\,\Boti{...} \cdots}$.
We denote by $\UU(\LL)$ the category obtained from the gluing category for a
\defectonemanifold\ $\LL$ by forgetting the balancings of the framed center, and by
$U_\LL\colon \ZZ(\LL) \,{\to}\, \UU(\LL)$ the corresponding forgetful functor, i.e.\
$\UU(\LL) \,{=}\, \calm_1^{\epsilon_1} \boti \calm_2^{\epsilon_2} \boti \cdots$.
Specifically, we consider the forgetful functor
$U\colon\ZZ(\partialg\Sigma)\,{\xrightarrow{~~}}\, \UU(\partialg\Sigma)$; $U$ is an exact
functor. Now note that, by construction, in the image of $U$ each of the categories $\calm_i$
assigned to an edge appears `twice' -- that is, once as the category itself and once as its
opposite -- namely once at either end of the connecting defect line. Thus we have
$\UU(\partialg\Sigma) \,{=}\, \bigboxtimes_i\calm_i^{}\,{\boxtimes}\,\bigboxtimes_i\CC{\calm_i}$.
We can then give

\begin{Definition}
The \emph{pre-block functor} assigned to a \defectsurface\ $\Sigma$ is the left exact functor
  \be
  \Zpre(\Sigma):\quad \ZZ(\partialg\Sigma) \to \vect 
  \label{eq:T(d)->vect}
  \ee
from the gluing category associated with the boundary of $\Sigma$ to the category of
vector spaces that is constructed in the following manner:
For each factor $\calm_i^{}$ in $\ZZ(\partialg\Sigma)$ we insert an object 
$m_i\,{\boxtimes}\, \cc{m_i} \,{\in}\, \calm_i^{} \,{\boxtimes}\, \CC{\calm_i}$
as a contravariant variable in a Hom functor and take the coend 
  \be
  \Zpre(\Sigma)(-) \,:= \int^{m_1,m_2,\ldots , m_n}\! \Hom
  (m_1\boti \cc m_1\boti\cdots\boti m_n\boti \cc m_n \,, U(-))
  \label{eq:def:Zpre}
  \ee
over these variables in the finite category of left exact functors;
here the contravariant and covariant arguments of the Hom functor are
matched according to the combinatorial configuration of $\Sigma$.
\end{Definition}
    
It follows directly from the definition that the pre-block functors depend on the
\defectsurface\ $\Sigma$ only via the incidence combinatorics of 
the gluing and free boundaries and defect lines of $\Sigma$ , and that they satisfy 
  \be
  \Zpre(\Sigma\,{\sqcup}\, \Sigma') = \Zpre(\Sigma) \boxtimes \Zpre(\Sigma') \,. 
  \label{eq:T=TbotiT}
  \ee

\begin{Remark}
We require the (pre-)block functors we are working with to be left exact. In 
particular, as explained in Appendix \ref{app:1}, ends and coends are taken in
categories of left exact functors. The reason for doing so is as follows.
Once we have constructed pre-block functors for connected surfaces, the pre-block
functors for a general surface have to be defined on Deligne products of finite 
categories. To achieve this we need to use the universal property of the Deligne 
product. This property, in turn, holds for left exact functors, and 
likewise for right exact functors. In contrast, it does not hold for general
functors, nor does it hold when both left exact and right exact functors are
admitted simultaneously.
 \\
We could thus alternatively  work with right exact functors only and their
Deligne product. For a finite category $\calm$, the Hom functor $\Hom\colon \calmopp \boti
\calm \,{\to}\, \vect$ is left exact, while the `dual Hom' functor
  \be
  \begin{array}{rl}
  \widetilde\Hom:\quad \calmopp \boti \calm &\!\!\! \xrightarrow{~~}\, \vect \,,
  \nxl4
  m\boxtimes n &\!\! \xmapsto{~~}\, {\Hom(n,m)}^*  \eear
  \ee
(with the star denoting the vector space dual) is right exact.
Accordingly, right exact pre-block functors can be defined by the end
  \be
  \ZpreRE(\Sigma)(-) := \int_{m_1,m_2,...\,, m_n}\!\!\!
  \widetilde\Hom(U(-),m_1\boti \cc m_1\boti \cdots \boti m_n\boti \cc m_n) \,.
  \ee
Thus, in agreement with the Eilenberg--Watts equivalences \eqref{Phile.Psile} and
\eqref{Phire...Psire}, working with right exact functors instead of left exact ones boils down to
replacing coends and Hom functors by ends and dual Hom functors.
\end{Remark}

\begin{Remark}
Taking the coend makes it legitimate to refer to these variables as \emph{state-sum variables}.
Indeed, in case that the categories are finitely semisimple, this reduces to a sum over isomorphism
classes of simple objects (or of ``spins'', in the parlance of part of the physics literature).
Accordingly we refer to coends of the form appearing in \eqref{eq:def:Zpre} also as
\emph{state-sum coends} and regard our prescription as a state-sum 
construction. This fits with the fact that in state-sum models, vector spaces
associated to closed surfaces, also called \emph{block spaces} or spaces of
\emph{conformal blocks}, are constructed as subspaces of auxiliary vector spaces
that need to be introduced first; we call the latter pre-block spaces. 
\\
Such a two-step procedure is also the basis of the use of state-sum models in the construction
of quantum codes where, however, typically a space bigger than the pre-block space is used.
In that case, the essential idea is to obtain the block space for a surface $\Sigma$ as the
image of the projector that the three-dimensional Turaev--Viro topological field theory 
assigns to the cylinder $\Sigma\,{\times}\,[-1,1]$. In our construction we impose instead
the condition of flat holonomy for every contractible \twocell\ of a \defectsurface, see Section \ref{sec:blocks}.
\end{Remark}

\begin{Example} \label{ex:4.4-part1}
Consider a \defectsurface\ $\DDsAMN$ whose underlying surface is a disk, with two free boundary 
intervals labeled by $\cala$-modules ${}_\cala\calm$ and ${}_\cala\caln$ for some finite tensor
category $\cala$ and two gluing intervals, and with the $2$-framing given by the constant vector 
field pointing in the direction of the two free boundary intervals (and thus,
as needed, parallel to them):
    \def\locpa  {0.88}   
    \def\locpb  {0.71}   
    \def\locpA  {2.63}   
    \def\locpc  {2.2pt}  
    \def\locpd  {0.12}   
  \be 
  \raisebox{-6.6em}{\begin{tikzpicture}
  \scopeArrow{0.29}{>}
  \filldraw[\colorCircle,fill=\colorSigma,line width=\locpc,postaction={decorate}]
           (0,\locpA) arc (90:270:\locpA cm)
	   node[left=46pt,yshift=14pt,color=black]{$\scriptstyle-1$};
  \end{scope}
  \scopeArrow{0.45}{>}
  \filldraw[\colorCircle,fill=white,line width=\locpc,postaction={decorate}]
           (0,-\locpa) arc (270:90:\locpa cm)
	   node[left=11pt,yshift=-16pt,color=black]{$\scriptstyle 1$};
  \end{scope}
  \begin{scope}[shift={(-0.61*\locpA,0.85*\locpa)}] \drawOrientation; \end{scope}
  \scopeArrow{0.77}{\arrowDefect}
  \draw[line width=\locpc,color=\colorFreebdy,postaction={decorate}]
           (0,\locpa) -- node[right=-2pt,yshift=-4pt]{${}_\cala\caln$} (0,\locpA);
  \draw[line width=\locpc,color=\colorFreebdy,postaction={decorate}]
           (0,-\locpA) -- node[right=-2pt,yshift=-4pt]{${}_\cala\calm$} (0,-\locpa);
  \end{scope}
  \filldraw[color=\colorDefect] (0,\locpa) circle (\locpd cm) node[right=1pt]{$+$}
                                (0,\locpA) circle (\locpd cm) node[right=1pt]{$-$}
                                (0,-\locpa) circle (\locpd cm) node[right=1pt]{$-$}
                                (0,-\locpA) circle (\locpd cm) node[right=1pt]{$+$};
  \begin{scope}[shift={(0.1,-0.4*\locpa)}]
  \foreach \xx in {110,133,...,250}
           \draw[thick,color=\colorFrame,->] (\xx:1.6*\locpa cm) -- +(0,\locpb);
  \end{scope}
  \begin{scope}[shift={(0,-0.45*\locpa)}]
  \foreach \xx in {112,135,...,254}
           \draw[thick,color=\colorFrame,->] (\xx:2.4*\locpa cm) -- +(0,\locpb);
  \end{scope}
  \node at (-\locpA-1.5,0) {$\DDsAMN ~=$};
  \end{tikzpicture}}
  \label{pic:ex:4.4}
  \ee
Here, and in similar pictures below, we draw the gluing segments as half-circles,
to remind of the fact that they amount to the presence of a `boundary insertion', which
often is indicated by removing a half-disk from a two-manifold with boundary.
(For instance, the picture above is `half' of the picture \eqref{eq:picture-closed-annulus}.
Also recall that we are allowing for smooth manifolds with corners.)
Denoting the gluing interval on the right hand side of \eqref{pic:ex:4.4} by $\LL_1$ and
the left one by $\LL_2$, the relevant gluing categories are
  \be
  \ZZ(\LL_1) = \cenf{ \calmopp \,{\stackrel{1}\boxtimes}\, \caln}\, \qquad \text{and}\qquad
  \ZZ(\LL_2) = \cenf{ \calnopp \,{\stackrel{-1}\boxtimes}\, \calm} \,.
  \label{eq:TL1-TL2}
  \ee
The resulting pre-block functor is given by
  \be
  x \boti y \,\xmapsto{~~}
  \int^{m\in\calm,n\in\caln}\!\! \Hom_{\CC{\calm}\boxtimes\caln\boxtimes\CC{\caln}\boxtimes\calm}
  (\cc m\boti n\boti \cc n\boti m, \cc x_m\boti x_n\boti \cc y_n\boti y_m)
  \label{eq:preblock440}
  \ee
with $ x\,{=}\, \cc x_m\boti x_n \,{\in}\, \ZZ(\LL_1)$ and
$y \,{=}\, \cc y_n\boti y_m\,{\in}\, \ZZ(\LL_2)$.
Using the convolution property \eqref{eq:convo} of the Hom functor, this can be simplified to
  \be
  \Zpre(\DDsAMN)(x \boti y) = \Hom_\calm(x_m,y_m) \otimes \Hom_\caln(y_n,x_n) \,.
  \label{eq:preblock441}
  \ee
Also note that according to Proposition \ref{Cor2prop:EW:coendzFz}(ii) there is a
distinguished equivalence $\ZZ(\LL_{1}) \,{\simeq}$
                                                  $\Funle_{\cala}(\calm,\caln)$
mapping $x \,{\in}\, \ZZ(\LL_1)$ to the left exact functor $\Phile(x)$, while by Lemma
\ref{lemma:Op-circle} we have an equivalence $\ZZ(\LL_{2}) \,{\simeq}\, \ZZ(\LL_{1})\Opp$.
Under these equivalences the pre-block space becomes
  \be
  \Zpre(\DDsAMN)(x \boti y) \cong \Nat(\Phile(\cc y),\Phile(x)) \,.
  \label{eq:preblock442}
  \ee
\end{Example}

As will become clear in Corollary \ref{cor:simpledisk} below, the \defectsurface\ 
$\DDsAMN$ \eqref{pic:ex:4.4} considered in Example \ref{ex:4.4-part1} is indeed
the most basic surface for us. In view of the particular form of the framing vector field on 
this surface, we will refer to $\DDsAMN$ as the \emph{\basicdisk}.

The results \eqref{eq:preblock441} and \eqref{eq:preblock442} express the pre-block
functor as a Deligne product of Hom functors and as natural transformations,
respectively. This is no coincidence, but is a generic feature of the construction, which is
a first hint at the power of the modular functor to produce algebraically interesting quantities.
To pinpoint this issue, let us also have a look at slightly more complicated surfaces.

It is worth pointing out that the framing enters the definition of the pre-block only via 
the gluing categories. Thus, as a direct consequence of 
the canonical equivalence of gluing categories in Proposition \ref{prop:lineflip}, we obtain

\begin{Lemma}\label{lem:lineflipPreblocks}
There is a canonical isomorphism between the pre-block functors for the configurations
    \def\locpa  {0.52}   
    \def\locpb  {1.4}    
    \def\locpc  {2.2pt}  
    \def\locpd  {0.12}   
    \def\locpe  {1.4}    
    \def\locpE  {1.5}    
  \be
  \raisebox{-7.1em}{\begin{tikzpicture}
  [decoration={random steps,segment length=2mm}]
  \fill[\colorSigma,decorate] (0,-\locpa-\locpb) rectangle (-\locpe,3*\locpa+2*\locpb);
  \fill[white] (0,0) circle (\locpa) (0,2*\locpa+\locpb) circle (\locpa);
  \scopeArrow{0.48}{>}
  \draw[line width=\locpc,\colorCircle,postaction={decorate}]
           (0,-\locpa) arc (270:90:\locpa)
	   node[left=11pt,yshift=-5pt,color=black] {$\kappa$};
  \draw[line width=\locpc,\colorCircle,postaction={decorate}]
           (0,\locpa+\locpb) arc (270:90:\locpa)
	   node[left=11pt,yshift=-5pt,color=black] {$\kappa'$};
  \end{scope}
  \scopeArrow{0.77}{\arrowDefect}
  \draw[line width=\locpc,color=\colorFreebdy,postaction={decorate}]
           (0,3*\locpa+\locpb) -- node[right=-2pt,yshift=-6pt] {$\calk$} +(0,\locpb);
  \draw[line width=\locpc,color=\colorFreebdy,postaction={decorate}]
           (0,-\locpa-\locpb) -- node[right=-2pt,yshift=-6pt] {$\calm$} +(0,\locpb);
  \end{scope}
  \scopeArrow{0.52}{\arrowDefectB}
  \draw[line width=\locpc,color=\colorFreebdy,postaction={decorate}]
           (0,\locpa) -- node[right=-2pt,yshift=8pt] {$\caln$} +(0,\locpb);
  \end{scope}
  \filldraw[color=\colorDefect] (0,\locpa) circle (\locpd) (0,-\locpa) circle (\locpd)
           (0,\locpa+\locpb) circle (\locpd) (0,3*\locpa+\locpb) circle (\locpd);
  \node at (2.22,\locpa+0.5*\locpb) {and};
  \begin{scope}[shift={(5.4,0)}]
  [decoration={random steps,segment length=2mm}]
  \fill[\colorSigma,decorate] (0,-\locpa-\locpb) rectangle (-\locpE,3*\locpa+2*\locpb);
  \fill[white] (0,0) circle (\locpa) (0,2*\locpa+\locpb) circle (\locpa);
  \scopeArrow{0.48}{>}
  \draw[line width=\locpc,\colorCircle,postaction={decorate}]
           (0,-\locpa) arc (270:90:\locpa)
	   node[left=10pt,yshift=-5pt,color=black] {$\kappa{+}1$};
  \draw[line width=\locpc,\colorCircle,postaction={decorate}]
           (0,\locpa+\locpb) arc (270:90:\locpa)
	   node[left=10pt,yshift=-5pt,color=black] {$\kappa'{-}1$};
  \end{scope}
  \scopeArrow{0.77}{\arrowDefect}
  \draw[line width=\locpc,color=\colorFreebdy,postaction={decorate}]
           (0,3*\locpa+\locpb) -- node[right=-2pt,yshift=-6pt]{$\calk$} +(0,\locpb);
  \draw[line width=\locpc,color=\colorFreebdy,postaction={decorate}] (0,\locpa) --
           node[right=-3pt,yshift=-8pt] {${}_{}^{1}\cc{\caln}$} +(0,\locpb);
  \draw[line width=\locpc,color=\colorFreebdy,postaction={decorate}]
           (0,-\locpa-\locpb) -- node[right=-2pt,yshift=-6pt]{$\calm$} +(0,\locpb);
  \end{scope}
  \filldraw[color=\colorDefect] (0,\locpa) circle (\locpd) (0,-\locpa) circle (\locpd)
           (0,\locpa+\locpb) circle (\locpd) (0,3*\locpa+\locpb) circle (\locpd);
  \end{scope}
  \end{tikzpicture}}
  \label{eq:lineflipPreblocks}
  \ee
involving a local change of the framing, as well as, more generally, for $z \,{\in}\, 2\Z{+}1$
an isomorphism with $\kappa{+}z$, $\kappa'{-}z$ and ${}_{}^{2-z}\cc{\caln}$
in place of $\kappa{+}1$, $\kappa'{-}1$ and ${}_{}^{1}\cc{\caln}$ on the right hand side.
Similarly, there is a canonical isomorphism between the pre-block functors for the configurations
    \def\locpa  {0.52}   
    \def\locpb  {1.4}    
    \def\locpc  {2.2pt}  
    \def\locpd  {0.12}   
    \def\locpe  {1.4}    
    \def\locpE  {1.5}    
  \be
  \raisebox{-7.1em}{\begin{tikzpicture}
  [decoration={random steps,segment length=2mm}]
  \fill[\colorSigma,decorate] (0,-\locpa-\locpb) rectangle (\locpe,3*\locpa+2*\locpb);
  \fill[white] (0,0) circle (\locpa) (0,2*\locpa+\locpb) circle (\locpa);
  \scopeArrow{0.48}{<}
  \draw[line width=\locpc,\colorCircle,postaction={decorate}]
           (0,-\locpa) arc (-90:90:\locpa)
	   node[right=11pt,yshift=-5pt,color=black] {$\kappa$};
  \draw[line width=\locpc,\colorCircle,postaction={decorate}]
           (0,\locpa+\locpb) arc (-90:90:\locpa)
	   node[right=11pt,yshift=-5pt,color=black] {$\kappa'$};
  \end{scope}
  \scopeArrow{0.77}{\arrowDefect}
  \draw[line width=\locpc,color=\colorFreebdy,postaction={decorate}]
           (0,3*\locpa+\locpb) -- node[left=-2pt,yshift=-6pt] {$\calk$} +(0,\locpb);
  \draw[line width=\locpc,color=\colorFreebdy,postaction={decorate}]
           (0,-\locpa-\locpb) -- node[left=-2pt,yshift=-6pt] {$\calm$} +(0,\locpb);
  \end{scope}
  \scopeArrow{0.52}{\arrowDefectB}
  \draw[line width=\locpc,color=\colorFreebdy,postaction={decorate}]
           (0,\locpa) -- node[left=-2pt,yshift=8pt] {$\caln$} +(0,\locpb);
  \end{scope}
  \filldraw[color=\colorDefect] (0,\locpa) circle (\locpd) (0,-\locpa) circle (\locpd)
           (0,\locpa+\locpb) circle (\locpd) (0,3*\locpa+\locpb) circle (\locpd);
  \node at (3.02,\locpa+0.5*\locpb) {and};
  \begin{scope}[shift={(5.4,0)}]
  [decoration={random steps,segment length=2mm}]
  \fill[\colorSigma,decorate] (0,-\locpa-\locpb) rectangle (\locpE,3*\locpa+2*\locpb);
  \fill[white] (0,0) circle (\locpa) (0,2*\locpa+\locpb) circle (\locpa);
  \scopeArrow{0.48}{<}
  \draw[line width=\locpc,\colorCircle,postaction={decorate}]
           (0,-\locpa) arc (-90:90:\locpa)
	   node[right=11pt,yshift=-5pt,color=black] {$\kappa{+}1$};
  \draw[line width=\locpc,\colorCircle,postaction={decorate}]
           (0,\locpa+\locpb) arc (-90:90:\locpa)
	   node[right=11pt,yshift=-5pt,color=black] {$\kappa'{-}1$};
  \end{scope}
  \scopeArrow{0.77}{\arrowDefect}
  \draw[line width=\locpc,color=\colorFreebdy,postaction={decorate}]
           (0,3*\locpa+\locpb) -- node[left=-2pt,yshift=-6pt]{$\calk$} +(0,\locpb);
  \draw[line width=\locpc,color=\colorFreebdy,postaction={decorate}] (0,\locpa) --
           node[left=-3.4pt,yshift=-9pt] {$\CC{\caln}_{}^{1}$} +(0,\locpb);
  \draw[line width=\locpc,color=\colorFreebdy,postaction={decorate}]
           (0,-\locpa-\locpb) -- node[left=-2pt,yshift=-6pt]{$\calm$} +(0,\locpb);
  \end{scope}
  \filldraw[color=\colorDefect] (0,\locpa) circle (\locpd) (0,-\locpa) circle (\locpd)
           (0,\locpa+\locpb) circle (\locpd) (0,3*\locpa+\locpb) circle (\locpd);
  \end{scope}
  \end{tikzpicture}}
  \label{eq:lineflipPreblocksR}
  \ee
\end{Lemma}

Our construction of pre-blocks is compatible with composition 
of left exact functors via the Eilenberg--Watts correspondence
  \be
  \Psile: \quad \Funle(\calm,\caln) \to \CC\calm\boxtimes\caln \,.
  \ee
Explicitly we have the following `fusion of boundary insertions':

\begin{Proposition}\label{prop:preblock-compose}
There is a canonical isomorphism between the pre-block functors for the configurations
    \def\locpa  {0.52}   
    \def\locpA  {0.84}   
    \def\locpb  {1.4}    
    \def\locpB  {2.3}    
    \def\locpc  {2.2pt}  
    \def\locpd  {0.12}   
    \def\locpe  {1.4}    
  \be
  \raisebox{-7.2em}{\begin{tikzpicture}
  [decoration={random steps,segment length=2mm}]
  \fill[\colorSigma,decorate] (0,-\locpa-\locpb) rectangle (-\locpe,3*\locpa+2*\locpb);
  \fill[white] (0,0) circle (\locpa);
  \fill[white] (0,2*\locpa+\locpb) circle (\locpa);
  \scopeArrow{0.52}{>}
  \draw[line width=\locpc,\colorCircle,postaction={decorate}]
           (0,-\locpa) arc (270:90:\locpa cm)
           node[left=9pt,yshift=-3pt,color=black]{$\Psile(F)$};
  \draw[line width=\locpc,\colorCircle,postaction={decorate}]
           (0,\locpa+\locpb) arc (270:90:\locpa cm)
           node[left=9pt,yshift=-3pt,color=black]{$\Psile(G)$};
  \end{scope}
  \scopeArrow{0.82}{\arrowDefect}
  \draw[line width=\locpc,color=\colorFreebdy,postaction={decorate}]
           (0,3*\locpa+\locpb) -- node[right=-2pt,yshift=-4pt]{$\calk$} +(0,\locpb);
  \draw[line width=\locpc,color=\colorFreebdy,postaction={decorate}]
           (0,\locpa) -- node[right=-2pt,yshift=-3pt]{$\caln$} +(0,\locpb);
  \draw[line width=\locpc,color=\colorFreebdy,postaction={decorate}]
           (0,-\locpa-\locpb) -- node[right=-2pt,yshift=-6pt]{$\calm$} +(0,\locpb);
  \end{scope}
  \filldraw[color=\colorDefect] (0,\locpa) circle (\locpd) (0,-\locpa) circle (\locpd)
           (0,\locpa+\locpb) circle (\locpd) (0,3*\locpa+\locpb) circle (\locpd);
  \begin{scope}[shift={(1.8,\locpa+0.5*\locpb)}] \node at (0,0) {and}; \end{scope}
  \begin{scope}[shift={(4.8,0)}]
  [decoration={random steps,segment length=2mm}]
  \fill[\colorSigma,decorate] (0,-\locpa-\locpb) rectangle (-\locpe,3*\locpa+2*\locpb);
  \end{scope}
  \begin{scope}[shift={(4.8,\locpA+\locpB-\locpa-\locpb)}]
  \fill[white] (0,0) circle (\locpA);
  \scopeArrow{0.48}{>}
  \draw[line width=\locpc,\colorCircle,postaction={decorate}]
           (0,-\locpA) arc (270:90:\locpA cm)
           node[left=14pt,yshift=-3pt,color=black]{$\Psile(G{\circ}F)$};
  \end{scope}
  \scopeArrow{0.82}{\arrowDefect}
  \draw[line width=\locpc,color=\colorFreebdy,postaction={decorate}]
           (0,\locpA) -- node[right=-2pt,yshift=-1pt]{$\calk$} +(0,\locpB);
  \end{scope}
  \scopeArrow{0.77}{\arrowDefect}
  \draw[line width=\locpc,color=\colorFreebdy,postaction={decorate}]
           (0,-\locpA-\locpB) -- node[right=-2pt,yshift=-9pt]{$\calm$} +(0,\locpB);
  \end{scope}
  \filldraw[color=\colorDefect] (0,\locpA) circle (\locpd) (0,-\locpA) circle (\locpd);
  \end{scope}
  \end{tikzpicture}}
  \label{eq:PsiG.PsiF=PsiGF}
  \ee
involving a local replacement of segments around a disk,
with the framing near the segments being along the positive $y$-axis.
\end{Proposition}
 
\begin{proof}
Consider left exact bimodule functors $F\colon \calm \,{\to}\, \caln$ and
$G\colon \caln \,{\to}\, \calk$ for $\calb$-$\cala$-bimodules $\calm,\caln, \calk$.
The composite of the isomorphisms (see \Cite{Cor.\,3.7}{fScS2})
  \be
  \bearll
  \multicolumn{2}{l}{\dsty \int^{n\in\caln}\!
  \Hom(\cc m\boxtimes n\boxtimes \cc n\boxtimes k, \Psile(F)\boxtimes \Psile(G)) }
  \Nxl1 \qquad\qquad &\dsty
  \cong \int^{n \in \caln}\! \Hom(\cc{m} \boxtimes n, \Psile(F)) \otimes \Hom(\cc{n} \boxtimes k, \Psile(G)) 
  \Nxl1 &\dsty
  \cong \int^{n \in \caln}\Hom(n, F(m)) \otimes \Hom(k,G(n)) 
  \stackrel{\eqref{eq:EW-xxx}}{\cong}
  \Hom(k, G{\circ}F(m)) 
  \label{eq:holonomy-fusion-line}

  \eear
  \ee
provides the desired isomorphism between the pre-block functors for the left and 
right hand sides of \eqref{eq:PsiG.PsiF=PsiGF}.
\end{proof}

Proposition \ref{prop:preblock-compose} illustrates how the modular functor realizes algebraic
structures: the fusion of boundary insertions provides the composition of functors and thus,
via the Eilenberg--Watts calculus, a composition on the Deligne product $\cc\calm\boti\caln$.
In particular, for any module category $\calm$ it yields a distinguished object in
$\cc\calm\boti\calm$, namely the one that acts like a unit for the type of local
replacement considered in \eqref{eq:PsiG.PsiF=PsiGF}.

\begin{Example} \label{exa:generaldisk}
The generalization of Example \ref{ex:4.4-part1} to a disk with any number $N$ of gluing and
free boundary segments, with the former oriented
as induced by the orientation of the disk and the latter of arbitrary orientation, and
with any indices is immediate.
First, by invoking Proposition \ref{prop:lineflip} we can restrict our attention to any
specific choice of orientations of the free boundary segments, say one of them oriented
counter-clockwise and all others clockwise, as indicated for $N \,{=}\, 4$ in the picture
    \def\locpa  {0.8}    
    \def\locpA  {2.5}    
    \def\locpb  {0.9}    
    \def\locpc  {2.2pt}  
    \def\locpd  {0.12}   
  \be
  \raisebox{-5.9em} 
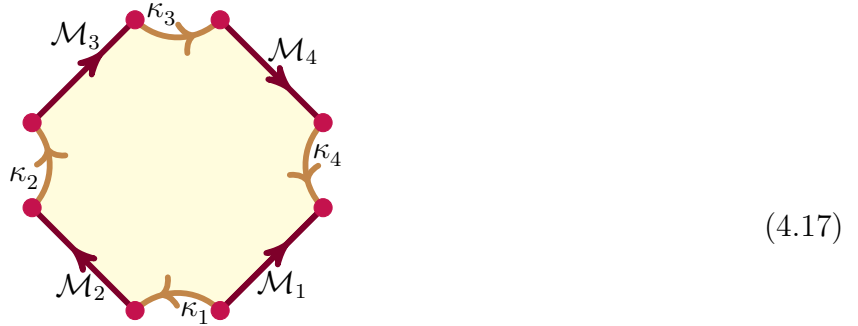
{\begin{tikzpicture}[scale=1.1]
  \scopeArrow{0.15}{\arrowDefect} \scopeArrow{0.64}{\arrowDefect} \scopeArrow{0.90}{\arrowDefect}
  \filldraw[\colorFreebdy,fill=\colorSigma,line width=\locpc,postaction={decorate}]
           (0,\locpA) node[right=29pt,yshift=-28pt,color=black] {$\calm_4$} -- (\locpA,0) 
           -- (0,-\locpA) node[left=22pt,yshift=25pt,color=black] {$\calm_2$}
           -- (-\locpA,0) node[right=18pt,yshift=49pt,color=black] {$\calm_3$} -- cycle;
  \end{scope} \end{scope} \end{scope} 
  \scopeArrow{0.51}{\arrowDefectB}
  \draw[\colorFreebdy,line width=\locpc,postaction={decorate}]
           (\locpA,0) -- (0,-\locpA) node[right=25pt,yshift=26pt,color=black] {$\calm_1$};
  \end{scope}
  \fill[white] (-\locpb,\locpA+\locpb-\locpa) arc (180:360:\locpb)
           (\locpb,\locpa-\locpb-\locpA) arc (0:180:\locpb)
           (\locpA,0) circle (\locpa) (-\locpA,0) circle (\locpa);
  \scopeArrow{0.70}{>}          
  \filldraw[\colorCircle,fill=white,line width=\locpc,postaction={decorate}]
           (0,-\locpA) +(45:\locpa) arc (45:135:\locpa)
           node[right=12pt,yshift=-1pt,color=black]{$\kappa_1^{}$};
  \filldraw[\colorCircle,fill=white,line width=\locpc,postaction={decorate}]
           (0,\locpA) +(225:\locpa) arc (225:315:\locpa)
           node[left=11.4pt,yshift=2pt,color=black]{$\kappa_3^{}$};
  \filldraw[\colorCircle,fill=white,line width=\locpc,postaction={decorate}]
           (-\locpA,0) +(315:\locpa) arc (-45:45:\locpa)
           node[below=12pt,xshift=-3pt,color=black]{$\kappa_2^{}$};
  \filldraw[\colorCircle,fill=white,line width=\locpc,postaction={decorate}]
           (\locpA,0) +(135:\locpa) arc (135:225:\locpa)
           node[above=11pt,xshift=1.7pt,color=black]{$\kappa_4^{}$};
  \end{scope}
  \filldraw[color=\colorDefect] 
           (\locpA,0) +(135:\locpa) circle (\locpd) (\locpA,0) +(225:\locpa) circle (\locpd)
           (-\locpA,0) +(45:\locpa) circle (\locpd) (-\locpA,0) +(315:\locpa) circle (\locpd)
           (0,\locpA) +(225:\locpa) circle (\locpd) (0,\locpA) +(315:\locpa) circle (\locpd)
           (0,-\locpA) +(45:\locpa) circle (\locpd) (0,-\locpA) +(135:\locpa) circle (\locpd);
  \end{tikzpicture}}
  \ee
(with a left $\cala$-module $\calm_1$ and right $\cala$-modules
$\calm_2$, $\calm_3$, $\calm_4$, for some finite tensor category $\cala$),
Next we can use Proposition \ref{prop:preblock-compose}
to reduce the number of gluing and free boundary segments by one, at the same time 
composing the functors that the Eilenberg--Watts equivalence assigns to the objects at the 
gluing segments. Doing so iteratively we end up with pre-blocks given by
  \be
  \Zpre(x_1 \boti x_2 \boti\cdots\boti x_N) \cong \Nat\big(\Phile(\cc{x_1}),
  \Phile(x_N)\,{\circ}\,\Phile(x_{N-1})\,{\circ}\cdots{\circ}\,\Phile(x_2) \big) \,,
  \label{eq:exa:generaldisk-Zpre}
  \ee
with $x_i \,{\in}\,\ZZ(\LL_i)$, such that $U(x_{i})\,{\in}\,\CC{\calm_{i+1}}\boti\calm_i$,
thereby generalizing (and using analogous notation as in) formula \eqref{eq:preblock442}.
\end{Example}
      
As this example indicates, by combining the Propositions and \ref{prop:preblock-compose} 
and \ref{prop:lineflip} one obtains

\begin{cor}\label{cor:simpledisk0}
The pre-block functor for any disk without defect lines and with an arbitrary number of free 
boundaries and
gluing segments can be reduced to the pre-block functor for the \basicdisk\ \eqref{pic:ex:4.4}.
\end{cor}

\begin{Example} \label{exa:gzlhs2-2Pachner}
For $\cala$ and $\calb$ finite tensor categories, $\calm_1$ and $\calm_4$
right $\calb$-modules, $\calm_2$ and $\calm_3$ right $\cala$-modules, and
$\calk$ an $\cala$-$\calb$-bimodule, consider the \defectsurface\
    \def\locpa  {0.95}   
    \def\locpA  {3.1}    
    \def\locpb  {1.1}    
    \def\locpc  {2.2pt}  
    \def\locpd  {0.12}   
  \be
  \raisebox{-7.3em} {\begin{tikzpicture}
  \scopeArrow{0.16}{\arrowDefect} \scopeArrow{0.41}{\arrowDefect}
  \scopeArrow{0.66}{\arrowDefect} \scopeArrow{0.91}{\arrowDefect}
  \filldraw[\colorDefect,fill=\colorSigma,line width=\locpc,postaction={decorate}]
           (0,\locpA) node[right=33pt,yshift=-32pt,color=black] {$\calm_4$}
           -- (\locpA,0) node[left=16pt,yshift=-44pt,color=black] {$\calm_1$}
           -- (0,-\locpA) node[left=31pt,yshift=33pt,color=black] {$\calm_2$}
           -- (-\locpA,0) node[right=15pt,yshift=46pt,color=black] {$\calm_3$} -- cycle;
  \end{scope} \end{scope} \end{scope} \end{scope}
  \scopeArrow{0.66}{\arrowDefect}
  \draw[\colorDefect,line width=\locpc,postaction={decorate}]
           (0,\locpa-\locpA) -- node[right=-3pt,yshift=-8pt,color=black] {$\calk$}
           node[left=20pt,yshift=24pt,color=black] {$\cala$}
           node[right=19pt,yshift=24pt,color=black] {$\calb$} (0,\locpA-\locpa);
  \end{scope}
  \fill[white] (-\locpb,\locpA+\locpb-\locpa) arc (180:360:\locpb)
           (\locpb,\locpa-\locpb-\locpA) arc (0:180:\locpb)
           (\locpA,0) circle (\locpa) (-\locpA,0) circle (\locpa);
  \scopeArrow{0.31}{>} \scopeArrow{0.76}{>}          
  \filldraw[\colorCircle,fill=white,line width=\locpc,postaction={decorate}]
           (0,-\locpA) +(45:\locpa) arc (45:135:\locpa)
           node[right=-8pt,yshift=13pt,color=black]{$\kappa_1'$}
           node[right=13pt,yshift=-1pt,color=black]{\begin{turn}{-33}$\kappa_1^{}$\end{turn}};
  \end{scope} \end{scope}
  \scopeArrow{0.39}{>} \scopeArrow{0.78}{>}          
  \filldraw[\colorCircle,fill=white,line width=\locpc,postaction={decorate}]
           (0,\locpA) +(225:\locpa) arc (225:315:\locpa)
	   node[left=26pt,yshift=-10pt,color=black]{$\kappa_3^{}$}
           node[left=-2pt,yshift=4pt,color=black]{\begin{turn}{33}$\kappa_3'$\end{turn}};
  \end{scope} \end{scope}
  \scopeArrow{0.76}{>}          
  \filldraw[\colorCircle,fill=white,line width=\locpc,postaction={decorate}]
           (-\locpA,0) +(315:\locpa) arc (-45:45:\locpa)
	   node[below=16pt,xshift=-1pt,color=black]{$\kappa_2^{}$};
  \filldraw[\colorCircle,fill=white,line width=\locpc,postaction={decorate}]
           (\locpA,0) +(135:\locpa) arc (135:225:\locpa)
	   node[above=16pt,xshift=1pt,color=black]{$\kappa_4^{}$};
  \end{scope}
  \filldraw[color=\colorDefect] 
           (\locpA,0) +(135:\locpa) circle (\locpd) (\locpA,0) +(225:\locpa) circle (\locpd)
           (-\locpA,0) +(45:\locpa) circle (\locpd) (-\locpA,0) +(315:\locpa) circle (\locpd)
           (0,\locpA) +(225:\locpa) circle (\locpd) (0,\locpA) +(315:\locpa) circle (\locpd)
           (0,-\locpA) +(45:\locpa) circle (\locpd) (0,-\locpA) +(135:\locpa) circle (\locpd)
           (0,\locpa-\locpA) circle (\locpd) (0,\locpA-\locpa) circle (\locpd);
  \node at (-\locpA-0.7,0) {$\Sigma ~=$};
  \end{tikzpicture}}
  \label{eq:gzlhs2-2Pachner}
  \ee
Applied to objects 
  \be
  \bearl
  z_1 \,{=}\, x_2 \boti u \boti \cc{y_1}
  \,\in \calm_2 \!\Boti{-\kappa_1'} \calk \Boti{-\kappa_1^{}\,} \CC{\calm_1} \,,
  \Nxl4
  z_2 \,{=}\, x_3 \boti \cc{y_2} \,\in \calm_3 \!\Boti{-\kappa_2^{}\,} \CC{\calm_2} \,,
  \Nxl4
  z_3 \,{=}\, x_4\boti\cc v\boti \cc{y_3}
  \,\in \calm_4 \!\Boti{-\kappa_3'\,} \CC\calk \Boti{-\kappa_3^{}\,} \CC{\calm_3} \,, \quad
  \Nxl4
  z_4 \,{=}\, x_1 \boti \cc{y_4} \,\in \calm_1 \!\Boti{-\kappa_4^{}\,} \CC{\calm_4}
  \eear
  \ee
of the gluing categories for the four gluing intervals, the pre-block functor gives
  \be
  \bearl\dsty
  \Zpre(z_1\boti z_2\boti z_3\boti z_4)\,
  = \int^{m_1\in\calm_1,m_2\in\calm_2,m_3\in\calm_3,m_4\in\calm_4,k\in\calk}
  \Nxl2
  \hspace*{9em} \Hom(m_1 \boti \cc{m_1} \boti m_2 \boti \cc{m_2} \boti m_3 \boti \cc{m_3}
  \boti m_4 \boti \cc{m_4} \boti k \boti \cc k \,,
  \Nxl3
  \hspace*{14.25em} x_1 \boti \cc{y_1} \boti x_2 \boti \cc{y_2} \boti
  x_3 \boti \cc{y_3} \boti x_4 \boti \cc{y_4} \boti u \boti \cc v) \,.
  \eear
  \ee
By a multiple application of the variant \eqref{eq:convo} of the Yoneda lemma, this reduces to
  \be
  \bearll
  \Zpre(z_1\boti z_2\boti z_3\boti z_4) = \!&
  \Hom_{\calm_1}(y_1,x_1) \oti \Hom_{\calm_2}(y_2,x_2)
  \Nxl4 &
  \oti \Hom_{\calm_3}(y_3,x_3) \oti \Hom_{\calm_4}(y_4,x_4) \oti \Hom_{\calk}(v,u) \,.
  \eear
  \label{eq:Zpre-gzlhs2-2}
  \ee
This result may again be written in terms of natural transformations. There are now two
distinguished ways to do so, one corresponding to fusing the $\calk$-defect to the
$\calm_2$- and $\calm_3$-de\-fects (this process of fusion will be discussed in detail in Section
\ref{sec:fusion} below), and one corresponding to fusing it to the $\calm_1$- and $\calm_4$-defects.
Let us write out the resulting expression for the former case: one gets the
natural transformation
  \be
  \Zpre(z_1\boti z_2\boti z_3\boti z_4) = \Nat(\Phile(z_4),G(z_1,z_2,z_3))
  \ee
of functors in $\Funle(\calm_1,\calm_4)$,
where $\Phile$ is the Eilenberg--Watts equivalence \eqref{Phile.Psile} and $G$ is the composition
  \be
  G(z_1,z_2,z_3) := \Phile(z_3) \big( [\Phile(z_2) \boti \Id_\calk] (\Phile(z_1)) \big) \,. 
  \ee
\end{Example}

In fact, the pre-block spaces for any arbitrary \defectsurface\ can be
expressed as tensor products of morphism spaces, analogously as in
\eqref{eq:Zpre-gzlhs2-2}. The so obtained expressions are, however, not particularly
illuminating. Expressing them through spaces of natural transformations
can be more informative, e.g.\ it often allows for a direct characterization of what
subspaces of pre-blocks furnish the block spaces.


\subsection{Holonomy} \label{sec:holonomy}

The pre-block functors do not see the framing of a \defectsurface\ and do not take the 
topology of the \twocells\ of a \defectsurface\ into account. The proper block functors
$\ZZ(\Sigma)$ that we are going to introduce will, on the other hand, depend on the framing,
and in their definition \twocells\ with the topology of a disk will play a crucial role.
To proceed from the pre-blocks to the block functors, we first define holonomy operations
on pre-blocks. For doing so we will restrict our attention to the subclass of surfaces
that can be patched together from disks. 

Recall from Definition \ref{def:fine}(ii) that
a \defectsurface\ $\Sigma$ is called \emph{fine} iff every \twocell\
(in the sense of Definition \ref{def:twocell}) is contractible.
{}From now on we assume that the \defectsurface\ $\Sigma$ under consideration
is fine. Then we have one holonomy operation for each \twocell\ of $\Sigma$.
The formal definition of these holonomy operations is notationally somewhat intricate.
Instead of spelling out the details, we explain these operations through concrete examples.

\begin{Example} \label{ex:4.4}
Let us give full details for the \basicdisk\ $\DDsAMN$ described in Example \ref{ex:4.4-part1}. 
In this case we have two gluing intervals $\LL_1$ and $\LL_2$, and the pre-block functor
on objects $x \,{\in}\, \ZZ(\LL_1)$ and $y \,{\in}\, \ZZ(\LL_2)$ in the gluing categories 
is given by formula \eqref{eq:preblock440}. For $a \,{\in}\, \cala$, the holonomy 
$\hol_{a,x}$ of $a$ starting at $x \,{=}\, \cc x_m\boti x_n$ is defined to be
the natural isomorphism between functors on $\cala \boti \ZZ(\LL_1) \boti \ZZ(\LL_2)$
that in terms of its components at $y \,{\in}\, \ZZ(\LL_2)$ is the composite
  \begin{eqnarray}
  \hspace*{-0.5em}\begin{array}{lll}
  \multicolumn{2}{l} {\dsty \int^{m\in\calm,n\in\caln}\!
  \Hom(\cc m\boti n\boti \cc n\boti m, \cc x_m\boti a.x_n\boti \cc y_n\boti y_m) }
  \Nxl1
  \quad & \tocong \Hom(\int_{m\in\calm,n\in\caln}\cc m\boti n \boti \cc n\boti m,
  \cc x_m\boti a.x_n\boti \cc y_n\boti y_m)
  & \text{Eq.\,[\eqref{eq:fScS2Prop3.4}]} 
  \Nxl7
  \quad & \tocong \Hom(\int_{m\in\calm,n\in\caln}\cc m\boti a^\vee\!.n \boti \cc n\boti m,
  \cc x_m\boti x_n\boti \cc y_n\boti y_m)
  & \text{[right duality]}
  \Nxl7
  & \tocong \Hom(\int_{m\in\calm,n\in\caln}\cc m\boti n \boti \cc {a.n}\boti m,
  \cc x_m\boti x_n\boti \cc y_n\boti y_m)
  & \text{[Lemma\,\ref{lemma:balancing-psil(F)}]}
  \Nxl7
  & \tocong \Hom(\int_{m\in\calm,n\in\caln}\cc m\boti n \boti \cc {n}\boti m,
  \cc x_m\boti x_n\boti \cc {a^\vee.y_n}\boti y_m)
  & \text{[right duality]}
  \Nxl7
  & \tocong \Hom(\int_{m\in\calm,n\in\caln}\cc m\boti n \boti \cc {n}\boti m,
  \cc x_m\boti x_n\boti \cc y_n \boti a^{\vee \vee}.y_m)
  & \text{[balancing of $y$]}
  \Nxl7
  & \tocong \Hom(\int_{m\in\calm,n\in\caln}\cc m\boti n \boti \cc {n}
  \boti a^{\vee\vee\vee}\!.m, \cc x_m\boti x_n\boti \cc y_n \boti y_m)
  & \text{[right duality]}
  \Nxl7
  & \tocong \Hom(\int_{m\in\calm,n\in\caln}\cc {a^{\vee\vee}.m}\boti n \boti \cc {n}\boti m,
  \cc x_m\boti x_n\boti \cc y_n \boti y_m)
  & \text{[Lemma \ref{lemma:balancing-psil(F)}]}
  \Nxl7
  & \tocong \Hom(\int_{m\in\calm,n\in\caln}\cc {m}\boti n \boti \cc {n}\boti m,
  \cc {a^{\vee\vee\vee}\!.x_m}\boti x_n\boti \cc y_n \boti y_m)
  \quad & \text{[right duality]} \,.
  \\[-1.4em]~
  \eear
  \nonumber \\[0.4em]~
  \label{eq:compositeHol}
  \end{eqnarray}
Thus the holonomy is a distinguished isomorphism $\hol_{a,x}$ with components
  \be
  \begin{array}{l}\dsty
  (\hol_{a,x})_y^{}: \quad \int^{m\in\calm,n\in\caln}\!\!
  \Hom(\cc m\boti n\boti \cc n\boti m, \cc x_m\boti a.x_n\boti \cc y_n\boti y_m)
  \Nxl1\dsty
  \hspace*{7.5 em} \tocong
  \int^{m\in\calm,n\in\caln}\!\! \Hom(\cc {m}\boti n \boti \cc {n}\boti m,
  \cc {a^{\vee\vee\vee}\!.x_m}\boti x_n\boti \cc y_n \boti y_m) 
  \eear
  \label{eq:iso1}
  \ee
of left exact functors. 
\\
The rationale behind this prescription is simple: we proceed -- counterclockwise, by convention
-- along the boundary of the \twocell\ and use alternatingly a duality to jump between a covariant 
and a contravariant argument of the Hom and a balancing in one of the two arguments. In the
contravariant state-sum variable, the balancing is given by Lemma \ref{lemma:balancing-psil(F)}
with $F \eq \Id$; 
in the covariant argument coming from a gluing boundary, it is part of the structure given by
the gluing categories. In the latter, the index resulting from the $2$-framing enters.
 \\
On the other hand, the balancing of the variable $x \,{\in}\, \ZZ(\LL_1)$ gives an isomorphism
  \be
  \begin{array}{l}\dsty
  \mu_{a,x}^{}:\quad \int^{m\in\calm,n\in\caln}\!\!
  \Hom(\cc m\boti n\boti \cc n\boti m, \cc x_m\boti a.x_n\boti \cc y_n\boti y_m)
  \Nxl3
  \hspace*{7.5 em} \tocong
  \Hom(\int_{m\in\calm,n\in\caln}\cc {m}\boti n \boti \cc {n}\boti m,
  \cc {a^\vee.x_m}\boti x_n\boti \cc y_n \boti y_m)
  \eear
  \label{eq:iso2}
  \ee
for any $y \,{\in}\, \ZZ(\LL_2)$.
In the following picture we indicate graphically how the holonomy 
\eqref{eq:iso1} and the isomorphism \eqref{eq:iso2} for the \basicdisk\ $\DDsAMN$ arise:
    \def\locpa  {0.88}   
    \def\locpA  {2.63}   
    \def\locpb  {1.03}   
    \def\locpB  {2.35}   
    \def\locpc  {2.2pt}  
    \def\locpd  {0.12}   
  \be
  \raisebox{-6.7em} {\begin{tikzpicture}
  \scopeArrow{0.29}{>}          
  \filldraw[\colorCircle,fill=\colorSigma,line width=\locpc,postaction={decorate}]
           (0,\locpA) arc (90:270:\locpA cm)
           node[left=66pt,yshift=40pt,color=black] {$\scriptstyle-1$}
           node[left=53pt,yshift=50pt,color=\colorHol,rotate=-77] {$\hol_{a,x}$};
  \end{scope}
  \scopeArrow{0.45}{>}
  \filldraw[\colorCircle,fill=white,line width=\locpc,postaction={decorate}]
           (0,-\locpa) arc (270:90:\locpa cm)
           node[right=-25pt,yshift=-17pt,color=black]{$\scriptstyle 1$};
  \end{scope}
  \scopeArrow{0.77}{\arrowDefect}
  \draw[line width=\locpc,color=\colorFreebdy,postaction={decorate}]
           (0,\locpa) -- node[right=-2pt,yshift=-2pt]{${}_\cala\caln$} (0,\locpA);
  \end{scope}
  \scopeArrow{0.51}{\arrowDefect}
  \draw[line width=\locpc,color=\colorFreebdy,postaction={decorate}]
           (0,-\locpA) -- node[right=-2pt,yshift=10pt]{${}_\cala\calm$} (0,-\locpa);
  \end{scope}
  \filldraw[color=\colorDefect] (0,\locpa) circle (\locpd cm) (0,\locpA) circle (\locpd cm)
                                (0,-\locpa) circle (\locpd cm) (0,-\locpA) circle (\locpd cm);
  \draw[\colorHol,thick,<-] plot[smooth,tension=0.55] coordinates{
          (260:\locpb) (260:\locpB) (230:\locpB) (205:\locpB) (180:\locpB) (155:\locpB) (130:\locpB)
          (100:\locpB) (100:\locpb) };
  \begin{scope}[shift={(\locpA+3.4,0)}]
  \scopeArrow{0.29}{>}          
  \filldraw[\colorCircle,fill=\colorSigma,line width=\locpc,postaction={decorate}]
           (0,\locpA) arc (90:270:\locpA cm)
           node[left=66pt,yshift=40pt,color=black]{$\scriptstyle-1$};
  \end{scope}
  \scopeArrow{0.45}{>}
  \filldraw[\colorCircle,fill=white,line width=\locpc,postaction={decorate}]
           (0,-\locpa) arc (270:90:\locpa cm)
           node[right=-25pt,yshift=-17pt,color=black]{$\scriptstyle 1$};
  \end{scope}
  \scopeArrow{0.77}{\arrowDefect}
  \draw[line width=\locpc,color=\colorFreebdy,postaction={decorate}]
           (0,\locpa) -- node[right=-2pt,yshift=-2pt]{${}_\cala\caln$} (0,\locpA);
  \end{scope}
  \scopeArrow{0.51}{\arrowDefect}
  \draw[line width=\locpc,color=\colorFreebdy,postaction={decorate}]
           (0,-\locpA) -- node[right=-2pt,yshift=10pt]{${}_\cala\calm$} (0,-\locpa);
  \end{scope}
  \filldraw[color=\colorDefect] (0,\locpa) circle (\locpd cm) (0,\locpA) circle (\locpd cm)
                                (0,-\locpa) circle (\locpd cm) (0,-\locpA) circle (\locpd cm);
  \draw[\colorHol,thick,<-] (265:1.17*\locpa)
	   node[left=33pt,yshift=16pt,rotate=-85] {$\mu_{a,x}^{}$} arc (290:72:1.23*\locpa);
  \end{scope}
  \end{tikzpicture}}
  \ee
This prescription generalizes as follows to disks with an arbitrary number $n$ of
gluing intervals, which denote by $\LL_j$ with $j \,{=}\, 1,2,...\,,n$ (and,
to be precise, with their free boundary segments oriented,
without loss of generality, in the way displayed in the picture \eqref{eq:disk-nbdies}
below). Select one of the gluing intervals, say $\LL_i$ for some
$i \,{\in}\, \{1,2,...\,,n\}$. Then for any $x \,{\in}\, \ZZ(\LL_i)$ and any
$a \,{\in}\, \cala$ we define a balancing isomorphism $\mu_{a,x}$ in the 
same way as in \eqref{eq:iso2}. Also, by iterating the procedure in 
\eqref{eq:compositeHol}, we can again define a holonomy $\hol_{a,x}$ which,
as a natural transformation, will now have components
$(\hol_{a,x})_{y_{i+1}^{},y_{i+2}^{},...,y_{i-1}^{}}$, where $y_j\,{\in}\, \ZZ(\LL_j)$
for each $j \,{=}\, i{+}1,i{+}2,...\,,i{+}n{-}1$ (each label taken modulo $n$).

As explained in the Introduction, we will obtain the block functor from the pre-block functor
by imposing ``flatness along the disk''. This amounts to considering an equalizer of two
morphisms that are built from the isomorphisms \eqref{eq:iso1} and \eqref{eq:iso2}. 
A precise description
will be given -- for general fine \defectsurface s -- in Definition \ref{Definition:block-func}.
\end{Example}

By Proposition \ref{prop:lineflip}, the \defectsurface\ of 
Example \ref{ex:4.4} gives, up to canonical equivalence, the same gluing category
as a disk with two free boundary segments having opposite orientations and indices $2$ and
$-2$. The following example generalizes this situation to the case of arbitrary even indices.

\begin{Example}\label{Example:k-fr-cylinder}
Consider the following disk having a framing with indices ${\pm}\kappa \,{\in}\, 2\Z$ along its
two gluing segments, continued as a cylinder along the direction of the free boundary intervals:
    \def\locpa  {2.8}    
    \def\locpA  {2.3}    
    \def\locpc  {2.2pt}  
    \def\locpd  {0.12}   
  \be
  \raisebox{-3.8em}{\begin{tikzpicture}
  \scopeArrow{0.36}{<}
  \scopeArrow{0.87}{<}
  \filldraw[\colorCircle,fill=\colorSigma,line width=\locpc,postaction={decorate}]
           (0,0) node[right=21pt,yshift=-7pt,color=black]{$\kappa$} rectangle
           (\locpa,\locpA) node[left=13pt,yshift=7pt,color=black]{$-\kappa$} ; 
  \end{scope} \end{scope}
  \scopeArrow{0.77}{\arrowDefect}
  \draw[line width=\locpc,color=\colorFreebdy,postaction={decorate}]
           (0,0) -- node[left=-2pt,yshift=-4pt]{$\calm_\Cala$} (0,\locpA);
  \draw[line width=\locpc,color=\colorFreebdy,postaction={decorate}]
           (\locpa,0) -- node[right=-2pt,yshift=-4pt]{${}_\cala\caln$} +(0,\locpA);
  \end{scope}
  \filldraw[color=\colorDefect] (0,0) circle (\locpd) (0,\locpA) circle (\locpd)
                                (\locpa,0) circle (\locpd) (\locpa,\locpA) circle (\locpd);
  \end{tikzpicture}}
  \label{eq:disk.k.-k}
  \ee
The gluing categories for the lower and upper segment are, respectively,
$ \ZZ_1 \,{:=}\, \cenf{\calm \,{\stackrel{\kappa}\boxtimes}\, \caln}$ and
$ \ZZ_2 \,{:=}\, \cenf{\calnopp \,{\stackrel{-\kappa}\boxtimes}\, \calmopp} $.
Hence we have again $\ZZ_{2} \,{\cong}\, \CC{\ZZ_{1}}$, and for objects
$x_{m} \boti x_{n} \,{\in}\, \ZZ_{1}$ and $\cc{y_{n}} \boti \cc{y_{m}} \,{\in}\, \ZZ_{2}$
the pre-block space is 
  \be
  \int^{m\in\calm, n\in\caln}\! \Hom(m \boti n \boti \cc{n} \boti \cc{m}, x_{m} \boti x_{n}
  \boti \cc{y_{n}} \boti \cc{y_{m}}) \,.
  \ee
Starting the holonomy of $a \,{\in}\, \cala$ at $x_{n}$ counterclockwise, and suppressing the 
balancing of the canonical objects $\int^{m} \cc{m} \boti m$, we arrive schematically
at the situation
    \def\locpa  {2.8}    
    \def\locpA  {2.3}    
    \def\locpc  {2.2pt}  
    \def\locpd  {0.12}   
  \be
  \raisebox{-3.8em}{\begin{tikzpicture}
  \scopeArrow{0.36}{<}
  \scopeArrow{0.87}{<}
  \filldraw[\colorCircle,fill=\colorSigma,line width=\locpc,postaction={decorate}]
           (0,0) node[right=21pt,yshift=-7pt,color=black]{$\kappa$} rectangle
           (\locpa,\locpA) node[left=13pt,yshift=7pt,color=black]{$-\kappa$} ; 
  \end{scope} \end{scope}
  \scopeArrow{0.77}{\arrowDefect}
  \draw[line width=\locpc,color=\colorFreebdy,postaction={decorate}]
           (0,0) -- node[left=-2pt,yshift=-4pt]{$\calm_\Cala$} (0,\locpA);
  \draw[line width=\locpc,color=\colorFreebdy,postaction={decorate}]
           (\locpa,0) -- node[right=-2pt,yshift=-4pt]{${}_\cala\caln$} +(0,\locpA);
  \end{scope}
  \filldraw[color=\colorDefect] (0,0) circle (\locpd) (0,\locpA) circle (\locpd)
                                (\locpa,0) circle (\locpd) (\locpa,\locpA) circle (\locpd);
  \node at (-0.77,-0.02) {$ x_m.a^{[\kappa]} $};
  \node at (\locpa+0.55,-0.1) {$ a.x_n $};
  \node at (-0.96,\locpA+0.1) {$ \cc{y_m.a^{[\kappa-1]}} $};
  \node at (\locpa+0.65,\locpA+0.1) {$ \cc{a^\vee\!.y_n} $};
  \end{tikzpicture}}
  \ee
This gives the holonomy
  \be
  \begin{array}{ll}
  \multicolumn{2}{l} {\dsty \int^{m\in\calm,n\in\caln}\!\!
  \Hom( m\boti n\boti \cc n\boti\cc m, x_m\boti a.x_n\boti \cc y_n\boti \cc y_m) }
  \nxl1
  \qquad\qquad &\dsty \tocong \int^{m\in\calm,n\in\caln}\!
  \Hom( m\boti n\boti \cc n\boti\cc m, x_m\boti x_n\boti \cc{ a^{\vee}.y_n}\boti \cc y_m) 
  \Nxl1
  & \tocong\dsty \int^{m\in\calm,n\in\caln}
  \Hom( m\boti n\boti \cc n\boti\cc m, x_m\boti x_n\boti \cc y_n\boti \cc{ y_m. a^{[\kappa-1]}})
  \Nxl1
  & \tocong\dsty \int^{m\in\calm,n\in\caln}
  \Hom( m\boti n\boti \cc n\boti\cc m,  x_m. a^{[\kappa]}\boti x_n\boti \cc y_n\boti \cc y_m) \,.
  \eear
  \ee
\end{Example}

The powers of duals appearing in the so obtained expressions for the holonomy turn out to
be significant: they will allow us to define the block functors as equalizers
(Definition \ref{Definition:block-func}).

\begin{Definition} \label{def:well-posed}
We say that the holonomy problem for a disk labeled by a finite tensor category $\cala$
is \emph{well-posed} iff for any $a \,{\in}\, \cala$ and any object $x$ in the gluing category
associated with a gluing segment of the boundary the holonomy $\hol_{a,x}$ and the
isomorphism $\mu_{a,x}$ provided by the balancing of $x$ 
coincide up to a double right dual in the same way as the isomorphisms 
\eqref{eq:iso1} and \eqref{eq:iso2} in the example of the \basicdisk\ $\DDsAMN$.
\\
The holonomy problem for a fine \defectsurface\ $\Sigma$ is said to be well-posed iff 
the holonomy problem for every \twocell\ of $\Sigma$ is well-posed.
\end{Definition}

By inspection, the holonomy is well-posed both in Example \ref{ex:4.4} and
in Example \ref{Example:k-fr-cylinder}.
Our next task will be to show that also in the general situation 
of a fine \defectsurface\ all holonomy problems are well-posed. To see this, we make use
of the following statement about framing indices.

\begin{Lemma}\label{lem:4.5}
Let $\DD$ be a $2$-framed disk with $N$ gluing boundary segments $s_i$ which are
oriented as induced by the orientation of $\DD$. Then
  \be
  \sum_{i=1}^N \ind(s_i) \,=\, N -2 \,. 
  \label{-2=sum-1}
  \ee
\end{Lemma}

\begin{proof} 
Consider the standard disk $ \DD_{\text{std}} \,{\subset}\, \R^{2}$ with marked points on
its boundary that divide $\partial\DD_{\text{std}}$ into segments. That the disk has Euler
characteristic 1 implies that for any non-vanishing continuous vector field on $\DD_{\text{std}}$
the sum of the indices of all segments on the boundary $\partial\DD_{\text{std}}$ is $-2$
when all segments are oriented counterclockwise. To be able to apply this fact to the situation 
at hand, we must smoothen out the framing at the corners of the defect disk $\DD$, at which 
a free boundary interval is adjacent to a gluing boundary. Thus we replace the two allowed
situations
    \def\locpa  {3.1}    
    \def\locpb  {2.1}    
    \def\locpc  {2.2pt}  
    \def\locpd  {0.12}   
    \def\locpe  {0.71}   
  \be
  \raisebox{-3.3em}
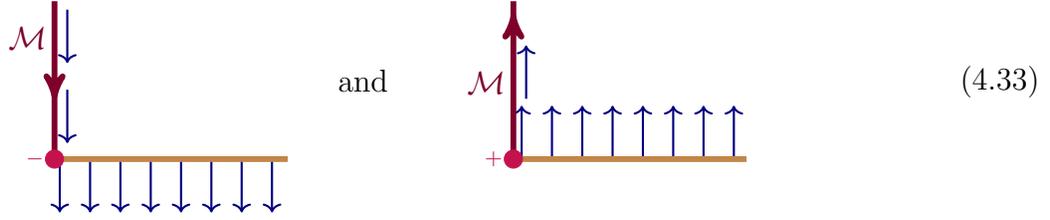
{\begin{tikzpicture}[scale=1.1]
  \foreach \xx in {0,1,...,7} \draw[thick,color=\colorFrame,->]
               (0.07+\xx*0.13*\locpa,0) -- +(0,-\locpe);
  \draw[thick,color=\colorFrame,->] (0.17,1.3*\locpe) -- +(0,-\locpe);
  \draw[thick,color=\colorFrame,->] (0.17,2.8*\locpe) -- +(0,-\locpe);
  \draw[line width=\locpc,draw=\colorCircle] (0,0) -- (\locpa,0);
  \scopeArrow{0.62}{\arrowDefect}
  \draw[line width=\locpc,color=\colorFreebdy,postaction={decorate}]
           (0,\locpb) -- node[left=-2pt,yshift=16pt]{$\calm$} (0,0);
  \end{scope}
  \filldraw[color=\colorDefect] (0,0) circle (\locpd) node[left]{$\scriptstyle -$};
  \begin{scope}[shift={(4.1,0.5*\locpb)}] \node at (0,0) {and}; \end{scope}
  \begin{scope}[shift={(6.1,0)}]
  \foreach \xx in {0,1,...,7} \draw[thick,color=\colorFrame,->]
               (0.11+\xx*0.13*\locpa,0) -- +(0,\locpe);
  \draw[thick,color=\colorFrame,->] (0.17,1.13*\locpe) -- +(0,\locpe);
  \draw[line width=\locpc,draw=\colorCircle] (0,0) -- (\locpa,0);
  \scopeArrow{0.92}{\arrowDefect}
  \draw[line width=\locpc,color=\colorFreebdy,postaction={decorate}]
           (0,0) -- node[left=-2pt,yshift=-1pt]{$\calm$} (0,\locpb);
  \end{scope}
  \filldraw[color=\colorDefect] (0,0) circle (\locpd) node[left]{$\scriptstyle +$};
  \end{scope}
  \end{tikzpicture}}
  \label{eq:smooth1}
  \ee
by the smoothened versions
    \def\locpa  {3.1}    
    \def\locpb  {2.1}    
    \def\locpc  {2.2pt}  
    \def\locpd  {0.12}   
    \def\locpe  {0.71}   
    \def\locpf  {0.9}    
    \def\locpg  {23}     
  \be
  \raisebox{-3.3em}
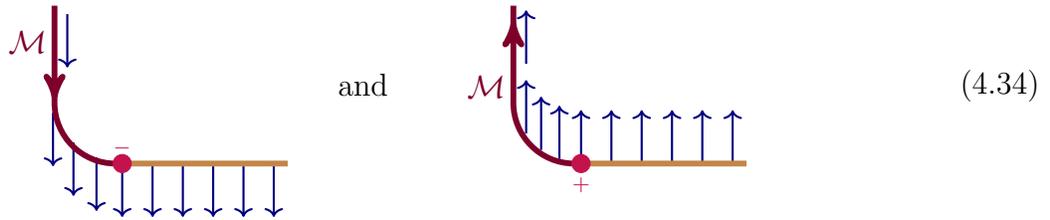
{\begin{tikzpicture}[scale=1.1]
  \foreach \xx in {0,1,...,5} \draw[thick,color=\colorFrame,->]
               (\locpf+\xx*0.13*\locpa,0) -- +(0,-\locpe);
  \draw[thick,color=\colorFrame,->] (0.17,2.8*\locpe) -- +(0,-\locpe);
  \draw[thick,color=\colorFrame,->] (-0.02,0.95*\locpe) -- +(0,-\locpe);
  \draw[thick,color=\colorFrame,->] (0.28*\locpf,0.28) -- +(0,-\locpe);
  \draw[thick,color=\colorFrame,->] (0.62*\locpf,0.08) -- +(0,-\locpe);
  \draw[line width=\locpc,draw=\colorCircle] (\locpf,0) -- (\locpa,0);
  \scopeArrow{0.47}{\arrowDefect}
  \draw[line width=\locpc,color=\colorFreebdy,rounded corners=\locpg pt,postaction={decorate}]
               (0,\locpb) -- node[left=-2pt,yshift=16pt]{$\calm$} (0,0) -- (\locpf,0);
  \end{scope}
  \filldraw[color=\colorDefect] (\locpf,0) circle (\locpd) node[above=-1pt]{$\scriptstyle -$};
  \begin{scope}[shift={(4.1,0.5*\locpb)}] \node at (0,0) {and}; \end{scope}
  \begin{scope}[shift={(6.1,0)}]
  \foreach \xx in {0,1,...,5} \draw[thick,color=\colorFrame,->]
               (\locpf+\xx*0.13*\locpa,0) -- +(0,\locpe);
  \draw[thick,color=\colorFrame,->] (0.17,1.87*\locpe) -- +(0,\locpe);
  \draw[thick,color=\colorFrame,->] (0.17,0.40) -- +(0,\locpe);
  \draw[thick,color=\colorFrame,->] (0.41*\locpf,0.18) -- +(0,\locpe);
  \draw[thick,color=\colorFrame,->] (0.68*\locpf,0.06) -- +(0,\locpe);
  \draw[line width=\locpc,draw=\colorCircle] (\locpf,0) -- (\locpa,0);
  \scopeArrow{0.92}{\arrowDefect}
  \draw[line width=\locpc,color=\colorFreebdy,rounded corners=\locpg pt,postaction={decorate}]
               (\locpf,0) -- (0,0) -- node[left=-2pt,yshift=-1pt]{$\calm$} (0,\locpb);
  \end{scope}
  \filldraw[color=\colorDefect] (\locpf,0) circle (\locpd) node[below=1pt]{$\scriptstyle +$};
  \end{scope}
  \end{tikzpicture}}
  \label{eq:smooth2}
  \ee
respectively, in which the free boundary segments are suitably deformed. Similarly,
in case we deal with a defect line, labeled by $\calm$, in the interior of $\DD$ rather
than a free boundary interval, we temporarily think of the defect line as a pair of
parallel free boundary intervals and apply the procedure above to both parts. 
Both of the situations \eqref{eq:smooth1} and \eqref{eq:smooth2} contribute,
in a counterclockwise sense, $-\frac12$. Since there are two arcs per
segment, each segment gets an additional contribution $-1$.  This way we obtain the equality
$-2 \,{=}\, \sum_{i=1}^N\! \big( \ind(s_i) - 1\big)$, thus proving \eqref{-2=sum-1}.
\end{proof}                     

\begin{Remark}\label{rem:framing}
In fact, vector fields on $\DD$ up to homotopy are in bijection with collections of indices that
obey the relation \eqref{-2=sum-1}. To see this, let $\DD$ be a disk whose boundary
is either a gluing circle (in case that $N \,{=}\, 0$) or consists of $N \,{\ge 1}\,$
free boundaries and $N$ gluing segments $s_i$. All gluing segments are
endowed with the induced orientation. The orientation of the free boundaries adds signs
to the points at which free boundaries and gluing intervals meet. Let $\kappa$ be a tuple
of framing indices for these signs such that $\sum_{i=1}^N \kappa_i - N \,{=}\, {-}2$. Then 
there exists a vector field on $\DD$ that is continuous on $\DD$,
parallel to the free boundaries and has the given index on the gluing segments, 
i.e.\ $\ind(s_i) \,{=}\, \kappa_i$ for each segment.
\end{Remark}

\begin{Proposition}\label{prop:wellposedD}
The holonomy problem of any labeled $2$-framed disk $\DD$ is well-posed. 
\end{Proposition} 

\begin{proof}
Consider a disk $\DD$. We invoke the line flip of Proposition \ref{prop:lineflip}
to assume without loss of generality that all free boundary segments on $\partial\DD$
are oriented as induced by the orientation of $\DD$.
It is straightforward to see that the holonomy problem is well-posed
for one collection of indices satisfying the sum rule \eqref{-2=sum-1} iff it is
well-posed for any other. It follows that we can restrict our attention to the situation
    \def\locpa  {0.4}    
    \def\locpb  {1.1}    
    \def\locpc  {2.2pt}  
    \def\locpd  {0.10}   
    \def\locpe  {2.9}    
    \def\locpf  {4.9}    
  \be
  \raisebox{-11.1em}
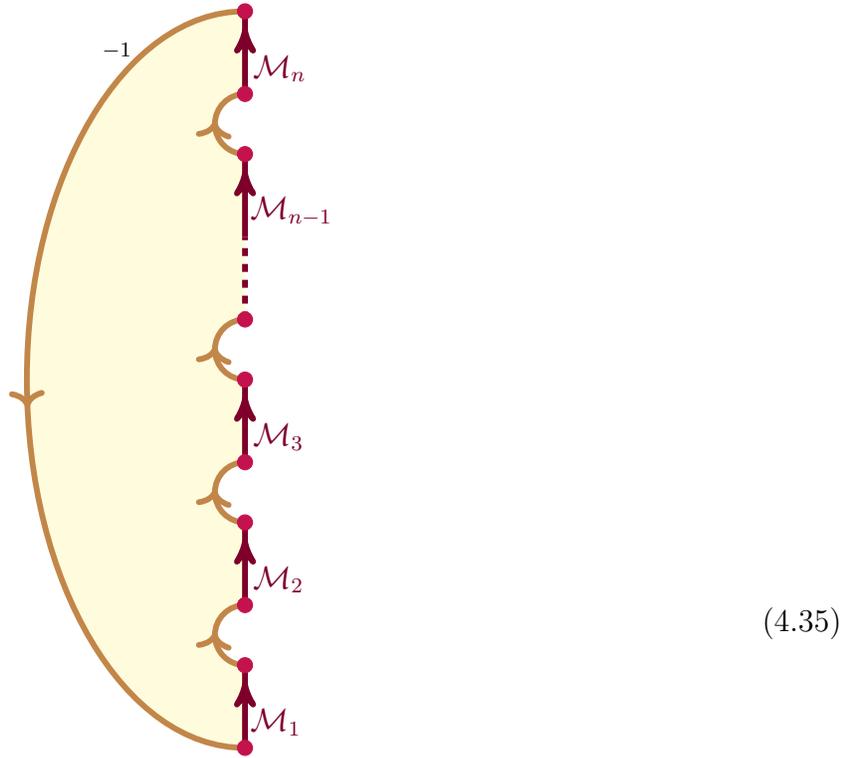
{\begin{tikzpicture}
  \scopeArrow{0.53}{>}
  \filldraw[line width=\locpc,\colorCircle,fill=\colorSigma,postaction={decorate}]
           (0,7*\locpa+5*\locpb)
           node[left=38pt,yshift=-15pt,color=black]{$\scriptstyle-1$}
           arc (90:270:\locpe cm and \locpf cm);
  \filldraw[\colorCircle,fill=white,line width=\locpc,postaction={decorate}]
           (0,-\locpa) arc (270:90:\locpa);
  \filldraw[\colorCircle,fill=white,line width=\locpc,postaction={decorate}]
           (0,\locpa+\locpb) arc (270:90:\locpa);
  \filldraw[\colorCircle,fill=white,line width=\locpc,postaction={decorate}]
           (0,3*\locpa+2*\locpb) arc (270:90:\locpa);
  \filldraw[\colorCircle,fill=white,line width=\locpc,postaction={decorate}]
           (0,5*\locpa+4*\locpb) arc (270:90:\locpa);
  \end{scope}
  \scopeArrow{0.80}{\arrowDefect}
  \draw[line width=\locpc,color=\colorFreebdy,postaction={decorate}]
           (0,-\locpa-\locpb) -- node[right=-3pt,yshift=-6pt]{$\calm_1$} +(0,\locpb);
  \draw[line width=\locpc,color=\colorFreebdy,postaction={decorate}]
           (0,\locpa) -- node[right=-2pt,yshift=-6pt]{$\calm_2$} +(0,\locpb);
  \draw[line width=\locpc,color=\colorFreebdy,postaction={decorate}]
           (0,3*\locpa+\locpb) -- node[right=-2pt,yshift=-6pt]{$\calm_3$} +(0,\locpb);
  \draw[line width=\locpc,color=\colorFreebdy,postaction={decorate}]
           (0,5*\locpa+3*\locpb) -- node[right=-3pt,yshift=-6pt]{$\calm_{n-1}$} +(0,\locpb);
  \draw[line width=\locpc,color=\colorFreebdy,postaction={decorate}]
           (0,7*\locpa+4*\locpb) -- node[right=-2pt,yshift=-6pt]{$\calm_{n}$} +(0,\locpb);
  \end{scope}
  \draw[dashed,line width=\locpc,color=\colorFreebdy] (0,5*\locpa+2*\locpb) -- +(0,\locpb);
  \filldraw[color=\colorDefect] 
           (0,-\locpa-\locpb) circle (\locpd)
           (0,-\locpa) circle (\locpd) (0,\locpa) circle (\locpd)
           (0,\locpa+\locpb) circle (\locpd) (0,3*\locpa+\locpb) circle (\locpd)
           (0,3*\locpa+2*\locpb) circle (\locpd) (0,5*\locpa+2*\locpb) circle (\locpd)
           (0,5*\locpa+4*\locpb) circle (\locpd) (0,7*\locpa+4*\locpb) circle (\locpd)
           (0,7*\locpa+5*\locpb) circle (\locpd);
  \end{tikzpicture}}
  \label{eq:disk-nbdies}
  \ee
corresponding to a framing given by the constant vector field pointing upwards
(recall that the index of all unlabeled oriented gluing intervals is $+1$).
Now for this situation well-posedness follows from Example \ref{ex:4.4} together with the 
observation that for each $i \,{\in}\, \{1,2,...\,,n{-}1\}$ we have
(abbreviating $\calm \,{=}\, \calm_i$ and $\caln \,{=}\, \calm_{i+1}$) 
    \def\locpa  {0.52}   
    \def\locpb  {1.4}    
    \def\locpB  {1.1}    
    \def\locpc  {2.2pt}  
    \def\locpd  {0.12}   
    \def\locpe  {1.4}    
  \be
  \raisebox{-7.1em}
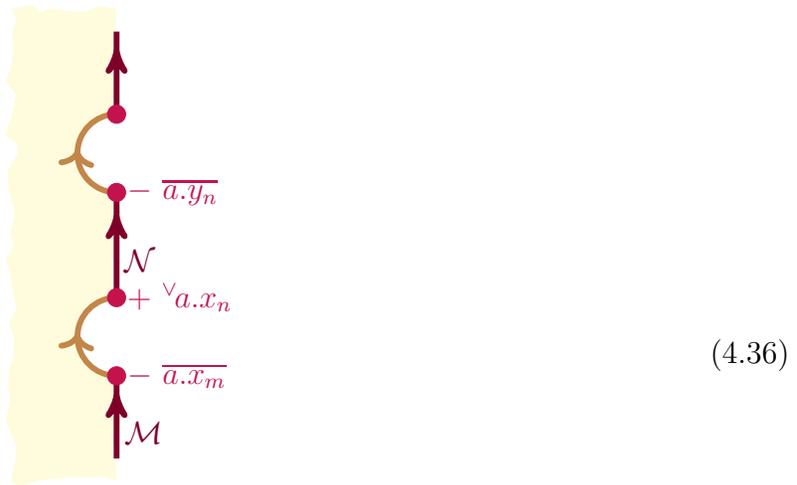
{\begin{tikzpicture}
  [decoration={random steps,segment length=2mm}]
  \fill[\colorSigma,decorate] (0,-\locpa-\locpb) rectangle (-\locpe,3*\locpa+2*\locpb);
  \fill[white] (0,0) circle (\locpa) (0,2*\locpa+\locpb) circle (\locpa);
  \scopeArrow{0.53}{>}
  \draw[line width=\locpc,\colorCircle,postaction={decorate}] (0,-\locpa) arc (270:90:\locpa);
  \draw[line width=\locpc,\colorCircle,postaction={decorate}] (0,\locpa+\locpb) arc (270:90:\locpa);
  \end{scope}
  \scopeArrow{0.80}{\arrowDefect}
  \draw[line width=\locpc,color=\colorFreebdy,postaction={decorate}]
           (0,3*\locpa+\locpb) -- +(0,\locpB);
  \draw[line width=\locpc,color=\colorFreebdy,postaction={decorate}]
           (0,\locpa) -- node[right=-2pt,yshift=-6pt]{$\caln$} +(0,\locpb);
  \draw[line width=\locpc,color=\colorFreebdy,postaction={decorate}]
           (0,-\locpa-\locpB) -- node[right=-2pt,yshift=-6pt]{$\calm$} +(0,\locpB);
  \end{scope}
  \filldraw[color=\colorDefect] 
           (0,-\locpa) circle (\locpd) node[right]{$-~\cc{a.x_m}$}
           (0,\locpa) circle (\locpd) node[right]{$+~\Vee a.x_n$}
           (0,\locpa+\locpb) circle (\locpd) node[right]{$-~\cc{a.y_n}$}
           (0,3*\locpa+\locpb) circle (\locpd);
  \end{tikzpicture}}
  \ee
so that no double duals arise in the part of the holonomy along the gluing segments with index 1.
\end{proof}

An immediate consequence is

\begin{cor} \label{cor:well-posed}
The holonomy problem of any disk $\DD$ in any \defectsurface\ $\Sigma$ is well-posed.
\end{cor}

In particular, the holonomy problems for the \basicdisk\ in Example \ref{ex:4.4} and for
the disks discussed in Example \ref{Example:k-fr-cylinder} are well-posed.


\subsection{Block functors for fine surfaces}\label{sec:blocks}

We are now in a position to set up functors that will eventually provide us with the 
block functors for fine \defectsurface s. Since the holonomy problem for the \basicdisk\
in Example \ref{ex:4.4} is well-posed, the holonomy operators
  \begin{eqnarray}
  \Zpre(\cc x_m\boti a.x_n\boti \cc y_n\boti y_m) \,= \int^{m\in\calm,n\in\caln}\!\!
  \Hom(\cc m\boti n\boti \cc n\boti m, \cc x_m\boti a.x_n\boti \cc y_n\boti y_m)~~
  \nonumber\\
  \xrightarrow[\cong]{~\hol_{a,x}^{}~}
  \int^{m\in\calm,n\in\caln}\!\! \Hom(\cc {m}\boti n \boti \cc {n}\boti m,
  \cc {a^{\vee\vee\vee}\!.\,x_m}\boti x_n\boti \cc y_n \boti y_m) \qquad
  \nonumber\\
  = \Zpre(\cc {a^{\vee\vee\vee}\!.\,x_m} \boti x_n\boti \cc y_n \boti y_m) \qquad
  \end{eqnarray}
allow us to formulate holonomy equations and thus to define blocks $\Zfine$ for
fine \defectsurface s as equalizers, see Definition \ref{Definition:block-func} below.
In contrast, for non-fine surfaces this is no longer possible.
As a consequence, to define the block functor $\ZZ$ for an arbitrary \defectsurface\
$\Sigma$ we will have to take a limit over suitable `fine refinements' of $\Sigma$; 
we relegate this to Definition \ref{def:THEdef}. In case $\Sigma$ is already fine,
the functor $\Zfine(\Sigma)$ from Definition \ref{Definition:block-func}
can be taken as a distinguished representative of the limit $\ZZ(\Sigma)$. 
This justifies an abuse of language: we
refer also to the functor $\Zfine(\Sigma)$ as the block functor for $\Sigma$.
	  
The following result will later allow us to
relate certain spaces of blocks to natural transformations.

\begin{Lemma}
\label{lemma:gen-equalizer}
Let $T\colon \calc \,{\to}\, \calc$ be a monad on a linear category, with category $\calc_T$
of modules, and denote by $U\colon \calc_{T}\,{\to}\, \calc$ the forgetful functor. Then for
any pair $m \,{=}\, (Um,\rho_m)$ and $n \,{=}\, (Un,\rho_n)$ of objects of $\calc_{T}$,
with actions $\rho_m$ and $\rho_n$, respectively, there is an equalizer diagram
  \be
  \Hom_{\calc_T}(m,n) \xrightarrow{~~~} \Hom_{\calc}(Um,Un)  
  \xrightrightarrows[~~\varphi_2~]{~~\varphi_1~} \Hom_{\calc}(T(Um),Un) \,,
  \label{eq:exSequ}
  \ee
where on a morphism $\gamma\colon U(m)\,{\to}\, U(n)$ the maps $\varphi_{1,2}$ are 
defined by $\varphi_1(\gamma) \,{:=}\, \gamma\,{\circ}\,\rho_m$
and $\varphi_2(\gamma) \,{:=}\, \rho_n \,{\circ}\, T(\gamma)$, respectively.
 \\
Similarly, for $S$ a comonad on $\calc$ and 
$x \,{=}\, (Ux,\delta_x)$ and $y \,{=}\, (Uy,\delta_y)$
objects in the category $\calc  ^{S}$ of $S$-comodules, there is an equalizer diagram
  \be
  \Hom_{\calc^S}(x,y) \xrightarrow{~~~} \Hom_{\calc}(Ux,Uy) 
  \xrightrightarrows[~~\psi_2~]{~~\psi_1~}\Hom_{\calc}(Ux,S(Uy))
  \label{eq:coexSequ}
  \ee
with $\psi_1(\gamma) \,{:=}\, \delta_y \,{\circ}\, \gamma$
and $\psi_2(\gamma) \,{:=}\, S(\gamma)\,{\circ}\, \delta_x$
for $\gamma\colon U(x)\,{\to}\, U(y)$.
\end{Lemma}

\begin{proof}
It follows from the definition that $\Hom_{\calc_{T}}(m,n)$ is the kernel of 
$\varphi_{1}\,{-}\,\varphi_{2}$: $\gamma$ is a module map if and only if
$\varphi_{1}(\gamma) \,{=}\, \varphi_{2}(\gamma)$. The statement for the comonad
follows by the same type of reasoning.
\end{proof}

We now consider, for a disk labeled by a finite tensor category $\cala$,
for any $a\iN\cala$ the composite
  \be
  \bearll
  \Zpre(...\, \cc{x_m} \boti x_n \boti y) \xrightarrow{~\coevl_*~} \!\! &
  \Zpre(...\, \cc{x_m} \boti (\Vee a \,{\otimes}\, a)\,.\,x_n \boti y)
  \Nxl2 &
  \xrightarrow[\cong]{~\hol~}
  \Zpre(...\, \cc{a^{\vee\vee}\!.\,x_{m}} \boti a.x_{n} \boti y)
  \label{eq:hola}
  \eear
  \ee
of the post-composition with the left coevaluation of $a$ and the holonomy
$\hol_{\Vee a,a.x_n}$ of $\Vee a$ around the disk. It is important to note that these maps
are dinatural in $a$ and thus factorize over the end.  We thus obtain a morphism of
left exact functors
  \be
  \hol_x: \quad \Zpre(...\, \cc{x_m} \boti x_n \boti y) 
  \,\xrightarrow{~~} \int_{a\in\cala}\!
  \Zpre(...\, \cc{ a^{\vee \vee} \!.\,x_m} \boti a\,.\,x_n \boti y) 
  \label{eq:hol}
  \ee
furnishing a well-posed holonomy problem. On the other hand, combining the right coevaluation
and the comodule structure of $\cc{x_{m}} \boxtimes x_{n}$ that expresses the balancing
provides us with a morphism $\cc{x_m} \boti x_n \,{\to} 
\int_{a\in\cala}\cc{a^{\vee \vee} \!.\,x_m} \boti a\,.\,x_n$, which defines another map 
  \be
  \begin{array}{r}\dsty
  (\mu_{x})_*^{} : \quad \Zpre(...\, \cc{x_m} \boti x_n \boti y)
  \xrightarrow{~\coevr_*~} 
  \int_{a \in \cala}\! \Zpre(...\, \cc{x_m} \boti (a\oti a^\vee)\,.\,x_n \boti y) \qquad
  \Nxl2 \dsty
  \xrightarrow{~\cong~} \int_{a \in \cala}\!
  \Zpre(...\, \cc{ a^{\vee\vee}\!.\,x_m} \boti a\,.\,x_n \boti y)
  \eear
  \ee
(the end over $\cala$ commutes with all coends in 
the pre-block functor, since those can be pulled in the first argument of the Hom).
We can thus define the \emph{block space} as the equalizer of the maps $\hol_x$ and 
$(\mu_{x})_{*}$. More precisely, we impose one such relation for each \twocell\
(which, as the surface is assumed to be fine, is a disk) and select for each 
\twocell\ $\PP_p$ a starting point $v_p$ among the defect points on $\partial\PP_p$.

\begin{Definition} \label{Definition:block-func}
Let $\Zprev$ be the pre-block functor associated with a fine \defectsurface.
The functor $\ZZ_{\!\text{fine}}$ associated with the surface is the equalizer
  \be
  \bearll
  \ZZ_{\!\text{fine} }(...\, \cc{x_m} \boti x_n \boti y) \,\xrightarrow{~~~} \!&
  \Zpre(...\, \cc{x_m} \boti x_n \boti y) 
  \Nxl3 &\dsty \hspace*{2.5em}
  \xrightrightarrows[~\prod({\mu_x)}_*] {~\prod\hol_x}\, \prod_p\! \int_{\!a\in\cala_p}
  \!\!\!\!  \Zpre(...\, \cc{a^{\vee\vee}\!.\, x_m} \boti a\,.\,x_n \boti y) \,,
  \eear
  \label{eq:defining-block}
  \ee
where the product is over all \twocell es $\PP_p$ of the \defectsurface\ and
$\cala_p$ is the finite tensor category labeling $\PP_p$.
\end{Definition}

As mentioned above, we call $\ZZ_{\!\text{fine}}$ the \emph{block functor} for
the \defectsurface, albeit a complete formulation that is valid for arbitrary
\defectsurface s will require to define the actual block functor $\ZZ$ as a limit.
While the pre-blocks for a \defectsurface\ $\Sigma$ only depend on the incidence combinatorics 
of boundary segments and defect lines of $\Sigma$, the blocks also depend, via the
holonomy, directly on the framing on the \twocells\ of $\Sigma$.

A priori the block functor $\ZZ_{\!\text{fine}} \,{=}\, \ZZ_{\!\text{fine},(v)}$
for a \defectsurface\ $\Sigma$ depends on the choice
$(v) \,{=}\, \{v_i\}$ of a starting point $v_i$ for each disk $\DD_i$ in $\Sigma$.
However, in fact it depends on these choices only up to canonical coherent natural 
isomorphism. To see this, we provide an alternative characterization of the equalizer: 
Composing Equation \eqref{eq:hola} with the balancing of $x$, we obtain a morphism 
  \be
  \widetilde{\hol_{a,x}}: \quad \Zpre(...\, \cc{x_m} \boti x_n \boti y) 
  \xrightarrow{~~} \Zpre(...\, \cc{x_m} \boti (\Vee a{\otimes}a) \,.\, x_n \boti y) \,,
  \label{eq:alternative-hol}
  \ee
which has as parallel morphism the composition with $\coevl_{a}$ at $x_{n}$. Again, these
morphisms factorize over the end, and by composing the defining equation \eqref{eq:defining-block}
with the balancing of $x$ we see that the block space is also the equalizer
  \be
  \bearl
  \ZZ_{\!\text{fine},(v)} (...\, \cc{x_{m}} \boti x_{n} \boti y) \xrightarrow{~~~} 
  \Zpre(...\, \cc{x_m} \boti x_n \boti y) \quad
  \Nxl1 \dsty \hspace*{11.9em}
  \xrightrightarrows[~\prod(\coevl)_{*}]{~\prod \widetilde{\hol_x}}
  \prod_x\!\int_{a\in\cala}\! \Zpre(...,\cc{x_m} \boti (\Vee a{\otimes}a)\,.\, x_n \boti y) \,.
  \eear
  \label{eq:alternative-block}
  \ee
We use the latter description of the block space to show 

\begin{Lemma} \label{lemma:independent-start}
The block functor depends on the choice of a starting point per disk only
up to canonical coherent natural isomorphism.
\end{Lemma}

\begin{proof}
We define canonical natural isomorphisms
  \be
  \Gamma_{\!v',v}: \quad \ZZ_{\!\text{fine},(v)}(...\,,\cc{x_{m}} \boti x_{n} \boti y)
  \xrightarrow{~} \ZZ_{\!\text{fine},(v')}(...\,,\cc{x_{m}} \boti x_{n} \boti y) \,,
  \label{eq:defGamma}
  \ee
for each pair of starting points $v,v'$ per disk, that satisfy the coherence
relation $ \Gamma_{\!v'',v'} \,{\circ}\, \Gamma_{\!v',v}$
	                                                ${=}\, \Gamma_{\!v'',v}$.
For a given disk with a choice of starting points $v,v'$ there is an isomorphism 
  \be
  \gamma_{v',v}:\quad \Zpre(...\,,\cc{x_m} \boti (\Vee a \,{\otimes}\, a) \,.\, x_n \boti y)
  \xrightarrow{~\cong~}
  \Zpre(...\,,\cc{x_m} \boti x_n \boti (\Vee a \,{\otimes}\, a) \,.\, y) \,,
  \ee
to which we refer as the \emph{parallel transport operation} from $v$ to $v'$. Here 
the first action of $\Vee a \,{\otimes}\, a$ is at the defect point $v$ and the second at
$v'$, which is constructed using the balancings precisely as in the holonomy operation following
a positive path along the boundary of the disk from $v$ to $v'$. The latter ensures
that the parallel transport operations $\gamma_{v',v}$ are coherent, and it also implies
that there are two commuting triangles of isomorphisms
  \be
  \begin{tikzcd}[row sep=3.4em,column sep=3.2em]
  \Zpre(...\,,\cc{x_m} \boti x_n \boti y)  
  \ar[shift left=0.8ex]{r}[above]{\widetilde{\hol_{a,v}}~}
  \ar[shift right=0.8ex]{r}[below]{\coevl~} 
  \ar[shift left=1ex,yshift=-3pt]{dr}[above,xshift=21pt,yshift=-7pt]{\widetilde{\hol_{a,v'}}}
  \ar[shift right=1.5ex,yshift=5pt]{dr}[below,xshift=-8pt]{\coevl}
  & \Zpre(...\,,\cc{x_m} \boti (\Vee a \,{\otimes}\, a) \,.\, x_n \boti y)
  \ar{d}{\gamma_{v',v}} 
  \\
  ~ & \Zpre(...\,,\cc{x_m} \boti x_n \boti (\Vee a \,{\otimes}\, a) \,.\, y 
  \end{tikzcd}
  \ee
Thus we obtain $\Gamma_{\!v',v}$ as the universal isomorphism between the corresponding 
equalizers and it inherits the coherence from $\gamma_{v',v}$.
\end{proof}

This result justifies to disregard the dependence of $\Zfine$ on the choice of starting points
in the sequel. More conceptionally, one can define $\Zfine(...\,,\cc{x_m} \boti x_n \boti y)$
as the direct limit over the isomorphisms $\Gamma_{\!v',v}$ defined in \eqref{eq:defGamma}.

The definition also implies directly that isomorphic \defectsurface s
give identical block functors:

\begin{Lemma} \label{lemma:naive-diffeo}
Let $\phi\colon \Sigma \,{\to}\, \Sigma'$ be an isomorphism of fine \defectsurface s. Then the
block functors for $\Sigma$ and $\Sigma'$ are equal on the nose,
  \be
  \Zfine(\Sigma) = \Zfine(\Sigma') \,.
  \ee
\end{Lemma}

\begin{proof}
We can identify the defects and the vector fields of $\Sigma$ with those of $\Sigma'$ via the
isomorphism $\phi$. Thus the pre-block functors as well as the holonomy operations for
$\Sigma$ and $\Sigma'$ coincide: they only depend on the incidence relations of the patches of 
various dimensions, and these are not changed by an isomorphism.
\end{proof}

\medskip

Next we use Lemma \ref{lemma:gen-equalizer} to show that block spaces produce spaces of
natural transformations in specific situations. We first note that the 
boundary-segment-flipping lemma \ref{lem:lineflipPreblocks} for pre-blocks extends to blocks:

\begin{Proposition}\label{prop:lineflipBlocks}
There is a canonical isomorphism between the block functors for the two configurations
\eqref{eq:lineflipPreblocks} that appear in Lemma $\ref{lem:lineflipPreblocks}$.
\end{Proposition}

\begin{proof}
We know from Lemma \ref{lem:lineflipPreblocks} that the two pre-block functors are 
canonically isomorphic. Moreover, the balancings on the segments labeled by $\caln$ and 
by ${}_{}^{1}\cc{\caln}$, respectively, agree by the definition of the respective actions. 
Thus the block functors are isomorphic as well.
\end{proof}

Similarly, our construction is compatible with composition of left exact functors via
the Eilenberg--Watts correspondence and thus with the fusion of boundary insertions
not only at the level of pre-blocks, but also for blocks:

\begin{Proposition}\label{prop:block-compose}
The isomorphism between the pre-block functors for the two configurations 
\eqref{eq:PsiG.PsiF=PsiGF} established in Proposition {\rm\ref{prop:preblock-compose}}
is compatible with the holonomy operators and hence induces a canonical
isomorphism of the corresponding block functors. 
\end{Proposition}

\begin{proof}
For left exact bimodule functors $F\colon \calm \,{\to}\, \caln$ and
$G\colon \caln \,{\to}\, \calk$ for $\calb$-$\cala$-bimodules $\calm$, $\caln$ and
$\calk$, on the level of pre-blocks the isomorphism is described by the composite
\eqref{eq:holonomy-fusion-line}.
By the definition of the balancings of $\Psile(F)$ and $\Psile(G)$ according to Equation 
\eqref{eq:def-bal-le}, the first two isomorphisms in this composite are compatible with the 
holonomy operation. 
We need to show that the isomorphism 
  \begin{eqnarray}
  \int^{n \in \caln}\! \Hom(k,G(n)) \otimes \Hom(n,F(m)) 
  \cong \int^{n \in \caln}\! \Hom(G^{\la}(k),n) \otimes \Hom(n,F(m)) \quad
  \nonumber\\[0.1em]
  \cong\, \Hom(G^{\la}(k),F(m)) \,\cong\, \Hom(k, G{\circ}F(m)) \qquad
  \label{eq:fusin-N}
  \end{eqnarray}
is compatible with the balancing structures. Therefore we consider the diagram 
  \be 
  \hspace*{-.6em}
  \begin{tikzcd}[column sep=0.2cm, row sep=0.7cm]
  \dsty\int^{n\in\caln}\!\!\! \Hom(k,G(n)) \otimes \Hom(n,F(m.a))
  \ar[end anchor={[xshift=-3.2em]}]{rrr}
  \ar[start anchor={[yshift=1.55ex]},end anchor={[yshift=-1.8ex]}]{d}
  & ~ & ~ & \hspace*{-3.2em}\Hom(k, GF(m.a)) 
  \ar[xshift=-0.8em]{ddd} \ar[xshift=-1.8em]{dl}
  \\
  \dsty\int^{n\in\caln}\!\!\! \Hom(G^{\la}(k),n) \otimes \Hom(n,F(m ) . a) \ar{rr}
  \ar[start anchor={[yshift=1.55ex]},end anchor={[yshift=-1.8ex]}]{d}
  & ~ & \Hom(G^{\la}(k),F(m) . a) \ar{d} & ~
  \\
  \dsty\int^{n\in\caln}\!\!\! \Hom(G^{\la}(k) . {}^{\vee\!}a,n)\otimes \Hom(n,F(m)) \ar{rr} 
  \ar[start anchor={[yshift=1.55ex]},end anchor={[yshift=-1.8ex]}]{d}
  & ~ & \Hom(G^{\la}(k) . \Vee a,F(m)) \ar[xshift=-1.8em]{dr} & ~
  \\
  \dsty\int^{n \in \caln}\!\!\! \Hom(k , \Vee a.G(n)) \otimes \Hom(n,F(m))
  \ar[end anchor={[xshift=-3.2em]}]{rrr}
  & ~ & ~ & \hspace*{-3.2em}\Hom(k . \Vee a, GF(m))
  \\[-2.5em] ~
  \end{tikzcd}
  \ee
  ~\\[-.2em]
of isomorphisms. Here the horizontal arrow in the top and bottom row are variants of the
isomorphism \eqref{eq:fusin-N}; the commutativity of the inner rectangle is the definition
in Equation \eqref{eq:def-bal-le} of the balancing structure of the coend over $\caln$, 
and all other arrows are composites of the dualities and the module structures of $F$ and
$G$. It thus follows directly that the whole diagram commutes and hence
\eqref{eq:holonomy-fusion-line} is compatible with the holonomy operation.  
\end{proof}

Analogously, the distinguished isomorphism used in Corollary \ref{cor:simpledisk0} is
compatible with the holonomy operators and hence induces an isomorphism of the block functors.
Thus we have

\begin{cor}\label{cor:simpledisk}
The block functor for any disk without defect points and with an arbitrary number of
free and gluing segments is isomorphic, by a distinguished isomorphism, to the block
functor for the \basicdisk\ \eqref{pic:ex:4.4}.
\end{cor}

Next we develop a conceptual formulation of the holonomy operations, which will in particular 
prove to be helpful later on, when we explore refinements of \defectsurface s. When doing so
we must account for the possibility that when performing a holonomy operation we move along one 
and the same defect line twice, in opposite directions. This is achieved by keeping track of 
normal directions, in the following manner. It suffices to
consider the case of a single disk $\DD$ in a \defectsurface\ $\Sigma$. We then consider
the set $\Edg$ consisting of all defect points and all free boundary segments of $\DD$,
where the latter may result from defect lines in $\Sigma$, and where each defect point
and segment is in addition equipped with the choice of a normal direction into the disk.
As an illustration, consider the annulus
    \def\locpa  {2.2}    
    \def\locpA  {0.5}    
    \def\locpb  {0.25}   
    \def\locpc  {2.2pt}  
    \def\locpC  {1.8pt}  
    \def\locpd  {0.12}   
    \def\locpD  {0.09}   
  \be
  \raisebox{-5.3em}
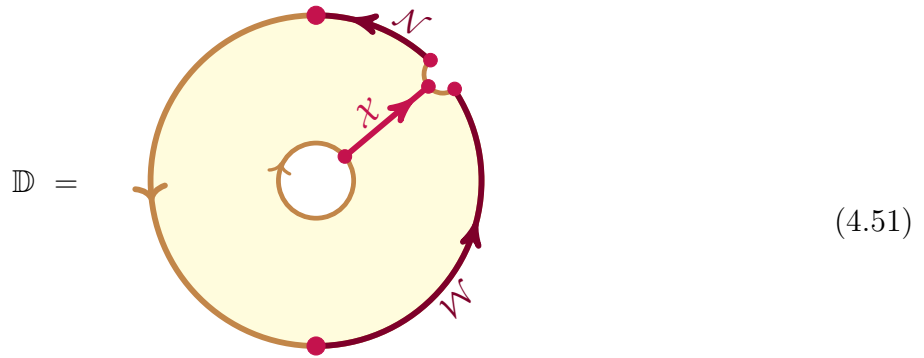
{\begin{tikzpicture}
  \scopeArrow{0.52}{>}
  \filldraw[\colorCircle,fill=\colorSigma,line width=\locpc,postaction={decorate}]
             (0,0) circle (\locpa);
  \end{scope}
  \scopeArrow{0.42}{\arrowDefect} \scopeArrow{0.93}{\arrowDefect}
  \draw[line width=\locpc,color=\colorFreebdy,postaction={decorate}]
	     (-90:\locpa) arc (-90:90:\locpa);
  \end{scope} \end{scope}
  \scopeArrow{0.76}{\arrowDefect}
  \draw[line width=\locpc,color=\colorDefect,postaction={decorate}] (0,0) --
             node[above,midway,sloped,yshift=-2pt,color=\colorDefect] {$\calx$} (40:\locpa);
  \end{scope}
  \filldraw[\colorCircle,fill=white,line width=\locpC] (0,0) circle (\locpA);
  \fill[white] (40:\locpa) circle (\locpb);
  \draw[\colorCircle,line width=\locpC] (40:\locpa)++(120:\locpb) arc (120:310:\locpb);
  \scopeArrow{0.47}{<}
  \draw[\colorCircle,line width=0.7*\locpC,postaction={decorate}] (0,0) circle (\locpA);
  \end{scope}
  \filldraw[color=\colorDefect] (90:\locpa) circle (\locpd) (270:\locpa) circle (\locpd)
             (40:\locpA) circle (\locpD) (40:\locpa-\locpb) circle (\locpD)
             (40:\locpa)++(131:\locpb) circle (\locpD)  
             (40:\locpa)++(-51:\locpb) circle (\locpD);
  \node[rotate=-131,color=\colorFreebdy] at (-40:1.12*\locpa) {$\calm$};
  \node[rotate=-32,color=\colorFreebdy] at (59:1.12*\locpa) {$\caln$};
  \node at (-\locpa-1.4,0) {$\DD ~=$};
  \end{tikzpicture}}
  \ee
  ~\\
In this example the elements of $\Edg$ are, besides the defect points, the two free boundaries 
labeled by $\calm$ and $\caln$, each with a single (namely, inward) normal direction, and twice
the defect line labeled by $\calx$, with two different choices of normal direction.

Each element of $\Edg$ is a possible start and end point of one of the parallel transport
operations on $\DD$. 
Recall that $\int_{a \in \cala} \cc{a} \boti a \iN \cc{\cala} \boti \cala$ is a coalgebra,
with comultiplication induced by the monoidal structure of $\cala$, and with counit 
$\epsilon\colon\! \int_{a \in \cala} \cc{a} \boti a \,{\to}\, \cc{\one} \boti \one$ given by 
the component at $\one$ of the universal dinatural transformation of the end. 
Recall further that forgetting the balancings provides a functor 
$\ZZ(\partial \DD) \,{\to}\, \UU(\partial \DD)$, with $\ZZ(\partial \DD)$ the Deligne product
of the gluing categories for the gluing segments on $\partial\DD$. Now for a disk $\DD$ the
pre-block functor is just the $\Hom$ functor on $\UU(\partial\DD)^{\opp} \boti \UU(\partial\DD)$,
and each $x \iN \Edg$ corresponds to a Deligne factor in 
$\UU(\partial\DD)^{\opp} \boti \UU(\partial\DD)$. 
We define a coaction of $\int_{a \in \cala} \cc{a} \boti a$ on 
$\UU(\partial\DD)^{\opp} \boti \UU(\partial\DD)$ as follows. In $\UU(\partial \DD)^{\opp}$ there 
are De\-lig\-ne products $\calm_{s} \boti \cc{\calm_{s}}$ and $\calm_{t} \boti \cc{\calm_{t}}$
corresponding to $x_{s}$ and  $x_{t}$. There are four cases to consider, depending on 
whether $\calm_{s}$ and $\calm_{t}$ are left or right $\cala$-modules. In each case the
coaction takes place on an object $\cc{m} \boti \cc{n} \iN \cc{\calm_{s}} \boti \cc{\calm_{t}}$,
and depends on two integers $\mu$ and $l$, where $\mu \,{=}\,\sum_{i=1}^{l}\mu_{i}$ is 
the sum of the framing indices counted clockwise from $x_{s}$ to $x_{t}$ and $l$ is the number
of gluing segments along that path. We define the coaction in the four cases as
  \be
  \cc{m} \boti \cc{n} \,\longmapsto \left\{ \bearll
  \int_{a\in\cala} \cc{a.m}  \boti \cc{a^{[\mu-l+1]}.n} & 
  \mbox{if $\calm_{s}$ and $\calm_{t}$ are left modules,}
  \Nxl5
  \int_{a\in\cala} \cc{a.m} \boti \cc{n.a^{[l-\nu]}} &
  \mbox{if $\calm_{s}$ is a left and $\calm_{t}$ a right module,}
  \Nxl5
  \int_{a\in\cala} \cc{m.a} \boti \cc{a^{[l-\mu]} n} &
  \mbox{if $\calm_{s}$ is a right and $\calm_{t}$ a left module,}
  \Nxl5
  \int_{a\in\cala} \cc{m.a} \boti \cc{n.a^{[l-\mu-1]}} &
  \mbox{if $\calm_{s}$ and $\calm_{t}$ are right modules,}
  \label{eq:cases-parallel-comon}
  \eear \right.
  \ee
respectively.

\begin{Proposition} \label{prop:structure-comodule-hol}
Let $\DD$ be a disk and $\cala$ the finite tensor category labeling its interior.
Denote by $\UU(\partial \DD)$ the category that is
obtained by taking the Deligne product over the labels at all defect points of $\DD$.
Let $x_{s}$ and $x_{t}$ be any two elements of the set $\Edg$ that correspond to defect lines. 
\Itemizeiii

\item[{\rm (i)}]
The $\cala$-coactions \eqref{eq:cases-parallel-comon} yield a canonical comonad 
$Z_{\DD,x_{s},x_{t}}$ on $\Funle( \UU(\partial\DD)^{\opp} \boti \UU(\partial\DD),$
                  $\vect)$.

\item[{\rm (ii)}]
The clock- and counterclockwise parallel transport operations from $x_{s}$ to $x_{t}$
provide two structures 
  \be
  \holc_{\DD,x_{s},x_{t}}, \holcc_{\DD,x_{s},x_{t}} :\quad \Zpre(\DD) \xrightarrow{~~} Z_{\DD,x_{s},x_{t}} (\Zpre(\DD))
  \label{eq:hol-comod}
  \ee
of a $Z_{\DD,x_{s},x_{t}}$-comodule on the functor $\Zpre(\DD)$.

\item[{\rm (iii)}]
The fine block functor $\Zfine(\DD)$ is the equalizer of $\holc_{\DD,x_{s},x_{t}}$ and $\holcc_{\DD,x_{s},x_{t}}$. It
depends on the choice of $x_{s},x_{t} \iN \Edg$ only up to a canonical isomorphism. 
\end{itemize}
\end{Proposition}

In case the start and end of the parallel transport are clear from the context, we just write
$\holc_{\DD}$ for $\holc_{\DD,x_{s},x_{t}}$, and analogously 
$\holcc_{\DD} \,{=}\, \holcc_{\DD,x_{s},x_{t}}$.

\begin{proof}
(i)\, By pre-composition (in the same way as in \eqref{eq:inherited-comonad}) 
we obtain from \eqref{eq:cases-parallel-comon} a coaction of the coalgebra
$\int_{a\in\cala} \cc{a} \boti a$ on the functor category 
$\Funle(\UU(\partial\DD)^{\opp} \boti \UU(\partial\DD),\vect)$. This way the coactions 
of $\cala$ on the defects $x_{s}$ and $x_{t}$ provide canonically a comonad
$Z_{\DD,x_{s},x_{t}}$ on $\Funle(\UU(\partial\DD)^{\opp} \boti \UU(\partial\DD), \vect)$.
 \\[.2em]
(ii)\,
We define the morphism $\holc_{\DD,x_{s},x_{t}}$ as follows. The pre-block functor on $\DD$ 
is given by
  \be
  \label{eq:prebl-DD}
  \Zpre(\DD)(-) = \int^{m_{s}\in\calm_s}\!\!\! \int^{m_{t}\in\calm_t}\!
  \cdots\, \Hom(\cc{m_{s}} \boti m_{s} \boti m_{t} \boti \cc{m_{t}} \boti \ldots, -) \,,
  \ee
where the ellipsis accounts for the additional defect lines on $\DD$. Consider the first case
in the list \eqref{eq:cases-parallel-comon}. Then the component 
$\holc_{\DD,x_{s},x_{t}}({}^{\vee\!}a)$ of the parallel transport $\holc_{\DD,x_{s},x_{t}}$ 
at ${}^{\vee\!} a \,{\in}\, \cala$ is defined as the composite
  \be
  \hspace*{-1.0em} \begin{array}{ll}
  \multicolumn2l{ \holc_{\DD,x_{s},x_{t}}({}^{\vee\!}a): \quad \Zpre(\DD)(-) }
  \Nxl5 \qquad & \dsty
  \xrightarrow{~(\evl)^{*}~}
  \int^{m_{s}\in\calm_s}\!\!\!\! \int^{m_{t}\in\calm_t}\!\! \cdots\,
  \Hom(\cc{m_s} \boti (a\,{\otimes}{}^{\vee\!} a).m_s \boti m_t \boti \cc{m_t} \boti \cdots\!,-) ~
  \Nxl2 & \dsty 
  \xrightarrow{\,\gamma_{x_{s},x_{t}}({}^{\vee}a)\,}
  \int^{m_{s}\in\calm_s}\!\!\!\! \int^{m_{t}\in\calm_t}\!\! \cdots\,
  \Hom(\cc{m_{s}} \boti a.m_{s} \boti a^{[\mu-l+1]}.m_{t} \boti \cc{m_{t}} \boti \cdots\!, -) \,.
  \eear
  \ee
Here $\gamma_{x_{s},x_{t}}^{\text{c}}({}^{\vee\!}a)$ are (a slight generalization of)
the parallel transport operations in Lemma \ref{lemma:independent-start}, which are now 
allowed to start and end at the variables corresponding to the defect lines; we consider
them in the clockwise version, i.e.\ the first isomorphism in
$\gamma_{x_{s},x_{t}}^{\text{c}}({}^{\vee\!}a)$ is induced by the isomorphism $\int^{m_{s}}\!
a.m_{s} \boti \cc{m_{s}} \,{\cong}\, \int^{m_{s}}\! a.m_{s} \boti \cc{a^{\vee}.m_{s}}$.
Then $\holc_{\DD,x_{s},x_{t}}$ is defined by taking the end over all 
$\holc_{\DD,x_{s},x_{t}}({}^{\vee\!}a)$. $\holcc_{\DD,x_{s},x_{t}}$ is defined analogously
by using the counterclockwise parallel transport operations instead.
 \\
Together with Proposition \ref{prop:wellposedD} these prescriptions imply that, in all 
four cases, both $\holc_{\DD}$ and $ \holcc_{\DD}$ are indeed natural transformations
from $\Zpre(\DD)$ to $Z_{\DD,x_{s},x_{t}}(\Zpre(\DD))$.
 \\[.2em]
To show that $\holc_{\DD}$ is a comodule structure, we first note that composing it with the
counit $\epsilon$ of $Z_{\DD,x_{s},x_{t}}$ yields the identity on $\Zpre(\DD)$, since the 
parallel transport of $\one$ on $\Zpre(\DD)$ is the identity. Next consider, for $a,b \iN \cala$,
the component $\holc_{\DD, a \otimes b}$ of the parallel transport from $x_{s}$ to $x_{t}$. In
each step involved in the parallel transport  operation we use either adjunction morphisms 
or a balancing, and both of these are compatible with the monoidal structure of $\cala$, 
i.e.\ can be split into the
composite of corresponding steps for $\holc_{\DD,a}$ and $\holc_{\DD,b}$. In both cases the
morphisms commute (up to changing their arguments accordingly), and after passing to the end
we obtain that indeed $\holc_{\DD}$ is a comodule structure for $\Zpre(\DD)$. The proof for 
$\holcc_{\DD}$ is analogous. 
 \\[.2em]
(iii)\, The statement follows from the use of the parallel transport operations, analogously
as in the proof of Lemma \ref{lemma:independent-start}.
\end{proof}

Analogous statements as in Proposition \ref{prop:structure-comodule-hol} still hold in the
situation that one or both of $x_{s}$ and $x_{t}$ in the Lemma is a general element of the
set $\Edg$: There is a comonad $Z_{\DD,x_{s},x_{t}}$, defined again case by case, 
on $\Funle(\UU(\partial\DD)^{\opp} \boti \UU(\partial\DD), \vect)$ such that the 
parallel transport  operations provide comodule structures for $\Zpre(\DD)$ over
$Z_{\DD,x_{s},x_{t}}$ and the block functor is the corresponding equalizer. 

As in Section \ref{sec:pre-blocks} we denote the Deligne product over labels for all 
boundary segments of a defect surface $\Sigma$ by $\UU(\partial \Sigma)$. The comonad
$Z_{\DD,x_{s},x_{t}}$ on $\Funle(\UU(\partial\DD)^{\opp} \boti \UU(\partial\DD), \vect)$
induces a corresponding comonad on the functor category
$\Funle(\UU(\partial\Sigma)^{\opp} \boti \UU(\partial\Sigma), \vect)$, which for simplicity
we denote again by $Z_{\DD,x_{s},x_{t}}$. 
Clearly, $\Zpre(\Sigma)$ becomes a comodule over the latter comonad $Z_{\DD,x_{s},x_{t}}$.
We call the comonads  $Z_{\DD,x_{s},x_{t}}$ the \emph{parallel transport comonads}.


\subsection{Blocks as module natural transformations}

In this subsection we express the block spaces for the prototypical example of the \basicdisk\
\eqref{pic:ex:4.4} considered in Example \ref{ex:4.4-part1} as spaces of module natural 
transformations. Recall that the pre-block functor in this example is given by formula
\eqref{eq:preblock441} or, equivalently, by the space \eqref{eq:preblock442} of
natural transformations, and that the holonomy is described in \eqref{eq:iso1}.

According to Lemma \ref{lemma:gen-equalizer}, for comodules $w,x \,{\in}\, \calm^T$ 
over a comonad $T$ on $\calm$, with coactions $\delta_w$ and $\delta_x$, the vector space
$\Hom_{\calm^T}(w,x)$ of comodule morphisms can be described as the equalizer of the two maps
  \be
  \psi_1,\psi_2:\quad \Hom_{\calm}(U(w),U(x)) \xrightarrow{~~} \Hom_{\calm}(U(w),T(U(x))) 
  \label{eq:psi1psi2}
  \ee
that are given by $\psi_1(f) \,{=}\, \delta_x\,{\circ}\,f$ and
$\psi_2(f) \,{=}\, T(f)\,{\circ}\,\delta_w$, respectively.
The properties of an adjunction readily imply

\begin{Lemma}
\label{lemma:gen-lemma-como-ad}
If the comonad $T$ on $\calm$ is left exact, then the map $\psi_2$ in
\eqref{eq:psi1psi2} equals the composite 
  \be
  \Hom_{\calm}(Uw,Ux) \xrightarrow{~(\delta_w^{\la})^*_{}~} \Hom_{\calm}(T^{\la}(Uw),Ux)
  \xrightarrow{~\cong~}  \Hom_{\calm}(Uw,T(Ux)) 
  \ee
of the pre-composition with the image $\delta_w^{\la}$ of $\delta_w$
under the adjunction $\Hom_\calm(Uw,T(Uw))$
          ${\cong}\, \Hom_\calm(T^{\la}(Uw),Uw)$ and the map provided by the adjunction.
\end{Lemma}

In the situation considered in Example \ref{ex:4.4-part1} we deal with the two 
categories \eqref{eq:TL1-TL2}, i.e.\
$\ZZ(\LL_1) \,{=}\, \cenf{\calmopp \,{\stackrel{1}\boxtimes \caln}\,}$ and
$\ZZ(\LL_2) \,{=}\, \cenf{\calnopp \,{\stackrel{-1}\boxtimes}\, \calm}
\,{\simeq}\, \ZZ(\LL_1)\Opp$.
By Lemma \ref{lemma:tw-center-cat-comod}, $\ZZ(\LL_1)$ is equivalent to the category of 
comodules over the comonad $Z_{[1]}$ on $\calmopp \boti \caln$, while $\ZZ(\LL_2)$ 
corresponds to the comodules over the comonad $Z_{[-1]}$ on $\calnopp \boti \calm$. Now 
note that $(\calnopp \boti \calm) \Opp \,{\simeq }\, \calmopp \boti \caln$.
By direct computation, using the dualities of the finite tensor category $\cala$
and the fact that taking the opposite interchanges coend and end, we then have 

\begin{Lemma}
\label{lemma:Ex-adj-comon}
The comonads on $Z_{[1]}$ on $\calmopp \boti \caln$ and $Z_{[-1]}$ on $\calnopp \boti \calm$
are related as 
  \be
  Z_{[-1]} \,\cong\, \cc{Z_{[1]}^{\la}} \,,
  \ee
where the overline indicates the opposite functor $($compare Eq.\ \eqref{eq:oppositefunctor}$)$.
\end{Lemma}

To highlight the power of our modular functor for producing higher algebra,
we now describe the block functor for the situation of Example \ref{ex:4.4}.

\begin{Proposition} \label{prop:modulenattrafo}
The block functor for the \basicdisk\ $\DDsAMN$ \eqref{pic:ex:4.4} is isomorphic
to the Hom functor for the category of comodules over the comonad $Z_{[1]}$, i.e.\ we have
  \be
  \Zfine(\DDsAMN)(x \boti y) \,\cong\, \Hom_{Z_{[1]}}(y,x)
  \label{eq:block4ex:4.4}
  \ee
for any pair of objects $x \,{\in}\, \ZZ(\LL_1)$ and $y \,{\in}\, \ZZ(\LL_2)$.
Hence the block space can be seen as a space of module natural transformations. 
\end{Proposition}
 
\begin{proof}
(i)~\,%
The main idea is to invoke Lemma \ref{lemma:gen-equalizer} which characterizes the
morphisms of comodules over a comonad as an equalizer. To this end we have to match the
two morphisms \eqref{eq:coexSequ} in that Lemma with the two parallel arrows in the  
Definition \ref{Definition:block-func} of the block functor. 
In the situation at hand, the block space is the equalizer
  \be
  \Zfine(x \boti y) \xrightarrow{~~~} \Zpre(x \boti y)
  \xrightrightarrows[~({\mu_x)}_*^{}] {~\hol_{x}} \, \Zpre( Z_{[1]}(x) \boti y) \,.
  \label{eq:blok-nat-trans}
  \ee
According to \eqref{eq:preblock441} the pre-blocks of our interest are given by
$\Zpre(x \boti y) \,{\cong}\, \Hom_{\calmopp \boxtimes \caln}(y,x)$. 
In this description the map $(\mu_x)_*^{}$ amounts to the map
  \be
  \begin{array}{c}
  ~ \Hom(y,x)\xrightarrow{~~\,} \Hom(y,Z_{[1]}(x)) \,,
  \Nxl3
  f \xmapsto{~~} \mu_{x} \circ f.
  \eear
  \ee
According to Lemma \ref{lemma:Ex-adj-comon}, the comodule structure $\delta_y$ on $y$ gives 
the module structure $Z_{[1]}^{\la} y \,{\to}\, y$ as the image $\delta_y^{\la}$ of $\delta_y$
under the adjunction $\Hom(y,Z_{[-1]}(y)) \,{\cong}\, \Hom(Z_{[-1]}^{\la}(y),y)$.  
\\[.2em]
(ii)~\,%
To obtain \eqref{eq:block4ex:4.4}, by Lemma \ref{lemma:gen-lemma-como-ad} we are thus left 
with showing that the morphism $\hol_x$ corresponds to the composite
  \be
  \Hom(y,x) \xrightarrow{\,(\delta_{y}^{\la})^*_{}} \Hom(Z^{\la}_{[1]}(y),x) 
  \xrightarrow{~\cong~} \Hom(y,Z_{[1]}(x)) \,,
  \label{eq:reform-hol}
  \ee
where the first map is pre-composition with $\delta_y^{\la}$ and the second one is provided 
by the adjunction. To show that this is the case, we consider for $a \,{\in}\, \cala$ the 
following two maps. First, the map
  \be
  \begin{array}{r}\dsty
  g_a:\quad \Zpre(\cc{x_m} \boti x_n \boti  \cc{y_{n}} \boti y_{m})
  \xrightarrow{~(\evr)^{*}~}
  \Zpre( \cc{x_m} \boti x_n \boti  \cc{(a^{\vee} \oti a) \,.\,y_{n}} \boti y_{m}) \quad
  \Nxl5 \dsty
  \xrightarrow{\,~\cong~\,}
  \Zpre(\cc{x_m} \boti x_n \boti  \cc{a^{\vee} .y_{n}} \boti a^{\vee}.y_{m}) \,,
  \eear
  \label{eq:mu-y}
  \ee
composed from an evaluation in the contravariant argument of the pre-block functor
and the balancing of $y$. And second, the isomorphism
  \be
  h_a:\quad \Zpre(\cc{x_{m}}\boti x_{n} \boti \cc{a^{\vee}.y_{n}} \boti a^{\vee}. y_{m}) 
  \xrightarrow{~\cong~}
  \Zpre(\cc{a^{\vee\vee}\!.\, x_m} \boti a\,.\,x_n \boti \cc{y_{n}} \boti y_{m})
  \ee
that is analogous to the composition of the second, third and fourth maps in
the chain \eqref{eq:compositeHol} of isomorphisms.
Together with the component $(\hol_{x})_a$ of the holonomy we thus have a diagram
  \be
  \hspace*{-0.3em}
  \begin{tikzcd}[column sep=-0.7cm,row sep=1.8em]
  \Zpre(\cc{x_m}\boti x_n \boti\cc{y_n} \boti y_m) \ar{rr}{(\hol_x)_a} 
  \ar{rd}{g_a}
  & ~ & \Zpre(\cc{a^{\vee\vee}\!. x_m} \boti a.x_n \boti \cc{y_n} \boti y_m)
  \\
  ~& \Zpre(\cc{x_m}\boti x_n \boti \cc{a^{\vee}\!.y_n} \boti a^{\vee}\!. y_m) \ar{ru}{h_a}
  \end{tikzcd}
  \label{eq:comm-diag-hol}
  \ee
By a straightforward, albeit lengthy, computation using the definition of the holonomy
$\hol_{x}$, the dualities of $\cala$ and the balancings \eqref{eq:mmbalancing} of the
object $\int_{m}\cc{m} \boti m $, it can be seen that this diagram commutes.
 \\
Next we take the end over $a \,{\in}\, \cala$ in the diagram \eqref{eq:comm-diag-hol}.
Then in the top row we get the holonomy, while
  \be
  {\dsty\int_a g_a}:\quad \Zpre(x \boti y) \xrightarrow{~~} 
  \Zpre(x \boti \!\int_{a} \cc{a^{\vee}\!. y_n} \boti a^{\vee}\!.  y_m)
  \,\cong\, \Zpre(x \boti \!\int_{a} \cc{a. y_n} \boti a. y_m) \,.
  \ee
Moreover, by the adjunction obtained in Lemma \ref{lemma:Ex-adj-comon} we have
  \be
  \Zpre(x \boti \!\int_{a} \cc{a. y_n} \boti a. y_m)
  =\, \Zpre(x \boti Z_{[-1]}(y)) \xrightarrow{~\cong~} \Hom(Z_{[1]}^{\la}(y),x) \,.
  \ee
It follows that we have a commuting diagram
  \be
  \begin{tikzcd}[row sep=2.2em]
  \Hom(Z_{[1]}^{\la}(y),x) \ar{r}{\cong} \ar{d}{\cong} & \Hom(y,Z_{[1]}(x)) \ar{d}{\cong}
  \\
  \int_{a \in \cala}\Zpre(\cc{x_{m}}\boti x_{n} \boti \cc{a^{\vee}.y_{n}} \boti a^{\vee}. y_{m})
  \ar{r}{\cong} &
  \int_{a \in \cala}\Zpre(\cc{a^{\vee\vee}\!.\, x_m} \boti a\,.\,x_n \boti \cc{y_{n}} \boti y_{m})
  \end{tikzcd}
  \label{eq:comm-diag-hol-adjun}
  \ee
This concludes the proof of the first statement of the Proposition. 
\\[.2em]
(iii)~\,%
Due to Corollary \ref{cor:lex-framed-app}, 
there is an equivalence $ \ZZ(\LL_1) \,{\simeq}\, \Funle_{\cala}(\calm,\caln)$.
\end{proof}

Combining the result \eqref{eq:block4ex:4.4} with
the equivalence $ \ZZ(\LL_1) \,{\simeq}\, \Funle_{\cala}(\calm,\caln)$, we get
(similarly as for the formula \eqref{eq:preblock442} for the pre-block functor):

\begin{cor} \label{cor:modulenattrafo}  
The block spaces for the \basicdisk\ $\DDsAMN$ \eqref{pic:ex:4.4} are given by the spaces
  \be
  \Zfine(\DDsAMN)(x \boti y) \cong \Nat_\cala(\Phile(\cc y),\Phile(x)) \,.
  \label{eq:block442}
  \ee
of module natural transformations between the module functors $\Phile(\cc y)$ and $\Phile(x)$. 
\end{cor}

\begin{Example} \label{exa:generaldisk2}
Consider a disk with any number $N$ of gluing and free boundary segments, as studied in 
Example \ref{exa:generaldisk}, imposing in addition the sum rule 
$\sum_i\! \kappa_i \,{=}\, N{-}2$ from Lemma \ref{lem:4.5} on the framing indices. Then
Proposition \ref{prop:modulenattrafo} implies that the block functor is given by
  \be
  \Zfine(x_1 \boti x_2 \boti\cdots\boti x_n) \cong \Nat_\cala\big(
  \Phile(\cc{x_1}),\Phile(x_n)\,{\circ}\,\Phile(x_n)\,{\circ}\cdots{\circ}\,\Phile(x_2) \big) \,,
  \ee
where the module structure on the module functors $\Phile(x_i)$ is determined by the
framing indices of the adjacent gluing intervals.
\end{Example}
           
In particular, the result of Proposition \ref{prop:modulenattrafo} for the block functor 
of the \basicdisk\ generalizes to the case of disks of the form \eqref{eq:disk.k.-k}. This, in 
turn, can be used to give the block functor for any cylinder that is built out of such disks,
i.e.\ for all \defectsurface s of the form
    \def\locpa  {0.88}   
    \def\locpA  {2.63}   
    \def\locpb  {0.59}   
    \def\locpB  {2.93}   
    \def\locpc  {2.2pt}  
    \def\locpd  {0.12}   
    \def\locpD  {0.04}   
    \def\locpe  {1.09}   
  \be
  \raisebox{-7.1em}{\begin{tikzpicture}[scale=1.1]
  \node[color=black] at (-\locpA-1.8,0) {$\Sigma ~=$};
  \scopeArrow{0.02}{>} \scopeArrow{0.15}{>} \scopeArrow{0.61}{>} \scopeArrow{0.74}{>}
  \scopeArrow{0.88}{>}
  \filldraw[\colorCircle,fill=\colorSigma,line width=\locpc,postaction={decorate}]
                       (0,0) circle (\locpA);
  \end{scope} \end{scope} \end{scope} \end{scope} \end{scope}
  \draw[dashed,line width=\locpc,color=white] (106:\locpA) arc (106:164:\locpA);
  \scopeArrow{0.02}{<} \scopeArrow{0.18}{<} \scopeArrow{0.62}{<} \scopeArrow{0.76}{<}
  \scopeArrow{0.90}{<}
  \filldraw[\colorCircle,fill=white,line width=\locpc,postaction={decorate}] (0,0) circle (\locpa);
  \end{scope} \end{scope} \end{scope} \end{scope} \end{scope}
  \draw[dashed,line width=\locpc,color=white] (109:\locpa) arc (109:161:\locpa);
  \scopeArrow{0.61}{\arrowDefect}
  \draw[line width=\locpc,color=\colorDefect,postaction={decorate}] (80:\locpa) -- (80:\locpA);
  \draw[line width=\locpc,color=\colorDefect,postaction={decorate}] (-20:\locpa) -- (-20:\locpA);
  \draw[line width=\locpc,color=\colorDefect,postaction={decorate}] (-120:\locpa) -- (-120:\locpA);
  \end{scope}
  \scopeArrow{0.61}{\arrowDefectB}
  \draw[line width=\locpc,color=\colorDefect,postaction={decorate}] (30:\locpa) -- (30:\locpA);
  \draw[line width=\locpc,color=\colorDefect,postaction={decorate}] (-70:\locpa) -- (-70:\locpA);
  \draw[line width=\locpc,color=\colorDefect,postaction={decorate}] (-170:\locpa) -- (-170:\locpA);
  \end{scope}
  \filldraw[color=\colorDefect] (80:\locpa) circle (\locpd) (80:\locpA) circle (\locpd)
                                (30:\locpa) circle (\locpd) (30:\locpA) circle (\locpd)
                               (-20:\locpa) circle (\locpd) (-20:\locpA) circle (\locpd)
                               (-70:\locpa) circle (\locpd) (-70:\locpA) circle (\locpd)
                              (-120:\locpa) circle (\locpd) (-120:\locpA) circle (\locpd)
                              (-170:\locpa) circle (\locpd) (-170:\locpA) circle (\locpd);
  \node at (45:1.03*\locpe) {$\kappa_1$}; \node at (55:\locpB) {${-}\kappa_1$};
  \node at (13:1.05*\locpe) {$\kappa_2$}; \node at (7:1.05*\locpB) {${-}\kappa_2$};
  \node at (-50:1.04*\locpe) {$\kappa_3$}; \node at (-38:1.03*\locpB) {${-}\kappa_3$};
  \node at (-103:\locpe) {$\kappa_4$}; \node at (-88:0.97*\locpB) {${-}\kappa_4$};
  \node at (-155:1.03*\locpe) {$\kappa_5$}; \node at (-141:1.03*\locpB) {${-}\kappa_5$};
  \node[color=\colorDefect] at (70:0.83*\locpA) {$\calm_1$};
  \node[color=\colorDefect] at (22:0.83*\locpA) {$\calm_2$};
  \node[color=\colorDefect] at (-29:0.75*\locpA) {$\calm_3$};
  \node[color=\colorDefect] at (-58:0.80*\locpA) {$\calm_4$};
  \node[color=\colorDefect] at (-110:0.74*\locpA) {$\calm_5$};
  \node[color=\colorDefect] at (-180:0.73*\locpA) {$\calm_6$};
  \node at (52:0.67*\locpA) {$\cala_1$}; \node at (2:0.70*\locpA) {$\cala_2$};
  \node at (-50:0.62*\locpA) {$\cala_3$}; \node at (-87:0.77*\locpA) {$\cala_4$};
  \node at (-145:0.69*\locpA) {$\cala_5$};
  \filldraw (115:0.5*\locpa+0.5*\locpA) circle (\locpD)
            (125:0.5*\locpa+0.5*\locpA) circle (\locpD)
            (135:0.5*\locpa+0.5*\locpA) circle (\locpD)
            (145:0.5*\locpa+0.5*\locpA) circle (\locpD)
            (155:0.5*\locpa+0.5*\locpA) circle (\locpD);
  \end{tikzpicture}
  \label{eq:refinedcyl}}
  \ee
We refer to surfaces of the form \eqref{eq:refinedcyl}, as well as to their counterparts
based on a gluing interval instead of a gluing circle, as \emph{defect cylinders}.

\begin{cor} \label{cor:gen-cylinder}
Let $\Sigma$ be defect cylinder as shown in picture \eqref{eq:refinedcyl}.
Denote the inner boundary circle of $\Sigma$ by $\SS_1$ and the outer one by $\SS_2$.
We have $\ZZ(\SS_1) \,{\cong}\, \ZZ(\SS_2)\opp_{}$. 
\\
Regarding $\Sigma$ as a bordism from
$\SS_1 \,{\sqcup}\, \SS_2$ to $\emptyset$, the block functor for $\Sigma$ takes the values 
  \be	
  \Zfine(\Sigma)(F\boti G) =\Hom_{\ZZ(\SS_2)}(F,G)
  \label{eq:ZZ-cylinder}	
  \ee
for $F \,{\in}\, \ZZ(\SS_1)$ and $G\,{\in}\, \ZZ(\SS_2)$. Regarding the surface
$\Sigma$ instead as a bordism from $\ZZ(\SS_1)$ to $\ZZ(\SS_2)\opp_{} \,{\cong}\, \ZZ(\SS_1)$,	
$\Zfine(\Sigma)$ is just the identity functor on $\ZZ(\SS_1)$.
\end{cor}

\begin{Example} \label{ex:3holedsphere}
For finite tensor categories $\cala$ and $\calb$, an $\cala$-$\calb$-bimodule $\calm_1$ and 
any finite number of $\cala$-$\calb$-bi\-mo\-dules $\calm_i$, $i \,{=}\, 2,3,...\,,n$,
consider a sphere with one circular defect line that is interrupted by gluing circles 
(with arbitrary admissible framing indices, which we suppress in the picture)
which cut it into intervals labeled by the bimodules $\calm_i$:
    \def\locpa  {3.0}    
    \def\locpA  {0.39}   
    \def\locpb  {1.25}   
    \def\locpC  {1.8pt}  
    \def\locpD  {0.08}   
    \def\locpf  {1.1pt}  
  \be
  \raisebox{-7.4em}
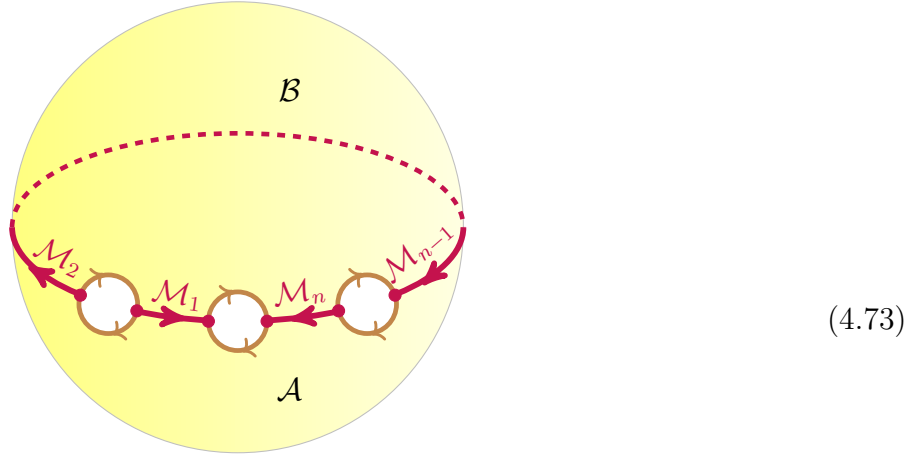
{\begin{tikzpicture}[scale=1.1]
  \coordinate (incirc1) at (0,-\locpb);
  \coordinate (incirc2) at ({\locpa*cos(235)},{\locpb*sin(235)}); 
  \coordinate (incirc3) at ({\locpa*cos(305)},{\locpb*sin(305)}); 
  \shade[left color=yellow!50,right color=yellow!11] (0,0) circle (\locpa);
  \draw[color=gray!50] (0,0) circle (\locpa);
  \scopeArrow{0.13}{\arrowDefectB} \scopeArrow{0.40}{\arrowDefect} \scopeArrow{0.65}{\arrowDefectB}
  \scopeArrow{0.91}{\arrowDefectB}
  \draw[line width=\locpc,color=\colorDefect,postaction={decorate}] 
           (-\locpa,0) arc (180:360:\locpa cm and \locpb cm);
  \end{scope} \end{scope} \end{scope} \end{scope}
  \draw[line width=0.8*\locpc,color=\colorDefect,dashed] 
           (-\locpa,0) arc (180:0:\locpa cm and \locpb cm);
  \filldraw[\colorCircle,fill=white,line width=\locpC]
           (incirc1) circle (\locpA) (incirc2) circle (\locpA) (incirc3) circle (\locpA);
  \scopeArrow{0.33}{<} \scopeArrow{0.83}{<}
  \draw[\colorCircle,line width=\locpf,postaction={decorate}]
	   (incirc1) node[color=\colorDefect,left=9pt,yshift=8pt,rotate=-8] {$\calm_1$}
	   node[color=\colorDefect,right=9pt,yshift=8pt,rotate=10] {$\calm_n$} circle (\locpA);
  \draw[\colorCircle,line width=\locpf,postaction={decorate}] (incirc2)
	   node[color=\colorDefect,left=6pt,yshift=12pt,rotate=-25] {$\calm_2$} circle (\locpA);
  \draw[\colorCircle,line width=\locpf,postaction={decorate}] (incirc3)
	   node[color=\colorDefect,right=4pt,yshift=7pt,rotate=41] {$\calm_{n{-}1}$} circle (\locpA);
  \end{scope} \end{scope}
  \node at (0.23*\locpa,0.61*\locpa) {$\calb$};
  \node at (0.23*\locpa,-0.72*\locpa) {$\cala$};
  \filldraw[color=\colorDefect]
           (incirc1)+(0:\locpA) circle (\locpD) (incirc1)+(180:\locpA) circle (\locpD) 
           (incirc2)+(-13:\locpA) circle (\locpD) (incirc2)+(162:\locpA) circle (\locpD) 
           (incirc3)+(193:\locpA) circle (\locpD) (incirc3)+(18:\locpA) circle (\locpD);
  \end{tikzpicture}}
  \ee
  ~\\[0.2em]
This can be recognized as the \defectsurface\ $\Sigma$ that is obtained when gluing a
disk of the type considered in Examples \ref{exa:generaldisk} and \ref{exa:generaldisk2} 
along its boundary to an oppositely oriented disk of the same type. Accordingly, 
the pre-block spaces -- which only depend on the incidence relations of $\vSigma \bs1$ 
(i.e.\ defect lines and the boundary of $\Sigma$) --
are given by the same expression \eqref{eq:exa:generaldisk-Zpre} as the pre-block spaces 
for one of those disks, albeit with their arguments now being objects in different 
categories, which are assigned to gluing circles rather than gluing intervals. Imposing
flat holonomy for the upper and lower disk, respectively, then amounts, in the formulation
with left exact functors, to requiring that the so obtained natural transformations are
module natural transformations with respect to the right $\cala$- and left $\calb$-action,
respectively. Imposing both of these (where the order does not matter) thus yields
$\cala$-$\calb$-\emph{bimodule} natural transformations:
  \be
  \Zfine(x_1 \boti x_2 \boti\cdots\boti x_n) \cong \Nat_{\cala,\calb}\big(\Phile(\cc{x_1}),
  \Phile(x_n)\,{\circ}\,\Phile(x_{n-1})\,{\circ}\cdots{\circ}\,\Phile(x_2) \big) \,.
  \label{eq:exa:3holedsphere}
  \ee
Of particular interest is the special case that all defects are transparent, i.e.\ that
$\calb \,{=}\, \cala$ and $\calm_i \,{=}\, \cala$ for each $i$, and that the framing indices
on the two segments of all gluing circles, except for the circle between the $\calm_1$- and
$\calm_{n}$-line, are equal to 1. Then the gluing categories for all other circles are given by
  \be
  \cala \,\Boti 1\, \CC{\cala} \,\Boti 1 \stackrel{\eqref{eq:Drinf-framed-cent-A}}\simeq\!
  \Funle_{\cala,\cala}(\cala,\cala) \simeq \calz(\cala) \,,
  \label{eq:Funle_AA(AA)=ZA}
  \ee
and the composition of functors in \eqref{eq:exa:3holedsphere} amounts to the tensor product in
the Drinfeld center $\calz(\cala)$. 
\end{Example}


\subsection{Contraction along a defect path}\label{subsection:Contrac}

Recall from Proposition \ref{prop:block-compose} the result about the fusion of 
boundary insertions, i.e.\ that there is a canonical
isomorphism between the block functors for the configurations shown in the picture
\eqref{eq:PsiG.PsiF=PsiGF}. These two configurations are related by contracting a 
single defect line between two gluing boundaries. In this subsection we obtain a vast
generalization of this procedure; we provide the behavior of the block functor
under ``contraction along a path of defect lines''. 

By a \emph{path of defect lines} in a \defectsurface\ $\Sigma$ we mean an ordered collection
  \be
  \gamma = (\delta_1,\delta_2,...\,,\delta_n)
  \ee
of defect lines such that each pair
$(\delta_i,\delta_{i+1})$ shares a common gluing boundary; we denote the latter by $\LL_{i,i+1}$. 
Such a path can be open or closed. In the following we first assume for simplicity that all the 
gluing boundaries $\LL_{i,i+1} \,{=}\, \SS_{i,i+1}$ are gluing \emph{circles}; the situation that 
also gluing intervals
are present will be briefly discussed afterwards. Given a path $\gamma$ of defect lines, we
consider a fattening of $\gamma$ to a tubular neighborhood $\nn\gamma \,{\subset}\, \Sigma$
which is sufficiently large such that its interior contains all gluing boundaries $\SS_{i,i+1}$ 
on the path and sufficiently small such that it does not meet any other gluing boundaries
of $\Sigma$. Paths $\gamma$ of defect lines come in several types, which are distinguished
by the form of the boundary $\partial \nn{\gamma}$ of their tubular neighborhood:
\Itemizeiii
\item[(i)] 
$\partial \nn{\gamma}$ consists of a single embedded circle in $\Sigma$.
\item[(ii)] 
$\partial \nn{\gamma}$ is the disjoint union of of two embedded circles
$\partial \nn{\gamma,1}$ and $\partial \nn{\gamma,2}$, both of which meet at least 
one defect line of $\Sigma$.
\item[(iii)] 
$\partial \nn{\gamma}$ is the disjoint union of of two embedded circles 
$\partial \nn{\gamma,1}$ and $\partial \nn{\gamma,2}$, and precisely one of them  
does not meet any defect line of $\Sigma$.
\item[(iv)] 
$\partial \nn{\gamma}$ is the disjoint union of of two embedded circles
$\partial \nn{\gamma,1}$ and $\partial \nn{\gamma,2}$, none of which meets 
any defect line of $\Sigma$.
\end{itemize}
If $\Sigma$ is fine, then in the cases (iii) and (iv) $\gamma$ is a loop that encloses a disk 
$\DD$ in $\Sigma$. The following pictures give examples for a path of type (i) 
and a path of type (ii), each consisting of two defect lines:
 \\[-0.5em]
    \def\locpa  {2.5}    
    \def\locpA  {0.34}   
    \def\locpb  {1.5}    
    \def\locpC  {1.8pt}  
    \def\locpD  {0.08}   
    \def\locpf  {1.1pt}  
  \be
  \raisebox{-7.4em}
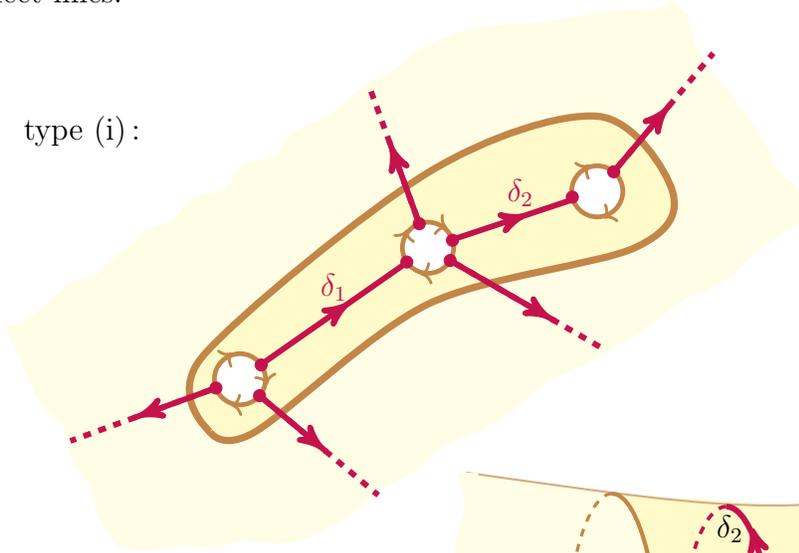
{\begin{tikzpicture}
  \coordinate (incirc1) at (0,0);
  \coordinate (incirc2) at (\locpa,0.7*\locpa);
  \coordinate (incirc3) at (1.9*\locpa,\locpa);
  \coordinate (oucirc0) at (0.2,-0.7);
  \coordinate (oucirc1) at (-0.4,-0.7);
  \coordinate (Oucirc1) at (-0.5,0.3);
  \coordinate (oucirc2) at (\locpa,1.1*\locpa);
  \coordinate (Oucirc2) at (\locpa,0.4*\locpa);
  \coordinate (oucirc3) at (1.9*\locpa,1.4*\locpa);
  \coordinate (Oucirc3) at (2.3*\locpa,\locpa);
  \coordinate (OUcirc3) at (2.1*\locpa,0.7*\locpa);
  \begin{scope}[decoration={random steps,segment length=4mm},rotate=32]
  \fill[\colorSigmaLight,decorate]
          (-1.7*\locpb,-0.7*\locpb) -- (-0.3*\locpb,-1.5*\locpb)
          -- (2*\locpa+1.7*\locpb,-1.5*\locpb) -- (2*\locpa+1.7*\locpb,0.8*\locpb)
	  node[right=2.5cm] {~}  
          -- (1.3*\locpa+\locpb,1.6*\locpb) -- (-1.5*\locpb,1.5*\locpb) -- cycle;
  \end{scope}
  \filldraw[line width=1.5*\locpC,color=\colorCircle,fill=\colorSigmaDark]
          plot[smooth,tension=0.55] coordinates{ (oucirc1) (Oucirc1)
          (oucirc2) (oucirc3) (Oucirc3) (OUcirc3) (Oucirc2) (oucirc0) (oucirc1) };
  \scopeArrow{0.57}{\arrowDefect}
  \draw[line width=\locpc,color=\colorDefect,postaction={decorate}]
          (incirc1) -- node[above] {$\delta_1$} (incirc2);
  \draw[line width=\locpc,color=\colorDefect,postaction={decorate}]
          (incirc2) -- node[above,xshift=3pt] {$\delta_2$} (incirc3);
  \end{scope}
  \scopeArrow{0.95}{\arrowDefect}
  \draw[line width=\locpc,color=\colorDefect,postaction={decorate}] (incirc1) -- ++(200:\locpb);
  \draw[line width=\locpc,color=\colorDefect,postaction={decorate}] (incirc1) -- ++(320:\locpb);
  \draw[line width=\locpc,color=\colorDefect,postaction={decorate}] (incirc2) -- ++(110:\locpb);
  \draw[line width=\locpc,color=\colorDefect,postaction={decorate}] (incirc2) -- ++(330:1.3*\locpb);
  \draw[line width=\locpc,color=\colorDefect,postaction={decorate}] (incirc3) -- ++(50:\locpb);
  \end{scope}
  \draw[line width=\locpc,color=\colorDefect,dashed] (incirc1) -- ++(200:1.6*\locpb)
           (incirc1) -- ++(320:1.6*\locpb) (incirc3) -- ++(50:1.6*\locpb)
           (incirc2) -- ++(110:1.5*\locpb) (incirc2) -- ++(330:1.8*\locpb);
  \filldraw[\colorCircle,fill=white,line width=\locpC]
           (incirc1) circle (\locpA) (incirc2) circle (\locpA) (incirc3) circle (\locpA);
  \scopeArrow{0.03}{<} \scopeArrow{0.36}{<} \scopeArrow{0.76}{<}
  \draw[\colorCircle,line width=\locpf,postaction={decorate}] (incirc1) circle (\locpA);
  \end{scope} \end{scope} \end{scope}
  \scopeArrow{0.19}{<} \scopeArrow{0.45}{<} \scopeArrow{0.77}{<}
  \draw[\colorCircle,line width=\locpf,postaction={decorate}] (incirc2) circle (\locpA);
  \end{scope} \end{scope} \end{scope}
  \scopeArrow{0.36}{<} \scopeArrow{0.85}{<}
  \draw[\colorCircle,line width=\locpf,postaction={decorate}] (incirc3) circle (\locpA);
  \end{scope} \end{scope}
  \filldraw[color=\colorDefect] (incirc1)+(200:\locpA) circle (\locpD) 
           (incirc1)+(320:\locpA) circle (\locpD) (incirc1)+(34:\locpA) circle (\locpD) 
           (incirc2)+(110:\locpA) circle (\locpD) (incirc2)+(330:\locpA) circle (\locpD)
           (incirc2)+(214:\locpA) circle (\locpD) (incirc2)+(17:\locpA) circle (\locpD) 
           (incirc3)+(193:\locpA) circle (\locpD) (incirc3)+(50:\locpA) circle (\locpD);
  \node at (-1.4*\locpb,2.2*\locpb) {type (i)\,:};
  \end{tikzpicture}}
  \ee
    \def\locpa  {2.5}    
    \def\locpA  {0.34}   
    \def\locpb  {4.1}    
    \def\locpB  {3.7}    
    \def\locpc  {2.2pt}  
    \def\locpC  {0.8pt}  
    \def\locpd  {0.09}   
    \def\locpD  {1.4pt}  
    \def\locpe  {0.14}   
    \def\locpu  {1.07}   
    \def\locpx  {0.66}   
    \def\locpX  {2.07}   
    \def\locpy  {2.45}   
    \def\locpY  {2.32}   
    \def\locpz  {2.9}    
  \be
  \hspace*{-2.5em}
  \raisebox{-7.1em}
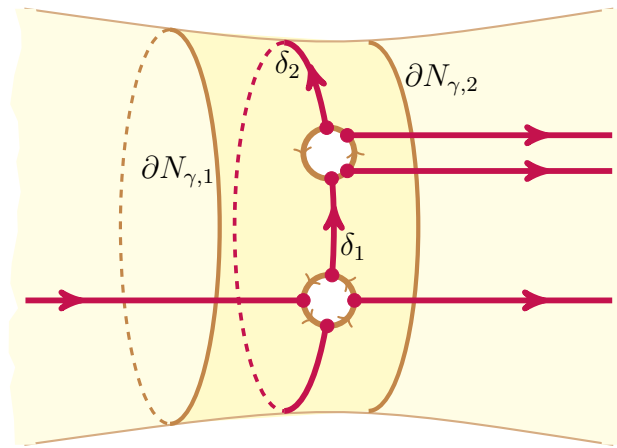
{\begin{tikzpicture}[scale=1.1]
  \fill[\colorSigmaLight] (-\locpb,-\locpz) .. controls (0,-\locpY) .. (\locpB,-\locpz)
           -- (\locpB,\locpz) .. controls (0,\locpY) .. (-\locpb,\locpz) -- cycle;
  \fill[\colorSigmaDark] (0.7*\locpx,\locpy) arc (90:-90:\locpx cm and \locpy cm)
           -- (-3.3*\locpx,-\locpu*\locpy) arc (-90:90:\locpx cm and \locpu*\locpy cm) -- cycle;
  \draw[line width=\locpC,draw=\colorCircleL]
           (-\locpb,-\locpz) .. controls (0,-\locpY) .. (\locpB,-\locpz)
           -- (\locpB,\locpz) .. controls (0,\locpY) .. (-\locpb,\locpz) -- cycle;
  \begin{scope}[decoration={random steps,segment length=3mm}]
           \fill[\colorSigmaLight,decorate]
           (-\locpb+\locpe,-\locpz) rectangle (-\locpb-\locpe,\locpz)
           (\locpB+\locpe,-\locpz) rectangle (\locpB-\locpe,\locpz);
  \end{scope}
  \draw[line width=2*\locpC,\colorCircle]
           (0.7*\locpx,\locpy) arc (90:-90:\locpx cm and \locpy cm)
           (-3.3*\locpx,\locpu*\locpy) arc (90:-90:\locpx cm and \locpu*\locpy cm);
  \draw[dashed,line width=1.5*\locpC,\colorCircle]
           (-3.3*\locpx,\locpu*\locpy) arc (90:270:\locpx cm and \locpu*\locpy cm);
  \scopeArrow{0.15}{\arrowDefectB} \scopeArrow{0.51}{\arrowDefectB} 
  \draw[line width=\locpc,color=\colorDefect,postaction={decorate}]
           (-\locpx,\locpy) arc (90:-90:\locpx cm and \locpy cm);
  \end{scope} \end{scope}
  \draw[dashed,line width=\locpD,color=\colorDefect,postaction={decorate}]
           (-\locpx,\locpy) arc (90:270:\locpx cm and \locpy cm);
  \scopeArrow{0.78}{\arrowDefect}
  \draw[line width=\locpc,color=\colorDefect,postaction={decorate}]
           (-0.1*\locpx,0.4*\locpy+0.7*\locpA) -- (\locpB,0.4*\locpy+0.7*\locpA);
  \draw[line width=\locpc,color=\colorDefect,postaction={decorate}]
           (-0.1*\locpx,0.4*\locpy-0.7*\locpA) -- (\locpB,0.4*\locpy-0.7*\locpA);
  \end{scope}
  \scopeArrow{0.10}{\arrowDefect} \scopeArrow{0.89}{\arrowDefect}
  \draw[line width=\locpc,color=\colorDefect,postaction={decorate}]
           (-\locpb,-0.4*\locpy) -- (\locpB,-0.4*\locpy);
  \end{scope} \end{scope}
  \filldraw[\colorCircle,fill=white,line width=\locpc]
           (-0.1*\locpx,0.4*\locpy) circle (\locpA) (-0.1*\locpx,-0.4*\locpy) circle (\locpA);
  \scopeArrow{0.03}{>} \scopeArrow{0.52}{>} 
  \draw[\colorCircle,line width=\locpC,postaction={decorate}]
           (-0.1*\locpx,0.4*\locpy) circle (\locpA);
  \end{scope} \end{scope} 
  \scopeArrow{0.15}{>} \scopeArrow{0.40}{>} \scopeArrow{0.65}{>} \scopeArrow{0.90}{>} 
  \draw[\colorCircle,line width=\locpC,postaction={decorate}]
           (-0.1*\locpx,-0.4*\locpy) circle (\locpA);
  \end{scope} \end{scope} \end{scope} \end{scope} 
  \filldraw[color=\colorDefect] (\locpA-0.1*\locpx,-0.4*\locpy) circle (\locpd)
           (-\locpA-0.1*\locpx,-0.4*\locpy) circle (\locpd)
           (-0.1*\locpx,-0.4*\locpy)++(83:\locpA) circle (\locpd)
           (-0.1*\locpx,-0.4*\locpy)++(266:\locpA) circle (\locpd)
           (-0.1*\locpx,0.4*\locpy)++(95:\locpA) circle (\locpd)
           (-0.1*\locpx,0.4*\locpy)++(44:\locpA) circle (\locpd)
           (-0.1*\locpx,0.4*\locpy)++(-44:\locpA) circle (\locpd)
           (-0.1*\locpx,0.4*\locpy)++(-84:\locpA) circle (\locpd);
  \node at (0.26,-0.26) {$\delta_1$};
  \node at (-\locpx+0.05,\locpy-0.29) {$\delta_2$};
  \node at (-3.13*\locpx,0.3*\locpy) {$\partial \nn{\gamma,1}$};
  \node at (2.18*\locpx,0.8*\locpy) {$\partial \nn{\gamma,2}$};
  \node at (-\locpb-2.1,0.1*\locpz) {type (ii)\,:};
  \end{tikzpicture}
  \label{eq:pics.type.i.ii}}
  \ee ~\\[0.1em]
For paths of type (i) and (ii), the boundary $\partial \nn\gamma$ has a natural structure 
of a \defectonemanifold, with framing inherited from the given vector field on $\Sigma$.
We denote by $\ZZ(\partial\nn\gamma)$ its gluing category; thus in case (ii),
$\ZZ(\partial \nn\gamma) \,{=}\,\ZZ(\partial \nn{\gamma,1}) \boti \ZZ(\partial \nn{\gamma,2})$.
In cases (iii) and (iv), in which $\partial \nn\gamma$ is not a proper \defectonemanifold,
we still define $\ZZ(\partial \nn\gamma)$ to be 
$\ZZ(\partial \nn{\gamma,1}) \boti \ZZ(\partial \nn{\gamma,2})$ where, however, now by
definition we assign to a component $\partial \nn{\gamma,i}$ that does not meet any defect 
line on $\Sigma$ the category $ \ZZ(\partial \nn{\gamma,i}) \,{:=}\, \vect$.

Let now $\gamma \,{=}\,(\delta_1,\delta_2,...\,,\delta_n)$ be a path of defect lines in 
a fine \defectsurface\ $\Sigma$, and denote the bimodule category labeling $\delta_i$ by
$\calm_i$. Denote by 
  \be
  \ZZ(\mathring{\nn\gamma}) := \bigboxtimes_{i=1}^{n-1}\, \ZZ(\SS_{i,i+1})
  \ee
the gluing category for the \defectonemanifold\ that is the disjoint union of all
gluing boundaries inside $\nn\gamma$, and (analogously as at the beginning of Section
\ref{sec:pre-blocks})) by $\UU(\mathring{\nn\gamma}) \,{=}\,
\bigboxtimes_{i=1}^{n-1} \UU(\SS_{i,i+1})$ the category obtained from  $\ZZ(\mathring{\nn\gamma})$
by forgetting all balancings of the framed centers. It follows from our conventions that
  \be
  \UU(\mathring{\nn\gamma}) \cong \Big( \bigboxtimes_{i=1}^n \, \calm_i \boti \CC{\calm_i} \Big)
  \boxtimes \UU(\partial \nn{\gamma}) \,.
  \ee
For an object $x \,{\in}\, \ZZ(\mathring{\nn\gamma})$ that has the form
$x \,{=}\, x_\circ \boti x_\partial$ with 
$U(x_\circ) \,{\in}\, \bigboxtimes_i \calm_i \boti \CC{\calm_i}$
and $U(x_\partial) \,{\in}\, \UU(\partial \nn{\gamma})$ we define
  \be
  \widehat x := \! \int^{m_1\in\calm_1,...,m_n\in\calm_n} \!\!
  \Hom_{\calm_1 \boxtimes \CC{\calm_1} \boxtimes \cdots \boti \calm_n \boxtimes \CC{\calm_n}}
  \big( m_1\boti \cc{m_{1}} \,{\boxtimes} \cdots {\boxtimes}\, m_n \boti \cc{m_n}\,,
  U(x_\circ) \big) \otimes U(x_\partial) \,.
  \label{eq:value-comp}
  \ee
This is by construction an object in $\UU(\partial \nn{\gamma})$.

\begin{Lemma}\label{lem:x-in-ZZ}
For $x \,{\in}\, \ZZ(\mathring{\nn\gamma})$, the object $\widehat x$ defined by
\eqref{eq:value-comp} comes naturally with balancings which endow it with the structure
of an object in the gluing category $\ZZ(\partial \nn{\gamma})$.
\end{Lemma} 

\begin{proof}
The required balancings are fully determined by the balancings of the object 
$x \,{\in}\, \ZZ(\mathring{\nn\gamma})$ and the parallel transport operations 
\eqref{eq:defGamma} for the \twocells\ in $\nn\gamma$. To see this, let us write
$x\,{=}\,\bigboxtimes_{i}\, x_{i,i+1}$ with $x_{i,i+1} \,{\in}\, \ZZ(\SS_{i,i+1})$. Consider
any two neighboring defect points $P_1$ and $P_2$ on a component $\partial \nn{\gamma,i}$
of $\partial \nn\gamma$. If $P_1$ and $P_2$ correspond to neighboring defect points 
on one and the same gluing boundary $\SS_{i,i+1}$, then the balancing on the object 
$\widehat x$ is directly provided 
by the balancing of $x_{i,i+1} \,{\in}\, \ZZ(\SS_{i,i+1})$. Otherwise $P_1$ and $P_2$ 
correspond to defect points that lie on different gluing boundaries which are connected by 
a sub-path $(\delta_{j},...\,,\delta_{j'})$ of $\gamma$. In this case the parallel transport 
operations \eqref{eq:defGamma} give the relevant balancing. 
That these balancings are indeed those for the category $\ZZ(\partial \nn{\gamma})$
follows directly from the definitions.
\end{proof}

It follows that the prescription \eqref{eq:value-comp} provides a functor from
$\ZZ(\mathring{\nn\gamma})$ to $\ZZ(\partial \nn\gamma)$, provided that
the path $\gamma$ is of the type (i) or (ii).
In case of the type (iii) or (iv), $\gamma$ is a closed path, and we apply after the 
prescription \eqref{eq:value-comp} in addition the block equalizers
(see \eqref{eq:defining-block}) for all disks that are enclosed by $\gamma$.
The balancings given in the proof of Lemma \ref{lem:x-in-ZZ} are not affected 
by applying these block equalizers, so that we end up again with a functor from 
$\ZZ(\mathring{\nn\gamma})$ to $\ZZ(\partial \nn\gamma)$:

\begin{Definition}\label{defi:excision}
Let $\gamma$ be a path of defect lines in a \defectsurface\ $\Sigma$.
\Itemizeiii

\item[{\rm (i)}] 
The \emph{excision functor} for the path $\gamma$ is the functor
  \be
  \ZE_\gamma :\quad
  \ZZ(\mathring{\nn\gamma}) \xrightarrow{~~} \ZZ(\partial \nn\gamma)
  \label{eq:gen-comp-fun}
  \ee
that is obtained by the prescription given above.

\item[{\rm (ii)}] 
The \emph{contraction of $\Sigma$ along $\gamma$}, denoted as $\Sigma_{\gamma}$, is the
following \defectsurface:
 \\
For $\gamma$ of the type {\rm (i)} or {\rm (ii)}, we take
$\Sigma_{\gamma} \,{:=}\, \Sigma \,{\setminus}\, \nn\gamma$ to be
the complement of $\nn\gamma$, with $\partial \nn\gamma$ as a new gluing boundary component
for type {\rm (i)}, respectively $\partial \nn{\gamma,1} \,{\sqcup}\, \partial \nn{\gamma,2}$ 
as two new gluing boundaries for type {\rm (ii)}. For $\gamma$ of the type {\rm (iii)},
$\Sigma_{\gamma}$ is the component of $\Sigma \,{\setminus} \nn\gamma$ whose boundary
contains the boundary component $\partial \nn{\gamma,i}$ of $\nn\gamma$ that meets at
least one defect line of $\Sigma$. Finally, in case $\gamma$ is of the type {\rm (iv)}
we set $\Sigma_{\gamma} \,{:=}\, \emptyset$.

\end{itemize}
\end{Definition}

The terminology `excision functor' suggests that this functor is related to locality properties
of our construction. Indeed it will play a role when describing a factorization structure for
the modular functor (see the proof of Theorem \ref{thm:T1T2=Ttr}).

We also call the object $\ZE_\gamma(x)$ the \emph{contraction of $x$ along $\gamma$}.
The following picture shows the contraction $\Sigma_{\gamma}$ for 
the path of type (ii) that is shown in the picture \eqref{eq:pics.type.i.ii}:
    \def\locpb  {4.8}    
    \def\locpB  {4.2}    
    \def\locpc  {2.2pt}  
    \def\locpC  {0.8pt}  
    \def\locpd  {0.12}   
    \def\locpD  {1.4pt}  
    \def\locpe  {0.14}   
    \def\locpu  {1.07}   
    \def\locpv  {3.5}    
    \def\locpx  {0.76}   
    \def\locpy  {2.85}   
    \def\locpY  {2.72}
    \def\locpz  {3.3}    
  \be
  \raisebox{-7.2em}{\begin{tikzpicture}
  \filldraw[fill=\colorSigma,line width=\locpC,draw=\colorCircleL]
           (-\locpb,-\locpz) .. controls (0,-\locpY) .. (\locpB,-\locpz)
           -- (\locpB,\locpz) .. controls (0,\locpY) .. (-\locpb,\locpz) -- cycle;
  \begin{scope}[decoration={random steps,segment length=3mm}]
           \fill[\colorSigmaLight,decorate]
           (-\locpb+\locpe,-\locpz) rectangle (-\locpb-\locpe,\locpz)
           (\locpB+\locpe,-\locpz) rectangle (\locpB-\locpe,\locpz);
  \end{scope}
  \scopeArrow{0.86}{\arrowDefect}
  \draw[line width=\locpc,color=\colorDefect,postaction={decorate}] (0,-0.4*\locpy) -- +(\locpB,0);
  \draw[line width=\locpc,color=\colorDefect,postaction={decorate}] (0,0.2*\locpy) -- +(\locpB,0);
  \draw[line width=\locpc,color=\colorDefect,postaction={decorate}] (0,0.5*\locpy) -- +(\locpB,0);
  \end{scope}
  \fill[white] (0.7*\locpx,1.026*\locpy) arc (90:-90:\locpx cm and 1.026*\locpy cm)
           -- (-\locpv*\locpx,-\locpu*\locpy) arc (-90:90:\locpx cm and \locpu*\locpy cm) -- cycle;
  \scopeArrow{0.095}{>} \scopeArrow{0.41}{>} \scopeArrow{0.99}{>}
  \draw[line width=\locpc,\colorCircle,postaction={decorate}]
           (0.7*\locpx,0) circle (\locpx cm and \locpy cm);
  \end{scope} \end{scope} \end{scope}
  \scopeArrow{0.13}{>}
  \draw[line width=\locpc,\colorCircle,postaction={decorate}]
           (-\locpv*\locpx,0) circle (\locpx cm and \locpu*\locpy cm);
  \end{scope}
  \draw[dashed,line width=1.1*\locpc,\colorSigma]
           (-\locpv*\locpx,\locpu*\locpy) arc (90:270:\locpx cm and \locpu*\locpy cm);
  \scopeArrow{0.33}{\arrowDefect}
  \draw[line width=\locpc,color=\colorDefect,postaction={decorate}]
           (-\locpb,-0.4*\locpy) -- (-\locpv*\locpx+0.93*\locpx,-0.4*\locpy);
  \end{scope}
  \filldraw[color=\colorDefect] 
           (-\locpv*\locpx+0.93*\locpx,-0.4*\locpy) circle (\locpd)
           (1.63*\locpx,-0.4*\locpy) circle (\locpd) 
           (1.68*\locpx,0.2*\locpy) circle (\locpd) (1.59*\locpx,0.5*\locpy) circle (\locpd);
  \draw[line width=\locpc,color=\colorDefect] (-4*\locpx,-0.4*\locpy) -- +(\locpx,0);
  \end{tikzpicture}}
  \ee
  ~\\[-1.2em]~

\begin{Lemma} \label{Lemma:Z-same}
Let $\gamma$ be a path of defect lines in a \defectsurface\ $\Sigma$.
There is a canonical isomorphism 
  \be
  \varphi_{\gamma} :\quad 
  \Zfine(\Sigma) \xrightarrow{~\cong~} \Zfine(\Sigma_\gamma) \circ (\ZE_\gamma \boti \Id) 
  \label{eq:can-iso-bl}
  \ee
of functors, where the identity functor is applied to the gluing categories for all
gluing boundaries of $\Sigma$ that are not met by $\gamma$ $($and where in the case of
$\gamma$ being of type {\rm (iv)}, the canonical equivalence $\vect\boti\calm \,{\simeq}\, \calm$
for any finite category $\calm$ is used implicitly on the right hand side$)$.  
\end{Lemma}

\begin{proof}
For types (i) and (ii) there is even a corresponding isomorphism involving pre-block functors
which holds by construction. Moreover,
by definition of the balancings of the objects \eqref{eq:value-comp}, all holonomy 
operations for $\Sigma$ and for $\Sigma_{\gamma}$ agree. 
Thus the isomorphism follows for all types (i)\,--\,(iv).
\end{proof}

If we allow $\gamma$ to contain also free boundaries, and thus some of the gluing boundaries
$\SS_{i,i+1}$ are gluing intervals, then the tubular neighborhood $\nn\gamma$ looks
as indicated in
    \def\locpa  {2.5}    
    \def\locpA  {0.34}   
    \def\locpb  {1.5}    
    \def\locpc  {2.2pt}  
    \def\locpC  {1.8pt}  
    \def\locpD  {0.08}   
    \def\locpf  {1.1pt}  
  \be
  \raisebox{-6.0em}
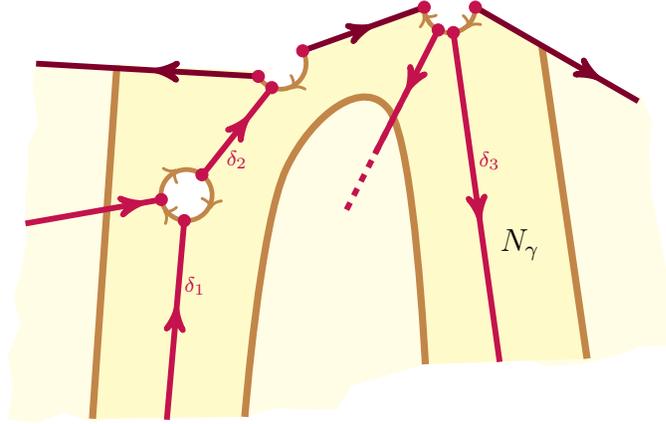
{\begin{tikzpicture}
  \coordinate (incirc1) at (0,0);
  \coordinate (incirc2) at (0.5*\locpa,0.7*\locpa);
  \coordinate (incirc2m) at (0.5*\locpa-0.86*\locpA,0.7*\locpa-0.5*\locpA);
  \coordinate (incirc3) at (1.4*\locpa,\locpa);
  \coordinate (incirc3m) at (1.4*\locpa-\locpA,\locpa);
  \begin{scope}[decoration={random steps,segment length=4mm}]
  \fill[\colorSigmaLight,decorate] (incirc2)+(-150:\locpA)+(-1.3*\locpa,0)
          -- (incirc2m) -- (incirc2) -- (incirc3) -- (2.42*\locpa,0.53*\locpa)
          -- (2.6*\locpa,-0.8*\locpa) -- (2.1*\locpa,-0.9*\locpa) -- (1.1*\locpa,-0.8*\locpa)
          -- (0.8*\locpa,-1.14*\locpa) -- (-0.95*\locpa,-1.19*\locpa) -- cycle;
  \end{scope}
  \fill[\colorSigmaDark] plot coordinates{ 
          (-0.5*\locpa,-1.19*\locpa) (-0.37*\locpa,0.69*\locpa) 
          (0.41*\locpa,0.66*\locpa) (1.4*\locpa,1.05*\locpa)
          (1.89*\locpa,0.79*\locpa) (2.13*\locpa,-0.87*\locpa)
          (1.27*\locpa,-0.90*\locpa) (1.08*\locpa,0.35*\locpa)
          (0.53*\locpa,0.04*\locpa) (0.31*\locpa,-1.18*\locpa)};
  \filldraw[line width=1.5*\locpC,color=\colorCircle,fill=\colorSigmaLight]
          plot[smooth,tension=0.65] coordinates{ 
          (1.27*\locpa,-0.90*\locpa) (1.08*\locpa,0.45*\locpa)
          (0.53*\locpa,0.14*\locpa) (0.31*\locpa,-1.18*\locpa)};
  \draw[line width=1.5*\locpC,color=\colorCircle]
          (-0.5*\locpa,-1.19*\locpa) -- (-0.37*\locpa,0.69*\locpa) 
          (1.89*\locpa,0.79*\locpa) -- (2.13*\locpa,-0.87*\locpa);
  \scopeArrow{0.47}{\arrowDefect}
  \draw[line width=\locpc,color=\colorFreebdy,postaction={decorate}]
          (incirc2)+(-150:\locpA) -- +(-1.3*\locpa,0);
  \end{scope}
  \scopeArrow{0.57}{\arrowDefect}
  \draw[line width=\locpc,color=\colorFreebdy,postaction={decorate}]
          (incirc2)+(30:\locpA) -- (incirc3m);
  \draw[line width=\locpc,color=\colorDefect,postaction={decorate}] (incirc3)++(-82:\locpA)
          -- node[right=-4pt,yshift=14pt] {$\scriptstyle \delta_3$}
          node[right=4pt,yshift=-18pt,color=black] {$\nn\gamma$} +(-82:1.77*\locpa);
  \end{scope}
  \scopeArrow{0.49}{\arrowDefect}
  \draw[line width=\locpc,color=\colorDefect,postaction={decorate}]
          (incirc3)++(243:\locpA) -- +(243:0.74*\locpa);
  \end{scope}
  \draw[line width=\locpc,color=\colorDefect,dashed]
          (incirc3)++(243:3*\locpA) -- +(243:0.81*\locpa);
  \scopeArrow{0.76}{\arrowDefect}
  \draw[line width=\locpc,color=\colorFreebdy,postaction={decorate}]
          (incirc3)++(\locpA,0) -- +(-30:\locpa);
  \end{scope}
  \scopeArrow{0.46}{\arrowDefectB}
  \draw[line width=\locpc,color=\colorDefect,postaction={decorate}]
          (incirc2)+(-110:\locpA) -- node[right=-6pt,yshift=-7pt] {$\scriptstyle \delta_2$} (incirc1);
  \end{scope}
  \scopeArrow{0.61}{\arrowDefectB}
  \draw[line width=\locpc,color=\colorDefect,postaction={decorate}]
          (incirc1) -- node[right=-2pt,yshift=8pt] {$\scriptstyle \delta_1$} +(-95:1.2*\locpa);
  \end{scope}
  \scopeArrow{0.42}{\arrowDefectB}
  \draw[line width=\locpc,color=\colorDefect,postaction={decorate}] (incirc1) -- +(190:0.87*\locpa);
  \end{scope}
  \fill[white] (incirc3) circle (\locpA);
  \filldraw[\colorCircle,fill=white,line width=\locpC] (incirc1) circle (\locpA)
           (incirc2)+(-150:\locpA) arc (-150:30:\locpA) (incirc3m) arc (-180:0:\locpA);
  \scopeArrow{0.37}{<} \scopeArrow{0.655}{<} \scopeArrow{0.98}{<}
  \draw[\colorCircle,line width=\locpf,postaction={decorate}] (incirc1) circle (\locpA);
  \end{scope} \end{scope} \end{scope}
  \scopeArrow{0.59}{<}
  \draw[\colorCircle,line width=\locpf,postaction={decorate}]
          (incirc2)+(-150:\locpA) arc (-150:30:\locpA);
  \end{scope}
  \scopeArrow{0.20}{<} \scopeArrow{0.85}{<}
  \draw[\colorCircle,line width=\locpf,postaction={decorate}] (incirc3m) arc (-180:0:\locpA);
  \end{scope} \end{scope}
  \filldraw[color=\colorDefect] (incirc1)+(-95:\locpA) circle (\locpD) 
          (incirc1)+(190:\locpA) circle (\locpD) (incirc1)+(53:\locpA) circle (\locpD) 
          (incirc2)+(-110:\locpA) circle (\locpD) (incirc2)+(30:\locpA) circle (\locpD)
          (incirc2)+(-150:\locpA) circle (\locpD) 
          (incirc3)+(243:\locpA) circle (\locpD) (incirc3)+(-82:\locpA) circle (\locpD)
          (incirc3)+(\locpA,0) circle (\locpD) (incirc3m) circle (\locpD);
  \end{tikzpicture}}
  \ee

The definition of contraction of $\Sigma$ along $\gamma$, as well as the statements of 
Lemma \ref{lem:x-in-ZZ} and Lemma \ref{Lemma:Z-same} generalize to this case accordingly.

\begin{Remark}
The construction can be further generalized to the case of an arbitrary graph $\Gamma$
formed by defect lines of $\Sigma$. Taking a tubular neighborhood $\nn{\Gamma}$ that 
encloses all gluing boundaries met by $\Gamma$ we obtain analogously categories
$\ZZ(\mathring{\nn{\Gamma}})$ and $\ZZ(\partial \nn{\Gamma})$, an excision functor 
$\ZE_{\Gamma} \colon \ZZ(\mathring{\nn{\Gamma}}) \,{\to}\, \ZZ(\partial \nn{\Gamma})$ 
between these, and a \defectsurface\ $\Sigma_{\Gamma}$ together with an isomorphism 
$\Zfine(\Sigma) \,{\to}\, \Zfine(\Sigma_{\Gamma}) \,{\circ}\, (\ZE_{\Gamma} \boti \Id)$. 
If we take the graph $\Gamma_{\!\text{tot}}$ formed by \emph{all} defect lines of $\Sigma$, we
obtain this way a canonical isomorphism
  \be
  \Zfine(\Sigma) \cong \ZE_{\Gamma_{\!\text{tot}}} 
  \ee
between the block functor for $\Sigma$ and the excision functor for the graph of
all defect lines of $\Sigma$.
\end{Remark}

\begin{Example}\label{ex:S1S2S3}
We use excision to compute the block functor for the following \defectsurface\ $\Sigma$:
 \\[-1.4em] 
    \def\locpa  {0.99}   
    \def\locpA  {0.98}   
    \def\locpc  {2.2pt}  
    \def\locpd  {0.12}   
    \def\locpD  {0.04}   
    \def\locpx  {3.5}    
    \def\locpy  {4.9}    
  \be
  \raisebox{-11.7em}{\begin{tikzpicture}
  \filldraw[fill=\colorSigma,line width=\locpc,draw=white]
           (0,0) circle (\locpx cm and \locpy cm);
  \scopeArrow{0.108}{>} \scopeArrow{0.223}{>} \scopeArrow{0.285}{>}
  \scopeArrow{0.396}{>} \scopeArrow{0.723}{>} \scopeArrow{0.795}{>}
  \draw[line width=\locpc,color=\colorCircle,postaction={decorate}]
           (0,0) circle (\locpx cm and \locpy cm);
  \end{scope} \end{scope} \end{scope} \end{scope} \end{scope} \end{scope}
  \coordinate (dashpt1) at ({\locpx*cos(81)},{\locpy*sin(81)});
  \coordinate (dashpt2) at ({\locpx*cos(261)},{\locpy*sin(261)});
  \draw[line width=\locpc,color=white,dashed]
           (dashpt1) arc (81:98:\locpx cm and \locpy cm)
           (dashpt2) arc (261:280:\locpx cm and \locpy cm);
  \coordinate (defpt1) at ({\locpx*cos(55)},{\locpy*sin(55)});
  \coordinate (defpt2) at ({\locpx*cos(68)},{\locpy*sin(68)});
  \coordinate (defpt3) at ({\locpx*cos(111)},{\locpy*sin(111)});
  \coordinate (defpt4) at ({\locpx*cos(130)},{\locpy*sin(130)});
  \coordinate (defpt5) at ({\locpx*cos(225)},{\locpy*sin(225)});
  \coordinate (defpt6) at ({\locpx*cos(237)},{\locpy*sin(237)});
  \coordinate (defpt7) at ({\locpx*cos(249)},{\locpy*sin(249)});
  \coordinate (defpt8) at ({\locpx*cos(296)},{\locpy*sin(296)});
  \coordinate (defpt9) at ({\locpx*cos(309)},{\locpy*sin(309)});
  \coordinate (Defpt1) at ({\locpa*cos(44)},{\locpA+\locpa+\locpa*sin(44)});
  \coordinate (Defpt2) at ({\locpa*cos(69)},{\locpA+\locpa+\locpa*sin(69)});
  \coordinate (Defpt3) at ({\locpa*cos(122)},{\locpA+\locpa+\locpa*sin(122)});
  \coordinate (Defpt4) at ({\locpa*cos(149)},{\locpA+\locpa+\locpa*sin(149)});
  \coordinate (Defpt5) at ({\locpa*cos(197)},{-\locpA-\locpa+\locpa*sin(197)});
  \coordinate (Defpt6) at ({\locpa*cos(226)},{-\locpA-\locpa+\locpa*sin(226)});
  \coordinate (Defpt7) at ({\locpa*cos(248)},{-\locpA-\locpa+\locpa*sin(248)});
  \coordinate (Defpt8) at ({\locpa*cos(294)},{-\locpA-\locpa+\locpa*sin(294)});
  \coordinate (Defpt9) at ({\locpa*cos(318)},{-\locpA-\locpa+\locpa*sin(318)});
  \scopeArrow{0.58}{\arrowDefectB}
  \draw[line width=\locpc,color=\colorDefect,postaction={decorate}]
           (defpt1) -- node[right=-3pt,yshift=-4pt] {$\caln_l$} (Defpt1);
  \draw[line width=\locpc,color=\colorDefect,postaction={decorate}] (defpt2) -- (Defpt2);
  \draw[line width=\locpc,color=\colorDefect,postaction={decorate}]
           (defpt3) -- node[left,yshift=5pt] {$\caln_2$} (Defpt3);
  \draw[line width=\locpc,color=\colorDefect,postaction={decorate}]
           (defpt4) -- node[left=-0.2pt,yshift=-2pt] {$\caln_1$} (Defpt4);
  \end{scope}
  \scopeArrow{0.54}{\arrowDefect}
  \draw[line width=\locpc,color=\colorDefect,postaction={decorate}]
           (defpt5) -- node[left=-6pt,yshift=8pt] {$\calm_n$} (Defpt5);
  \draw[line width=\locpc,color=\colorDefect,postaction={decorate}] (defpt6) -- (Defpt6);
  \draw[line width=\locpc,color=\colorDefect,postaction={decorate}] (defpt7) -- (Defpt7);
  \draw[line width=\locpc,color=\colorDefect,postaction={decorate}]
           (defpt8) -- node[left=-2pt,yshift=1pt] {$\calm_2$} (Defpt8);
  \draw[line width=\locpc,color=\colorDefect,postaction={decorate}]
           (defpt9) -- node[right=-3pt,yshift=5pt] {$\calm_1$} (Defpt9);
  \draw[line width=\locpc,color=\colorDefect,postaction={decorate}]
           (0,-\locpA) -- node[right=-2pt,yshift=7pt] {$\calx$} (0,\locpA);
  \end{scope}
  \begin{scope}[shift={(0,-\locpa-\locpA)}]
  \filldraw[draw=\colorSigma,line width=\locpc,fill=white] (0,0) circle (\locpa);
  \scopeArrow{0.04}{<} \scopeArrow{0.48}{<} \scopeArrow{0.76}{<}
  \draw[line width=\locpc,color=\colorCircle,postaction={decorate}] (0,0) circle (\locpa);
  \end{scope} \end{scope} \end{scope}
  \foreach \xx in {259,268,276,285} \filldraw (\xx:2.1*\locpa) circle (\locpD);
  \node[rotate=38] at (-52:0.69*\locpa) {$\rho_1$};
  \node[rotate=-49] at (207:0.67*\locpa) {$\rho_{n-1}$};
  \node at (-8:0.69*\locpa) {$\kappa_2$};
  \node at (145:0.69*\locpa) {$\kappa_1$};
  \node[\colorCircle] at (55:1.29*\locpa) {$\SS_1$};
  \end{scope}
  \begin{scope}[shift={(0,\locpa+\locpA)}]
  \filldraw[draw=\colorSigma,line width=\locpc,fill=white] (0,0) circle (\locpa);
  \scopeArrow{0.27}{<} \scopeArrow{0.60}{<} \scopeArrow{0.97}{<}
  \draw[line width=\locpc,color=\colorCircle,postaction={decorate}] (0,0) circle (\locpa);
  \end{scope} \end{scope} \end{scope}
  \foreach \xx in {77,85,93,101} \filldraw (\xx:2.1*\locpa) circle (\locpD);
  \node[rotate=40] at (133:0.69*\locpa) {$\mu_1$};
  \node[rotate=-30] at (57:0.67*\locpa) {$\mu_{l-1}$};
  \node at (-35:0.69*\locpa) {$\kappa_4$};
  \node at (195:0.69*\locpa) {$\kappa_3$};
  \node[\colorCircle] at (233:1.24*\locpa) {$\SS_2$};
  \end{scope}
  \filldraw[color=\colorDefect]
           (defpt1) circle (\locpd) (defpt2) circle (\locpd) (defpt3) circle (\locpd)
           (defpt4) circle (\locpd) (defpt5) circle (\locpd) (defpt6) circle (\locpd)
           (defpt7) circle (\locpd) (defpt8) circle (\locpd) (defpt9) circle (\locpd) 
           (Defpt1) circle (\locpd) (Defpt2) circle (\locpd) (Defpt3) circle (\locpd)
           (Defpt4) circle (\locpd) (Defpt5) circle (\locpd) (Defpt6) circle (\locpd)
           (Defpt7) circle (\locpd) (Defpt8) circle (\locpd) (Defpt9) circle (\locpd)
           (0,\locpA) circle (\locpd) (0,-\locpA) circle (\locpd); 
  \node at (-\locpx-1.5,0) {$\Sigma ~=$};
  \node at (-0.5*\locpx,0) {$\cala$};
  \node at (0.5*\locpx,0) {$\calb$};
  \node[rotate=-73] at (\locpx+0.05,0.52*\locpx) {$\kappa_2{+}\kappa_4{-}1$};
  \node[rotate=71] at (-\locpx+0.05,0.59*\locpx) {$\kappa_1{+}\kappa_3{-}1$};
  \end{tikzpicture}}
  \ee
  ~\\
We regard $\Sigma$ as a pair of pants with $\SS_{1}$ and $\SS_{2}$ as incoming boundary
circles. Consider first the case that 
$\kappa_{1} \,{=}\, \kappa_{2} \,{=}\, \kappa_{3} \,{=}\, \kappa_{4} \,{=}\,1$. 
Then  using the notation $\vec\calm\,{\Boti{\vec\rho}}\,\,{:=}\,
\cenf{\calm_1 \,\Boti{\rho^{}_1}\,\cdots\!{\Boti{\rho_{n-1}}} \! \calm_n\, }$, 
by Corollary \ref{cor:lex-framed} we obtain
  \be
  \ZZ(\SS_{1}) \cong \Funle_{\cala,\calb}( \vec\calm\,{\Boti{\vec\rho}}\,, \calx) 
  \qquad\text{and}\qquad  
  \ZZ(\SS_{2} ) \cong \Funle_{\cala,\calb}(\calx, \vec\caln\,{\Boti{\vec\mu}}\,) \,.
  \label{eq:explicit-state}
  \ee
For the path $\gamma$ of defect lines in $\Sigma$ that consists just of the defect line 
with label $\calx$, the neighborhood $\nn{\gamma}$ encloses the gluing circles $\SS_{1}$
and $\SS_{2}$. For $F \iN \ZZ(\SS_{1})$ and $G \iN \ZZ(\SS_{2})$ we then compute 
$\ZE_{\gamma}(F \boti G) \iN \ZZ(\partial \nn{\gamma}) \,{\cong}\, \ZZ(\SS_{3})^{\opp}$ 
by invoking the formula \eqref{eq:value-comp}:
  \be
  \begin{array}{ll}
  \ZE_{\gamma}(F{\boxtimes} G) \!\! &\dsty
  = \int^{x\in\calx}\!\! \Hom\big( x\boti\cc{x}\,, \big( \int^{z\in\vec\calm\,{\Boti{\vec\rho}}}
  \! \cc{z} \boti F(z) \big) \boti  \big( \int^{y \in \calx}\!\! \cc{y} \boti G(y) \big) \big)
  \Nxl1 & \dsty
  \cong \int^{z\in\vec\calm\,{\Boti{\vec\rho}}}\!\!\! \int^{y\in\calx}\!\!\! \int^{x\in\calx} 
  \! \Hom(x,F(z)) \otik \Hom(y,x) \oti \cc{z} \,\boxtimes G(y)
  \Nxl1 & \dsty
  \cong \int^{z \in \vec\calm\,{\Boti{\vec\rho}}}\! \cc{z} \boxtimes G \cir F(z) \,.
  \eear
  \ee
Under the Eilenberg--Watts equivalence, this amounts to the functor $G\,{\circ}\,F \iN 
\Funle_{\cala,\calb}( \vec\calm\,{\Boti{\vec\rho}}\,,$
                                                     $\vec\caln\,{\Boti{\vec\mu}}\,)$. 
Invoking Corollary \ref{cor:gen-cylinder} we thus we conclude that
  \be
  \Zfine(\Sigma) :\quad \ZZ(\SS_{1}) \boti \ZZ(\SS_{2}) \xrightarrow{~~~} \ZZ(\SS_{3})
  \ee
is the functor that corresponds to the composition of bimodule functors. (This is a simple
instance of the way in which algebraic structure, here composition of left exact functors,
can be extracted from our framed modular functor.) The 
case of general values of $\kappa_{1}, \kappa_{2}, \kappa_{3}$ and $\kappa_{4}$ is reduced
to the one considered above by the equivalence obtained in Remark \ref{Remark:Kappa-twist}.
When all defects involved are transparent, then via the Eilenberg--Watts isomorphisms the
composition of functors gives the tensor product in the Drinfeld center (compare Example
\ref{ex:3holedsphere}); note that this functor is the one that also the standard
Turaev--Viro constructions assigns to a pair of pants.
\end{Example}


\subsection{Fusion of defect lines along \twocells} \label{sec:fusion}

We will now show that the block functor for a \defectsurface\ remains unchanged if
two parallel defect lines are \emph{fused}, meaning that two neighboring defect lines 
are replaced by a single one with appropriate label. To see this,
consider two parallel defect lines labeled by $\calm_{\Cala}$ and ${}_{\cala}\caln$ that
end on two gluing segments with framing of index ${\pm}\kappa$ for $\kappa \,{\in}\, 2\Z$.
For our discussion only the right and left module structure, respectively, matters
(and is displayed), so that we are in the situation of Example \ref{Example:k-fr-cylinder}. 
We will see that locally we can replace this combinatorial configuration by a single
defect line labeled by the framed center $\cenf{\calm_{\Cala} \,{\Boti\kappa}
{}_{\cala}\caln}$; schematically,
    \def\locpa  {2.8}    
    \def\locpA  {2.3}    
    \def\locpc  {2.2pt}  
    \def\locpd  {0.12}   
  \be
  \raisebox{-3.6em}{\begin{tikzpicture}
  \scopeArrow{0.36}{<} \scopeArrow{0.87}{<}
  \filldraw[\colorCircle,fill=\colorSigma,line width=\locpc,postaction={decorate}]
           (0,0) node[right=21pt,yshift=-7pt,color=black]{$\kappa$} rectangle
           (\locpa,\locpA) node[left=13pt,yshift=7pt,color=black]{$-\kappa$} ; 
  \end{scope} \end{scope}
  \scopeArrow{0.73}{\arrowDefect}
  \draw[line width=\locpc,color=\colorFreebdy,postaction={decorate}]
           (0,0) -- node[left=-2pt,yshift=-7pt]{$\calm_\Cala$} (0,\locpA);
  \draw[line width=\locpc,color=\colorFreebdy,postaction={decorate}]
           (\locpa,0) -- node[right=-2pt,yshift=-7pt]{${}_\cala\caln$} +(0,\locpA);
  \end{scope}
  \filldraw[color=\colorDefect] (0,0) circle (\locpd) (0,\locpA) circle (\locpd)
                                (\locpa,0) circle (\locpd) (\locpa,\locpA) circle (\locpd);
  \begin{scope}[shift={(\locpa+1.81,0.5*\locpA)}]
  \node at (0,0) {{\Large$\rightsquigarrow$}};
  \end{scope}
  \begin{scope}[shift={(\locpa+3.2,0)}];
  \scopeArrow{0.73}{\arrowDefect}
  \draw[line width=\locpc,color=\colorDefect,postaction={decorate}] (0,0) --
  node[right=-1pt,yshift=-5pt]{$\cenf{\calm_{\Cala} \,{\Boti\kappa} {}_{\cala}\caln}$} (0,\locpA);
  \end{scope}
  \filldraw[color=\colorDefect] (0,0) circle (\locpd) (0,\locpA) circle (\locpd);
  \end{scope}
  \end{tikzpicture}
  \label{eq:M.kk.N=MkN}}
  \ee
This replacement is meant to happen locally within a generic \defectsurface,
with all data in the parts not involved in the replacement remaining unchanged.
In particular, the gluing segments on which the defect lines end will typically be
part of gluing circles or gluing intervals that contain further defect points. This is
indicated in the following picture:
 \\[-2.2em]~   
    \def\locpa  {1.9}    
    \def\locpA  {2.5}    
    \def\locpb  {1.9}    
    \def\locpc  {2.2pt}  
    \def\locpd  {0.12}   
  \be
  \raisebox{-7.2em}{\begin{tikzpicture}
  [decoration={random steps,segment length=2mm}]
  \scopeArrow{0.36}{<} \scopeArrow{0.87}{<}
  \filldraw[\colorCircle,fill=\colorSigma,line width=\locpc,postaction={decorate}]
           (0,0) node[right=15pt,yshift=-7pt,color=black]{$\kappa$} rectangle
           (\locpa,\locpA) node[left=9pt,yshift=7pt,color=black]{$-\kappa$}; 
  \end{scope} \end{scope}
  \coordinate (rectpt1) at ({-\locpb*sin(70)},{\locpb*cos(70)-\locpb});
  \filldraw[fill=\colorSigma,draw=white,line width=\locpc,postaction={decorate}]
           (0,\locpA) arc (270:200:\locpb) -- (rectpt1) arc (160:90:\locpb);
  \coordinate (rectpt2) at ({\locpb*sin(70)+\locpa},{\locpb*cos(70)+\locpA+0.7});  
  \filldraw[fill=\colorSigma,draw=white,line width=\locpc,postaction={decorate}]
           (\locpa,0) arc (90:20:\locpb) -- (rectpt2) arc (340:270:\locpb);
  \draw[\colorCircle,line width=\locpc]
           (0,\locpA) arc (270:200:\locpb) (rectpt1) arc (160:90:\locpb)
           (\locpa,0) arc (90:20:\locpb) (\locpa,\locpA) arc (270:340:\locpb);
  \scopeArrow{0.77}{\arrowDefect}
  \draw[line width=\locpc,color=\colorDefect,postaction={decorate}]
           (0,0) -- node[left=-2pt,yshift=-7pt]{$\calm_\Cala$} (0,\locpA);
  \draw[line width=\locpc,color=\colorDefect,postaction={decorate}]
           (\locpa,0) -- node[right=-2pt,yshift=-7pt]{${}_\cala\caln$} +(0,\locpA);
  \end{scope}
  \draw[line width=\locpc,color=\colorDefect]
				(0,-\locpb)++(120:\locpb) -- +(105:0.18*\locpA)
				(0,\locpA+\locpb)++(240:\locpb) -- +(265:0.18*\locpA)
				(\locpa,\locpA+\locpb)++(290:\locpb) -- +(275:0.18*\locpA)
				(\locpa,\locpA+\locpb)++(310:\locpb) -- +(285:0.18*\locpA);
  \draw[line width=0.8*\locpc,color=\colorDefect,dashed]
				(0,-\locpb)++(120:\locpb) -- +(105:0.47*\locpA)
				(0,\locpA+\locpb)++(240:\locpb) -- +(265:0.47*\locpA)
				(\locpa,\locpA+\locpb)++(290:\locpb) -- +(275:0.47*\locpA)
				(\locpa,\locpA+\locpb)++(310:\locpb) -- +(285:0.47*\locpA);
  \filldraw[color=\colorDefect] (0,0) circle (\locpd) (0,\locpA) circle (\locpd)
                                (\locpa,0) circle (\locpd) (\locpa,\locpA) circle (\locpd)
                                (0,-\locpb)+(120:\locpb) circle (\locpd)
                                (0,\locpA+\locpb)+(240:\locpb) circle (\locpd)
                                (\locpa,\locpA+\locpb)+(290:\locpb) circle (\locpd)
                                (\locpa,\locpA+\locpb)+(310:\locpb) circle (\locpd);
  \begin{scope}[shift={(\locpa+\locpb+1.5,0.5*\locpA)}]
  \node at (-0.42,0) {{\Large$\rightsquigarrow$}};
  \end{scope}
  \begin{scope}[shift={(\locpa+2*\locpb+2.0,0)}];
  \filldraw[white,fill=\colorSigma,line width=\locpc,postaction={decorate}]
           (0,\locpA+\locpb)+(200:\locpb) arc (200:270:\locpb) -- (0,0) arc (90:160:\locpb)
           -- cycle
           (0,\locpA+\locpb)+(340:\locpb) arc (340:270:\locpb) -- (0,0) arc (90:20:\locpb)
           -- cycle;
  \scopeArrow{0.79}{\arrowDefect}
  \draw[line width=\locpc,color=\colorDefect,postaction={decorate}] (0,0) --
           node[right=-1pt,yshift=-5pt]{\begin{turn}{-90}
           $\cenf{\calm_{\Cala} \,{\Boti\kappa} {}_{\cala}\caln}$\end{turn}} (0,\locpA);
  \end{scope}
  \draw[\colorCircle,line width=\locpc,postaction={decorate}]
           (0,0) arc (90:160:\locpb) (0,\locpA) arc (270:200:\locpb) 
           (0,0) arc (90:20:\locpb) (0,\locpA) arc (270:340:\locpb);
  \draw[line width=\locpc,color=\colorDefect]
				(0,-\locpb)++(120:\locpb) -- +(105:0.18*\locpA)
				(0,\locpA+\locpb)++(240:\locpb) -- +(265:0.18*\locpA)
				(0,\locpA+\locpb)++(290:\locpb) -- +(275:0.18*\locpA)
				(0,\locpA+\locpb)++(310:\locpb) -- +(285:0.18*\locpA);
  \draw[line width=0.8*\locpc,color=\colorDefect,dashed]
				(0,-\locpb)++(120:\locpb) -- +(105:0.47*\locpA)
				(0,\locpA+\locpb)++(240:\locpb) -- +(265:0.47*\locpA)
				(0,\locpA+\locpb)++(290:\locpb) -- +(275:0.47*\locpA)
				(0,\locpA+\locpb)++(310:\locpb) -- +(285:0.47*\locpA);
  \filldraw[color=\colorDefect] (0,0) circle (\locpd) (0,\locpA) circle (\locpd)
                                (0,-\locpb)+(120:\locpb) circle (\locpd)
                                (0,\locpA+\locpb)+(240:\locpb) circle (\locpd)
                                (0,\locpA+\locpb)+(290:\locpb) circle (\locpd)
                                (0,\locpA+\locpb)+(310:\locpb) circle (\locpd);
  \end{scope}
  \end{tikzpicture}
  \label{eq:M.kk.N=MkN-2}}
  \ee
We refer to this procedure as the \emph{fusion of the defect lines} along a 
constantly framed \twocell.

\begin{thm} \label{thm:fusion-rel-Del}
There is a canonical isomorphism between the block functors associated with
the two combinatorial configurations \eqref{eq:M.kk.N=MkN-2}.
\end{thm}
                                   
\begin{proof}
First note that there is a canonical equivalence between the gluing categories
for the gluing segments on which the two defect lines on the left hand side of
\eqref{eq:M.kk.N=MkN} end. Consider objects $y  \,{\in}\,\calm \,{\Boti\kappa} \caln$ 
and $\cc{x} \,{\in}\, \cc{\caln} \,{\Boti{-\kappa}} \cc{\calm}$ 
in the relevant gluing categories. After fusion, the disk $\DD$
between the two lines has disappeared and there is just one defect line left; the latter
gives rise to the block functor in the upper left corner of the following diagram:
  \be
  \hspace*{-0.4em}
  \begin{tikzcd}[column sep=2.2em]
  \int^{z \in \calm \,{\Boti\kappa} \caln  }\! 
  \Hom_{\cc{\caln} \,{\Boti{-\kappa}} \cc{\calm} \boxtimes \calm \,{\Boti\kappa} \caln }
  (\cc{z} \boti z, \cc{x} \boti y) \ar{r} \ar[,end anchor={[yshift=-0.9ex]}]{d}{f} &
  \Hom_{\calm \,{\Boti\kappa} \caln} (x,y) \ar{d}
  \\[0.4em]	 
  \int^{m \in \calm}\!\! \int^{n \in \caln}\!
  \Hom_{\cc\caln{\boxtimes}\cc\calm{\boxtimes}\calm{\boxtimes}\caln}
  (\cc{n}{\boxtimes}\cc{m}{\boxtimes}m{\boxtimes}n, U(\cc{x}) \boti U(y))
  \ar[shift right=1.3]{d} \arrow[shift left=1.3]{d} \ar{r}
  & \Hom_{\calm{\boxtimes}\caln}(U(x),U(y) )
  \ar[shift right=1.3]{dd} \arrow[shift left=1.3]{dd}
  \\[0.4em] 
  \int^{m}\!\! \int^{n}\! \Hom_{\cc\caln{\boxtimes}\cc\calm{\boxtimes}\calm{\boxtimes}\caln}
  (\cc{n}{\boxtimes}\cc{m}{\boxtimes}m{\boxtimes}n, \hspace*{8.35em}
  \\[-2.3em]	 
  \hspace*{10.5em} \int_a \cc{a .\, U(x) \,. {}^{[\kappa+3]} a} \boti U(y)) \ar{r} 
  & \Hom_{\calm{\boxtimes}\caln}(Z_{[\kappa]}(U(x)),U(y))
  \end{tikzcd}
  \label{eq:cd4fusion} 
  \ee
  ~\\[-.6em]
In this diagram each of the horizontal morphisms is a variant of the convolution property
\eqref{eq:convo} of the Hom functor. The equalizer of the two morphisms on the bottom left
of this diagram is by definition the block functor before fusion. The left column of
the diagram is of the form of \eqref{eq:coexSequ} and is thus an equalizer diagram.
\\
This diagram \eqref{eq:cd4fusion} commutes: That the upper square commutes is (after using 
that the coend over $\calm\boti\caln$ can be expressed as a double coend over $\calm$ and
$\caln$, see \Cite{Cor.\,3.12}{fScS2}) just the definition of the morphism $f$, while
commutativity of the two lower squares follows from the definition of the holonomy operation.
\\
Together it follows that the right column is an equalizer diagram as well, and thus 
that indeed also the situation with the fused defect line describes the block functor.
\end{proof}

Iterating the procedure \eqref{eq:M.kk.N=MkN-2}, one can analogously fuse any number
of neighboring defect lines. Thus in particular there is a
distinguished isomorphism between the block functors associated with the two situations 
    \def\locpa  {0.72}   
    \def\locpA  {0.95}   
    \def\locpb  {1.17}   
    \def\locpB  {2.1}    
    \def\locpc  {2.2pt}  
    \def\locpC  {1.9pt}  
    \def\locpd  {0.10}   
    \def\locpD  {0.06}   
    \def\locpe  {1.95}   
    \def\locpE  {1.31}   
    \def\locpf  {0.04}   
    \def\locpF  {1.31}   
    \def\locps  {0.19}   
    \def\locpy  {1.96}   
  \be
  \raisebox{-6.6em}{\begin{tikzpicture}
  [decoration={random steps,segment length=2mm}]
  \fill[\colorSigma,decorate] (\locpe,-\locpa-\locpB) rectangle (-\locpe,\locpa+\locpB);
  \fill[white] (0,0) circle (\locpa);
  \scopeArrow{0.53}{>}
  \draw[line width=\locpc,\colorCircle,postaction={decorate}]
           (0,-\locpa) arc (270:90:\locpa cm) node[left=15pt,yshift=-9pt,color=black]{$\mu$}
           node[left=3pt,yshift=6.6pt,color=black]{$\scriptstyle\kappa_1^{}$}
           node[left=-10pt,yshift=5pt,color=black]{$\scriptstyle\kappa_2^{}$};
  \draw[line width=\locpc,\colorCircle,postaction={decorate}]
           (0,\locpa) arc (90:-90:\locpa cm) node[right=16pt,yshift=13pt,color=black]{$\lambda$}
           node[right=4pt,yshift=-6pt,color=black]{$\scriptstyle\nu_1^{}$};
  \end{scope}
  \scopeArrow{0.69}{\arrowDefect}
  \draw[line width=\locpc,color=\colorDefect,postaction={decorate}]
           (125:\locpa) -- node[left=5pt,yshift=11pt]{$\calk_1$} (125:\locpb*\locpa+\locpb*\locpB);
  \draw[line width=\locpc,color=\colorDefect,postaction={decorate}]
           (105:\locpa) -- node[left,yshift=12pt]{$\calk_2$} (105:\locpa+\locpB);
  \draw[line width=\locpc,color=\colorDefect,postaction={decorate}]
           (55:\locpa) -- node[right=7pt,yshift=11pt]{$\calk_k$} (55:\locpb*\locpa+\locpb*\locpB);
  \end{scope}
  \scopeArrow{0.66}{\arrowDefect}
  \draw[line width=\locpc,color=\colorDefect,postaction={decorate}]
           (235:\locpb*\locpa+\locpb*\locpB) -- node[left=6pt,yshift=-9pt]{$\caln_n$} (235:\locpa);
  \draw[line width=\locpc,color=\colorDefect,postaction={decorate}]
           (285:\locpa+\locpB) -- node[right=2pt,yshift=-12pt]{$\caln_2$} (285:\locpa);
  \draw[line width=\locpc,color=\colorDefect,postaction={decorate}]
           (305:\locpb*\locpa+\locpb*\locpB) -- node[right=7pt,yshift=-10pt]{$\caln_1$} (305:\locpa);
  \end{scope}
  \filldraw[color=\colorDefect] (55:\locpa) circle (\locpd) (105:\locpa) circle (\locpd)
                         (125:\locpa) circle (\locpd) (235:\locpa) circle (\locpd)
                         (285:\locpa) circle (\locpd) (305:\locpa) circle (\locpd);
  \filldraw[color=black] (68:\locpF)  circle (\locpf) (80:\locpF)  circle (\locpf)
                         (92:\locpF)  circle (\locpf) (248:\locpF) circle (\locpf)
                         (260:\locpF) circle (\locpf) (272:\locpF) circle (\locpf);
  \begin{scope}[shift={(3.2,0)}] \node at (0,0) {and}; \end{scope}
  \begin{scope}[shift={(5.8,0)}]
  \fill[\colorSigma,decorate] (\locpE,-\locpa-\locpB) rectangle (-\locpE,\locpa+\locpB);
  \fill[white] (0,0) circle (\locpa);
  \scopeArrow{0.53}{>}
  \draw[line width=\locpc,\colorCircle,postaction={decorate}]
           (0,-\locpa) arc (270:90:\locpa cm) node[left=14pt,yshift=-9pt,color=black]{$\mu$};
  \draw[line width=\locpc,\colorCircle,postaction={decorate}]
           (0,\locpa) arc (90:-90:\locpa cm) node[right=15pt,yshift=12pt,color=black]{$\lambda$};
  \end{scope}
  \scopeArrow{0.57}{\arrowDefect}
  \draw[line width=\locpc,color=\colorDefect,postaction={decorate}]
           (0,\locpa) -- node[right=-2pt,yshift=5pt]{$\calk$} (0,\locpy);
  \end{scope}
  \scopeArrow{0.79}{\arrowDefect}
  \draw[line width=\locpc,color=\colorDefect,postaction={decorate}]
           (0,-\locpy) -- node[right=-1pt]{$\caln$} (0,-\locpa);
  \end{scope}
  \draw[line width=\locpc,color=\colorDefect]
           (0,\locps-\locpy) -- node[left=3pt,yshift=-3pt]{$\caln_n$} (-0.7,-\locpB-\locpa)
           (0,-\locpy) -- (0.3,-\locpB-\locpa)
           (0,\locps-\locpy) -- node[right=3pt,yshift=-3pt]{$\caln_1$} (0.7,-\locpB-\locpa)
           (0,\locpy-\locps) -- node[left=1pt,yshift=4pt]{$\calk_1$} (-0.7,\locpB+\locpa)
           (0,\locpy) -- (-0.3,\locpB+\locpa)
	   (0,\locpy-\locps) -- node[right=2pt,yshift=4pt]{$\calk_k$} (0.7,\locpB+\locpa);
  \filldraw[line width=\locpC,draw=\colorCircle,fill=white]
	   (0,\locpy) node {$*$} circle (\locps) (0,-\locpy) node {$*$} circle (\locps);
  \filldraw[color=\colorDefect] (90:\locpa) circle (\locpd) (270:\locpa) circle (\locpd)
           (0,\locpy-\locps) circle (\locpD) (0,-\locpy+\locps) circle (\locpD)  
           (0,\locpy)+(25.2:\locps) circle (\locpD) (0,-\locpy)+(-25.2:\locps) circle (\locpD)  
           (0,\locpy)+(110.7:\locps) circle (\locpD) (0,-\locpy)+(-70:\locps) circle (\locpD)  
           (0,\locpy)+(157:\locps) circle (\locpD) (0,-\locpy)+(-157:\locps) circle (\locpD);
  \filldraw[color=black]
           (0.18,\locpy+2.1*\locps) circle (\locpf) (0.02,\locpy+2.3*\locps) circle (\locpf)
           (-0.18,-\locpy-2.1*\locps) circle (\locpf) (-0.02,-\locpy-2.3*\locps) circle (\locpf);
  \end{scope}
  \end{tikzpicture}}
  \label{eq:fusion}
  \ee
for any numbers $n$ of incoming and $k$ of outgoing defect points, with 
$\calk \,{=}\, \calk_1 \,\Boti{\kappa_1} \,\,{\cdots} \!{\Boti{\kappa_{k-1}}}\! \calk_k$ and
$\caln \,{=}\, \caln_1 \,\Boti{\nu_1} \,\,{\cdots} \!{\Boti{\nu_{n-1}}}\! \caln_n$.
Here the asterisk at the additional gluing circles on the right hand side
(for which, for better readability, we omit the orientation and labelings) indicates that
when calculating a block functor they have to be evaluated on the
relevant distinguished fusion objects (as introduced after Example \ref{ex:MNKfusion}).


\section{The modular functor} \label{sec:5}

\subsection{Factorization}\label{sec:facto}

We now study how our description of block functors fits together with the gluing of
\defectsurface s, assuming for now that both the initial and the glued surface are fine. In
the literature, this issue is often formulated as the behavior of blocks under
\emph{factorization}, and we adopt this term here. Our bicategorical setting makes it
manifest that factorization amounts to a \emph{structure}, rather than being a property. We
will show that such a structure 
exists, but will not investigate its uniqueness in the present paper. The general case of
factorization for not necessarily fine \defectsurface s will follow from the definition of
the block functor of a non-fine surface as a limit, see Theorem \ref{thm:Facto}.

According to Definition \ref{def:Bordcats}, the horizontal composition of morphisms in
the bicategory  $\Bfrdefdec$ is given by the gluing of surfaces. A 
2-functor comes with additional data, implementing the compatibility of
horizontal composition. The purpose of the present subsection is to exhibit
these data for the  block functor $\ZZ$. We refer to the resulting structure as the
\emph{factorization structure} of the 2-functor.

To define this factorization structure, consider the situation that the boundary of a 
\gfine\ \defectsurface\ $\Sigma$ contains gluing boundary components $\LL_1 \,{=}\, \LL$
and $\LL_2 \,{=}\, \cc{\LL}$. Denote by $\Tr_\LL(\Sigma)$ the (fine)
\defectsurface\ that results from gluing $\LL$ and $\cc{\LL}$, that is, by identifying
the boundary components corresponding to their parametrizations. In case that the
\defectsurface\ is a disjoint sum $\Sigma \,{=}\, \Sigma_1 \,{\sqcup}\, \Sigma_2$ 
and that $\LL \,{\subseteq}\, \partial\Sigma_1$ and $\cc{\LL} \,{\subseteq}\, \partial\Sigma_2$,
the surface $\Tr_{\LL}(\Sigma) \,{=}\, \Sigma_{2} \,{\circ}\, \Sigma_{1}$ is the (partial)
horizontal composition of the 1-morphisms $\Sigma_1$ and $\Sigma_2$ in $\Bfrdefdec$. 

We will show below that there is a canonical isomorphism that expresses the block functor
for the surface $\Tr_\LL(\Sigma)$, by taking a coend in the gluing category $\ZZ(\LL)$:
  \be
  \Zfine({\Tr_\LL}{(\Sigma)}) \,\cong \int^{z\in\ZZ(\LL)}\! \Zfine(\Sigma)({-}\boti z\boti \cc{z}) \,.
  \label{eq:fact-expl}
  \ee
Here on the right hand side we evaluate the block functor on objects $z$ and $\cc{z}$ 
in the gluing categories for the gluing boundaries $\LL$ and $\cc{\LL}$, respectively.  
We select this canonical isomorphism as the definition of the factorization structure of 
the modular functor. 

To identify this canonical isomorphism, we
first recall the notion of a defect cylinder, see e.g.\ the picture \eqref{eq:refinedcyl}. 
It follows directly from the definition that any defect cylinder $\Sigma$ satisfies
  \be
  {\Tr_\LL}{(\Sigma \,{\sqcup}\, \Sigma)} = \Sigma \,.
  \ee

\begin{Proposition}\label{prop:CylCyl=Cyl}
Let $\Sigma$ be a defect cylinder over a \defectonemanifold. Then the endofunctor $\Zfine(\Sigma)$
is a strict idempotent under horizontal composition:
  \be
  \Zfine(\Sigma) \hoco \Zfine(\Sigma) = \Zfine(\Sigma) \,.
  \ee
\end{Proposition}

\begin{proof}
This follows directly from the Eilenberg--Watts calculus: according to Corollary
\ref{cor:gen-cylinder}, the  block functor (with values in $\vect$) for a defect cylinder is a
Hom functor, which via the Eilenberg--Watts equivalences corresponds to the identity functor.
When working instead with Deligne products, the  isomorphism is provided by the variant 
\eqref{eq:convo} of the Yoneda lemma: Denoting the two copies of $\Sigma$ by $\Sigma_1$ and 
$\Sigma_2$, such that $\partial \Sigma_1 \,{=}\,\cc{\LL} \,{\sqcup}\, \LL$ and 
$\partial \Sigma_2 \,{=}\, \LL \,{\sqcup}\, \cc{\LL}$, and fixing objects
$\cc{y} \,{\in}\,  \ZZ(\cc{\LL}) \,{\cong}\, \ZZ(\LL)\opp$ and $w \,{\in}\, \ZZ(\LL)$ in 
the gluing categories for $\cc{\LL} \,{\subset}\,\Sigma_{1}$ and 
$\LL \,{\subset}\, \Sigma_2$, respectively, we have
  \be
  \bearll
  \big( \Zfine(\Sigma_1) \circ \Zfine(\Sigma_2) \big) (\cc y \boti w) \!\!&\dsty
  \stackrel{\eqref{eq:EW-xxx}}\cong
  \int^{x \in \ZZ(\LL)}\! \Zfine(\Sigma{\sqcup}\Sigma)(\cc y \boti x \boti \cc x \boti w) 
  \Nxl1 & \dsty
  \! \stackrel{\eqref{eq:ZZ-cylinder}}\cong
  \int^{x \in \ZZ(\LL)}\! \Hom_{\ZZ(\LL)}(y,x) \oti \Hom_{\ZZ(\LL)}(x,w)
  \Nxl1 &
  \stackrel{\eqref{eq:convo}}\cong \Hom_{\ZZ(\LL)}(y,w) \,=\, \Zfine(\Sigma)(\cc y \boti w) \,.
  \\[-2.7em]~
  \eear
  \ee
\end{proof}

When combined with the results of Section \ref{subsection:Contrac}, we obtain an analogous
canonical isomorphism also in more general situations:

\begin{thm}\label{thm:T1T2=Ttr}
Let $\Sigma$ be a fine \defectsurface\ with gluing boundary components $\LL$ and $\cc{\LL}$
such that the glued surface $\Tr_\LL(\Sigma)$ is again fine. There is a canonical isomorphism 
  \be
  \int^{x \in \ZZ(\LL)} \! \Zfine(\Sigma)({-} \boti x \boti \cc x)
  \,\cong\, \Zfine(\Tr_{\LL}(\Sigma))(-) 
  \ee
of functors.
\end{thm}

\begin{proof}
Consider the subset of defect lines of $\Sigma$, to be called the \emph{link $\gamma_{\LL}$
of the boundary component $\LL$}, that arises as follows: Consider all \twocells\ that have
a common boundary with $\LL$ and take all the defect lines on any of these \twocells\ 
that do not have an end point on $\LL$. The following picture shows the link $\gamma_{\LL}$,
indicated by the thickened defect lines, for a sample boundary component $\LL$:
    \def\locpa  {2.9}    
    \def\locpA  {0.8}    
    \def\locpb  {0.45}   
    \def\locpc  {2.2pt}  
    \def\locpd  {0.12}   
    \def\locpD  {0.10}   
  \be
  \raisebox{-8.8em}
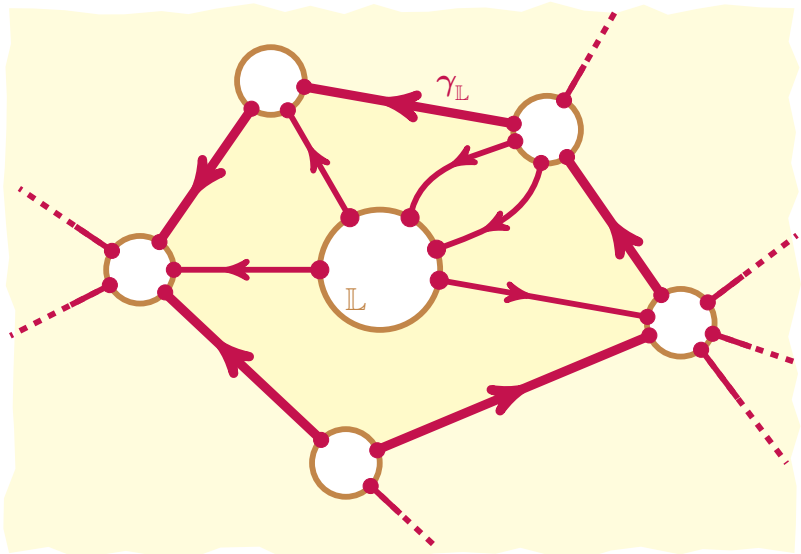
{\begin{tikzpicture}
  [decoration={random steps,segment length=4mm}]
  \fill[\colorSigma,decorate] (-1.7*\locpa,-1.3*\locpa) rectangle (1.9*\locpa,1.2*\locpa);
  \coordinate (oucirc1) at (40:\locpa);
  \coordinate (oucirc2) at (120:\locpa);
  \coordinate (oucirc3) at (180:1.1*\locpa);
  \coordinate (oucirc4) at (260:0.9*\locpa);
  \coordinate (oucirc5) at (350:1.4*\locpa);
  \scopeArrow{0.12}{\arrowDefect} \scopeArrow{0.30}{\arrowDefect}
  \scopeArrow{0.47}{\arrowDefectB} \scopeArrow{0.71}{\arrowDefect} \scopeArrow{0.92}{\arrowDefect}
  \filldraw[line width=1.6*\locpc,draw=\colorDefect,fill=\colorSigmaDark,postaction={decorate}]
             (oucirc1) -- (oucirc2) -- (oucirc3) -- (oucirc4) -- (oucirc5) -- cycle;
  \end{scope} \end{scope} \end{scope} \end{scope} \end{scope}
  \scopeArrow{0.50}{\arrowDefect}
  \draw[line width=\locpc,color=\colorDefect,postaction={decorate}]
	     (oucirc1)++(200:\locpb) [out=200,in=80] to (60:\locpA);
  \draw[line width=\locpc,color=\colorDefect,postaction={decorate}] (0,0) -- (oucirc5);
  \end{scope} 
  \scopeArrow{0.65}{\arrowDefect}
  \draw[line width=\locpc,color=\colorDefect,postaction={decorate}]
	     (oucirc1)++(260:\locpb) [out=260,in=20] to (20:\locpA);
  \draw[line width=\locpc,color=\colorDefect,postaction={decorate}] (0,0) -- (oucirc2);
  \draw[line width=\locpc,color=\colorDefect,postaction={decorate}] (0,0) -- (oucirc3);
  \end{scope} 
  \draw[line width=\locpc,color=\colorDefect,dashed] (oucirc1) -- ++(60:4.2*\locpb)
             (oucirc3) -- ++(146:4.3*\locpb)  (oucirc3) -- ++(207:4.3*\locpb)
             (oucirc4) -- ++(317:3.6*\locpb)  (oucirc5) -- ++(37:4.3*\locpb)
             (oucirc5) -- ++(307:5.2*\locpb)  (oucirc5) -- ++(341:3.7*\locpb);
  \draw[line width=\locpc,color=\colorDefect] (oucirc1) -- ++(60:2.2*\locpb)
             (oucirc3) -- ++(146:2.3*\locpb)  (oucirc3) -- ++(207:2.3*\locpb)
             (oucirc4) -- ++(317:2.2*\locpb)  (oucirc5) -- ++(37:2.3*\locpb)
             (oucirc5) -- ++(307:2.8*\locpb)  (oucirc5) -- ++(341:2.2*\locpb);
  \filldraw[\colorCircle,fill=white,line width=\locpc]
	     (0,0) circle (\locpA) node[left,yshift=-11pt] {$\LL$}
	     (oucirc1) circle (\locpb) (oucirc2) circle (\locpb)
             (oucirc3) circle (\locpb) (oucirc4) circle (\locpb) (oucirc5) circle (\locpb);
  \filldraw[color=\colorDefect]
             (20:\locpA) circle (\locpd)  (oucirc1)++(260:\locpb) circle (\locpD)
             (60:\locpA) circle (\locpd)  (oucirc1)++(200:\locpb) circle (\locpD)
             (120:\locpA) circle (\locpd) (oucirc2)++(300:\locpb) circle (\locpD)
             (180:\locpA) circle (\locpd) (oucirc3)++(\locpb,0) circle (\locpD)
             (350:\locpA) circle (\locpd) (oucirc5)++(170:\locpb) circle (\locpD)
             (oucirc1)++(60:\locpb) circle (\locpD)  (oucirc1)++(170:\locpb) circle (\locpD)
             (oucirc1)++(306:\locpb) circle (\locpD)
             (oucirc2)++(235:\locpb) circle (\locpD) (oucirc2)++(350:\locpb) circle (\locpD)
             (oucirc3)++(55:\locpb) circle (\locpD)  (oucirc3)++(146:\locpb) circle (\locpD)
             (oucirc3)++(207:\locpb) circle (\locpD) (oucirc3)++(318:\locpb) circle (\locpD)
             (oucirc4)++(22:\locpb) circle (\locpD)  (oucirc4)++(137:\locpb) circle (\locpD)
             (oucirc4)++(317:\locpb) circle (\locpD)
             (oucirc5)++(37:\locpb) circle (\locpD)  (oucirc5)++(125:\locpb) circle (\locpD)
             (oucirc5)++(202:\locpb) circle (\locpD) (oucirc5)++(307:\locpb) circle (\locpD)
             (oucirc5)++(341:\locpb) circle (\locpD);
  \node[color=\colorDefect] at (68:0.9*\locpa) {{\boldmath $\gamma_\LL^{}$}}; 
  \end{tikzpicture}}
  \ee
Contracting $\Sigma$ along the the graph $\gamma_{\LL}$ in the way described in Section
\ref{subsection:Contrac} results in a disconnected \defectsurface\ $\Sigma_{\gamma_\LL^{}}$
such that the component of $\Sigma_{\gamma_\LL^{}}$ that contains $\LL$ is a defect cylinder
over $\LL$. Performing the same construction also over the gluing boundary $\cc{\LL}$ of
$\Sigma_{\gamma_\LL^{}}$ then gives another \defectsurface\ 
$(\Sigma_{\gamma_\LL^{}})_{\gamma_{\cc\LL^{}}}$, among the components of which there are
the defect cylinders over $\LL$ and $\cc{\LL}$, which are
identical as \defectsurface s (each having boundary $\LL \,{\sqcup}\, \cc{\LL}$). Applying
Proposition \ref{prop:CylCyl=Cyl} to the so obtained defect cylinders, it follows that
  \be
  \bearll \dsty
  \int^{x\in\ZZ(\LL)} \! \Zfine(\Sigma)({-} \boti x \boti \cc x) \!\!& \dsty
  \cong \int^{x\in\ZZ(\LL)}\! \Zfine\big((\Sigma_{\gamma_\LL^{}})_{\gamma_{\cc{\LL}^{}}}\big)
  ({-} \boti x \boti \cc x)
  \Nxl4 &
  \cong\, \Zfine(\Tr_\LL\big( (\Sigma_{\gamma_\LL^{}}){}_{\gamma_{\cc\LL}^{}} \big) )(-) 
  \,\cong\, \Zfine(\Tr_\LL(\Sigma))(-) \,. 
  \eear
  \ee
Here the first and last isomorphisms hold by the canonical isomorphism involving 
the excision functor as obtained in Lemma \ref{Lemma:Z-same}, and we also use that
$\Tr_\LL \big( (\Sigma_{\gamma_\LL^{}})_{\gamma_{\cc\LL^{}}} \big) \,{\cong}\, 
\big( (\Tr_\LL^{}(\Sigma))_{\gamma_\LL^{}} \big){}_{\gamma_{\cc\LL^{}}}$.
\end{proof}

It follows directly from the locality of the factorization in Theorem \ref{thm:T1T2=Ttr} that
factorizations performed on several boundary components commute:

\begin{cor} \label{cor:locality-Fact}
Let $(\Sigma_{1}, \Sigma_{2},\Sigma_{3})$ be an ordered triple of fine
\defectsurface s that can composed in the given order. Then the two isomorphisms 
  \be
  \Zfine(\Sigma_3) \circ \Zfine(\Sigma_2) \circ \Zfine(\Sigma_1)
  \xrightrightarrows{~\cong~} \Zfine(\Sigma_3 \,{\circ}\, \Sigma_2 \,{\circ}\, \Sigma_1)
  \label{eq:two-facto}
  \ee
that correspond to either first gluing the surfaces $\Sigma_{3}$ and $\Sigma_{2}$ and then the
resulting surface with $\Sigma_{1}$, or else first gluing $\Sigma_{2}$ and $\Sigma_{1}$ and then 
the result with $\Sigma_{3}$, are equal. 
\end{cor}


\subsection{Transparent defects and \filldisk s} \label{sec:transparent}

To be able to define the modular functor also on \defectsurface s that are not fine, we need
a suitable notion of \emph{refinement} of surfaces. This will be introduced in Section
\ref{sec:refin-defect-surf}. As a preparation, in the present subsection we provide two
notions that will be convenient in that context: a subclass of defect lines that we will 
call transparent defects, and a subclass of \defectsurface s that we will call \filldisk s.

We start with

\begin{Definition} \label{def:filling}
Let $\LL$ be a gluing circle or gluing interval of a \defectsurface\ $\Sigma$, and regard
$\Sigma$ as a bordism with domain $\LL$.
$($In case $\LL$ is an interval, it is convenient to think of it concretely as a half-circle,
as we have done before, e.g.\ in Proposition $\ref{prop:preblock-compose}.)$
We say that $\LL$ is \emph{fillable} iff there exists a \defectsurface\ 
$\DD_\LL\colon \emptyset \,{\to}\, \LL$ such that its underlying surface $D_\LL$ has the topology 
of a disk $($not containing any additional gluing circles$)$. 
\end{Definition}

If a gluing boundary $\LL$ is fillable in this way, we call the \defectsurface\ 
$\Sigma \,{\circ}\, \DD_\LL$ a \emph{filling} of $\Sigma$ by the disk $\DD_\LL$.
The filling of a \defectsurface\ by any finite number of disks is defined analogously.

If a gluing circle or gluing interval is fillable, then it has necessarily an even number of
defect points, coming in pairs which carry the same label and have opposite orientation. The
following picture shows an example of a fillable circle $\LL$ and a fillable interval
$\LL'$ and of corresponding disks $\DD_\LL$ and $\DD_{\LL'}$ that fill them:
    \def\locpa  {1.3}    
    \def\locpA  {1.6}    
    \def\locpb  {1.5}    
    \def\locpc  {2.2pt}  
    \def\locpd  {0.12}   
    \def\locpe  {0.7}    
  $$
  \begin{tikzpicture}
  \begin{scope}[decoration={random steps,segment length=2mm}]
  \fill[\colorSigma,decorate] (0,-\locpa-\locpb) rectangle (-\locpa-\locpe,\locpa+\locpb);
  \end{scope}
  \fill[white] (0,0) circle (\locpa);
  \scopeArrow{0.17}{>} \scopeArrow{0.51}{>} \scopeArrow{0.91}{>}
  \draw[line width=\locpc,\colorCircle,postaction={decorate}] (0,-\locpa) arc (270:90:\locpa);
  \end{scope} \end{scope} \end{scope}
  \scopeArrow{0.82}{\arrowDefect}
  \draw[line width=\locpc,color=\colorFreebdy,postaction={decorate}]
           (0,\locpa) -- node[right=-2pt,yshift=-1pt]{$\calm$} +(0,\locpb);
  \end{scope}
  \scopeArrow{0.70}{\arrowDefect}
  \draw[line width=\locpc,color=\colorFreebdy,postaction={decorate}]
           (0,-\locpa-\locpb) -- node[right=-2pt,yshift=-6pt]{$\calm$} +(0,\locpb);
  \end{scope}
  \filldraw[color=\colorDefect] (0,\locpa) circle (\locpd) (0,-\locpa) circle (\locpd)
           (133:\locpa) circle (\locpd) node[left=-2pt,yshift=3pt] {$\caln$}
           node[right=-3pt,yshift=-5pt] {$+$}
           (227:\locpa) circle (\locpd) node[left,yshift=-4pt] {$\caln$}
           node[right=-3pt,yshift=5pt] {$-$};
  \node at (-\locpa-\locpe-1.1,0) {$\LL ~=$};
  \begin{scope}[shift={(7.0,0)}]
  \fill[\colorSigma] (0,-\locpA) arc (270:90:\locpA);
  \scopeArrow{0.17}{>} \scopeArrow{0.51}{>} \scopeArrow{0.91}{>}
  \draw[line width=\locpc,\colorCircle,postaction={decorate}] (0,-\locpA) arc (270:90:\locpA);
  \end{scope} \end{scope} \end{scope}
  \scopeArrow{0.59}{\arrowDefect}
  \draw[line width=\locpc,color=\colorFreebdy,postaction={decorate}]
           (0,-\locpA) -- node[right=-2pt,yshift=-7pt]{$\calm$} (0,\locpA);
  \end{scope}
  \scopeArrow{0.52}{\arrowDefect}
  \draw[line width=\locpc,color=\colorDefect,postaction={decorate}]
           (133:\locpA) to [out=330,in=40] (227:\locpA) node[right=-6pt,yshift=17pt] {$\caln$};
  \end{scope}
  \filldraw[color=\colorDefect] (0,\locpA) circle (\locpd) (0,-\locpA) circle (\locpd)
	   (133:\locpA) circle (\locpd) node[left=-3pt,yshift=4pt] {$-$}
           (227:\locpA) circle (\locpd) node[left=-3pt,yshift=-5pt] {$+$};
  \node at (-\locpA-1.3,0) {$\DD_\LL ~=$};
  \end{scope}
  \end{tikzpicture}
  $$
  \be
  \raisebox{-4.9em}{\begin{tikzpicture}
  \begin{scope}[decoration={random steps,segment length=4mm}]
  \fill[\colorSigma,decorate] (-\locpa-\locpe,-\locpa-\locpe) rectangle (\locpa+\locpe,\locpa+\locpe);
  \end{scope}
  \fill[white] (0,0) circle (\locpa);
  \scopeArrow{0.51}{>}
  \draw[line width=\locpc,\colorCircle,postaction={decorate}] (0,0) circle (\locpa);
  \end{scope}
  \filldraw[color=\colorDefect]
     (50:\locpa) circle (\locpd) node[left=-3pt,yshift=-6pt] {$+$} node[above=1pt,xshift=2pt] {$\calk$}
    (100:\locpa) circle (\locpd) node[below=-1pt,xshift=2pt] {$+$} node[above=1pt] {$\calm$}
    (140:\locpa) circle (\locpd) node[right=-3pt,yshift=-5pt] {$-$} node[left=-2pt,yshift=2pt] {$\caln$}
    (220:\locpa) circle (\locpd) node[right=-4pt,yshift=6pt] {$+$} node[left=-1pt,yshift=-1pt] {$\caln$}
    (260:\locpa) circle (\locpd) node[above=-2pt,xshift=1pt] {$-$} node[below,xshift=-2pt] {$\calm$}
    (310:\locpa) circle (\locpd) node[above=-3pt,xshift=-5pt] {$-$} node[below=1pt] {$\calk$};
  \node at (-\locpa-\locpe-1.1,0) {$\LL' ~=$};
  \begin{scope}[shift={(7.3,0)}]
  \scopeArrow{0.51}{>}
  \filldraw[line width=\locpc,\colorCircle,fill=\colorSigma,postaction={decorate}] (0,0) circle (\locpA);
  \end{scope}
  \scopeArrow{0.52}{\arrowDefect}
  \draw[line width=\locpc,color=\colorDefect,postaction={decorate}]
           (50:\locpA) to [out=225,in=130] (310:\locpA) node[left=-3.5pt,yshift=36pt] {$\calk$};
  \draw[line width=\locpc,color=\colorDefect,postaction={decorate}]
           (100:\locpA) to [out=285,in=75] (260:\locpA) node[right=-1pt,yshift=14pt] {$\calm$};
  \draw[line width=\locpc,color=\colorDefect,postaction={decorate}]
           (220:\locpA) to [out=30,in=320] (140:\locpA) node[right=7pt,yshift=-47pt] {$\caln$};
  \end{scope}
  \filldraw[color=\colorDefect]
            (50:\locpA) circle (\locpd) node[right=-2pt,yshift=4pt] {$-$}
	   (100:\locpA) circle (\locpd) node[above=-2pt] {$-$}
	   (140:\locpA) circle (\locpd) node[left=-3pt,yshift=5pt] {$+$}
	   (220:\locpA) circle (\locpd) node[left=-2pt,yshift=-4pt] {$-$}
	   (260:\locpA) circle (\locpd) node[below,xshift=-1pt] {$+$}
	   (310:\locpA) circle (\locpd) node[below=-3pt,xshift=5pt] {$+$};
  \node at (-\locpA-1.3,0) {$\DD_{\LL'_{}}^{} ~=$};
  \end{scope}
  \end{tikzpicture}}
  \ee
By abuse of terminology, we apply the notion of fillability also to circles that do not
contain any defect point and thus are not proper \defectonemanifold s.

Next we introduce a convention which will be used later on:

\begin{Definition} \label{def:IIkappa}
For any finite tensor category $\cala$ and any even integer $\kappa$ we call the
$\cala$-bi\-mo\-dule
  \be
  \II_\kappa := {}^{-\kappa\!}\!\cala 
  \ee
the \emph{$\kappa$-trans\-parent defect label} $($of type $\cala)$. A defect line labeled 
by $\II_\kappa$ is called a \emph{$\kappa$-transparent defect} $($of type $\cala)$.
\\
For brevity, in case that $\kappa \,{=}\, 0$, we call $\II_0 \,{=}\, \cala$ just the
\emph{transparent defect label} $($of type $\cala)$, and a defect line labeled by 
$\II_0$ just a \emph{transparently labeled} defect.
\end{Definition}

This notation and terminology is justified by the fact that $\II_\kappa$ behaves like 
a unit for the $\kappa$-twisted center: according to formula \eqref{eq:kappaA=Akappa}
and Lemma \ref{lem:ZMA=M} for an $\cala$-bimodule category $\calm$ there are specified 
equivalences $\II_\kappa \,{\simeq}\, \cala^{-\kappa}$ and
  \be
  \calm\stackrel \kappa\boxtimes \II_\kappa \,\simeq\, \calm \,\simeq\,
  \cenf{\II_{\kappa} {\stackrel{\kappa}{\boti}} \calm} \,.
  \label{eq:M.I_n=M}
  \ee

\begin{Remark}
For any $\kappa,\kappa' \,{\in}\, 2\Z$ the double dual induces an equivalence 
  \be
  {}^{\kappa}_{\phantom0}\II_0^{\kappa'} \simeq \II_0^{\kappa+\kappa'} \,,
  \ee
in particular 
  \be
  {}^{\kappa}_{}\II_{0}^{-\kappa} \,\simeq\, \II_0 \,.
  \label{eq:kappaII-kappa=II}
  \ee
Together with formula \eqref{eq:[k]l=k+l}, this implies that for any $\kappa \,{\in}\, 2\Z$ we have
  \be
  \calm {\stackrel{\kappa}{\boti}} \II_0 {\stackrel{-\kappa}{\boti}} \caln
  \simeq \calm {\stackrel{0}{\boti}} {}^{\kappa}_{}\II_{0}^{-\kappa} {\stackrel{0}{\boti}} \caln
  \simeq \calm {\stackrel{0}{\boti}} \II_0 {\stackrel{0}{\boti}} \caln 
  \simeq \calm {\stackrel{0}{\boti}} \caln.
  \label{eq:M.II0.N=M.N}
  \ee
\end{Remark}

\medskip

\begin{Example}
The equivalences \eqref{eq:M.I_n=M} lead to distinguished objects in certain gluing categories:
Let $\cala$ be a finite tensor category and $\calm$ an $\cala$-bimodule. Via the Eilenberg--Watts
equivalences \eqref{Phile.Psile} the gluing categories for the \defectonemanifold s
    \def\locpa  {0.73}   
    \def\locpA  {1.9}    
    \def\locpc  {2.2pt}  
    \def\locpd  {0.12}   
    \def\locpx  {270}    
    \def\locpX  {90}     
    \def\locpy  {40}     
    \def\locpY  {-40}    
    \def\locpz  {140}    
    \def\locpZ  {-140}   
  $$  
  \begin{tikzpicture}
  \scopeArrow{0.31}{<} \scopeArrow{0.55}{<} \scopeArrow{0.95}{<}
  \filldraw[line width=\locpc,\colorCircle,fill=white,postaction={decorate}]
           node[left=1.9cm,color=black] {$\LII_\kappa^{\sss\nearrow}(\calm) ~:= $}
           (0,0) circle (\locpa)
           node[above=0.66cm,xshift=7pt,color=black]{$\scriptstyle \kappa$};
  \end{scope} \end{scope} \end{scope}
  \scopeArrow{0.61}{\arrowDefectB}
  \draw[line width=\locpc,color=\colorDefect,postaction={decorate}]
           (\locpx:\locpa) -- (\locpx:0.9*\locpA)
           node[right=-3pt,yshift=2pt,color=\colorDefect] {$\calm$};
  \end{scope}
  \scopeArrow{0.81}{\arrowDefect}
  \draw[line width=\locpc,color=\colorTransp,postaction={decorate}]
           (\locpy:\locpa) -- (\locpy:\locpA)
           node[right=-2pt,color=\colorDefect] {$\II_\kappa$};
  \draw[line width=\locpc,color=\colorDefect,postaction={decorate}]
           (\locpz:\locpa) -- (\locpz:\locpA) node[left=-3pt,color=\colorDefect] {$\calm$};
  \end{scope}
  \filldraw[color=\colorDefect] (\locpx:\locpa) circle (\locpd)
           (\locpy:\locpa) circle (\locpd) (\locpz:\locpa) circle (\locpd);
  \begin{scope}[shift={(7.9,0)}] 
  \scopeArrow{0.02}{<} \scopeArrow{0.43}{<} \scopeArrow{0.73}{<}
  \filldraw[line width=\locpc,\colorCircle,fill=white,postaction={decorate}]
           node[left=1.9cm,color=black] {$\LII_\kappa^{\sss\swarrow}(\calm) ~:= $}
           (0,0) circle (\locpa)
           node[below=0.6cm,xshift=9pt,color=black]{$\scriptstyle-\kappa$};
  \end{scope} \end{scope} \end{scope}
  \scopeArrow{0.81}{\arrowDefect}
  \draw[line width=\locpc,color=\colorDefect,postaction={decorate}]
           (\locpX:\locpa) -- (\locpX:0.9*\locpA)
           node[left=-2pt,yshift=-2pt,color=\colorDefect] {$\calm$};
  \end{scope}
  \scopeArrow{0.65}{\arrowDefectB}
  \draw[line width=\locpc,color=\colorDefect,postaction={decorate}]
           (\locpY:\locpa) -- (\locpY:\locpA) node[right=-5pt,color=\colorDefect] {$\calm$}
           node[right=28pt] {~};    
  \draw[line width=\locpc,color=\colorTransp,postaction={decorate}]
           (\locpZ:\locpa) -- (\locpZ:\locpA)
           node[left=-7pt,yshift=-1pt,color=\colorDefect] {$\II_\kappa$};
  \end{scope}
  \filldraw[color=\colorDefect] (\locpX:\locpa) circle (\locpd)
           (\locpY:\locpa) circle (\locpd) (\locpZ:\locpa) circle (\locpd);
  \end{scope}
  \end{tikzpicture}
  $$  
~\\[-4.0em] 
  \be
  \raisebox{-4.4em}{\begin{tikzpicture}
  \scopeArrow{0.31}{<} \scopeArrow{0.55}{<} \scopeArrow{0.95}{<}
  \filldraw[line width=\locpc,\colorCircle,fill=white,postaction={decorate}]
           node[left=1.9cm,color=black] {$\LII_\kappa^{\sss\nwarrow}(\calm) ~:= $}
           (0,0) circle (\locpa)
           node[above=0.66cm,xshift=7pt,color=black]{$\scriptstyle \kappa$};
  \end{scope} \end{scope} \end{scope}
  \scopeArrow{0.61}{\arrowDefectB}
  \draw[line width=\locpc,color=\colorDefect,postaction={decorate}]
           (\locpx:\locpa) -- (\locpx:0.9*\locpA)
           node[right=-3pt,yshift=2pt,color=\colorDefect] {$\calm\,$};
  \end{scope}
  \scopeArrow{0.81}{\arrowDefect}
  \draw[line width=\locpc,color=\colorDefect,postaction={decorate}]
           (\locpy:\locpa) -- (\locpy:\locpA)
           node[right=-5pt,color=\colorDefect] {$\calm$};
  \draw[line width=\locpc,color=\colorTransp,postaction={decorate}]
           (\locpz:\locpa) -- (\locpz:\locpA) node[left=-5pt,color=\colorDefect] {$\II_\kappa$};
  \end{scope}
  \filldraw[color=\colorDefect] (\locpx:\locpa) circle (\locpd)
           (\locpy:\locpa) circle (\locpd) (\locpz:\locpa) circle (\locpd);
  \begin{scope}[shift={(7.9,0)}] 
  \scopeArrow{0.02}{<} \scopeArrow{0.43}{<} \scopeArrow{0.73}{<}
  \filldraw[line width=\locpc,\colorCircle,fill=white,postaction={decorate}]
           node[left=1.9cm,color=black] {$\LII_\kappa^{\sss\searrow}(\calm) ~:= $}
           (0,0) circle (\locpa)
           node[below=0.6cm,xshift=9pt,color=black]{$\scriptstyle-\kappa$};
  \end{scope} \end{scope} \end{scope}
  \scopeArrow{0.81}{\arrowDefect}
  \draw[line width=\locpc,color=\colorDefect,postaction={decorate}]
           (\locpX:\locpa) -- (\locpX:0.9*\locpA)
           node[left=-2pt,yshift=-2pt,color=\colorDefect] {$\calm$};
  \end{scope}
  \scopeArrow{0.65}{\arrowDefectB}
  \draw[line width=\locpc,color=\colorTransp,postaction={decorate}]
           (\locpY:\locpa) -- (\locpY:\locpA) node[right=-4pt,color=\colorDefect] {$\II_\kappa$};
  \draw[line width=\locpc,color=\colorDefect,postaction={decorate}]
           (\locpZ:\locpa) -- (\locpZ:\locpA) node[left=-4pt,color=\colorDefect] {$\calm$};
  \end{scope}
  \filldraw[color=\colorDefect] (\locpX:\locpa) circle (\locpd)
           (\locpY:\locpa) circle (\locpd) (\locpZ:\locpa) circle (\locpd);
  \end{scope}
  \end{tikzpicture}}
  \label{eq:circle.IIn-M-M}
  \ee
are canonically equivalent to the category $\Funle_{\Cala,\cala}(\calm,\calm)$ of 
bimodule endofunctors: we have
  \be
  \bearll
  \ZZ(\LII_\kappa^{\sss\nearrow}(\calm)) \!\!\!&
  = \cenf{ \CC\calm \,{\stackrel1\boxtimes}\, \calm \,{\stackrel \kappa\boxtimes}\, \II_\kappa
  \,{\stackrel1\boxtimes}} 
  \Nxl4 &\!\!
  \!\stackrel{\eqref{eq:equ-lex-framed-cent}}\simeq\!
  \Funle_{\Cala,\cala}(\calm,\cenf{\calm{\stackrel \kappa\boxtimes}\II_\kappa})
  \!\stackrel{\eqref{eq:M.I_n=M}}\simeq\! \Funle_{\Cala,\cala}(\calm,\calm) \qquad\text{and}
  \Nxl6
  \ZZ(\LII_\kappa^{\sss\swarrow}(\calm)) \!\!\!&
  = \cenf{ \CC\calm \,{\stackrel{-\kappa}\boxtimes}\, \CC{\II_\kappa} \,{\stackrel1\boxtimes}\,
  \calm \,{\stackrel1\boxtimes}}
  \Nxl4 &\!\!
  \stackrel{\eqref{eq:MkNopp=N-kM}}\simeq\!
  \cenf{ \CC{\II_\kappa \,{\stackrel{\kappa}\boxtimes}\,\calm} \,{\stackrel1\boxtimes}\,
  \calm \,{\stackrel1\boxtimes}} 
  \,\simeq\, \cenf{ \CC\calm \,{\stackrel1\boxtimes}\, \calm \,{\stackrel1\boxtimes}}
  \stackrel{\eqref{eq:equ-lex-framed-cent}}\simeq\! \Funle_{\Cala,\cala}(\calm,\calm) \,,
  \eear
  \label{eq:ZZIkappa=Funle-2}
  \ee
and similarly for the other two gluing circles.
The so obtained endofunctor categories have the identity functor as a distinguished object.
\end{Example}

\begin{Remark}
The pre-images of the identity functor under the equivalences \eqref{eq:ZZIkappa=Funle-2}
constitute distinguished objects in the respective gluing categories.
Using the precise form of the equivalences \eqref{eq:M.I_n=M} (see the proof of Lemma
\ref{lem:ZMA=M}), one can express these distinguished objects as the subtle compounds
  \be
  \begin{array}{lll} \dsty
  \int^{m\in\calm}\! \cc m\boti Z_{[\kappa]}(m \boti \one) \!\!&\dsty
  = \int^{m\in\calm}\!\!\! \int_{\!a\in\cala}\cc m \boti m. a \boti {}^{[\kappa-1]}a
  & \in \ZZ(\LII_\kappa^{\sss\nearrow}(\calm))
  \Nxl2
  \text{and}
  \Nxl2 \dsty
  \int^{m\in\calm}\! Z_{[\kappa]}(\cc{m} \boti \cc{\DA} )\boti m
  \!\!& \multicolumn{2}{l}{ \dsty
  = \int^{m\in\calm}\!\!\! \int_{\!a\in\cala} \cc{a. m} \boti \cc{\DA.{}^{[\kappa+3]}a} \boti m }
  \Nxl2 & \dsty
  = \int^{a\in\cala}\!\!\! \int^{m\in\calm}\! \cc{ m} \boti \cc{a} \boti {}^{[\kappa+1]}a. m 
  & \in \ZZ(\LII_\kappa^{\sss\swarrow}(\calm))
  \eear
  \ee
of ends and coends. Here the last equality follows from $\int^{a}\! a \boti
\cc{\DA\oti{}^{[\kappa+3]}a} \,{\cong} \int_{a} a \boti \cc{{}^{[\kappa+1]}a} \,{\in}$
                            $\cala \boti \CC\cala$ (compare \eqref{eq:disting-iso}).
\end{Remark}

Next we recall the so-called \emph{braided induction}. Given a left $\cala$-module $\calm$ and
a right $\cala$-mo\-du\-le $\caln$, consider, for $z\iN \calz(\cala)$, the endofunctors 
  \be
  F_{\!z}(n) := n. z \qquad\text{and}\qquad {}_{z}F(m) \,{=}\, z.m 
  \label{eq:braidedind}
  \ee
of $\caln$ and $\calm$, respectively. The half-braiding on $z$ 
endows the functors $F_{\!z}$ and ${}_{z}F$ with the structure of a 
module functor. By assigning to $z$ these functors we obtain two monoidal functors 
  \be
  F_{\!\bullet} :\quad \calz(\cala) \xrightarrow{~~} \Funle_{\Cala}(\caln,\caln)
  \qquad\text{and}\qquad
  {}_{\bullet}F :\quad \calz(\cala) \xrightarrow{~~} \Funle_{\Cala}(\calm,\calm) \,.
  \ee
These are the functors of braided induction, also termed $\alpha$-\emph{induction},
see e.g.\ \Cite{Sect.\,5.1}{ostr} for the categorical formulation used here. (Alternatively
one may use the inverse half-braiding, which gives another pair of such functors.)
 
\begin{Example} \label{exa:braidedinduction}
For $\cala$ a finite tensor category, $\calm$ a left $\cala$-module and $\caln$
a right $\cala$-module, consider the following \defectsurface s $\DD_1$ and $\DD_2$:
    \def\locpa  {2.2}    
    \def\locpA  {0.5}    
    \def\locpb  {0.25}   
    \def\locpc  {2.2pt}  
    \def\locpC  {1.8pt}  
    \def\locpd  {0.12}   
    \def\locpD  {0.09}   
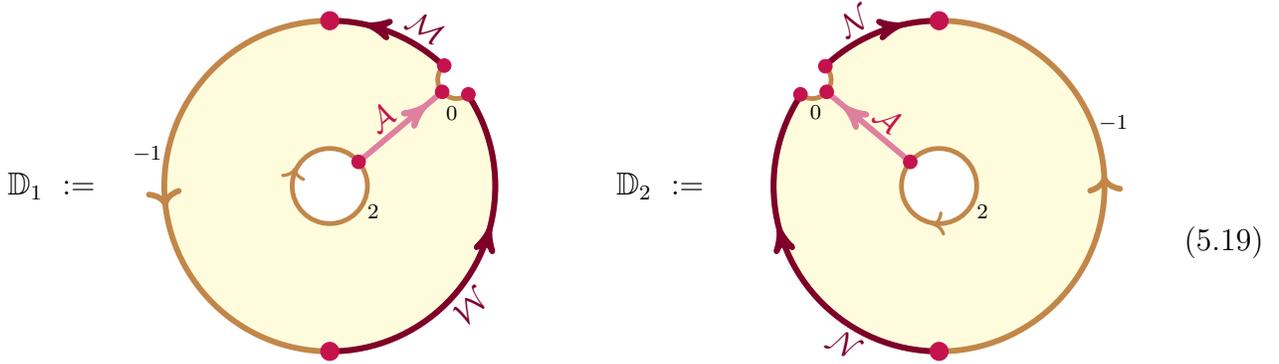
  \begin{eqnarray}
  \begin{tikzpicture}
  \scopeArrow{0.52}{>}
  \filldraw[\colorCircle,fill=\colorSigma,line width=\locpc,postaction={decorate}]
             (0,0) circle (\locpa);
  \end{scope}
  \scopeArrow{0.42}{\arrowDefect} \scopeArrow{0.93}{\arrowDefect}
  \draw[line width=\locpc,color=\colorFreebdy,postaction={decorate}]
	     (-90:\locpa) arc (-90:90:\locpa);
  \end{scope} \end{scope}
  \scopeArrow{0.76}{\arrowDefect}
  \draw[line width=\locpc,color=\colorTransp,postaction={decorate}] (0,0) --
             node[above,midway,sloped,yshift=-2pt,color=\colorDefect] {$\cala$} (40:\locpa);
  \end{scope}
  \filldraw[\colorCircle,fill=white,line width=\locpC] (0,0) circle (\locpA);
  \fill[white] (40:\locpa) circle (\locpb);
  \draw[\colorCircle,line width=\locpC] (40:\locpa)++(120:\locpb) arc (120:310:\locpb);
  \scopeArrow{0.47}{<}
  \draw[\colorCircle,line width=0.7*\locpC,postaction={decorate}] (0,0) circle (\locpA);
  \end{scope}
  \filldraw[color=\colorDefect] (90:\locpa) circle (\locpd) (270:\locpa) circle (\locpd)
             (40:\locpA) circle (\locpD) (40:\locpa-\locpb) circle (\locpD)
             (40:\locpa)++(131:\locpb) circle (\locpD)
             (40:\locpa)++(-51:\locpb) circle (\locpD);
  \node[rotate=-131,color=\colorFreebdy] at (-40:1.12*\locpa) {$\calm$};
  \node[rotate=-32,color=\colorFreebdy] at (59:1.12*\locpa) {$\calm$};
  \node at (170:1.12*\locpa) {$\scriptstyle {-}1$};
  \node at (330:1.33*\locpA) {$\scriptstyle 2$};
  \node at (31:0.86*\locpa) {$\scriptstyle 0$};
  \node at (-\locpa-1.5,0) {$\DD_1 ~:=$};
  \begin{scope}[shift={(3*\locpa+1.5,0)}]
  \scopeArrow{0.01}{>}
  \filldraw[\colorCircle,fill=\colorSigma,line width=\locpc,postaction={decorate}]
             (0,0) circle (\locpa);
  \end{scope}
  \scopeArrow{0.42}{\arrowDefect} \scopeArrow{0.93}{\arrowDefect}
  \draw[line width=\locpc,color=\colorFreebdy,postaction={decorate}]
	     (-90:\locpa) arc (270:90:\locpa);
  \end{scope} \end{scope}
  \scopeArrow{0.75}{\arrowDefect}
  \draw[line width=\locpc,color=\colorTransp,postaction={decorate}] (0,0) --
             node[above,midway,sloped,yshift=-2pt,color=\colorDefect] {$\cala$} (140:\locpa);
  \end{scope}
  \filldraw[\colorCircle,fill=white,line width=\locpC] (0,0) circle (\locpA);
  \fill[white] (140:\locpa) circle (\locpb);
  \draw[\colorCircle,line width=\locpC] (140:\locpa)++(60:\locpb) arc (60:-130:\locpb);
  \scopeArrow{0.77}{<}
  \draw[\colorCircle,line width=0.7*\locpC,postaction={decorate}] (0,0) circle (\locpA);
  \end{scope}
  \filldraw[color=\colorDefect] (90:\locpa) circle (\locpd) (270:\locpa) circle (\locpd)
             (140:\locpA) circle (\locpD) (140:\locpa-\locpb) circle (\locpD)
             (140:\locpa)++(46:\locpb) circle (\locpD)
             (140:\locpa)++(231:\locpb) circle (\locpD);
  \node[rotate=31,color=\colorFreebdy] at (117:1.12*\locpa) {$\caln$};
  \node[rotate=-212,color=\colorFreebdy] at (239:1.11*\locpa) {$\caln$};
  \node at (20:1.12*\locpa) {$\scriptstyle {-}1$};
  \node at (330:1.33*\locpA) {$\scriptstyle 2$};
  \node at (149:0.87*\locpa) {$\scriptstyle 0$};
  \node at (-\locpa-1.51,0) {$\DD_2 ~:=$};
  \end{scope}
  \end{tikzpicture}
  \nonumber \\[-3.1em] ~
  \label{eq:picbraidedinduction}
  \\[-0.4em] ~ \nonumber
  \end{eqnarray}
Here, and in all pictures below, defect lines that are transparently labeled 
are drawn in a lighter color.
Both surfaces have the same inner boundary circle $\LL$; we denote the outer boundary circle
of $\DD_{i}$ by $\LL_i$, for $i\,{=}\,1,2$. We have $\ZZ(\LL) \,{=}\, \calz(\cala)$ and
$\ZZ(\LL_1) \,{=}\, \cc{\calm} \,\Boti{-1}\, {\calm} \,{\stackrel{\eqref{eq:MbotikappaN=Lex}}\simeq}$
           $\cc{ \Funle_{\Cala}(\calm,\calm) }$, while 
$\ZZ(\LL_2) \,{=}\, {\caln} \,\Boti{-1}\, \cc{\caln} \,{\stackrel{\eqref{eq:MbotikappaN=Lex}}\simeq}
\cc{ \Funle_{\Cala}(\caln,\caln) }$.
The block functors for the disks $\DD_{1}$ and $\DD_{2}$ are given by 
  \be
  \bearl
  \Zfine(\DD_{1})(\cc{G} \boti z)=\Funle_{\cala}(G, F_{z}) \qquad\text{and}
  \Nxl8
  \Zfine(\DD_{2})(\cc{G} \boti z)=\Funle_{\cala}(G, {}_{z}F) \,,
  \eear
  \label{eq:blocks-D1D2}
  \ee
respectively, for $G \iN \Funle_{\Cala,\cala}(\cala,\cala) $ and $z \iN \calz(\cala)$, 
and thus, by the Yoneda lemma, represent the braided inductions.
In the case of $\DD_{1}$ this is seen as follows (the computation of $\ZZ(\DD_{2})$ 
is analogous). Consider the object $\widetilde z$ in $ \cc{\calm} \,\Boti{-1}\, {\calm}$ that
for $z \iN \ZZ(\LL_1)$ results from contracting the defect line that is connected to
$\LL_1$. Using the formulas given in Example \ref{ex:filldisk} we obtain
  \be
  \widetilde z \,= \int^{m\in\calm}\!\!\! \int^{a\in\cala}\! \Hom_{\cala}(a, z)
  \otimes \cc{m} \boti a.m  \, \cong \int^{m\in\calm}\! \cc{m} \boti z.m \,.
  \label{eq:result-contract}
  \ee
Under the Eilenberg--Watts correspondence this object gives the bimodule functor $F_z$. 
The claimed result for $\Zfine(\DD_{1})$ then follows from Corollary \ref{cor:modulenattrafo}.
\end{Example}

Combining the notion of transparent defect with the one of fillable circles and intervals
allows us to introduce a particular type of \defectsurface s. Topologically these specific
surfaces are disks with holes, but since the holes are fillable and their gluing categories
contain distinguished objects, for brevity we abuse terminology and refer to these
\defectsurface s just as ``disks''. 

\begin{Definition} \label{def:filldisk}
Let $\XX$ be a \defectsurface\ such that the underlying surface $X$ is a disk and
$\partial\XX$ contains at most one free boundary segment. 
\Itemizeiii

\item[{\rm (i)}]
A \emph{\filldisk\ $(\DD_{\XX},\delta_{\mathrm{tr}},\partial_{\mathrm{outer}})$
of type $\XX$} is a \defectsurface\ $\DD_{\XX}$ with having following structure:
~\\[-1.65em]\def\leftmarginii{1.24em}
\begin{itemize} \addtolength\itemsep{-0.3em}

\item[{\rm 1.}]
Every boundary circle of $\DD_{\XX}$ is oriented as induced by the orientation of the surface. 

\item[{\rm 2.}]
$\DD_{\XX}$ contains a set $\delta_\text{tr}$ of
distinguished transparent defect lines on $\DD_{\XX}$, as well as one distinguished 
gluing boundary component $\partial_\text{outer} \,{\subseteq}\, \partialg\DD_{\XX}$,
which we call the \emph{outer boundary} of $\DD_{\XX}$.

\item[{\rm 3.}]
Deletion of all transparent defect lines that belong to $\delta_\text{tr}$ results 
in a \defectsurface\ for which all gluing circles and gluing intervals, except for
$\partial_\text{outer}$, are fillable by a disk in the sense of Definition 
$\ref{def:filling}$, and the so obtained filling is isomorphic to the \defectsurface\ $\XX$.
\end{itemize}

\item[{\rm (ii)}]
We depict a \filldisk\ $\DD_{\XX}$ as a subset of the plane, with the non-fillable
boundary circle $($or interval, in case $\XX$ has  a free boundary$)$ forming 
the outer boundary, while the fillable gluing circles and intervals are referred to as
\emph{inner boundary circles} or \emph{intervals}. 
\end{itemize}
\end{Definition}

It follows that every defect line on $\DD_\XX$ is either a distinguished transparent defect or 
corresponds to a defect of $\XX$, possibly interrupted by fillable gluing boundaries.
Note that a \filldisk\ is in general not topologically a disk. 
As an illustration, the following is an example of a \filldisk\ (the non-fillable gluing
boundary of $\DD_{\XX}$ is the interval labeled as $\LL$, various other labels are omitted):
    \def\locpa  {2.7}    
    \def\locpA  {0.34}   
    \def\locpb  {2.0}    
    \def\locpc  {2.2pt}  
    \def\locpC  {1.8pt}  
    \def\locpd  {0.12}   
    \def\locpD  {0.08}   
    \def\locpf  {1.1pt}  
    \def\locpm  {225}    
    \def\locpM  {75}     
    \def\locpx  {20}     
    \def\locpX  {-70}    
    \def\locpz  {36.5}   
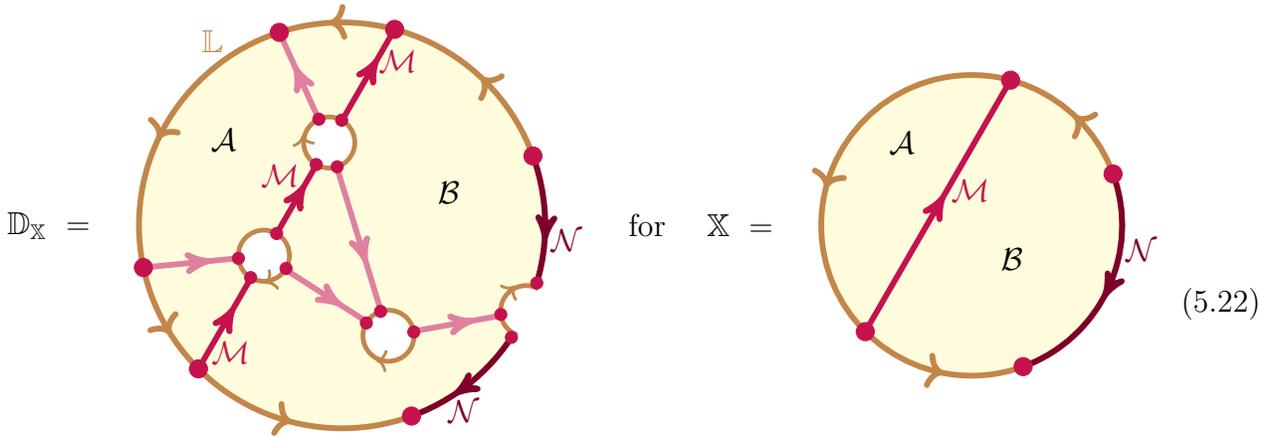
  \begin{eqnarray}
  \begin{tikzpicture}
  \coordinate (incirc1) at ({\locpa*cos(\locpm)/3+2*\locpa*cos(\locpM)/3},
                            {\locpa*sin(\locpm)/3+2*\locpa*sin(\locpM)/3});
  \coordinate (incirc2) at ({2*\locpa*cos(\locpm)/3+\locpa*cos(\locpM)/3},
                            {2*\locpa*sin(\locpm)/3+\locpa*sin(\locpM)/3});
  \coordinate (incirc3) at (0.96*\locpX:0.59*\locpa);
  \coordinate (oucirc1) at (\locpM+33:\locpa);
  \coordinate (oucirc2) at (\locpm-33:\locpa);
  \coordinate (bdcirc1) at (0.5*\locpx+0.5*\locpX:\locpa);
  \node at (-\locpa-1.2,0) {$\DD_\XX ~=$};
  \scopeArrow{0.13}{>} \scopeArrow{0.26}{>} \scopeArrow{0.43}{>} \scopeArrow{0.59}{>}
  \scopeArrow{0.71}{>}
  \filldraw[\colorCircle,fill=\colorSigma,line width=\locpc,postaction={decorate}] (0,0) circle (\locpa);
  \end{scope} \end{scope} \end{scope} \end{scope} \end{scope}
  \scopeArrow{0.18}{\arrowDefect} \scopeArrow{0.54}{\arrowDefect} \scopeArrow{0.91}{\arrowDefect}
  \draw[line width=\locpc,color=\colorDefect,postaction={decorate}]
                  (\locpm:\locpa) node[right,yshift=6pt] {$\calm$}
		  -- node[left=-6pt,yshift=9pt] {$\calm$}
		  (\locpM:\locpa) node[right=-11pt,yshift=-12pt] {$\calm$};
  \end{scope} \end{scope} \end{scope}
  \scopeArrow{0.27}{\arrowDefect} \scopeArrow{0.85}{\arrowDefect}
  \draw[line width=\locpc,color=\colorFreebdy,postaction={decorate}]
		  (\locpx:\locpa) node[right=2pt,yshift=-32pt] {$\caln$} 
                  arc (\locpx:\locpX:\locpa) node[right=9pt,yshift=2pt] {$\caln$};
  \end{scope} \end{scope}
  \scopeArrow{0.56}{\arrowDefect}
  \draw[line width=\locpc,color=\colorTransp,postaction={decorate}] (oucirc2) -- (incirc2);
  \end{scope}
  \scopeArrow{0.62}{\arrowDefect}
  \draw[line width=\locpc,color=\colorTransp,postaction={decorate}] (incirc1) -- (incirc3);
  \draw[line width=\locpc,color=\colorTransp,postaction={decorate}] (incirc2) -- (incirc3);
  \draw[line width=\locpc,color=\colorTransp,postaction={decorate}] (incirc3) -- (bdcirc1);
  \end{scope}
  \scopeArrow{0.69}{\arrowDefect}
  \draw[line width=\locpc,color=\colorTransp,postaction={decorate}] (incirc1) -- (oucirc1);
  \end{scope}
  \filldraw[\colorCircle,fill=white,line width=\locpC]
                  (incirc1) circle (\locpA) (incirc2) circle (\locpA) (incirc3) circle (\locpA)
                  (bdcirc1) circle (\locpA);
  \scopeArrow{0.38}{<}
  \draw[\colorCircle,line width=\locpf,postaction={decorate}] (bdcirc1) circle (\locpA);
  \end{scope}
  \scopeArrow{0.52}{<}
  \draw[\colorCircle,line width=\locpf,postaction={decorate}] (incirc1) circle (\locpA);
  \end{scope}
  \scopeArrow{0.82}{<}
  \draw[\colorCircle,line width=\locpf,postaction={decorate}] (incirc2) circle (\locpA);
  \end{scope}
  \scopeArrow{0.74}{<}
  \draw[\colorCircle,line width=\locpf,postaction={decorate}] (incirc3) circle (\locpA);
  \end{scope}
  \fill[white] (\locpX+\locpz:1.02*\locpa) -- (\locpx-\locpz:1.02*\locpa) --
               (\locpx-\locpz:1.2*\locpa) -- (\locpX+\locpz:1.2*\locpa) -- cycle;
  \filldraw[color=\colorDefect] (\locpm:\locpa) circle (\locpd) (\locpM:\locpa) circle (\locpd)
               (\locpx:\locpa) circle (\locpd) (\locpX:\locpa) circle (\locpd)
               (oucirc1) circle (\locpd) (oucirc2) circle (\locpd)
               (\locpx-\locpz:\locpa) circle (\locpD) (\locpX+\locpz:\locpa) circle (\locpD)
               (bdcirc1)+(187.5:\locpA) circle (\locpD) 
               (incirc1)+(61:\locpA) circle (\locpD) (incirc1)+(114:\locpA) circle (\locpD)
               (incirc1)+(239.5:\locpA) circle (\locpD) (incirc1)+(288:\locpA) circle (\locpD)
               (incirc2)+(-30:\locpA) circle (\locpD) (incirc2)+(60:\locpA) circle (\locpD)
               (incirc2)+(186:\locpA) circle (\locpD) (incirc2)+(239:\locpA) circle (\locpD)
               (incirc3)+(8.5:\locpA) circle (\locpD) (incirc3)+(108:\locpA) circle (\locpD)
               (incirc3)+(150:\locpA) circle (\locpD);
  \node[color=\colorCircle] at (125:1.11*\locpa) {$\LL$};
  \node at (144:0.72*\locpa) {$\cala$};
  \node at (17:0.55*\locpa) {$\calb$};
  \begin{scope}[shift={(3.1*\locpa,0)}]
  \node at (-\locpa-0.9,0) {for $\quad \XX ~=$};
  \scopeArrow{0.13}{>} \scopeArrow{0.46}{>} \scopeArrow{0.72}{>}
  \filldraw[\colorCircle,fill=\colorSigma,line width=\locpc,postaction={decorate}] (0,0) circle (\locpb);
  \end{scope} \end{scope} \end{scope}
  \scopeArrow{0.54}{\arrowDefect}
  \draw[line width=\locpc,color=\colorDefect,postaction={decorate}]
                  (\locpm:\locpb) -- node[right,yshift=5pt] {$\calm$} (\locpM:\locpb);
  \draw[line width=\locpc,color=\colorFreebdy,postaction={decorate}]
                  (\locpx:\locpb) node[right,yshift=-30pt] {$\caln$} arc (\locpx:\locpX:\locpb);
  \end{scope}
  \filldraw[color=\colorDefect] (\locpm:\locpb) circle (\locpd) (\locpM:\locpb) circle (\locpd)
                                (\locpx:\locpb) circle (\locpd) (\locpX:\locpb) circle (\locpd);
  \node at (130:0.72*\locpb) {$\cala$};
  \node at (320:0.35*\locpb) {$\calb$};
  \end{scope}
  \end{tikzpicture}
  \nonumber \\[-3.5em] ~
  \label{eq:exampletranspdisk}
  \\[0.4em] ~ \nonumber
  \end{eqnarray}
We also consider \filldisk s for which, while they are proper \defectsurface s themselves,
the corresponding
surface $\XX$ is improper in the sense that it does not contain any defect lines or free boundaries.
In this case we call the disk a \emph{\transdisk} and denote it just by $\DDt$; the following
is an example of a \transdisk:
    \def\locpa  {2.7}    
    \def\locpA  {0.34}   
    \def\locpc  {2.2pt}  
    \def\locpC  {1.8pt}  
    \def\locpd  {0.12}   
    \def\locpD  {0.08}   
    \def\locpf  {1.1pt}  
  \begin{eqnarray}
  \begin{tikzpicture}
  \coordinate (incirc1) at (-55:0.66*\locpa);
  \coordinate (incirc2) at (3:0.71*\locpa);
  \coordinate (incirc3) at (65:0.55*\locpa);
  \coordinate (incirc4) at (145:0.66*\locpa);
  \coordinate (incirc5) at (225:0.35*\locpa);
  \coordinate (oucirc1) at (40:\locpa);
  \coordinate (oucirc2) at (90:\locpa);
  \coordinate (oucirc3) at (180:\locpa);
  \coordinate (oucirc4) at (225:\locpa);
  \coordinate (oucirc5) at (257:\locpa);
  \node at (-\locpa-1.1,0) {$\DDt ~=$};
  \scopeArrow{0.19}{>} \scopeArrow{0.38}{>} \scopeArrow{0.57}{>} \scopeArrow{0.68}{>}
  \scopeArrow{0.94}{>}
  \filldraw[\colorCircle,fill=\colorSigma,line width=\locpc,postaction={decorate}] (0,0) circle (\locpa);
  \end{scope} \end{scope} \end{scope} \end{scope} \end{scope}
  \scopeArrow{0.49}{\arrowDefectB}
  \draw[line width=\locpc,color=\colorTransp,postaction={decorate}] (oucirc4) -- (incirc5);
  \draw[line width=\locpc,color=\colorTransp,postaction={decorate}] (oucirc5) -- (incirc5);
  \end{scope}
  \scopeArrow{0.56}{\arrowDefect}
  \draw[line width=\locpc,color=\colorTransp,postaction={decorate}] (oucirc1) -- (incirc3);
  \draw[line width=\locpc,color=\colorTransp,postaction={decorate}] (oucirc2) -- (incirc3);
  \draw[line width=\locpc,color=\colorTransp,postaction={decorate}] (oucirc3) -- (incirc5);
  \end{scope}
  \scopeArrow{0.63}{\arrowDefect}
  \draw[line width=\locpc,color=\colorTransp,postaction={decorate}] (incirc4) -- (incirc3);
  \draw[line width=\locpc,color=\colorTransp,postaction={decorate}] (incirc3) -- (incirc2);
  \draw[line width=\locpc,color=\colorTransp,postaction={decorate}] (incirc2) -- (incirc1);
  \draw[line width=\locpc,color=\colorTransp,postaction={decorate}] (incirc3) -- (incirc5);
  \draw[line width=\locpc,color=\colorTransp,postaction={decorate}] (incirc5) -- (incirc1);
  \end{scope}
  \filldraw[\colorCircle,fill=white,line width=\locpC]
           (incirc1) circle (\locpA) (incirc2) circle (\locpA) (incirc3) circle (\locpA)
           (incirc4) circle (\locpA) (incirc5) circle (\locpA);
  \scopeArrow{0.82}{<}
  \draw[\colorCircle,line width=\locpf,postaction={decorate}] (incirc1) circle (\locpA);
  \draw[\colorCircle,line width=\locpf,postaction={decorate}] (incirc2) circle (\locpA);
  \draw[\colorCircle,line width=\locpf,postaction={decorate}] (incirc3) circle (\locpA);
  \draw[\colorCircle,line width=\locpf,postaction={decorate}] (incirc4) circle (\locpA);
  \end{scope}
  \scopeArrow{0.09}{<}
  \draw[\colorCircle,line width=\locpf,postaction={decorate}] (incirc5) circle (\locpA);
  \end{scope}
  \filldraw[color=\colorDefect] (oucirc1) circle (\locpd) (oucirc2) circle (\locpd)
           (oucirc3) circle (\locpd) (oucirc4) circle (\locpd) (oucirc5) circle (\locpd)
           (incirc1)+(60:\locpA) circle (\locpD) (incirc1)+(156:\locpA) circle (\locpD)
           (incirc2)+(136:\locpA) circle (\locpD) (incirc2)+(239.3:\locpA) circle (\locpD)
           (incirc3)+(14:\locpA) circle (\locpD) (incirc3)+(115:\locpA) circle (\locpD)
           (incirc3)+(188:\locpA) circle (\locpD) (incirc3)+(238:\locpA) circle (\locpD)
           (incirc3)+(316:\locpA) circle (\locpD) (incirc4)+(8:\locpA) circle (\locpD)
           (incirc5)+(57:\locpA) circle (\locpD) (incirc5)+(161:\locpA) circle (\locpD)
           (incirc5)+(225:\locpA) circle (\locpD) (incirc5)+(272:\locpA) circle (\locpD)
           (incirc5)+(337:\locpA) circle (\locpD);
  \begin{scope}[shift={(3.1*\locpa,0)}]
  \node at (-\locpa-0.9,0) {for $\quad \XX ~=$};
  \scopeArrow{0.46}{>}
  \filldraw[\colorCircle,fill=\colorSigma,line width=\locpc,postaction={decorate}]
           (0,0) circle (\locpb) node[right=1.8em] {~};   
  \end{scope}
  \end{scope}
  \end{tikzpicture}
  \nonumber \\[-3.2em] ~
  \label{eq:exampleAtranspdisk}
  \\[0.3em] ~ \nonumber
  \end{eqnarray}

By definition, all defect lines of a \transdisk\ $\DDt$ of type $\cala$ are labeled by
the transparent defect label $\II_0 \eq \cala$ for one and the same tensor category $\cala$.
The framing of a \filldisk\ $\DD_{\XX}$ has the following property: since
by definition every inner gluing segment is fillable by a disk, the indices on the outer boundary
$\partial_\text{outer}\DD_{\XX}$ must add up to the value 2 
(thus $\partial_\text{outer}\DD_{\XX}$ is not fillable). Moreover, the label of every defect
line of $\XX$ appears an even number of times as a label of a defect point on $\partial \DD_{\XX}$.
(The latter criteria are, however, not sufficient for the existence of a \filldisk\ with a given 
boundary: If the cyclically ordered string of bimodules for the boundary $\partial\XX$ contains
a string $(\calm, \caln, \cc{\calm}, \cc{\caln})$ 
with generic bimodules $\calm$ and $\caln$, then there does not exist a corresponding \filldisk.)

\medskip

Let us also mention that the gluing categories for fillable circles and intervals can often
conveniently be described as functor categories, whereby one can in particular identify
certain distinguished objects in such gluing categories. Also, in Lemma \ref{Lemma:rho-part1}
in the appendix, we exhibit distinguished objects for all fillable circles and intervals.  

\begin{Example}\label{exa:JJ}
\Itemizeiii

\item[(i)]
The gluing categories in \eqref{eq:circle.IIn-M-M} are fillable if and only if $\kappa \,{=}\, 0$.

\item[(ii)]
Consider, for $\cala$ a finite tensor category, the \defectonemanifold
    \def\locpa  {0.73}   
    \def\locpA  {1.9}    
    \def\locpb  {0.83}   
    \def\locpc  {2.2pt}  
    \def\locpd  {0.12}   
  \be 
  \raisebox{-4.3em}{\begin{tikzpicture}
  \scopeArrow{0.47}{<} \scopeArrow{0.97}{<}
  \filldraw[line width=\locpc,\colorCircle,fill=white,postaction={decorate}]
       node[left=1.88cm,color=black]{$\JJ_\kappa ~:=$} (0,0) circle (\locpa)
       node[left=0.64cm,yshift=-6pt,color=black]{$\scriptstyle 1+ \kappa$}
       node[right=0.6cm,yshift=7pt,color=black]{$\scriptstyle 1-\kappa$};
  \end{scope} \end{scope}
  \scopeArrow{0.81}{\arrowDefect}
  \draw[line width=\locpc,color=\colorTransp,postaction={decorate}] (0,\locpa) -- +(0,\locpb)
       node[right=-3pt,yshift=1pt,color=\colorDefect] {$\II_{0}$};
  \end{scope}
  \scopeArrow{0.71}{\arrowDefectB}
  \draw[line width=\locpc,color=\colorTransp,postaction={decorate}] (0,-\locpa) -- +(0,-\locpb)
       node[right=-4pt,yshift=-1pt,color=\colorDefect] {$\,\II_{0}$};
  \end{scope}
  \filldraw[color=\colorDefect] (0,\locpa) circle (\locpd) (0,-\locpa) circle (\locpd);
  \end{tikzpicture}}
  \label{eq:transp-Serre-def-points}
  \ee
with $\kappa \,{\in}\, 2\Z$.
It follows from \eqref{eq:M.II0.N=M.N} that there is a distinguished equivalence 
  \be
  \ZZ(\JJ_\kappa) =
  \cc{\II_{0}} \,{\stackrel{1+\kappa}{\boti}} \II_{0} {\stackrel{1-\kappa}{\boti}}
  \stackrel{\eqref{eq:M.II0.N=M.N}} \simeq
  \cc{\II_{0}} {\stackrel{1}{\boti}} {}^{[\kappa]}_{}\II_{0}^{[-\kappa]} {\stackrel{1}{\boti}}
  \stackrel{\eqref{eq:kappaII-kappa=II}} \simeq
  \cc{\II_{0}} {\stackrel{1}{\boti}} \II_{0} {\stackrel{1}{\boti}}
  \stackrel{\eqref{eq:MbotikappaN=Lex}}\simeq \Funle_{\Cala,\cala}(\II_0,\II_0) \,.~
  \label{eq:ZZ(JJ)=...}
  \ee
Thus in $\ZZ(\JJ_\kappa)$ there is a distinguished object, given by the pre-image of the identity
functor in $\Funle_{\Cala,\cala}(\II_0,\II_0)$ under the equivalence \eqref{eq:ZZ(JJ)=...}.
It is given explicitly by $\int^{a\in\cala} \cc{a}\boti{}^{[\kappa]}a \,{\in}\, \ZZ(\JJ_\kappa)$,
with balancing determined by the one of $\int^{a\in\cala}\cc{a} \boti a$.

\item[(iii)]
Similarly, for the \defectonemanifold s that coincide with $\JJ_\kappa$ as manifolds, but
have general framings, i.e.\ with general indices $\kappa$ and $\kappa'$, the gluing category
is canonically equivalent to $\Funle_{\Cala,\cala}(\II_0^\kappa,\II_0^{\kappa'})$.
If $\kappa \,{=}\, \kappa'$, then this functor category is equivalent to the
Drinfeld center and is thus monoidal; if $\kappa \,{-}\, \kappa' \,{=}\, 2$, then it is
equivalent to the category of objects $x \iN \cala$ together with coherent natural isomorphisms
$a \oti x \,{\cong}\, x \oti a^{\vee\vee}$ for all $a \iN \cala$, which is in general not
monoidal (e.g., the monoidal unit of $\cala$ might not have the structure of an object 
in this category). The categories for all other cases are equivalent to one obtained for
the latter two cases, determined by $\kappa \,{-}\, \kappa'$ mod 4, using the
distinguished invertible object in $\cala$ as in Equation \eqref{eq:DA:kappa+4}.
 \\[.2em]
We take this observation as an opportunity to remark that
in \cite{doSs3} the framing on a circle (without defect points) is described with the help
of a corona instead of an index. For instance, the three framed manifolds shown in Table 3 of
\cite{doSs3} correspond to $\kappa{+}\kappa' \,{=}\, 2$, $0$ and $-2$, respectively.
\end{itemize}
\end{Example}


\subsection{Refinement of defect surfaces}\label{sec:refin-defect-surf}

In Section \ref{sec:pre-blocks} we have assigned a pre-block functor $\Zpre(\Sigma)$ to any
arbitrary \defectsurface\ $\Sigma$. In contrast, the block functor $\ZZ(\Sigma)$ could so far
be defined only for surfaces that are \emph{fine} in the sense of Definition \ref{def:fine}(ii).
To complete the definition of the modular 2-functor, this restriction must be removed.
As a first step towards this end we are going to construct, for any \defectsurface\ $\Sigma$, a
family of fine \defectsurface s $\Sigma\refi$, to be called \emph{fine refinements} of $\Sigma$.
In the present subsection we introduce the notion of refinement and show that fine refinements
exist for any surface.
Later on we will use the family of fine refinements of a given surface $\Sigma$ to define the
functor $\ZZ(\Sigma)$; in case $\Sigma$ is already fine itself, this must coincide with the
functor given by Definition \ref{Definition:block-func}. 

We need a precise notion of refinement of \defectsurface s:

\begin{Definition} \label{def:refinement}
Let $\Sigma$ and $\Sigma'$ be two \defectsurface s, both not necessarily fine.
A \emph{refinement} from $\Sigma$ to $\Sigma'$ is a triple $(\Sigma;\Sigma';\varphi)$, 
with distinguished sets $\mathbf I'$ of transparent defect lines, $\mathbf S'$  of gluing 
circles and $\mathbf L'$ of gluing intervals in $\Sigma'$, such that the following holds:
\Itemizeiii

\item[{\rm 1.}]
Both end points of each of the defect lines 
in $\mathbf I'$ are contained in $\mathbf S' \,{\cup}\, \mathbf L'$.
\item[{\rm 2.}]
When removing all defect lines in $\mathbf I'$ from $\Sigma'$, each of the modified gluing 
circles $\SS_j''$ that result from some $\SS_j' \iN \mathbf S'$ and modified gluing intervals 
$\LL_k''$  that result from some $\LL_j' \iN \mathbf L'$ is fillable.
\item[{\rm 3.}]
The resulting filling $\Sigma''$ of $\Sigma'$ is isomorphic as a \defectsurface\ to 
$\Sigma$, with isomorphism $\varphi\colon \Sigma \,{\to}\, \Sigma''$.
\end{itemize}
 
\noindent
We denote the subset of the gluing boundary $\partialg\Sigma\refi$ of the surface
$\Sigma\refi$ consisting of those gluing circles and gluing intervals that are not inherited
from the gluing boundary $\partialg\Sigma$ of $\Sigma$ by $\partialr\Sigma\refi$. 
We usually suppress the isomorphism $\varphi$ from the notation and denote the refinement just
by $(\Sigma;\Sigma')$; if $(\Sigma;\Sigma')$ is a refinement, then we also say that $\Sigma$
can be \emph{refined} to $\Sigma'$, and that $\Sigma'$ \emph{refines} $\Sigma$.
\end{Definition}

The definition implies that each defect line $\delta$ of $\Sigma$ corresponds either to a single
defect line of $\Sigma'$ with the same label as $\delta$, or else splits into several defect lines 
of $\Sigma'$ that all carry the same label as $\delta$ and which are interrupted 
by fillable gluing circles that are not present in $\Sigma$.
Moreover, the gluing boundaries of $\Sigma$ correspond to identical gluing boundaries of $\Sigma'$.
As an illustration, the following picture shows the refinement $(\DD;\DD')$ of a 
one-holed disk $\DD$ for which $\DD'$ has three additional gluing circles:
    \def\locpa  {3.1}    
    \def\locpA  {0.34}   
    \def\locpc  {2.2pt}  
    \def\locpC  {1.8pt}  
    \def\locpd  {0.12}   
    \def\locpD  {0.08}   
    \def\locpf  {1.1pt}  
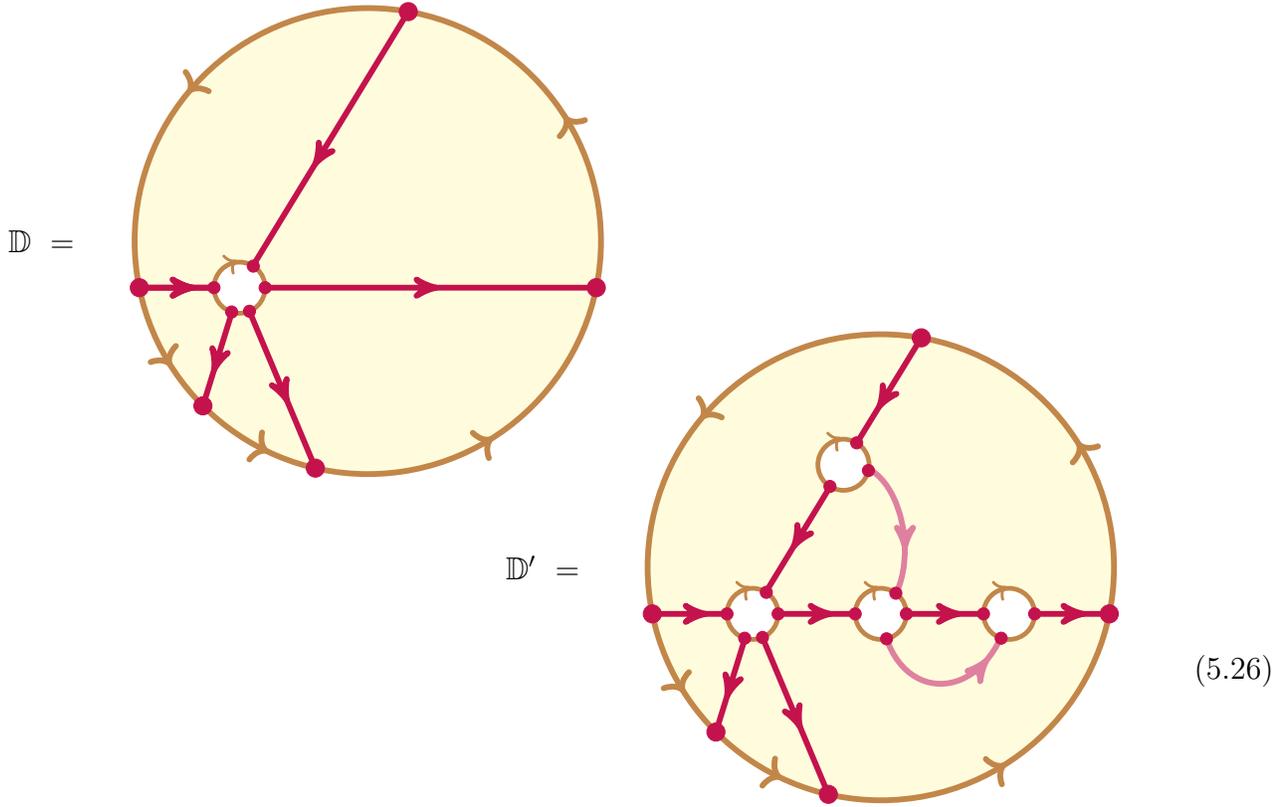
  $$
  \begin{tikzpicture}
  \coordinate (oucirc1) at (80:\locpa);
  \coordinate (oucirc3) at ({-sqrt(0.96)*\locpa},{-0.2*\locpa});
  \coordinate (oucirc4) at (225:\locpa);
  \coordinate (oucirc5) at (257:\locpa);
  \coordinate (oucirc6) at ({sqrt(0.96)*\locpa},{-0.2*\locpa});
  \coordinate (incirc5) at ({-0.55*\locpa},{-0.2*\locpa});
  \node at (-1.4*\locpa,0) {$\DD ~=$};
  \scopeArrow{0.09}{>} \scopeArrow{0.39}{>} \scopeArrow{0.59}{>} \scopeArrow{0.68}{>}
  \scopeArrow{0.84}{>}
  \filldraw[\colorCircle,fill=\colorSigma,line width=\locpc,postaction={decorate}] (0,0) circle (\locpa);
  \end{scope} \end{scope} \end{scope} \end{scope} \end{scope}
  \scopeArrow{0.49}{\arrowDefectB}
  \draw[line width=\locpc,color=\colorDefect,postaction={decorate}] (oucirc4) -- (incirc5);
  \draw[line width=\locpc,color=\colorDefect,postaction={decorate}] (oucirc5) -- (incirc5);
  \end{scope}
  \scopeArrow{0.56}{\arrowDefect}
  \draw[line width=\locpc,color=\colorDefect,postaction={decorate}] (oucirc3) -- (incirc5);
  \end{scope}
  \scopeArrow{0.56}{\arrowDefect}
  \draw[line width=\locpc,color=\colorDefect,postaction={decorate}] (oucirc1) -- (incirc5);
  \draw[line width=\locpc,color=\colorDefect,postaction={decorate}] (oucirc3) -- (incirc5);
  \draw[line width=\locpc,color=\colorDefect,postaction={decorate}] (incirc5) -- (oucirc6);
  \end{scope}
  \filldraw[\colorCircle,fill=white,line width=\locpC] (incirc5) circle (\locpA);
  \scopeArrow{0.33}{<}
  \draw[\colorCircle,line width=\locpf,postaction={decorate}] (incirc5) circle (\locpA);
  \end{scope}
  \filldraw[color=\colorDefect] (oucirc1) circle (\locpd) (oucirc3) circle (\locpd)
           (oucirc4) circle (\locpd) (oucirc5) circle (\locpd) (oucirc6) circle (\locpd)
           (incirc5)+(0:\locpA) circle (\locpD)
           (incirc5)+(57.5:\locpA) circle (\locpD) (incirc5)+(180:\locpA) circle (\locpD)
           (incirc5)+(252:\locpA) circle (\locpD) (incirc5)+(292.8:\locpA) circle (\locpD);
  \node at (2*\locpa+2.0,0) {~};
  \end{tikzpicture}
  $$
  \be
  \raisebox{-7.1em}{\begin{tikzpicture}
  \begin{scope}[shift={(2.2*\locpa,-1.4*\locpa)}]
  \coordinate (oucirc1) at (80:\locpa);
  \coordinate (oucirc3) at ({-sqrt(0.96)*\locpa},{-0.2*\locpa});
  \coordinate (oucirc4) at (225:\locpa);
  \coordinate (oucirc5) at (257:\locpa);
  \coordinate (oucirc6) at ({sqrt(0.96)*\locpa},{-0.2*\locpa});
  \coordinate (incirc1) at (0,{-0.2*\locpa});
  \coordinate (incirc2) at ({0.55*\locpa},{-0.2*\locpa});
  \coordinate (incirc3) at ({-0.16*\locpa},{0.44*\locpa});
  \coordinate (incirc5) at ({-0.55*\locpa},{-0.2*\locpa});
  \node at (-1.45*\locpa,0) {$\DD' ~=$};
  \node at (-2.4*\locpa,0) {~};
  \scopeArrow{0.09}{>} \scopeArrow{0.39}{>} \scopeArrow{0.59}{>} \scopeArrow{0.68}{>}
  \scopeArrow{0.84}{>}
  \filldraw[\colorCircle,fill=\colorSigma,line width=\locpc,postaction={decorate}] (0,0) circle (\locpa);
  \end{scope} \end{scope} \end{scope} \end{scope} \end{scope}
  \scopeArrow{0.49}{\arrowDefectB}
  \draw[line width=\locpc,color=\colorDefect,postaction={decorate}] (oucirc4) -- (incirc5);
  \draw[line width=\locpc,color=\colorDefect,postaction={decorate}] (oucirc5) -- (incirc5);
  \end{scope}
  \scopeArrow{0.56}{\arrowDefect}
  \draw[line width=\locpc,color=\colorDefect,postaction={decorate}] (oucirc1) -- (incirc3);
  \draw[line width=\locpc,color=\colorDefect,postaction={decorate}] (oucirc3) -- (incirc5);
  \draw[line width=\locpc,color=\colorDefect,postaction={decorate}] (incirc3) -- (incirc5);
  \end{scope}
  \scopeArrow{0.63}{\arrowDefect}
  \draw[line width=\locpc,color=\colorDefect,postaction={decorate}] (incirc1) -- (incirc2);
  \draw[line width=\locpc,color=\colorDefect,postaction={decorate}] (incirc5) -- (incirc1);
  \draw[line width=\locpc,color=\colorTransp,postaction={decorate}]
              (incirc3) to [out=10,in=40] (incirc1);
  \end{scope}
  \scopeArrow{0.76}{\arrowDefect}
  \draw[line width=\locpc,color=\colorTransp,postaction={decorate}] (incirc1)
              .. controls (0.05*\locpa,-0.6*\locpa) and (0.45*\locpa,-0.6*\locpa) .. (incirc2);
  \end{scope}
  \scopeArrow{0.77}{\arrowDefect}
  \draw[line width=\locpc,color=\colorDefect,postaction={decorate}] (incirc2) -- (oucirc6);
  \end{scope}
  \filldraw[\colorCircle,fill=white,line width=\locpC] (incirc1) circle (\locpA)
           (incirc2) circle (\locpA) (incirc3) circle (\locpA) (incirc5) circle (\locpA);
  \scopeArrow{0.33}{<}
  \draw[\colorCircle,line width=\locpf,postaction={decorate}] (incirc1) circle (\locpA);
  \draw[\colorCircle,line width=\locpf,postaction={decorate}] (incirc2) circle (\locpA);
  \draw[\colorCircle,line width=\locpf,postaction={decorate}] (incirc3) circle (\locpA);
  \draw[\colorCircle,line width=\locpf,postaction={decorate}] (incirc5) circle (\locpA);
  \end{scope}
  \filldraw[color=\colorDefect] (oucirc1) circle (\locpd) (oucirc3) circle (\locpd)
           (oucirc4) circle (\locpd) (oucirc5) circle (\locpd) (oucirc6) circle (\locpd)
           (incirc1)+(0:\locpA) circle (\locpD) (incirc1)+(54:\locpA) circle (\locpD)
           (incirc1)+(180:\locpA) circle (\locpD) (incirc1)+(283:\locpA) circle (\locpD)
           (incirc2)+(0:\locpA) circle (\locpD) (incirc2)+(180:\locpA) circle (\locpD)
           (incirc2)+(252:\locpA) circle (\locpD)
           (incirc3)+(59:\locpA) circle (\locpD) (incirc3)+(238:\locpA) circle (\locpD)
           (incirc3)+(347:\locpA) circle (\locpD) (incirc5)+(0:\locpA) circle (\locpD)
           (incirc5)+(57.5:\locpA) circle (\locpD) (incirc5)+(180:\locpA) circle (\locpD)
           (incirc5)+(252:\locpA) circle (\locpD) (incirc5)+(292.8:\locpA) circle (\locpD);
  \end{scope}
  \end{tikzpicture}}
  \label{eq:picFeb1}
  \ee
 
\begin{Remark}\label{refi-vs-transdisk}
\Itemizeiii

\item[(i)]
For any refinement $(\Sigma;\Sigma\refi)$, by definition every gluing boundary of $\Sigma\refi$
that is not a gluing boundary of $\Sigma$ is fillable.

\item[(ii)]
It is worth comparing the notion of refinement with Definition \ref{def:filldisk} of a
\filldisk\ of type $\XX$: It is easily seen that for every refinement $(\XX,\XX\refi)$
with $\XX$ a \defectsurface\ for which the underlying surface $X$ is a disk and which has
at most one free boundary segment, the refined 
surface $\XX\refi$ is a \filldisk\ of type $\XX$. On the other hand, the converse is not true,
because in general the outer boundary of a \filldisk\ of
type $\XX$ contains more defect points (transparently labeled) than the boundary of $\XX$.

\item[(iii)]
It follows directly that  the gluing of any two refinements is again a refinement, i.e.\ that
for any two refinements $(\Sigma_1^{};\Sigma_1')$ and $(\Sigma_2^{};\Sigma_2')$ for which
$\Sigma_1$ and $\Sigma_2$ can be glued, the same holds for $\Sigma_1'$ and $\Sigma_2'$ and the
pair $(\Sigma_1^{}\,{\circ}\,\Sigma_2^{};\Sigma_1'\,{\circ}\,\Sigma_2')$ is a refinement as well.
Moreover, if $\Sigma_{1}'$ and $\Sigma_{2}'$ are \gfine\ refinements, then 
$\Sigma_1'\,{\circ}\,\Sigma_2'$ is \gfine, too.

\item[(iv)]
If $(\Sigma;\Sigma\refi)$ is a refinement, then for any \defectsurface\ $\Sigma\refi'$ that
is isomorphic to $\Sigma\refi$, $(\Sigma;\Sigma\refi')$ is a refinement of $\Sigma$ as well.
\end{itemize}
\end{Remark}

We proceed to show that fine refinements exist. From Definition \ref{def:refinement}
it follows immediately that any refinement $\Sigma'$ of $\Sigma$ can be obtained by 
separately refining every \twocell\ of $\Sigma$, in any order.
Now like any two-manifold with boundary, a \twocell\
admits pair-of-pants decompositions, i.e.\ as a manifold it can be obtained by
gluing a finite number of disks, annuli, and pairs of pants (or trinions).
Since in the situation of our interest we are working with framed surfaces, we must 
in addition account for the framing. Moreover,
unlike for manifolds without defects and boundaries we also have to consider
\twocells\ whose boundary contains free boundary segments. These still admit
generalized pair-of-pants decompositions with a larger number of building blocks
(compare \Cite{Prop.\,3.8}{laPf}).
It is not hard to see that indeed any \twocell\ of any \defectsurface\ can be obtained by
gluing  a finite number of the following specific two-manifolds:
\Itemize
\item
a disk with boundary being a gluing circle without defect points and with framing index $-2$;
\item 
an annulus with boundary components 
being gluing circles without defect points and with framing indices $\pm\kappa \iN 2\mathbb Z$;
\item
a pair of pants with boundary components being gluing circles without defect points and with
framing indices $\kappa,\kappa' \iN 2\mathbb Z$ and $2{-}\kappa{-}\kappa'$, respectively; 
\item
a disk whose boundary is the union of one free boundary and of one gluing interval that does
not have any defect points in its interior and has framing index $-1$;
\item
a disk whose boundary is the union of three free boundaries and of three gluing intervals 
that do not have any defect points in their interior and have framing indices
$\kappa,\kappa' \iN \mathbb Z$ and $1{-}\kappa{-}\kappa'$, respectively;
\item
an \emph{open-closed pipe}: an annulus such that one boundary component is
the union of one free boundary and of one gluing interval without 
defect points in its interior and with framing index $\kappa \iN \mathbb Z$,
while the other boundary component is a gluing circle 
without defect points and with framing index $-1{-}\kappa$.
\end{itemize}

It is worth recalling that a \defectonemanifold\ in the sense of Definition \ref{def:Bordcats}
contains at least one defect point. Thus those of the building blocks in this list which
have a gluing circle without defect points as a boundary component are themselves not
\defectsurface s. However, our prescription implies that any fine surface that results from
the refinement of any (proper) \defectsurface\ \emph{is} again a proper \defectsurface.

We are now in a position to show

\begin{thm}\label{thm:no-name - so far not refrred to!}
Every \defectsurface\ can be refined to a fine \defectsurface.
\end{thm}

\begin{proof}
We will construct a fine refinement for an arbitrary \twocell; combining the refinements of 
all \twocells\ of a \defectsurface\ $\Sigma$ then provides a fine refinement of $\Sigma$.
Thus consider a \twocell\ $\PP$ of a \defectsurface. Take any (generalized)
pair-of-pants decomposition $\Ppop$ of $\PP$, with building blocks of the form just described,
making sure that the curves on $\Sigma$ that define $\Ppop$ intersect all defect lines
on $\Sigma$ transversally (this can be achieved by applying, if necessary, a small isotopy 
to any chosen set of curves).
The resulting building blocks come in two kinds: either they contain a gluing boundary of 
$\Sigma$ or they do not. We first retrict our attention to building blocks of the first kind.
We would now like to define standard refinements of these building blocks of $\Ppop$
that glue together to a refinement of $\Sigma$. However, to do so, we need to alter also gluing
boundaries in the building blocks of $\Ppop$, namely those that do not correspond to gluing 
boundaries of $\Sigma$. We account for this by allowing for \emph{generalized refinements},
which are defined in the same way as refinements, except that there is no restriction on the
end points of the defect lines $\LII_i'$. 
(This is not in conflict with the rationale of our construction, because the
block functor is only applied after gluing such generalized refinements to genuine ones.)
Thus for each of the building blocks of $\Ppop$ we now provide a standard fine generalized
refinement that is chosen in such a manner that gluing together any two of the thereby 
obtained fine \defectsurface s results in a \defectsurface\ that is fine as well.
 \\[.2em]
Let us first present our prescription in much detail for the case that the building block 
in question is an annulus $\AA$ (without defect points on its boundary).
Denote by $\pm\kappa \iN 2\mathbb Z$ the framing indices of the boundary circles of $\AA$.
We refine $\AA$ in two steps. In the first step we add a single transparent defect line
if $\kappa \,{=}\, 0$, while for $\kappa \,{\ne}\, 0$ we add $|\kappa|$
transparent defect lines of alternating orientation
that connect the two boundary circles, in such a way that we obtain a fine \defectsurface\
$\AA_\text{ref,0}$ each of whose \twocells\ is a \basicdisk. This is illustrated in the
following picture: 
    \def\locpa  {0.93}   
    \def\locpA  {2.7}    
    \def\locpb  {0.64}   
    \def\locpB  {3.0}    
    \def\locpc  {2.2pt}  
    \def\locpd  {0.12}   
    \def\locpD  {0.04}   
    \def\locpe  {1.15}   
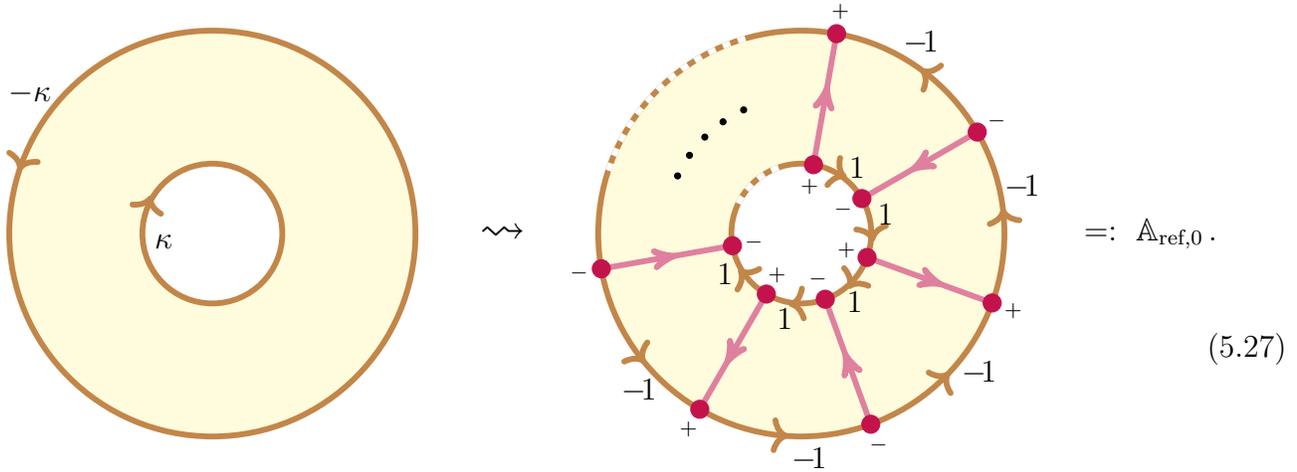
  \begin{eqnarray} 
  \begin{tikzpicture}
  \scopeArrow{0.45}{>}
  \filldraw[\colorCircle,fill=\colorSigma,line width=\locpc,postaction={decorate}]
           (0,0) circle (\locpA) node[left=56pt,yshift=53pt,color=black]{${-}\kappa$};
  \end{scope}
  \scopeArrow{0.45}{<}
  \filldraw[\colorCircle,fill=white,line width=\locpc,postaction={decorate}]
           (0,0) circle (\locpa) node[left=10pt,yshift=-3pt,color=black]{$\kappa$};
  \end{scope}
  \end{tikzpicture} \hspace*{16.5em}
  \nonumber \\[-2.3em]
  \hspace*{5em} \raisebox{-7.6em} {\begin{tikzpicture}
  \begin{scope}[shift={(2.84*\locpA,0)}]
  \scopeArrow{0.02}{>} \scopeArrow{0.15}{>} \scopeArrow{0.61}{>} \scopeArrow{0.74}{>}
  \scopeArrow{0.88}{>}
  \filldraw[\colorCircle,fill=\colorSigma,line width=\locpc,postaction={decorate}]
                       (0,0) circle (\locpA);
  \end{scope} \end{scope} \end{scope} \end{scope} \end{scope}
  \draw[dashed,line width=\locpc,color=white] (106:\locpA) arc (106:164:\locpA);
  \scopeArrow{0.02}{<} \scopeArrow{0.18}{<} \scopeArrow{0.62}{<} \scopeArrow{0.76}{<}
  \scopeArrow{0.90}{<}
  \filldraw[\colorCircle,fill=white,line width=\locpc,postaction={decorate}] (0,0) circle (\locpa);
  \end{scope} \end{scope} \end{scope} \end{scope} \end{scope}
  \draw[dashed,line width=\locpc,color=white] (109:\locpa) arc (109:161:\locpa);
  \scopeArrow{0.61}{\arrowDefect}
  \draw[line width=\locpc,color=\colorTransp,postaction={decorate}] (80:\locpa) -- (80:\locpA);
  \draw[line width=\locpc,color=\colorTransp,postaction={decorate}] (-20:\locpa) -- (-20:\locpA);
  \draw[line width=\locpc,color=\colorTransp,postaction={decorate}] (-120:\locpa) -- (-120:\locpA);
  \end{scope}
  \scopeArrow{0.61}{\arrowDefectB}
  \draw[line width=\locpc,color=\colorTransp,postaction={decorate}] (30:\locpa) -- (30:\locpA);
  \draw[line width=\locpc,color=\colorTransp,postaction={decorate}] (-70:\locpa) -- (-70:\locpA);
  \draw[line width=\locpc,color=\colorTransp,postaction={decorate}] (-170:\locpa) -- (-170:\locpA);
  \end{scope}
  \filldraw[color=\colorDefect] (80:\locpa) circle (\locpd) (80:\locpA) circle (\locpd)
                                (30:\locpa) circle (\locpd) (30:\locpA) circle (\locpd)
                               (-20:\locpa) circle (\locpd) (-20:\locpA) circle (\locpd)
                               (-70:\locpa) circle (\locpd) (-70:\locpA) circle (\locpd)
                              (-120:\locpa) circle (\locpd) (-120:\locpA) circle (\locpd)
                              (-170:\locpa) circle (\locpd) (-170:\locpA) circle (\locpd);
  \node at (80:\locpb)   {$\scriptstyle +$}; \node at (80:\locpB)   {$\scriptstyle +$};
  \node at (30:\locpb)   {$\scriptstyle -$}; \node at (30:\locpB)   {$\scriptstyle -$};
  \node at (-20:\locpb)  {$\scriptstyle +$}; \node at (-20:\locpB)  {$\scriptstyle +$};
  \node at (-70:\locpb)  {$\scriptstyle -$}; \node at (-70:\locpB)  {$\scriptstyle -$};
  \node at (-120:\locpb) {$\scriptstyle +$}; \node at (-120:\locpB) {$\scriptstyle +$};
  \node at (-170:\locpb) {$\scriptstyle -$}; \node at (-170:\locpB) {$\scriptstyle -$};
  \node at (50:\locpe)   {$1$}; \node at (58:\locpB)   {${-}\!1$};
  \node at (13:\locpe)   {$1$}; \node at (12:\locpB)   {${-}\!1$};
  \node at (-52:\locpe)  {$1$}; \node at (-38:\locpB)  {${-}\!1$};
  \node at (-101:\locpe) {$1$}; \node at (-88:\locpB)  {${-}\!1$};
  \node at (-152:\locpe) {$1$}; \node at (-136:\locpB) {${-}\!1$};
  \filldraw (115:0.5*\locpa+0.5*\locpA) circle (\locpD) (125:0.5*\locpa+0.5*\locpA) circle (\locpD)
            (135:0.5*\locpa+0.5*\locpA) circle (\locpD) (145:0.5*\locpa+0.5*\locpA) circle (\locpD)
            (155:0.5*\locpa+0.5*\locpA) circle (\locpD);
  \end{scope}
  \node at (1.4*\locpA,0) {{\Large$\rightsquigarrow$}};
  \node at (4.57*\locpA,0) {$=:~ \AA_\text{ref,0} \,.$};
  \end{tikzpicture}}
  \end{eqnarray}
While the so obtained surface $\AA_\text{ref,0}$ is already fine, gluing it to another
fine surface can still result in a \defectsurface\ that is no longer fine. The second step 
of the prescription eliminates this unwanted feature: We refine $\AA_\text{ref,0}$ further 
by suitably inserting additional transparent defects together with transparent gluing circles
with three defect points that are of the form shown in the picture \eqref{eq:circle.IIn-M-M}
with $\kappa \,{=}\, 0$ and $\calm \,{=}\, \II_0$; thereby we obtain a fine (generalized)
refinement $(\AA;\AA\refi)$ of the form indicated in in the following picture
(for better readability we omit the defect points and orientation of the gluing intervals of the
transparent gluing circles as well as the orientation of some of the defect lines):
    \def\locpa  {0.88}   
    \def\locpA  {3.15}   
    \def\locpb  {0.59}   
    \def\locpB  {3.45}   
    \def\locpc  {2.2pt}  
    \def\locpd  {0.10}   
    \def\locpD  {0.04}   
    \def\locpe  {1.09}   
    \def\locpf  {1.75}   
    \def\locpF  {2.26}   
  \be 
  \raisebox{-9.1em}{\begin{tikzpicture}[scale=1.2]
  \node at (-5,0) {$\AA\refi ~=$};
  \scopeArrow{0.02}{>} \scopeArrow{0.16}{>} \scopeArrow{0.60}{>} \scopeArrow{0.75}{>}
  \scopeArrow{0.88}{>}
  \filldraw[\colorCircle,fill=\colorSigma,line width=\locpc,postaction={decorate}]
           (0,0) circle (\locpA);
  \end{scope} \end{scope} \end{scope} \end{scope} \end{scope}
  \scopeArrow{0.04}{<} \scopeArrow{0.17}{<} \scopeArrow{0.62}{<} \scopeArrow{0.76}{<}
  \scopeArrow{0.90}{<}
  \filldraw[\colorCircle,fill=white,line width=\locpc,postaction={decorate}] (0,0) circle (\locpa);
  \end{scope} \end{scope} \end{scope} \end{scope} \end{scope}
  \scopeArrow{0.27}{\arrowDefect} \scopeArrow{0.90}{\arrowDefect}
  \draw[line width=\locpc,color=\colorTransp,postaction={decorate}] (80:\locpa) -- (80:\locpA);
  \draw[line width=\locpc,color=\colorTransp,postaction={decorate}] (-20:\locpa) -- (-20:\locpA);
  \draw[line width=\locpc,color=\colorTransp,postaction={decorate}] (-120:\locpa) -- (-120:\locpA);
  \end{scope} \end{scope}
  \scopeArrow{0.24}{\arrowDefectB} \scopeArrow{0.88}{\arrowDefectB}
  \draw[line width=\locpc,color=\colorTransp,postaction={decorate}] (30:\locpa) -- (30:\locpA);
  \draw[line width=\locpc,color=\colorTransp,postaction={decorate}] (-70:\locpa) -- (-70:\locpA);
  \draw[line width=\locpc,color=\colorTransp,postaction={decorate}] (-170:\locpa) -- (-170:\locpA);
  \end{scope} \end{scope}
  \scopeArrow{0.55}{\arrowDefectB}
  \draw[line width=\locpc,color=\colorTransp,postaction={decorate}] (80:\locpf) arc (80:30:\locpf);
  \draw[line width=\locpc,color=\colorTransp,postaction={decorate}] (-20:\locpF) arc (-20:30:\locpF);
  \draw[line width=\locpc,color=\colorTransp,postaction={decorate}] (-20:\locpf) arc (-20:-70:\locpf);
  \draw[line width=\locpc,color=\colorTransp,postaction={decorate}] (-120:\locpF) arc (-120:-70:\locpF);
  \draw[line width=\locpc,color=\colorTransp,postaction={decorate}] (-120:\locpf) arc (-120:-170:\locpf);
  \end{scope}
  \draw[line width=\locpc,color=\colorTransp] (80:\locpF) arc (80:105:\locpF);
  \draw[line width=\locpc,color=\colorTransp] (-170:\locpF) arc (-170:-195:\locpF);
  \transv {-170:\locpf} node[right=-1.1pt,yshift=-4pt,color=black]{$\scriptstyle 0$};
  \transv {-170:\locpF} node[right=-3pt,yshift=7pt,color=black]{$\scriptstyle 0$};
  \transv {-70:\locpf}  node[right=-4pt,yshift=7pt,color=black]{$\scriptstyle 0$};
  \transv {-70:\locpF}  node[left,yshift=4pt,color=black]{$\scriptstyle 0$};
  \transv {-120:\locpf} node[above,xshift=-1pt,color=black]{$\scriptstyle 0$};
  \transv {-120:\locpF} node[above=-4pt,xshift=7pt,color=black]{$\scriptstyle 0$};
  \transv {-20:\locpf}; \transv {-20:\locpF}; 
  \transv {30:\locpf}; \transv {30:\locpF};
  \transv {80:\locpf}; \transv {80:\locpF}; 
  \filldraw[color=\colorDefect] (80:\locpa) circle (\locpd) (80:\locpA) circle (\locpd)
                                (30:\locpa) circle (\locpd) (30:\locpA) circle (\locpd)
                               (-20:\locpa) circle (\locpd) (-20:\locpA) circle (\locpd)
                               (-70:\locpa) circle (\locpd) (-70:\locpA) circle (\locpd)
                              (-120:\locpa) circle (\locpd) (-120:\locpA) circle (\locpd)
                              (-170:\locpa) circle (\locpd) (-170:\locpA) circle (\locpd);
  \filldraw (115:0.5*\locpa+0.5*\locpA) circle (\locpD)
            (125:0.5*\locpa+0.5*\locpA) circle (\locpD)
            (135:0.5*\locpa+0.5*\locpA) circle (\locpD)
            (145:0.5*\locpa+0.5*\locpA) circle (\locpD)
            (155:0.5*\locpa+0.5*\locpA) circle (\locpD);
  \draw[dashed,line width=\locpc,color=white] (106:\locpA) arc (106:164:\locpA);
  \draw[dashed,line width=\locpc,color=white] (109:\locpa) arc (109:161:\locpa);
  \end{tikzpicture}}
  \ee
The framing on the so obtained fine surface $\AA\refi$ is uniquely specified by the orientation
of the various defect lines together with the framing on $\AA_\text{ref,0}$. Specifically, the
framing indices of the three segments of each of the new transparent gluing circles
are $0$, $1$ and $1$. We indicate the position of some of those segments that have framing index
$0$ in the picture; the position of the index-$0$ segments for the other transparent
gluing circles is analogous.
 \\[.2em]
The standard fine refinements for the other building blocks of a pair-of-pants decomposition
are obtained in a similar manner as in the case of the annulus $\AA$. We content ourselves to
display the resulting fine refinements $(\BB;\BB\refi)$ and $(\mathbb U;\mathbb U\refi)$ for
the cases of a pair of pants $\BB$ and of an open-closed pipe $\mathbb U$. The standard fine
refinement for $\BB$ looks as follows:
 \\[-1.4em]
    \def\locpa  {0.88}   
    \def\locpA  {1.1}    
    \def\locpc  {2.2pt}  
    \def\locpd  {0.12}   
    \def\locpD  {0.04}   
    \def\locpx  {5.3}    
    \def\locpy  {3.0}    
  \begin{eqnarray}
  \begin{tikzpicture}[scale=1.1]
  \node at (-7.1,0) {$\BB\refi ~=$};
  \filldraw[fill=\colorSigma,line width=\locpc,draw=white]
           (0,0) circle (\locpx cm and \locpy cm);
  \node at (0.15*\locpx,0.92*\locpy) {$\scriptstyle 0$};
  \node at (-0.15*\locpx,-0.92*\locpy) {$\scriptstyle 0$};
  \scopeArrow{0.09}{>} \scopeArrow{0.25}{>} \scopeArrow{0.35}{>}
  \scopeArrow{0.57}{>} \scopeArrow{0.75}{>} \scopeArrow{0.85}{>}
  \draw[line width=\locpc,color=\colorCircle,postaction={decorate}]
           (0,0) circle (\locpx cm and \locpy cm);
  \end{scope} \end{scope} \end{scope} \end{scope} \end{scope} \end{scope}
  \coordinate (dashpt1) at ({\locpx*cos(154)},{\locpy*sin(154)});
  \coordinate (dashpt2) at ({\locpx*cos(-24)},{\locpy*sin(-24)});
  \draw[line width=\locpc,color=white,dashed]
           (dashpt1) arc (154:204:\locpx cm and \locpy cm)
           (dashpt2) arc (-24:26:\locpx cm and \locpy cm);
  \begin{scope}[shift={(-\locpa-\locpA,0)}]
  \filldraw[draw=\colorSigma,line width=\locpc,fill=white] (0,0) circle (\locpa);
  \scopeArrow{0.30}{<} \scopeArrow{0.64}{<} \scopeArrow{0.99}{<}
  \draw[line width=\locpc,color=\colorCircle,postaction={decorate}] (0,0) circle (\locpa);
  \end{scope} \end{scope} \end{scope}
  \draw[line width=\locpc,color=white,dashed] (145:\locpa) arc (145:210:\locpa);
  \filldraw (180:2.26*\locpa) circle (\locpD)
            (191:2.2*\locpa) circle (\locpD) (169:2.2*\locpa) circle (\locpD)
            (202:2.0*\locpa) circle (\locpD) (158:2.0*\locpa) circle (\locpD);
  \end{scope}
  \begin{scope}[shift={(\locpa+\locpA,0)}]
  \filldraw[draw=\colorSigma,line width=\locpc,fill=white] (0,0) circle (\locpa);
  \scopeArrow{0.13}{<} \scopeArrow{0.51}{<} \scopeArrow{0.79}{<}
  \draw[line width=\locpc,color=\colorCircle,postaction={decorate}] (0,0) circle (\locpa);
  \end{scope} \end{scope} \end{scope}
  \draw[line width=\locpc,color=white,dashed] (-30:\locpa) arc (-30:35:\locpa);
  \filldraw (0:2.26*\locpa) circle (\locpD)
            (11:2.2*\locpa) circle (\locpD) (-11:2.2*\locpa) circle (\locpD)
            (22:2.0*\locpa) circle (\locpD) (-22:2.0*\locpa) circle (\locpD);
  \end{scope}
  \coordinate (defpt1) at ({\locpx*cos(55)},{\locpy*sin(55)});
  \coordinate (defpt2) at ({\locpx*cos(104)},{\locpy*sin(104)});
  \coordinate (defpt3) at ({\locpx*cos(136)},{\locpy*sin(136)});
  \coordinate (defpt4) at ({\locpx*cos(235)},{\locpy*sin(235)});
  \coordinate (defpt5) at ({\locpx*cos(284)},{\locpy*sin(284)});
  \coordinate (defpt6) at ({\locpx*cos(316)},{\locpy*sin(316)});
  \coordinate (Defpt1) at ({\locpA+\locpa+\locpa*cos(69)},{\locpa*sin(69)});
  \coordinate (Defpt2) at ({-\locpA-\locpa+\locpa*cos(76)},{\locpa*sin(76)});
  \coordinate (Defpt3) at ({-\locpA-\locpa+\locpa*cos(122)},{\locpa*sin(122)});
  \coordinate (Defpt4) at ({-\locpA-\locpa+\locpa*cos(249)},{\locpa*sin(249)});
  \coordinate (Defpt5) at ({\locpA+\locpa+\locpa*cos(250)},{\locpa*sin(250)});
  \coordinate (Defpt6) at ({\locpA+\locpa+\locpa*cos(302)},{\locpa*sin(302)});
  \scopeArrow{0.34}{\arrowDefect} \scopeArrow{0.85}{\arrowDefect}
  \draw[line width=\locpc,color=\colorTransp,postaction={decorate}] (defpt1) -- (Defpt1);
  \draw[line width=\locpc,color=\colorTransp,postaction={decorate}] (defpt2) -- (Defpt2);
  \draw[line width=\locpc,color=\colorTransp,postaction={decorate}] (defpt6) -- (Defpt6);
  \end{scope} \end{scope}
  \scopeArrow{0.34}{\arrowDefectB} \scopeArrow{0.85}{\arrowDefectB}
  \draw[line width=\locpc,color=\colorTransp,postaction={decorate}] (defpt3) -- (Defpt3);
  \draw[line width=\locpc,color=\colorTransp,postaction={decorate}] (defpt4) -- (Defpt4);
  \draw[line width=\locpc,color=\colorTransp,postaction={decorate}] (defpt5) -- (Defpt5);
  \end{scope} \end{scope}
  \scopeArrow{0.62}{\arrowDefect}
  \draw[line width=\locpc,color=\colorTransp,postaction={decorate}]
          ($(defpt2)!0.5!(Defpt2)$) to [out=170,in=46] ($(defpt3)!0.5!(Defpt3)$);
  \end{scope}
  \scopeArrow{0.62}{\arrowDefectB}
  \draw[line width=\locpc,color=\colorTransp,postaction={decorate}]
          ($(defpt5)!0.5!(Defpt5)$) to [out=-10,in=216] ($(defpt6)!0.5!(Defpt6)$);
  \draw[line width=\locpc,color=\colorTransp,postaction={decorate}]
          ($(defpt4)!0.5!(Defpt4)$) to [out=-25,in=155] ($(defpt1)!0.5!(Defpt1)$);
  \end{scope}
  \filldraw[color=\colorDefect] (defpt1) circle (\locpd) (defpt2) circle (\locpd)
                                (defpt3) circle (\locpd) (defpt4) circle (\locpd)
                                (defpt5) circle (\locpd) (defpt6) circle (\locpd)
                                (Defpt1) circle (\locpd) (Defpt2) circle (\locpd)
                                (Defpt3) circle (\locpd) (Defpt4) circle (\locpd)
                                (Defpt5) circle (\locpd) (Defpt6) circle (\locpd);
  \transv {$(defpt1)!0.5!(Defpt1)$} node[left=-0.6pt,yshift=-2.5pt,color=black]{$\scriptstyle 0$};
  \transv {$(defpt2)!0.5!(Defpt2)$} node[left=-1.1pt,yshift=-4.8pt,color=black]{$\scriptstyle 0$};
  \transv {$(defpt3)!0.5!(Defpt3)$} node[right=.1pt,yshift=-0.1pt,color=black]{$\scriptstyle 0$};
  \transv {$(defpt4)!0.5!(Defpt4)$} node[right=-.5pt,yshift=3.3pt,color=black]{$\scriptstyle 0$};
  \transv {$(defpt5)!0.5!(Defpt5)$} node[right=-1.6pt,yshift=5.2pt,color=black]{$\scriptstyle 0$};
  \transv {$(defpt6)!0.5!(Defpt6)$} node[left=.5pt,yshift=0.8pt,color=black]{$\scriptstyle 0$};
  \end{tikzpicture}
  \nonumber \\[-3.3em] ~
  \end{eqnarray}
Drawing the open-closed pipe $\mathbb U$ as on the right hand side of
    \def\locpa  {2.8}    
    \def\locpA  {2.3}    
    \def\locpb  {0.42}   
    \def\locpc  {2.2pt}  
    \def\locpC  {1.8pt}  
    \def\locpd  {0.12}   
    \def\locpD  {0.08}   
    \def\locpe  {6.8}    
    \def\locpf  {1.1pt}  
    \def\locpt  {0.34}   
    \def\locpv  {50}     
    \def\locpw  {130}    
    \def\locpx  {70}     
    \def\locpy  {110}    
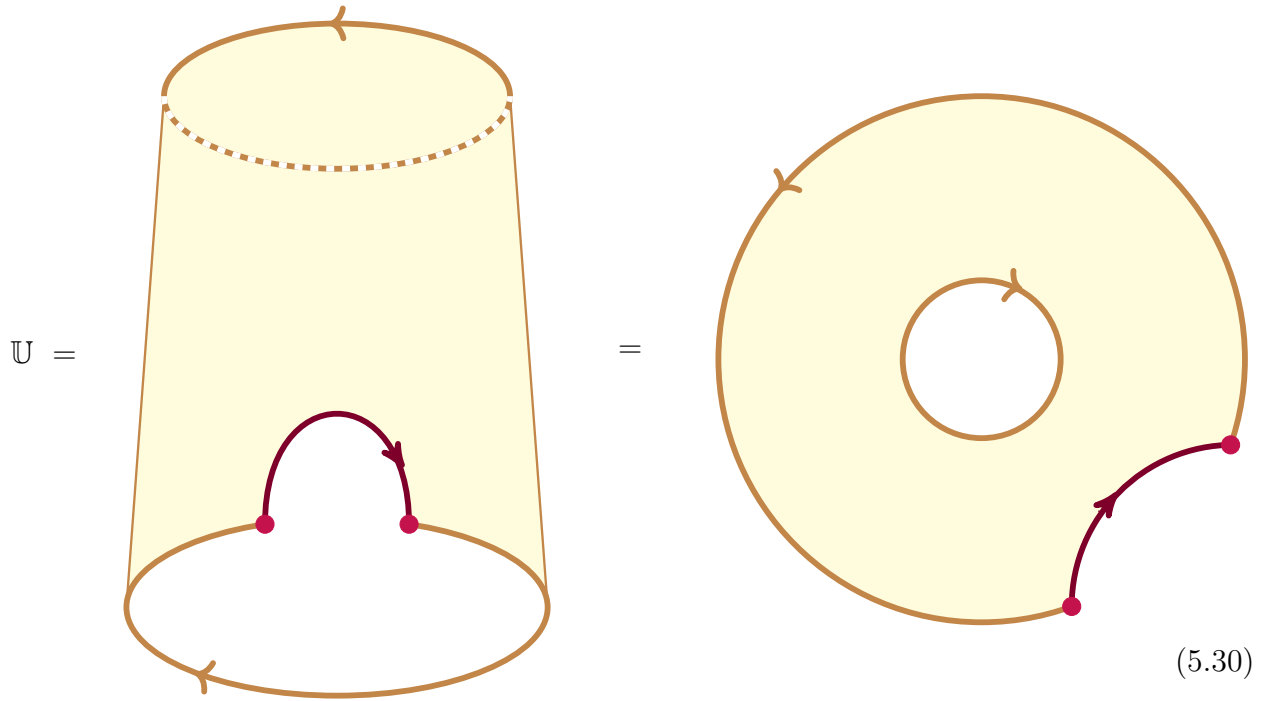
  \begin{eqnarray}
  \begin{tikzpicture}
  \coordinate (bpoint1b) at ({\locpa*cos(\locpv)},{\locpb*\locpa*sin(\locpv)});
  \coordinate (bpoint2b) at ({\locpa*cos(\locpw)},{\locpb*\locpa*sin(\locpw)});
  \coordinate (bpoint1N) at ({\locpa*cos(\locpx)},{\locpb*\locpa*sin(\locpx)});
  \coordinate (bpoint2N) at ({\locpa*cos(\locpy)},{\locpb*\locpa*sin(\locpy)});
  \filldraw[line width=0.4*\locpc,\colorCircle,fill=\colorSigma]
           (-\locpA,\locpe) -- (-\locpa,0) -- (\locpa,0) -- (\locpA,\locpe) -- cycle;
  \scopeArrow{0.62}{<}
  \filldraw[line width=\locpc,\colorCircle,fill=white,postaction={decorate}]
           (0,0) circle (\locpa cm and \locpb*\locpa cm);
  \end{scope}
  \scopeArrow{0.26}{>}
  \filldraw[line width=\locpc,\colorCircle,fill=\colorSigma,postaction={decorate}]
           (0,\locpe) circle (\locpA cm and \locpb*\locpA cm);
  \end{scope}
  \draw[line width=1.3*\locpc,color=white]
           (bpoint1N) arc (\locpx:\locpy:\locpa cm and \locpb*\locpa cm);
  \draw[line width=\locpc,color=white,dashed]
           (-\locpA,\locpe) arc (180:360:\locpA cm and \locpb*\locpA cm);
  \scopeArrow{0.28}{\arrowDefectB}
  \filldraw[line width=\locpc,draw=\colorFreebdy,fill=white,postaction={decorate}]
           (bpoint1N) .. controls ++(0,0.7*\locpa) and ++(0,0.7*\locpa) .. (bpoint2N);
  \end{scope}
  \filldraw[color=\colorDefect] (bpoint1N) circle (\locpd) (bpoint2N) circle (\locpd);
  \node at (-\locpa-0.9,0.5*\locpe) {$\mathbb U ~=$};
  \node at (\locpa+0.7,0.5*\locpe) {$=$};
    \def\locpa  {3.5}    
    \def\locpb  {1.05}   
    \def\locpc  {2.2pt}  
    \def\locpd  {0.12}   
    \def\locpe  {2.15}   
    \def\locpy  {290}    
    \def\locpY  {180}    
    \def\locpz  {341}    
    \def\locpZ  {90}     
  \begin{scope}[shift={(2*\locpa+1.0,\locpa-0.2)}]
  \scopeArrow{0.39}{>}
  \filldraw[\colorCircle,fill=\colorSigma,line width=\locpc,postaction={decorate}] (0,0) circle (\locpa);
  \end{scope}
  \scopeArrow{0.20}{<}
  \filldraw[\colorCircle,fill=white,line width=\locpc,postaction={decorate}] (0,0) circle (\locpb);
  \end{scope}
  \fill[white] (0.5*\locpy+0.5*\locpz:\locpa+0.55*\locpe) circle (\locpe);
  \draw[\colorFreebdy,line width=\locpc] (\locpy:\locpa) arc (\locpY:\locpZ:\locpe);
  \filldraw[color=\colorDefect] (\locpy:\locpa) circle (\locpd) (\locpz:\locpa) circle (\locpd);
  \scopeArrow{0.5}{\arrowDefectB}
  \draw[line width=\locpc,\colorFreebdy,postaction={decorate}] (308:0.737*\locpa) -- +(230:0.02);
  \end{scope}
  \end{scope}
  \end{tikzpicture} \hspace*{-6.7em}
  \nonumber\\[-5.6em]~\\[-0.1em]~\nonumber
  \end{eqnarray}
its standard fine refinement looks as follows:
    \def\locpa  {3.9}    
    \def\locpb  {0.85}   
    \def\locpc  {2.2pt}  
    \def\locpd  {0.12}   
    \def\locpD  {0.104}  
    \def\locpe  {2.38}   
    \def\locpf  {0.51}   
    \def\locpx  {9.9}    
    \def\locpy  {290}    
    \def\locpY  {180}    
    \def\locpz  {341}    
    \def\locpZ  {90}     
  \be
  \raisebox{-14.4em}
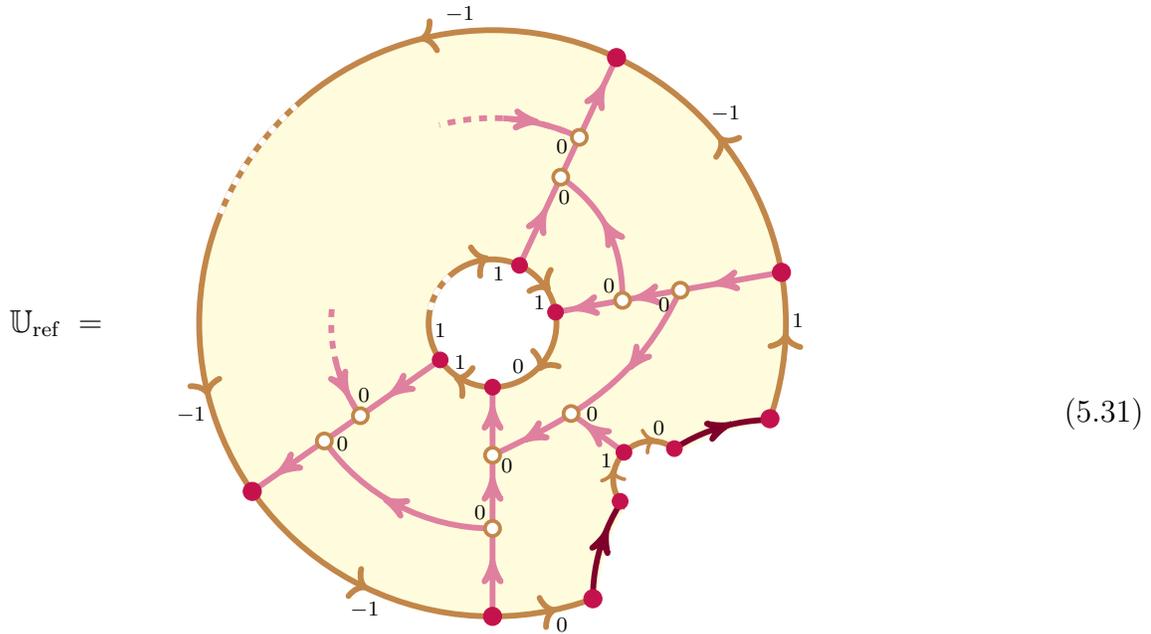
{\begin{tikzpicture}
  \coordinate (cooxvx) at (311:0.408*\locpa);   
  \node at (-\locpa-1.9,0) {$\mathbb U\refi ~=$};
  \scopeArrow{0.11}{>} \scopeArrow{0.29}{>} \scopeArrow{0.54}{>} \scopeArrow{0.68}{>}
  \scopeArrow{0.785}{>} \scopeArrow{0.995}{>}
  \filldraw[\colorCircle,fill=\colorSigma,line width=\locpc,postaction={decorate}]
           (0,0) circle (\locpa);
  \end{scope} \end{scope} \end{scope} \end{scope} \end{scope} \end{scope}
  \scopeArrow{0.124}{<} \scopeArrow{0.30}{<} \scopeArrow{0.693}{<} \scopeArrow{0.914}{<}
  \filldraw[\colorCircle,fill=white,line width=\locpc,postaction={decorate}] (0,0) circle (\locpb);
  \end{scope} \end{scope} \end{scope} \end{scope}
  \draw[white,line width=\locpc,dashed]
           (158:\locpa) arc (158:132:\locpa) (168:\locpb) arc (168:130:\locpb);
  \scopeArrow{0.27}{\arrowDefect} \scopeArrow{0.87}{\arrowDefect}
  \draw[line width=\locpc,\colorTransp,postaction={decorate}] (65:\locpb) -- (65:\locpa);
  \draw[line width=\locpc,\colorTransp,postaction={decorate}] (215:\locpb) -- (215:\locpa);
  \end{scope} \end{scope}
  \scopeArrow{0.19}{\arrowDefectB} \scopeArrow{0.46}{\arrowDefectB}
  \scopeArrow{0.81}{\arrowDefectB}
  \draw[line width=\locpc,\colorTransp,postaction={decorate}] (10:\locpb) -- (10:\locpa);
  \end{scope} \end{scope} \end{scope}
  \scopeArrow{0.19}{\arrowDefectB} \scopeArrow{0.50}{\arrowDefectB}
  \scopeArrow{0.86}{\arrowDefectB}
  \draw[line width=\locpc,\colorTransp,postaction={decorate}] (270:\locpb) -- (270:\locpa);
  \end{scope} \end{scope} \end{scope}
  \scopeArrow{0.67}{\arrowDefectB} 
  \draw[line width=\locpc,\colorTransp,postaction={decorate}] 
	   (cooxvx) -- (0.5*\locpy+0.5*\locpz:0.628*\locpa);
  \end{scope} 
  \scopeArrow{0.23}{\arrowDefectB} \scopeArrow{0.73}{\arrowDefectB} 
  \draw[line width=\locpc,\colorTransp,postaction={decorate}]
           (-90:0.45*\locpa) [out=25,in=245] to (10:0.65*\locpa);
  \end{scope} \end{scope}
  \scopeArrow{0.57}{\arrowDefect} 
  \draw[line width=\locpc,\colorTransp,postaction={decorate}]
           (10:0.45*\locpa) [out=90,in=325] to (65:0.55*\locpa);
  \end{scope} 
  \scopeArrow{0.87}{\arrowDefectB} 
  \draw[line width=\locpc,\colorTransp,postaction={decorate}]
           (65:0.7*\locpa) arc (65:87:0.7*\locpa);
  \end{scope} 
  \scopeArrow{0.62}{\arrowDefect} 
  \draw[line width=\locpc,\colorTransp,postaction={decorate}]
           (193:0.55*\locpa) arc (193:215:0.55*\locpa);
  \end{scope} 
  \draw[line width=\locpc,\colorTransp,dashed] (65:0.7*\locpa) arc (65:105:0.7*\locpa);
  \draw[line width=\locpc,\colorTransp,dashed] (175:0.55*\locpa) arc (175:215:0.55*\locpa);
  \scopeArrow{0.55}{\arrowDefectB} 
  \draw[line width=\locpc,\colorTransp,postaction={decorate}]
           (215:0.7*\locpa) arc (215:270:0.7*\locpa);
  \end{scope} 
  \draw[\colorFreebdy,line width=\locpc] (\locpy:\locpa) arc (\locpY:\locpZ:\locpe);
  \filldraw[\colorCircle,fill=white,line width=\locpc]
           (0.5*\locpy+0.5*\locpz:0.76*\locpa) circle (\locpf);
  \scopeArrow{0.285}{<} \scopeArrow{0.51}{<}
  \draw[\colorCircle,line width=0.65*\locpc,postaction={decorate}]
           (0.5*\locpy+0.5*\locpz:0.76*\locpa) circle (\locpf);
  \end{scope} \end{scope} 
  \fill[white] (0.5*\locpy+0.5*\locpz:\locpa+0.59*\locpe) circle (\locpe);
  \filldraw[color=\colorDefect] (10:\locpb) circle (\locpD) (65:\locpb) circle (\locpD)
           (215:\locpb) circle (\locpD) (270:\locpb) circle (\locpD)
	   (10:\locpa) circle (\locpd) (65:\locpa) circle (\locpd)
           (215:\locpa) circle (\locpd) (270:\locpa) circle (\locpd)
           (\locpy:\locpa) circle (\locpd) (\locpz:\locpa) circle (\locpd)
           (0.5*\locpy+0.5*\locpz-\locpx:0.747*\locpa) circle (\locpD)
           (0.5*\locpy+0.5*\locpz+\locpx:0.753*\locpa) circle (\locpD)
           (0.5*\locpy+0.5*\locpz:0.628*\locpa) circle (\locpD) 
                            node[left=0.4pt,yshift=-3.1pt,color=black]{$\scriptstyle 1$}
                            node[right=6.9pt,yshift=9.3pt,color=black]{$\scriptstyle 0$};
  \transv {10:0.45*\locpa}  node[above=-1.5pt,xshift=-5.1pt,color=black]{$\scriptstyle 0$}; 
  \transv {10:0.65*\locpa}  node[below=-1.8pt,xshift=-6.1pt,color=black]{$\scriptstyle 0$}; 
  \transv {65:0.55*\locpa}  node[below=0.4pt,xshift=1.4pt,color=black]{$\scriptstyle 0$};
  \transv {65:0.7*\locpa}   node[below=-3.9pt,xshift=-6.5pt,color=black]{$\scriptstyle 0$};
  \transv {215:0.55*\locpa} node[above=0.5pt,xshift=1.4pt,color=black]{$\scriptstyle 0$};
  \transv {215:0.7*\locpa}  node[right=0.1pt,yshift=-1.1pt,color=black]{$\scriptstyle 0$}; 
  \transv {270:0.45*\locpa} node[below=-3.3pt,xshift=5.5pt,color=black]{$\scriptstyle 0$};
  \transv {270:0.7*\locpa}  node[above=-1.1pt,xshift=-4.7pt,color=black]{$\scriptstyle 0$};
  \transv {cooxvx} node[right=1.3pt,color=black]{$\scriptstyle 0$};
  \node at (24:\locpb-0.17) {$\scriptstyle 1$};
  \node at (83:\locpb-0.19) {$\scriptstyle 1$};
  \node at (188:\locpb-0.15) {$\scriptstyle 1$};
  \node at (230:\locpb-0.18) {$\scriptstyle 1$};
  \node at (301:\locpb-0.19) {$\scriptstyle 0$};
  \node at (0.5:\locpa+0.16) {$\scriptstyle 1$};
  \node at (42:\locpa+0.27) {$\scriptstyle {-}1$};
  \node at (96:\locpa+0.24) {$\scriptstyle {-}1$};
  \node at (197:\locpa+0.29) {$\scriptstyle {-}1$};
  \node at (246:\locpa+0.28) {$\scriptstyle {-}1$};
  \node at (283:\locpa+0.21) {$\scriptstyle 0$};
  \scopeArrow{0.5}{\arrowDefect}
  \draw[line width=\locpc,\colorFreebdy,postaction={decorate}] (299:0.805*\locpa) -- +(71:0.02);
  \draw[line width=\locpc,\colorFreebdy,postaction={decorate}] (337:0.881*\locpa) -- +(19:0.02);
  \end{scope}
  \end{tikzpicture}}
  \ee
  ~\\[-5.8em]

Finally we consider a building block $\PP$ of the $\Ppop$ of the second kind, i.e.\ one
that contains a gluing boundary of $\Sigma$. By picking, if necessary, a finer 
pair-of-pants decomposition $\Ppop'$, we can assume that $\PP$ is a cylinder. 
Consider first the case that the gluing boundary of $\Sigma$ in question is a circle; then
we deal with a \defectsurface\ of the following form:
    \def\locpa  {3.4}    
    \def\locpA  {0.96}   
    \def\locpb  {0.5*\locpa+0.55*\locpA}
    \def\locpc  {2.2pt}  
    \def\locpC  {1.8pt}  
    \def\locpd  {0.12}   
    \def\locpD  {0.09}   
    \def\locpf  {1.1pt}  
    \def\locpt  {0.34}   
    \def\locpu  {105}    
    \def\locpv  {66}     
    \def\locpw  {191}    
    \def\locpx  {227}    
    \def\locpy  {267}    
    \def\locpz  {36}     
  \be
  \raisebox{-8.5em}
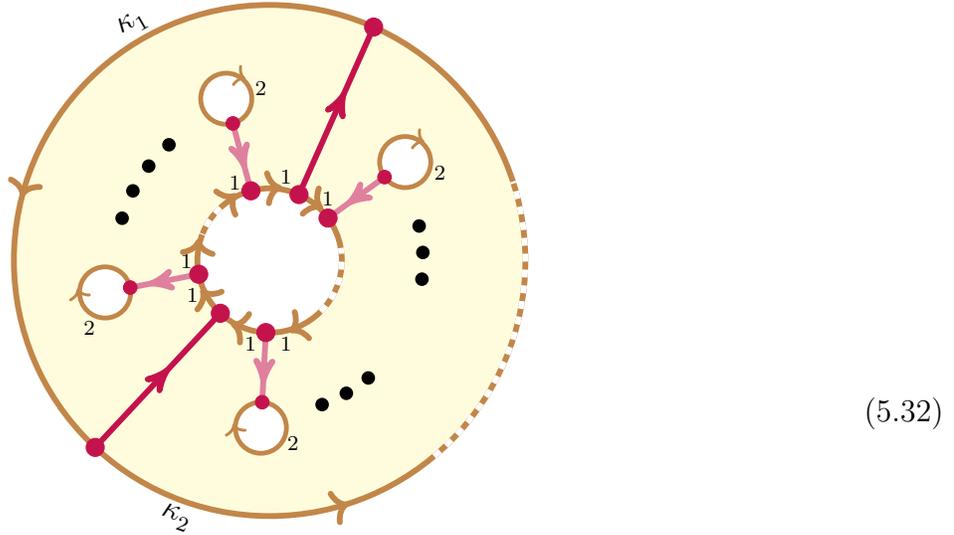
{\begin{tikzpicture}
  \scopeArrow{0.46}{>} \scopeArrow{0.80}{>}
  \filldraw[\colorCircle,fill=\colorSigma,line width=\locpc,postaction={decorate}] (0,0) circle (\locpa);
  \end{scope} \end{scope}
  \draw[line width=\locpc,color=white,dashed] (-50:\locpa) arc (-50:18:\locpa);
  \scopeArrow{0.61}{\arrowDefect}
  \draw[line width=\locpc,color=\colorDefect,postaction={decorate}] (\locpv:\locpA) -- (\locpv:\locpa);
  \end{scope}
  \scopeArrow{0.55}{\arrowDefectB}
  \draw[line width=\locpc,color=\colorDefect,postaction={decorate}] (\locpx:\locpA) -- (\locpx:\locpa);
  \end{scope}
  \scopeArrow{0.51}{\arrowDefectB}
  \draw[line width=\locpc,color=\colorTransp,postaction={decorate}]
       (\locpu:\locpA) -- (\locpu:\locpb) node[right=6pt,yshift=4pt,color=black] {$\scriptstyle 2$};
  \draw[line width=\locpc,color=\colorTransp,postaction={decorate}]
       (\locpz:\locpA) -- (\locpz:\locpb) node[right=6pt,yshift=-4pt,color=black] {$\scriptstyle 2$};
  \end{scope}
  \scopeArrow{0.54}{\arrowDefect}
  \draw[line width=\locpc,color=\colorTransp,postaction={decorate}]
       (\locpw:\locpA) -- (\locpw:\locpb) node[below=6pt,xshift=-6pt,color=black] {$\scriptstyle 2$};
  \draw[line width=\locpc,color=\colorTransp,postaction={decorate}]
       (\locpy:\locpA) -- (\locpy:\locpb) node[right=5pt,yshift=-6pt,color=black] {$\scriptstyle 2$};
  \end{scope}
  \scopeArrow{0.16}{<} \scopeArrow{0.26}{<} \scopeArrow{0.36}{<} \scopeArrow{0.49}{<}
  \scopeArrow{0.60}{<} \scopeArrow{0.70}{<} \scopeArrow{0.83}{<}
  \filldraw[\colorCircle,fill=white,line width=\locpc,postaction={decorate}] (0,0) circle (\locpA);
  \end{scope} \end{scope} \end{scope} \end{scope} \end{scope} \end{scope} \end{scope}
  \draw[line width=\locpc,color=white,dashed]
             (-45:\locpA) arc (-45:10:\locpA) (133:\locpA) arc (133:159:\locpA);
  \filldraw[\colorCircle,fill=white,line width=\locpC]
             (\locpu:\locpb) circle (\locpt) (\locpw:\locpb) circle (\locpt)
             (\locpy:\locpb) circle (\locpt) (\locpz:\locpb) circle (\locpt);
  \node at (47:1.18*\locpA) {$\scriptstyle 1$};
  \node at (79:1.18*\locpA) {$\scriptstyle 1$};
  \node at (114:1.18*\locpA) {$\scriptstyle 1$};
  \node at (180.8:1.15*\locpA) {$\scriptstyle 1$};
  \node at (204:1.17*\locpA) {$\scriptstyle 1$};
  \node at (257:1.18*\locpA) {$\scriptstyle 1$};
  \node at (281:1.18*\locpA) {$\scriptstyle 1$};
  \scopeArrow{0.54}{<}
  \draw[\colorCircle,line width=\locpf,postaction={decorate}] (\locpw:\locpb) circle (\locpt);
  \draw[\colorCircle,line width=\locpf,postaction={decorate}] (\locpy:\locpb) circle (\locpt);
  \end{scope}
  \scopeArrow{0.19}{<}
  \draw[\colorCircle,line width=\locpf,postaction={decorate}] (\locpu:\locpb) circle (\locpt);
  \draw[\colorCircle,line width=\locpf,postaction={decorate}] (\locpz:\locpb) circle (\locpt);
  \end{scope}
  \filldraw[color=\colorDefect] (\locpv:\locpa) circle (\locpd) (\locpv:\locpA) circle (\locpd)
             (\locpx:\locpa) circle (\locpd) (\locpx:\locpA) circle (\locpd)
             (\locpu:\locpA) circle (\locpd) (\locpu:\locpb-\locpt) circle (\locpD)
             (\locpw:\locpA) circle (\locpd) (\locpw:\locpb-\locpt) circle (\locpD)
             (\locpy:\locpA) circle (\locpd) (\locpy:\locpb-\locpt) circle (\locpD)
             (\locpz:\locpA) circle (\locpd) (\locpz:\locpb-\locpt) circle (\locpD);
  \foreach \xx in {131,142,153,164,290,300,310,353,3,13} \fill (\xx:0.6*\locpa) circle (\locpD);
  \node[rotate=30]  at (120:1.07*\locpa) {$\kappa_1^{}$};
  \node[rotate=-20] at (250:1.07*\locpa) {$\kappa_2^{}$};
  \end{tikzpicture}}
  \ee
Here for each gluing segment with framing index $\kappa_i$ of the outer circle there are
$|\kappa_i^{}|\,{-}\,1$ additional boundary circles that are connected by a single transparent defect 
line to the inner gluing circle in such a way that the framing index of each of the resulting
new gluing segments on the inner circle is $+1$ (as shown in the picture) if $\kappa_i$
is positive and $-1$ if $\kappa_i$ is negative. This prescription is compatible with 
the previous treatment of building blocks that do not meet the boundary of $\Sigma$,
in such a way that the resulting refinement of $\Sigma$ is still fine.
The case that the relevant gluing boundary of $\Sigma$ is an interval is treated similarly.
 \\[.2em]
To summarize: Given a \defectsurface\ $\Sigma$ we choose for each \twocell\ $\PP$ of 
$\Sigma$ any pair-of-pants decomposition $\Ppop$ and refine 
the building pieces of $\Ppop$ in
the way described above. Gluing the so obtained refined building pieces 
back together gives a 
fine refinement $(\PP;\PP\refi)$ of $\PP$. Combining the fine refinements
$(\PP;\PP_{\ell;\text{ref}})$ of all \twocells\ $\PP_\ell$ of $\Sigma$ then provides
a fine refinement $(\Sigma;\Sigma\refi)$ of $\Sigma$, provided that the new transparent
defect lines in neighboring \twocells\ are inserted in a coordinated fashion, which is easily
accomplished.  Thus in short, each collection of pair-of-pants decompositions of the 
\twocells\ of $\Sigma$ defines a specific fine refinement of $\Sigma$.
\end{proof}

\begin{Remark}
Consider the standard two-torus $\TT \eq [0,1]^{2}/{\sim}$. It inherits a framing 
from the standard framing of the plane $\R^{2}$. We denote the \defectsurface\ with this
standard framing by $(\TT, \chi^\text{std})$. Any other framing $\chi$ on $\TT$ defines an 
isomorphic \defectsurface. Indeed, for any non-vanishing vector field on $\TT$ there exists a
pair of closed curves with zero winding number for the framing that generate the fundamental 
group $\pi_1(\TT)$ of the torus. If 
we place transparent defect lines on $\TT$ that follow these closed
curves, except that their intersection is resolved into two gluing circles each having
three defect points, we arrive at a valid \defectsurface\ whose gluing circles are of type
\eqref{eq:circle.IIn-M-M} with $\kappa \,{=}\, 0$, and thereby at a fine refinement
$(\TT;\TT\refi)_\chi^{}$ of $\TT$;
identifying $(\TT;\TT\refi)_\chi^{}$ with a refinement of the standard framed torus 
yields the desired isomorphism. 
 \\
Moreover, every automorphism $\varphi$ of $(\TT, \chi^\text{std})$ is isotopic to the 
identity: We can assume that $\varphi$ preserves $(0,0)$. Since $\varphi$ must preserve
flow lines, it maps the generator corresponding to $(1,0)$ to itself, and since it
must preserves angles between tangent vectors, it maps the $(0,1)$ generator to itself
as well. 
\end{Remark}


\subsection{Block functors for non-fine surfaces}\label{sec:nonfine}

In Section \ref{sec:pre-blocks} we have constructed the pre-block functor $\Zpre(\Sigma)$
for any \defectsurface\ $\Sigma$, not necessarily fine. In contrast, the block functor
cannot be defined for $\Sigma$ as in Section \ref{sec:blocks} unless the surface $\Sigma$
is fine, since the topology of the disk enters via the well-posedness, in the sense of
Definition \ref{def:well-posed}, for disks (see Corollary \ref{cor:well-posed}).
Our goal is now to define the block functor $\ZZ(\Sigma)$ for an
arbitrary \defectsurface\ $\Sigma$ by making use of the existence
of fine refinements of $\Sigma$ that was shown in Section \ref{sec:refin-defect-surf}.

Let $(\Sigma;\Sigma\refi)$ be a fine refinement from $\Sigma$ to $\Sigma\refi$. In the
first place, $\Sigma\refi$ is just a specific fine \defectsurface, which comes with its
own pre-block functor $\Zpre(\Sigma\refi)$ and block functor $\ZZ(\Sigma\refi)$. 
Typically $\ZZ(\Sigma\refi)$ lies in a different functor category than the
desired functor $\ZZ(\Sigma)$. In order to obtain a functor in the correct functor category,
we are going to assign to the refinement a left exact functor $\Zprime(\Sigma;\Sigma\refi)$ 
that is obtained from the block functor $\Zfine(\Sigma\refi)$ for the refined surface 
$\Sigma\refi$ by evaluating it on a distinguished object in the gluing category for the part
$\partialr\Sigma\refi$ of the gluing boundary of $\Sigma\refi$ that is not inherited from 
$\Sigma$. We call this (still to be defined) distinguished object in $\ZZ(\partialr\Sigma\refi)$
the \emph{\mute} object of the refinement $(\Sigma;\Sigma\refi)$ and denote it by 
$\om_{\Sigma;\Sigma\refi}$, and call $\Zprime(\Sigma;\Sigma\refi)$ the \emph{\csfun} for the 
fine refinement $(\Sigma;\Sigma\refi)$. Thus we set
  \be
  \Zprime(\Sigma;\Sigma\refi)(-) := \ZZ(\Sigma\refi) (- \boti \om_{\Sigma;\Sigma\refi})
  :\quad \ZZ(\partialg\Sigma\refi{\setminus}\partialr\Sigma\refi) \to \vect \,.
  \label{eq:def:Zprime}
  \ee
Further, we can identify the complement of $\partialr\Sigma\refi$ in $\partialg\Sigma\refi$
with the gluing boundary of $\Sigma$, whereby the \csfun\ becomes a left exact functor 
\be
\label{eq:lexZprime}
  \Zprime(\Sigma;\Sigma\refi) :\quad \ZZ(\partialg\Sigma) \to \vect \,,
  \ee
and thus an object in the same functor category in which also the block functor for $\Sigma$ 
should be an object.
Of course, $\Zprime(\Sigma;\Sigma\refi)$ depends on both $\Sigma$ and $\Sigma\refi$, whereas
the block functor $\ZZ(\Sigma)$ must only depend on $\Sigma$. We will therefore define 
$\ZZ(\Sigma)$ as a limit of the collection $\{ \Zprime(\Sigma;\Sigma\refi) \}$ of functors
over the refinements that refine the surface $\Sigma$ which, as we will show,
can be endowed with the structure of a diagram with values in the category of left exact
functors. Addressing this issue, it still remains to specify the \mute\ 
object $\om_{\Sigma;\Sigma\refi}$ for any refinement. Since according to Remark
\ref{refi-vs-transdisk}(i) every component of $\partialr\Sigma\refi$ is fillable by a disk, it 
is sufficient to define a \mute\ object $\om(\LL) \,{\in}\, \TT(\LL)$ for every fillable 
gluing circle or gluing interval $\LL$. Given those objects, we set
  \be
  \om_{\Sigma;\Sigma\refi} := \bigboxtimes_i\, \om(\LL_i) \,,
  \ee
where the Deligne product is taken over all components $\LL_i$ of $\partialr\Sigma\refi$.

To define, in turn, the objects $\om(\LL)$ we proceed in two steps. We first consider the case
that $\LL$ is the simplest possible type of a fillable circle: for $\cala$ a finite tensor
category and $\epsilon \,{\in}\, \{\pm1\}$, a circle with a single defect point that is 
transparently labeled and has orientation $\epsilon$. We denote such a circle by
$\QQ^\cala_\epsilon$ and refer to it as a \emph{\tadpolecircle}. Pictorially, 
    \def\locpa  {0.75}   
    \def\locpb  {0.81}   
    \def\locpc  {2.2pt}  
    \def\locpd  {0.12}   
  \be
  \raisebox{-3.9em}{\begin{tikzpicture}
  \scopeArrow{0.11}{<}
  \filldraw[line width=\locpc,\colorCircle,fill=white,postaction={decorate}]
           (0,0) circle (\locpa) node[right=0.66cm,color=black]{$2$};
  \end{scope}
  \scopeArrow{0.81}{\arrowDefect}
  \draw[line width=\locpc,color=\colorTransp,postaction={decorate}] (0,\locpa) -- +(0,\locpb)
           node[above=-2pt,xshift=1pt,color=\colorDefect] {$\II_0$};
  \end{scope}
  \filldraw[color=\colorDefect] (0,\locpa) circle (\locpd);
  \node at (-1.9,0.5) {$\QQ^{\cala}_+ ~=$};
  \begin{scope}[shift={(6.1,1.0)}]
  \scopeArrow{0.11}{<}
  \filldraw[line width=\locpc,\colorCircle,fill=white,postaction={decorate}]
           (0,0) circle (\locpa) node[right=0.66cm,color=black]{$2$};
  \end{scope}
  \scopeArrow{0.71}{\arrowDefectB}
  \draw[line width=\locpc,color=\colorTransp,postaction={decorate}] (0,-\locpa) -- +(0,-\locpb)
           node[below=-1pt,xshift=1pt,color=\colorDefect] {$\II_0$};
  \end{scope}
  \filldraw[color=\colorDefect] (0,-\locpa) circle (\locpd);
  \node at (-1.9,-0.5) {$\QQ^{\cala}_- ~=$};
  \end{scope}
  \label{eq:pic:tadpolecirc}
  \end{tikzpicture}}
  \ee

The relevant \mute\ objects are obtained via a remarkable interplay
between algebraic structures in finite tensor categories and the geometry of framings.
The gluing categories for \tadpolecircle s are canonically equivalent to (twisted) centers:
  \be
  \ZZ(\QQ^{\cala}_+) = \II_0\, \Boti{2} \, \cong \cent(\cala) \qquad \text{and} \qquad 
  \ZZ(\QQ^{\cala}_-) = \cc{\II_0}\, \Boti{2} \, \cong \cc{\cent^{-4}_{}(\cala)}
  \ee
(recall Definition \ref{def:Znkappa} of the twisted center $\cent_{}^\kappa$). Thus in particular
they contain canonical objects, namely the monoidal unit in $\cent(\cala)$ and the distinguished
invertible object (see page \pageref{eq:end-coend-D}) in $\cc{\cent^{-4}(\cala)}$,
respectively. We define the respective \mute\ objects to be these two specific objects:
  \be
  \om(\QQ^{\cala}_+) := \one ~\in \cent(\cala) \,{\cong}\, \ZZ(\QQ^{\cala}_+)
  \qquad \text{and} \qquad
  \om(\QQ^{\cala}_-) := \cc{\DA} ~\in \cc{\cent^{-4}(\cala)} \,{\cong}\, \ZZ(\QQ^{\cala}_-) \,.
  \label{eq:tadpole-mute}
  \ee

In the second step we consider an arbitrary fillable gluing circle or interval $\LL$. There
exists (uniquely up to isomorphism) a \defectsurface\ $\DDtad_\LL$ which is a fillable disk (in
the sense of Definition \ref{def:filldisk}) such that $\LL$ is its outer boundary and each of
its inner boundary circles is a \tadpolecircle\ $\QQ_i$ (for a suitable finite tensor category)
whose single defect point is connected by a single (transparent) defect line with $\LL$.
The following picture shows an example of such a \defectsurface\ $\DDtad_\LL$
(involving a left $\cala$-module $\caln$ and an $\cala$-$\calb$-bimodule $\calx$):
 \\[-2.2em] 
    \def\locpa  {2.7}    
    \def\locpA  {1.5}    
    \def\locpc  {2.2pt}  
    \def\locpC  {1.8pt}  
    \def\locpd  {0.12}   
    \def\locpD  {0.08}   
    \def\locpe  {0.09}   
    \def\locpf  {1.1pt}  
    \def\locpr  {1.21}   
    \def\locps  {1.38}   
    \def\locpt  {0.43}   
    \def\locpu  {144}    
    \def\locpv  {65}     
    \def\locpw  {225}    
    \def\locpx  {172}    
    \def\locpX  {331}    
    \def\locpy  {-6}     
    \def\locpY  {40}     
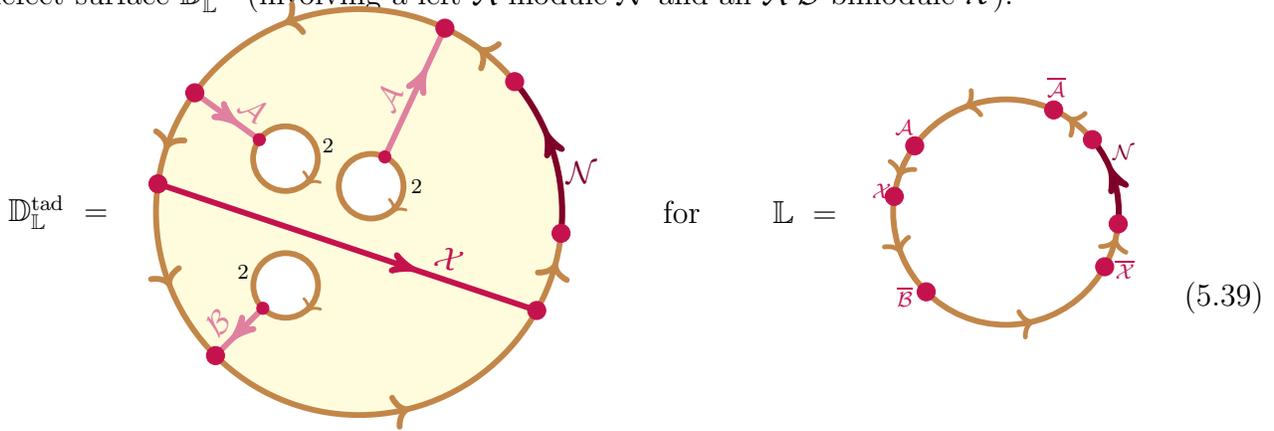
  \begin{eqnarray}
  \begin{tikzpicture}
  \scopeArrow{0.15}{>} \scopeArrow{0.31}{>} \scopeArrow{0.45}{>} \scopeArrow{0.56}{>}
  \scopeArrow{0.79}{>} \scopeArrow{0.96}{>}
  \filldraw[\colorCircle,fill=\colorSigma,line width=\locpc,postaction={decorate}]
             (0,0) circle (\locpa);
  \end{scope} \end{scope} \end{scope} \end{scope} \end{scope} \end{scope}
  \scopeArrow{0.74}{\arrowDefect}
  \draw[line width=\locpc,color=\colorTransp,postaction={decorate}]
             (\locpv:0.38) -- node[sloped,above,yshift=-1pt]{$\cala$} (\locpv:\locpa);
  \end{scope}
  \scopeArrow{0.47}{\arrowDefect}
  \draw[line width=\locpc,color=\colorTransp,postaction={decorate}]
             (\locpu:\locpa) -- node[sloped,above,yshift=-1pt]{$\cala$} (\locpu:\locpr);
  \end{scope}
  \scopeArrow{0.77}{\arrowDefect}
  \draw[line width=\locpc,color=\colorTransp,postaction={decorate}]
             (\locpw:\locps) -- node[sloped,above,near end, yshift=-1pt]{$\calb$} (\locpw:\locpa);
  \end{scope}
  \scopeArrow{0.64}{\arrowDefect}
  \draw[line width=\locpc,color=\colorFreebdy,postaction={decorate}] (\locpy:\locpa) 
             node[right=-3pt,yshift=23pt] {$\caln$} arc (\locpy:\locpY:\locpa);
  \end{scope}
  \filldraw[\colorCircle,fill=white,line width=\locpc]
	     (\locpv:0.38) circle (\locpt) node[right=10pt,color=black] {$\scriptstyle 2$}
             (\locpu:\locpr) circle (\locpt) node[right=9pt,color=black,yshift=5pt] {$\scriptstyle 2$}
             (\locpw:\locps) circle (\locpt) node[left=9pt,color=black,yshift=5pt] {$\scriptstyle 2$};
  \scopeArrow{0.91}{<}
  \draw[\colorCircle,line width=\locpf,postaction={decorate}] (\locpu:\locpr) circle (\locpt);
  \draw[\colorCircle,line width=\locpf,postaction={decorate}] (\locpw:\locps) circle (\locpt);
  \draw[\colorCircle,line width=\locpf,postaction={decorate}] (\locpv:0.38) circle (\locpt);
  \end{scope}
  \scopeArrow{0.68}{\arrowDefect}
  \draw[line width=\locpc,color=\colorDefect,postaction={decorate}]
             (\locpx:\locpa) -- node[sloped,above,near end,yshift=-1pt]{$\calx$} (\locpX:\locpa);
  \end{scope}
  \filldraw[color=\colorDefect] (\locpv:\locpa) circle (\locpd) (\locpv:\locpt+0.38) circle (\locpD)
             (\locpu:\locpa) circle (\locpd) (\locpu:\locpr+\locpt) circle (\locpD)
             (\locpw:\locpa) circle (\locpd) (\locpw:\locps+\locpt) circle (\locpD)
             (\locpx:\locpa) circle (\locpd) (\locpX:\locpa) circle (\locpd)
             (\locpy:\locpa) circle (\locpd) (\locpY:\locpa) circle (\locpd);
  \node at (-\locpa-1.3,0) {$\DDtad_\LL ~=$};
  \node at (\locpa+2.5,0) {for $\qquad \LL ~=$};
  \begin{scope}[shift={(\locpa+\locpA+4.4,0)}]
  \draw[\colorCircle,line width=\locpc] (0,0) circle (\locpA);
  \scopeArrow{0.158}{>} \scopeArrow{0.31}{>} \scopeArrow{0.45}{>} \scopeArrow{0.56}{>}
  \scopeArrow{0.79}{>} \scopeArrow{0.964}{>}
  \draw[\colorCircle,line width=\locpC,postaction={decorate}] (0,0) circle (\locpA);
  \end{scope} \end{scope} \end{scope} \end{scope} \end{scope} \end{scope}
  \scopeArrow{0.64}{\arrowDefect}
  \draw[line width=\locpc,color=\colorFreebdy,postaction={decorate}]
             node[right=-7pt,yshift=27pt] {$\scriptstyle\caln$} arc (\locpy:\locpY:\locpA);
  \end{scope}
  \filldraw[color=\colorDefect]
             (\locpu:\locpA) circle (\locpd) node[above,xshift=-4pt] {$\scriptstyle \cala$} 
	     (\locpv:\locpA) circle (\locpd) node[above=1pt,xshift=1pt] {$\scriptstyle \CC\cala$}
             (\locpw:\locpA) circle (\locpd) node[left=1pt,yshift=-2pt] {$\scriptstyle \CC\calb$}
             (\locpy:\locpA) circle (\locpd) (\locpY:\locpA) circle (\locpd)
             (\locpx:\locpA) circle (\locpd) node[left=-3pt,yshift=2pt] {$\scriptstyle \calx$}
             (\locpX:\locpA) circle (\locpd) node[right,yshift=-1pt] {$\scriptstyle \CC\calx$};
  \end{scope}
  \end{tikzpicture}
  \label{eq:exa:DDtadLL}
  \nonumber\\[-3.8em]~\\[0.5em]~\nonumber
  \end{eqnarray}
We regard $\DDtad_\LL$ as a bordism
$\DDtad_\LL \colon \bigsqcup_{i=1}^n \QQ_i \,{\xrightarrow{~~}}\, \LL$ and define the 
\mute\ object for $\LL$ to be the image
  \be
  \om(\LL) := \ZZ(\DDtad_\LL)(\om(\QQ_1) \boti \cdots \boti \om(\QQ_n)) ~ \in \ZZ(\LL) 
  \label{eq:def:om(LL)}
  \ee
of the \mute\ objects for the inner boundary circles under the block functor for the fine
\defectsurface\ $\DDtad_\LL$.

To summarize: For any fine refinement $(\Sigma;\Sigma\refi)$ the prescription
\eqref{eq:def:Zprime} provides us with a functor $\Zprime(\Sigma;\Sigma\refi)$ in the 
category $\FunleSigmavect$ in which also the pre-block functor $\Zpre(\Sigma)$ is an object.
By construction, the \csfun\ $\Zprime(\Sigma;\Sigma\refi)$ depends both on $\Sigma\refi$ and on 
$\Sigma$. Now recall from Section \ref{sec:labeling} that we denote by $\FUNle$ the symmetric 
monoidal bicategory of finite \ko-linear categories, left exact functors and natural transformations.
We can thus collect the data of the \csfun\ by considering the following 
symmetric monoidal bicategory $\widetilde{\Gamma}$: Objects of $\widetilde{\Gamma}$
are \defectonemanifold s.  For \defectonemanifold s $\LL$ and $\LL'$, a $1$-morphism
$(\Sigma;\Sigma')$ from $\LL$ to $\LL'$ is a pair consisting of a $1$-morphism
$\Sigma \colon \LL\,{\to}\,\LL'$ in $\Bfrdefdec$, i.e.\ a \defectsurface , and of
a fine \defectsurface\ $\Sigma'$ that refines $\Sigma$. The  composition is as in
$\Bfrdefdec$, the only 2-morphisms are identity morphisms, and the symmetric 
monoidal structure is disjoint union.
Then from Equation \eqref{eq:lexZprime} and Theorem \ref{thm:T1T2=Ttr} we see that
  \be
  \Zprime :\quad \widetilde{\Gamma} \xrightarrow{~~} \FUNle
  \label{eq:Zprime-bifunctor}
  \ee
is a symmetric monoidal functor, since by construction the disjoint union of objects and 
1-morphisms is mapped to the Deligne product of the corresponding objects in   $\FUNle$. 
   
We want the block functor of $\Sigma$ ultimately to be a functor in $\FunleSigmavect$ 
as well. Accordingly we add more 2-morphisms to $\widetilde{\Gamma}$:
     
\begin{Definition}\label{def:Gamma-Pi}
\Itemizeiii

\item[{\rm (i)}]
The \emph{bicategory} $\Gamma$ \emph{of fine refinements} of \defectsurface s is the 
following symmetric monoidal bicategory: Objects of $\,\Gamma$ are \defectonemanifold s.
For \defectonemanifold s $\LL$ and $\LL'$, a $1$-morphism
$(\Sigma;\Sigma')$ from $\LL$ to $\LL'$ is a pair consisting of a $1$-morphism
$\Sigma \colon \LL\,{\to}\,\LL'$ in $\Bfrdefdec$, i.e.\ a \defectsurface , and of
a fine \defectsurface\ $\Sigma'$ that refines $\Sigma$.
The only $2$-mor\-phisms of $\,\Gamma$ are one single isomorphism 
$(\Sigma;\Sigma')\,{\xrightarrow{~\cong~}}\,(\Sigma;\Sigma'')$ 
for each pair of fine refinements that refine the same \defectsurface\ $\Sigma$. 
The composition of $1$-morphisms is given by the gluing of surfaces and of refinements.
The symmetric monoidal structure of $\,\Gamma$ is given by the disjoint union of
objects and of $1$-morphisms. Since the $2$-morphisms are either $1$-element sets or
empty there is a unique way to define the symmetric monoidal product of $2$-morphisms.

\item[{\rm (ii)}]
A \emph{parallelization} $\Pi$ for the collection of \csfun s for all fine refinements
is a symmetric monoidal bifunctor
  \be
  \Pi:\quad \Gamma \xrightarrow{~~} \FUNle
  \ee
that maps objects and $1$-morphisms according to $\eqref{eq:Zprime-bifunctor}$, i.e.\
$\Pi((\Sigma;\Sigma')) \,{:=}\, \Zprime(\Sigma;\Sigma')$.
\end{itemize}
\end{Definition}

Spelled out, a parallelization $\Pi$  is a collection of natural isomorphisms
  \be
  \Pi_{\Sigma;\Sigma',\Sigma''}^{}: \quad
  \Zprime(\Sigma;\Sigma') \Longrightarrow \Zprime(\Sigma;\Sigma'')
  \label{eq:spelledout}
  \ee
of functors that is coherent in the sense that every diagram involving three refinements
commutes, $\Pi_{\Sigma;\Sigma_2,\Sigma_3}^{} \,{\veco}\,\Pi_{\Sigma;\Sigma_1,\Sigma_2}^{}
\,{=}\, \Pi_{\Sigma;\Sigma_1,\Sigma_3}^{}$ 
    and that is compatible with factorization in the
following sense: if $\Sigma \,{=}\, \Sigma' \,{\circ}\, \Sigma''$ is a \defectsurface\
obtained by gluing and we are given refinements $(\Sigma';\Sigma'_i)$ and
$(\Sigma'';\Sigma''_{i})$, for $i\,{=}\,1,2$, that start at $\Sigma'$ and at $\Sigma''$,
respectively, then gluing the refined surfaces results in refinements $(\Sigma;\Sigma_{i})$
that start at $\Sigma$, and the associated isomorphisms satisfy
  \be
  \Pi_{\Sigma;\Sigma_1^{},\Sigma_2^{}}^{}
  = \Pi_{\Sigma';\Sigma'_1,\Sigma'_2}^{} \hoco \Pi_{\Sigma'';\Sigma''_1,\Sigma''_2}^{} \,.
  \label{eq:facto-pi}
  \ee

Our goal is now to define the block functor $\ZZ(\Sigma)$ as the limit
of a parallelization of $\Sigma$. For this to make sense, we have to show that for any
\defectsurface\ $\Sigma$ at least one parallelization exists.
The rest of this subsection is devoted to the construction of such a parallelization.

\begin{Remarks}
\Itemizeiii

\item[(i)]
The terminology `parallelization' is chosen because this structure may be thought of 
as a bicategorical analogue of the parallelization of a tangent bundle.

\item[(ii)]
A limit construction similar to ours is also used in the standard Turaev--Viro 
construction. In that case, parallelizations are provided by the unique natural 
transformations assigned to 3-manifolds that are cylinders over 2-manifolds (see e.g.\
\cite{balKi}). In our framework, in which we do not have a three-dimensional topological
field theory at our disposal, we are instead going to construct parallelizations 
in a purely two-dimensional setting.

\item[(iii)]
For any diagram $\alpha\colon \Xi \,{\to}\, \calc$ in which every morphism
$\alpha(\xi{\stackrel g\to} \xi')$ is an isomorphism, the limit and colimit of $\alpha$
can be identified. (To see this, just note that for any representative $\xi_0$,
$\alpha(\xi_0)$ can be endowed with the structure of a colimit by the structure morphisms
$\alpha(\xi{\stackrel g\to} \xi_0)$, and with the structure of a limit by the structure morphisms
$\alpha(\xi_0{\stackrel g\to} \xi)$.)
Thus the block functor $\ZZ(\Sigma)$ is both a limit and a colimit
of a parallelization of the collection of \csfun s.
\end{itemize}
\end{Remarks}


\subsection{Construction of a parallelization}

We now outline the construction of a parallelization $\Pi$; the details of the arguments 
are deferred to Appendix \ref{section:Refind}. 
Recall again from Definition \ref{def:filldisk} the notion of a
\filldisk\ $\DD \,{=}\, \DD_\XX$ of type $\XX$ and its outer boundary.
We proceed in two steps. The first step addresses a local aspect of the construction,
while the second step deals with common subrefinements of \defectsurface s so as
to globalize the local construction. We start by constructing
a distinguished isomorphism between the block functors for 
any two \filldisk s of the same type. This amounts to a parallelization for the subcategory
of refinements $(\Sigma;\Sigma')$ for which $\Sigma$ is a \filldisk. To achieve this goal, fix a
\defectonemanifold\ $\LL \,{=}\, \LL_\XX$ that can appear as the outer boundary of a \filldisk\
$\DD \,{=}\, \DD_\XX$ of some type $\XX$, as introduced in Definition \ref{def:filldisk}.
We have already seen that there exists a unique
\filldisk\ $\DDtad_\XX$ with outer boundary $\LL$ such that each inner boundary circle of
$\DDtad_\XX$ is a \tadpolecircle\ (see the picture \eqref{eq:exa:DDtadLL}).
In Appendix \ref{sec:B2}\,--\,\ref{lem:Xirho=rhoXi}
we construct explicitly a canonical isomorphism
  \be
  \varphi_{\DD}^{} :\quad
  \Zfine(\DD)(\om(\DD)) \tocong \Zfine(\DDtad_\XX)(\om(\DDtad_\XX)) \,.
  \label{eq:varphiDD}
  \ee
Our construction of this isomorphism makes use of a combinatorial datum, namely a spanning
tree for the graph $\Gamma_{\!\DD}$ that has as edges those defect lines of $\DD$ that come
from $\XX$, and as vertices the end points of these defect lines on $\partial\DD$ together
with the inner boundary circles of $\DD$. But as we show in Lemma \ref{lem:varphiDD-indep}
in the appendix, the isomorphism \eqref{eq:varphiDD} does not depend on this datum.

Having obtained the isomorphisms \eqref{eq:varphiDD} we can define the desired isomorphism
between the block functors for any two \filldisk s $\DD$ and $\DD'$ of the same type $\XX$
(and thus in particular with the same \mute\ object $\om(\DD)$) as the vertical composition 
  \be
  \varphi_{\DD,\DD'} := \varphi_{\DD'}^{-1} \veco \varphi_{\DD}^{}:\quad
  \Zfine(\DD)(-\boti\om(\DD)) \to \Zfine(\DD')(-\boti\om(\DD)) \,.
  \ee
We can then further show -- see Proposition \ref{Proposition:main-disk-repl} -- that
the so obtained family of natural isomorphisms labeled by pairs of \filldisk s 
of any given type $\XX$ is coherent in the sense that 
  \be
  \varphi_{\DD,\DD''}^{} = \varphi_{\DD',\DD''}^{} \veco \varphi_{\DD,\DD'}^{} \,,
  \ee
and that it
satisfies the following \emph{factorization property}: For a \filldisk\ $\DD$ of the form
$\DD \,{=}\, \mathbb Y \,{\circ}\, (\DD_1 \,{\sqcup} \cdots {\sqcup}\, \DD_n)$, with
$\DD_1, ...\,, \DD_n$ non-intersecting \filldisk s in $\DD$ and $\mathbb Y$ the
\defectsurface\ that is obtained by removing all the $\DD_i$ from $\DD$,
for any $n$-tuple of replacements of the \filldisk s $\DD_i^{}$ by \filldisk s $\DD_i'$
of the same type $\XX_i$ and with the same outer boundary, the equality
  \be
  \varphi_{\DD,\DD'}^{} = (\varphi_{\DD_1^{},\DD_1'} \boti \cdots \boti
  \varphi_{\DD_n^{},\DD_n'}) \hoco {\Zfine(\mathbb Y)}^{}
  \ee
of natural transformations holds, where $\DD'$ is the \defectsurface\
$\DD' \,{:=}\, \mathbb Y \,{\circ}\, (\DD_1' \,{\sqcup} \cdots {\sqcup}\, \DD_n')$ and the symbol `$\circ$'
stands for the horizontal composition (whiskering) of the natural transformation
$\bigboxtimes_i \varphi_{\DD_i^{},\DD_i'}$ with the functor $\Zfine(\mathbb Y)$. (Also, here
$\mathbb Y$ is regarded as a bordism from $\partialgg \DD$ to
$\partialgg\DD_1 {\sqcup}\cdots{\sqcup} \partialgg\DD_n$, with $\partialgg$ denoting the
gluing part of the outer boundary, so that we deal with a functor
$\Zfine(\mathbb Y)\colon \ZZ(\partial \DD) \,{\to}\, \bigboxtimes_i\ZZ(\partial \DD_i)$
and thus tacitly invoke the equivalence \eqref{eq:LexMN-LexMNbvect}.)

Note that the factorization property involves the manipulation of replacing a \filldisk\
$\DD_i^{}$ inside the \defectsurface\ $\Sigma \,{=}\,\mathbb Y$ by another \filldisk\ $\DD_i'$
with the same boundary. (That such a manipulation is possible rests on the factorization result
of Theorem \ref{thm:Facto} below.) We call the corresponding manipulation
for a generic \defectsurface\ $\Sigma$ a \emph{\filldiskrep} of $\DD$ by $\DD'$ in $\Sigma$ 
and denote the resulting \defectsurface\ by $\DR(\Sigma) \,{\equiv}\, \DR_{\DD,\DD'}(\Sigma)$.
(For an example of a \filldiskrep, see the picture \eqref{eq:exa-filldiskrep}.)

By a \emph{refinement replacement} from $(\Sigma;\Sigma')$ to $(\Sigma;\Sigma'')$ we then mean
a sequence of (possibly intersecting) \filldiskrep s $(\DR_1,...\,,\DR_n)$ such that
$\DR_n( \cdots \DR_1(\Sigma') \cdots) \,{=}\, \Sigma''$. Given any two refinements
$(\Sigma;\Sigma_1)$ to $(\Sigma;\Sigma_2)$ one can construct a \emph{common subrefinement}
$(\Sigma;\Sigma_{1,2})$ and specific \emph{standard} refinement replacements from
$(\Sigma;\Sigma_i)$ for $i\,{=}\,1,2$ to $(\Sigma;\Sigma_{1,2})$ as a sequence of
\filldiskrep s of a very restricted type (compare Definition \ref{def:resolve} and the
pictures \eqref{eq:resolve} and \eqref{eq:repl:creation}). This implies in particular
(see Lemma \ref{lem:refrepl}) that a refinement replacement
exists between any two refinements that refine the same \defectsurface. Moreover, it
follows from the results about \filldisk s in Proposition \ref{Proposition:main-disk-repl} that
any refinement replacement $(\DR_1,...\,,\DR_n)$ from $(\Sigma;\Sigma')$ to $(\Sigma;\Sigma'')$
provides a distinguished natural isomorphism
  \be
  \varphi^{}_{ \Sigma', \DR_n( \cdots\, \DR_1(\Sigma') \,\cdots) }
  :\quad \Zfine(\Sigma')(-\boti\om') \xrightarrow{~~} \Zfine(\Sigma'')(-\boti\om'') \,.
  \label{eq:ref-repl-M-0}
  \ee

The second step of the construction of a parallelization uses these isomorphisms to
provide the natural isomorphisms \eqref{eq:spelledout}.
To this end we first show, with the help of a suitable notion of common subrefinement,
that for any two refinements $(\Sigma;\Sigma')$ and $(\Sigma;\Sigma'')$
a refinement replacement $(\DR_1,...\,,\DR_n)$ from $(\Sigma;\Sigma')$ to $(\Sigma;\Sigma'')$
exists (Lemma \ref{lem:refrepl}).
Then we show in Lemma \ref{Lemma:equiv-comm-sub} that for any refinement replacement 
$(\DR_1,...\,,\DR_n)$ in a fine \defectsurface\ there exists a refinement replacement
$(\DR_1',...\,,\DR_{n'}')$ that induces the same natural isomorphism and that involves only
very specific types of \filldiskrep s (which are introduced in Definition \ref{def:resolve}).
This finally allows us to conclude, in Proposition \ref{Proposition:main-generic-ref},
that any two refinement replacements $(\DR_1,...\,,\DR_n)$ and $(\DR_1',...\,,\DR_{n'}')$
between any two given refinements $(\Sigma;\Sigma_{\rm ref_1})$ and $(\Sigma;\Sigma_{\rm ref_2})$
of an arbitrary \defectsurface\ $\Sigma$ induce the same natural isomorphism,
  \be
  \varphi^{}_{ \Sigma_{\rm ref_1}, \DR_{n'}'( \cdots\, \DR_1'(\Sigma_{\rm ref_1}) \,\cdots) }
  = \varphi^{}_{ \Sigma_{\rm ref_1}, \DR_n^{}( \cdots\, \DR_1^{}(\Sigma_{\rm ref_1}) \,\cdots) }
  \,.
  \ee
Thus in the situation considered in Equation \eqref{eq:spelledout} we define a
parallelization $\Pi$ by picking a refinement replacement between two given
refinements $\Sigma'$ and $\Sigma''$ of a \defectsurface\ $\Sigma$ and setting
  \be
  \Pi_{\Sigma;\Sigma',\Sigma''}^{}:=  \varphi^{}_{ \Sigma', \DR_n( \cdots\, \DR_1(\Sigma') \,\cdots) }
  \label{eq:PI-explicit}
  \ee
as in Equation \eqref{eq:ref-repl-M-0}. By the uniqueness result of Proposition
\ref{Proposition:main-generic-ref}, it now follows directly that $\Pi$ is compatible
with the composition of 2- and 1-morphisms of $\Gamma$ in accordance with the 
coherence condition \ref{eq:facto-pi}. We summarize our findings in the

\begin{thm} \label{thm:exists-pi}
There exists a parallelization $\Pi$ for the collection of \csfun s for all fine refinements of
\defectsurface s.
\end{thm}

In fact, a parallelization not only exists, but it also satisfies a universal property
with respect to parallel transport operations, and is therefore uniquely characterized.
Before showing this universality, we have a look at consequences of the existence of $\Pi$.
We start by giving the
	
\begin{Definition}\label{def:THEdef}
The block functor $\ZZ(\Sigma)$ for a $($not necessarily fine$)$ \defectsurface\ $\Sigma$ is the
limit 
  \be
  \ZZ(\Sigma) := \lim_{\Sigma'} \Pi_{\Sigma;\Sigma'}
  \ee
$($or, equivalently, colimit$)$ of the parallelization of $\Sigma$ described above.
\end{Definition}

\begin{Remarks}
\Itemizeiii

\item[(i)]
The proof of Proposition \ref{Proposition:main-generic-ref} in Appendix \ref{section:Refind}
relies on the details of the parallelization only through Proposition 
\ref{Proposition:main-disk-repl} which covers the case of \filldisk s. This shows that the
datum of a modular functor can be obtained from a functor $\Bfrdefdecgfine \To \Funle$, 
that is, from a functor that is defined on \emph{fine} defect surfaces,
and that is equipped with a family of `transparent' defects such that
Proposition \ref{Proposition:main-disk-repl} holds.

\item[(ii)]
Thus the crucial input for the equivalence of refinement replacements
stated in Proposition \ref{Proposition:main-generic-ref} is the corresponding assertion
for disks in Proposition \ref{Proposition:main-disk-repl}. This aspect of our
construction is reminiscent of the way \Cite{Sect.\,5.5}{luri7}
in which factorization homology allows one to extend locally constant factorization algebras
from open disks to arbitrary surfaces. Accordingly one may suspect that
Proposition \ref{Proposition:main-disk-repl} is related to the structure
of a factorization algebra on disks that is locally constant with respect to the
stratification given by the defects on fillable disks. A detailed analysis of this
idea is beyond the scope of the present paper. But it is reassuring that
analogs of defects can be treated in the setting of locally constant factorization 
algebras for stratified manifolds \Cite{Sect.\,6}{ginot} and that 
factorization homology can be extended to stratifications \cite{ayfT3}.
Moreover, transparent defects (in our language) have appeared in this
context as well, e.g.\ as a tool for defining a topological field theory 
associated with an $E_n$-algebra \cite{scheC2}.
\end{itemize} 
\end{Remarks}

It follows directly from Proposition \ref{Proposition:Equivalence-disks} that 
as a particular case of the parallelization we have

\begin{Lemma} \label{Lemma:excision-pi}
Let $(\Sigma;\Sigma\refi)$ be a refinement of a \defectsurface\ $\Sigma$ such that the 
defect lines in $\Sigma\refi{\setminus}\Sigma$ together with $\partialr\Sigma\refi$
form a tree $\gamma_{\Sigma;\Sigma\refi}$ in $\Sigma$. Then the excision isomorphism from 
Lemma $\ref{Lemma:Z-same}$ applied to $\gamma_{\Sigma;\Sigma\refi}$ provides an isomorphism
  \be
  \Zpre( \Sigma\refi ) (- \boti \om_{\Sigma;\Sigma\refi}) \cong \Zpre(\Sigma)
  \label{eq:tree-iso}
  \ee
already for the pre-block functors, and this isomorphism 
induces the parallelization isomorphism for the \csfun s.
\end{Lemma}

Let us now come back to the factorization issue considered in
Section \ref{sec:facto}, which so far could be discussed only for \emph{fine} surfaces. 
Our results on refinements allow us to extend the results of Section \ref{sec:facto}
directly to the factorization of general \defectsurface s.

\begin{thm} \label{thm:Facto}
\emph{Factorization}:
Let $\Sigma$ be a \defectsurface\ with gluing boundary components $\LL$ and $\cc{\LL}$. There
is a canonical isomorphism 
  \be
  \int^{x \in \ZZ(\LL)} \! \ZZ(\Sigma)({-} \boti x \boti \cc x)
  \,\cong\, \ZZ(\Tr_{\LL}(\Sigma))(-) 
  \label{eq:thm:fact-fine}
  \ee
of left exact functors. Moreover, the canonical isomorphisms obtained this way with any two distinct 
pairs $(\LL_1,\cc{\LL_1}) $ and  $(\LL_2,\cc{\LL_2})$ of gluing boundaries of $\Sigma$ commute. 
\end{thm}

\begin{proof}
For notational simplicity we restrict our attention to the case that 
$\Sigma \,{=}\, \Sigma_{1} \,{\sqcup}\, \Sigma_{2}$ is the disjoint union of two \defectsurface s
$\Sigma_{1}$ and $\Sigma_{2}$ such that $\Tr_{\LL}(\Sigma) \,{=}\,\Sigma_2 \,{\circ}\, \Sigma_1$
is the gluing of $\Sigma_1$ and $\Sigma_2$; the general case is covered by the same arguments. 
Let $(\Sigma_1^{};\Sigma_1')$ and $(\Sigma_2^{};\Sigma_2')$ be fine refinements that refine
$\Sigma_{1}$ and $\Sigma_{2}$, respectively. Then 
$(\Sigma_1 \,{\circ}\, \Sigma_2;\Sigma_1' \,{\circ}\, \Sigma_2')$ is a fine refinement that
refines $\Sigma_{1} \,{\circ}\, \Sigma_{2}$. Theorem \ref{thm:T1T2=Ttr} 
thus provides us with an isomorphism 
  \be
  \widetilde\varphi_{\Sigma_2',\Sigma_1'}^{}:\quad
  \Zprime(\Sigma_2^{};\Sigma_2') \circ \Zprime(\Sigma_1^{}; \Sigma_1')
  \xrightarrow{\,\cong\,} \Zprime(\Sigma_2^{}\,{\circ}\,\Sigma_1^{}; \Sigma_2'\,{\circ}\,\Sigma_1') \,.
  \label{eq:fact-fine}
  \ee
By taking the limit over the fine refinements that refine $\Sigma_1$ and $\Sigma_2$
we then also obtain an isomorphism
  \be
  \varphi_{\Sigma_2',\Sigma_1'}^{}:\quad \ZZ(\Sigma_2^{}) \circ \ZZ(\Sigma_1^{})
  \xrightarrow{\,\cong\,}
  \Zprime(\Sigma_2^{}\,{\circ}\,\Sigma_1^{}; \Sigma_2' \,{\circ}\, \Sigma_1') \,.
  \label{eq:fact-fineII}
  \ee
Now not every refinement $(\Sigma_2\,{\circ}\,\Sigma_1;\Sigma')$ of $\Sigma_2\,{\circ}\,\Sigma_1$
is of the form $(\Sigma_2^{}\,{\circ}\,\Sigma_1^{}; \Sigma_2'\,{\circ}\,\Sigma_1')$, since
$\Sigma'$ might also refine the gluing boundary $\LL$. However, we can compose
$ \varphi_{\Sigma_2',\Sigma_1'}$ with  the parallelization
$\Pi_{\Sigma_2^{}\circ\Sigma_1^{}; \Sigma', \Sigma_2'\circ\Sigma_1'}$ of
$\Sigma_2 \,{\circ}\, \Sigma_1$ so as to obtain a coherent family of isomorphisms 
  \be
  \varphi_{\Sigma'}^{} :\quad \ZZ(\Sigma_2) \circ \ZZ(\Sigma_1)
  \xrightarrow{\,\cong\,} \Zprime(\Sigma_2 \,{\circ}\, \Sigma_1; \Sigma') \,, 
  \label{eq:fact-fineIIi}
  \ee
which does not any longer depend on the choice of $\Sigma_{1}'$ and $\Sigma_{2}'$,
owing to the coherence property \eqref{eq:facto-pi} of $\Pi$. 
Thus by the universal property of the limit we get an isomorphism 
  \be
  \ZZ(\Sigma_2) \,{\circ}\, \ZZ(\Sigma_1) \xrightarrow{~\cong~}
  \ZZ(\Sigma_2{\circ}\Sigma_1) \,.
  \ee
Note that the block functors are a priori functors to $\vect$; thus in \eqref{eq:fact-fineIIi}
we implicitly use the Eilenberg--Watts calculus to turn the functors $\ZZ(\Sigma_i)$ into
functors from the gluing category of the incoming to the one of the outgoing boundary, 
compare the discussion around \eqref{eq:LexMN-LexMNbvect}. The desired isomorphism
\eqref{eq:thm:fact-fine} then follows from \eqref{eq:fact-fineIIi} by invoking the 
compatibility of the Eilenberg--Watts calculus with the composition of functors, as
expressed by the isomorphism \eqref{eq:EW-xxx}.
 \\[.2em]
Moreover, given three composable \defectsurface s $\Sigma_1, \Sigma_2, \Sigma_3$ and
three corresponding fine refinements $(\Sigma_1^{};\Sigma_1')$, $(\Sigma_2^{};\Sigma_2')$ and 
$(\Sigma_3^{};\Sigma_3')$ it follows from Corollary \ref{cor:locality-Fact} that the order 
of the factorizations does not matter. By standard arguments this property passes to the limit. 
\end{proof}

We have thus defined the behavior of block functors under the horizontal composition 
of 1-morphisms in $\Bfrdefdec$, i.e.\ a factorization structure of the modular functor.

\begin{Remark} \label{rem:TV}
Let us comment on the relation of our construction to the standard Turaev--Viro construction.
The latter takes as an input a fusion category $\cala$ with a spherical structure, i.e.\ a pivotal
structure such that left and right traces coincide. In its standard incarnation 
\cite{bawe2,tuVi,balKi}, the Turaev--Viro construction defines a symmetric monoidal 2-functor
   \be
   \mathrm{Bord}_{3,2,1}\xrightarrow{~~} 2\mbox-\vect \,,
   \ee
where $2\mbox-\vect$ is the symmetric monoidal bicategory of finitely
semisimple $\complex$-linear categories. $\mathrm{Bord}_{3,2,1}$ is a bicategory
of bordisms, usually without defects.
The manifolds in this approach are \emph{oriented}, rather than framed.
  \\
Since we do not assign natural transformations to arbitrary three-manifolds with corners,
a comparison with the Turaev--Viro construction can only be made at the level of
the gluing categories, the functors for two-manifolds with boundary, and the representations of
mapping class groups. Concerning the level of categories, we remark that the
pivotal structure allows one to canonically identify all twisted centers
with the Drinfeld center $\cent(\cala)$; this is indeed the category associated
to a circle in the standard Turaev--Viro construction.
Similarly, in the presence of defects all twisted balancings (as introduced in
Definition \ref{def:framedcenter}) get identified and the $\kappa$-framed centers
only depend on the orientation of the defect lines.
 \\[.2em]
To compare the functors assigned to surfaces without (visible) defects, we first
sketch how to assign in our setting a functor to an oriented surface $\Sigma$, possibly
with boundary. Again we must choose auxiliary data:
a refinement $(\Sigma;\Sigma')$ of $\Sigma$ as a surface together with a framing on $\Sigma'$
(it is not hard to convince oneself that this can indeed be found), and then a refinement of 
the resulting framed surface. Given these data we get fine block functors by the prescription
in Section \ref{sec:preblock}. The relevant index category $\widetilde\Gamma$ has the same
objects as $\Gamma$, but more 1-morphisms than the one in Definition \ref{def:Gamma-Pi}, since
now framings are included as auxiliary data. However, $\widetilde\Gamma$ also has more 2-morphisms
than $\Gamma$, because it is defined to have precisely one morphisms between any two choices of
framings. Thus the bicategories $\widetilde\Gamma$ and $\Gamma$ are equivalent by construction. 
 \\[.2em]
Next one observes that
any two framings with the same underlying orientation differ by a suitable application of
double duals. The pivotal structure should therefore give us enough natural transformations
to construct a parallelization also for the index category $\widetilde\Gamma$. 
The blocks for the oriented bordism category can then again be defined, as in Definition 
\ref{def:THEdef}, as (co)limits. Let us call them the \emph{oriented blocks}.
This prescription has the additional advantage of exhibiting the oriented blocks as 
(co)limits over the framed blocks. Since the orientation amounts to a tangential
structure, this is in line with general expectations within the framework of the 
cobordism hypthesis, compare \cite[Thm.\,2.4.26]{luri4}.  
 \\[.2em]
By picking a framing we can represent an oriented block by the block obtained in our approach.
(This is completely analogous to representing a block using a specific refinement.)
Thus in particular we know that the oriented blocks obey factorization. As a consequence
it suffices to compare the oriented conformal blocks with the standard Turaev--Viro blocks 
for the case of a three-punctured sphere. The standard Turaev--Viro blocks (see e.g.\  
\Cite{Ex.\,8.6}{balKi}) are invariants; that the same is true for our oriented blocks 
is a consequence of Proposition \ref{prop:modulenattrafo}.
\end{Remark}

\begin{Example} \label{exa:ENOM}
Consider for an $\cala$-$\calb$-bimodule ${}_{\cala}\caln_{\calb}$ the following 
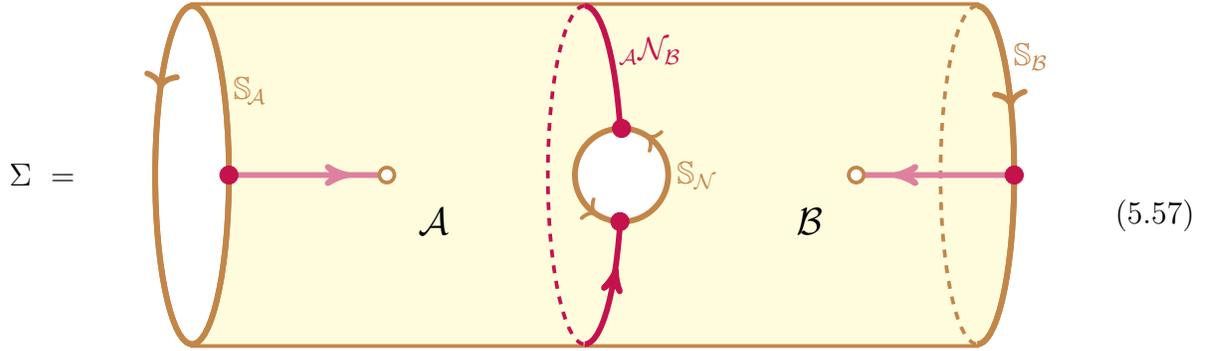
\defectsurface\ $\Sigma$:
    \def\locpa  {0.62}   
    \def\locpb  {2.1}    
    \def\locpB  {4.6}    
    \def\locpc  {2.2pt}  
    \def\locpC  {1.4pt}  
    \def\locpd  {0.12}   
    \def\locpD  {1.4pt}  
    \def\locpx  {0.49}   
    \def\locpy  {2.25}   
    \def\locpY  {2.275}  
  \be
  \raisebox{-5.4em} {\begin{tikzpicture}
  \filldraw[fill=\colorSigma,line width=\locpC,draw=\colorCircle]
           (-\locpa-\locpB-\locpx,-\locpY) rectangle (\locpa+\locpB-\locpx,\locpY)
           node[above=4pt] {~};  
  \scopeArrow{0.83}{\arrowDefectB}
  \draw[line width=\locpc,color=\colorDefect,postaction={decorate}]
           (-\locpx,\locpy) node[right=8pt,yshift=-17pt]{${}_\cala\caln_\calb$}
	   arc (90:-90:\locpx cm and \locpy cm);
  \end{scope}
  \draw[dashed,line width=\locpD,color=\colorDefect,postaction={decorate}]
           (-\locpx,\locpy) arc (90:270:\locpx cm and \locpy cm);
  \fill[\colorSigma] (\locpa+\locpB-\locpx,0) circle (\locpx cm and \locpy cm);
  \filldraw[fill=white,line width=\locpc,draw=\colorCircle] (0,0) circle (\locpa);
  \scopeArrow{0.16}{>} \scopeArrow{0.66}{>}
  \draw[line width=\locpC,\colorCircle,postaction={decorate}] (0,0) circle (\locpa);
  \end{scope} \end{scope}
  \filldraw[fill=white,line width=\locpc,draw=\colorCircle]
           (-\locpa-\locpB-\locpx,0) circle (\locpx cm and \locpy cm);
  \scopeArrow{0.38}{>}
  \draw[line width=\locpc,\colorCircle,postaction={decorate}]
           (-\locpa-\locpB-\locpx,0) circle (\locpx cm and \locpy cm);
  \end{scope}
  \filldraw[fill=\colorSigma,line width=\locpc,draw=\colorCircle]
           (\locpa+\locpB-\locpx,\locpy) arc (90:-90:\locpx cm and \locpy cm);
  \scopeArrow{0.31}{>}
  \draw[line width=\locpc,\colorCircle,postaction={decorate}]
           (\locpa+\locpB-\locpx,\locpy) arc (90:-90:\locpx cm and \locpy cm);
  \end{scope}
  \draw[dashed,line width=\locpD,draw=\colorCircle]
           (\locpa+\locpB-\locpx,\locpy) arc (90:270:\locpx cm and \locpy cm);
  \scopeArrow{0.77}{\arrowDefect}
  \draw[line width=\locpc,color=\colorTransp,postaction={decorate}]
           (-\locpa-\locpB,0) -- +(\locpb,0);
  \draw[line width=\locpc,color=\colorTransp,postaction={decorate}]
           (\locpa+\locpB,0) -- +(-\locpb,0);
  \end{scope}
  \transv{\locpb-\locpa-\locpB,0};
  \transv{\locpa+\locpB-\locpb,0};
  \filldraw[color=\colorDefect] (0,\locpa) circle (\locpd) (268:\locpa) circle (\locpd)
           (\locpa+\locpB,0) circle (\locpd) (-\locpa-\locpB,0) circle (\locpd);
  \node at (\locpb-\locpB,-\locpa) {\large $\cala$};
  \node at (\locpB-\locpb,-\locpa) {\large $\calb$};
  \node[color=\colorCircle] at (\locpx cm+\locpa cm-3pt,0) {$\SS_\caln$};
  \node[color=\colorCircle] at (-\locpB cm-9.6pt,0.5*\locpy) {$\SS_\Cala$};
  \node[color=\colorCircle] at (\locpx cm+\locpB cm+10pt,0.7*\locpy) {$\SS_\calb$};
  \node at (-\locpa-\locpB-2*\locpx-1.5,0) {$\Sigma ~=$};
  \end{tikzpicture}
  \label{eq:Fig5.8-1}}~
  \ee
  ~\\[0.2em]
We may think of $\Sigma$ as describing how the situation in the region
labeled by $\cala$ is `transmitted' to the situation in the region
labeled by $\calb$, and accordingly refer to the block functor
$\ZZ(\Sigma)$ as a \emph{transmission functor}. To calculate the
transmission functor, we suitably redraw $\Sigma$
and pick a refinement $(\Sigma;\Sigma')$ with $\Sigma'$ having two additional
transparent gluing circles (to be evaluated at their respective \mute\ objects),
as shown in the following picture:        
    \def\locpa  {0.55}   
    \def\locpA  {0.31}   
    \def\locpc  {2.2pt}  
    \def\locpC  {1.4pt}  
    \def\locpd  {0.12}   
    \def\locpD  {0.09}   
    \def\locpe  {0.9pt}  
    \def\locpx  {4.6}    
    \def\locpy  {3.4}    
  \be
  \raisebox{-8.0em}
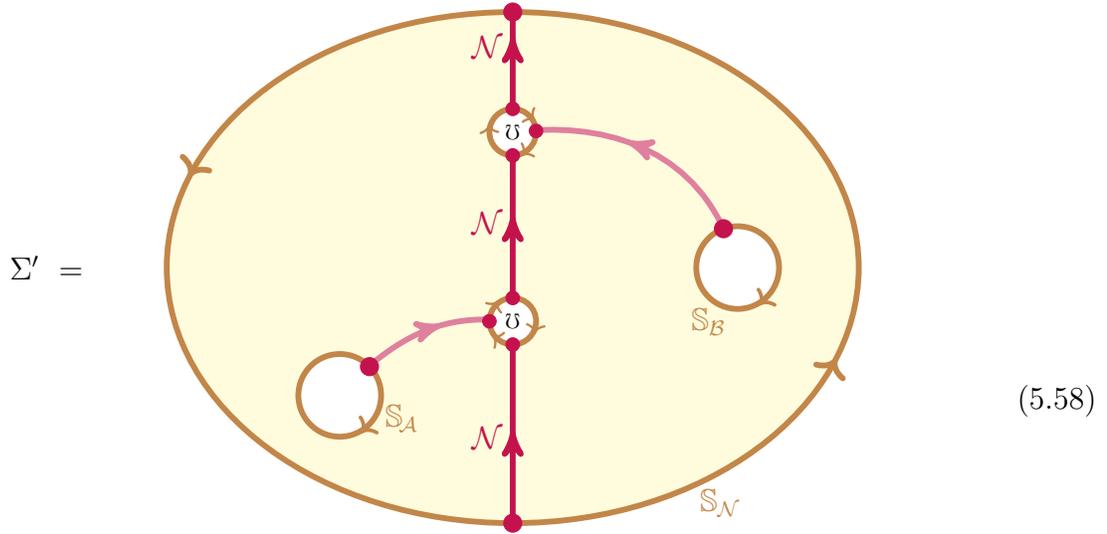
{\begin{tikzpicture}
  \scopeArrow{0.45}{>} \scopeArrow{0.95}{>}
  \filldraw[\colorCircle,fill=\colorSigma,line width=\locpc,postaction={decorate}]
           (0,0) circle (\locpx cm and \locpy cm);
  \end{scope} \end{scope}
  \scopeArrow{0.18}{\arrowDefect} \scopeArrow{0.60}{\arrowDefect} \scopeArrow{0.95}{\arrowDefect}
  \draw[line width=\locpc,color=\colorDefect,postaction={decorate}]
           (0,-\locpy) -- node[left=-1.2pt,yshift=17pt] {$\caln$} 
           node[left=-1.2pt,very near start,yshift=8pt] {$\caln$} 
           node[left=-1.2pt,very near end,yshift=11pt] {$\caln$} (0,\locpy);
  \end{scope} \end{scope} \end{scope}
  \scopeArrow{0.53}{\arrowDefect}
  \draw[line width=\locpc,color=\colorTransp,postaction={decorate}]
           (-0.5*\locpx,-0.5*\locpy)++(45:\locpa) [out=40,in=173] to (0,-0.21*\locpy);
  \draw[line width=\locpc,color=\colorTransp,postaction={decorate}]
           (0.65*\locpx,0)++(110:\locpa) [out=115,in=7] to (0,0.53*\locpy);
  \end{scope}
  \filldraw[\colorCircle,fill=white,line width=\locpc]
           (-0.5*\locpx,-0.5*\locpy) circle (\locpa) (0.65*\locpx,0) circle (\locpa)
           (0,-0.21*\locpy) circle (\locpA) (0,0.53*\locpy) circle (\locpA);
  \scopeArrow{0.88}{<}
  \draw[line width=\locpC,\colorCircle,postaction={decorate}]
           (-0.5*\locpx,-0.5*\locpy) circle (\locpa);
  \draw[line width=\locpC,\colorCircle,postaction={decorate}] (0.65*\locpx,0) circle (\locpa);
  \end{scope}
  \scopeArrow{0.41}{<} \scopeArrow{0.66}{<} \scopeArrow{0.98}{<}
  \draw[line width=\locpe,\colorCircle,postaction={decorate}] (0,-0.21*\locpy) circle (\locpA);
  \end{scope} \end{scope} \end{scope}
  \scopeArrow{0.14}{<} \scopeArrow{0.51}{<} \scopeArrow{0.87}{<}
  \draw[line width=\locpe,\colorCircle,postaction={decorate}] (0,0.53*\locpy) circle (\locpA);
  \end{scope} \end{scope} \end{scope}
  \filldraw[color=\colorDefect]
           (-\locpA,-0.21*\locpy) circle (\locpD) (0,\locpA-0.21*\locpy) circle (\locpD)
           (0,-\locpA-0.21*\locpy) circle (\locpD)
	   (0,0.53*\locpy)++(3:\locpA) circle (\locpD) (0,\locpA+0.53*\locpy) circle (\locpD)
           (0,-\locpA+0.53*\locpy) circle (\locpD)
           (0,-\locpy) circle (\locpd) (0,\locpy) circle (\locpd)
           (-0.5*\locpx,-0.5*\locpy)++(44:\locpa) circle (\locpd)
           (0.65*\locpx,0)++(110:\locpa) circle (\locpd);
  \node at (0,-0.21*\locpy) {$\scriptstyle\om$};
  \node at (0,0.53*\locpy) {$\scriptstyle\om$};
  \node[color=\colorCircle] at (0.6*\locpx,-.92*\locpy) {$\SS_\caln$};
  \node[color=\colorCircle] at (-0.5*\locpx+1.5*\locpa,-0.54*\locpy-0.3*\locpa) {$\SS_\Cala$};
  \node[color=\colorCircle] at (0.65*\locpx-0.7*\locpa,-1.3*\locpa) {$\SS_\calb$};
  \node at (-\locpx-1.6,0) {$\Sigma' ~=$};
  \end{tikzpicture}}~
  \ee
By using the isomorphisms from Lemma \ref{Lemma:excision-pi} to take care of the tadpole
circles and applying Example \ref{exa:braidedinduction} to each of the two \twocells\ 
of $\Sigma'$ we arrive at the description
  \be
  \begin{array}{rcl}
  \ZZ(\Sigma') :\quad 
  \ZZ(\SS_\Cala) \boti \ZZ(\SS_\calb) = \calz(\cala) \boti \calz(\calb)
  &\!\! \xrightarrow{~~~} \!\!& \Funle_{\cala,\calb}(\caln,\caln) = \ZZ(\SS_\caln) \,, 
  \Nxl4
  x \boti y &\!\! \xmapsto{~~~} \!\!& F_{\!y} \cir {}_{x}F
  \eear
  \ee
of the (fine) block functor of $\Sigma'$; where we use the braided induction functors 
introduced in \eqref{eq:braidedind}, i.e.\ $\ZZ(\Sigma')(x\boti y)(n)\eq x.n.y$. In 
particular, when evaluating at the identity functor $\Id \iN \Funle_{\cala,\calb}(\caln,\caln)$
in the gluing category for the outer circle $\SS_\caln$ we obtain the functor 
  \be
  \calz(\cala) \boti \calz(\calb) \to  \vect \,, \quad~
  x \boti y \mapsto \Nat_{\cala,\calb}(x.(-).y, \Id) \,.
  \label{eq:insert-id}
  \ee
Then by using the adjunctions we obtain a functor 
$\calz(\cala) \,{\to}\, \calz(\calb)^{\opp} \,{\cong}\, \calz(\calb)$. In case that the
bimodule $\caln$ is invertible, it follows from \Cite{Sect.\,5.1}{enoM} 
that the so obtained transmission functor is an equivalence between the Drinfeld centers. 
\end{Example}

\begin{Remark}
Ideally one would also like to compute the block functor for the \defectsurface\ that is
obtained from the surface \eqref{eq:Fig5.8-1} by omitting the gluing circle on the defect 
line. This requires in particular to compute the block functor for the non-fine surface
$\widehat{\Sigma}$ shown in the following picture:
    \def\locpa  {0.59}   
    \def\locpA  {2.6}    
    \def\locpb  {1.7}    
    \def\locpc  {2.2pt}  
    \def\locpd  {0.12}   
  \be
  \raisebox{-5.4em}
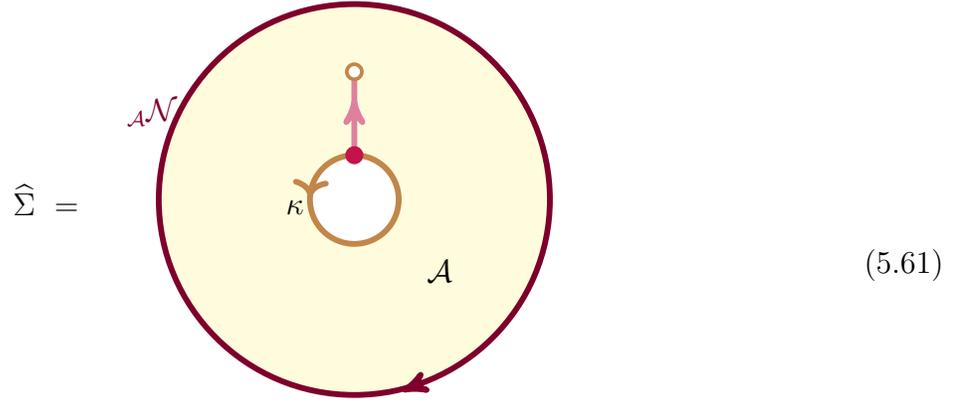
{\begin{tikzpicture}
  \scopeArrow{0.81}{\arrowDefectB}
  \filldraw[\colorFreebdy,fill=\colorSigma,line width=\locpc,postaction={decorate}]
	   (0,0) circle (\locpA) node[left=62pt,yshift=33pt]{${}_\cala\caln$};
  \end{scope}
  \scopeArrow{0.49}{>}
  \filldraw[\colorCircle,fill=white,line width=\locpc,postaction={decorate}]
           (0,0) circle (\locpa) node[left=14pt,yshift=-3pt,color=black]{$\kappa$};
  \end{scope}
  \scopeArrow{0.67}{\arrowDefect}
  \draw[line width=\locpc,\colorTransp,postaction={decorate}] (0,\locpa) -- (0,\locpb);
  \end{scope}
  \transv {0,\locpb};
  \fill[\colorDefect] (0,\locpa,0) circle (\locpd);
  \node at (-40:0.57*\locpA) {$\cala$};
  \node at (-\locpA-1.5,0) {$\widehat{\Sigma} ~=$};
  \end{tikzpicture}}
  \ee
Here the framing index at the gluing circle is necessarily  $\kappa \eq 0$, and the 
block functor is a functor 
  \be
  \ZZ(\widehat{\Sigma}) :\quad \cala \, \Boti 0 \xrightarrow{~~} \vect \,.
  \ee
Now $\cala\, {\Boti 0}$ is a \emph{twisted} Drinfeld center, with objects consisting of an
object $z \iN \cala$ together with a balancing $z\oti a \,{\cong}\, {}^{\vee \vee\!}a \oti z$
for $a \iN \cala$. Consider the refinement $(\Sigma;\widetilde{\Sigma}')$, with
$\widetilde{\Sigma}'$ the following surface:
    \def\locpa  {0.59}   
    \def\locpA  {2.6}    
    \def\locpb  {0.43}   
    \def\locpc  {2.2pt}  
    \def\locpC  {1.3pt}  
    \def\locpd  {0.12}   
  \be 
  \raisebox{-5.3em}
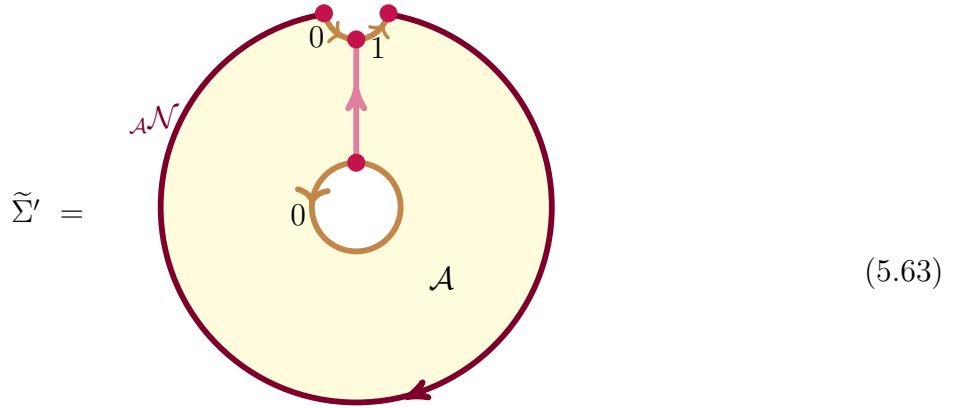
{\begin{tikzpicture}
  \scopeArrow{0.81}{\arrowDefectB}
  \filldraw[\colorFreebdy,fill=\colorSigma,line width=\locpc,postaction={decorate}]
           (0,0) circle (\locpA) node[left=62pt,yshift=33pt]{${}_\cala\caln$};
  \end{scope}
  \scopeArrow{0.49}{>}
  \filldraw[\colorCircle,fill=white,line width=\locpc,postaction={decorate}]
           (0,0) circle (\locpa) node[left=14pt,yshift=-3pt,color=black]{$0$};
  \end{scope}
  \scopeArrow{0.51}{\arrowDefect}
  \draw[line width=\locpc,\colorTransp,postaction={decorate}] (0,\locpa) -- (0,\locpA);
  \end{scope}
  \filldraw[\colorCircle,fill=white,line width=\locpc]
           (-\locpb,\locpA+0.06) arc (180:360:\locpb)
           node[left=19.3pt,yshift=-11.3pt,color=black]{$0$}
           node[right=-12pt,yshift=-15.5pt,color=black]{$1$};
  \scopeArrow{0.36}{>} \scopeArrow{0.83}{>}
  \draw[line width=\locpC,\colorCircle,postaction={decorate}]
           (-\locpb,\locpA+0.06) arc (180:360:\locpb);
  \end{scope} \end{scope}
  \fill[\colorDefect] (0,\locpa) circle (\locpd) (0,\locpA+0.06-\locpb) circle (\locpd)
	   (0,\locpA)++(-3:\locpb) circle (\locpd) (0,\locpA)++(183:\locpb) circle (\locpd);
  \node at (-40:0.57*\locpA) {$\cala$};
  \node at (-\locpA-1.5,0) {$\widetilde{\Sigma}' ~=$};
  \end{tikzpicture}}
  \ee
  ~\\[0.4em]
As pre-block functor of this surface we get
$\Zpreprime(\widetilde{\Sigma}')(z) \,{=} \int^{n\in\caln}\! \Hom(n,z.n)$. 
After invoking the identity \Cite{Eq.\,(3.52)}{fScS2}
$\int^{n\in\caln}\! \cc{n} \boti n \,{\cong} \int_{n\in\caln} \cc{n} \boti \pift(n)$
(with $\pift$ the Nakayama functor \eqref{eq:Nr}) this can be recognized as a space of
natural transformations. The block functor is then given by the corresponding
module natural transformations,
  \be
  \ZZ(\widetilde{\Sigma})(z)=\Nat_{\cala}(\pift,z.{-}) \,,  
  \label{eq:3}
  \ee
where we also the fact that both the Nakayama functor $\pift$ and $z.{-}$ are twisted module
functors. This shows that, in general, when omitting the gluing circle $\SS_{\caln}$ from 
the \defectsurface\ in \eqref{eq:Fig5.8-1} one does \emph{not} obtain the transmission functor.
\end{Remark}


\subsection{Universality of the parallelization}

We are now going to show that the parallelization $\Pi$ whose existence was
established in Theorem \ref{thm:exists-pi} can be uniquely characterized, namely by a
universal property with respect to the parallel transport operations that were introduced in
Section \ref{sec:holonomy}. This observation demonstrates that the choice of parallelization
$\Pi$ obtained by our construction is distinguished. Moreover, it will allow for concrete
computations of the parallelization isomorphisms by working with pre-block functors.

Recall the parallel transport comonad $Z_{\DD,x_{s},x_{t}}$ on
$\Funle(\UU(\partial\DD)^{\opp} \boti \UU(\partial\DD), \vect)$ 
that was considered in Proposition \ref{prop:structure-comodule-hol}, for $\DD$ a disk
in a \defectsurface\ $\Sigma$ and $x_s,x_t$ elements of the set $\Edg$ defined there.
Also, just as the \csfun\ $\Zprime(\Sigma;\Sigma') \colon \ZZ(\Sigma) \,{\to}\, \vect$ 
introduced in \eqref{eq:def:Zprime} is the block functor for a refinement $\Sigma'$ of 
$\Sigma$ evaluated on \mute\ objects in $\partialr(\Sigma')$, we can define a \cspfun\ 
$\Zpreprime(\Sigma;\Sigma') \colon \Zpre(\Sigma) \,{\to}\, \vect$ as the pre-block functor 
on $\Sigma'$ with \mute\ objects inserted in $\partialr(\Sigma')$, according to
  \be
  \Zpreprime(\Sigma;\Sigma') := \Zpre(\Sigma') (- \boti \om_{\Sigma;\Sigma'}) \,.
  \ee

The following result shows that applying the parallel transport  comonad to the pre-blocks for
$\Sigma$ provides the pre-blocks for a refined surface $\Sigma'$:

\begin{Lemma}\label{lem:addnormaldefect}
Let $\DD$ be a disk in a \defectsurface\ $\Sigma$ and $x_s,x_t \iN \Edg$ defect lines of $\DD$. 
\Itemizeiii

\item[{\rm (i)}]
There is $($unique up to isomorphism$)$ a refinement $(\Sigma;\Sigma')$ such that $\Sigma'$
differs from $\Sigma$ by a new transparent defect line that connects $x_{s}$ to $x_{t}$ 
in their normal directions and such that one framing index at $x_{s}$ in $\Sigma'$ is $0$. 

\item[{\rm (ii)}] 
With $\Sigma'$ defined this way, there is a canonical isomorphism 
  \be
  Z_{\DD,x_{s},x_{t}}(\Zpre(\Sigma)) \tocong \Zpreprime(\Sigma')
  \label{eq:can-prebl}
  \ee
of functors.

\item[{\rm (iii)}]
Conversely, let $(\Sigma;\Sigma')$ be a refinement and $\delta$ a transparent defect line
in $\Sigma'$ that is not a defect line of $\Sigma$ such that $(\Sigma;\Sigma'')$, with $\Sigma''$
the \defectsurface\ obtained by deleting $\delta$, is still a refinement of $\Sigma$. Then for
$x_{s},x_{t} \,{\in}\, E_{\Sigma''}$ corresponding to $($choices of$)$ two defect points 
on the gluing boundaries at the start and end of $\delta$, as indicated in
    \def\locpa  {0.7}    
    \def\locpc  {2.2pt}  
    \def\locpd  {0.10}   
    \def\locpD  {0.08}   
    \def\locpe  {0.53}   
    \def\locpx  {2.2}    
    \def\locpy  {2.8}    
  \begin{eqnarray}
  \begin{tikzpicture}
  \begin{scope}[shift={(-3.7,0.8)}]
  \begin{scope}[decoration={random steps,segment length=3mm}] \fill[\colorSigma,decorate] 
           (0,\locpy)++(120:\locpx) arc (120:180:\locpx) -- (-\locpx,0) arc (180:230:\locpx)
           -- (310:\locpx) arc (310:360:\locpx) -- ++(0,\locpy) arc (0:60:\locpx) -- cycle;
  \end{scope}
  \scopeArrow{0.51}{\arrowDefect}
  \draw[line width=\locpc,\colorTransp,postaction={decorate}]
           (0,0) -- node[right=-2pt,yshift=6pt,color=black]{$\delta$} (0,\locpy);
  \end{scope}
  \scopeArrow{0.83}{\arrowDefect}
  \foreach \xx in {0,180,220,340} 
  \draw[line width=\locpc,\colorDefect,postaction={decorate}] (0,0) -- (\xx:\locpx);
  \foreach \xx in {252,271,289,307} \fill (\xx:1.9*\locpa) circle (\locpD);
  \begin{scope}[shift={(0,\locpy)}]
  \foreach \xx in {0,25,50,155,180} 
  \draw[line width=\locpc,\colorDefect,postaction={decorate}] (0,0) -- (\xx:\locpx);
  \foreach \xx in {80,95,110,125} \fill (\xx:1.9*\locpa) circle (\locpD);
  \end{scope}
  \end{scope}
  \filldraw[draw=\colorCircle,fill=white,line width=\locpc]
           (0,0) circle (\locpa) (0,\locpy) circle (\locpa);
  \scopeArrow{0.16}{>} \scopeArrow{0.41}{>} \scopeArrow{0.80}{>}
  \draw[\colorCircle,line width=0.8*\locpc,postaction={decorate}] (0,0) circle (\locpa);
  \end{scope} \end{scope} \end{scope}
  \scopeArrow{0.31}{>} \scopeArrow{0.66}{>} \scopeArrow{0.91}{>}
  \draw[\colorCircle,line width=0.8*\locpc,postaction={decorate}] (0,\locpy) circle (\locpa);
  \end{scope} \end{scope} \end{scope}
  \foreach \xx in {0,90,180,220,340} \filldraw[color=\colorDefect] (\xx:\locpa) circle (\locpd);
  \foreach \xx in {0,25,50,155,180,270}
           \filldraw[color=\colorDefect] (0,\locpy)++(\xx:\locpa) circle (\locpd);
  \node at (-\locpx-1.32,0.5*\locpy) {$\Sigma' ~=$};
  \end{scope}
  \begin{scope}[shift={(\locpx-1.6,-2.0*\locpy)}]
  \begin{scope}[decoration={random steps,segment length=3mm}] \fill[\colorSigma,decorate] 
           (0,\locpy)++(120:\locpx) arc (120:180:\locpx) -- (-\locpx,0) arc (180:230:\locpx)
           -- (310:\locpx) arc (310:360:\locpx) -- ++(0,\locpy) arc (0:60:\locpx) -- cycle;
  \end{scope}
  \scopeArrow{0.83}{\arrowDefect}
  \foreach \xx in {0,180,220,340} 
  \draw[line width=\locpc,\colorDefect,postaction={decorate}] (0,0) -- (\xx:\locpx);
  \foreach \xx in {252,271,289,307} \fill (\xx:1.9*\locpa) circle (\locpD);
  \begin{scope}[shift={(0,\locpy)}]
  \foreach \xx in {0,25,50,155,180} 
  \draw[line width=\locpc,\colorDefect,postaction={decorate}] (0,0) -- (\xx:\locpx);
  \foreach \xx in {80,95,110,125} \fill (\xx:1.9*\locpa) circle (\locpD);
  \end{scope}
  \end{scope}
  \filldraw[draw=\colorCircle,fill=white,line width=\locpc]
           (0,0) circle (\locpa) (0,\locpy) circle (\locpa);
  \scopeArrow{0.30}{>} \scopeArrow{0.83}{>}
  \draw[\colorCircle,line width=0.8*\locpc,postaction={decorate}] (0,0) circle (\locpa);
  \draw[\colorCircle,line width=0.8*\locpc,postaction={decorate}] (0,\locpy) circle (\locpa);
  \end{scope} \end{scope}
  \foreach \xx in {0,180,220,340} \filldraw[color=\colorDefect] (\xx:\locpa) circle (\locpd);
  \foreach \xx in {0,25,50,155,180}
           \filldraw[color=\colorDefect] (0,\locpy)++(\xx:\locpa) circle (\locpd);
  \node at (-\locpx-1.32,0.5*\locpy) {$\Sigma'' ~=$};
  \node[\colorCircle] at (46:1.44*\locpa) {$\SS_1$};
  \node[\colorCircle] at (1.12*\locpa,0.77*\locpy) {$\SS_2$};
  \node (NodeS1) at (168:0.91*\locpa) {};
  \node (NodeS2) at (12:0.91*\locpa) {};
  \node (NodeT1) at (-0.87*\locpa,\locpy-0.15) {};
  \node (NodeT2) at (0.87*\locpa,\locpy-0.15) {};
  \node (NodeSS) [draw=black,thick,rectangle,rounded corners,inner sep=0pt]
	   at (2.2*\locpx,0.32*\locpy) {\begin{tabular}c \rm\footnotesize possible\\[-4pt]
           \rm\footnotesize choices for $x_s$ \end{tabular}};
  \node (NodeTT) [draw=black,thick,rectangle,rounded corners,inner sep=0pt]
           at (2*\locpx,0.76*\locpy) {\begin{tabular}c \rm\footnotesize possible\\[-4pt] 
           \rm\footnotesize choices for $x_t$ \end{tabular}};
  \draw[->] (NodeSS) [out=165,in=160] to (NodeS1);
  \draw[->] (NodeSS) [out=195,in=20] to (NodeS2);
  \draw[->] (NodeTT) [out=195,in=200] to (NodeT1);
  \draw[->] (NodeTT) [out=165,in=-20] to (NodeT2);
  \end{scope}
  \end{tikzpicture} \hspace*{-3em} 
  \nonumber\\[-3.3em]~\\[-0.1em]~\nonumber
  \end{eqnarray}
there is a canonical isomorphism 
  \be
  Z_{\DD,x_{s},x_{t}}(\Zpreprime(\Sigma'')) \tocong \Zpreprime(\Sigma') \,.
  \label{eq:can-prebl-delete}
  \ee
\end{itemize}
\end{Lemma}

\begin{proof}
(i)\, $\DD$ is a disk having $x_{s}$ and $x_{t}$ as parts of its boundary. Connecting
$x_{s}$ and $x_{t}$ by an additional defect line leads to four new gluing boundaries;
setting one of the corresponding new framing indices to $0$ uniquely fixes the other
three. On the disk any other defect line from $x_{s}$ to $x_{t}$ would be isotopic to
the chosen one, hence all refinements that arise this way are isomorphic, 
compare Remark \ref{Remark:structure-def-bord}(iii).
 \\[.2em]
(ii)\, We compute $\Zpre(\Sigma')$ for the first case in the proof of Proposition
\ref{prop:structure-comodule-hol}: The local situation at the new defect line in $\Sigma'$ is
as indicated in the following picture: 
 \\[-1.9em]~ 
    \def\locpa  {0.46}   
    \def\locpc  {2.2pt}  
    \def\locpd  {0.13}   
    \def\locpe  {0.53}   
    \def\locpx  {4.9}    
    \def\locpy  {4.3}    
  \be
  \raisebox{-7.7em}
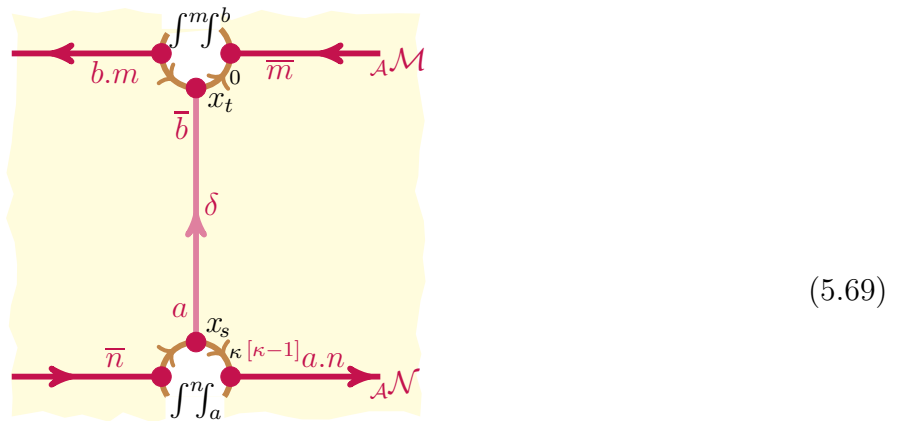
{\begin{tikzpicture}
  \begin{scope}[decoration={random steps,segment length=3mm}] \fill[\colorSigma,decorate] 
	   (0,-\locpe) rectangle (\locpx+\locpe,\locpy+\locpe);
  \end{scope}
  \scopeArrow{0.09}{\arrowDefectB} \scopeArrow{0.90}{\arrowDefectB}
  \draw[line width=\locpc,\colorDefect,postaction={decorate}]
           (\locpx,0) node[above=-1.5pt,xshift=-32pt] {${}^{[\kappa-1]}a.n$} 
           node[right=-9pt,yshift=-3pt] {${}_\cala\caln$}
	   -- (0,0) node[above=-1.5pt,xshift=39pt] {$\cc n$};
  \end{scope} \end{scope}
  \scopeArrow{0.17}{\arrowDefectB} \scopeArrow{0.90}{\arrowDefectB}
  \draw[line width=\locpc,\colorDefect,postaction={decorate}]
           (0,\locpy) node[below=-2.2pt,xshift=39pt] {$b.m$} -- (\locpx,\locpy)
	   node[below=-3pt,xshift=-38pt] {$\cc m$} node[right=-9pt,yshift=-3pt] {${}_\cala\calm$} ;
  \end{scope} \end{scope}
  \scopeArrow{0.51}{\arrowDefect}
  \draw[line width=\locpc,\colorTransp,postaction={decorate}]
	   (0.5*\locpx,0) node[left=-2pt,yshift=25pt,color=\colorDefect] {$a$}
           -- node[right=-2pt,yshift=4pt,color=\colorDefect] {$\delta$}
           +(0,\locpy) node[left=-2pt,yshift=-27pt,color=\colorDefect] {$\cc b$};
  \end{scope}
  \filldraw[draw=\colorCircle,fill=white,line width=\locpc]
           (0.5*\locpx,0)++(-35:\locpa) arc (-35:215:\locpa)
           (0.5*\locpx,\locpy)++(145:\locpa) arc (145:405:\locpa);
  \fill[white] (0.5*\locpx,0)++(-35:\locpa) rectangle +(-1.7*\locpa,-1.8*\locpe)
           (0.5*\locpx,\locpy)++(45:\locpa) rectangle +(-1.7*\locpa,1.8*\locpe);
  \scopeArrow{0.29}{>} \scopeArrow{0.81}{>}
  \draw[\colorCircle,line width=0.7*\locpc,postaction={decorate}]
           (0.5*\locpx-\locpa,0) arc (180:0:\locpa);
  \end{scope} \end{scope}
  \scopeArrow{0.29}{>} \scopeArrow{0.81}{>}
  \draw[\colorCircle,line width=0.7*\locpc,postaction={decorate}]
           (0.5*\locpx-\locpa,\locpy) arc (180:360:\locpa);
  \end{scope} \end{scope}
  \foreach \xx in {0,90,180} \filldraw[color=\colorDefect]
           (0.5*\locpx,0)++(\xx:\locpa) circle (\locpd);
  \foreach \xx in {0,180,270} \filldraw[color=\colorDefect]
           (0.5*\locpx,\locpy)++(\xx:\locpa) circle (\locpd);
  \node at (0.5*\locpx+1.14*\locpa,0.68*\locpa) {$\scriptstyle \kappa$};
  \node at (0.5*\locpx+1.13*\locpa,\locpy-0.68*\locpa) {$\scriptstyle 0$};
  \node at (0.5*\locpx+0.65*\locpa,1.42*\locpa) {$x_{\!s}$};
  \node at (0.5*\locpx,-0.69*\locpa) {$\int^n\!\!\! \int_a$};
  \node at (0.5*\locpx+0.72*\locpa,\locpy-1.42*\locpa) {$x_t$};
  \node at (0.5*\locpx+0.05,\locpy+0.74*\locpa) {$\int^m\!\!\! \int^b$};
  \end{tikzpicture}}
  \ee
Here the labels at the defect lines together with the (co)ends indicate the relevant \mute\ 
objects. In the pre-block functor we can use the Yoneda lemma for $a,b$, so as to obtain
the \cspfun
  \be
  \Zpreprime(\Sigma')(-) \,= \int^{n \in \caln}\!\!\!\!\int^{m \in \calm}\!\!\!\!\int_{a\in\cala}
  \!  \Hom (\cc{n} \boti {}^{[\kappa-1]}a.n  \boti a.m \boti \cc{m} \,\cdots\!, -) \,, 
  \label{eq:ZpreprimeSimgap} 
  \ee
where the ellipsis indicates the remaining defect lines  of $\Sigma'$. With $\mu$
the sum of the framing indices counted clockwise from $x_{s}$ to $x_{t}$ and $l$ the number
of the gluing segments along that path, it follows that $\kappa \,{=}\,l{-}\mu$ and thus
$\Zpreprime(\Sigma')$ is isomorphic to $Z_{\DD,x_{s},x_{t}}(\Zpre(\Sigma))$ according to the
definition of the comonad in Proposition \ref{prop:structure-comodule-hol}.
 \\[.2em]
(iii)\, Consider the \mute\ objects $\om(\SS_{i})$ in the gluing categories for $\SS_{1}$ 
and $\SS_{2}$ for $\Sigma''$. The only case that is not already covered by (ii) above is 
the situation that the relevant framings on $\SS_{i}$ are not both $\pm 1$ (this can happen 
if the adjacent defect lines are transparent, see the case of the defect surfaces considered
in Proposition \ref{Proposition:dehntwist}). However, the transparent objects are explicitly 
known (see Example \ref{exa:JJ}(ii)), so that we can proceed as in the proof of statement (ii). 
\end{proof}

Denote, as above, the new transparent defect line on $\Sigma'$ featuring in Lemma
\ref{lem:addnormaldefect} by $\delta$. Consider the situation shown in the following picture:
    \def\locpa  {0.5}    
    \def\locpc  {2.2pt}  
    \def\locpC  {1.4pt}  
    \def\locpd  {0.10}   
    \def\locpe  {0.53}   
    \def\locpx  {5.2}    
    \def\locpy  {3.1}    
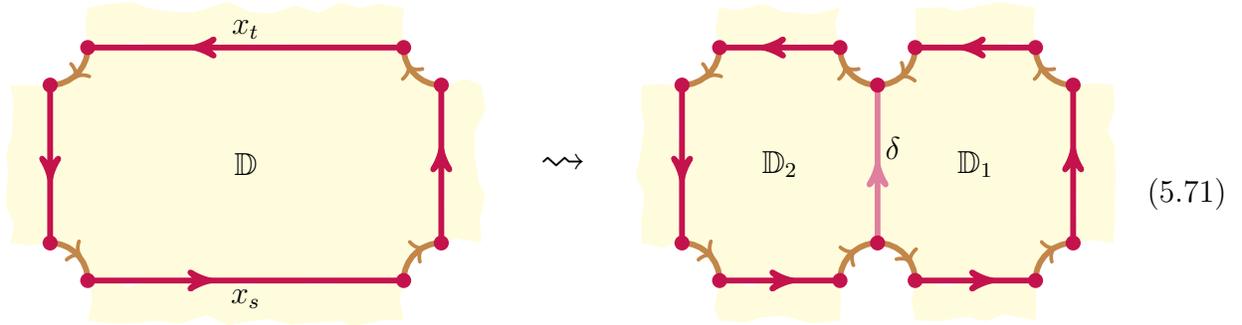
  \begin{eqnarray}
  \begin{tikzpicture}
  \begin{scope}[decoration={random steps,segment length=3mm}]
	  \fill[\colorSigma,decorate] (\locpa,0) rectangle (\locpx-\locpa,-\locpe);
	  \fill[\colorSigma,decorate] (\locpa,\locpy) rectangle (\locpx-\locpa,\locpy+\locpe);
	  \fill[\colorSigma,decorate] (0,\locpa) rectangle (-\locpe,\locpy-\locpa);
	  \fill[\colorSigma,decorate] (\locpx,\locpa) rectangle (\locpx+\locpe,\locpy-\locpa);
  \end{scope}
  \scopeArrow{0.10}{\arrowDefectB} \scopeArrow{0.32}{\arrowDefectB}
  \scopeArrow{0.60}{\arrowDefectB} \scopeArrow{0.89}{\arrowDefectB}
  \filldraw[line width=\locpc,draw=\colorDefect,fill=\colorSigma,postaction={decorate}]
           (0,0) rectangle (\locpx,\locpy);
  \end{scope} \end{scope} \end{scope} \end{scope}
  \fill[white] (0,0) circle (\locpa) (\locpx,0) circle (\locpa)
           (0,\locpy) circle (\locpa) (\locpx,\locpy) circle (\locpa);
  \draw[line width=\locpc,\colorCircle]
           (\locpa,0) arc (0:90:\locpa) (\locpx-\locpa,0) arc (180:90:\locpa)
           (0,\locpy-\locpa) arc (270:360:\locpa) (\locpx,\locpy-\locpa) arc (270:180:\locpa);
  \scopeArrow{0.62}{<}
  \draw[line width=\locpC,\colorCircle,postaction={decorate}] (\locpa,0) arc (0:90:\locpa);
  \draw[line width=\locpC,\colorCircle,postaction={decorate}]
           (0,\locpy-\locpa) arc (270:360:\locpa);
  \end{scope}
  \scopeArrow{0.70}{>}
  \draw[line width=\locpC,\colorCircle,postaction={decorate}]
           (\locpx-\locpa,0) arc (180:90:\locpa);
  \draw[line width=\locpC,\colorCircle,postaction={decorate}]
           (\locpx,\locpy-\locpa) arc (270:180:\locpa);
  \end{scope}
  \fill[\colorDefect] (\locpa,0) circle (\locpd) (0,\locpa) circle (\locpd)
           (\locpx-\locpa,0) circle (\locpd) (\locpx,\locpa) circle (\locpd)
           (\locpa,\locpy) circle (\locpd) (0,\locpy-\locpa) circle (\locpd)
           (\locpx-\locpa,\locpy) circle (\locpd) (\locpx,\locpy-\locpa) circle (\locpd);
  \node at (0.5*\locpx,-7pt) {$x_s$};
  \node at (0.5*\locpx,\locpy cm+7pt) {$x_t$};
  \node at (0.5*\locpx,0.5*\locpy) {$\DD$};
  \node at (\locpx+1.62,0.5*\locpy) {{\Large$\rightsquigarrow$}};
  \begin{scope}[shift={(\locpx+3.2,0)}]
  \begin{scope}[decoration={random steps,segment length=3mm}]
	  \fill[\colorSigma,decorate] (\locpa,0) rectangle (0.5*\locpx-\locpa,-\locpe);
	  \fill[\colorSigma,decorate] (0.5*\locpx+\locpa,0) rectangle (\locpx-\locpa,-\locpe);
	  \fill[\colorSigma,decorate] (\locpa,\locpy) rectangle (0.5*\locpx-\locpa,\locpy+\locpe);
	  \fill[\colorSigma,decorate] (0.5*\locpx+\locpa,\locpy) rectangle (\locpx-\locpa,\locpy+\locpe);
	  \fill[\colorSigma,decorate] (0,\locpa) rectangle (-\locpe,\locpy-\locpa);
	  \fill[\colorSigma,decorate] (\locpx,\locpa) rectangle (\locpx+\locpe,\locpy-\locpa);
  \end{scope}
  \scopeArrow{0.10}{\arrowDefectB} \scopeArrow{0.27}{\arrowDefectB}
  \scopeArrow{0.43}{\arrowDefectB} \scopeArrow{0.60}{\arrowDefectB}
  \scopeArrow{0.77}{\arrowDefectB} \scopeArrow{0.93}{\arrowDefectB}
  \filldraw[line width=\locpc,draw=\colorDefect,fill=\colorSigma,postaction={decorate}]
           (0,0) rectangle (\locpx,\locpy);
  \end{scope} \end{scope} \end{scope} \end{scope} \end{scope} \end{scope}
  \scopeArrow{0.51}{\arrowDefect}
  \draw[line width=\locpc,\colorTransp,postaction={decorate}]
           (0.5*\locpx,0) -- node[right=-2pt,yshift=6pt,color=black]{$\delta$} ++(0,\locpy);
  \end{scope}
  \fill[white] (0,0) circle (\locpa) (\locpx,0) circle (\locpa)
           (0,\locpy) circle (\locpa) (\locpx,\locpy) circle (\locpa)
           (0.5*\locpx,0) circle (\locpa) (0.5*\locpx,\locpy) circle (\locpa);
  \draw[line width=\locpc,\colorCircle]
           (\locpa,0) arc (0:90:\locpa) (\locpx-\locpa,0) arc (180:90:\locpa)
           (0,\locpy-\locpa) arc (270:360:\locpa) (\locpx,\locpy-\locpa) arc (270:180:\locpa)
           (0.5*\locpx+\locpa,0) arc (0:180:\locpa) (0.5*\locpx-\locpa,\locpy) arc (180:360:\locpa);
  \scopeArrow{0.62}{<}
  \draw[line width=\locpC,\colorCircle,postaction={decorate}] (\locpa,0) arc (0:90:\locpa);
  \draw[line width=\locpC,\colorCircle,postaction={decorate}]
           (0,\locpy-\locpa) arc (270:360:\locpa);
  \end{scope}
  \scopeArrow{0.70}{>}
  \draw[line width=\locpC,\colorCircle,postaction={decorate}]
           (\locpx-\locpa,0) arc (180:90:\locpa);
  \draw[line width=\locpC,\colorCircle,postaction={decorate}]
           (\locpx,\locpy-\locpa) arc (270:180:\locpa);
  \end{scope}
  \scopeArrow{0.30}{>}
  \scopeArrow{0.85}{>}
  \draw[line width=\locpC,\colorCircle,postaction={decorate}]
           (0.5*\locpx-\locpa,0) arc (180:0:\locpa);
  \draw[line width=\locpC,\colorCircle,postaction={decorate}]
           (0.5*\locpx+\locpa,\locpy) arc (360:180:\locpa);
  \end{scope} \end{scope}
  \fill[\colorDefect] (\locpa,0) circle (\locpd) (0,\locpa) circle (\locpd)
           (\locpx-\locpa,0) circle (\locpd) (\locpx,\locpa) circle (\locpd)
           (\locpa,\locpy) circle (\locpd) (0,\locpy-\locpa) circle (\locpd)
           (\locpx-\locpa,\locpy) circle (\locpd) (\locpx,\locpy-\locpa) circle (\locpd)
           (0.5*\locpx-\locpa,0) circle (\locpd) (0.5*\locpx+\locpa,0) circle (\locpd)
           (0.5*\locpx-\locpa,\locpy) circle (\locpd) (0.5*\locpx+\locpa,\locpy) circle (\locpd)
           (0.5*\locpx,\locpa) circle (\locpd) (0.5*\locpx,\locpy-\locpa) circle (\locpd);
  \node at (0.25*\locpx,0.5*\locpy) {$\DD_2$};
  \node at (0.75*\locpx,0.5*\locpy) {$\DD_1$};
  \end{scope}
  \end{tikzpicture}
  \nonumber \\[-2.5em]~
  \label{eq:BildSigmaprime-ref}
  \\[1.1em]~ \nonumber
  \end{eqnarray}
The defect line $\delta$ in $\Sigma'$ splits the disk $\DD$ in two disks $\DD_{1}$ 
and $\DD_{2}$, which we label in such a way that the counterclockwise 
path from $x_{s}$ to $x_{t}$ lies in $\DD_{1}$. We start with the following 

\begin{Lemma} \label{Lemma:techn-parallel-transp}
Let $\DD$, $\DD_1$ and $\DD_2$ be disks as in the situation shown in
\eqref{eq:BildSigmaprime-ref}. Then $\holcc_{\DD_{2}} \,{=}\, \holc_{\DD_{1}}$.
\end{Lemma}

\begin{proof}
We show that under the isomorphism 
$\Zpreprime(\Sigma') \,{\cong}\,  Z_{\DD_{1},x_{s},x_{t}}(\Zpre(\Sigma'))$
(with $x_{s}$ and $x_{t}$ the defect lines appearing in \eqref{eq:BildSigmaprime-ref}),
the morphisms $\holcc_{\DD_{2}}$ and $\holc_{\DD_{1}}$ both coincide with the comultiplication
$Z_{\DD_{1},x_{s},x_{t}}(\Zpre(\Sigma')) \,{\to}\, Z_{\DD_{1},x_{s},x_{t}}^{2}(\Zpre(\Sigma'))$. 
To this end we first note that the comultiplication of the coalgebra
$\int_{a \in \cala} \cc{a} \boti a \iN \cc{\cala} \boti \cala$ is given by the end over 
$b \iN \cala $ of the composite morphisms
  \be
  \int_{a\in\cala} \cc{a} \boti a \xrightarrow{~\coevl_{b}~}
  \int_{a\in\cala} \cc{a} \boti (a \otimes {}^{\vee\!} b \otimes b)
  \xrightarrow{~\cong~} \int_{a\in\cala} \cc{b \oti a}\, \boti (a \oti b) \,, 
  \ee
where the second morphism is obtained from the canonical central structure of the end. 
The end over these morphisms equals the end over the morphisms 
  \be
  \int_{a\in\cala} \cc{a} \boti a \xrightarrow{~\coevr_{b}~}
  \int_{a\in\cala} \cc{a} \boti (b \oti b^{\vee} {\otimes}\, a )
  \xrightarrow{~\cong~} \int_{a\in\cala} \cc{b \oti a}\, \boti (a \oti b) \,. 
  \ee
According to Proposition \ref{prop:structure-comodule-hol}, the comonad 
$Z_{\DD_{1},x_{s},x_{t}}(\Zpre(\Sigma'))$ is given by acting with $\int_{a} \cc{a} \boti a$ 
on $\UU(\partial\DD)^{\opp} \boti \UU(\partial\DD)$. It follows that the two morphisms above 
are mapped to the morphisms $\holc_{\DD_{1}}$ and $\holcc_{\DD_{1}}$, respectively.
Thus the statement follows. 
\end{proof}

We now define canonical morphisms between the pre-block functors on $\Sigma$ and $\Sigma'$
that corresponding to the creation and deletion of $\delta$, respectively. To this end we make
uses of the following fact (which is weaker than the statement that morphisms between diagrams
induce morphisms between their limits, but is still elementary):

\begin{Lemma} \label{Lemma:tech-equali}
For $\calc$ an abelian category, let $E_1 \,{\xrightarrow{~\iota_1~}}\,A_1$ be the equalizer
of a pair of morphisms $f_{1},f_{2}\colon A_{1} \,{\to}\, A_{2}$ in $\calc$, and 
$E_2 \,{\xrightarrow{~\iota_2~}}\,A_2$ the equalizer of morphisms
$g_{1},g_{2}\colon B_{1} \,{\to}\, B_{2}$.  Let $a\colon A_{1} \,{\to}\, B_{1}$ and
$b\colon A_{2} \,{\to}\, B_{2}$ be such that for every $i\,{\in}\, \{1,2\}$ there is a
$j \,{\in}\, \{1,2\}$ so that $g_{i} \cir a  \cir \iota_{1}\,{=}\, b \cir f_{j} \cir \iota_{1}$.
Then by restriction $a$ induces a unique morphism from $E_1$ to $E_2$.
\end{Lemma}

\begin{proof}
It follows directly that the morphism $a \cir \iota_{1}$ equalizes the pair 
$(g_{1},g_{2})$, and thus the morphism from $E_{1}$ to $E_{2}$ exists  by the universal property of the equalizer $E_{2}$. In particular, the induced morphism does not depend on the choice of $b$ that fulfills the assumption. 
\end{proof}

Note that, in terms of diagrams, in 
  \be
  \begin{tikzcd}[column sep=1.4cm, row sep=1.1cm]
  E_{1} \ar{r}{\iota_{1}} \ar[dashed]{d}{} & A_{1}
  \ar[yshift=-2]{r}[below]{f_{2}} \ar[yshift=2]{r}{f_{1}} \ar{d}{a} & A_{2} \ar{d}{b}
  \\
  E_{2} \ar{r}{\iota_{2}} & B_{1}\ar[yshift=-2]{r}[below]{g_{2}} \ar[yshift=2]{r}{g_{1}} & B_{2}
  \end{tikzcd}
  \label{eq:komm-diag-equa}
  \ee
it is not required that all four right squares commute. As a special case, we see that if 
in the situation of the Lemma for $i\,{\in}\, \{1,2\}$ there is a $j \,{\in}\, \{1,2\}$ so that 
$g_i \cir a \,{=}\, b \cir f_j$, then there is a unique induced morphism between the equalizers.

If for a morphism $a\colon A_{1} \,{\to}\, B_{1}$ as in this lemma there exists a
morphism $b\colon A_{2} \,{\to}\, B_{2}$ that fulfills the condition stated in the lemma,
we say that $a$ is \emph{compatible with the equalizers}. It turns out that such a 
compatibility is present for the parallel transports $\holc_{\DD}$ and $\holcc_{\DD}$ on 
the pre-blocks in the situation treated in Lemma \ref{lem:addnormaldefect}:

\begin{Proposition} \label{prop:can-Pi}
Let $\Sigma$ be a \defectsurface\ and $(\Sigma;\Sigma')$ the refinement 
described in Lemma $\ref{lem:addnormaldefect}$.
\Itemizeiii

\item[{\rm (i)}]
The clock- and counterclockwise parallel transports 
  \be
  \holc_{\DD}, \holcc_{\DD} :\quad
  \Zpre(\Sigma) \rightarrow  Z_{\DD,x_{s},x_{t}} (\Zpre(\Sigma)) \cong \Zpreprime(\Sigma')
  \label{eq:comp-hol}
  \ee
are compatible with the equalizers to the block spaces for all parallel transport operations
on $\Sigma$ and on $\Sigma'$. Both morphisms 
induce the same morphism between the corresponding fine block functors. 

\item[{\rm (ii)}]
Let now $\Sigma'$ and $\Sigma''$ be as in Lemma {\rm \ref{lem:addnormaldefect}(iii)}. 
The counit $\epsilon$ of the comonad $Z_{\DD,x_{s},x_{t}}$ defines a morphism
  \be
  \epsilon_{\DD}^{} :\quad
  \Zpreprime(\Sigma') \cong  Z_{\DD,x_{s},x_{t}} (\Zpre(\Sigma)) \xrightarrow{~~} \Zpre(\Sigma)
  \label{eq:counit-eps}
  \ee
which is compatible with the equalizers to the block spaces for all parallel transport
operations on $\Sigma$ and on $\Sigma'$.

\end{itemize}
\end{Proposition}

\begin{proof}
For brevity of the exposition, we identify $Z_{\DD,x_{s},x_{t}}(\Zpre(\Sigma))$ with
$\Zpreprime(\Sigma')$ via the canonical isomorphism \eqref{eq:can-prebl}.
 \\[.2em]
(i)\, Consider again the situation shown in Figure \eqref{eq:BildSigmaprime-ref}. 
We treat the case of $\holcc_{\DD} \colon \Zpre(\Sigma) \,{\to}$
	                                                       $\Zpreprime(\Sigma')$;
the proof for $\holc_{\DD}$ is analogous. We consider the parallel transports on the 
disks $\DD_{1}$ and $\DD_{2}$. For $\DD_{1}$ we deal with
the morphisms $\holc_{\DD_{1}}$ and $\holcc_{\DD_{1}}$ and find corresponding parallel
transports for $\Sigma$ as in Lemma \ref{Lemma:tech-equali}. This gives rise to the diagrams
  \be
  \hspace*{-1.2em}
  \begin{tikzcd}[column sep=1.1cm,row sep=1.1cm]
  \Zpre(\Sigma)\ar{r}{\holcc_{\DD}} \ar{d}[swap]{\holcc_{\DD}}
  & Z_{\DD,x_{s},x_{t}}(\Zpre(\Sigma)) \ar[xshift=-19pt]{d}{Z_{\DD,x_{s},x_{t}}(\holcc_{\DD})}
  \\ 
  \Zpreprime(\Sigma')\ar{r}{\holcc_{\DD_{1}}} & Z_{\DD,x_{s},x_{t}}(\Zpreprime(\Sigma')) 
  \end{tikzcd}  
  \quad \text{~and~} \quad
  \begin{tikzcd}[column sep=1.1cm,row sep=1.1cm]
  \Zpre(\Sigma)\ar{r}{\holcc_{\DD}} \ar{d}[swap]{\holcc_{\DD}}
  & Z_{\DD,x_{s},x_{t}}(\Zpre(\Sigma)) \ar[xshift=-19pt]{d}{Z_{\DD,x_{s},x_{t}}(\holcc_{\DD})}
  \\
  \Zpreprime(\Sigma')\ar{r}{\holc_{\DD_{1}}} & Z_{\DD,x_{s},x_{t}}(\Zpreprime(\Sigma')) 
  \end{tikzcd}  
  \ee
The first of these diagrams commutes directly. The second diagram commutes as well,
owing to the fact that, according to Proposition \ref{prop:structure-comodule-hol},
$\holcc_{\DD}$ provides a comodule structure. 
 \\
The corresponding diagrams for the disk $\DD_{2}$ are
  \be
  \hspace*{-1.2em}
  \begin{tikzcd}[column sep=1.1cm,row sep=1.1cm]
  \Zpre(\Sigma)\ar{r}{\holcc_{\DD}} \ar{d}[swap]{\holcc_{\DD}}
  & Z_{\DD,x_{s},x_{t}}(\Zpre(\Sigma)) \ar[xshift=-19pt]{d}{Z_{\DD,x_{s},x_{t}}(\holcc_{\DD})}
  \\
  \Zpreprime(\Sigma')\ar{r}{\holcc_{\DD_{2}}} & Z_{\DD,x_{s},x_{t}}(\Zpreprime(\Sigma')) 
  \end{tikzcd}  
  \quad \text{~and~} \quad
  \begin{tikzcd}[column sep=1.1cm,row sep=1.1cm]
  \Zpre(\Sigma)\ar{r}{\holc_{\DD}} \ar{d}[swap]{\holcc_{\DD}}
  & Z_{\DD,x_{s},x_{t}}(\Zpre(\Sigma)) \ar[xshift=-19pt]{d}{Z_{\DD,x_{s},x_{t}}(\holcc_{\DD})}
  \\
  \Zpreprime(\Sigma')\ar{r}{\holc_{\DD_{2}}} & Z_{\DD,x_{s},x_{t}}(\Zpreprime(\Sigma')) 
  \end{tikzcd}  
  \ee
Here the first diagram commutes because $\holcc_{\DD}$ is a comodule structure, while the
second diagram commutes because the parallel transports take place in different disks. 
It follows analogously that $\holc_{\DD}$ is compatible with the equalizers, and this
induces another morphism from $\ZZ(\Sigma)$ to $\Zprime(\Sigma')$. Since, by Proposition 
\ref{prop:structure-comodule-hol}(iii), $\ZZ(\Sigma)$ is the equalizer of both parallel 
transports, it follows that the two induced morphisms between the block functors agree. 
 \\[.2em]
(ii)\, We use the counit $\epsilon$ of $Z_{\DD, x_{s},x_{t}}$, which by definition of the 
comonad is given by the component at $\one$ of the dinatural transformation of the end. 
Then, with $\Sigma'$ as in Figure \eqref{eq:BildSigmaprime-ref}, we have  
$\Zpreprime(\Sigma') \,{\cong}\, Z_{\DD,x_{s},x_{t}}(\Zpre(\Sigma))$, and it follows that
$\epsilon$ defines a natural transformation 
$\epsilon_{\DD} \colon \Zpreprime(\Sigma') \,{\to}\, \Zpre(\Sigma)$.
For the parallel transports $\holc_{\DD}$ and $\holcc_{\DD}$ we need to provide
corresponding parallel transports for one disk of $\Sigma'$ such that the condition
of Lemma \ref{Lemma:tech-equali} is satisfied. Consider the  diagrams 
  \be
  \hspace*{-1.2em}
  \begin{tikzcd}[column sep=1.1cm,row sep=1.1cm]
  \Zpreprime(\Sigma')\ar{r}{\holc_{\DD_{2}}} \ar{d}[swap]{\epsilon_{\DD}}
  & Z_{\DD,x_{s},x_{t}}(\Zpreprime(\Sigma')) \ar[xshift=-19pt]{d}{Z_{\DD,x_{s},x_{t}}(\epsilon_{\DD})}
  \\
  \Zpre(\Sigma)\ar{r}{\holc_{\DD}} & Z_{\DD,x_{s},x_{t}}(\Zpre(\Sigma)) 
  \end{tikzcd}  
  \quad \text{~and~} \quad
  \begin{tikzcd}[column sep=1.1cm,row sep=1.1cm]
  \Zpreprime(\Sigma')\ar{r}{\holcc_{\DD_{2}}} \ar{d}[swap]{\epsilon_{\DD}}
  & Z_{\DD,x_{s},x_{t}}(\Zpreprime(\Sigma')) \ar[xshift=-19pt]{d}{Z_{\DD,x_{s},x_{t}}(\epsilon_{\DD})}
  \\
  \Zpre(\Sigma)\ar{r}{\holcc_{\DD}} & Z_{\DD,x_{s},x_{t}}(\Zpre(\Sigma)). 
  \end{tikzcd}  
  \ee
The left diagram commutes directly. In contrast, the right diagram does not, in general,
commute. Let us compose the right diagram with the equalizer of the block functor for $\Sigma'$
to obtain
  \be
  \begin{tikzcd}[column sep=1.1cm,row sep=0.9cm]
  \Zprime(\Sigma') \ar{dr}{\iota}  \ar[bend left]{drr} \ar[bend right]{rdd} && 
  \\[-0.4cm]
  & \Zpreprime(\Sigma')\ar{r}{\holcc_{\DD_{2}}} \ar{d}[swap]{\epsilon_{\DD}}
  & Z_{\DD,x_{s},x_{t}}(\Zpreprime(\Sigma')) \ar[xshift=-19pt]{d}{Z_{\DD,x_{s},x_{t}}(\epsilon_{\DD})}
  \\
  & \Zpre(\Sigma)\ar{r}{\holcc_{\DD}} & Z_{\DD,x_{s},x_{t}}(\Zpre(\Sigma)). 
  \end{tikzcd}     
  \label{eq:outer-diag}
  \ee
where the unlabeled arrows are defined as the corresponding composites. 
Consider the outer square in this diagram. From Lemma \ref{Lemma:techn-parallel-transp}
we conclude that $\holcc_{\DD_{2}} \,{=}\, \holc_{\DD_{1}}$, while the fact that $\iota$
is the equalizer of the parallel transports implies that 
$\holc_{\DD_{1}} \cir \iota \,{=}\, \holcc_{\DD_{1}} \cir \iota$. Since the diagram 
  \be
  \begin{tikzcd}[column sep=1.1cm,row sep=1.1cm]
  \Zpreprime(\Sigma')\ar{r}{\holcc_{\DD_{1}}} \ar{d}[swap]{\epsilon_{\DD}}
  & Z_{\DD,x_{s},x_{t}}(\Zpreprime(\Sigma')) \ar[xshift=-19pt]{d}{Z_{\DD,x_{s},x_{t}}(\epsilon_{\DD})}
  \\
  \Zpre(\Sigma)\ar{r}{\holcc_{\DD}} & Z_{\DD,x_{s},x_{t}}(\Zpre(\Sigma))
  \end{tikzcd}  
  \ee
commutes directly, it follows directly that the outer square in \eqref{eq:outer-diag} commutes.
We can now invoke Lemma \ref{Lemma:tech-equali}(ii) we obtain a morphism between the block functors. 
\end{proof}

We denote the induced morphisms between the block functors from part (i) and (ii) of
Proposition \ref{prop:can-Pi} by $\widehathol_{\Sigma} \colon \ZZ(\Sigma) \,{\to}\, \Zprime(\Sigma')$ and 
$\widehat{\epsilon}_{\Sigma} \colon \Zprime(\Sigma') \,{\to}\, \ZZ(\Sigma)$.  
To show that these morphisms agree with the parallelizations we require a compatibility
with factorization: Let $\Sigma \,{=}\, \Sigma_{1} \cir \Sigma_{2}$ be a fine \defectsurface\
that is the composite of two fine \defectsurface s along a common boundary component
$\SS$. Assume that $x_{s},x_{t} \,{\in}\, \Edg$ are two defect lines that intersect $\SS$,
i.e.\ they lie in  both $\Sigma_{1}$ and $\Sigma_{2}$. We then obtain morphisms 
  \be
  \widehathol_{\!\Sigma_1} :\quad \ZZ(\Sigma_{1}) \rightarrow \Zprime(\Sigma_{1}')
  \qquad \text{and} \qquad
  \widehathol_{\!\Sigma_2} :\quad \ZZ(\Sigma_{2}) \rightarrow \Zprime(\Sigma_{2}') \,, 
  \ee
with $\Sigma_{i}'$ the respective refinements. 

\begin{Lemma} \label{Lemma:hol-factor}
In the situation just described we have isomorphisms 
  \be
  \Sigma'_{1} \cir \Sigma_2^{} \tocong \Sigma_1^{} \cir \Sigma_{2}' \tocong \Sigma'
  \ee
of \defectsurface s, where $(\Sigma;\Sigma')$ is the refinement of $\Sigma$ that corresponds 
to to the pair $x_{s},x_{t}$. Furthermore, the morphisms on the block functors satisfy
  \be
  \widehathol_{1} \circ \ZZ(\Sigma_{2}) = \ZZ(\Sigma_{1}) \circ \widehathol_{2}
  = \widehathol_{\Sigma} \,.   
  \ee
\end{Lemma}

\begin{proof}
Consider the following situation:
    \def\locpa  {0.44}   
    \def\locpA  {1.22}   
    \def\locpc  {2.2pt}  
    \def\locpC  {1.4pt}  
    \def\locpd  {0.10}   
    \def\locpe  {0.53}   
    \def\locpx  {6.5}    
    \def\locpy  {5.3}    
  \be
  \raisebox{-7.5em}
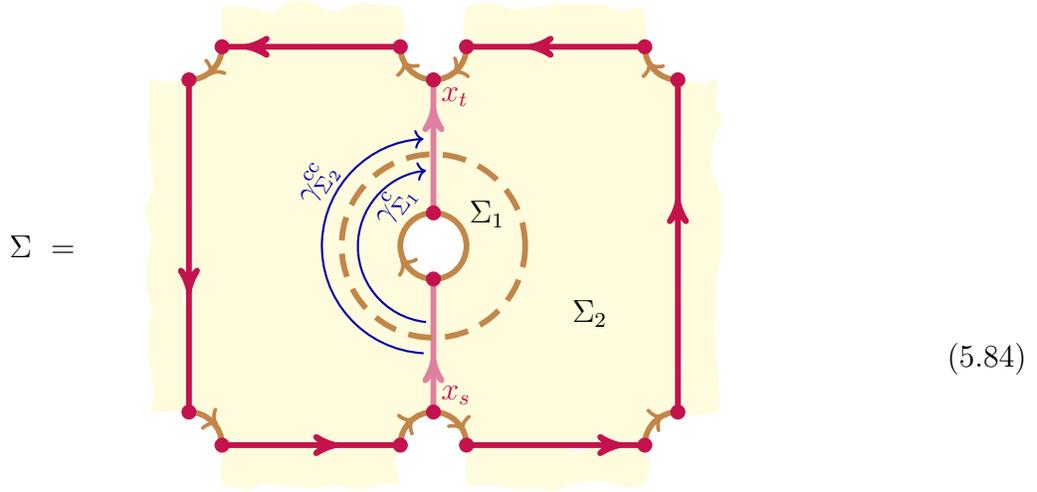
{\begin{tikzpicture}
  \begin{scope}[decoration={random steps,segment length=3mm}]
	  \fill[\colorSigma,decorate] (\locpa,0) rectangle (0.5*\locpx-\locpa,-\locpe);
	  \fill[\colorSigma,decorate] (0.5*\locpx+\locpa,0) rectangle (\locpx-\locpa,-\locpe);
	  \fill[\colorSigma,decorate] (\locpa,\locpy) rectangle (0.5*\locpx-\locpa,\locpy+\locpe);
	  \fill[\colorSigma,decorate] (0.5*\locpx+\locpa,\locpy) rectangle (\locpx-\locpa,\locpy+\locpe);
	  \fill[\colorSigma,decorate] (0,\locpa) rectangle (-\locpe,\locpy-\locpa);
	  \fill[\colorSigma,decorate] (\locpx,\locpa) rectangle (\locpx+\locpe,\locpy-\locpa);
  \end{scope}
  \scopeArrow{0.10}{\arrowDefectB} \scopeArrow{0.27}{\arrowDefectB}
  \scopeArrow{0.43}{\arrowDefectB} \scopeArrow{0.60}{\arrowDefectB}
  \scopeArrow{0.77}{\arrowDefectB} \scopeArrow{0.93}{\arrowDefectB}
  \filldraw[line width=\locpc,draw=\colorDefect,fill=\colorSigma,postaction={decorate}]
           (0,0) rectangle (\locpx,\locpy);
  \end{scope} \end{scope} \end{scope} \end{scope} \end{scope} \end{scope}
  \scopeArrow{0.22}{\arrowDefect} \scopeArrow{0.85}{\arrowDefect}
  \draw[line width=\locpc,\colorTransp,postaction={decorate}] (0.5*\locpx,0) --
           node[right=-2pt,very near start,color=\colorDefect] {$x_s$}
           node[right=-2pt,very near end,color=\colorDefect] {$x_t$} ++(0,\locpy);
  \end{scope} \end{scope}
  \draw[\colorCircle,line width=\locpc,dash pattern=on 9pt off 5pt]
           (0.5*\locpx,0.5*\locpy) circle (\locpA);
  \fill[white] (0,0) circle (\locpa) (\locpx,0) circle (\locpa)
           (0,\locpy) circle (\locpa) (\locpx,\locpy) circle (\locpa)
           (0.5*\locpx,0) circle (\locpa) (0.5*\locpx,\locpy) circle (\locpa);
  \filldraw[line width=\locpc,\colorCircle,fill=white] (0.5*\locpx,0.5*\locpy) circle (\locpa);
  \draw[line width=\locpc,\colorCircle]
           (\locpa,0) arc (0:90:\locpa) (\locpx-\locpa,0) arc (180:90:\locpa)
           (0,\locpy-\locpa) arc (270:360:\locpa) (\locpx,\locpy-\locpa) arc (270:180:\locpa)
           (0.5*\locpx+\locpa,0) arc (0:180:\locpa) (0.5*\locpx-\locpa,\locpy) arc (180:360:\locpa);
  \scopeArrow{0.62}{<}
  \draw[line width=\locpC,\colorCircle,postaction={decorate}]
           (0.5*\locpx,0.5*\locpy) circle (\locpa);
  \draw[line width=\locpC,\colorCircle,postaction={decorate}] (\locpa,0) arc (0:90:\locpa);
  \draw[line width=\locpC,\colorCircle,postaction={decorate}]
           (0,\locpy-\locpa) arc (270:360:\locpa);
  \end{scope}
  \scopeArrow{0.70}{>}
  \draw[line width=\locpC,\colorCircle,postaction={decorate}]
           (\locpx-\locpa,0) arc (180:90:\locpa);
  \draw[line width=\locpC,\colorCircle,postaction={decorate}]
           (\locpx,\locpy-\locpa) arc (270:180:\locpa);
  \end{scope}
  \scopeArrow{0.30}{>} \scopeArrow{0.86}{>}
  \draw[line width=\locpC,\colorCircle,postaction={decorate}]
           (0.5*\locpx-\locpa,0) arc (180:0:\locpa);
  \draw[line width=\locpC,\colorCircle,postaction={decorate}]
           (0.5*\locpx+\locpa,\locpy) arc (360:180:\locpa);
  \end{scope} \end{scope}
  \fill[\colorDefect] (\locpa,0) circle (\locpd) (0,\locpa) circle (\locpd)
           (\locpx-\locpa,0) circle (\locpd) (\locpx,\locpa) circle (\locpd)
           (\locpa,\locpy) circle (\locpd) (0,\locpy-\locpa) circle (\locpd)
           (\locpx-\locpa,\locpy) circle (\locpd) (\locpx,\locpy-\locpa) circle (\locpd)
           (0.5*\locpx-\locpa,0) circle (\locpd) (0.5*\locpx+\locpa,0) circle (\locpd)
           (0.5*\locpx-\locpa,\locpy) circle (\locpd) (0.5*\locpx+\locpa,\locpy) circle (\locpd)
           (0.5*\locpx,\locpa) circle (\locpd) (0.5*\locpx,\locpy-\locpa) circle (\locpd)
           (0.5*\locpx,0.5*\locpy-\locpa) circle (\locpd)
           (0.5*\locpx,0.5*\locpy+\locpa) circle (\locpd);
  \draw[\colorHol,thick,->] (0.485*\locpx,0.5*\locpy-0.83*\locpA) arc (264:96:0.83*\locpA)
           node[left=4pt,yshift=-2.5pt,rotate=53] {$\holc_{\Sigma_1}$};
  \draw[\colorHol,thick,->] (0.48*\locpx,0.5*\locpy-1.17*\locpA) arc (267:93:1.17*\locpA)
           node[left=33pt,yshift=-6pt,rotate=57] {$\holcc_{\Sigma_2}$};
  \node at (-1.95,0.5*\locpy) {$\Sigma ~=$};
  \node at (0.61*\locpx,0.58*\locpy) {$\Sigma_1$};
  \node at (0.82*\locpx,0.33*\locpy) {$\Sigma_2$};
  \end{tikzpicture}}
  \ee
Since by Proposition \ref{prop:can-Pi} the morphisms $\widehathol_{i}$ do not depend on
the choice of parallel transport in $\Sigma_{i}$, we may just consider the parallel transports
$\holcc_{\Sigma_{2}}$ and $\holc_{\Sigma_{1}}$. By factorization for pre-blocks we have 
  \be
  \Zpre(\Sigma_{1}) \circ \Zpre(\Sigma_{2}) \,\cong\, \Zpre(\Sigma_1) \boti \Zpre(\Sigma_2)\big(\mbox{\large$\int^{z \in \ZZ(\SS)}$} \cc{z} \boti z\big)  ,
  \ee
and the parallel transports $\holcc_{\Sigma_{2}}$ and $\holc_{\Sigma_{1}}$ correspond to post-composition
with the canonical morphism $\int^{z \in \ZZ(\SS)} \cc{z} \boti z \,{\to} 
\int^{z\in\ZZ(\SS)}\!\!\int_{a\in\cala} \cc{a.z} \boti a.z$.
It follows that $\widehathol_1 \cir \ZZ(\Sigma_2) \,{=}\, \ZZ(\Sigma_1) \cir \widehathol_2$.
Using again that $\widehathol_{\Sigma}$ is independent of the choice of clock- or
counterclockwise parallel transport, it follows that both of these are also equal to
$\widehathol_{\Sigma}$.
\end{proof}

Next we show that
the induced morphisms between the block functors are identical to the parallelization 
morphisms from Theorem \ref{thm:exists-pi}, which implies in particular that the
morphisms $\widehathol_{\Sigma}$ and $\widehat{\epsilon}_{\Sigma}$ are indeed isomorphisms. 

\begin{Proposition}
\label{Proposition:Hol-Pi}
The morphisms $\widehathol_{\Sigma} \colon \Zfine(\Sigma) \,{\to}\, \Zfine(\Sigma')$ 
and $\widehat{\epsilon}_{\Sigma} \colon \Zfine(\Sigma') \,{\to}\, \Zfine(\Sigma)$ that are
induced by the morphisms 
\eqref{eq:comp-hol} and \eqref{eq:counit-eps} are the same as the parallelization isomorphisms.
\end{Proposition}

\begin{proof}
Consider first the situation that the \defectsurface s $\Sigma_1^{}$ and $\Sigma_{1}'$ are 
as in the following picture:
    \def\locpa  {1.6}    
    \def\locpb  {0.36}   
    \def\locpc  {2.2pt}  
    \def\locpd  {0.12}   
  \be
  \raisebox{-3.9em}
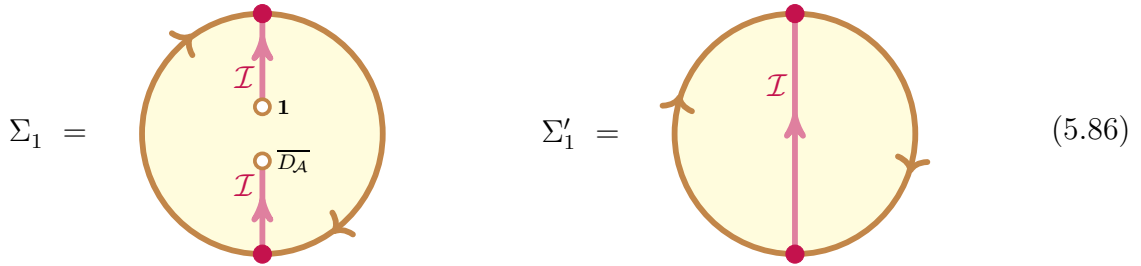
{\begin{tikzpicture}
  \scopeArrow{0.37}{<} \scopeArrow{0.87}{<}
  \filldraw[\colorCircle,fill=\colorSigma,line width=\locpc,postaction={decorate}] (0,0) circle (\locpa);
  \end{scope} \end{scope}
  \scopeArrow{0.51}{\arrowDefectB}
  \draw[line width=\locpc,color=\colorTransp,postaction={decorate}]
             (90:\locpa) node[left=-2pt,yshift=-24pt,color=\colorDefect] {$\II$} -- (90:\locpb);
  \end{scope}
  \scopeArrow{0.66}{\arrowDefectB}
  \draw[line width=\locpc,color=\colorTransp,postaction={decorate}]
             (270:\locpb) node[left=-2pt,yshift=-9pt,color=\colorDefect] {$\II$} -- (270:\locpa);
  \end{scope}
  \transv {90:\locpb} node[right=1pt,color=black] {$\scriptstyle \one$};
  \transv {270:\locpb} node[right=1pt,color=black] {$\scriptstyle \cc{D_\Cala}$};
  \filldraw[color=\colorDefect] (90:\locpa) circle (\locpd) (270:\locpa) circle (\locpd);
  \node at (-\locpa-1.25,0) {$\Sigma_1^{} ~=$};
  \begin{scope}[shift={(2*\locpa+3.88,0)}]
  \scopeArrow{0.47}{<} \scopeArrow{0.97}{<}
  \filldraw[\colorCircle,fill=\colorSigma,line width=\locpc,postaction={decorate}] (0,0) circle (\locpa);
  \end{scope} \end{scope}
  \scopeArrow{0.52}{\arrowDefectB}
  \draw[line width=\locpc,color=\colorTransp,postaction={decorate}]
             (90:\locpa) node[left=-2pt,yshift=-27pt,color=\colorDefect] {$\II$} -- (270:\locpa);
  \end{scope}
  \filldraw[color=\colorDefect] (90:\locpa) circle (\locpd) (270:\locpa) circle (\locpd);
  \node at (-\locpa-1.25,0) {$\Sigma_1' ~=$};
  \end{scope}
  \end{tikzpicture}}
  \ee
As start $x_{s}$ and end $x_{t}$ of the parallel transport we select the two vertical defect
lines. Then by Lemma \ref{lem:addnormaldefect} we identify 
$Z_{\Sigma_{1},x_{s},x_{t}}(\Zpre(\Sigma_{1}))$ with  $\Zpreprime(\Sigma_{1}')$. Then the claim
follows by observing that the diagram 
  \be
  \begin{tikzcd}
  \ZZ(\Sigma_{1}) \ar{r} \ar{d}[swap]{\Pi_{\Sigma_{1};\Sigma_{1},\Sigma_{1}'}}
  & \Zpre(\Sigma_{1}) \ar{d}{\holc_{\Sigma_{1}}(x_{s},x_{t})}
  \\
  \Zprime(\Sigma_{1}') \ar{r} & \Zpreprime(\Sigma_{1}')
  \end{tikzcd}
  \label{eq:from-LemmC3}
  \ee
commutes by Lemma \ref{Lemma:univ-blockiso}. 
\\
For the case of a general \defectsurface\ $\Sigma$, with start $x_{s}$ and end $x_{t}$
of the parallel transport, take the \defectsurface\ $\Sigma'$ as on the right hand side
of \eqref{eq:BildSigmaprime-ref}, and the \defectsurface\ $\widetilde{\Sigma}$, with
factorization into $\Sigma_{1}$ and $\Sigma_{2}$, as indicated in 
    \def\locpa  {0.44}   
    \def\locpA  {1.17}   
    \def\locpc  {2.2pt}  
    \def\locpC  {1.4pt}  
    \def\locpd  {0.10}   
    \def\locpe  {0.53}   
    \def\locpx  {6.5}    
    \def\locpy  {5.3}    
  \be
  \raisebox{-7.5em}{\begin{tikzpicture}
  \begin{scope}[decoration={random steps,segment length=3mm}]
	  \fill[\colorSigma,decorate] (\locpa,0) rectangle (0.5*\locpx-\locpa,-\locpe);
	  \fill[\colorSigma,decorate] (0.5*\locpx+\locpa,0) rectangle (\locpx-\locpa,-\locpe);
	  \fill[\colorSigma,decorate] (\locpa,\locpy) rectangle (0.5*\locpx-\locpa,\locpy+\locpe);
	  \fill[\colorSigma,decorate] (0.5*\locpx+\locpa,\locpy) rectangle (\locpx-\locpa,\locpy+\locpe);
	  \fill[\colorSigma,decorate] (0,\locpa) rectangle (-\locpe,\locpy-\locpa);
	  \fill[\colorSigma,decorate] (\locpx,\locpa) rectangle (\locpx+\locpe,\locpy-\locpa);
  \end{scope}
  \scopeArrow{0.10}{\arrowDefectB} \scopeArrow{0.27}{\arrowDefectB}
  \scopeArrow{0.43}{\arrowDefectB} \scopeArrow{0.60}{\arrowDefectB}
  \scopeArrow{0.77}{\arrowDefectB} \scopeArrow{0.93}{\arrowDefectB}
  \filldraw[line width=\locpc,draw=\colorDefect,fill=\colorSigma,postaction={decorate}]
           (0,0) rectangle (\locpx,\locpy);
  \end{scope} \end{scope} \end{scope} \end{scope} \end{scope} \end{scope}
  \scopeArrow{0.56}{\arrowDefect}
  \draw[line width=\locpc,\colorTransp,postaction={decorate}] (0.5*\locpx,0) -- ++(0,0.4*\locpy);
  \end{scope}
  \scopeArrow{0.54}{\arrowDefectB}
  \draw[line width=\locpc,\colorTransp,postaction={decorate}]
           (0.5*\locpx,\locpy) -- ++(0,-0.4*\locpy);
  \end{scope}
  \transv {0.5*\locpx,0.4*\locpy};
  \transv {0.5*\locpx,0.6*\locpy};
  \draw[\colorCircle,line width=\locpc,dash pattern=on 9pt off 5pt]
           (0.5*\locpx,0.5*\locpy) circle (\locpA);
  \fill[white] (0,0) circle (\locpa) (\locpx,0) circle (\locpa)
           (0,\locpy) circle (\locpa) (\locpx,\locpy) circle (\locpa)
           (0.5*\locpx,0) circle (\locpa) (0.5*\locpx,\locpy) circle (\locpa);
  \draw[line width=\locpc,\colorCircle]
           (\locpa,0) arc (0:90:\locpa) (\locpx-\locpa,0) arc (180:90:\locpa)
           (0,\locpy-\locpa) arc (270:360:\locpa) (\locpx,\locpy-\locpa) arc (270:180:\locpa)
           (0.5*\locpx+\locpa,0) arc (0:180:\locpa) (0.5*\locpx-\locpa,\locpy) arc (180:360:\locpa);
  \scopeArrow{0.62}{<}
  \draw[line width=\locpC,\colorCircle,postaction={decorate}] (\locpa,0) arc (0:90:\locpa);
  \draw[line width=\locpC,\colorCircle,postaction={decorate}]
           (0,\locpy-\locpa) arc (270:360:\locpa);
  \end{scope}
  \scopeArrow{0.70}{>}
  \draw[line width=\locpC,\colorCircle,postaction={decorate}]
           (\locpx-\locpa,0) arc (180:90:\locpa);
  \draw[line width=\locpC,\colorCircle,postaction={decorate}]
           (\locpx,\locpy-\locpa) arc (270:180:\locpa);
  \end{scope}
  \scopeArrow{0.30}{>}
  \scopeArrow{0.85}{>}
  \draw[line width=\locpC,\colorCircle,postaction={decorate}]
           (0.5*\locpx-\locpa,0) arc (180:0:\locpa);
  \draw[line width=\locpC,\colorCircle,postaction={decorate}]
           (0.5*\locpx+\locpa,\locpy) arc (360:180:\locpa);
  \end{scope} \end{scope}
  \fill[\colorDefect] (\locpa,0) circle (\locpd) (0,\locpa) circle (\locpd)
           (\locpx-\locpa,0) circle (\locpd) (\locpx,\locpa) circle (\locpd)
           (\locpa,\locpy) circle (\locpd) (0,\locpy-\locpa) circle (\locpd)
           (\locpx-\locpa,\locpy) circle (\locpd) (\locpx,\locpy-\locpa) circle (\locpd)
           (0.5*\locpx-\locpa,0) circle (\locpd) (0.5*\locpx+\locpa,0) circle (\locpd)
           (0.5*\locpx-\locpa,\locpy) circle (\locpd) (0.5*\locpx+\locpa,\locpy) circle (\locpd)
           (0.5*\locpx,\locpa) circle (\locpd) (0.5*\locpx,\locpy-\locpa) circle (\locpd);
  \node at (-1.95,0.5*\locpy) {$\widetilde\Sigma ~=$};
  \node at (0.43*\locpx,0.5*\locpy) {$\Sigma_1$};
  \node at (0.82*\locpx,0.33*\locpy) {$\Sigma_2$};
  \end{tikzpicture}}
  \ee
Then consider the diagram 
  \be
  \label{eq:comm-diag-facto}
  \begin{tikzcd}
  \ZZ(\Sigma) \ar{rr} \ar[bend right=75]{dd}[left]{\Pi_{_{\Sigma;\Sigma,\Sigma'}}} \ar{d}
  & ~ & \Zpre(\Sigma) \ar{d} \ar[bend left=75]{dd}[right]{\holcc_{\Sigma}}
  \\
  \ZZ(\widetilde{\Sigma}) \ar{d}{\Pi_{\Sigma_{1};\Sigma_{1},\Sigma_{1}'}}^{} \ar{rr}
  & ~ & \Zpre(\widetilde{\Sigma}) \ar{d}[left]{\holc_{\Sigma_{1}}}
  \\
  \Zprime(\Sigma') \ar{rr} & ~ & \Zpreprime(\Sigma')
  \end{tikzcd}
  \ee
where all horizontal arrows are the canonical morphisms from the block to the pre-block 
functor, the vertical arrows in the upper square are obtained from Lemma \ref{Lemma:excision-pi}.
This diagram commutes: Commutativity of the upper square follows
from excision, while commutativity of the lower square is a consequence of factorization 
and commutativity of the diagram \ref{eq:from-LemmC3} above. The subdiagram to the left
commutes by definition of the parallelization, and the subdiagram to the right by Lemma 
\ref{Lemma:hol-factor}. It follows that $\widehathol_{\Sigma} \,{=}\, \Pi_{\Sigma;\Sigma,\Sigma'}$
is  in particular invertible.
 \\
The morphism $\widehat{\epsilon}_{\Sigma}$ is by construction right inverse to
$\holcc_{\Sigma}(x_{s},x_{t})$: For the relative pre-preblock functor $\Zpreprime(\Sigma')$,
the equality $\holcc_{\Sigma}(x_{s},x_{t}) \circ \epsilon=\id$ holds by construction. 
Thus it follows that $\widehathol_{\Sigma} \cir \widehat{\epsilon} \eq \id$, 
and since $\widehathol_{\Sigma}$ is invertible, $\widehat{\epsilon}$ its two-sided inverse,
which thus agrees with the inverse of the parallelization morphism. 
\end{proof}

\begin{Remark} 
Proposition \ref{Proposition:Hol-Pi} provides a more conceptual understanding of the
parallelization $\Pi$. Our specific construction of $\Pi$ has, in contrast, the virtue that it
is more local and thereby allows one to establish the coherence properties of the parallelization. 
\end{Remark}


\subsection{Actions of mapping class groups} \label{sec:mcg}

The structures we have defined provide us directly with a representation of the mapping class
group of a \defectsurface\ $\Sigma$ by isomorphisms of the block functor $\ZZ(\Sigma)$. To
make this precise we define the value of the block functor on the 2-morphisms of the bicategory
$\Bfrdefdec$ (recall the description of the latter in Definition \ref{def:Bordcats}).

For the present purposes it is convenient to record all the structure of a \defectsurface\
$\Sigma$ in the notation. Thus we write $\Sigma \,{=}\, (\vSigma, \rho, \delta, \chi)$, where 
$\vSigma$ is the underlying surface, whose boundary is parametrized according to the
parametrization $\rho$, $\delta$ is the set of defect lines on $\Sigma$, and $\chi$ is the 
framing. Recall from Section \ref{sec:Bfrdef} that a 2-morphism 
  \be
  \varphi :\quad (\vSigma_{1}, \rho_{1}, \delta_{1}, \chi_{1}) \xrightarrow{~~}
  (\vSigma_{2}, \rho_{2}, \delta_{2}, \chi_{2})
  \label{eq:2-morphi-Brfr}
  \ee
in $\Bfrdefdec$ is an isotopy class of isomorphisms from $\vSigma_{1}$ to $\vSigma_{2}$
preserving all the structure, i.e.\ it is a diffeomorphism 
relative to the boundary parametrizations, it respects the defect lines,
and the push-forward vector field $\varphi_{*} \chi_{1}$ is equal to $\chi_{2}$. In
particular, given an isomorphism $\varphi \colon \Sigma_{1} \,{\to}\, \Sigma_{2}$,
with $(\vSigma_{1}, \rho_{1}, \delta_{1}, \chi_{1})$ a \defectsurface, there is the
induced structure $(\vSigma_2, \varphi(\rho_1), \varphi_*(\delta_1), \varphi_*(\chi_1))$
of a \defectsurface\ on $\Sigma_{2}$ such that $\varphi$ represents a morphism of
\defectsurface s. In case that  $\varphi \colon \vSigma \,{\to}\, \vSigma$ is an automorphism
of the underlying surface $\vSigma$ of $ \Sigma \,{=}\, (\vSigma, \rho, \delta, \chi)$,
it can happen that $(\vSigma, \varphi(\rho), \varphi_*(\delta),
\varphi_{*}(\chi))$ is not the same object as $\Sigma$ in $\Bfrdefdec$: the boundary
parametrization might have changed and/or the induced vector field may not be
homotopic to the original one. We will exhibit examples of both phenomena later. 

Based on the parallelization and on the definition of the block functor $\ZZ$ as a limit 
(see Definition \ref{def:THEdef})
we can directly specify the value of $\ZZ$ on the 2-morphisms in $\Bfrdefdec$:
  
\begin{Lemma} \label{lemma:repMCG}
Let $ \varphi \colon (\Sigma_{1}, \rho_{1}, \delta_{1}, \chi_{1}) \,{\to}\, 
(\Sigma_{2}, \rho_{2}, \delta_{2}, \chi_{2})$ be a morphism of \defectsurface s. For every
fine refinement $(\Sigma_1^{};\Sigma_1')$, the image $\varphi_*(\delta_1')$ of the defects in
$\Sigma_1'$ provides a fine refinement of $\Sigma_2$. As a consequence,
$\varphi$ induces a unique isomorphism 
  \be
  \ZZ(\varphi) :\quad \lim \Pi_{\Sigma_1^{};\Sigma_1'} \xrightarrow{~\cong~}
  \lim \Pi_{\Sigma_2^{}; \Sigma_2'} 
  \label{eq:T(varphi)}
  \ee
of block functors.
Moreover, if $\varphi_1$ and $\varphi_2$ are composable morphisms of \defectsurface s,
then $\ZZ(\varphi_2 \,{\veco}\, \varphi_1 ) = \ZZ(\varphi_2) \,{\veco}\, \ZZ(\varphi_1)$.
  \end{Lemma}

\begin{proof}
The first statement is geometrically obvious. As in Definition \ref{def:Gamma-Pi} we
denote by $\Gamma(\Sigma_{i})$ the category of fine refinements of the surface $\Sigma_i$,
so that we have 
  \be
  \ZZ(\Sigma_i) = \lim_{\Sigma_i'}\Pi_{\Sigma_i}(\Sigma_i') \,,
  \ee
where
$\Pi_{\Sigma_i} \colon \Gamma(\Sigma_i) \,{\to}\, \Funle(\ZZ(\partial\Sigma_i), \vect)$ with
$\Pi_{\Sigma_i}(\Sigma') \,{=}\, \Pi_{\Sigma_i;\Sigma'}$ the parallelization on $\Sigma_{i}$.
The push-forward along $\varphi$ provides an equivalence 
$\varphi_{*}\colon \Gamma(\Sigma_1) \,{\xrightarrow{\,\simeq\,}}\, \Gamma(\Sigma_2)$. Since
the boundaries agree, we have $\ZZ(\partial\Sigma_1)\,{=}\,\ZZ(\partial\Sigma_2)$. Next we will
show that there is a strict equality $\Pi_{\Sigma_1} \,{=}\, \Pi_{\Sigma_2}\cir\varphi_{*}
\colon \Gamma(\Sigma_{1})
\,{\to}\, \Funle(\ZZ(\partial\Sigma_1), \vect)$ of functors. It then 
follows by standard arguments for limits that there is a unique
induced isomorphism $\ZZ(\varphi)$ as in \ref{eq:T(varphi)}.
 \\[.2em]
To verify our claim, consider a refinement $(\Sigma_1^{},\rho_1^{},\chi_1^{};\Sigma_1')$
of $\Sigma_{1}$, i.e.\ an object in $\Gamma(\Sigma_{1})$. Since the block functor is
determined entirely by the combinatorial data of the framing indices and the incidence
relations of the defect lines, and since these data agree on the two \defectsurface s
$(\Sigma_2^{},\rho_2^{},\chi_2^{};\varphi(\Sigma_1'))$ and
$(\Sigma_1^{},\rho_{1}^{},\chi_{1}^{};\Sigma_{1}')$, it follows that the functors
$\Pi_{\Sigma_1}(\Sigma_1^{};\Sigma_1') $ and $\Pi_{\Sigma_2}(\Sigma_2^{};\varphi(\Sigma_{1}'))$
are equal. Thus the functors $\Pi_{\Sigma_1}$ and $\Pi_{\Sigma_2} \cir \varphi_*$ agree on 
objects. Moreover, every local \filldiskrep\ on $\Sigma_{1}'$ corresponds to a local 
\filldiskrep\ under $\varphi_{*}$. Therefore the functors agree on morphisms as well. 
Finally, the uniqueness of the induced morphisms between the limits directly implies
that $\ZZ$ is compatible with the composition of 2-morphisms.
\end{proof}

\begin{Definition} \label{def:MF2-morph}
The value of the modular functor $\ZZ$ on a 2-morphism $\varphi$ in $\Bfrdefdec$
is the isomorphism \eqref{eq:T(varphi)} described in Lemma \ref{lemma:repMCG}.
\end{Definition}

To determine the isomorphism $\ZZ(\varphi)$ in practice, one picks representatives for the
block functors on $\Sigma_1$ and $\Sigma_2$ by choosing fine refinements 
$(\Sigma_1^{};\Sigma_1')$ and $(\Sigma_2^{};\Sigma_2')$ and uses the cone isomorphisms  
$\ZZ(\Sigma_1) \,{\cong}\, \Zprime(\Sigma_1, \rho_1, \delta_1, \chi_1 ; \Sigma_1')$ and 
$\ZZ(\Sigma_2) \,{\cong}\, \Zprime(\Sigma_2, \rho_2, \delta_2, \chi_2 ; \Sigma_2')$ to 
identify the block functors with the latter representatives. Then $\ZZ(\varphi)$ can be
computed as follows:

\begin{Lemma} \label{lem:T2=T1}
The functor $\ZZ(\Sigma_{2}, \rho_{2}, \chi_{2}; \varphi(\Sigma_{1}'))$ is equal to the functor 
$\ZZ(\Sigma_{1},\rho_{1}, \chi_{1}; \Sigma_{1}')$. The isomorphism $\ZZ(\varphi)$ is given by
  \be
  \bearl
  \ZZ(\varphi) = \Pi_{\Sigma_{2}; \varphi(\Sigma_{1}'), \Sigma_{2}'} : 
  \Nxl6 \hspace*{5em}
  \Zprime(\Sigma_{1}, \rho_{1}, \delta_{1}, \chi_{1} ; \Sigma_{1}')
  = \Zprime(\Sigma_{2}, \rho_{2}, \chi_{2}; \varphi(\Sigma_{1}')) 
  \xrightarrow{~~} \Zprime(\Sigma_{2}, \rho_{2}, \chi_{2}; \Sigma_{2}') 
  \label{eq:rep-varphi}
  \eear
  \ee
on the chosen representatives for the block functors.
\end{Lemma}

\begin{proof}
The first statement follows from the proof of Lemma \ref{lemma:repMCG}.
The second statement follows  from the definition of $\ZZ(\varphi)$ as the universal arrow
between the limits, given that we chose representatives for the limits. 
\end{proof}

With Lemma \ref{lem:T2=T1} at hand we can finally collect our results to complete 
the proof of our main Theorem about a modular functor.

\begin{Proofmain}\label{prf:Proofmain}
We collect the data and axioms of a modular functor
$\ZZ\colon \Bfrdefdec \,{\xrightarrow{~}}\, \Funle $: The gluing category assigned to
a \defectonemanifold\ is given in Definition \ref{def:T(L,e,k)}, the block functor for 
a \defectsurface\ is in Definition \ref{def:THEdef}, and the value on a 2-morphism,
i.e.\ an isomorphism of \defectsurface s, is given in Definition \ref{def:MF2-morph}.
The coherence data of $\ZZ$ with respect to the symmetric monoidal product are trivial:
By construction, the disjoint union of objects, 1- or 2-morphisms in $\Bfrdefdec$ is
mapped to the Deligne product of the respective structures. The coherence isomorphism
with respect to the composition of 1-morphisms is provided by factorization as proven
in Theorem \ref{thm:Facto}. For the compositions with units we have again
trivial coherence data for the symmetric monoidal product (to the disjoint union
with the empty set there is assigned an identity morphism), and as shown in
Corollary \ref{cor:gen-cylinder}, the assignment for a cylinder is an identity
1-morphisms in $\Funle$.
 \\[.2em]
Concerning the axioms of a symmetric monoidal functor, the only non-trivial 
statement to be checked is the pentagon axiom for four
composable 1-morphisms. That this is satisfied is shown in Theorem \ref{thm:Facto}.
\end{Proofmain}


\subsection{Explicit actions of mapping class group elements}

For the computation of mapping class group representations we need to derive 
explicit descriptions of the parallelization isomorphisms. 
As a preparation we introduce the following graphical conventions. First, a \tadpolecircle\
$\QQ_\pm$ inside a \defectsurface\ is drawn as an unlabeled small circle, according to
    \def\locpa  {1.3}    
    \def\locpc  {2.2pt}  
  \be
  \raisebox{-0.6em}{\begin{tikzpicture}
  \scopeArrow{0.61}{\arrowDefect}
  \draw[line width=\locpc,color=\colorTransp,postaction={decorate}]
             (0,0) -- node[near end,above=-2pt,color=\colorDefect] {$\II$} (\locpa,0);
  \end{scope}
  \transv {0,0};
  \node at (-1.28,0) {$\QQ_+ ~=$};
  \begin{scope}[shift={(\locpa+1.8,0)}] \node at (0,0) {and}; \end{scope}
  \begin{scope}[shift={(\locpa+5.4,0)}]
  \scopeArrow{0.76}{\arrowDefect}
  \draw[line width=\locpc,color=\colorTransp,postaction={decorate}]
             (0,0) -- node[above=-2pt,near start,color=\colorDefect] {$\II$} (\locpa,0);
  \end{scope}
  \transv {\locpa,0};
  \node at (-1.28,0) {$\QQ_- ~=$};
  \end{scope}
  \end{tikzpicture}}
  \ee
respectively, whenever this facilitates the description. Further, in case we deal with a
block functor for some surface $\Sigma$ that is to be evaluated at the 
\mute\ object $\om(\QQ_+) \,{=}\, \one \,{\in}\, \cent(\cala) \,{\cong}$
                                                                       $\ZZ(\QQ_+)$ or
$\om(\QQ_-) \,{=}\, \cc{\DA} \,{\in}\, \cc{\cent^{-4}(\cala)} \,{\cong}\, \ZZ(\QQ_-)$ of a
\tadpolecircle, we introduce a further convention that allows us to omit the argument $\om(\QQ_\pm)$
of the block functor: we attach instead the \mute\ object as a label to that \tadpolecircle:
    \def\locpa  {1.3}    
    \def\locpb  {0.65}   
    \def\locpc  {2.2pt}  
    \def\locpd  {0.33}   
  \be
  \raisebox{-3.6em}{\begin{tikzpicture}
  \begin{scope}[decoration={random steps,segment length=4mm}]
	  \fill[\colorSigma,decorate] (-\locpd,-\locpb) rectangle (\locpa+\locpd,\locpb);
  \end{scope}
  \scopeArrow{0.61}{\arrowDefect}
  \draw[line width=\locpc,color=\colorTransp,postaction={decorate}]
             (0,0) -- node[near end,above=-2pt,color=\colorDefect] {$\II$} (\locpa,0);
  \end{scope}
  \transv {0,0} node[below=1pt,color=black] {$\scriptstyle \one$};
  \node at (-\locpd-0.5,0) {$\ZZ\,\Big($};
  \node at (\locpa+\locpd+0.77,0) {$\Big)\quad := $};
  \begin{scope}[shift={(\locpa+3.56,0)}]
  \begin{scope}[decoration={random steps,segment length=4mm}]
	  \fill[\colorSigma,decorate] (-\locpd,-\locpb) rectangle (\locpa+\locpd,\locpb);
  \end{scope}
  \scopeArrow{0.61}{\arrowDefect}
  \draw[line width=\locpc,color=\colorTransp,postaction={decorate}]
             (0,0) -- node[near end,above=-2pt,color=\colorDefect] {$\II$} (\locpa,0);
  \end{scope}
  \transv {0,0};
  \node at (-\locpd-0.5,0) {$\ZZ\,\Big($};
  \node at (\locpa+\locpd+1.48,0) {$\Big)(-\boti \om(\QQ_+))$};
  \end{scope}
  \begin{scope}[shift={(0,-2*\locpb-0.5)}]
  \begin{scope}[decoration={random steps,segment length=4mm}]
          \fill[\colorSigma,decorate] (-\locpd,-\locpb) rectangle (\locpa+\locpd,\locpb);
  \end{scope}
  \scopeArrow{0.76}{\arrowDefect}
  \draw[line width=\locpc,color=\colorTransp,postaction={decorate}]
             (0,0) -- node[above=-2pt,near start,color=\colorDefect] {$\II$} (\locpa,0);
  \end{scope}
  \transv {\locpa,0} node[below=1pt,color=black] {$\scriptstyle \cc{D_\Cala}$};
  \node at (-\locpd-2.5,0) {and};
  \node at (-\locpd-0.5,0) {$\ZZ\,\Big($};
  \node at (\locpa+\locpd+0.77,0) {$\Big)\quad := $};
  \end{scope}
  \begin{scope}[shift={(\locpa+3.56,-2*\locpb-0.5)}]
  \begin{scope}[decoration={random steps,segment length=4mm}]
          \fill[\colorSigma,decorate] (-\locpd,-\locpb) rectangle (\locpa+\locpd,\locpb);
  \end{scope}
  \scopeArrow{0.76}{\arrowDefect}
  \draw[line width=\locpc,color=\colorTransp,postaction={decorate}]
             (0,0) -- node[above=-2pt,near start,color=\colorDefect] {$\II$} (\locpa,0);
  \end{scope}
  \transv {\locpa,0};
  \node at (-\locpd-0.5,0) {$\ZZ\,\Big($};
  \node at (\locpa+\locpd+1.48,0) {$\Big)(-\boti \om(\QQ_-))$};
  \end{scope}
  \end{tikzpicture}
  \label{eq:insert1orDA}}
  \ee
(Moreover, later on 
we will sometimes want to refrain from specifying whether a tadpole circle is of the form $\QQ_+$
or $\QQ_-$. We then just use $\om$ as a generic symbol for the \mute\ object of such a circle.
Likewise we treat the appearance of \mute\ objects for other gluing circles, such as for the 
trivalent gluing circles in the picture \eqref{eq:canisoXi} below.)

 \medskip
    
In order to discuss the braid group 
representation, we now consider the following two situations:
    \def\locpa  {2.5}    
    \def\locpA  {0.46}   
    \def\locpB  {0.66}   
    \def\locpc  {2.2pt}  
    \def\locpd  {0.12}   
    \def\locpe  {0.26}   
    \def\locpE  {0.51}   
    \def\locpf  {0.33}   
    \def\locpy  {3.7}    
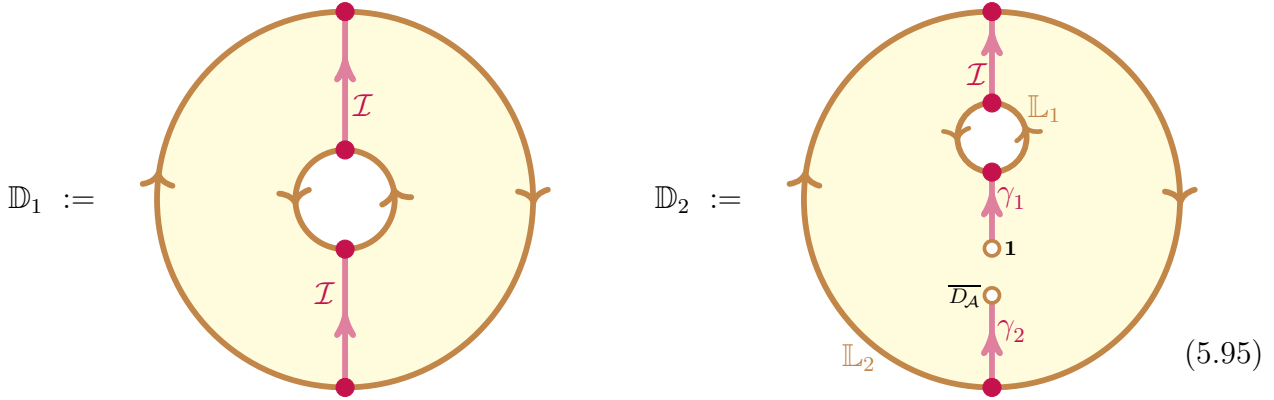
  \begin{eqnarray}
  \begin{tikzpicture}
  \scopeArrow{0.01}{<} \scopeArrow{0.49}{<}
  \filldraw[\colorCircle,fill=\colorSigma,line width=\locpc,postaction={decorate}] (0,0) circle (\locpa);
  \end{scope} \end{scope}
  \scopeArrow{0.20}{\arrowDefect} \scopeArrow{0.88}{\arrowDefect}
  \draw[line width=\locpc,color=\colorTransp,postaction={decorate}]
             (270:\locpa) -- node[left=-1.1pt,near start,color=\colorDefect] {$\II$}
             node[right=-2pt,near end,color=\colorDefect] {$\II$} (90:\locpa);
  \end{scope} \end{scope}
  \scopeArrow{0.04}{>} \scopeArrow{0.52}{>}
  \filldraw[\colorCircle,fill=white,line width=\locpc,postaction={decorate}] (0,0) circle (\locpB);
  \end{scope} \end{scope}
  \filldraw[color=\colorDefect] (90:\locpa) circle (\locpd) (270:\locpa) circle (\locpd)
             (90:\locpB) circle (\locpd) (270:\locpB) circle (\locpd);
  \node at (-\locpa-1.3,0) {$\DD_1~:=$};
  \begin{scope}[shift={(2*\locpa+3.5,0)}]
  \scopeArrow{0.01}{<} \scopeArrow{0.49}{<}
  \filldraw[\colorCircle,fill=\colorSigma,line width=\locpc,postaction={decorate}] (0,0) circle (\locpa);
  \end{scope} \end{scope}
  \scopeArrow{0.85}{\arrowDefect}
  \draw[line width=\locpc,color=\colorTransp,postaction={decorate}]
             (90:\locpf*\locpa) -- node[left=-3.2pt,color=\colorDefect] {$\II$} (90:\locpa);
  \end{scope}
  \scopeArrow{0.72}{\arrowDefectB}
  \draw[line width=\locpc,color=\colorTransp,postaction={decorate}] (90:\locpf*\locpa) -- 
             node[right=-3pt,yshift=-2pt,color=\colorDefect] {$\gamma_1^{}$} (270:\locpe*\locpa);
  \end{scope}
  \scopeArrow{0.66}{\arrowDefectB}
  \draw[line width=\locpc,color=\colorTransp,postaction={decorate}] (270:\locpE*\locpa) --
               node[right=-3pt,yshift=4pt,color=\colorDefect] {$\gamma_2^{}$} (270:\locpa);
  \end{scope}
  \filldraw[\colorCircle,fill=white,line width=\locpc] (90:\locpf*\locpa) circle (\locpA);
  \scopeArrow{0.06}{>} \scopeArrow{0.52}{>}
  \draw[\colorCircle,line width=0.8*\locpc,postaction={decorate}]
             (90:\locpf*\locpa) circle (\locpA);
  \end{scope} \end{scope}
  \transv {270:\locpe*\locpa} node[right,color=black] {$\scriptstyle \one$};
  \transv {270:\locpE*\locpa} node[left=-0.4pt,color=black] {$\scriptstyle \cc{D_\Cala}$};
  \filldraw[color=\colorDefect] (90:\locpa) circle (\locpd) (270:\locpa) circle (\locpd)
             (90:\locpf*\locpa+\locpA) circle (\locpd) (90:\locpf*\locpa-\locpA) circle (\locpd);
  \node at (-\locpa-1.3,0) {$\DD_2~:=$};
  \node[color=\colorCircle] at (60:0.56*\locpa) {$\LL_1$};
  \node[color=\colorCircle] at (230:1.1*\locpa) {$\LL_2$};
  \end{scope}
  \end{tikzpicture} \hspace*{-2.0em}
  \nonumber \\[-1.7em] ~ 
  \label{eq:Fig6.6-1}
  \\[-1.4em] ~ \nonumber
  \end{eqnarray}
As indicated in the picture for $\DD_2$, we denote by $\LL_1$ and $\LL_2$, respectively, the
inner and outer gluing circle (which are the same for the two disks $\DD_1$ and $\DD_2$), and 
by $\gamma_i$, for $i\,{=}\,1,2$, the defect line that connects $\LL_{i}$ to a \tadpoledisk. 
Invoking Corollary \ref{cor:gen-cylinder} we see that the block functor for 
$\DD_1$ is given by $\ZZ(\DD_1)(\cc{G} \boti F) \,{=}\, \Nat_{\cala,\cala}(G,F)$
for $F \iN \ZZ(\LL_{1}) \,{\simeq}\, \Funle_{\cala,\cala}(\cala,\cala)$ and 
$\cc{G} \iN \ZZ(\LL_{2}) \,{\simeq}\, \cc{\ZZ(\LL_1)}$. 
Recall that for $F \iN  \Funle_{\cala,\cala}(\cala,\cala)$, the object $F(\one) \in \cala$
carries a canonical structure of an object in $\calz(\cala)$, using the structure of $F$ as
a bimodule functor. This structure is used in

\begin{Lemma} \label{Lemma:block-iso-disk}
\Itemizeiii

\item[{\rm (i)}]
The block functor $\ZZ(\DD_{2})$ is canonically  isomorphic to the functor 
  \be
  \ZZ(\LL_{1}) \boti \ZZ(\LL_{2}) \ni~ F \boti \cc{G}
  \,\longmapsto\, \Hom_{\calz(\cala)}(G(\one),F(\one)) ~\in \vect \,.
  \label{eq:block-D2}
  \ee

\item[{\rm (ii)}]
The parallelization isomorphism $\Pi_{\DD_{1}; \DD_{1},\DD_{2}}$ is given by 
  \be
  \begin{array}{rcl}
  \Nat_{\cala,\cala}(G,F) &\!\!\! \xrightarrow{~~} \!\!\!& \Hom_{\calz(\cala)}(G(\one),F(\one))
  \Nxl3
  \mu &\!\!\! \xmapsto{~~} \!\!\!& \mu_{\one} \,,
  \end{array}
  \label{eq:parallel-D2}
  \ee
and its inverse is provided by applying the braided induction \eqref{eq:braidedind} to morphisms.
\end{itemize}
\end{Lemma}

\begin{proof}
(i)\, We denote by 
$\ZE_{\gamma_{1}^{}} \colon \Funle_{\Cala,\cala}(\cala,\cala) \,{\to}\, \calz(\cala)$ and by
$\ZE_{\gamma_{2}^{}} \colon \cc{\Funle_{\Cala,\cala}(\cala,\cala)} \,{\to}\, \cc{\calz(\cala)}$
the excision functors (as introduced in Definition \ref{defi:excision}) that are associated
to the defect lines $\gamma_1$ and $\gamma_2$ in the disk $\DD_2$, respectively. 
Lemma \ref{Lemma:Z-same} provides us with a canonical isomorphism $\ZZ(\DD_2)(\cc{G}\boti F)
\,{\tocong}\, \Hom_{\calz(\cala)}(\ZE_{\gamma_2^{}}(G),\ZE_{\gamma_1^{}}(F))$.
Moreover, by the Eilen\-berg--Watts equivalences we have $\ZE_{\gamma_1^{}}(F) \,{=}
\int^{a\in\cala}\!\Hom_{\cala}(a, \one) \boti F(a) \,{\cong}\, F(\one)$, while for 
$\ZE_{\gamma_2^{}}$ we can use that $\cc{\int^{a\in\cala}\! G(a) \boti \cc{a}}\,
{=}\, \int_{a\in\cala} \cc{G(a)} \boti a$ in $\ZZ(\LL_{2})$ and obtain 
  \be
  \begin{array}{ll} \dsty
  \ZE_{\gamma_2^{}}(\int_{a\in\cala}\! \cc{G(a)} \boti a) \!\!\!&\dsty
  = \int_{a\in\cala}\! \Hom_{\cala}(D, a) \otimes \cc{G(a)}
  \Nxl3 &\dsty
  \cong \int^{a\in\cala}\! \Hom_{\cala}(D, D\,{\otimes} {}^{\vee\vee\!}a) \otimes \cc{G(a)}
  \Nxl3 &\dsty
  \cong \int^{a\in\cala}\! \Hom_{\cala}(\one,\! {}^{\vee\vee\!}a) \otimes \cc{G(a)}
  \Nxl3 &\dsty
  \cong \int^{a\in\cala}\! \Hom_{\cala}(\one, a) \otimes \cc{G(a)} \,\cong\, \cc{G(\one)} \,.
  \eear
  \ee
Here we use in the first step that $G$ is exact and then make use of the isomorphism 
\eqref{eq:end-coend-D}.
\\[.2em]
(ii)\, The prescription \eqref{eq:parallel-D2} indeed defines an isomorphism. The local \filldiskrep\
from $\DD_{1}$ to $\DD_{2}$ is covered by the situation analyzed in Lemma \ref{Lemma:univ-blockiso}:
We find sub-disks of $\DD_{1}$ and $\DD_{2}$ which are of the same type as the ones on 
the left and on the right of \ref{eq:canisoXi}. Thus in view of Lemma \ref{Lemma:univ-blockiso},
what we need to show is that \eqref{eq:parallel-D2} comes from the 
corresponding morphism between the pre-block functors. 
The latter is given by the dinatural morphism from an end and in our case is given by
  \be
  \Nat(G,F) \,\cong \int_{\!a\in\cala}\Hom_{\cala}(G(a), F(a))
  \xrightarrow{~~} \Hom_{\cala}(G(\one),F(\one)) \,.
  \label{eq:dinat-end}
  \ee
This proves the claim. 
\end{proof}

We now consider in detail two specific types of isomorphisms $\varphi$: the exchange of 
two boundary circles on a three-holed sphere and the Dehn twist on a cylinder over a circle. 
These are of particular interest, because Dehn twists generate the mapping class groups,
while manipulations analogous to the one for the three-holed sphere generate 
a braid group, which is a prominent example of (a subgroup of) a mapping class group.

 \medskip 

To deduce the action of the first of these two types of isomorphisms
we start with the following \filldisk\ $\DD$:
    \def\locpa  {2.2}    
    \def\locpA  {2.4*\transvR}    
    \def\locpb  {0.56}   
    \def\locpc  {2.2pt}  
    \def\locpC  {1.8pt}  
    \def\locpd  {0.12}   
    \def\locpD  {0.09}   
  \be
  \raisebox{-5.4em}
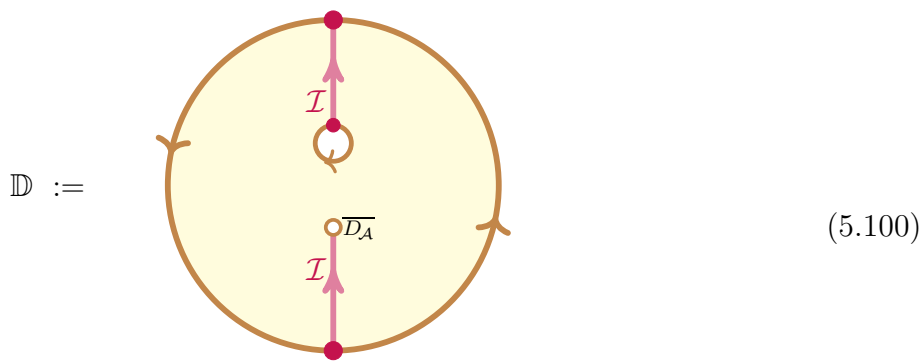
{\begin{tikzpicture}
  \scopeArrow{0.47}{>} \scopeArrow{0.97}{>}
  \filldraw[\colorCircle,fill=\colorSigma,line width=\locpc,postaction={decorate}] (0,0) circle (\locpa);
  \end{scope} \end{scope}
  \scopeArrow{0.67}{\arrowDefect}
  \draw[line width=\locpc,color=\colorTransp,postaction={decorate}]
             (270:\locpa) node[left=-2.3pt,yshift=31pt,color=\colorDefect] {$\II$} -- (270:\locpb);
  \end{scope}
  \scopeArrow{0.70}{\arrowDefect}
  \draw[line width=\locpc,color=\colorTransp,postaction={decorate}]
             (90:\locpb) node[left=-2.3pt,yshift=16pt,color=\colorDefect] {$\II$} -- (90:\locpa);
  \end{scope}
  \transv {270:\locpb} node[right=1pt,color=black] {$\!\scriptstyle \cc{D_\Cala}$};
  \filldraw[\colorCircle,fill=white,line width=\locpC] (90:\locpb) circle (\locpA);
  \scopeArrow{0.77}{<}
  \draw[\colorCircle,line width=0.7*\locpC,postaction={decorate}] (90:\locpb) circle (\locpA);
  \end{scope}
  \filldraw[color=\colorDefect] (90:\locpa) circle (\locpd) (270:\locpa) circle (\locpd)
                  (90:\locpb+\locpA) circle (\locpD);
  \node at (-\locpa-1.6,0) {$\DD ~:=$};
  \end{tikzpicture}
  \label{eq:canisoXi}}
  \ee
Now consider the two refinements $(\DD;\DD_1)$ and $(\DD;\DD_2)$ with
    \def\locpa  {2.8}    
    \def\locpA  {2.4*\transvR}    
    \def\locpb  {0.56}   
    \def\locpB  {0.5*\locpa+0.5*\locpb}
    \def\locpc  {2.2pt}  
    \def\locpC  {1.8pt}  
    \def\locpd  {0.12}   
    \def\locpD  {0.09}   
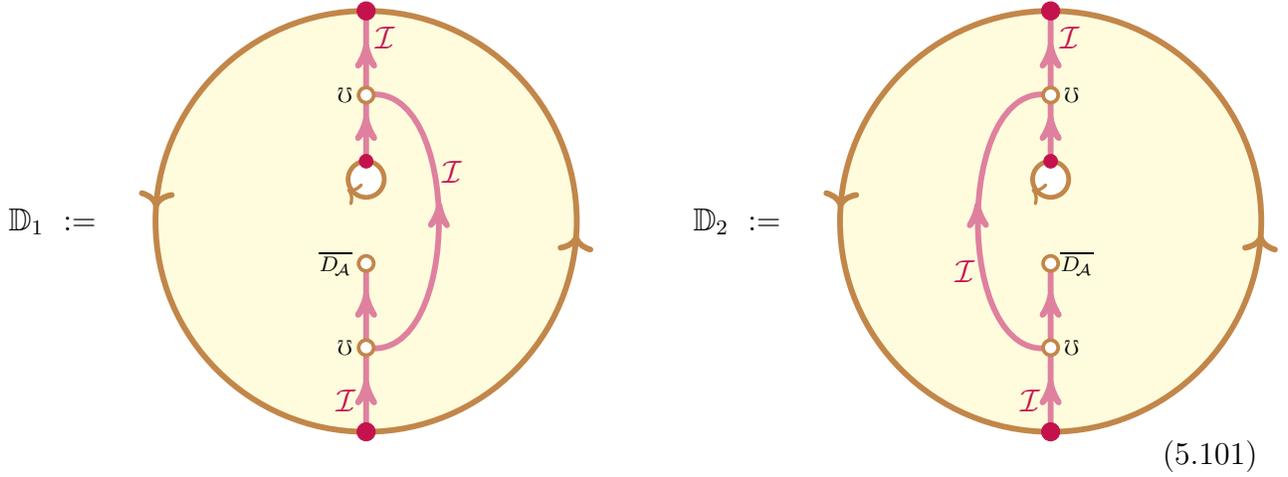
  \begin{eqnarray}
  \begin{tikzpicture}
  \begin{scope}[shift={(-0.6,0)}]
  \scopeArrow{0.49}{>} \scopeArrow{0.99}{>}
  \filldraw[\colorCircle,fill=\colorSigma,line width=\locpc,postaction={decorate}] (0,0) circle (\locpa);
  \end{scope} \end{scope}
  \scopeArrow{0.38}{\arrowDefect} \scopeArrow{0.81}{\arrowDefect}
  \draw[line width=\locpc,color=\colorTransp,postaction={decorate}] (90:\locpb)
             -- node[right=-1.6pt,yshift=6pt,near end,color=\colorDefect] {$\II$} (90:\locpa);
  \end{scope} \end{scope}
  \scopeArrow{0.31}{\arrowDefect} \scopeArrow{0.83}{\arrowDefect}
  \draw[line width=\locpc,color=\colorTransp,postaction={decorate}] (270:\locpa) 
             -- node[left=-1pt,near start,yshift=-4pt,color=\colorDefect] {$\II$} (270:\locpb);
  \end{scope} \end{scope}
  \scopeArrow{0.52}{\arrowDefectB}
  \draw[line width=\locpc,color=\colorTransp,postaction={decorate}]
             (90:\locpB) [out=10,in=-10] to (270:\locpB);
  \end{scope}
  \transv {270:\locpb} node[left=0.9pt,color=black] {$\!\scriptstyle \cc{D_\Cala}$};
  \transv {270:\locpB} node[left=0.5pt,color=black] {$\scriptstyle \om$};
  \transv {90:\locpB} node[left=0.5pt,color=black] {$\scriptstyle \om$};
  \filldraw[\colorCircle,fill=white,line width=\locpC] (90:\locpb) circle (\locpA);
  \scopeArrow{0.67}{<}
  \draw[\colorCircle,line width=0.7*\locpC,postaction={decorate}] (90:\locpb) circle (\locpA);
  \end{scope}
  \filldraw[color=\colorDefect] (90:\locpa) circle (\locpd) (270:\locpa) circle (\locpd)
                  (90:\locpb+\locpA) circle (\locpD);
  \node [color=\colorDefect] at (30:0.47*\locpa) {$\II$};
  \node at (-\locpa-1.37,0) {$\DD_1 ~:=$};
  \node at (0,-1.4*\locpa) {~};  
  \end{scope}
  \begin{scope}[shift={(2*\locpa+0.6,-3.1)}]
  \scopeArrow{0.49}{>} \scopeArrow{0.99}{>}
  \filldraw[\colorCircle,fill=\colorSigma,line width=\locpc,postaction={decorate}] (0,0) circle (\locpa);
  \end{scope} \end{scope}
  \scopeArrow{0.31}{\arrowDefect} \scopeArrow{0.83}{\arrowDefect}
  \draw[line width=\locpc,color=\colorTransp,postaction={decorate}] (270:\locpa)
             -- node[left=-1pt,near start,yshift=-4pt,color=\colorDefect] {$\II$} (270:\locpb);
  \end{scope} \end{scope}
  \scopeArrow{0.38}{\arrowDefect} \scopeArrow{0.81}{\arrowDefect}
  \draw[line width=\locpc,color=\colorTransp,postaction={decorate}] (90:\locpb) 
             -- node[right=-1.6pt,yshift=6pt,near end,color=\colorDefect] {$\II$} (90:\locpa);
  \end{scope} \end{scope}
  \scopeArrow{0.52}{\arrowDefectB}
  \draw[line width=\locpc,color=\colorTransp,postaction={decorate}]
             (90:\locpB) [out=170,in=190] to (270:\locpB);
  \end{scope}
  \transv {270:\locpb} node[right=1.3pt,color=black] {$\!\scriptstyle \cc{D_\Cala}$};
  \transv {270:\locpB} node[right=0.5pt,color=black] {$\scriptstyle \om$};
  \transv {90:\locpB} node[right=0.5pt,color=black] {$\scriptstyle \om$};
  \filldraw[\colorCircle,fill=white,line width=\locpC] (90:\locpb) circle (\locpA);
  \scopeArrow{0.67}{<}
  \draw[\colorCircle,line width=0.7*\locpC,postaction={decorate}] (90:\locpb) circle (\locpA);
  \end{scope}
  \filldraw[color=\colorDefect] (90:\locpa) circle (\locpd) (270:\locpa) circle (\locpd)
                  (90:\locpb+\locpA) circle (\locpD);
  \node [color=\colorDefect] at (210:0.47*\locpa) {$\II$};
  \node at (-\locpa-2.0,-0.7) {$\DD_2 ~:=$};
  \end{scope}
  \end{tikzpicture} \hspace*{-2.0em}
  \nonumber \\[-1.9em] ~ 
  \label{eq:Fig2-14.5}
  \\[-1.3em] ~ \nonumber
  \end{eqnarray}

We denote the outer and inner gluing circles (of both $\DD_1$ and $\DD_2$)
by $\LL$ and $\LL'$, respectively, and regard the block functors for $\DD_1$ and $\DD_2$
as functors from $\ZZ(\cc\LL)\boti \ZZ(\LL')$ to $\vect$.

\begin{Lemma} \label{Lemma:MCG-br}
\Itemizeiii

\item[{\rm (i)}]
The block functors for the \defectsurface s \eqref{eq:Fig2-14.5} are canonically isomorphic to 
  \be
  \bearll
  \Zfine(\DD_{1})(\cc{G} \boti z)=\Funle_{\cala}(G, F_{z}) \qquad\text{and}\qquad
  \Nxl6
  \Zfine(\DD_{2})(\cc{G} \boti z)=\Funle_{\cala}(G, {}_{z}F) \,,
  \label{eq:blocks-D1D2-a}
  \eear
  \ee
respectively, for $G \iN \Funle_{\Cala,\cala}(\cala,\cala) \,{=}\, \Zfine(\LL)$ and
$z \iN \calz(\cala)\,{=}\, \Zfine(\LL')$, with $F_{z}$ and ${}_{z}F$ the braided-induced
functors in $\Funle_{\Cala,\cala}(\cala,\cala)$, see \eqref{eq:braidedind}.

\item[{\rm (ii)}]
The braiding $F_{z}(y) \,{=}\, y \oti z \xrightarrow{c_{y,z}} z \oti y$
provides a bimodule  isomorphism $F_{z} \,{\cong}\, {}_{z}F$ and the parallelization
isomorphism $\Pi_{\DD; \DD_{1}, \DD_{2}}$ is given by post-composing with this  isomorphism.
\end{itemize}
\end{Lemma}

\begin{proof}
(i)\, Similarly to Example \ref{exa:braidedinduction} we use for 
$z \iN \calz(\cala)\,{=}\, \Zfine(\LL')$ the excision functors to obtain the situation 
of the left hand side of \eqref{eq:Fig6.6-1}
with object $F_{z}$ in the inner boundary. Thus the block functors are given by 
  \be
  \ZZ(\DD_{1})(\cc{G} \boti z)=\Funle_{\cala,\cala}(G, F_{z}) \qquad\text{and}\qquad
  \ZZ(\DD_{2})(\cc{G} \boti z)=\Funle_{\cala,\cala}(G, {}_{z}F) 
  \label{eq:9}
  \ee
for $G \iN \ZZ(\LL) \,{=}\, \Funle_{\cala,\cala}(\cala,\cala)$.
and $z \iN \ZZ(\LL') \,{=}\, \calz(\cala)$.
 \\[.2em]
(ii)\, To compute the natural isomorphism $\Pi_{\DD;\DD_1,\DD_2}\colon \ZZ(\DD_1) \,{\to}\, \ZZ(\DD_2)$,
we factor it through the disk $\DD$, which is possible because $\DD$ is already fine. This is
achieved by the \filldiskrep\ that replaces the disk $\DD_1$ by the left hand side of
    \def\locpa  {2.8}    
    \def\locpA  {2.4}    
    \def\locpb  {1.1}    
    \def\locpB  {0.56}   
    \def\locpc  {2.2pt}  
    \def\locpC  {1.8pt}  
    \def\locpd  {0.12}   
    \def\locpD  {0.09}   
    \def\locpe  {0.105}  
    \def\locpf  {1.1pt}  
    \def\locpr  {1.92}   
    \def\locpR  {0.52}   
    \def\locpt  {0.33}   
    \def\locpu  {208}    
    \def\locpv  {280}    
  \begin{eqnarray}
  \begin{tikzpicture}
  \scopeArrow{0.02}{>} \scopeArrow{0.52}{>}
  \filldraw[\colorCircle,fill=\colorSigma,line width=\locpc,postaction={decorate}]
             (0,0) circle (\locpa);
  \end{scope} \end{scope}
  \draw[\colorCircle,line width=\locpc,dash pattern=on 9pt off 5pt] (0,0) circle (\locpb);
  \scopeArrow{0.24}{\arrowDefectB} \scopeArrow{0.75}{\arrowDefectB}
  \draw[line width=\locpc,color=\colorTransp,postaction={decorate}] (\locpv:0.1) -- (\locpv:\locpa);
  \end{scope} \end{scope}
  \scopeArrow{0.82}{\arrowDefect}
  \draw[line width=\locpc,color=\colorTransp,postaction={decorate}]
             (\locpu:\locpr) [out=20,in=270] to (90:\locpa);
  \end{scope}
  \transv {\locpv:0.1} node[right,yshift=-1pt,color=black] {$\scriptstyle \cc{D_\Cala}$};
  \filldraw[\colorCircle,fill=white,line width=\locpc] (\locpu:\locpr) circle (\locpt);
  \scopeArrow{0.64}{<}
  \draw[\colorCircle,line width=\locpf,postaction={decorate}] (\locpu:\locpr) circle (\locpt);
  \end{scope}
  \filldraw[color=\colorDefect] (\locpv:\locpa) circle (\locpd) (90:\locpa) circle (\locpd)
             (\locpu:\locpr-\locpt) circle (\locpD) (\locpv:\locpb) circle (\locpe)
	     (206.3:\locpb) circle (\locpe) (97.3:\locpb) circle (\locpe);
  \node[color=\colorDefect] at (97:0.79*\locpa) {$\II$};
  \begin{scope}[shift={(2*\locpa+2.2,0)}]
  \scopeArrow{0.02}{>} \scopeArrow{0.52}{>}
  \filldraw[\colorCircle,fill=\colorSigma,line width=\locpc,postaction={decorate}]
             (0,0) circle (\locpA);
  \end{scope} \end{scope}
  \scopeArrow{0.58}{\arrowDefect}
  \draw[line width=\locpc,color=\colorTransp,postaction={decorate}] (280:\locpA) -- (280:\locpB);
  \end{scope}
  \scopeArrow{0.70}{\arrowDefect}
  \draw[line width=\locpc,color=\colorTransp,postaction={decorate}]
             (90:\locpR) -- node[right=-2.1pt,yshift=12pt,color=\colorDefect] {$\II$} (90:\locpA);
  \end{scope}
  \transv {280:\locpB} node[right=1.8pt,color=black] {$\!\scriptstyle \cc{D_\Cala}$};
  \filldraw[\colorCircle,fill=white,line width=\locpC] (90:\locpR) circle (\locpt);
  \scopeArrow{0.74}{<}
  \draw[\colorCircle,line width=0.7*\locpC,postaction={decorate}] (90:\locpR) circle (\locpt);
  \end{scope}
  \filldraw[color=\colorDefect] (280:\locpA) circle (\locpd) (90:\locpA) circle (\locpd)
                  (90:\locpR+\locpt) circle (\locpD);
  \node at (-\locpA-1.2,0) {{\Large$\simeq$}};
  \node at (\locpA+1.2,0) {$=~\DD \,.$};
  \end{scope}
  \end{tikzpicture}
  \nonumber \\[-4.9em]~ \\[1.7em] ~ \nonumber
  \end{eqnarray}
By Lemma \ref{Lemma:block-iso-disk} there is a canonical isomorphism $\ZZ(\DD)(\cc{G} \boti z)
\,{\tocong}\, \Hom_{\calz(\cala)}(G(\one), z)$ and the parallelization isomorphism
$\Pi_{\DD; \DD_{1}, \DD} \colon \ZZ(\DD_{1}) \,{\to}\, \ZZ(\DD)$ is given by 
  \be
  \Funle_{\cala,\cala}(G, F_{z}) \,\ni \mu
  \,\longmapsto\, \mu_\one \in\, \Hom_{\calz(\cala)}(G(\one), z) \,.
  \label{eq:nat-iso-to-D}
  \ee
Analogously $\Pi_{\DD; \DD_{2}, \DD} \colon \ZZ(\DD_{2}) \,{\tocong}\, \ZZ(\DD)$ is 
given by evaluation on $\one \iN \cala$ which has the braided induction ${}_zF$
as a quasi-inverse functor. The claimed expression for $\Pi_{\DD; \DD_{1}, \DD_{2}}$
now follows directly from
$ \Pi_{\DD; \DD_1, \DD_2}  \,{=}\, \Pi_{\DD; \DD_2, \DD}^{-1} \cir \Pi_{\DD; \DD_1, \DD}$.
\end{proof}

Consider now a braiding move on a three-punctured sphere $(\Sigma,\varrho)$ with boundary 
parame\-tri\-za\-tion $\varrho$, i.e.\ a diffeomorphism §$\varphi \colon (\Sigma,\varrho,\chi) \,{\to}\,
(\Sigma, \tau\cir\varrho, \chi)$ that changes the boundary parame\-tri\-za\-tion to $\tau \cir \varrho$,
where $\tau$ is the symmetric monoidal braiding on $\Bfrdefdec$ as  indicated in
 \\[-1.5em]
    \def\locpa  {3.5}    
    \def\locpA  {0.49}   
    \def\locpb  {0.37}   
    \def\locpB  {55}     
    \def\locpc  {2.2pt}  
    \def\locpC  {1.6pt}  
    \def\locpd  {0.12}   
    \def\locpD  {0.10}   
    \def\locpe  {1.27}   
    \def\locpE  {1.04}   
    \def\locpr  {green!52!black}   
    \def\locpy  {0.7*\locpa}  
    \def\locpz  {1.11}   
  \begin{eqnarray}
  \begin{tikzpicture}
  \coordinate (circ1) at (270-\locpB:\locpb*\locpa);
  \coordinate (circ2) at (270+\locpB:\locpb*\locpa);
  \coordinate (dashpt1) at ({\locpa*cos(143)},{-\locpy*sin(143)});
  \coordinate (dashpT1) at ({\locpa*cos(143)},{\locpy*sin(143)});
  \coordinate (dashpt2) at ({\locpa*cos(125)},{-\locpy*sin(125)});
  \coordinate (dashpT2) at ({\locpa*cos(125)},{\locpy*sin(125)});
  \coordinate (dashpt3) at ({\locpa*cos(110)},{-\locpy*sin(110)});
  \coordinate (dashpT3) at ({\locpa*cos(110)},{\locpy*sin(110)});
  \coordinate (dashpt4) at ({\locpa*cos(96)},{-\locpy*sin(96)});
  \coordinate (dashpT4) at ({\locpa*cos(96)},{\locpy*sin(96)});
  \coordinate (dashpt5) at ({\locpa*cos(82)},{-\locpy*sin(82)});
  \coordinate (dashpT5) at ({\locpa*cos(82)},{\locpy*sin(82)});
  \coordinate (dashpt6) at ({\locpa*cos(68)},{-\locpy*sin(68)});
  \coordinate (dashpT6) at ({\locpa*cos(68)},{\locpy*sin(68)});
  \coordinate (dashpt7) at ({\locpa*cos(53)},{-\locpy*sin(53)});
  \coordinate (dashpT7) at ({\locpa*cos(53)},{\locpy*sin(53)});
  \coordinate (dashpt8) at ({\locpa*cos(35)},{-\locpy*sin(35)});
  \coordinate (dashpT8) at ({\locpa*cos(35)},{\locpy*sin(35)});
  \coordinate (circ3ou) at ({\locpa*cos(44)},{-1.08*\locpy*sin(44)});
  \scopeArrow{0.01}{>}
  \filldraw[\colorCircle,fill=\colorSigma,line width=\locpc,postaction={decorate}]
           (0,0) circle (\locpa cm and \locpy cm);
  \end{scope}
  \scopeArrow{0.75}{>}
  \draw[blue!50!black,dashed,line width=0.5*\locpc,postaction={decorate}] (dashpt1) -- (dashpT1);
  \draw[blue!50!black,dashed,line width=0.5*\locpc,postaction={decorate}] (dashpt2) -- (dashpT2);
  \draw[blue!50!black,dashed,line width=0.5*\locpc,postaction={decorate}] (dashpt3) -- (dashpT3);
  \draw[blue!50!black,dashed,line width=0.5*\locpc,postaction={decorate}] (dashpt4) -- (dashpT4);
  \draw[blue!50!black,dashed,line width=0.5*\locpc,postaction={decorate}] (dashpt5) -- (dashpT5);
  \draw[blue!50!black,dashed,line width=0.5*\locpc,postaction={decorate}] (dashpt6) -- (dashpT6);
  \draw[blue!50!black,dashed,line width=0.5*\locpc,postaction={decorate}] (dashpt7) -- (dashpT7);
  \draw[blue!50!black,dashed,line width=0.5*\locpc,postaction={decorate}] (dashpt8) -- (dashpT8);
  \end{scope}
  \scopeArrow{0.80}{\arrowDefect}
  \draw[line width=\locpc,color=\colorTransp,postaction={decorate}] (circ1) -- ++(0,\locpe);
  \draw[line width=\locpc,color=\colorTransp,postaction={decorate}] (circ2) -- ++(0,\locpe);
  \end{scope}
  \scopeArrow{0.72}{\arrowDefect}
  \draw[line width=\locpc,color=\colorTransp,postaction={decorate}]
           (0,\locpy-\locpE) -- ++(0,\locpE);
  \end{scope}
  \scopeArrow{0.80}{<}
  \filldraw[\colorCircle,fill=white,line width=\locpC,postaction={decorate}] (circ1) circle (\locpA);
  \filldraw[\colorCircle,fill=white,line width=\locpC,postaction={decorate}] (circ2) circle (\locpA);
  \end{scope}
  \draw[\colorCircle,line width=\locpc] (circ1) circle (\locpA) (circ2) circle (\locpA);
  \transv {0,\locpy-\locpE};
  \begin{scope}[shift={(circ1)}] \transv {0,\locpe}; \end{scope}
  \begin{scope}[shift={(circ2)}] \transv {0,\locpe}; \end{scope}
  \filldraw[color=\colorDefect] (0,\locpy) circle (\locpd) 
	   (circ1)++(0,\locpA) circle (\locpD) (circ2)++(0,\locpA) circle (\locpD);
  \scopeArrow{0.86}{triangle 45}
  \draw[\locpr,line width=1.3*\locpc,postaction={decorate}] (circ1)++(230:1.2*\locpA)
           [out=240,in=110] to node[left=-2.4pt,color=\locpr] {{$\varrho$}}
           (-3.5*\locpA,-\locpa-\locpz);
  \end{scope}
  \scopeArrow{0.86}{triangle 45}
  \draw[\locpr,line width=1.3*\locpc,postaction={decorate}] (circ2)++(230:1.2*\locpA)
           [out=220,in=110] to node[left=-1pt,yshift=-7pt,color=\locpr] {{$\varrho$}}
           (0,-\locpa-\locpz);
  \end{scope}
  \scopeArrow{0.82}{triangle 45}
  \draw[\locpr,line width=1.3*\locpc,postaction={decorate}] (circ3ou)
           [out=290,in=70] to node[left=-6pt,yshift=13pt,color=\locpr] {{$\varrho$}}
           (3.5*\locpA,-\locpa-\locpz);
  \end{scope}
  \begin{scope}[shift={(0,-\locpa-\locpz)}]
  \fill[white] (-3.5*\locpA,0) circle (1.2*\locpA)
           (0,0) circle (1.2*\locpA) (3.5*\locpA,0) circle (1.2*\locpA);
  \scopeArrow{0.80}{<}
  \draw[\colorCircle,line width=\locpC,postaction={decorate}] (-3.5*\locpA,0) circle (\locpA);
  \draw[\colorCircle,line width=\locpC,postaction={decorate}] (0,0) circle (\locpA);
  \end{scope}
  \scopeArrow{0.01}{>}
  \draw[\colorCircle,line width=\locpC,postaction={decorate}] (3.5*\locpA,0) circle (\locpA);
  \end{scope}
  \draw[\colorCircle,line width=\locpc] (-3.5*\locpA,0) circle (\locpA)
           (0,0) circle (\locpA) (3.5*\locpA,0) circle (\locpA);
  \end{scope}
  \node at (\locpa+1.0,0) {\large{\boldmath $\xrightarrow{~\,\varphi\,~}$}};
  \begin{scope}[shift={(2*\locpa+1.9,0)}]
  \coordinate (circ1) at (270-\locpB:\locpb*\locpa);
  \coordinate (circ2) at (270+\locpB:\locpb*\locpa);
  \coordinate (dashpt1) at ({\locpa*cos(143)},{-\locpy*sin(143)});
  \coordinate (dashpT1) at ({\locpa*cos(143)},{\locpy*sin(143)});
  \coordinate (dashpt2) at ({\locpa*cos(125)},{-\locpy*sin(125)});
  \coordinate (dashpT2) at ({\locpa*cos(125)},{\locpy*sin(125)});
  \coordinate (dashpt3) at ({\locpa*cos(110)},{-\locpy*sin(110)});
  \coordinate (dashpT3) at ({\locpa*cos(110)},{\locpy*sin(110)});
  \coordinate (dashpt4) at ({\locpa*cos(96)},{-\locpy*sin(96)});
  \coordinate (dashpT4) at ({\locpa*cos(96)},{\locpy*sin(96)});
  \coordinate (dashpt5) at ({\locpa*cos(82)},{-\locpy*sin(82)});
  \coordinate (dashpT5) at ({\locpa*cos(82)},{\locpy*sin(82)});
  \coordinate (dashpt6) at ({\locpa*cos(68)},{-\locpy*sin(68)});
  \coordinate (dashpT6) at ({\locpa*cos(68)},{\locpy*sin(68)});
  \coordinate (dashpt7) at ({\locpa*cos(53)},{-\locpy*sin(53)});
  \coordinate (dashpT7) at ({\locpa*cos(53)},{\locpy*sin(53)});
  \coordinate (dashpt8) at ({\locpa*cos(35)},{-\locpy*sin(35)});
  \coordinate (dashpT8) at ({\locpa*cos(35)},{\locpy*sin(35)});
  \coordinate (circ3ou) at ({\locpa*cos(44)},{-1.08*\locpy*sin(44)});
  \scopeArrow{0.01}{>}
  \filldraw[\colorCircle,fill=\colorSigma,line width=\locpc,postaction={decorate}]
           (0,0) circle (\locpa cm and \locpy cm);
  \end{scope}
  \scopeArrow{0.85}{>}
  \draw[blue!50!black,dashed,line width=0.5*\locpc,postaction={decorate}]
       plot[smooth,tension=0.65] coordinates{
           (dashpt1) (196:0.91*\locpa) (180:0.90*\locpa) (160:0.85*\locpa) (dashpT1) };
  \draw[blue!50!black,dashed,line width=0.5*\locpc,postaction={decorate}]
       plot[smooth,tension=0.65] coordinates{
           (dashpt2) (218:0.80*\locpa) (200:0.83*\locpa) 
           (180:0.80*\locpa) (160:0.75*\locpa) (145:0.73*\locpa) (dashpT2) };
  \draw[blue!50!black,dashed,line width=0.5*\locpc,postaction={decorate}]
       plot[smooth,tension=0.65] coordinates{
           (dashpt3) (250:0.66*\locpa) (250:0.62*\locpa) (240:0.67*\locpa) (230:0.73*\locpa)
           (210:0.74*\locpa) (187:0.70*\locpa) (155:0.63*\locpa) (120:0.55*\locpa) (dashpT3) };
  \draw[blue!50!black,dashed,line width=0.5*\locpc,postaction={decorate}]
       plot[smooth,tension=0.65] coordinates{
           (dashpt4) (267:0.57*\locpa) (300:0.57*\locpa) (325:0.59*\locpa) (340:0.57*\locpa)
           (6:0.49*\locpa) (50:0.31*\locpa) (262:0.44*\locpa) (225:0.68*\locpa) (190:0.62*\locpa)
           (150:0.53*\locpa) (121:0.42*\locpa) (102:0.55*\locpa) (dashpT4) };
  \draw[blue!50!black,dashed,line width=0.5*\locpc,postaction={decorate}]
       plot[smooth,tension=0.65] coordinates{
           (dashpt5) (289:0.62*\locpa) (320:0.67*\locpa) (1:0.61*\locpa) (57:0.39*\locpa)
           (255:0.35*\locpa) (225:0.61*\locpa) (190:0.54*\locpa) (155:0.44*\locpa)
           (90:0.28*\locpa) (77:0.54*\locpa) (dashpT5) };
  \draw[blue!50!black,dashed,line width=0.5*\locpc,postaction={decorate}]
       plot[smooth,tension=0.65] coordinates{
           (dashpt7) (327:0.83*\locpa) (342:0.81*\locpa)
           (9:0.78*\locpa) (43:0.70*\locpa) (dashpT7) };
  \end{scope}
  \scopeArrow{0.75}{>}
  \draw[blue!50!black,dashed,line width=0.5*\locpc,postaction={decorate}]
       plot[smooth,tension=0.65] coordinates{
           (dashpt6) (307:0.71*\locpa) (325:0.75*\locpa) (0:0.71*\locpa) (27:0.61*\locpa)
           (49:0.54*\locpa) (65:0.51*\locpa) (65:0.60*\locpa) (60:0.67*\locpa) (dashpT6) };
  \draw[blue!50!black,dashed,line width=0.5*\locpc,postaction={decorate}]
       plot[smooth,tension=0.65] coordinates{
           (dashpt8) (339:0.89*\locpa) (1:0.88*\locpa) (19:0.83*\locpa) (dashpT8) };
  \end{scope}
  \scopeArrow{0.80}{\arrowDefect}
  \draw[line width=\locpc,color=\colorTransp,postaction={decorate}] (circ1) -- ++(0,\locpe);
  \draw[line width=\locpc,color=\colorTransp,postaction={decorate}] (circ2) -- ++(0,\locpe);
  \end{scope}
  \scopeArrow{0.72}{\arrowDefect}
  \draw[line width=\locpc,color=\colorTransp,postaction={decorate}]
           (0,\locpy-\locpE) -- ++(0,\locpE);
  \end{scope}
  \scopeArrow{0.80}{<}
  \filldraw[\colorCircle,fill=white,line width=\locpC,postaction={decorate}] (circ1) circle (\locpA);
  \filldraw[\colorCircle,fill=white,line width=\locpC,postaction={decorate}] (circ2) circle (\locpA);
  \end{scope}
  \draw[\colorCircle,line width=\locpc] (circ1) circle (\locpA) (circ2) circle (\locpA);
  \transv {0,\locpy-\locpE};
  \begin{scope}[shift={(circ1)}] \transv {0,\locpe}; \end{scope}
  \begin{scope}[shift={(circ2)}] \transv {0,\locpe}; \end{scope}
  \filldraw[color=\colorDefect] (0,\locpy) circle (\locpd) 
           (circ1)++(0,\locpA) circle (\locpD) (circ2)++(0,\locpA) circle (\locpD);
  \scopeArrow{0.87}{triangle 45}
  \draw[\locpr,line width=1.3*\locpc,postaction={decorate}] (circ2)++(240:1.2*\locpA) [out=250,in=84]
           to node[right=-1pt,yshift=-3pt,color=\locpr] {{$\tau{\circ}\varrho$}}
           (-3.5*\locpA,-\locpa-\locpz);
  \end{scope}
  \scopeArrow{0.86}{triangle 45}
  \draw[\locpr,line width=1.3*\locpc,postaction={decorate}] (circ1)++(250:1.2*\locpA) [out=250,in=110]
           to node[left=1pt,yshift=2pt,color=\locpr] {{$\tau{\circ}\varrho$}}
           (0,-\locpa-\locpz);
  \end{scope}
  \scopeArrow{0.82}{triangle 45}
  \draw[\locpr,line width=1.3*\locpc,postaction={decorate}] (circ3ou) [out=290,in=70]
           to node[right,yshift=6pt,color=\locpr] {{$\tau{\circ}\varrho$}}
           (3.5*\locpA,-\locpa-\locpz);
  \end{scope}
  \end{scope}
  \begin{scope}[shift={(2*\locpa+1.9,-\locpa-\locpz)}]
  \fill[white] (-3.5*\locpA,0) circle (1.2*\locpA)
           (0,0) circle (1.2*\locpA) (3.5*\locpA,0) circle (1.2*\locpA);
  \scopeArrow{0.80}{<}
  \draw[\colorCircle,line width=\locpC,postaction={decorate}] (-3.5*\locpA,0) circle (\locpA);
  \draw[\colorCircle,line width=\locpC,postaction={decorate}] (0,0) circle (\locpA);
  \end{scope}
  \scopeArrow{0.01}{>}
  \draw[\colorCircle,line width=\locpC,postaction={decorate}] (3.5*\locpA,0) circle (\locpA);
  \end{scope}
  \draw[\colorCircle,line width=\locpc] (-3.5*\locpA,0) circle (\locpA)
           (0,0) circle (\locpA) (3.5*\locpA,0) circle (\locpA);
  \end{scope}
  \end{tikzpicture} \hspace*{-2.4em}
  \nonumber\\[-0.7em]~\\[-1.9em]~\nonumber
  \end{eqnarray}
Note that the framing $\varphi_{*}(\chi)$ is homotopic to $\chi$. Thus we can consider 
the two refinements of $(\Sigma, \varrho, \chi; \Sigma_{1})$ and 
$(\Sigma, \tau \circ \varrho, \chi; \Sigma_{2})$ depicted in 
    \def\locpa  {3.7}    
    \def\locpA  {0.49}   
    \def\locpb  {0.61}   
    \def\locpB  {50}     
    \def\locpc  {2.2pt}  
    \def\locpC  {1.6pt}  
    \def\locpd  {0.12}   
    \def\locpD  {0.10}   
    \def\locpe  {2.06}   
    \def\locpE  {1.84}   
    \def\locpy  {0.85*\locpa}  
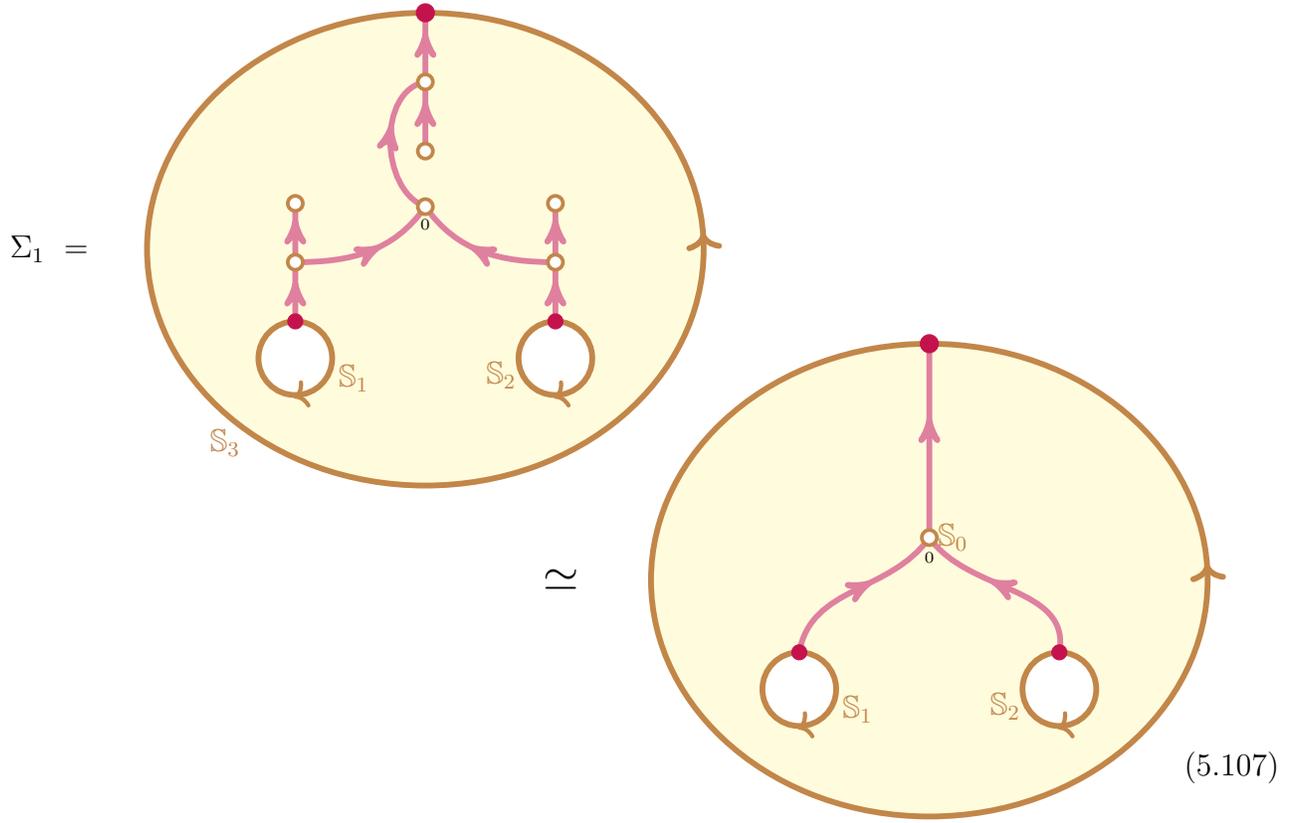
  \begin{eqnarray}
  \begin{tikzpicture}
  \coordinate (circ1) at (270-\locpB:\locpb*\locpa);
  \coordinate (circ2) at (270+\locpB:\locpb*\locpa);
  \scopeArrow{0.01}{>}
  \filldraw[\colorCircle,fill=\colorSigma,line width=\locpc,postaction={decorate}]
           (0,0) circle (\locpa cm and \locpy cm);
  \end{scope}
  \scopeArrow{0.50}{\arrowDefect} \scopeArrow{0.90}{\arrowDefect}
  \draw[line width=\locpc,color=\colorTransp,postaction={decorate}] (circ1) -- ++(0,\locpe);
  \draw[line width=\locpc,color=\colorTransp,postaction={decorate}] (circ2) -- ++(0,\locpe);
  \end{scope} \end{scope}
  \scopeArrow{0.36}{\arrowDefect} \scopeArrow{0.86}{\arrowDefect}
  \draw[line width=\locpc,color=\colorTransp,postaction={decorate}]
           (0,\locpy-\locpE) -- ++(0,\locpE);
  \end{scope} \end{scope}
  \scopeArrow{0.59}{\arrowDefect}
  \draw[line width=\locpc,color=\colorTransp,postaction={decorate}] 
	   (circ1)++(0,0.62*\locpe) [out=0,in=235] to (0,0.18*\locpy);
  \draw[line width=\locpc,color=\colorTransp,postaction={decorate}] 
	   (circ2)++(0,0.62*\locpe) [out=180,in=305] to (0,0.18*\locpy);
  \end{scope}
  \scopeArrow{0.62}{\arrowDefect}
  \draw[line width=\locpc,color=\colorTransp,postaction={decorate}] 
	   (0,0.18*\locpy) [out=160,in=190] to (0,\locpy-0.5*\locpE);
  \end{scope}
  \scopeArrow{0.8}{<}
  \filldraw[\colorCircle,fill=white,line width=\locpC,postaction={decorate}] (circ1) circle (\locpA);
  \filldraw[\colorCircle,fill=white,line width=\locpC,postaction={decorate}] (circ2) circle (\locpA);
  \end{scope}
  \draw[\colorCircle,line width=\locpc] (circ1) circle (\locpA) (circ2) circle (\locpA);
  \transv {0,0.18*\locpy} node[below,color=black] {$\scriptscriptstyle 0$};
  \transv {0,\locpy-\locpE};
  \transv {0,\locpy-0.5*\locpE};
  \begin{scope}[shift={(circ1)}] \transv {0,\locpe}; \transv {0,0.62*\locpe}; \end{scope}
  \begin{scope}[shift={(circ2)}] \transv {0,\locpe}; \transv {0,0.62*\locpe}; \end{scope}
  \filldraw[color=\colorDefect] (0,\locpy) circle (\locpd) 
	   (circ1)++(0,\locpA) circle (\locpD) (circ2)++(0,\locpA) circle (\locpD);
  \node[\colorCircle] at (241:0.53*\locpa) {$\SS_1$};
  \node[\colorCircle] at (301:0.53*\locpa) {$\SS_2$};
  \node[\colorCircle] at (224:\locpa) {$\SS_3$};
  \node at (-1.3-\locpa,0) {$\Sigma_1 ~=$};
  \begin{scope}[shift={(2*\locpa-1.5,-1.4*\locpy)}]
  \coordinate (circ1) at (270-\locpB:\locpb*\locpa);
  \coordinate (circ2) at (270+\locpB:\locpb*\locpa);
  \scopeArrow{0.01}{>}
  \filldraw[\colorCircle,fill=\colorSigma,line width=\locpc,postaction={decorate}]
           (0,0) circle (\locpa cm and \locpy cm);
  \end{scope}
  \scopeArrow{0.59}{\arrowDefect}
  \draw[line width=\locpc,color=\colorTransp,postaction={decorate}] 
           (circ1)++(0,\locpA) [out=80,in=235] to (0,0.18*\locpy);
  \draw[line width=\locpc,color=\colorTransp,postaction={decorate}] 
           (circ2)++(0,\locpA) [out=80,in=305] to (0,0.18*\locpy);
  \end{scope}
  \scopeArrow{0.62}{\arrowDefect}
  \draw[line width=\locpc,color=\colorTransp,postaction={decorate}] 
	   (0,0.18*\locpy) -- (0,\locpy);
  \end{scope}
  \scopeArrow{0.8}{<}
  \filldraw[\colorCircle,fill=white,line width=\locpC,postaction={decorate}] (circ1) circle (\locpA);
  \filldraw[\colorCircle,fill=white,line width=\locpC,postaction={decorate}] (circ2) circle (\locpA);
  \end{scope}
  \draw[\colorCircle,line width=\locpc] (circ1) circle (\locpA) (circ2) circle (\locpA);
  \transv {0,0.18*\locpy} node[below=1pt,color=black] {$\scriptscriptstyle 0$};
  \filldraw[color=\colorDefect] (0,\locpy) circle (\locpd) 
	   (circ1)++(0,\locpA) circle (\locpD) (circ2)++(0,\locpA) circle (\locpD);
  \node at (-1.2-\locpa,0) {{\Large $\simeq$}};
  \node[\colorCircle,xshift=9pt] at (0,0.18*\locpy) {$\SS_0$};
  \node[\colorCircle] at (241:0.53*\locpa) {$\SS_1$};
  \node[\colorCircle] at (301:0.53*\locpa) {$\SS_2$};
  \end{scope}
  \end{tikzpicture}
  \nonumber \\[-2.7em]~
  \label{eq:reduce-brd}
  \\[-5.3em] ~\nonumber
  \end{eqnarray}
and
    \def\locpa  {3.7}    
    \def\locpA  {0.49}   
    \def\locpb  {0.61}   
    \def\locpB  {50}     
    \def\locpc  {2.2pt}  
    \def\locpC  {1.6pt}  
    \def\locpd  {0.12}   
    \def\locpD  {0.10}   
    \def\locpe  {2.06}   
    \def\locpE  {1.84}   
    \def\locpr  {green!52!black}   
    \def\locpy  {0.85*\locpa}  
    \def\locpz  {1.22}   
  \be
  \raisebox{-12.9em}
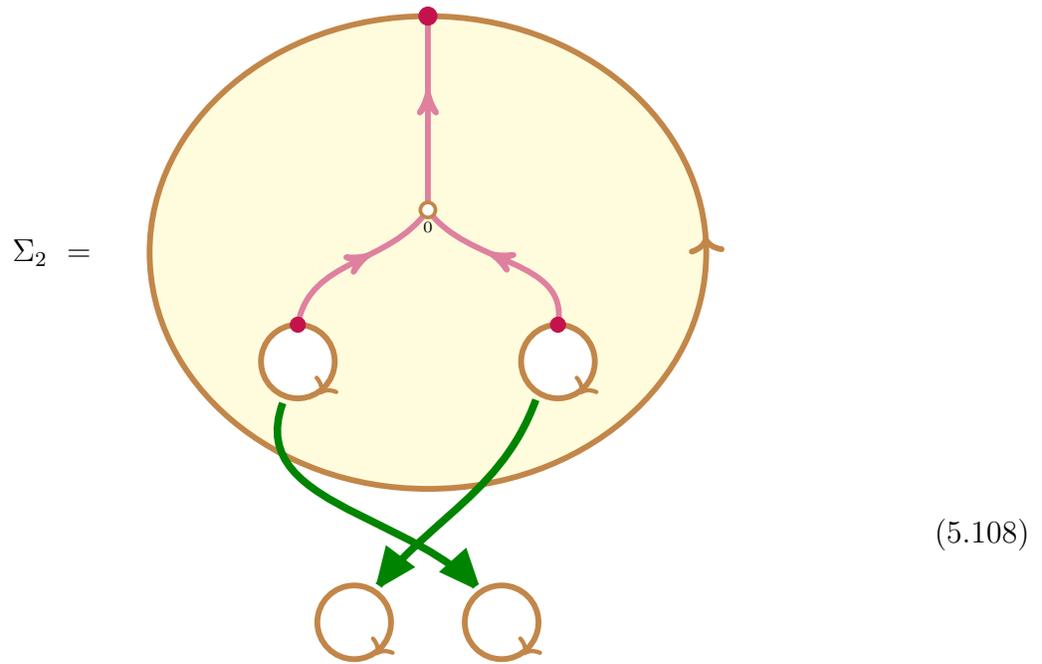
{\begin{tikzpicture}
  \coordinate (circ1) at (270-\locpB:\locpb*\locpa);
  \coordinate (circ2) at (270+\locpB:\locpb*\locpa);
  \scopeArrow{0.01}{>}
  \filldraw[\colorCircle,fill=\colorSigma,line width=\locpc,postaction={decorate}]
           (0,0) circle (\locpa cm and \locpy cm);
  \end{scope}
  \scopeArrow{0.59}{\arrowDefect}
  \draw[line width=\locpc,color=\colorTransp,postaction={decorate}] 
           (circ1)++(0,\locpA) [out=80,in=235] to (0,0.18*\locpy);
  \draw[line width=\locpc,color=\colorTransp,postaction={decorate}] 
           (circ2)++(0,\locpA) [out=80,in=305] to (0,0.18*\locpy);
  \end{scope}
  \scopeArrow{0.62}{\arrowDefect}
  \draw[line width=\locpc,color=\colorTransp,postaction={decorate}] 
	   (0,0.18*\locpy) -- (0,\locpy);
  \end{scope}
  \scopeArrow{0.90}{<}
  \filldraw[\colorCircle,fill=white,line width=\locpC,postaction={decorate}] (circ1) circle (\locpA);
  \filldraw[\colorCircle,fill=white,line width=\locpC,postaction={decorate}] (circ2) circle (\locpA);
  \end{scope}
  \draw[\colorCircle,line width=\locpc] (circ1) circle (\locpA) (circ2) circle (\locpA);
  \transv {0,0.18*\locpy} node[below,color=black] {$\scriptscriptstyle 0$};
  \filldraw[color=\colorDefect] (0,\locpy) circle (\locpd) 
	   (circ1)++(0,\locpA) circle (\locpD) (circ2)++(0,\locpA) circle (\locpD);
  \node at (-1.3-\locpa,0) {$\Sigma_2 ~=$};
  \scopeArrow{0.87}{triangle 45}
  \draw[\locpr,line width=1.3*\locpc,postaction={decorate}] (circ2)++(240:1.2*\locpA) [out=250,in=61]
           to (-2*\locpA,-\locpa-\locpz);
  \end{scope}
  \scopeArrow{0.89}{triangle 45}
  \draw[\locpr,line width=1.3*\locpc,postaction={decorate}] (circ1)++(250:1.2*\locpA) [out=250,in=115]
           to (2*\locpA,-\locpa-\locpz);
  \end{scope}
  \begin{scope}[shift={(0,-\locpa-\locpz)}]
  \fill[white] (-2*\locpA,0) circle (1.2*\locpA) (2*\locpA,0) circle (1.2*\locpA);
  \scopeArrow{0.90}{<}
  \draw[\colorCircle,line width=\locpC,postaction={decorate}] (-2*\locpA,0) circle (\locpA);
  \draw[\colorCircle,line width=\locpC,postaction={decorate}] (2*\locpA,0) circle (\locpA);
  \end{scope}
  \draw[\colorCircle,line width=\locpc] (-2*\locpA,0) circle (\locpA) (2*\locpA,0) circle (\locpA);
  \end{scope}
  \end{tikzpicture}}~~
  \ee
respectively. The gluing categories for the inner gluing circles $\SS_1$ and $\SS_2$ are 
$\calz(\cala)$, while the one for the outer gluing circle $\SS_3$ is $\cc{\calz(\cala)}$.

\begin{Proposition} \label{prop:braid}
The block functors for $\Sigma_{1}$ and $\Sigma_{2}$ evaluated at the objects
$x \iN \ZZ(\SS_{1}) \,{=}$
                         $\calz(\cala)$, $y \iN \ZZ(\SS_{2}) \,{=}\, \calz(\cala)$ and
$\cc{z} \iN \ZZ(\SS_{3}) \,{=}\, \cc{\calz(\cala)}$ are
  \be
  \bearl
  \ZZ(\Sigma_1)(x \boti y \boti \cc{z}) = \Hom_{\calz(\cala)}(z,x \oti y) \qquad \text{and}
  \Nxl6
  \ZZ(\Sigma_2)(x \boti y \boti \cc{z}) = \Hom_{\calz(\cala)}(z, y \oti x) \,,
  \label{eq:blocks-sigmai}
  \eear
  \ee
respectively. The isomorphism corresponding to the value $\ZZ(\varphi)$ of the block functor
on the diffeomorphism $\varphi$ is given by composition with the braiding of $x$ and $y$. 
\end{Proposition}

\begin{proof}
First we use again the excision functors to reduce the situation to the one on the 
right hand side of Equation \eqref{eq:reduce-brd}. We then contract the objects in the gluing
categories for the circles $\SS_{1}$ and $\SS_{2}$ in $\Sigma_1$ with the \mute\ object 
for the transparent gluing circle $\SS_0$ in the diagram on the right hand side of 
\eqref{eq:reduce-brd}; the latter circle is the one denoted by $\LII_0^{\sss\swarrow}(\calm)$,
with $\calm$ set to $\cala$, in \eqref{eq:circle.IIn-M-M}, and its \mute\ object is recorded
in Equation \eqref{eq:mute-expl}. Hereby we obtain the object 
  \be
  \int^{b\in\cala}\!\!\! \int^{a\in\cala}\! \Hom_{\cala}(a,x) \otik \Hom_{\cala}(b,y) \otimes a.b
  \,\cong\, x \oti y ~\iN \calz(\cala) \,. 
  \label{eq:contract} 
  \ee
This gives the claimed expression for the block functor on $\Sigma_{1}$. Analogously we compute
the block functor for $\Sigma_{2}$. Next consider the result of applying the diffeomorphism 
$\varphi$ to the refined surface $\Sigma_{1}$, which results in
the refinement $\varphi(\Sigma_{1})$ in the first of the following pictures:
    \def\locpa  {3.7}    
    \def\locpA  {0.49}   
    \def\locpb  {0.39}   
    \def\locpB  {78}     
    \def\locpc  {2.2pt}  
    \def\locpC  {1.6pt}  
    \def\locpd  {0.12}   
    \def\locpD  {0.10}   
    \def\locpe  {1.85}   
    \def\locpy  {0.85*\locpa}  
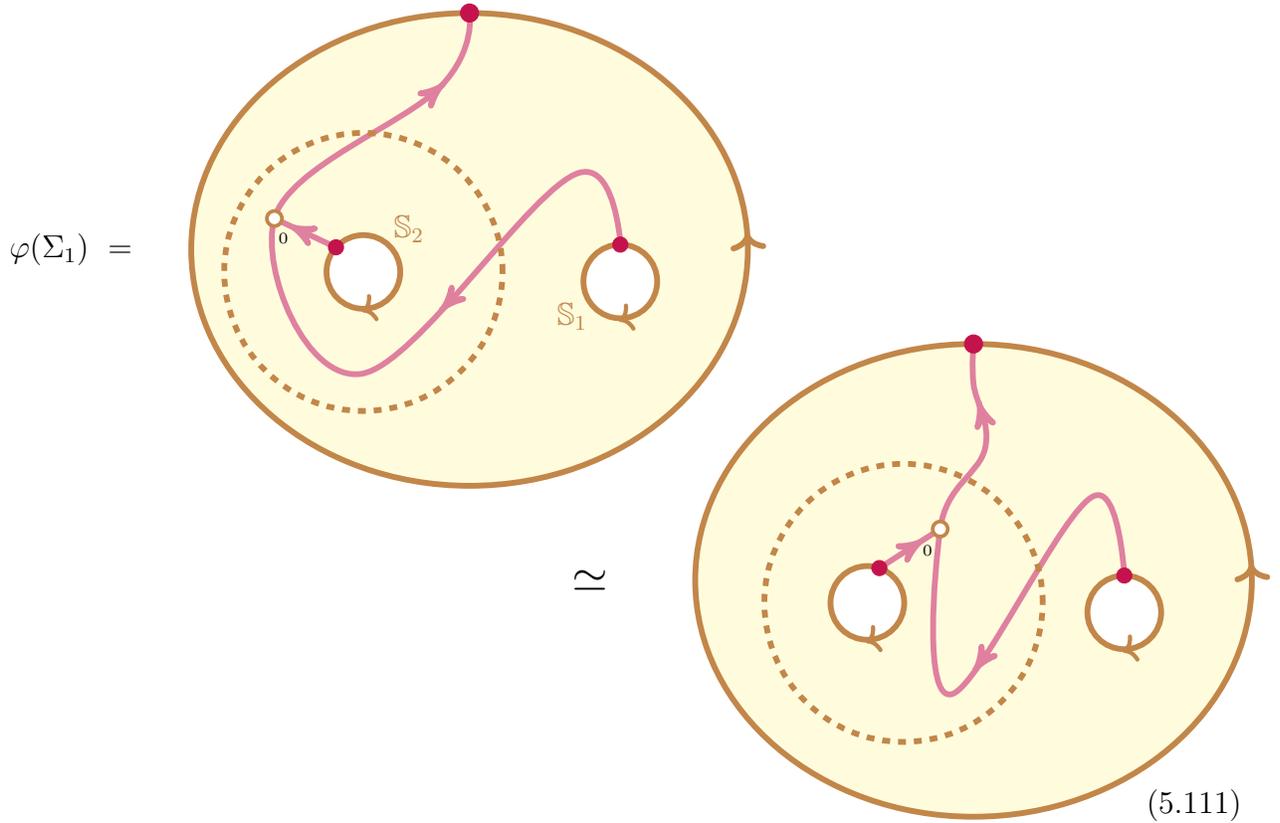
  \begin{eqnarray}
  \begin{tikzpicture}
  \coordinate (circ1) at (270-\locpB:\locpb*\locpa);
  \coordinate (Circ1) at (307-\locpB:1.52*\locpb*\locpa);
  \coordinate (circ2) at (270+\locpB:1.42*\locpb*\locpa);
  \coordinate (Circ2) at (265+\locpB:1.48*\locpb*\locpa);
  \coordinate (postr) at (171:0.71*\locpa);
  \scopeArrow{0.01}{>}
  \filldraw[\colorCircle,fill=\colorSigma,line width=\locpc,postaction={decorate}]
           (0,0) circle (\locpa cm and \locpy cm);
  \end{scope}
  \scopeArrow{0.29}{\arrowDefect}
  \draw[line width=\locpc,color=\colorTransp] 
           plot[smooth,tension=0.65] coordinates{
           (Circ2) (35:0.48*\locpa) (Circ1) (postr) (100:0.7*\locpy) (0,\locpy) };
  \end{scope}
  \scopeArrow{0.73}{\arrowDefect}
  \draw[line width=\locpc,color=\colorTransp,postaction={decorate}] (circ1)++(138:\locpA) -- (postr);
  \draw[line width=\locpc,color=\colorTransp,postaction={decorate}] (100:0.7*\locpy) -- ++ (41:0.05);
  \draw[line width=\locpc,color=\colorTransp,postaction={decorate}] (246:0.27*\locpy) -- ++ (226:0.05);
  \end{scope}
  \scopeArrow{0.8}{<}
  \filldraw[\colorCircle,fill=white,line width=\locpC,postaction={decorate}] (circ1) circle (\locpA);
  \filldraw[\colorCircle,fill=white,line width=\locpC,postaction={decorate}] (circ2) circle (\locpA);
  \end{scope}
  \draw[\colorCircle,line width=\locpc] (circ1) circle (\locpA) (circ2) circle (\locpA);
  \draw[\colorCircle,line width=\locpc,dashed] (circ1) circle (\locpe);
  \transv {postr} node[below=1pt,xshift=3.4pt,color=black] {$\scriptscriptstyle 0$};
  \filldraw[color=\colorDefect] (0,\locpy) circle (\locpd) 
          (circ1)++(138:\locpA) circle (\locpD) (circ2)++(0,\locpA) circle (\locpD);
  \node at (-1.6-\locpa,0) {$\varphi(\Sigma_1)~=$};
  \node[\colorCircle] at (121:0.32*\locpa) {$\DD$};
  \node[\colorCircle] at (161:0.23*\locpa) {$\SS_2$};
  \node[\colorCircle] at (327:0.44*\locpa) {$\SS_1$};
  \node at (1.6+\locpa,0) {~};   
  \begin{scope}[shift={(2*\locpa-2.1,-1.7*\locpy)}]
  \coordinate (circ1) at (270-\locpB:\locpb*\locpa);
  \coordinate (Circ1) at (337-\locpB:1.07*\locpb*\locpa);
  \coordinate (circ2) at (270+\locpB:1.42*\locpb*\locpa);
  \coordinate (Circ2) at (265+\locpB:1.48*\locpb*\locpa);
  \coordinate (postr) at (123:0.22*\locpa);
  \scopeArrow{0.01}{>}
  \filldraw[\colorCircle,fill=\colorSigma,line width=\locpc,postaction={decorate}]
           (0,0) circle (\locpa cm and \locpy cm);
  \end{scope}
  \scopeArrow{0.29}{\arrowDefect}
  \draw[line width=\locpc,color=\colorTransp] 
           plot[smooth,tension=0.65] coordinates{
           (Circ2) (35:0.53*\locpa) (Circ1) (postr) (85:0.56*\locpy) (90:0.83*\locpy) (0,\locpy) };
  \end{scope}
  \scopeArrow{0.70}{\arrowDefect}
  \draw[line width=\locpc,color=\colorTransp,postaction={decorate}]
           (circ1)++(71:\locpA) -- (postr);
  \draw[line width=\locpc,color=\colorTransp,postaction={decorate}]
           (88:0.74*\locpy) -- ++ (106:0.05);
  \draw[line width=\locpc,color=\colorTransp,postaction={decorate}]
           (272:0.37*\locpy) -- ++ (235:0.05);
  \end{scope}
  \scopeArrow{0.8}{<}
  \filldraw[\colorCircle,fill=white,line width=\locpC,postaction={decorate}] (circ1) circle (\locpA);
  \filldraw[\colorCircle,fill=white,line width=\locpC,postaction={decorate}] (circ2) circle (\locpA);
  \end{scope}
  \draw[\colorCircle,line width=\locpc] (circ1) circle (\locpA) (circ2) circle (\locpA);
  \draw[\colorCircle,line width=\locpc,dashed] (circ1) ++(0.13*\locpa,0) circle (\locpe);
  \transv {postr} node[below=1.5pt,xshift=-4.8pt,color=black] {$\scriptscriptstyle 0$};
  \filldraw[color=\colorDefect] (0,\locpy) circle (\locpd) 
          (circ1)++(71:\locpA) circle (\locpD) (circ2)++(0,\locpA) circle (\locpD);
  \node at (-1.4-\locpa,0) {{\Large $\simeq$}};
  \end{scope}
  \end{tikzpicture}
  \hspace*{-2.8em}
  \nonumber\\[-1.4em]~
  \label{eq:diffeo-refinchange}
  \\[-1.3em]~\nonumber
  \end{eqnarray}
We factor the computation of $\varphi(\Sigma_{1})$ through the disk $\DD$ to obtain the 
isomorphism (compare the proof of Lemma \ref{lemma:repMCG})
  \be
  \bearll
  \ZZ(\varphi(\Sigma_{1}))(x \boti y  \boti  \cc{z}) \!\!\!&
  = \ZZ(\Sigma_{1})(x \boti y \boti \cc{z})
  \Nxl6 &
  = \Hom_{\calz(\cala)}(z,x \oti y ) \cong \Hom_{\calz(\cala)}(z, {}_{x}F(y)) \,.
  \eear
  \ee
Next we perform the refinement change from the first to the second picture in 
\eqref{eq:diffeo-refinchange} and apply Lemma \ref{Lemma:MCG-br} so as to obtain the isomorphism
${}_{x}F(y) \,{\cong}\, F_{x}(y)$, which is given by the braiding of $x$ and $y$. From 
the pictures it is evident that the resulting refinement is isotopic to $\Sigma_{2}$. Thus the
functor we have computed is indeed $\ZZ(\varphi)$. 
\end{proof}

\begin{Remark} \label{rem:TV-2}
Recall the comparison with the standard Turaev--Viro construction in Remark \ref{rem:TV}.
It remains to compare the representations of the mapping class groups. To this end we 
observe that an oriented modular functor is completely determined by a Lego-Teichm\"uller
game (compare e.g.\ \cite{BAki,fuSc22}). As a consequence it suffices to compare the modular
functors on $n$-holed spheres, which has been discussed in Remark \ref{rem:TV}, and on
5 elementary moves between specific surfaces. A detailed comparison of all these moves
in the present and standard Turaev--Viro approaches is beyond the scope of this paper. But
note that the braiding move is covered by Proposition \ref{prop:braid}. We leave a detailed 
analysis of all moves, in particular of the S-move, which involves surfaces of genus one,
to future work.
\end{Remark}

Instead, let us now compute, besides the braiding move, also the value of the modular functor
on a Dehn twist. In this case we start with a non-fine \defectsurface\ $(\Sigma, \chi)$ given
by a two-punctured sphere, as shown in the following picture:
    \def\locpa  {0.59}   
    \def\locpA  {2.6}    
    \def\locpb  {0.71}   
    \def\locpB  {1.35}   
    \def\locpc  {2.2pt}  
    \def\locpC  {1.3pt}  
    \def\locpd  {0.12}   
    \def\locpe  {1.1}    
    \def\locpf  {0.77}   
  \be
  \raisebox{-6.8em}
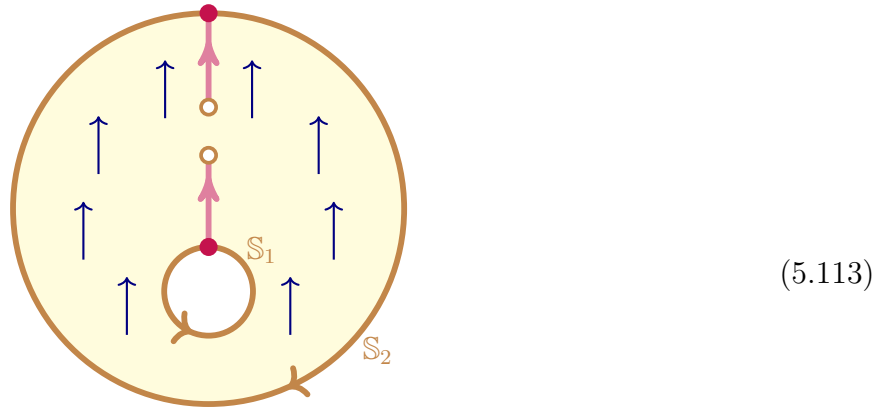
{\begin{tikzpicture}[scale=1.1]
  \scopeArrow{0.83}{<}
  \filldraw[\colorCircle,fill=\colorSigma,line width=\locpc,postaction={decorate}]
           (0,0) circle (\locpA);
  \end{scope}
  \scopeArrow{0.69}{>}
  \filldraw[\colorCircle,fill=white,line width=\locpc,postaction={decorate}]
           (0,-\locpe) circle (\locpa);
  \end{scope}
  \scopeArrow{0.71}{\arrowDefect}
  \draw[line width=\locpc,\colorTransp,postaction={decorate}] (0,\locpB) -- (0,\locpA);
  \end{scope}
  \scopeArrow{0.77}{\arrowDefect}
  \draw[line width=\locpc,\colorTransp,postaction={decorate}] (0,\locpa-\locpe) -- (0,\locpb);
  \end{scope}
  \transv {0,\locpb};
  \transv {0,\locpB};
  \foreach \xx in {-50,-10,...,230} \draw[thick,color=\colorFrame,->]
       (\xx:0.65*\locpA)++(0,-0.5*\locpf) -- +(0,\locpf);
  \fill[\colorDefect] (0,\locpa-\locpe) circle (\locpd) (0,\locpA) circle (\locpd);
  \node[\colorCircle] at (1.2*\locpa,-0.5*\locpe) {$\SS_1$};
  \node[\colorCircle] at (-40:1.13*\locpA) {$\SS_2$};
  \end{tikzpicture}}
  \ee
Here the framing is the one that corresponds to a straight cylinder. The gluing categories are 
$\ZZ(\SS_{1}) \eq \calz(\cala)$ and $\ZZ(\SS_{2}) \eq \cc{\calz(\cala)}$. Consider the 
diffeomorphism 
$\varphi \colon (\Sigma, \chi) \rightarrow (\Sigma, \varphi_{*}(\chi))$ indicated in
    \def\locpa  {0.47}   
    \def\locpA  {2.6}    
    \def\locpb  {0.43}   
    \def\locpB  {1.43}   
    \def\locpc  {2.2pt}  
    \def\locpC  {1.3pt}  
    \def\locpd  {0.12}   
    \def\locpe  {1.1}    
    \def\locpf  {0.77}   
  \be
  \raisebox{-5.7em}
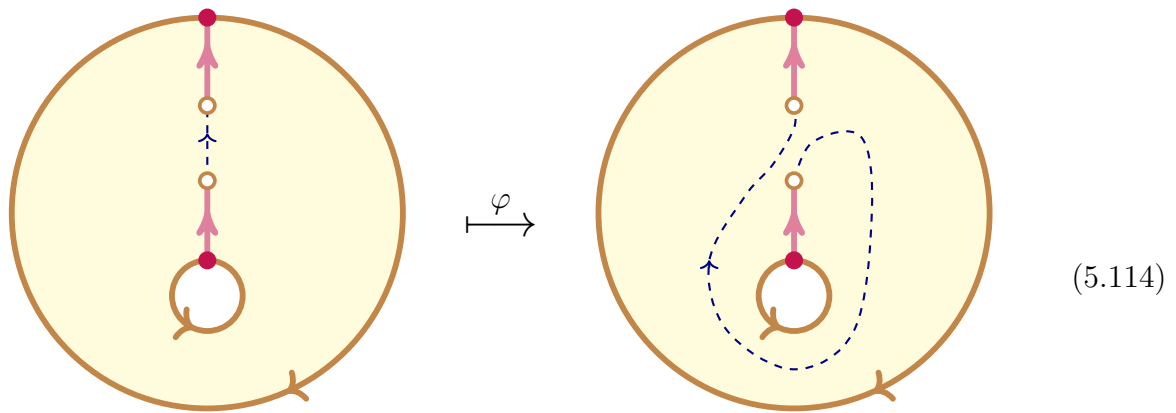
{\begin{tikzpicture}
  \scopeArrow{0.83}{<}
  \filldraw[\colorCircle,fill=\colorSigma,line width=\locpc,postaction={decorate}]
           (0,0) circle (\locpA);
  \end{scope}
  \scopeArrow{0.69}{>}
  \filldraw[\colorCircle,fill=white,line width=\locpc,postaction={decorate}]
           (0,-\locpe) circle (\locpa);
  \end{scope}
  \scopeArrow{0.71}{\arrowDefect}
  \draw[line width=\locpc,\colorTransp,postaction={decorate}] (0,\locpB) -- (0,\locpA);
  \end{scope}
  \scopeArrow{0.63}{\arrowDefect}
  \draw[line width=\locpc,\colorTransp,postaction={decorate}] (0,\locpa-\locpe) -- (0,\locpb);
  \end{scope}
  \scopeArrow{0.63}{>}
  \draw[thick,\colorFrame,dashed,postaction={decorate}] (0,\locpb) -- (0,\locpB);
  \end{scope}
  \transv {0,\locpb};
  \transv {0,\locpB};
  \fill[\colorDefect] (0,\locpa-\locpe) circle (\locpd) (0,\locpA) circle (\locpd);
  \node at (1.5*\locpA,0) {{\Large$\stackrel\varphi\longmapsto$}};
  \begin{scope}[shift={(3*\locpA,0)}]
  \scopeArrow{0.83}{<}
  \filldraw[\colorCircle,fill=\colorSigma,line width=\locpc,postaction={decorate}]
           (0,0) circle (\locpA);
  \end{scope}
  \scopeArrow{0.69}{>}
  \filldraw[\colorCircle,fill=white,line width=\locpc,postaction={decorate}]
           (0,-\locpe) circle (\locpa);
  \end{scope}
  \scopeArrow{0.71}{\arrowDefect}
  \draw[line width=\locpc,\colorTransp,postaction={decorate}] (0,\locpB) -- (0,\locpA);
  \end{scope}
  \scopeArrow{0.63}{\arrowDefect}
  \draw[line width=\locpc,\colorTransp,postaction={decorate}] (0,\locpa-\locpe) -- (0,\locpb);
  \end{scope}
  \draw[thick,\colorFrame,dashed] plot[smooth,tension=0.61] coordinates{
           (0,\locpb) (0.25,2.2*\locpb) (1.9*\locpa,2.35*\locpb) (2.2*\locpa,0)
           (1.7*\locpa,-0.6*\locpA) (0,-0.8*\locpA) (-1.9*\locpa,-0.6*\locpA)
           (-2.4*\locpa,-0.25*\locpA) (-1.4*\locpa,0.5*\locpb) (-0.6*\locpa,1.5*\locpb)
           (0,\locpB-0.3) (0,\locpB) };
  \draw[thick,\colorFrame,->] (-2.4*\locpa,-0.28*\locpA) -- +(0,0.1); 
  \transv {0,\locpb};
  \transv {0,\locpB};
  \fill[\colorDefect] (0,\locpa-\locpe) circle (\locpd) (0,\locpA) circle (\locpd);
  \end{scope}
  \end{tikzpicture}~~
  \label{eq:Fig9-14.5.}}
  \ee
The framings $\chi$ and $\varphi_{*}(\chi)$ are non-homotopic, as can e.g.\ be seen 
by realizing that the the winding number along the dashed line on the right hand side of 
\eqref{eq:Fig9-14.5.} is non-zero.

To compute the value of $\ZZ$ on $\varphi$, we make use of the universal property of the
parallelization that were developed in Proposition \ref{prop:can-Pi}, which allows us to work
entirely on the level of pre-block spaces. Consider the refinements $(\Sigma;\Sigma_{i})$, 
with $i \iN \{1,2,3\}$, of $\varphi_{*}(\Sigma)$ shown in the following picture:
    \def\locpa  {0.47}   
    \def\locpA  {2.6}    
    \def\locpc  {2.2pt}  
    \def\locpC  {1.3pt}  
    \def\locpd  {0.12}   
    \def\locpe  {0.5}    
    \def\locpE  {0.99}   
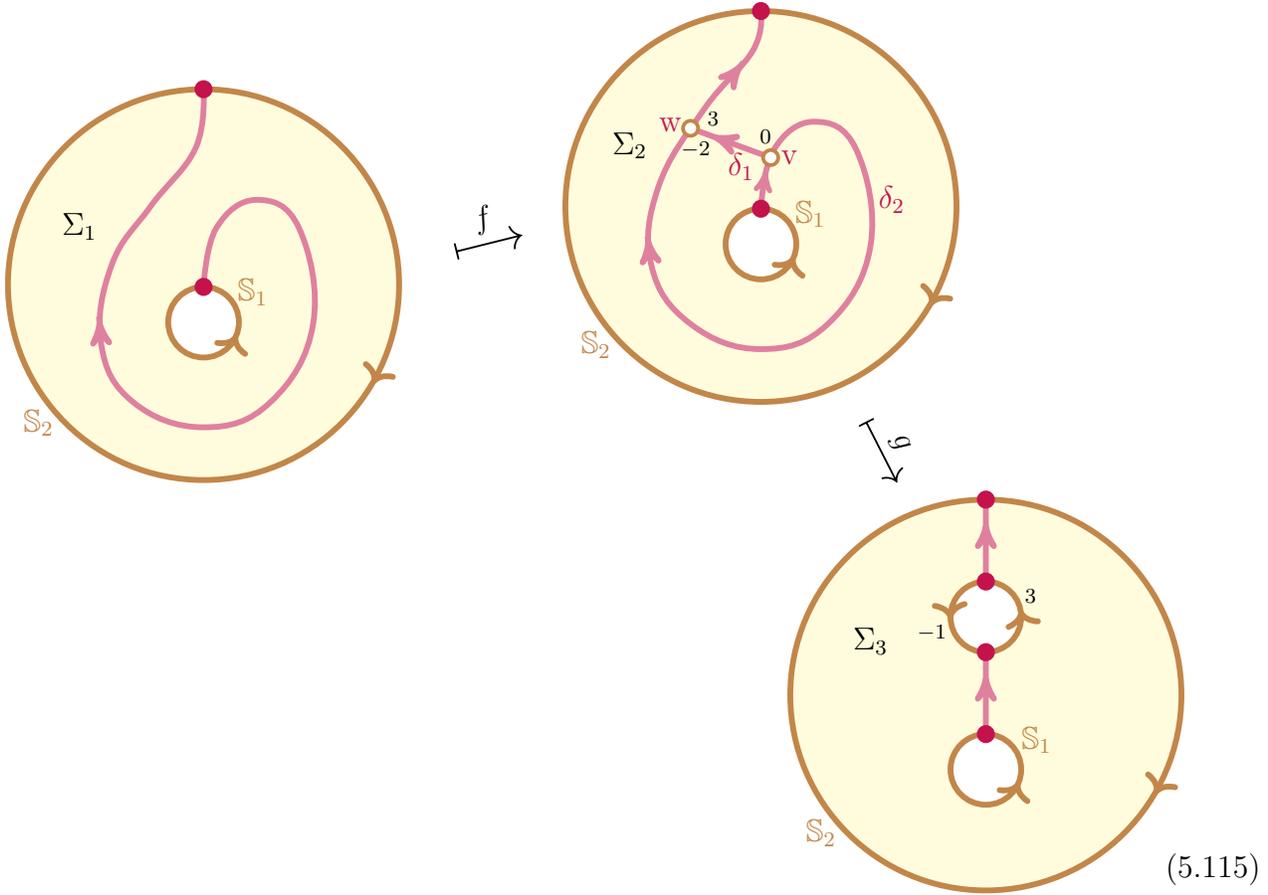
  \begin{eqnarray}
  \begin{tikzpicture}
  \scopeArrow{0.93}{<}
  \filldraw[\colorCircle,fill=\colorSigma,line width=\locpc,postaction={decorate}]
           (0,0) circle (\locpA);
  \end{scope}
  \scopeArrow{0.94}{>}
  \filldraw[\colorCircle,fill=white,line width=\locpc,postaction={decorate}]
           (0,-\locpe) circle (\locpa);
  \end{scope}
  \draw[line width=\locpc,\colorTransp] plot[smooth,tension=0.71] coordinates{
           (0,\locpa-\locpe) (0.36*\locpa,1.6*\locpa) (1.5*\locpa,2.4*\locpa)
	   (2.7*\locpa,1.6*\locpa) (3.1*\locpa,-0.2*\locpA)
           (1.9*\locpa,-0.6*\locpA) (0,-0.73*\locpA) (-2*\locpa,-0.6*\locpA)
           (-2.9*\locpa,-0.3*\locpA) (-2.6*\locpa,0.1*\locpA) (-1.5*\locpa,0.4*\locpA)
           (-0.04*\locpA,0.7*\locpA) (0,\locpA) };
  \scopeArrow{0.5}{\arrowDefect}
  \draw[line width=\locpc,\colorTransp,postaction={decorate}]
           (-2.94*\locpa,-1.05*\locpe) -- +(-0.008,0.1); 
  \end{scope}
  \fill[\colorDefect] (0,\locpa-\locpe) circle (\locpd) (0,\locpA) circle (\locpd);
  \node[\colorCircle] at (1.4*\locpa,-0.1) {$\SS_1$};
  \node[\colorCircle] at (220:1.1*\locpA) {$\SS_2$};
  \node at (155:0.70*\locpA) {$\Sigma_1$};
  \begin{scope}[shift={(2.85*\locpA,0.4*\locpA)}]
  \coordinate (trnv1) at (0.05*\locpA,0.25*\locpA);
  \coordinate (trnv2) at (-0.36*\locpA,0.4*\locpA);
  \scopeArrow{0.93}{<}
  \filldraw[\colorCircle,fill=\colorSigma,line width=\locpc,postaction={decorate}]
           (0,0) circle (\locpA);
  \end{scope}
  \scopeArrow{0.94}{>}
  \filldraw[\colorCircle,fill=white,line width=\locpc,postaction={decorate}]
           (0,-\locpe) circle (\locpa);
  \end{scope}
  \draw[line width=\locpc,\colorTransp] plot[smooth,tension=0.71] coordinates{
           (0,\locpa-\locpe) (0.36*\locpa,1.6*\locpa) (1.5*\locpa,2.4*\locpa)
	   (2.7*\locpa,1.6*\locpa) (3.1*\locpa,-0.2*\locpA)
           (1.9*\locpa,-0.6*\locpA) (0,-0.73*\locpA) (-2*\locpa,-0.6*\locpA)
           (-3.1*\locpa,-0.3*\locpA) (-2.9*\locpa,0.1*\locpA) (-1.6*\locpa,0.5*\locpA)
           (-0.04*\locpA,0.8*\locpA) (0,\locpA) };
  \scopeArrow{0.68}{\arrowDefect}
  \draw[line width=\locpc,\colorTransp,postaction={decorate}]
           (-0.12*\locpA,0.7*\locpA) -- +(0.08,0.1); 
  \draw[line width=0.8*\locpc,\colorTransp,postaction={decorate}]
           (0.12*\locpa,1.93*\locpa-\locpe) -- +(0.016,0.1); 
  \draw[line width=\locpc,\colorTransp,postaction={decorate}]
           (-3.18*\locpa,-\locpe) -- +(-0.012,0.1); 
  \end{scope}
  \scopeArrow{0.72}{\arrowDefect}
  \draw[line width=\locpc,\colorTransp,postaction={decorate}]
           (trnv1) -- node[midway,xshift=4pt,yshift=-9pt,color=\colorDefect] {$\delta_1$} (trnv2);
  \end{scope}
  \transv {trnv1};
  \transv {trnv2};
  \node[xshift=-2pt,yshift=8pt] at (trnv1) {$\scriptstyle 0$};
  \node[xshift=7.2pt,color=\colorDefect] at (trnv1) {$\mathrm v$};
  \node[xshift=8.5pt,yshift=3.5pt] at (trnv2) {$\scriptstyle 3$};
  \node[xshift=2pt,yshift=-8pt] at (trnv2) {$\scriptstyle {-}2$};
  \node[xshift=-7.5pt,yshift=1.5pt,color=\colorDefect] at (trnv2) {$\mathrm w$};
  \fill[\colorDefect] (0,\locpa-\locpe) circle (\locpd) (0,\locpA) circle (\locpd);
  \node[\colorCircle] at (1.4*\locpa,-0.1) {$\SS_1$};
  \node[\colorCircle] at (220:1.1*\locpA) {$\SS_2$};
  \node[\colorDefect] at (3:0.67*\locpA) {$\delta_2$};
  \node at (155:0.74*\locpA) {$\Sigma_2$};
  \end{scope}
  \begin{scope}[shift={(4.0*\locpA,-2.1*\locpA)}]
  \scopeArrow{0.93}{<}
  \filldraw[\colorCircle,fill=\colorSigma,line width=\locpc,postaction={decorate}]
           (0,0) circle (\locpA);
  \end{scope}
  \scopeArrow{0.94}{>}
  \filldraw[\colorCircle,fill=white,line width=\locpc,postaction={decorate}]
           (0,-\locpE) circle (\locpa);
  \end{scope}
  \scopeArrow{0.03}{>} \scopeArrow{0.51}{>}
  \filldraw[\colorCircle,fill=white,line width=\locpc,postaction={decorate}]
           (0,0.5*\locpA-0.5*\locpE+0.5*\locpa) circle (\locpa);
  \end{scope} \end{scope}
  \scopeArrow{0.71}{\arrowDefect}
  \draw[line width=\locpc,\colorTransp,postaction={decorate}]
           (0,\locpa-\locpE) -- (0,0.5*\locpA-0.5*\locpE-0.5*\locpa);
  \draw[line width=\locpc,\colorTransp,postaction={decorate}]
	   (0,0.5*\locpA-0.5*\locpE+1.5*\locpa) -- (0,\locpA);
  \end{scope}
  \node[xshift=3.5pt,yshift=8pt] at (\locpa,0.5*\locpA-0.5*\locpE+0.5*\locpa) {$\scriptstyle 3$};
  \node[xshift=-7pt,yshift=-6pt] at (-\locpa,0.5*\locpA-0.5*\locpE+0.5*\locpa) {$\scriptstyle {-}1$};
  \fill[\colorDefect] (0,\locpa-\locpE) circle (\locpd) (0,\locpA) circle (\locpd)
           (0,0.5*\locpA-0.5*\locpE-0.5*\locpa) circle (\locpd)
           (0,0.5*\locpA-0.5*\locpE+1.5*\locpa) circle (\locpd);
  \node[\colorCircle] at (1.44*\locpa,-0.23*\locpA) {$\SS_1$};
  \node[\colorCircle] at (220:1.1*\locpA) {$\SS_2$};
  \node at (155:0.65*\locpA) {$\Sigma_3$};
  \end{scope}
  \node[rotate=14] at (1.44*\locpA,0.29*\locpA) {{\Large$\stackrel f\longmapsto$}};
  \node[rotate=-63] at (3.52*\locpA,-.83*\locpA) {{\Large$\stackrel g\longmapsto$}};
  \end{tikzpicture}
  \hspace*{-2.2em} \nonumber\\[-1.1em]~
  \label{Sigma i=123}
  ~\\[-1.2em] ~\nonumber
  \end{eqnarray}
The change $(\Sigma;\Sigma_1) \,{\xmapsto{~f\,}}\, (\Sigma;\Sigma_2)$ 
of refinements consists of adding the defect line $\delta_{1}$, while the 
refinement change $(\Sigma;\Sigma_2) \,{\xmapsto{~g\,}}\, (\Sigma;\Sigma_3)$ 
is the deletion of $\delta_{2}$. For obtaining a representative for the value of the
block functor on $\varphi_{*}(\Sigma)$ we pick the refinement $(\Sigma;\Sigma_{3})$. With
the help of Lemma \ref{lemma:repMCG} we then get
  \be
  \ZZ(\varphi) = \Pi_{\Sigma; \varphi(\Sigma),\Sigma_{3}} :\quad
  \Zprime(\Sigma) \xrightarrow{~~} \Zprime(\varphi_{*}(\Sigma); \Sigma_{3}) \,.
  \label{eq:expl-dehn}
  \ee
We first consider the pre-block spaces. The proof of the following statements follows directly
by combining the canonical isomorphisms from the Yoneda lemma with the explicit expressions
for the \mute\ objects of the two relevant gluing circles. As in \eqref{Sigma i=123}
we denote the latter circles by $\mathrm v$ and $\mathrm w$, and choose the conventions
    \def\locpa  {0.53}   
    \def\locpA  {1.4}    
    \def\locpc  {2.2pt}  
    \def\locpd  {0.11}   
    \def\locpx  {270}    
    \def\locpX  {90}     
    \def\locpy  {40}     
    \def\locpY  {-40}    
    \def\locpz  {140}    
    \def\locpZ  {-140}   
  \be
  \raisebox{-3.3em}{\begin{tikzpicture}
  \scopeArrow{0.28}{<} \scopeArrow{0.58}{<} \scopeArrow{0.95}{<}
  \filldraw[line width=\locpc,\colorCircle,fill=white,postaction={decorate}]
           node[left=1.5cm,color=black,yshift=-5pt] {$\mathrm v ~=$} (0,0) circle (\locpa)
           node[above=12pt,xshift=6pt,color=black]{$\scriptstyle 0$};
  \end{scope} \end{scope} \end{scope}
  \scopeArrow{0.77}{\arrowDefectB}
  \draw[line width=\locpc,color=\colorTransp,postaction={decorate}]
           (\locpx:\locpa) -- (\locpx:\locpA) node[below=-2pt,color=\colorDefect] {$\cc b$};
  \end{scope}
  \scopeArrow{0.75}{\arrowDefect}
  \draw[line width=\locpc,color=\colorTransp,postaction={decorate}]
           (\locpy:\locpa) -- (\locpy:\locpA)
           node[right=-3pt,yshift=1pt,color=\colorDefect] {$c^\vee{\otimes}b$};
  \draw[line width=\locpc,color=\colorTransp,postaction={decorate}]
           (\locpz:\locpa) -- (\locpz:\locpA) node[left=-2pt,yshift=1pt,color=\colorDefect] {$c$};
  \end{scope}
  \filldraw[color=\colorDefect] (\locpx:\locpa) circle (\locpd)
           (\locpy:\locpa) circle (\locpd) (\locpz:\locpa) circle (\locpd);
  \begin{scope}[shift={(2*\locpA+3.2,-0.73)}]
  \scopeArrow{0.06}{<} \scopeArrow{0.44}{<} \scopeArrow{0.75}{<}
  \filldraw[line width=\locpc,\colorCircle,fill=white,postaction={decorate}]
           node[left=1.5cm,color=black] {$\mathrm w ~=$} (0,0) circle (\locpa)
           node[below=10.5pt,xshift=6.2pt,color=black]{$\scriptstyle {-}2$}
           node[left=12.5pt,yshift=1pt,color=black]{$\scriptstyle 3$};
  \end{scope} \end{scope} \end{scope}
  \scopeArrow{0.81}{\arrowDefect}
  \draw[line width=\locpc,color=\colorTransp,postaction={decorate}]
           (\locpX:\locpa) -- (\locpX:\locpA)
           node[above=-2.5pt,xshift=3pt,color=\colorDefect] {$d{\otimes}a^{\vee\vee}$};
  \end{scope}
  \scopeArrow{0.65}{\arrowDefectB}
  \draw[line width=\locpc,color=\colorTransp,postaction={decorate}]
           (\locpY:\locpa) -- (\locpY:\locpA) node[right=-3pt,yshift=1pt,color=\colorDefect] {$\cc a$};
  \draw[line width=\locpc,color=\colorTransp,postaction={decorate}]
           (\locpZ:\locpa) -- (\locpZ:\locpA) node[left=-2.2pt,yshift=1pt,color=\colorDefect] {$\cc d$};
  \end{scope}
  \filldraw[color=\colorDefect] (\locpX:\locpa) circle (\locpd)
           (\locpY:\locpa) circle (\locpd) (\locpZ:\locpa) circle (\locpd);
  \end{scope}
  \end{tikzpicture}} \qquad
  \ee
for the labels of these \defectonemanifold s; hereby their \mute\ objects take the form
  \be
  \om_{\mathrm v}
  = \int^{b\in\cala}\!\!\!\! \int_{\!c\in\cala} c \boti c^{\vee} {\otimes} b \boti \cc{b}
  \qquad\text{and}\qquad
  \om_{\mathrm w}
  = \int^{a\in\cala}\!\!\!\! \int^{m\in\cala} m {\otimes} a^{\vee\vee} \boti \cc a \boti \cc m \,.
  \ee

\begin{Lemma} \label{Lemma:prebl-sigmai}
The $($relative$)$
pre-block spaces for the \defectsurface s $\Sigma_{i}$ shown in \eqref{Sigma i=123} are
  \be
  \bearll
  \Zpre(\Sigma_{1})(\cc{y} \boti x) \!\!&\dsty
  = \int^{a \in \cala}\! \Hom(y,a) \otik \Hom(a,x) \,\cong\, \Hom(y,x) \,,
  \Nxl5
  \Zpreprime(\Sigma_{2})(\cc{y} \boti x) \!\!&\dsty
  = \int^{b\in\cala}\!\!\!\! \int_{\!c\in\cala} \int^{a\in\cala}\!\!\!\! \int^{d\in\cala}\!
  \Hom(b,x) \otik \Hom(d, c^{\vee} \oti b)
  \nxl1 &\dsty \hspace*{11em}
  \otik \Hom(a,c) \otik \Hom(y, d \oti a^{\vee\vee})
  \Nxl1 &\dsty 
  \cong \int^{d\in\cala}\!\!\!\! \int_{\!c\in\cala}\! \Hom(y,d \oti c^{\vee \vee}) \otik \Hom(c^{\vee \vee} \oti d,x)
  \Nxl4 &\dsty 
  \cong \int^{d\in\cala}\!\!\!\! \int_{\!c\in\cala}\! \Hom(y,d \oti c) \otik \Hom(c \oti d,x) \,,
  \Nxl5
  \Zpreprime(\Sigma_{3})(\cc{y} \boti x) \!\!&\dsty
  = \int^{d\in\cala}\! \Hom(y, d^{\vee\vee}) \otik \Hom(d,x) \,\cong\, \Hom(y,x^{\vee\vee}) 
  \label{eq:ZpreSigma2}
  \eear
  \ee
for $x \iN \ZZ(\SS_1)$	and $y \iN \ZZ(\SS_2)$
(here all Hom spaces are morphism spaces in $\cala$).
\end{Lemma}

With the results on the block functors given in Corollary \ref{cor:gen-cylinder} we then obtain

\begin{Lemma} \label{Lemma:bl-sigma13}
The block functors for the \defectsurface s
$\Sigma_{1}$ and $\Sigma_{3}$ are, up to canonical isomorphism, given by 
  \be
  \ZZ(\Sigma_{1})(\cc{y} \boti x) = \Hom_{\calz(\cala)}(y,x) \qquad\text{and}\qquad
  \ZZ(\Sigma_{3})(\cc{y} \boti x) = \Hom_{\calz(\cala)}(y,x^{\vee \vee}) \,,~
  \label{eq:blcokSigmai}
  \ee
respectively.
\end{Lemma}

We will not need an explicit expression for the block functor on $\Sigma_{2}$.
 
The following notion is a direct generalization of the one of a ribbon twist in a 
ribbon category:

\begin{Definition} \label{def:twist}
The \emph{twist} on an object $x$ in a braided monoidal category $\calc$ with right duality 
is the isomorphism $x \,{\tocong}\, x^{\vee\vee}$ given by
  \be
  x \xrightarrow{~\coevr_{\!x^{\vee}}\otimes\id_{x}\,} x^{\vee} \oti x^{\vee\vee} \oti x   
  \xrightarrow{~\id_{x^{\vee}} \otimes c_{x^{\vee\vee},x}\,} x^{\vee} \oti x \oti x^{\vee\vee}
  \xrightarrow{~\evr_{x} \otimes \id_{x^{\vee \vee}}\,}\, x^{\vee \vee} ,
  \label{eq:twist}
  \ee
with $c$ the braiding in $\calc$. 
\end{Definition}

In case $\calc$ is actually ribbon, it is in particular pivotal, and the ribbon twist is the
composition of the twist \eqref{eq:twist} and the pivotal structure \Cite{Eq.\,(2.2.27)}{BAki}.
In the situation at hand, we deal with the case that $\calc \eq \calz(\cala)$,
with $\cala$ not assumed to be pivotal.
An ordinary ribbon twist is not monoidal, but rather obeys a compatibility
condition with the double braiding. Since this condition comes from a relation
in the mapping class group of a three-holed sphere, the twist in the sense of
Definition \ref{def:twist} obeys a similar compatibility condition, which can
also be shown algebraically.

We proceed to compute $\ZZ(\varphi)$. 

\begin{Proposition} \label{Proposition:dehntwist}
The value of $\ZZ$ on the diffeomorphism $\varphi$ from \eqref{eq:Fig9-14.5.}
is the natural isomorphism 
  \be
  \ZZ(\varphi)(\cc{y} \boti x) :\quad
  \Hom_{\calz(\cala)}(y,x)\rightarrow  \Hom_{\calz(\cala)}(y,x^{\vee \vee})
  \ee
that is given by post-composition with the twist \eqref{eq:twist} on $x$. 
\end{Proposition}

\begin{proof}
We are going to show that the morphism between the pre-block spaces for $\Sigma_{1}$ and 
$\Sigma_{3}$ that is induced by the refinement changes in \eqref{Sigma i=123} is given by
post-composition with the twist on $x$. By the universal property of 
$\Pi_{\Sigma; \varphi(\Sigma),\Sigma_{3}} \eq \ZZ(\varphi)$ obtained in Proposition
\ref{prop:can-Pi} we can then conclude that $\ZZ(\varphi)$ is the post-composition with 
the twist as well. 
 \\[.2em]
We first compute the morphism $f^{\pre} \colon \Zpre(\Sigma_1) \,{\to}\, \Zpreprime(\Sigma_2)$ 
that corresponds to the first refinement change in \eqref{Sigma i=123}. With the help of the
parallel transport operations on $\Sigma_2$ we see that $f^{\pre}$ is the following composite:
  \be
  \hspace*{-1.1em}\bearll
  \multicolumn2l {\dsty
  \int^{a\in\cala}\! \Hom(y,a) \otik \Hom(a,x) \,\xrightarrow{~~~}
  \int_{\!c\in\cala}\!\int^{a\in\cala}\! \Hom(y,{}^{\vee\vee\!}c\oti a) \otik \Hom(c\oti a,x) }
  \Nxl2 \hspace*{4.5em} &\dsty
  \tocong \int_{\!c\in\cala}\! \int^{a\in\cala}\! \Hom({}^{\vee\!}c\oti y,a) \otik \Hom(c\oti a,x)
  \Nxl4 &\dsty
  \tocong \int^{a\in\cala}\!\!\!\! \int_{\!c\in\cala}\! \Hom(y\oti {}^{\vee\!}c,a)
  \otik \Hom(c\oti a,x)
  \Nxl5 &\dsty
  \tocong \int^{d\in\cala}\!\!\!\! \int_{\!c\in\cala}\! \Hom(y,d\oti c) \otik \Hom(c\oti d,x)
  \,\tocong\, \Zpreprime(\Sigma_{2})(\cc{y} \boti x) \,.
  \eear
  \ee
Here the first morphism is the comonad structure of $\int^{m\in\calm} \cc{m} \boti m \iN
\CC{\calm} \,{\Boti 1}\, \calm \,{\Boti 1}$, the second and fourth are adjunction isomorphisms,
and the third is obtained from the braiding of $y \iN \calz(\cala)$; the last isomorphism 
holds by Lemma \ref{Lemma:prebl-sigmai}. 
 \\[.2em]
Next we compute the morphism $g^{\pre} \colon \Zpreprime(\Sigma_2) \,{\to}\, \Zpreprime(\Sigma_3)$ 
that corresponds to the second refinement change in \eqref{Sigma i=123}, given by deletion
of the defect line $\delta_{2}$. To this end we would like to make use of the dinatural
component at $\one \iN \cala$ of the end. Before we can do so we must, however,
first use the isomorphism
  \be
  \int^{b\in\cala}\!\!\!\! \int_{\!c\in\cala} c \boti c^{\vee} {\otimes} b \boti \cc{b}
  \,\cong \int^{b\in\cala}\!\!\!\! \int_{\!c\in\cala} b{\otimes}c \boti c^{\vee} \boti \cc{b} 
  \label{eq:mute-iso}
  \ee
for the \mute\ object at the gluing circle $v$ in \eqref{Sigma i=123}, 
which arises from the canonical isomorphism $\int_{\!c\in\cala} c^{\vee}{\otimes}b \boti c
\,{\cong} \int_{\!c\in\cala} c^{\vee} \boti b {\otimes} c$.
When using the first expression for $\Zpreprime(\Sigma_{2})$ in \eqref{eq:ZpreSigma2}
we obtain the first isomorphism in the following composite, which is $g^{\pre}$:
  \be
  \hspace*{-0.8em}
  \begin{array}{l} \dsty
  \int^{b\in\cala}\!\!\! \int_{\!c\in\cala} \int^{a\in\cala}\!\!\! \int^{d\in\cala}\! 
  \Hom(b,x) \otik \Hom(d,c^{\vee}{\otimes}\,b) \otik \Hom(a,c) \otik \Hom(y,d \oti a^{\vee\vee})
  \Nxl4 \dsty \qquad
  \tocong \int^{b}\!\!\! \int_{\!c} \int^{a}\!\!\! \int^{d}
  \Hom(b,x) \otik \Hom(d, c^{\vee}) \otik \Hom(a,b \oti c) \otik \Hom(y,d \oti a^{\vee\vee})~~
  \Nxl2 
  \multicolumn1r {\dsty
  \tocong \int_{\!c\in\cala}\Hom(y,c^{\vee} \oti (x\oti c)^{\vee\vee})
  \,\xrightarrow{~~~}\, \Hom(y,x^{\vee\vee}) \,.}
  \eear
  \ee
Here the last morphism is the dinatural component of the end at $\one$. 
 \\[.2em]
Using repeatedly the Yoneda isomorphism, one can express the morphism $g^{\pre}{\circ}f^{\pre}$
as the composite along the upper-right path from $\Hom(y,x)$ to $\Hom(y \oti x^{\vee}, \one)$
in the diagram
  \be
  \begin{tikzcd}[column sep=1.1cm]
  \Hom(y,x) \ar{r} \ar{ddd} & \int_{c} \Hom(y, {}^{\vee\!}c \oti c \oti x) \ar{r}{\cong}
  & \int_{c} \Hom(c \oti y , c \oti x) \ar{d}{\cong} 
  \\
  ~ &  ~ &  \int_{c} \Hom(y \oti c, c \oti x) \ar{d}{\cong} 
  \\
  ~ & ~ & \int_{c} \Hom(y \oti c \oti x^{\vee},c) \ar{d}{}
  \\
  \Hom({x}^{\vee} \oti y, \one) \ar{rr} & ~   & \Hom(y \oti x^{\vee}, \one)
  \end{tikzcd} \qquad
  \ee
Here the first morphism in the upper row is given by the canonical morphism
$\one \,{\to} \int_{c\in\cala}\! {}^{\vee\!}c \oti c$ that is obtained 
from the coevaluation of $\cala$. The subsequent morphisms are obtained from  the adjunction
and the braiding of $y$, followed by the morphism induced from 
$\int_{c} \cc{c} \boti c \oti x \,{\cong} \int_{c} \cc{c \oti x^{\vee}} \boti c$, while the
final vertical morphism in this path is again the dinatural component at $\one$. The left-lower
path from $\Hom(y,x)$ to $\Hom(y \oti x^{\vee}, \one)$ is provided directly by the 
adjunction and the braiding of $y$. It follows from the dinaturality of the involved morphisms 
that the diagram commutes. Further, the naturality of the braiding implies that 
composing the left-lower path with the adjunction isomorphism
$\Hom(y \oti x^{\vee}, \one) \,{\cong}\, \Hom(y,x^{\vee\vee})$ is the same as 
post-com\-po\-si\-tion with the twist on $x$. Thus we have shown that $g^{\pre} \,{\circ} f^{\pre}$
is given by post-composition with the twist on $x$. Now post-composing with the twist is even
an isomorphism $\Hom_{\calz(\cala)}(y,x)\,{\tocong}\, \Hom_{\calz(\cala)}(y,x^{\vee\vee})$,
and thus we have a commuting diagram  
  \be
  \begin{tikzcd}
  \Hom_{\calz(\cala)}(y,x) \ar{r} \ar{d}[swap]{U}
  & \Hom_{\calz(\cala)}(y,x^{\vee \vee}) \ar{d}{U}
  \\
  \Hom_{\cala}(y,x) \ar{r}
  & \Hom_{\cala}(y,x^{\vee\vee})
  \end{tikzcd}
  \ee
where the horizontal arrows are post-composition with the twist on $x$ and the
vertical arrows forget the half-braidings. By the universal property of $\ZZ(\varphi)$, as obtained
in Proposition \ref{prop:can-Pi}, we conclude that $\ZZ(\varphi)$ 
is the post-composition with the twist. 
 \end{proof}


\appendix

\section{Framing shifts} \label{app:shift}

In this appendix we examine the operation of shifting the framing of a \defectsurface\ 
along a defect line. As we will see, performing this operation simultaneously for all
defect lines yields canonical autoequivalences of the bicategory of \undefectsurface s. 

Recall the bicategory $\Bfrdef$ of \undefectsurface\ defined in Section \ref{sec:Bfrdef}.
A  homotopy $\chi_{t}$ between two framings $\chi_{1}$ and $\chi_{2}$ on an \undefectsurface\ 
is by definition required to be in particular a framing for each $t$, so the corresponding 
vector field has to be parallel to the defect lines of $\Sigma$ for all $t$. By Remark
\ref{Remark:structure-def-bord}(iii) such homotopies yield isomorphic \defectsurface s,
$(\Sigma, \chi_{1}) \,{\cong}\, (\Sigma, \chi_{2})$. If one would relax 
the condition, so that $\chi_{t}$ is only required to be a non-vanishing 
vector field for all $t$, 
but not necessarily parallel to the defect lines for $t\iN (0,1)$, many more framings would 
be related, but the resulting 1-morphisms in $\Bfrdef$ would, in general, be non-isomorphic. 

There is a local prototype for this kind of more general homotopy, which allows one to relate
different framings on the same \undefectsurface\ $\Sigma$: Consider a defect line $\delta$
on a \defectsurface\ $\Sigma$ with a local neigbourhood that looks like the left hand side 
of the following picture (using the framing index to specify the vector field locally):  
 \\[-2.1em]~
    \def\locpA  {2.5}    
    \def\locpb  {1.9}    
    \def\locpc  {2.2pt}  
    \def\locpd  {0.12}   
  \be
  \raisebox{-7.5em}
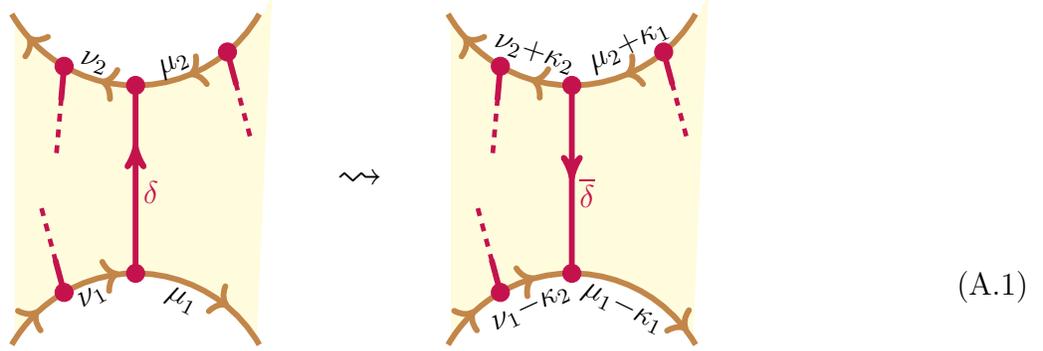
{\begin{tikzpicture}
  [decoration={random steps,segment length=2mm}]
  \filldraw[white,fill=\colorSigma,line width=\locpc,postaction={decorate}]
           (0,\locpA+\locpb)+(210:\locpb) arc (210:270:\locpb) -- (0,0) arc (90:150:\locpb)
           -- cycle
           (0,\locpA+\locpb)+(350:\locpb) arc (350:270:\locpb) -- (0,0) arc (90:30:\locpb)
           -- cycle;
  \scopeArrow{0.69}{\arrowDefect}
  \draw[line width=\locpc,color=\colorDefect,postaction={decorate}] (0,0) --
	   node[right=-2pt,yshift=-5pt] {$\delta$} (0,\locpA);
  \end{scope}
  \scopeArrow{0.22}{<} \scopeArrow{0.61}{<} \scopeArrow{0.91}{<}
  \draw[\colorCircle,line width=\locpc,postaction={decorate}]
	   (30:\locpb)+(0,-\locpb) arc (30:150:\locpb);
  \end{scope} \end{scope} \end{scope}
  \scopeArrow{0.32}{>} \scopeArrow{0.62}{>} \scopeArrow{0.91}{>}
  \draw[\colorCircle,line width=\locpc,postaction={decorate}]
	   (0,\locpA+\locpb)+(330:\locpb) arc (330:210:\locpb);
  \end{scope} \end{scope} \end{scope}
  \draw[line width=\locpc,color=\colorDefect]
				(0,-\locpb)++(120:\locpb) -- +(105:0.18*\locpA)
				(0,\locpA+\locpb)++(240:\locpb) -- +(265:0.18*\locpA)
				(0,\locpA+\locpb)++(310:\locpb) -- +(285:0.18*\locpA);
  \draw[line width=0.8*\locpc,color=\colorDefect,dashed]
				(0,-\locpb)++(120:\locpb) -- +(105:0.47*\locpA)
				(0,\locpA+\locpb)++(240:\locpb) -- +(265:0.47*\locpA)
				(0,\locpA+\locpb)++(310:\locpb) -- +(285:0.47*\locpA);
  \filldraw[color=\colorDefect] (0,0) circle (\locpd) (0,\locpA) circle (\locpd)
                                (0,-\locpb)+(120:\locpb) circle (\locpd)
                                (0,\locpA+\locpb)+(240:\locpb) circle (\locpd)
                                (0,\locpA+\locpb)+(310:\locpb) circle (\locpd);
  \node[rotate=-28] at (0.32*\locpb,-0.37) {$\mu_1$};
  \node[rotate=29]  at (-0.30*\locpb,-0.32) {$\nu_1$};
  \begin{scope}[shift={(0,\locpA)}]
  \node[rotate=26]  at (0.27*\locpb,0.28) {$\mu_2$};
  \node[rotate=-18] at (-0.29*\locpb,0.30) {$\nu_2$};
  \end{scope}
  \begin{scope}[shift={(\locpb+1.5,0.5*\locpA)}]
  \node at (-0.42,0) {{\Large$\rightsquigarrow$}};
  \end{scope}
  \begin{scope}[shift={(2*\locpb+2.0,0)}];
  \filldraw[white,fill=\colorSigma,line width=\locpc,postaction={decorate}]
           (0,\locpA+\locpb)+(210:\locpb) arc (210:270:\locpb) -- (0,0) arc (90:150:\locpb)
           -- cycle
           (0,\locpA+\locpb)+(350:\locpb) arc (350:270:\locpb) -- (0,0) arc (90:30:\locpb)
           -- cycle;
  \scopeArrow{0.63}{\arrowDefectB}
  \draw[line width=\locpc,color=\colorDefect,postaction={decorate}] (0,0) --
           node[right=-2pt,yshift=-5pt] {$\overline\delta$} (0,\locpA);
  \end{scope}
  \scopeArrow{0.11}{<} \scopeArrow{0.68}{<} \scopeArrow{0.91}{<}
  \draw[\colorCircle,line width=\locpc,postaction={decorate}]
           (30:\locpb)+(0,-\locpb) arc (30:150:\locpb);
  \end{scope} \end{scope} \end{scope}
  \scopeArrow{0.32}{>} \scopeArrow{0.64}{>} \scopeArrow{0.91}{>}
  \draw[\colorCircle,line width=\locpc,postaction={decorate}]
           (0,\locpA+\locpb)+(330:\locpb) arc (330:210:\locpb);
  \end{scope} \end{scope} \end{scope}
  \draw[line width=\locpc,color=\colorDefect]
				(0,-\locpb)++(120:\locpb) -- +(105:0.18*\locpA)
				(0,\locpA+\locpb)++(240:\locpb) -- +(265:0.18*\locpA)
				(0,\locpA+\locpb)++(310:\locpb) -- +(285:0.18*\locpA);
  \draw[line width=0.8*\locpc,color=\colorDefect,dashed]
				(0,-\locpb)++(120:\locpb) -- +(105:0.47*\locpA)
				(0,\locpA+\locpb)++(240:\locpb) -- +(265:0.47*\locpA)
				(0,\locpA+\locpb)++(310:\locpb) -- +(285:0.47*\locpA);
  \filldraw[color=\colorDefect] (0,0) circle (\locpd) (0,\locpA) circle (\locpd)
                                (0,-\locpb)+(120:\locpb) circle (\locpd)
                                (0,\locpA+\locpb)+(240:\locpb) circle (\locpd)
                                (0,\locpA+\locpb)+(310:\locpb) circle (\locpd);
  \node[rotate=-26] at (0.36*\locpb,-0.45) {$\mu_1{-}\kappa_1$};
  \node[rotate=25]  at (-0.29*\locpb,-0.44) {$\nu_1{-}\kappa_2$};
  \begin{scope}[shift={(0,\locpA)}]
  \node[rotate=26]  at (0.40*\locpb,0.53) {$\mu_2{+}\kappa_1$};
  \node[rotate=-18] at (-0.26*\locpb,0.44) {$\nu_2{+}\kappa_2$};
  \end{scope}
  \end{scope}
  \end{tikzpicture}
  \label{eq:delta-deltabar}}
  \ee
  ~\\[-2.2em]
For any $(\kappa_{1},\kappa_{2}) \iN (2 \mathbb{Z}{+}1)^{2}$ the surface shown on the
right hand side of \eqref{eq:delta-deltabar} has a unique framing that differs only
locally from the framing on the left hand side. An instance of this operation in the 
analogous situation of a free boundary occurs in Proposition \ref{prop:lineflip}. 
We call such operations on an \undefectsurface\ $\Sigma$ the \emph{odd framing shifts} 
at $\delta_{i}$ on $\Sigma$ and denote them by $S_{\kappa_{1},\kappa_{2}}(\delta_{i})(\Sigma)$.
Analogously there are \emph{even framing shifts} $S_{\kappa_{1},\kappa_{2}}(\delta_{i})(\Sigma)$
at $\delta_{i}$ for  $(\kappa_{1},\kappa_{2}) \iN (2 \mathbb{Z})^{2}$, for which the orientation
of $\delta_{i}$ does not change.
 
The so defined framing shifts amount to the following modification of the framings of
\undefectonemanifold s $\LL$: 
For $p$ a defect point on $\LL$ and $(\kappa_{1},\kappa_{2}) \iN (2 \mathbb{Z}{+}1)^{2}$,
the one-manifold $S_{\kappa_{1},\kappa_{2}}(p)(\LL)$ is the \undefectonemanifold\ for which
the sign of $p$ is flipped and, if $p$ is positive, the framing index on the segment preceding
$p$ is decreased by $\kappa_{2}$ while the one after $p$ is decreased by $\kappa_{1}$; if
instead $p$ is negative, the framing index on the segment preceding $p$ is increased by
$\kappa_1$ and the one after $p$ is increased by $\kappa_2$.
If $(\kappa_{1},\kappa_{2}) \iN (2\mathbb{Z})^{2}$, the orientation of $p$ is kept and the 
change on the framing indices is, with the same notation as in Figure \eqref{eq:delta-deltabar},
  \be
  \nu_1 \mapsto \nu_1 - \kappa_2 \,,\quad~ \mu_1 \mapsto \mu_1 - \kappa_1 \,,\quad~
  \nu_2 \mapsto \nu_2 + \kappa_2 \quad~ \text{and} \quad~ \mu_2 \mapsto \mu_2 + \kappa_1 \,.
  \ee

For $(\kappa_{1},\kappa_{2}) \iN (2 \mathbb{Z})^{2}$ consider, for a given
an \undefectsurface\ $\Sigma$, the \undefectsurface\
  \be
  S_{\kappa_{1},\kappa_{2}}(\Sigma):=  \prod_{i}S_{\kappa_{1},\kappa_{2}}(\delta_{i})(\Sigma) \,,
  \label{eq:prod-all}
  \ee
where the product ranges over all defect lines of $\Sigma$. For an \undefectonemanifold\ $\LL$ 
we define analogously 
$S_{\kappa_{1},\kappa_{2}}(\LL) \,{:=}\, \prod_{i}S_{\kappa_{1},\kappa_{2}}(p_{i})(\LL)$,
with product over all defect points of $\LL$.
We tacitly extend these prescriptions for defect lines and defect points
to apply also to free boundary segments and their end points, in which case only one instead
of two framing indices are shifted, without modifying the notation.
The following statements then follow directly from the definitions.

\begin{Lemma} \label{Lemma:framing-shift-func}
Let $(\kappa_{1},\kappa_{2}) \iN (2 \mathbb{Z})^{2}$. The assignments $\LL\,{\mapsto}\,
S_{\kappa_{1},\kappa_{2}}(\LL)$ and $\Sigma\,{\mapsto}\,S_{\kappa_{1},\kappa_{2}}(\Sigma)$ for
\undefectonemanifold s $\LL$ and \undefectsurface s $\Sigma$ define a symmetric monoidal functor
$S_{\kappa_{1},\kappa_{2}}\colon \Bfrdef \,{\to}\, \Bfrdef$. This functor is an autoequivalence.
\end{Lemma}

We call the symmetric monoidal functor $S_{\kappa_1,\kappa_2}\colon \Bfrdef \,{\to}\, \Bfrdef$
the ($\kappa_1,\kappa_2)$\emph{-framing shift functor}. For a \defectsurface\ $\Sigma$ with
labels we define $S_{\kappa_{1},\kappa_{2}}(\Sigma)$ by applying $S_{\kappa_{1},\kappa_{2}}$
to the underlying \undefectsurface\ and keeping the labels, and analogously for
labeled \defectonemanifold s. Thereby we obtain, by composing with
framing shifts, a whole family of modular functors from a given one:

\begin{Proposition}
Let $\ZZ\colon \Bfrdefdec \,{\xrightarrow{~}}\, \mathcal S$ be a modular functor. Then for 
every pair $(\kappa_{1},\kappa_{2}) \iN (2 \mathbb{Z})^{2}$, the functor 
  \be
  {}^{\kappa_{1}}\ZZ^{\kappa_{2}} := \ZZ \circ S_{\kappa_{1},\kappa_{2}}
  \ee
is a modular functor as well. 
\end{Proposition}

We call the so obtained functor ${}^{\kappa_{1}}\ZZ^{\kappa_{2}}$ the 
($\kappa_1,\kappa_2)$\emph{-shift of $\ZZ$}. This functor can also be described in terms of 
$\ZZ$ by shifting the labels instead of the framings: For each 
$(\kappa_{1},\kappa_{2}) \iN (2 \mathbb{Z})^{2}$ and bimodule $\calm_{i}$ we define the bimodule
  \be
  \widetilde{S}_{\kappa_{1},\kappa_{2}}(\calm_{i})
  := {}_{}^{-\kappa_{1}\!\!}\calm{}^{-\kappa_{2}}_{} .
  \ee
For a \defectsurface\ $\Sigma$ the \defectsurface\ $\widetilde{S}_{\kappa_1,\kappa_2}(\Sigma)$
is now defined as the \emph{same} \undefectsurface\ $\Sigma$, but with
$\widetilde{S}_{\kappa_{1},\kappa_{2}}$ applied to the labels of each defect line of $\Sigma$. 
It follows that $\widetilde{S}_{\kappa_{1},\kappa_{2}}(\Sigma)$ is again a \defectsurface,
and extending the operation $\widetilde{S}_{\kappa_{1},\kappa_{2}}$ in the obvious way to
\defectonemanifold s yields a symmetric monoidal functor 
  \be
  \widetilde{S}_{\kappa_{1},\kappa_{2}}:\quad \Bfrdefdec \xrightarrow {~~} \Bfrdefdec . 
  \ee
It is straightforward to see that the analogue of  Proposition \ref{prop:lineflip} for defect 
lines holds. Thus we conclude that  

\begin{Proposition}
Let $\,\ZZ \colon \Bfrdefdec \,{\xrightarrow{~~}}\, \mathcal S$ be a modular functor. For
any pair $(\kappa_{1},\kappa_{2}) \,{\in}$
	                                 $(2 \mathbb{Z})^{2}$ the modular functors 
${}^{\kappa_{1}}\ZZ^{\kappa_{2}} \,{=}\, \ZZ \cir S_{\kappa_{1},\kappa_{2}}$ and 
$\ZZ \cir \widetilde{S}_{\kappa_{1},\kappa_{2}} $ are canonically isomorphic. 
\end{Proposition}


\section{Categorical constructions}

In this appendix we briefly mention pertinent concepts and constructions on categories, 
such as (co)monads, (co)ends and the Eilenberg--Watts calculus on finite linear categories.

\subsection{General concepts}\label{app:1}

\paragraph*{Finite tensor categories.}

Throughout this paper, all categories appearing as labels for defect manifolds 
are assumed to be enriched over \findim\ \ko-vector spaces, with \ko\ a fixed 
algebraically closed field.  A (\ko-linear) \emph{finite category} is an abelian category 
for which the set of isomorphism classes of simple objects is finite and
for which every object has finite length and a projective cover.
A (\ko-linear) \emph{finite tensor category} is a finite \ko-linear monoidal
category with simple monoidal unit and with a left and a right duality.
For $\cala$ and $\calb$ finite categories we denote by $\Funle(\cala,\calb)$ and
$\Funre(\cala,\calb)$ the finite categories of (\ko-linear) left exact and right exact
functors from $\cala$ to $\calb$, respectively.

\paragraph*{Module categories.}

A right \emph{module category} over a monoidal category $\cala$, or right 
$\cala$-\emph{mo\-du\-le}, for short, is a category $\calm$ together with a bilinear bifunctor
$\calm \,{\times}\, \cala \,{\to}\, \cala$ (which we denote by ``.'') that  is exact in each
argument, and with a natural family of isomorphisms $(m.a)\,.\,b \,{\to}\, m\,.\,(a{\otimes}b)$ 
for $a,b \,{\in}\,\cala$ and $m\,{\in}\,\calm$ obeying obvious coherence axioms. 
Left $\cala$-modules are defined analogously. For $\cala$ and $\calb$ monoidal categories,
an $\cala$-$\calb$-\emph{bimodule} is a category that is a left $\cala$-module and a right
$\calb$-module together with natural coherence isomorphisms connecting the left and right actions.

A (right) $\cala$-\emph{module functor} $(F,\phi)$ between right $\cala$-modules $\calm$ and
$\caln$ is a linear functor $F\colon \calm \,{\to}\, \caln$ together with a natural family $\phi$
of isomorphisms $F(m.a) \,{\to}\, F(m).a$, for $a \,{\in}\,\cala$ and $m\,{\in}\,\calm$,
such that the obvious pentagon and triangle diagrams commute.
A \emph{module natural transformation} $\eta$ between right $\cala$-module functors $(F,\phi^F)$
and $(G,\phi^G)$ is a natural transformation $\eta\colon F \,{\Rightarrow}\, G$ such that
$\phi^G_{m,a}\,{\circ}\,\eta_{m.a}^{} \,{=}\, (\eta^{}_m \,.\,\id^{}_a) \,{\circ}\, \phi^F_{m,a}$
for all $a \,{\in}\,\cala$ and $m\,{\in}\,\calm$,

\paragraph*{Ends and coends.}

A \emph{dinatural transformation} from a functor 
$F\colon\, \calc \,{\times}\, \calc\opp \,{\to}\, \cald$ to an object $d\,{\in}\,\cald$
is a family of morphisms $\varphi_c\colon F(c,c) \,{\to}\, d$ for $c \,{\in}\, \calc$ such that
  \be
  \varphi_ {c'} \,{\circ}\, F(g,c') = \varphi_c \,{\circ}\, F(c,g)
  \ee
for all $g\,{\in}\,\Hom_\calc(c,c')$.
A \emph{coend} $(C,\iota)$ for $F\colon\, \calc \,{\times}\, \calc\opp \,{\to}\, \cald$ 
is an object $C \,{\in}\, \cald$ together with a dinatural transformation 
$\iota\colon F \,{\to}\, C$ that is universal among all dinatural transformations $\varphi$
from $F$ to an object of $\cald$, i.e.\ for any such dinatural transformation
there is a unique morphism 
$\varpi$ obeying $\varphi_c \,{=}\, \varpi \,{\circ}\, \iota_c$ for all $c \,{\in}\, \calc$.
The coend $(C,\iota)$, as well as the underlying object $C$, is denoted by 
$\int^{c\in\calc}\! F(c,c)$.
Dually, an \emph{end} $(E,\jmath) \,{=}\, \int_{c\in\calc} F(c,c)$ for $F$ is a 
universal dinatural transformation from a constant to $F$.
If a coend or end exists, then it is unique up to unique isomorphism.
When the categories in question are functor categories, one often wants 
specific properties of the functors to be preserved under taking ends and coends.
In the case of interest to us, the relevant property is representability.
Accordingly, we are considering left exact functors and impose the
universal property for ends and coends within categories of left exact functors. (Unlike e.g.\
in \cite{lyub11}, we do not use a separate symbol $\oint$ for such `left exact (co)ends'.)
For more information on ends and coends see e.g.\ \cite{fuSc23} and the literature cited there.

The Hom functor preserves and reverses ends in the sense that
  \be
  \bearl
  \Hom_\cald(-,\int_{c\in\calc}F(c,c)) \,\cong \dsty\int_{c\in\calc}\! \Hom_\cald(-,F(c,c))
  \qquad{\rm and}
  \Nxl4
  \Hom_\cald(\int^{c\in\calc}\! F(c,c),-) \,\cong \dsty\int_{c\in\calc}\! \Hom_\cald(F(c,c),-)
  \eear
  \ee
for any functor $F\colon \, \calc \,{\times}\, \calc\opp \,{\to}\, \cald$.
Further, it follows from the (co-)Yoneda lemma
(see e.g.\ \Cite{Prop.\,2.7}{fScS2}) that for any linear functor $F$ between
finite linear categories $\cala$ and $\cala'$ there are natural isomorphisms
  \be
  \int^{a \in \cala}\! \Hom_\cala(a,-) \otimes F(a) \,\cong\, F \,\cong\,
  \int_{a \in \cala} \Hom_\cala(-,a)^*_{} \otimes F(a)
  \label{eq:convo}
  \ee
of linear functors.

By considering the particular functor $F \eq \Hom_\calc(c,G(-))$, it follows for instance 
that for any pair $H \colon \calb \,{\to}\, \cala$ and $G \colon \cala \,{\to}\, \calc$ 
of composable linear functors between finite linear categories and for any $c \iN \calc$
there is a canonical isomorphism
  \be
  \Hom_\calc(c, G \cir H(-)) \,\cong \int^{a\in\cala}\! \Hom_\calc(c, G(a))
  \otimes_\ko \Hom_\cala(a, H(-)) \,.
  \label{eq:EW-xxx}
  \ee
            
When dealing with Deligne products of finite linear categories one has in addition
  \be
  \bearl
  \Hom_{\CC\calc\boxtimes\cald} \big(-\,, \int^{c\in\calc}\cc c \boti F(c) \big) \,\cong 
  \dsty \int^{c\in\calc}\!\! \Hom_{\CC\calc\boxtimes\cald} \big(-\,, \cc c \boti F(c) \big)
  \eear
  \label{eq:fScS2Prop3.4}
  \ee
for any \emph{left exact} functor $F \,{\in}\, \Funle(\calc,\cald)$, 
and a similar reversed isomorphism for right exact functors \Cite{Prop.\,3.4}{fScS2}.
These isomorphisms
can e.g.\ be used, in conjunction with the Eilenberg--Watts equivalences \eqref{Phile.Psile},
to show that for any two finite categories $\calm$ and $\caln$ the mapping
  \be
  F \longmapsto \Hom_\caln(-,F(-))
  \ee
defines an equivalence
  \be
  \Funle(\calm,\caln) \,\xrightarrow{~~\simeq~~}\, \Funle(\calm\boti\CC\caln,\vect) \,.
  \label{eq:MN->MccNvect}
  \ee
An inverse equivalence
$\Funle(\calm\boti\CC\caln,\vect) \,{\xrightarrow{~\simeq~}}\, \Funle(\calm,\caln)$
is given by
  \be
  G \,\longmapsto \int^{n\in\caln}\!\! G(-\boti\cc n) \otimes n \,.
  \ee

\paragraph*{Monads and comonads.}

A \emph{monad} $M \,{=}\, (M,\mu,\eta)$ on a category $\calc$ is an algebra (or monoid) in the
monoidal category of endofunctors of $\calc$, that is, an endofunctor $M$ together with natural
transformations $\mu \,{=}\, (\mu_c)_{c\in\calc}\colon M \,{\circ}\, M \,{\Rightarrow}\, M$
(product) and $\eta \,{=}\, (\eta_c)_{c\in\calc}\colon \id_\calc \,{\Rightarrow}\, M$ (unit) which
satisfy associativity and unit properties, i.e.\
$\mu_c \circ M(\mu_c) \,{=}\, \mu_c \,{\circ}\, \mu_{M(c)}$ and
$\mu_c \circ M(\eta_c) \,{=}\, \id_{M(c)} \,{=}\, \mu_c \,{\circ}\, \eta_{M(c)}$.
Analogously, a \emph{comonad} on $\calc$ is a coalgebra in the category of endofunctors
of $\calc$, i.e.\ an endofunctor $W$ together with a coproduct
$W \,{\Rightarrow}\, W \,{\circ}\, W$ and counit $W \,{\Rightarrow}\, \id_\calc$
satisfying coassociativity and counit properties.
Every adjoint pair of functors $F$ and $G$ ($G$ right adjoint to $F$) gives rise to a monad
structure on the endofunctor $G\,{\circ}\,F$ and a comonad structure on $F\,{\circ}\,G$. 

A (left) \emph{comodule} over a comonad $W$ on $\calc$ (also called a $W$-coalgebra) is 
an object $x \,{\in}\, \calc$ together with a morphism $\delta\colon x \,{\to}\, W(x)$
satisfying analogous compatibility conditions with the coproduct $\Delta$ and counit $\eps$
of $W$ as a comodule over a comonoid, i.e.\
$\Delta_x \,{\circ}\, \delta \,{=}\, W(\delta)\,{\circ}\, \delta$ and
$\eps_x \,{\circ}\, \delta \,{=}\, \id_x$. Right comodules and left and right modules over
a monad are defined analogously.

If $W$ is a comonad on a category $\calc$, then the opposite functor
  \be
  \cc{W}:\quad \cc{\calc} \rightarrow \cc{\calc}
  \label{eq:oppositefunctor}
  \ee
is a monad on $\CC\calc$. Further, if the adjoint functors $W^{\la}$ and $W^{\ra}$ exist, 
then they are monads on $\calc$, while
the functors $\cc{W^{\la}}, \cc{W^{\ra}}\colon \cc{\calc} \,{\to}\, \cc{\calc}$ are comonads
on $\CC\calc$.
Moreover, if $x \,{\in}\, \calc$ is a comodule over $W$, then it is also naturally a module 
over the monads $W^{\la}$ and $W^{\ra}$, while the object $\cc{x} \,{\in}\, \CC{\calc}$ 
is a module over the monad $\cc W$ and a comodule over the comonads $\cc{W^{\la}}$ and
$\cc{W^{\ra}}$. Analogous statements hold for monads on $\calc$.


\subsection{The canonical $\kappa$-twisted (co)monads and their (co)modules}

Here we define various versions of $\kappa$-balanced functors and exhibit the framed center 
as a universal category for $\kappa$-balanced functors, 
making use of corresponding $\kappa$-twisted (co)monads.
 
We consider a cyclically composable string
$(\calm_1^{\epsilon_1}, \calm_2^{\epsilon_2},...\,, \calm_n^{\epsilon_n})$ of bimodules 
with corresponding balancings $\kappa_i$ between $\calm_i$ and $\calm_{i+1}$. 
There are four different situations to be distinguished, indexed by the sequence 
$\{ (\epsilon_{i}, \epsilon_{i+1}) \}$. The following definition captures all cases:

\begin{Definition}
  \label{definition:bal-functor}
Let $(\calm_1^{-},\calm_2^{-},\calm_3^{+},\calm_4^+)$ be a cyclically composable string of
bimodules with  corresponding sequence $\vec{\kappa} \,{=}\, (\kappa_1,\kappa_2,\kappa_3,
\kappa_4)$ of balancings $\kappa_i$ between $\calm_i$ and $\calm_{i+1}$. The category 
$\Funlebal_{\vec{\kappa}}((\calm_1^{-},\calm_2^{-},\calm_3^{+},\calm_4^{+}),\calx)$
of $\kappa$-balanced functors to a finite category $\calx$ consists of functors 
$F \colon \calm_{1}^- \boti \cdots \boti \calm_4^+ \,{\to}\, \calx$ with coherent isomorphisms
  \be
  \begin{array}{l}
  F(\cc{a.m_1} \boti \cc{m_2} \boti m_3  \boti m_4) \,\xrightarrow{~\cong\,}\,
  F(\cc{m_1} \boti \cc{m_2. {}^{[\kappa_1+4]}a} \boti m_3  \boti m_4) \,,
  \Nxl8
  F(\cc{m_1} \boti \cc{b.m_2} \boti m_3  \boti m_4) \,\xrightarrow{~\cong\,}\,
  F(\cc{m_1} \boti \cc{m_2} \boti {}^{[\kappa_{2}+2]}b.m_3  \boti m_4) \,,
  \Nxl8
  F(\cc{m_1} \boti \cc{m_2} \boti m_3.c  \boti m_4) \,\xrightarrow{~\cong\,}\,
  F(\cc{m_1} \boti \cc{m_2} \boti m_3  \boti {}^{[\kappa_{3}]}c.m_4) \qquad\text{and}
  \Nxl8
  F(\cc{m_1} \boti \cc{m_2} \boti m_3  \boti m_4.d) \,\xrightarrow{~\cong\,}\,
  F(\cc{m_1.{}^{[\kappa_4+2]}}d \boti \cc{m_2} \boti m_3  \boti m_4) \,.
  \eear
  \ee
\end{Definition}

In accordance with the various situations, for each value of $i$ there are four (related) comonads
$Z_{[\kappa_{i}]}(\calm_{i}^{\epsilon_{i}}\boti \calm_{i+1}^{\epsilon_{i+1}})$ on the category
$\calm_i^{\epsilon_{i}}\boti \calm_{i+1}^{\epsilon_{i+1}}$, defined by
  \be
  \begin{array}{l}\dsty
  Z_{[\kappa_{1}]}(\calm_{1}^{-} \boti \calm_{2}^{-})(\cc{m_1} \boti \cc{m_2})
  := \int_{a \in \cala} \cc{a.m_{1}} \boti \cc{m_{2}.{}^{[\kappa_{1}+3]}a}
  \Nxl7\dsty
  Z_{[\kappa_{2}]}(\calm_{2}^{-} \boti \calm_{3})(\cc{m_2} \boti {m_3})
  := \int_{b \in \calb} \cc{b.m_{2}} \boti {}^{[\kappa_{2}+1]}b.m_{3}
  \Nxl7\dsty
  Z_{[\kappa_{3}]}(\calm_{3} \boti \calm_{4})(m_{3} \boti {m_4})
  := \int_{c \in \calc} \cc{m_{3}.c} \boti {}^{[\kappa_{3}-1]}c.m_{4} \qquad\text{and}
  \Nxl7\dsty
  Z_{[\kappa_{4}]}(\calm_{4} \boti \calm_{1}^{-})(m_{4} \boti \cc{m_1})
  := \int_{d \in \cald} \cc{m_{4}.d} \boti \cc{ m_{1}.{}^{[\kappa_{4}+1]}d} \,,
  \eear
  \label{eq:4comonads}
  \ee
respectively. We can combine these comonads to a comonad
on $\vec{\calm} \,{=}\, (\calm_1^{-}, \calm_2^{-}, \calm_3^{+}, \calm_4^{+})$. We denote
this comonad by $Z_{\vec\kappa} \,{\eq}\, Z_{\vec\kappa}(\vec{\calm})$ and refer to it as
the \emph{canonical $\vec\kappa$-twisted comonad} on $\vec{\calm}$. 
This terminology is justified by the following result:

\begin{Proposition} \label{proposition:framed-center-univ}
Let $\vec\calm \,{=}\, (\calm_1^{\epsilon_1},...\,,\calm_n^{\epsilon_n})$ be a string
of cyclically composable bimodules with balancings $\vec\kappa$. Denote by $\vec\calm
\,{\boxtimes} \,{:=}\, \calm_1^{\epsilon_1}\,{\boxtimes}\cdots{\boxtimes}\,\calm_n^{\epsilon_n}$
the Deligne product of the bimodules and by
$\vec\calm\,{\Boti{\vec\kappa}}\,{:=}\, \cenf{\calm_1^{\epsilon_1} \,\Boti{\kappa^{}_1}\,
\calm_2^{\epsilon_2} \,\Boti{\kappa^{}_2}\, \calm_3^{\epsilon_3} \,\Boti{\kappa_3^{}}\, \cdots
\!\!{\Boti{\kappa_{n-1}}} \! \calm_n^{\epsilon_n}\, {\Boti{\kappa^{}_n}}\, }$ the 
corresponding framed center, as in Definition {\rm\ref{def:framedcenter}}. 
\Itemizeiii

\item[{\rm (i)}]
For any object $x \iN \vec\calm\,{\boxtimes}\,$ the object $Z_{[\vec{\kappa}]}(x)$ has
a canonical structure of an object in the category $\vec\calm{\Boti{\vec\kappa}}$, to be 
denoted by $I_{[\vec{\kappa}]}(x)$. 
 \\
This naturally defines a functor $I_{[\vec{\kappa}]}
\iN \Funlebal_{\vec\kappa}(\vec\calm\,{\boxtimes},\vec\calm{\Boti{\vec\kappa}})$ such that
$Z_{[\vec{\kappa}]} \,{=}\, U \,{\circ}\, I_{[\vec{\kappa}]} $, with 
$U\colon \vec\calm{\Boti{\vec\kappa}} \To \vec\calm$ the functor that forgets the balancing. 

\item[{\rm (ii)}]
The framed center $\vec\calm\,{\Boti{\vec\kappa}}$ together with the $\vec\kappa$-balanced functor 
$I_{[\vec\kappa]}$ is universal for $\vec\kappa$-balanced functors: For any finite category 
$\calx$, pre-composition with $I_{[\vec\kappa]}$ is an equivalence
  \be
  \Funle(\vec\calm{\Boti{\vec\kappa}}\,, \calx) \xrightarrow{~\simeq~}\,
  \Funlebal_{\vec\kappa}(\vec\calm\,{\boxtimes}\,, \calx) \,.
  \ee
\end{itemize}
\end{Proposition}
 
\begin{proof}
Consider part (i) in the special case of a single bimodule $\calm$. For $\kappa \,{\in}\, 2\Z$
the comonad $Z_{[\kappa]}$ on $\calm$ is given by the endofunctor
  \be
  Z_{[\kappa]}:\quad  m \,\xmapsto{~~} \int_{a\in\cala} a\,.\,m\,.\,a^{[\kappa-1]} .
  \label{eq:def:comonadZn2}
  \ee
$Z_{[\kappa]}$ can be viewed as acting with the object $\int_{a\in\cala} a\boti \cc a
\,{\in}\, \cala\boti\CC\cala \,\,{\cong}\, \Psire(\id_\Cala)$ on $\calm$
(with $\Psire$ the Eilenberg--Watts equivalence \eqref{Phire...Psire}),
after applying a suitable power of the double dual functor to it.
The balancings of $I_{[\kappa]}(m) \iN \vec\calm\,{\Boti{\kappa}}$ with underlying object 
$Z_{[\kappa]}(m)$ are determined by invoking the coherent isomorphisms
  \be
  \int_{a\in\cala}\! b \oti a \boxtimes \cc a
  \,\cong\! \int_{a\in\cala}\! a \boxtimes \cc{{}^{\vee_{}} b \oti a} \qquad \text{and} \qquad
  \int_{a\in\cala}\! a \oti b \boxtimes \cc a 
  \,\cong\! \int_{a\in\cala}\! a \boxtimes \cc{a \oti b^{\vee}}~
  \ee
for $b \,{\in}\, \cala$; these are obtained by setting $\calm \eq \cala$ in the isomorphisms
\eqref{eq:mmbalancing}.
It follows that this way we have defined a functor $I_{[\kappa]}\colon \calm \To 
\vec\calm\,{\Boti{\kappa}}$ that satisfies $U \,{\circ}\, I_{[\kappa]} \,{=}\, Z_{[\kappa]}$.
The proof of the general case follows by the same reasoning. 
 \\[.2em]
(ii)\, We first show that the framed center is equivalent to the category of comodules over 
$Z_{\vec{\kappa}}$: This follows in the case of a single bimodule $\calm$ 
with the help of the linear isomorphisms
  \be
  \Hom_\calm(m.a,{}^{[\kappa-2]}a.m) \,\cong\, \Hom_\calm(m, {}^{[\kappa-2]}a.m.a^{\vee})
  \ee
after taking the end. The general case is treated analogously. The universal property of the
framed center now follows, analogously as the universal property of the center of 
a bimodule category \cite{genN}, by observing that
the left adjoint of a $\kappa$-balanced functor takes values in the framed center.  
\end{proof}

Analogous considerations apply to monads: On a bimodule with framing $\kappa$ we define the 
monad $Z^{[\kappa]}$ by $Z^{[\kappa]}(m) \,{:=}\, \int^{a\in\cala}\! a\,.\,m\,.\,a^{[\kappa-3]}$.
Similarly there are monads in each of the four types of situations for the framed center;
explicit expressions for these are obtained by replacing in the formulas \eqref{eq:4comonads}
the end by an coend and $\kappa$ by $\kappa{-}2$. The monads define corresponding
induction functors $\IZ^{\vec{\kappa}}$; these are the universal functors for
$\kappa{-}2$-balanced right exact functors. 

\begin{cor} \label{cor:ind-co-left}
\Itemizeiii

\item[{\rm (i)}]
The forgetful functor $U \colon \cent^{\vec\kappa}(\vec{\calm}) \,{\to}\, \vec{\calm}$ is left 
adjoint to the co-induction functor $\IZ_{[\vec\kappa]}$ and right adjoint to the induction functor
$\IZ^{[\vec\kappa]}$ that correspond to the endofunctors $Z_{[\vec\kappa]}$ and $Z^{[\vec\kappa]}$,
respectively.

\item[{\rm (ii)}]
The category $\cent^\kappa(\calm)$ is equivalent to the category of modules over the monad
$Z^{[\kappa]}$, and equivalent to the category of comodules over the comonad $Z_{[\kappa]}$.
\end{itemize}
\end{cor}


\subsection{Extensions of the Eilenberg--Watts calculus}\label{sec:moduleEW}

We now collect a few useful results which extend the Eilenberg--Watts calculus of
\cite{fScS2} that is recapitulated in Section \ref{sec:balanc-pair-bimod}.
We first present a mild generalization of the Eilenberg--Watts equivalences
\eqref{Phire...Psire} and \eqref{Phile.Psile}:

\begin{Lemma} \label{Lemma:EW-adjun}
For finite categories $\calm,\calk$ and $\caln$ there are adjoint equivalences 
  \be
  \Funle(\calm \boti \calk, \caln) \simeq \Funle(\calm, \caln \boti \CC\calk)
  \label{eq:EW-inLEX}
  \ee
and
  \be
  \Funre(\calm \boti \calk, \caln) \simeq \Funre(\calm, \caln \boti \CC\calk) \,.
  \ee
\end{Lemma}

\begin{proof}
This is a direct consequence of the equivalences
  \be
  \Funle(\calm \boti \calk, \caln) \,\simeq\, \cc{\calm \boti \calk} \boxtimes \caln
  \,\simeq\, \cc{\calm} \boxtimes (\cc{\calk} \boti \caln)
  \,\simeq\, \Funle(\calm, \CC\calk \boti \caln) \,,
  \ee
and of the corresponding chain of equivalences of categories of right exact functors. Here
the first and third equivalences are given by the Eilenberg--Watts functors \eqref{Phile.Psile}.
The one in the middle follows by elementary associativity properties of the Deligne product.
\end{proof}

The ordinary Eilenberg--Watts equivalences are recovered from this statement by taking $\calm$ 
to be $\vect$. 

We next collect without proof a few statements that involve a lift of the Eilenberg--Watts
calculus to categories of (co)modules over a (co)monad on a functor category.
A proof of these statements is given in \cite{fScS3}, where a module Eilenberg--Watts calculus
is set up which e.g.\ allows for a novel perspective on the center of module categories.
The proof in \cite{fScS3} makes use of the fact that 
for $\Phi\colon\,\calx \,{\leftrightarrows}\, \caly \,\colon \Psi$ an adjoint equivalence
between categories and  $\TX$ a (co)monad on $\calx$, the functor
$\Phi \,{\circ}\, \TX \,{\circ}\, \Psi \,{=:}\, \TY$ is canonically a (co)monad on $\caly$.
Also note that given a (co)monad $\TZ \colon \calm \,{\to}\, \calm$ on a category $\calm$,
for any category $\calx$ the functor category $\Fun(\calm,\calx)$ inherits a (co)monad 
$\TZ^*_{}$ by pre-composition with $\TZ$, i.e.
  \be
  \TZ^*_{}(F) = F \circ \TZ
  \label{eq:inherited-comonad}
  \ee
  for $F \iN \Fun(\calm,\calx)$.

\begin{Proposition} \label{Cor2prop:EW:coendzFz}
Let $\calm \,{\eq}\, {}_\cala \calm_\cala$ be a finite bimodule category over a finite tensor category
$\cala$ and $\cent^\kappa(\calm)$ its $\kappa$-twisted center,
and let $\calx$ be a finite linear category.
\Itemizeiii
\item[{\rm (i)}]
The Eilenberg--Watts calculus provides explicit equivalences
  \be
  \begin{tikzcd}
  \Funle^{\kappa}(\calm,\calx ) \ar{rr}{\simeq} \ar{dr}[left,yshift=-4pt]{\Psile}
  &~& \Funle(\cent^\kappa(\calm),\calx)
  \\
  ~& \CC{\cent^\kappa(\calm)} \boti \calx  \ar{ur}[right,yshift=-3pt]{\Phile}  &~
  \end{tikzcd}
  \label{equation:factor-EW-module}
  \ee
of linear categories. 
Moreover, the category $\Funle^{\kappa}(\calm,\calx)$ is equivalent to the category of 
comodules over the comonad ${(Z_{[\kappa]})}_{\!*}$ on $\Funle(\calm,\calx)$,
  \be
  \Funle^{\kappa}(\calm,\calx) \,\simeq\, {(Z_{[\kappa]})}_{\!*}\text{-comod}(\Funle(\calm,\calx)) \,.
  \ee

\item[{\rm (ii)}]
For any left exact $\kappa{+}2$-balanced functor $F\colon \calm \,{\to}\, \calx$ there is 
an isomorphism 
  \be
  \int^{z \in \cent^\kappa(\calm)}\!\! \cc{z} \boti \widehat{F}(z) \,\cong
  \int^{m \in \calm}\! \cc{m} \boti F(m) \,
  \label{eq:Main-ModuleEW-concret}
  \ee
of objects in $\CC{\cent^\kappa(\calm)} \boti \calx$, where
$\widehat F \,{:=}\, \Phile\cir\Psile(F)$.

\item[{\rm (iii)}]
Specifically, for the co-induction functor $\IZ_{[\kappa]}\colon \calm \,{\to}\,
\cent^\kappa_{}(\calm)$ that corresponds to the co\-mo\-nad $Z_{[\kappa]}$, the corresponding
functor $\widehat{I_{[\kappa]}}\colon \cent^\kappa(\calm) \,{\to}\, \cent^\kappa_{}(\calm) $ 
is the identity functor, whereby the isomorphism \eqref{eq:Main-ModuleEW-concret} reduces to 
  \be
  \int^{z \in \cent^\kappa_{}(\calm)}\! \cc{z} \boti \Id(z)
  \,\cong \int^{m\in\calm}\! \cc{m} \boti Z_{[\kappa]}(m)
  \,\cong \int^{m\in\calm}\! \cc{m} \boti \int_{a\in\cala} a .m. a^{[\kappa-1]} .
  \label{eq:EW:coendzz=...}
  \ee
An analogous isomorphism holds for ends: $ \int_{z\in \cent^\kappa(\calm)} \cc z\boxtimes z
\;{\cong} \int^{m\in\calm} \cc m\boxtimes Z^{[\kappa]}(m)$.
\end{itemize}
\end{Proposition}
 
Combining these assertions with Lemma \ref{Lemma:EW-adjun} we arrive at

\begin{Lemma} \label{Lemma:Module-ew-adj}
The equivalences \eqref{eq:EW-inLEX} induce equivalences 
  \be
  \Funle(\caln_{\cala} {\stackrel{\kappa}{\boti}} {}_{\cala}\calm,\calk)
  \cong \Funle_{\cala}(\caln_{\cala}, \calk \boti \cc{\calm}^{-\kappa+1})
  \label{eq:module-=EW-adj}
  \ee
for any $\kappa \,{\in}\, 2\Z$. 
\end{Lemma}

\begin{proof}
Using that $\caln_{\cala} {\stackrel{\kappa}{\boti}} {}_{\cala}\calm \,{\cong}\,
\Comod_{\cent_{[\kappa]}}(\caln \boti \calm)$, Proposition \ref{Cor2prop:EW:coendzFz}(i)
implies that
  \be
  \begin{array}{rl}
  \Funle(\caln_{\cala}{\stackrel{\kappa}{\boti}}{}_{\cala}\calm,\calk) \!\!\! &
  \cong \Funle( \Comod_{\cent_{[\kappa]}}(\caln \boti \calm), \calk)
  \Nxl6 &
  \cong \Comod_{(\cent_{[\kappa]})^{*}} \Funle(\caln \boti \calm, \calk)
  \Nxl6 &
  \cong \Comod_{\widetilde{\cent_{[\kappa]}}}\Funle(\caln, \calk \boti \cc{\calm}) \,,
  \eear
  \label{eq:combined-equ}
  \ee
where $\widetilde{\cent_{[\kappa]}}$ is the comonad on $\Funle(\caln, \calk \boti \cc{\calm})$
that is induced by the equivalence from Lemma \ref{Lemma:EW-adjun}. We proceed to compute 
$\widetilde{\cent_{[\kappa]}}$. For $F \,{\in}\, \Funle(\caln, \calk \boti \cc{\calm})$ and
$n \,{\in}\, \caln$ we have
  \be
  \begin{array}{rl}
  \widetilde{\cent_{[\kappa]}}(F)(n) \!\! & \dsty
  = \big( \hatPsile(\cent_{[\kappa]})^{*}( \hatPhile(F) \big) (n) 
  \Nxl6 &\dsty
  = \int^{m\in\calm} \!\!\! \int_{a\in\cala} \Hom_{\CC{\calm}}(a.m, F(n.a^{[\kappa+1]})) \boti \cc{m}
  \Nxl6 &\dsty
  = \int^{m\in\calm} \!\!\! \int_{\!a\in\cala} \Hom_{\CC{\calm}}(m, F(n.a^{[\kappa+1]}))
  \boti \cc{a^{\vee}\!.m}
  \Nxl6 &\dsty
  \cong  \int_{\!a\in\cala} \int^{m\in\calm}\!\! \Hom_{\CC{\calm}}(m, F(n.a^{[\kappa+1]}) )
  \boti \cc{a^{\vee}\!.m}
  \Nxl6 &\dsty
  \cong \int_{\!a \in \cala} a^{\vee}\!.\,F(n. a^{[\kappa+1]})
  \,\cong \int_{\!a\in\cala} a\,.\,F(n.a^{[\kappa]}) \,.
  \eear
  \label{eq:ind-comon}
  \ee
As a consequence, a $\widetilde{\cent_{[\kappa]}}$-comodule structure on $F$ consists of
a coherent family of morphisms $F(n) \,{\to}\, a. F(n.a^{[\kappa]})$ for $n \,{\in}\, \caln$. By
adjunction, this is equivalent to a family of coherent natural isomorphisms
$F(n.a) \,{\xrightarrow{\,\cong\,}}\, a^{[\kappa-1]}.F(n)$ which, in turn, is equivalent to 
$F$ belonging to the category $\Funle_\Cala(\caln_\Cala, \calk \boti \CC{\calm}^{-\kappa+1})$.
\end{proof}

As a particular case we obtain 

\begin{cor} \label{cor:lex-framed-app}
For bimodule categories ${}_{\cala}\calm_{\calb}$ and ${}_{\cala}\caln_{\calb}$,
the Eilenberg--Watts equivalences induce an equivalence 
  \be
  \cenf{ \CC\calm \,{\stackrel1\boxtimes}\, \caln \,{\stackrel1\boxtimes}}\,
  \,\simeq\, \Funle_{\Cala,\calb}(\calm,\caln)
  \label{eq:equ-lex-framed-cent-0}
  \ee
of categories.
\end{cor}

\begin{proof}
For the $\cala$-action we have 
  \be
  \begin{array}{rl}
  {}_{\cala}\caln \,{\stackrel1\boxtimes}\, {}_{\cala}\CC\calm \!\!\! &
  \simeq \Funle\big(\CC{\caln \,{\stackrel1\boxtimes}\, \CC\calm }, \vect \big)
  \,\simeq\, \Funle(\calm \,{\stackrel{-1}\boxtimes}\, \CC\caln,\vect )
  \Nxl4 &
  =\, \Funle(\calm \,{\stackrel{0}\boxtimes}\, {}^{[-1]\!}\caln,\vect)
  \,\simeq\, \Funle_\Cala(\calm_\Cala,\caln_\Cala) \,,
  \eear
  \label{eq:Lex-tw}
  \ee
where the last step uses Lemma \ref{Lemma:Module-ew-adj} as well as the canonical 
equivalence $\CC{ \!{\phantom|}^{[-1]}_{} {\CC\caln} }^{[1]} {\simeq}\, \caln_\Cala$.
The $\calb$-action is treated analogously. 
\end{proof}


\subsection{Twisted identity bimodules}

Recall from Section \ref{sec:duals} the twisted variant
${}_{}^{\kappa_1\!\!}\calm{}^{\kappa_2}_{}$ of a bimodule $\calm$, for any pair $\kappa_1$,
$\kappa_2$ of even integers. If $\calm \,{\eq}\, \cala$ is the regular $\cala$-bimodule,
the double duality functor $(-)^{\vee\vee}$ (which is monoidal) provides a distinguished
equivalence
  \be
  {}_{}^{\kappa_1\!\!}\cala{}^{\kappa_2}_{}
  \simeq {}_{}^{\kappa_1-2\!\!}\cala{}^{\kappa_2+2}_{} ,
  \label{eq:kappaA=Akappa}
  \ee
of bimodules, and thus by iteration we get in particular 
${}_{}^{\kappa\!\!}\calm{}^{0} \,{\simeq}\, {}_{}^{0\!\!}\calm{}^{\kappa}$, i.e.\
${}_{\phantom|}^{\kappa\!\!\!}\cala \,{\simeq}\, \cala{}^{\kappa}_{\phantom|}$. 
Furthermore, we have

\begin{Lemma}\label{lem:ZMA=M}
Let $\cala$ be a finite tensor category, $\calm$ a right and $\caln$ a left $\cala$-module,
and let $\kappa \,{\in}\, 2\Z$. The functors
  \be
  \begin{array}{rl}
  \widetilde{\rho} \colon & \calm \boti \cala \,\xrightarrow{~~~} \calm \,,
  \Nxl4
  & ~~m \boti a \,\xmapsto{~~~} m\,.\,a^{[\kappa]}
  \eear\qquad
  \label{eq:rho-abst}
  \ee
and
  \be
  \qquad\begin{array}{rl}
  \widetilde{\lambda}: \quad \cala\boti \caln \,\xrightarrow{~~~} &\!\! \caln \,,
  \Nxl3
  a\boti n \,\xmapsto{~~~} &\!\! a^{[-\kappa]} .\,n \hspace*{4.2em}
  \eear
  \label{equation:EquANA}
  \ee
furnish distinguished equivalences
  \be
  \cenf{\calm\,{\stackrel \kappa\boxtimes}\,\cala^{-\kappa}_{}} \,\tosim\, \calm
  \qquad\text{and}\qquad
  \cenf{{}^{-\kappa\!\!\!}_{}\cala \,{\stackrel \kappa\boxtimes}\, \caln} \,\tosim\, \caln
  \label{eq:ZMA=M}
  \ee
of module categories, given explicitly by
  \be
  \begin{array}{rl}
  \rho \colon & \cenf{\calm\,{\stackrel \kappa\boxtimes}\, \cala} \,\xrightarrow{~~~} \calm \,,
  \Nxl4
  & ~~m \boti a \,\xmapsto{~~~} \Hom_{\cala}(D_\Cala,a) \otik m 
  \eear
  \label{eq:rho-expl}
  \ee
and
  \be
  \begin{array}{rl}
  \lambda \colon & \cenf{\cala\,{\stackrel \kappa\boxtimes}\, \caln} \,\xrightarrow{~~~} \caln \,,
  \Nxl3
  & ~~a \boti n \,\xmapsto{~~~} \Hom_{\cala}(D_\Cala,a) \otik n 
  \eear~~
  \label{eq:lambda-expl}
  \ee     
\end{Lemma}

\begin{proof}
The linear functor \eqref{eq:rho-abst}
has a natural structure of a $\kappa$-ba\-lan\-ced module functor, i.e.\ we have
$\widetilde{\rho}(m.x \boti a) \,{\cong}\, \widetilde{\rho}(m \boti {}^{[\kappa]}x.a)$.
By the universal property of $Z_{[\kappa]}$ from Proposition \ref{proposition:framed-center-univ},
$\widetilde\rho$ therefore induces a module functor
$\rho'\colon \cenf{\calm{\stackrel \kappa\boxtimes}\cala^{-\kappa}_{}} \,{\to}\,\calm$; we claim 
that $\rho'\,{=}\,\rho$ is the functor defined in \eqref{eq:rho-expl}. To this end, set
$\widetilde\rho^{-1}(m)\,{:=}\,m \boti \one$ and 
$\rho^{-1}\,{:=}\,Z_{[\kappa]} \,{\circ}\, \widetilde\rho^{-1}$ and consider the diagram 
  \be
  \begin{tikzcd}[row sep=3.2em]
  \calm \boti \cala \ar{dr}[xshift=-2pt]{\widetilde\rho}
  \ar{rr}{\widetilde\rho^{-1} \circ\, \widetilde\rho} \ar{d}[swap]{Z_{[\kappa]}}
  & ~ & \calm \boti \cala \ar{d}{Z_{[\kappa]}}
  \\
  \cenf{\calm{\stackrel \kappa\boxtimes}\cala^{-\kappa}_{}} \ar{r}[swap]{\rho'}
  & \calm \ar{r}[swap,xshift=-1pt,yshift=1pt]{\rho^{-1}}
  \ar{ur}[xshift=3pt,yshift=-4pt]{\widetilde\rho^{-1}}
  & \cenf{\calm{\stackrel \kappa\boxtimes}\cala^{-\kappa}_{}}    
  \end{tikzcd}
  \label{eq:diag-rho}
  \ee
Since $Z_{[\kappa]}$ is $\kappa$-balanced, the functor
  \be
  \begin{array}{rl}
  \rho^{-1} \colon & \calm\, \xrightarrow{~~~} \cenf{\calm\,{\stackrel \kappa\boxtimes}\, \cala}
  \Nxl5
  & ~m \,\xmapsto{~~~} Z_{[\kappa]}(m\boti\one) = \int_{a\in\cala}m.a^{[\kappa-1]}\boti a
  \eear
  \label{eq:M-to-M.kappa.A}
  \ee
satisfies
  \be
  m.a \,\xmapsto{~~~} Z_{[\kappa]}(m.a\boti \one) \cong
  Z_{[\kappa]}(m\boti a^{[-\kappa]}) = Z_{[\kappa]}(m\boti \one) \otimes a^{[-\kappa]} 
  \ee
and is thus a module functor. Moreover, the functors $\rho'$ and $\rho^{-1}$ are
quasi-inverse, hence $\rho'$ is an equivalence: By the balancing of $Z_{[\kappa]}$, the functor
$Z_{[\kappa]} \,{\circ}\, \widetilde\rho^{-1} \,{\circ}\, \widetilde{\rho}$ is isomorphic to
$Z_{[\kappa]}$ as a right $\cala$-module functor, and thus by the universal property of
$Z_{[\kappa]}$ we have $\rho^{-1} \,{\circ}\, \rho' \,{\cong}\, \Id$ as module functors. It is
even more direct to see that 
$\widetilde{\rho} \,{\circ}\, \widetilde\rho^{-1} \,{\cong}\, \id_{\calm}$
as module functors, so that we also have $\rho' \,{\circ}\, \rho^{-1} \,{\cong}\, \Id$. 
\\[.2em]
We now compute the functor $\rho'$ explicitly. Denote by
$U\colon \calm\,{\stackrel \kappa\boxtimes}\, \cala \,{\to}\, \calm \boti \cala$ the
forgetful functor and consider the diagram 
  \be
  \begin{tikzcd}[row sep=2.5em]
  \calm \boti \cala^{\kappa} \ar{d}[swap,xshift=2pt]{Z_{[\kappa]}} \ar{r}{\widetilde\rho}
  & \calm
  \\ \calm\,{\stackrel \kappa\boxtimes}\, \cala \ar{r}{U}
  & \calm\boti\cala \ar{u}[swap]{\Hom_{\cala}(D_\Cala,\reflectbox{$\scriptstyle?$}) \otimes_\ko ?}
  \end{tikzcd}
  \label{eq:diag-rho-2}
  \ee
This diagram commutes up to a module natural isomorphism:
By \eqref{eq:end-coend-D} there is a distinguished isomorphism $\int_{a\in\cala}
\cc{a^{\vee \vee}{\otimes}\,D_\Cala} \boti a \;{\cong} \int^{a\in\cala} \cc{a} \boti a$,
and hence we obtain a distinguished isomorphism
  \be
  \begin{array}{rl}
  \big(\Hom(D_\Cala,-)\oti{-}\big) \circ U \circ Z_{[\kappa]}(m \boti b) \!\!& \dsty
  = \int_{\!a\in\cala}\! \Hom_{\cala}(D_\Cala,a \oti b) \otik m. a^{[\kappa-1]}
  \Nxl5 & \dsty
  \cong \int_{\!a\in\cala}\! \Hom(a^{\vee\vee}{\otimes}\, D_\Cala,b) \otik m.a^{[\kappa]}
  \Nxl5 & \dsty
  \cong \int^{a\in\cala}\!\! \Hom(a,b) \otik a^{[\kappa]} 
  \,\cong\, m.b^{[\kappa]}
  \eear
  \label{eq:comp-thre}
  \ee
for any $m \boti b \,{\in}\, \calm \boti \cala$.
This shows that $\rho' \,{=}\, \Hom_{\cala}(D_\Cala,\reflectbox?) \otik ? \,{\equiv}\, \rho$,
as claimed.
\\[3pt]
To show the second of the equivalences \eqref{eq:ZMA=M}, consider the functor
  \be
  \begin{array}{rl}
  \caln \,\xrightarrow{~~~} &\!\!
  \cenf{{}^{-\kappa\!\!\!}_{\phantom|}\cala \,{\stackrel \kappa\boxtimes}\,\caln} \,,
  \Nxl6
  n \,\xmapsto{~~~} &\!\!
  Z_{[\kappa]}(\one\boti n) = \int_{a\in\cala} a^{[\kappa-1]}\boti a. n \,.
  \eear
  \ee
This is an $\cala$-module functor, and it is straightforward to check that it is an 
equivalence with a quasi-inverse given by the functor that is induced by the functor
$\widetilde{\lambda}$ as given in \eqref{equation:EquANA}, which is $\kappa{}$-balanced 
and is an $\cala$-module functor.
The explicit form of $\lambda$ follows by analogous considerations as those leading 
to \eqref{eq:comp-thre}.
\end{proof}

Next we consider the case that $\kappa$ is odd. Recall that $\II_0 \,{=}\, \cala$. 
We will construct an equivalence $\calm  \,{\stackrel \kappa\boxtimes}\, \cc{\II_{0}} 
\,{\simeq}\, \calm$, where $\calm_{\cala}$ is a right module as in the case of even $\kappa$, 
but now the action on $\calm$ gets balanced with the right action on $\cc{\II_{0}}$, whereby
the remaining action on $\calm \,{\stackrel \kappa\boxtimes}\, \cc{\II_{0}}$ is a
left $\cc{\cala}$-action. Accordingly also on the right hand side we need to work with a 
left $\cc{\cala}$-action, which we obtain in the form $\cc{a}.m \,{:=}\, m.a^{[\kappa]}$.
In accordance with the notation in \eqref{eq:opp-module} we denote the resulting module by
$~{}_{\cc{\cala}}^{\!\!\!\!\![-\kappa]}\!\calm$. The analogous notation for left modules
${}_{\cala}\caln$ is ${\caln}_{\cc{\cala}}^{[-\kappa]}$ with $n.\cc{a} \,{=}\, {}^{[\kappa]}a.n$.

\begin{Lemma}
Let $\cala$ be a finite tensor category, $\calm$ a right and $\caln$ a left $\cala$-module, and
let $\kappa \,{\in}\, 2\Z{+}1$. There are canonical equivalences
  \be
  \cenf{\calm\,{\stackrel \kappa\boxtimes}\,\cc{\II_0}}
  \,\simeq ~~{}_{\cc{\cala}}^{\!\!\!\!\![-\kappa]}\!\calm \qquad\text{and}\qquad
  \cenf{\cc{\II_0} \,{\stackrel \kappa\boxtimes}\, \caln}
  \,\simeq\, {\caln}_{\cc{\cala}}^{[\kappa-2]}
  \ee
of module categories.
\end{Lemma}

\begin{proof}
We treat explicitly the case of a right module $\calm$, which is analogous to the proof of 
Lemma \ref{lem:ZMA=M}. Again we define a diagram 
  \be
  \begin{tikzcd}[row sep=3.2em]
  \calm \boti \cc{\II_0} \ar{dr}[xshift=-2pt]{\widetilde\rho}
  \ar{rr}[xshift=-1pt]{\widetilde\rho^{-1} \circ\, \widetilde\rho} \ar{d}[swap]{Z_{[\kappa]}}
  & ~ & \calm \boti \cc{\II_0} \ar{d}{Z_{[\kappa]}}
  \\
  \cenf{\calm\,{\stackrel \kappa\boxtimes}\,\cc{\II_0}} \ar{r}[swap]{\rho}
  & \calm \ar{r}[swap,xshift=-1pt,yshift=1pt]{\rho^{-1}}
  \ar{ur}[xshift=3pt,yshift=-4pt]{\widetilde\rho^{-1}}
  & \cenf{\calm\,{\stackrel \kappa\boxtimes}\,\cc{\II_0}}
  \end{tikzcd}
  \ee
of functors, where this time we set $\widetilde\rho(m \boti \cc{a}) \,{:=}\, 
m.(D_\Cala \oti a^{[\kappa-2]})$. The functor $\widetilde\rho$ is $\kappa{+}2$-ba\-lan\-ced:
  \be
  \bearll
  \widetilde\rho(m.b \boti \cc{a}) \cong m.(b \oti D_\Cala \oti a^{[\kappa-2}])
  \!\!& \cong m.(D_\Cala \oti {}^{[4]}b \oti a^{[\kappa-2]})
  \Nxl5 &
  \cong \widetilde\rho(m \boti \cc{a. {}^{[\kappa+2]}b}) \,.
  \eear
  \label{eq:bal-rhotilde}
  \ee
The functor $\widetilde\rho$ also obeys
$\widetilde\rho(m\boti\cc{b.a})\,{\cong}\, \widetilde\rho(m\boti a).b^{[\kappa-2]}$ and thus is
a left $\cc{\cala}$-module functor $\widetilde\rho \colon \calm \,\Boti\kappa\,
\cc{\II_0}\,{\to}\, {}_{~~~\cc{\cala}}^{[\kappa-2]}\!\calm$. Hence
by the universal property of the twisted center, $\rho$ is a well defined module functor.  
\\[2pt]
Further, we define 
  \be
  \widetilde\rho^{-1}(m) := m \boti \cc{D_\Cala} \,{\in}\, \calm \boti \cc{\II_0}
  \ee
and set $\rho^{-1}\,{:=}\,Z_{[\kappa]} \,{\circ}\, \widetilde\rho^{-1}$.
This is indeed a module functor, as 
  \be
  \begin{array}{ll}
  \rho^{-1}(m.x) \!\!& \dsty = \int_{a\in\cala} m.x.a \boti \cc{D_\Cala \oti {}^{[\kappa+1]}a}
  \Nxl6 & \dsty
  \cong \int_{a\in\cala} m.a \boti \cc{D_\Cala \oti {}^{[\kappa+1]}({}^{\vee\!}x \oti a)}
  \Nxl6 & \dsty
  \cong \int_{a\in\cala} m.a \boti \cc{{}^{[\kappa-2]}x \oti D_\Cala \oti {}^{[\kappa+1]}a}
  \,=\, {}^{[\kappa-2]}x. \rho^{-1}(m) \,.
  \eear
  \label{eq:module-rho-inv}
  \ee
Analogously as in the proof of Lemma \ref{lem:ZMA=M} it then follows that $\rho$ and $\rho^{-1}$
furnish an equivalence of module categories, as required. Moreover, again analogously as above
we see that $\rho$ is given explicitly by
  \be
  \begin{array}{rl}
  \rho \colon & \cenf{\calm\,{\stackrel \kappa\boxtimes}\, \cc{\II_{0}}} \,\xrightarrow{~~~} \calm
  \Nxl4
  & ~~\,m \boti \cc{a}  \,\xmapsto{~~~} \Hom_{\cala}(a, \one) \otimes m \,.
  \eear
  \label{eq:rho-expl-cc}
  \ee
The case of a left module follows directly by taking opposite categories. 
\end{proof}


\section{Construction of a parallelization} \label{section:Refind}

In this appendix we provide details of the construction of a parallelization $\Pi$ for the
collection of \csfun s for all fine refinements, as defined in Definition \ref{def:Gamma-Pi}(ii).

\subsection{Isomorphisms among block functors of \filldisk s}

As a preparatory step we restrict our attention to a specific class of \defectsurface s,
namely \emph{\filldisk s}, and construct a distinguished isomorphism between the block functors
for any two such disks of the same type. Recall from Definition \ref{def:filldisk} that a
\filldisk\ of type $\XX$ is a \defectsurface\ $\DD_\XX$ together with a set $\delta_\text{tr}$
of transparent defect lines and with a distinguished boundary segment $\partial_\text{outer}$,
called the \emph{outer boundary} of $\DD_\XX$. The removal of $\delta_\text{tr}$ from
$\DD_\XX$ gives a \defectsurface\ for which every gluing circle and gluing interval except
for $\partial_\text{outer}$ is fillable by a disk in the sense of Definition \ref{def:filling},
and the corresponding filling of $\DD_\XX{\setminus}\delta_\text{tr}$ yields a
\defectsurface\ $\XX$ with underlying surface being a disk and with $\partial\XX$
containing at most one free boundary segment.

For a \filldisk\ $\DD \,{=}\, \DD_\XX$ of arbitrary type $\XX$, denote by
$\LL \,{:=}\, \partialg\DD{\setminus}\partialr\DD$ the gluing part of the outer boundary of 
$\DD$. The prescription \eqref{eq:def:om(LL)} provides a distinguished object
  \be
  \om(\LL) \in \ZZ(\LL) \,,
  \ee
called the \emph{\mute} object for $\LL$, in the gluing category for $\LL$.
We will use these \mute\ objects to specify 
particular isomorphisms between functors associated to \filldisk s.   
For the construction of these isomorphisms we restrict our attention temporarily to the
situation that the outer boundary of $\DD$ is a gluing circle $\SS$ all of whose defect points
are transparently labeled. The assumption that $\DD$ is fillable means that $\SS$ is fillable
in the sense of Definition \ref{def:filling}; also, each of its defect points is labeled either
by $\II$ for one and the same finite tensor category $\cala$ or by $\CC\cala$. In the sequel
we write $\II$ for $\cala$ in order to remind us that it appears as a transparent label.
 
We denote such a gluing circle with $n \,{>}\, 1$ defect points and
$n$-tuple $\kappa \,{=}\, (\kappa_n, \kappa_{n-1}, ...\,, \kappa_1)$ of framing indices (and
with corresponding orientations $\epsilon_i \,{\in}\, \{1,-1\}$ of the defect points)
by $\SS_{n,\kappa}$, i.e.
    \def\locpa  {1.8}    
    \def\locpc  {2.2pt}  
    \def\locpd  {0.12}   
  \be
  \raisebox{-6.3em}{\begin{tikzpicture}
  \scopeArrow{0.34}{>} \scopeArrow{0.58}{>} \scopeArrow{0.85}{>}
  \draw[line width=\locpc,color=\colorCircle,postaction={decorate}]
           node[left=2.3cm,color=black] { } (0,0) circle (\locpa);
  \end{scope} \end{scope} \end{scope}
  \draw[line width=\locpc,color=white,dashed] (0.866*\locpa,-0.5*\locpa) arc (-30:35:\locpa);
  \filldraw[color=\colorDefect]
           (0,\locpa) circle (\locpd) node[above=-5pt,color=black]{$~ \II^{\epsilon_2}_{\phantom|}$}
           (-\locpa,0) circle
                   (\locpd) node[right=1pt,yshift=-1pt,color=black]{$\II^{\epsilon_1}_{\phantom|}$}
           (0,-\locpa) circle (\locpd) node[below=2pt,color=black]{$~\II^{\epsilon_n}_{\phantom|}$}
           (0.7*\locpa,0.79*\locpa) node[above,xshift=-3pt,color=black]{$\kappa_2$}
           (-\locpa,-\locpa) node[above,xshift=12pt,color=black]{$\kappa_n$}
           (-\locpa,\locpa) node[below=7pt,xshift=5pt,color=black]{$\kappa_1$};
  \node at (-\locpa-1.6,0) {$\SS_{n,\kappa} ~=$};
  \end{tikzpicture}}
  \label{eq:def:Snkappa}
  \ee  
Further, denote by $\SS_{n,\kappa}\bs i$, for $i \,{\in}\,\{1,2,...\,,n\}$, 
the fillable circle that is obtained by removing the $i$th defect point from the circle
$\SS_{n,\kappa}$. Thus $\SS_{n,\kappa}\bs i$ is a circle of type $\SS_{n-1,\kappa'}$ with
framing indices $\kappa' \,{=}\,
(\kappa_n, ...\,, \kappa_{i+2},\kappa_{i+1}{+}\kappa_{i},\kappa_{i-1}, ...\,, \kappa_1)$.
Then for any possible choice of $n$, $\kappa$ and $i$ we consider the two \defectsurface s
    \def\locpa  {3.0}    
    \def\locpA  {0.96}   
    \def\locpc  {2.2pt}  
    \def\locpC  {1.8pt}  
    \def\locpd  {0.12}   
    \def\locpD  {0.08}   
    \def\locpf  {1.1pt}  
    \def\locpt  {0.34}   
    \def\locpu  {105}    
    \def\locpv  {62}     
    \def\locpw  {151}    
    \def\locpx  {233}    
    \def\locpy  {281}    
  \be
  \raisebox{-8.0em}{\begin{tikzpicture}[scale=1.1]
  \scopeArrow{0.26}{>} \scopeArrow{0.52}{>} \scopeArrow{0.72}{>} \scopeArrow{0.85}{>}
  \filldraw[\colorCircle,fill=\colorSigma,line width=\locpc,postaction={decorate}] (0,0) circle (\locpa);
  \end{scope} \end{scope} \end{scope} \end{scope}
  \draw[line width=\locpc,color=white,dashed] (-30:\locpa) arc (-30:35:\locpa);
  \scopeArrow{0.66}{\arrowDefect}
  \end{scope}
  \scopeArrow{0.61}{\arrowDefectB}
  \draw[line width=\locpc,color=\colorTransp,postaction={decorate}] (\locpv:\locpA) -- (\locpv:\locpa);
  \draw[line width=\locpc,color=\colorTransp,postaction={decorate}] (\locpw:\locpA) -- (\locpw:\locpa);
  \draw[line width=\locpc,color=\colorTransp,postaction={decorate}] (\locpx:\locpA) -- (\locpx:\locpa);
  \draw[line width=\locpc,color=\colorTransp,postaction={decorate}] (\locpy:\locpA) -- (\locpy:\locpa);
  \end{scope}
  \scopeArrow{0.51}{\arrowDefectB}
  \draw[line width=\locpc,color=\colorTransp,postaction={decorate}] (\locpu:\locpA)
             node[left=2pt,yshift=10pt,color=black,rotate=43] {$\kappa_i$}
	     node[right=-4pt,yshift=6pt,color=black] {$\kappa_{i+1}$}
	     -- (\locpu:0.6*\locpa+0.5*\locpA);
  \end{scope}
  \node[rotate=27] at (115:\locpa+0.23) {$\kappa_i{+}\kappa_{i+1}$};
  \scopeArrow{0.39}{>} \scopeArrow{0.55}{>} \scopeArrow{0.73}{>}
  \filldraw[\colorCircle,fill=white,line width=\locpc,postaction={decorate}] (0,0) circle (\locpA);
  \end{scope} \end{scope} \end{scope}
  \draw[line width=\locpc,color=white,dashed] (-25:\locpA) arc (-25:30:\locpA);
  \filldraw[\colorCircle,fill=white,line width=\locpC] (\locpu:0.6*\locpa+0.5*\locpA) circle (\locpt);
  \scopeArrow{0.54}{<}
  \draw[\colorCircle,line width=\locpf,postaction={decorate}]
             (\locpu:0.6*\locpa+0.5*\locpA) circle (\locpt);
  \end{scope}
  \filldraw[color=\colorDefect] (\locpv:\locpa) circle (\locpd) (\locpv:\locpA) circle (\locpd)
             (\locpw:\locpa) circle (\locpd) (\locpw:\locpA) circle (\locpd)
             (\locpx:\locpa) circle (\locpd) (\locpx:\locpA) circle (\locpd)
             (\locpy:\locpa) circle (\locpd) (\locpy:\locpA) circle (\locpd)
             (\locpu:\locpA) circle (\locpd) (\locpu:0.6*\locpa+0.5*\locpA-\locpt) circle (\locpD);
  \node at (-\locpa-1.8,0) {$\Sigma_{n,\kappa}\bs i ~:=$};
  \end{tikzpicture}
  \label{eq:Sigma-fig2}}
  \ee
~\\[-1.8em]
and
~\\[-1.9em]
    \def\locpa  {3.0}    
    \def\locpA  {0.96}   
    \def\locpc  {2.2pt}  
    \def\locpC  {1.8pt}  
    \def\locpd  {0.12}   
    \def\locpD  {0.08}   
    \def\locpf  {1.1pt}  
    \def\locpt  {0.34}   
    \def\locpu  {105}    
    \def\locpv  {62}     
    \def\locpw  {151}    
    \def\locpx  {233}    
    \def\locpy  {281}    
  \be
  \raisebox{-7.3em}{\begin{tikzpicture}
  \scopeArrow{0.22}{<} \scopeArrow{0.37}{<} \scopeArrow{0.54}{<} \scopeArrow{0.72}{<}
  \scopeArrow{0.86}{<}
  \filldraw[\colorCircle,fill=\colorSigma,line width=\locpc,postaction={decorate}] (0,0) circle (\locpa);
  \end{scope} \end{scope} \end{scope} \end{scope} \end{scope}
  \draw[line width=\locpc,color=white,dashed] (-30:\locpa) arc (-30:35:\locpa);
  \scopeArrow{0.66}{\arrowDefect}
  \end{scope}
  \scopeArrow{0.61}{\arrowDefectB}
  \draw[line width=\locpc,color=\colorTransp,postaction={decorate}] (\locpv:\locpA) -- (\locpv:\locpa);
  \draw[line width=\locpc,color=\colorTransp,postaction={decorate}] (\locpw:\locpA) -- (\locpw:\locpa);
  \draw[line width=\locpc,color=\colorTransp,postaction={decorate}] (\locpx:\locpA) -- (\locpx:\locpa);
  \draw[line width=\locpc,color=\colorTransp,postaction={decorate}] (\locpy:\locpA) -- (\locpy:\locpa);
  \end{scope}
  \scopeArrow{0.76}{\arrowDefectB}
  \draw[line width=\locpc,color=\colorTransp,postaction={decorate}]
	     (\locpu:0.44*\locpa+0.5*\locpA) -- (\locpu:\locpa)
	     node[left=13pt,yshift=6pt,color=black,rotate=43] {$\kappa_i^{}$}
	     node[right=6pt,yshift=9pt,color=black] {$\kappa_{i+1}$};
  \end{scope}
  \node[rotate=27] at (113:\locpA+0.22) {$\kappa_i{+}\kappa_{i+1}$};
  \scopeArrow{0.24}{<} \scopeArrow{0.54}{<} \scopeArrow{0.73}{<}
  \filldraw[\colorCircle,fill=white,line width=\locpc,postaction={decorate}] (0,0) circle (\locpA);
  \end{scope} \end{scope} \end{scope}
  \draw[line width=\locpc,color=white,dashed] (-25:\locpA) arc (-25:30:\locpA);
  \filldraw[\colorCircle,fill=white,line width=\locpC] (\locpu:0.44*\locpa+0.5*\locpA) circle (\locpt);
  \scopeArrow{0.54}{<}
  \draw[\colorCircle,line width=\locpf,postaction={decorate}]
             (\locpu:0.44*\locpa+0.5*\locpA) circle (\locpt);
  \end{scope}
  \filldraw[color=\colorDefect] (\locpv:\locpa) circle (\locpd) (\locpv:\locpA) circle (\locpd)
             (\locpw:\locpa) circle (\locpd) (\locpw:\locpA) circle (\locpd)
             (\locpx:\locpa) circle (\locpd) (\locpx:\locpA) circle (\locpd)
             (\locpy:\locpa) circle (\locpd) (\locpy:\locpA) circle (\locpd)
             (\locpu:\locpa) circle (\locpd) (\locpu:0.44*\locpa+0.5*\locpA+\locpt) circle (\locpD);
  \node at (-\locpa-1.8,0) {$\widetilde\Sigma_{n,\kappa}\bs i ~:=$};
  \end{tikzpicture}}
  \ee
The boundary of $\Sigma_{n,\kappa}\bs i$ and of $\widetilde\Sigma_{n,\kappa}\bs i$
is the union of the circles $\SS_{n,\kappa}$ and $\SS_{n,\kappa}\bs i$ (respectively their
opposites) and of a \emph{tadpole circle} $\QQ_\pm$, i.e.\ a fillable circle with a single
transparently labeled defect point, as indicated in \eqref{eq:pic:tadpolecirc}.
We regard these surfaces as bordisms 
  \be
  \Sigma_{n,\kappa}\bs i:\quad \SS_{n,\kappa} \,{\sqcup}\, \QQ_+ \,{\to}\, \SS_{n,\kappa}\bs i
  \qquad\text{and}\qquad
  \widetilde\Sigma_{n,\kappa}\bs i:\quad \SS_{n,\kappa}\bs i \,{\sqcup}\, \QQ_-
  \,{\to}\, \Sigma_{n,\kappa}\bs i ,
  \ee
respectively. We then consider the functors 
  \be
  \bearll
  G_{n,\kappa}\bs i := \ZZ(\Sigma_{n,\kappa}\bs i)(-\boti \om(\QQ_{-\epsilon_i})) :~ &
  \ZZ(\SS_{n,\kappa}) \,{\xrightarrow{~~}}\, \ZZ(\SS_{n,\kappa}\bs i) \qquad\text{and}
  \Nxl7
  \widetilde G_{n,\kappa}\bs i := \ZZ(\widetilde\Sigma_{n,\kappa}\bs i)(-\boti\om(\QQ_{\epsilon_i}))
  : & \ZZ(\SS_{n,\kappa}\bs i) \,{\xrightarrow{~~}}\, \ZZ(\SS_{n,\kappa}) \,,
  \eear
  \label{eq:def:G,widetildeG}
  \ee
respectively, where $\om(\QQ_\pm)$ are the \mute\ objects \eqref{eq:tadpole-mute} for the tadpole 
circles. These functors which may be viewed as \csfun s
of the form described in \eqref{eq:def:Zprime} for two fine refinements $(\Sigma;\Sigma')$ for
which $\Sigma$ is a cylinder over a circle $\SS_{n-1,\kappa'}$.

Let us describe the functor $G_{n,\kappa}\bs i$ in detail. First note that for 
$z_1^{\epsilon_1} \,{\boxtimes}\cdots{\boxtimes}\, z_n^{\epsilon_n} \,{\in}\,\ZZ(\SS_{n,\kappa})$
and $\epsilon_i\,{=}\,{-}1$ the object
  \be
  \Hom_{\cala}(z_{i}, \one) \otimes z_1^{\epsilon_1} \boti \cdots
  \widehat{z_i^{\epsilon_i}} \cdots \boti z_n^{\epsilon_n} ~\in \ZZ(\SS_{n,\kappa}\bs i)
  \label{eq:objectGink()}
  \ee
(with the symbol $\widehat{z}$ indicating that the factor $z$ is to be removed from the 
expression) comes canonically with the following balancings:
For $a \,{\in}\, \cala$ the balancing between $z_{i-1}$  and $z_{i+1}$ is, in 
case the orientations are as in the picture \eqref{eq:Sigma-fig2},
  \be
  \hspace*{-1.2em}\begin{array}{ll}
  \multicolumn{2}{l} {\Hom_{\cala}(z_{i}, \one) \otimes z_{i-1}.a \boti \cc{z_{i+1}}
  = \Hom_{\cala}(z_{i}, \one) \otimes (z_{i-1} \otimes a) \boti \cc{z_{i+1}} }
  \Nxl7 \quad&
  \cong \Hom_{\cala}(z_{i} \otimes {}^{[\kappa_{i}]}a, \one) \otimes z_{i-1} \boti \cc{z_{i+1}}
  \,\cong \Hom_{\cala}(z_{i} ,  {}^{[\kappa_{i}-1]}a ) \otimes z_{i-1} \boti \cc{z_{i+1}}
  \Nxl7 &
  \cong \Hom_{\cala}( {}^{[\kappa_{i}-2]}a \otimes z_i,\one) \otimes z_{i-1} \boti \cc{z_{i+1}}
  \,\cong \Hom_{\cala}(z_i,\one) \otimes z_{i-1} \boti \cc{z_{i+1}.{}^{[\kappa_i+\kappa_{i+1}]}a} \,, 
  \eear
  \label{eq:bal-G}
  \ee
and similarly for other combinations of orientations. In these expressions we display
only the relevant part of the object, and we use the balancings of the object
$z_1^{\epsilon_1} \,{\boxtimes} \cdots {\boxtimes}\, z_n^{\epsilon_n}$ in steps two and five, 
and the definition of the module structures on $\II$ in steps one and five.

\begin{Proposition} \label{Proposition:Equivalence-disks}
\Itemizeiii

\item[{\rm (i)}]
The \csfun\ $G_{n,\kappa}\bs i$
assigned to the \defectsurface\ \eqref{eq:Sigma-fig2} is an equivalence 
  \be
  G_{n,\kappa}\bs i :\quad \ZZ(\SS_{n,\kappa}) \xrightarrow{~\simeq~} \ZZ(\SS_{n,\kappa}\bs i) \,.
  \ee
For $\epsilon_i\,{=}\,{-}1$ this functor is given explicitly by
  \be
  G_{n,\kappa}\bs i(z_1^{\epsilon_1} \,{\boxtimes} \cdots {\boxtimes}\, z_n^{\epsilon_n})
  = \Hom_{\cala}(z_{i}, \one) \otimes z_1^{\epsilon_1} \boti \cdots 
  \widehat{z_i^{\epsilon_i}} \cdots \boti z_n^{\epsilon_n}
  \label{eq:Gi-eps-1}
  \ee
for $z_1^{\epsilon_1} \,{\boxtimes}\cdots{\boxtimes}\, z_n^{\epsilon_n} \,{\in}\,\ZZ(\SS_{n,\kappa})$
$($with balancings as described in \eqref{eq:bal-G}$)$, while for $\epsilon_{i}\,{=}\,1$ it is 
  \be
  G_{n,\kappa}\bs i(z_1^{\epsilon_1} \,{\boxtimes} \cdots {\boxtimes}\, z_n^{\epsilon_n})
  = \Hom_\cala(D_\Cala,z_i) \otimes z_1^{\epsilon_1} \boti \cdots 
  \widehat{z_i^{\epsilon_i}} \cdots \boti z_n^{\epsilon_n} .
  \label{eq:Gi-eps+1}
  \ee

\item[{\rm (ii)}]
Similarly, the functor $\widetilde G_{n,\kappa}\bs i$ is an equivalence as well, and for
$\epsilon_{i}\,{=}\,1$ it is given by 
  \be
  \widetilde G_{n,\kappa}\bs i (z_1^{\epsilon_1} \,{\boxtimes} \cdots {\boxtimes}\, z_n^{\epsilon_n})
  = z_1 \boti \cdots Z_{\kappa_i}(z_i^{\epsilon_1} \boti \one)\boti \cdots \boti z_n^{\epsilon_n},
  \label{eq:tGi-eps+1}
  \ee
while for $\epsilon_{i}\,{=}\,{-}1$ it is
  \be
  \widetilde G_{n,\kappa}\bs i (z_1^{\epsilon_1} \,{\boxtimes} \cdots {\boxtimes}\, z_n^{\epsilon_n})
  = z_1 \boti \cdots Z_{\kappa_i}(z_i^{\epsilon_1} \boti \cc{D_\Cala} )\boti \cdots
  \boti z_n^{\epsilon_n} .
  \ee
\end{itemize}
\end{Proposition}

\begin{proof}
(i)\,
Clearly the object \eqref{eq:objectGink()} with balancings \eqref{eq:bal-G} is an object
in the gluing category $\ZZ(\SS_{n,\kappa}\bs i)$, and prescribing
$G_{n,\kappa}\bs i(z_1^{\epsilon_1} \,{\boxtimes} \cdots {\boxtimes}\, z_n^{\epsilon_n})
\,{\in}\, \ZZ(\SS_{n,\kappa}\bs i)$ as in \eqref{eq:Gi-eps-1} defines a functor
$H_{n,\kappa}\bs i \colon \ZZ(\SS_{n,\kappa}) \,{\to}\, \ZZ(\SS_{n,\kappa}\bs i)$. 
That $H_{n,\kappa}\bs i$ indeed coincides with $G_{n,\kappa}\bs i$ as defined in
\eqref{eq:def:G,widetildeG} is seen as follows. Recall that we denote by 
  \be
  U: \quad \ZZ(\SS_{n,\kappa}\bs i)\,{\to}\, \UU(\SS_{n,\kappa}\bs i)
  = (\II^{\epsilon_1}_{})_1 \boti (\II^{\epsilon_2}_{})_2 \boti\cdots
  \widehat{(\II^{\epsilon_i}_{})_i} \cdots\boti (\II^{\epsilon_n}_{})_n
  \ee
the forgetful functor to the category that is obtained from the 
gluing category for $\SS_{n,\kappa}\bs i$ by ignoring the balancings. With the help of $U$ the 
pre-block functor for $\Sigma_{n,\kappa}\bs i$, when evaluated on the \mute\ object of $\QQ_+$,
can be written as 
  \be
  \Zpre(\Sigma_{n,\kappa}\bs i)(z\boti \om(\QQ_+))
  = \Hom_{\UU(\SS_{n,\kappa}\bs i)}(U(-),UH_{n,\kappa}\bs i(z))
  \ee
for $z \,{\in}\, \ZZ(\SS_{n,\kappa})$.
Moreover, it follows directly from the definition of the balancing of $H_{n,\kappa}\bs i(z)$ that
the forgetful functor provides the equalizer
  \be
  \Hom_{\ZZ(\SS_{n,\kappa}\bs i)}(-, H_{n,\kappa}\bs i (z))
  \to \Hom_{\UU(\SS_{n,\kappa}\bs i)}(U(-),UG_{n,\kappa}\bs i(z))
  \ee
for the parallel transport  equations on $\Sigma_{n,\kappa}\bs i$. As a consequence 
we also have
  \be
  G_{n,\kappa}\bs i(z) = \Hom_{\ZZ(\SS_{n,\kappa}\bs i)} (-,H_{n,\kappa}\bs i(z)) \,,
  \ee
and thus $G_{n,\kappa}\bs i \,{=}\, H_{n,\kappa}\bs i$, as claimed.
 \\[.2em]
Further, the explicit form of the functor $\rho$ in \eqref{eq:rho-expl-cc} (applied here to the
case that $\calm \,{=}\, \II$) tells us that $G_{n,\kappa}\bs i$ can also be seen as coming from
the equivalence $\cenf{\calm\,{\stackrel\kappa\boxtimes}\, \cc{\II}} \,{\xrightarrow{~~}}\, \calm$,
with $\calm$ the bimodule labeling the defect point that is adjacent to the one labeled by
$\II$ on $\SS_{n,\kappa}$; thus in particular $G_{n,\kappa}\bs i$ is an equivalence. 
 \\[.2em]
The case $\epsilon_i\,{=}\,1$ is treated analogously.  
 \\[4pt]
(ii)\,
For analyzing the functor $\widetilde G_{n,\kappa}\bs i$ in the case $\epsilon_i\,{=}\,1$, 
we consider an object of the form $z \,{=}$
           $z_1^{\epsilon_{1}} \,{\boxtimes} \cdots z_{i-1}^{\epsilon_{i-1}} {\boxtimes} 
z_{i+1}^{\epsilon_{i+1}}\cdots{\boxtimes}\, z_n^{\epsilon_{n}} \,{\in}\, \ZZ(\SS_{n,\kappa}\bs i)$.
It is straightforward to see that via the balancing of the comonad $Z_{[\kappa]}$ the object
$\widetilde G_{n,\kappa} \bs i(z)$ has a canonical structure of an object in $\ZZ(\SS_{n,\kappa})$.
The pre-block functor for $\widetilde\Sigma_{n,\kappa}\bs i$ takes the values
  \be
  \bearll
  \Zpre(\widetilde\Sigma_{n,\kappa}\bs i)(x,z\boti \om(\QQ_-))
  \Nxl5 \hspace*{5.7em}
  = \Hom_{\UU(\SS_{n,\kappa})}(U(x), 
  z_{1}^{\epsilon_1} \boti {\cdots} \boti z_{i-1}^{\epsilon_{i-1}} \boti \one
  \boti z_{i+1}^{\epsilon_{i+1}} \boti {\cdots} \boti z_n^{\epsilon_n})
  \eear
  \ee
on objects $x \,{\in}\, \ZZ(\SS_{n,\kappa})$, where $\UU(\SS_{n,\kappa}) \,{=}\,
(\II^{\epsilon_1}_{})_1 \boti (\II^{\epsilon_2}_{})_2 \boti\cdots\boti (\II^{\epsilon_n}_{})_n$.
Consider now the category $\widetilde{\UU}(\SS_{n,\kappa}) \,{:=}\, (\II^{\epsilon_1}_{})_1
\boti {\cdots} \boti \big((\II^{\epsilon_{i-1}}_{})_{i-1}{\Boti\kappa}(\II^{\epsilon_i}_{})_i\big)
\boti (\II^{\epsilon_{i+1}}_{})_{i+1} \boti\cdots\boti (\II^{\epsilon_n}_{})_n$ with
corresponding forgetful functor $\widetilde{U}\colon \ZZ(\SS_{n,\kappa})\,{\to}\,
\widetilde{\UU}(\SS_{n,\kappa})$. Using co-induction gives an isomorphism
  \be
  \bearl
  \Hom_{\UU(\SS_{n,\kappa})}(U(x), z_1^{\epsilon_1} \boti \cdots \boti z_{i-1}^{\epsilon_{i-1}}
  \boti \one \boti \cdots z_n^{\epsilon_{n}})
  \Nxl7
  \hspace*{6.5em}
  \cong \Hom_{\widetilde{\UU}(\SS_{n,\kappa})}(\widetilde{U}(x), 
  z_1^{\epsilon_1} \boti \cdots \boti Z_{[\kappa_i]}(z_{i-1}^{\epsilon_{i-1}}\boti\one)
  \boti\cdots\boti z_n^{\epsilon_n}) \,.
  \eear
  \label{eq:once-coind}
  \ee
Next note that the forgetful functor $\widetilde{U}$ may be regarded as a composition
of forgetful functors $\widetilde{U}_{j,j+1}$, each of which is applied to the twisted center 
associated with two adjacent defect points $(j,j+1)$ on $\SS_{n,\kappa}\bs i$; thus
for every \twocell\ of $\widetilde\Sigma_{n,\kappa}\bs i$ there is a corresponding pair of
forgetful functors in the Hom functor on the right hand side of \eqref{eq:once-coind}.
It therefore follows in the same way as in the case of $G_{n,\kappa} \bs i$ that the
block functor is given by $\ZZ(\widetilde\Sigma_{n,\kappa}\bs i)(x,z\boti \om(\QQ_-)) \,{=}\,
\Hom_{\ZZ(\SS_{n,\kappa})}(x,\widetilde G_{n,\kappa} \bs i(z))$, thus proving the explicit 
form of $\widetilde G_{n,\kappa} \bs i$ given in \eqref{eq:tGi-eps+1}. Moreover, again as 
in the case of $G_{n,\kappa} \bs i(z)$, we see from the expression \eqref{eq:M-to-M.kappa.A}
for the functor $\rho^{-1}$ that $\widetilde G_{n,\kappa}\bs i$ is an equivalence. 
 \\	  
The case $\epsilon_i\,{=}\,{-}1$ is treated analogously.  
\end{proof}

The boundary circles $\SS_{n,\kappa}^{}$ and $\SS_{n,\kappa}\bs i$ 
of the \defectsurface s $\Sigma_{n,\kappa} \bs i$ and $\widetilde{\Sigma}_{n,\kappa} \bs i$,
for which all defect points are transparently labeled, can play the role of the outer boundary
of a \emph{\transdisk}, i.e.\ a \filldisk\ of the type shown in \eqref{eq:exampleAtranspdisk}.
We now allow for general fillable circles as well as for fillable intervals, which 
can play the role of the outer boundary of a \filldisk\ of arbitrary type $\XX$, like
e.g.\ the one shown in \eqref{eq:exampletranspdisk}. There are then
obvious analogues $\Sigma_{\XX} \bs i$ and $\widetilde{\Sigma}_{\XX} \bs i$ of the
surfaces $\Sigma_{n,\kappa} \bs i$ and $\widetilde{\Sigma}_{n,\kappa} \bs i$. For instance,
for $\XX$ as in \eqref{eq:exampletranspdisk}, an example of a defect surface
$\Sigma_{\XX}^{}\colon \LL_{\XX}^{\phantom i} \,{\to}\, \LL_{\XX} \bs i$ is given by
    \def\locpa  {3.5}    
    \def\locpA  {1.2}    
    \def\locpc  {2.2pt}  
    \def\locpC  {1.8pt}  
    \def\locpd  {0.12}   
    \def\locpD  {0.08}   
    \def\locpf  {1.1pt}  
    \def\locpt  {0.34}   
    \def\locpu  {150}    
    \def\locpv  {39}     
    \def\locpw  {115}    
    \def\locpx  {76}     
  \be
  \raisebox{-5.6em} 
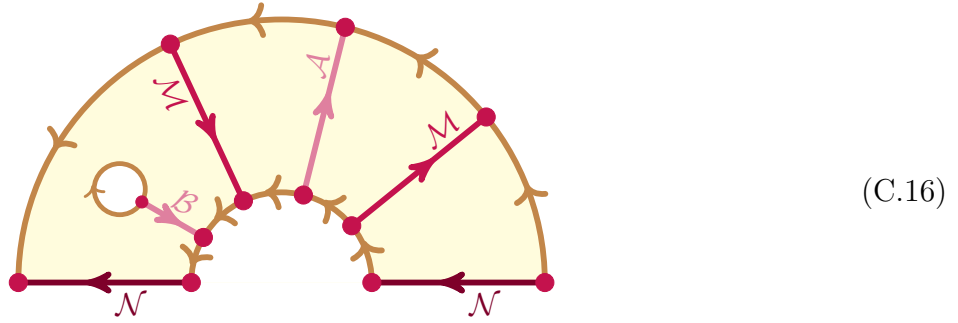
{\begin{tikzpicture}[scale=1.1]
  \scopeArrow{0.19}{<} \scopeArrow{0.48}{<} \scopeArrow{0.69}{<} \scopeArrow{0.90}{<}
  \filldraw[\colorCircle,fill=\colorSigma,line width=\locpc,postaction={decorate}]
             (-\locpa,0) arc (180:0:\locpa);
  \end{scope} \end{scope} \end{scope} \end{scope}
  \scopeArrow{0.12}{<} \scopeArrow{0.30}{<} \scopeArrow{0.50}{<} \scopeArrow{0.71}{<}
  \scopeArrow{0.91}{<}
  \filldraw[\colorCircle,fill=white,line width=\locpc,postaction={decorate}]
             (-\locpA,0) arc (180:0:\locpA);
  \end{scope} \end{scope} \end{scope} \end{scope} \end{scope}
  \scopeArrow{0.61}{\arrowDefect}
  \draw[line width=\locpc,color=\colorFreebdy,postaction={decorate}]
             (-\locpA,0) -- node[below=-1pt,xshift=10pt] {$\caln$} (-\locpa,0);
  \draw[line width=\locpc,color=\colorFreebdy,postaction={decorate}]
             (\locpa,0) -- node[below=-1pt,xshift=10pt] {$\caln$} (\locpA,0);
  \end{scope}
  \scopeArrow{0.61}{\arrowDefect}
  \draw[line width=\locpc,color=\colorDefect,postaction={decorate}]
             (\locpv:\locpA) -- node[sloped,above,near end,yshift=-2pt] {$\calm$} (\locpv:\locpa);
  \draw[line width=\locpc,color=\colorDefect,postaction={decorate}]
             (\locpw:\locpa) -- node[sloped,below,near start,yshift=2pt] {$\calm$} (\locpw:\locpA);
  \draw[line width=\locpc,color=\colorTransp,postaction={decorate}]
             (\locpx:\locpA) -- node[sloped,above,near end,yshift=-2pt] {$\cala$} (\locpx:\locpa);
  \draw[line width=\locpc,color=\colorTransp,postaction={decorate}]
             (\locpu:0.6*\locpa) -- node[sloped,above,yshift=-1pt] {$\calb$} (\locpu:\locpA);
  \end{scope}
  \filldraw[\colorCircle,fill=white,line width=\locpC]
             (\locpu:0.71*\locpa) circle (\locpt);
  \scopeArrow{0.54}{<}
  \draw[\colorCircle,line width=\locpf,postaction={decorate}]
             (\locpu:0.71*\locpa) circle (\locpt);
  \end{scope}
  \filldraw[color=\colorDefect] (\locpa,0) circle (\locpd) (\locpA,0) circle (\locpd)
             (-\locpa,0) circle (\locpd) (-\locpA,0) circle (\locpd)
             (\locpv:\locpa) circle (\locpd) (\locpv:\locpA) circle (\locpd)
             (\locpw:\locpa) circle (\locpd) (\locpw:\locpA) circle (\locpd)
             (\locpx:\locpa) circle (\locpd) (\locpx:\locpA) circle (\locpd)
             (\locpu:\locpA) circle (\locpd) (\locpu:0.71*\locpa-\locpt) circle (\locpD);
  \end{tikzpicture}}
  \label{eq:Bild14.3}
  \ee
  ~\\[-0.7em]
where $\LL_{\XX}$ appears as the inner and $\LL_{\XX} \bs i$ as the outer boundary interval.

The obvious analogue of Proposition \ref{Proposition:Equivalence-disks} holds in this generic 
case, too. That is, the functors 
  \be
  \bearl
  G_{\XX}\bs i := \ZZ(\Sigma_{\XX} \bs i)(-\boti \om(\QQ_{-\epsilon_i})) \qquad\text{and}
  \Nxl8
  \widetilde G_{\XX}\bs i := \ZZ(\widetilde\Sigma_{\XX} \bs i)(-\boti \om(\QQ_{\epsilon_i}))
  \label{eq:def:GXi}
  \eear
  \ee
are equivalences and have similar expressions as in the transparent case. For example,
for an object $z \,{=}\, \cc{x_{\caln}} \boti \cc{a} \boti \cc{x_{\calm}} \boti b
\boti y_{\calm} \boti y_{\caln}$ in the gluing category 
  \be
  \ZZ(\SS_{\XX}) = \CC\caln \boti \CC\cala \boti \CC\calm \boti \cala
  {\boxtimes}\, \calm \boti \caln
  \ee
for the inner boundary interval of the surface \eqref{eq:Bild14.3} one has
  \be
  G_{\XX}\bs 2(z) = \Hom_{\cala}(a, \one) \otimes
  \cc{x_{\caln}} \boti \cc{x_{\calm}} \boti b \boti y_{\calm} \boti y_{\caln} \,. 
  \label{eq:G-genneric}
  \ee

We now establish further properties of the functors \eqref{eq:def:GXi}.
Recall the convention \eqref{eq:insert1orDA} for block functors evaluated at \mute\ objects.
We use this convention in the following statement:

\begin{Lemma}\label{lem:canisoXi1}
There is a canonical isomorphism
    \def\locpa  {1.6}    
    \def\locpb  {0.36}   
    \def\locpc  {2.2pt}  
    \def\locpd  {0.12}   
  \be
  \raisebox{-3.6em}{\begin{tikzpicture}
  \scopeArrow{0.37}{<} \scopeArrow{0.87}{<}
  \filldraw[\colorCircle,fill=\colorSigma,line width=\locpc,postaction={decorate}] (0,0) circle (\locpa);
  \end{scope} \end{scope}
  \scopeArrow{0.52}{\arrowDefectB}
  \draw[line width=\locpc,color=\colorTransp,postaction={decorate}]
             (0:\locpa) node[above=-2pt,xshift=-27pt,color=\colorDefect] {$\II$} -- (180:\locpa);
  \end{scope}
  \filldraw[color=\colorDefect] (0:\locpa) circle (\locpd) (180:\locpa) circle (\locpd);
  \node at (-\locpa-1.25,0) {$\Xi: \quad \ZZ\,\Big($};
  \node at (\locpa+1.39,0) {$\Big) \quad\xrightarrow{~~\cong~~}$};
  \begin{scope}[shift={(2*\locpa+3.88,0)}]
  \scopeArrow{0.37}{<} \scopeArrow{0.87}{<}
  \filldraw[\colorCircle,fill=\colorSigma,line width=\locpc,postaction={decorate}] (0,0) circle (\locpa);
  \end{scope} \end{scope}
  \scopeArrow{0.51}{\arrowDefectB}
  \draw[line width=\locpc,color=\colorTransp,postaction={decorate}]
             (0:\locpa) node[above=-2pt,xshift=-24pt,color=\colorDefect] {$\II$} -- (0:\locpb);
  \end{scope}
  \scopeArrow{0.66}{\arrowDefectB}
  \draw[line width=\locpc,color=\colorTransp,postaction={decorate}]
             (180:\locpb) node[above=-2pt,xshift=-9pt,color=\colorDefect] {$\II$} -- (180:\locpa);
  \end{scope}
  \transv {0:\locpb} node[below=1pt,color=black] {$\scriptstyle \one$};
  \transv {180:\locpb} node[below=1pt,color=black] {$\scriptstyle \cc{D_\Cala}$};
  \filldraw[color=\colorDefect] (0:\locpa) circle (\locpd) (180:\locpa) circle (\locpd);
  \node at (-\locpa-0.74,0) {$\ZZ\,\Big($};
  \node at (\locpa+0.61,0) {$\Big)$};
  \end{scope}
  \end{tikzpicture}}
  \label{eq:canisoXi-a}
  \ee
  ~\\[-1.8em]
of functors.
\end{Lemma}

\begin{proof}
Let us factorize the disk $\DD$ that appears on the right hand side of \eqref{eq:canisoXi-a}
in the way indicated by the dashed circle in
    \def\locpa  {2.7}    
    \def\locpb  {0.96}   
    \def\locpc  {2.2pt}  
    \def\locpC  {1.8pt}  
    \def\locpd  {0.12}   
    \def\locpD  {0.08}   
    \def\locpe  {0.105}  
    \def\locpf  {1.1pt}  
    \def\locpr  {1.58}   
    \def\locpt  {0.23}   
    \def\locpu  {125}    
    \def\locpv  {56}     
  \be
  \raisebox{-7.0em}
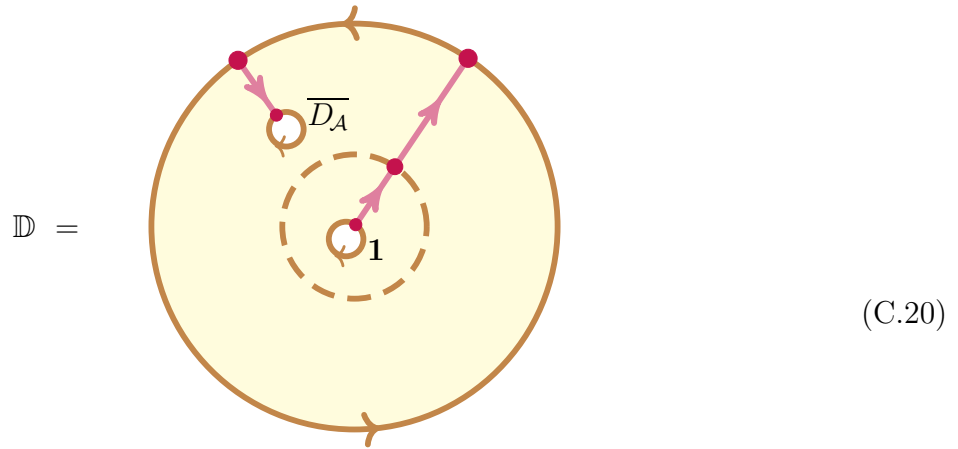
{\begin{tikzpicture}
  \scopeArrow{0.26}{>} \scopeArrow{0.77}{>}
  \filldraw[\colorCircle,fill=\colorSigma,line width=\locpc,postaction={decorate}] (0,0) circle (\locpa);
  \end{scope} \end{scope}
  \draw[\colorCircle,line width=\locpc,dash pattern=on 9pt off 5pt] (0,0) circle (\locpb);
  \scopeArrow{0.29}{\arrowDefect} \scopeArrow{0.77}{\arrowDefect}
  \draw[line width=\locpc,color=\colorTransp,postaction={decorate}]
             (\locpv:-0.2) -- (\locpv:\locpa);
  \end{scope} \end{scope}
  \scopeArrow{0.59}{\arrowDefect}
  \draw[line width=\locpc,color=\colorTransp,postaction={decorate}]
             (\locpu:\locpa) -- (\locpu:\locpr);
  \end{scope}
  \filldraw[\colorCircle,fill=white,line width=\locpc] (\locpv:-0.2) circle (\locpt)
             node[right=3pt,yshift=-4pt,color=black] {$\one$};
  \filldraw[\colorCircle,fill=white,line width=\locpc] (\locpu:\locpr) circle (\locpt)
             node[right=3pt,yshift=6pt,color=black] {$\cc{D_\Cala}$};
  \scopeArrow{0.74}{<}
  \draw[\colorCircle,line width=\locpf,postaction={decorate}] (\locpu:\locpr) circle (\locpt);
  \draw[\colorCircle,line width=\locpf,postaction={decorate}] (\locpv:-0.2) circle (\locpt);
  \end{scope}
  \filldraw[color=\colorDefect] (\locpv:\locpa) circle (\locpd)
             (\locpv:\locpt-0.2) circle (\locpD) (\locpv:\locpb) circle (\locpe)
             (\locpu:\locpa) circle (\locpd) (\locpu:\locpr+\locpt) circle (\locpD);
  \node at (-\locpa-1.4,0) {$\DD~=$};
  \label{eq:pic2tadpfact}
  \end{tikzpicture}}
  \ee
With the help of the description of the functors $G_{}\bs i$ and $\widetilde G_{} \bs i$
in Proposition \ref{Proposition:Equivalence-disks} we see that $\ZZ(\DD) \,{=}\,
Z_{[1]}(\one \boti \cc{D_\Cala}) \,{=} \int_{a\in\cala} a \boti \cc{D_\Cala.{}^{[2]}a}
\,{\cong} \int^{a\in\cala}\! a \boti \cc{a} \,{\in}\, \ZZ(\partial\DD)$, 
using the canonical isomorphism \eqref{eq:disting-iso} of objects in $\ZZ(\partial \DD)$.
Noticing that $\int^{a}\! a \boti \cc{a}$ is the value of the block functor on the 
left hand side of \eqref{eq:canisoXi-a} then establishes the isomorphism $\Xi$.
\end{proof}  

The so obtained isomorphism is a universal morphism in the following sense. Denote by
$\DD_{1}$ and $\DD_{2}$ the \defectsurface s on the left and right hand sides of 
\eqref{eq:canisoXi-a}, respectively, and by $\SS$ their common gluing boundary, with gluing category 
  \be
  \ZZ(\SS) = \cc{\cala} \Boti{-1} \cala \Boti{-1}\,.
  \ee
Then the pre-block functors are given by
  \be
  \bearl \dsty
  \Zpre(\DD_1)(G)= \int^{a\in\cala}\! \Hom_{\cala}(G_{1},a) \otik \Hom_{\cala}(a,G_2)
  \qquad\text{and}
  \Nxl6
  \Zpre(\DD_2)(G)= \Hom_{\cala}(G_{2},\one) \otik \Hom_{\cala}(D_\Cala,G_{\one}) \,,
  \eear
  \ee
respectively, for $G \,{=}\, \cc{G_1} \boti G_2 \iN \ZZ(\SS)$. Using first the isomorphism
between coend and end that follows from the isomorphism \eqref{eq:end-coend-D}, and then the
dinatural transformation of the end, we obtain a canonical morphism
  \be
  \bearll
  \Zpre(\DD_1)(G) \!\!& \dsty
  = \int^{a\in\cala}\! \Hom_{\cala}(G_1,a) \otik \Hom_{\cala}(a,G_2)
  \Nxl2 &\dsty
  \cong \int_{a\in\cala}\! \Hom_{\cala}(G_1,a^{\vee\vee}) \otik \Hom_{\cala}(D_\Cala \oti a, G_2)
  \Nxl6 & \hspace*{4.5em}
  \xrightarrow{~~~} \Hom_{\cala}(G_2,\one) \otik \Hom_{\cala}(D_\Cala,G_1) \,=\, \Zpre(\DD_2)(G) \,.
  \eear
  \label{eq:can-uni-iso}
  \ee

\begin{Lemma} \label{Lemma:univ-blockiso} {}
\Itemizeiii

\item[{\rm (i)}]
The canonical morphism \eqref{eq:can-uni-iso} from $\Zpre(\DD_1)$ to $\Zpre(\DD_2)$
is compatible with all parallel transport  operations on $\DD_1$ and $\DD_2$, i.e.\ it 
commutes with $\hol_{a,x}$ for all objects $a\iN \cala$ and starting points $x$.

\item[{\rm (ii)}]
The morphism \eqref{eq:can-uni-iso} induces a morphism between the parallel transport equalizers,
i.e.\ between the corresponding block functors. The so obtained morphism between the block 
functors is the isomorphism \eqref{eq:canisoXi-a} in Lemma {\rm\ref{lem:canisoXi1}}.
\end{itemize}
\end{Lemma}

\begin{proof}
(i)\, The compatibility with the parallel transport  operations follows by direct computation.
\\[.2em]
(ii)\, We obtain a morphism between the block functors by the universal property of the
equalizer. It is straightforward to check that applying the forgetful functor from blocks to
pre-blocks to the morphism \eqref{eq:canisoXi-a} reproduces the morphism \eqref{eq:can-uni-iso}. 
\end{proof}

As a consequence of Lemma \ref{lem:canisoXi1} we have

\begin{Lemma}\label{lem:canisoXi2}
The isomorphism $\Xi$ in \eqref{eq:canisoXi-a} provides a distinguished adjoint equivalence
between the functors $G_{\XX}\bs i$ and $\widetilde G_{\XX} \bs i$ for any type $\XX$.
\end{Lemma}

\begin{proof}
We factorize the bordism $\widetilde\Sigma_{\XX}\bs i \,{\circ}\,\Sigma_{\XX}\bs i$ in such a way
that one of the factors is the disk $\DD$ on the right hand side of \eqref{eq:canisoXi-a}. The
isomorphism $\Xi$ can then be used to define a natural isomorphism 
$\widetilde G_{\XX}\bs i \,{\circ}\, G_{\XX}\bs i \,{\cong}\, \id_{\SS_{\XX}}$. By Lemma
\ref{lem:ZMA=M}, the two functors are inverse equivalences; as a consequence there is a unique
way to define the isomorphism $G_{\XX}\bs i \,{\circ}\, \widetilde G_{\XX}\bs i 
\,{\xrightarrow{\cong}}\, \id_{\SS_{\XX}\bs i}$ 
in such a way that the equivalence is an adjoint equivalence. 
\end{proof}  

Next we note that by successively applying the functors $\widetilde G_{n,\kappa}\bs i$
with all possible values of $i$ to the transparently labeled gluing circle $\SS_{n,\kappa}$
\eqref{eq:def:Snkappa} we obtain a fillable disk $\DD_{n,\kappa}$ each of whose inner boundaries
is a \tadpolecircle. This allows for the following description of the \mute\ object
$\om(\SS_{n,\kappa})$, as defined according to \eqref{eq:def:om(LL)}:

\begin{Lemma} \label{Lemma:rho-part1}
The \mute\ object $\om_{n,\kappa} \,{:=}\, \om(\SS_{n,\kappa})$ for any transparent gluing
circle $\SS_{n,\kappa}$ can be recovered from the functors $\widetilde G_{\ell,\kappa}\bs{i_\ell}$
with $\ell \eq 2,3,...\,,n$ and with $(i_1,i_2,...\,,i_n)$ any permutation of $(1,2,...\,,n)$
as follows
$($for brevity we abuse notation by writing the same generic label $\kappa$ for all the 
tuples of framing indices involved$){:}$ there is a canonical isomorphism
$\rho \colon \om_{n,\kappa} \,{\tocong}\, \widetilde G(\om(\QQ_{\epsilon_{i_1}}))$ with
$\widetilde G$ the composite
  \be
  \ZZ(\QQ_{\epsilon_{i_1}}) \xrightarrow{~\widetilde G_{2,\kappa}\bs{i_2}~}
  \ZZ(\SS_{2,\kappa}) \xrightarrow{~\widetilde G_{3,\kappa}\bs{i_3}~}
  \ZZ(\SS_{3,\kappa}) \xrightarrow{~\widetilde G_{4,\kappa}\bs{i_4}~} \cdots
  \cdots \xrightarrow{~\widetilde G_{n,\kappa}\bs{i_n}~} \ZZ(\SS_{n,\kappa}) \,.
  \ee
An analogous statement holds for the \mute\ object for the outer boundary of any \filldisk\
of arbitrary type $\XX$. 
\end{Lemma}

\begin{proof}
This statement follows directly by applying the block functor to a situation involving
consecutive gluings each of which involves a single
\tadpolecircle, as indicated in the following picture
    \def\locpa  {3.1}    
    \def\locpb  {0.77}   
    \def\locpB  {1.9}    
    \def\locpc  {2.2pt}  
    \def\locpC  {1.8pt}  
    \def\locpd  {0.12}   
    \def\locpD  {0.08}   
    \def\locpe  {0.105}  
    \def\locpf  {1.1pt}  
    \def\locpr  {2.23}   
    \def\locps  {1.09}   
    \def\locpu  {118}    
    \def\locpv  {140}    
    \def\locpw  {170}    
  \be
  \raisebox{-7.3em}
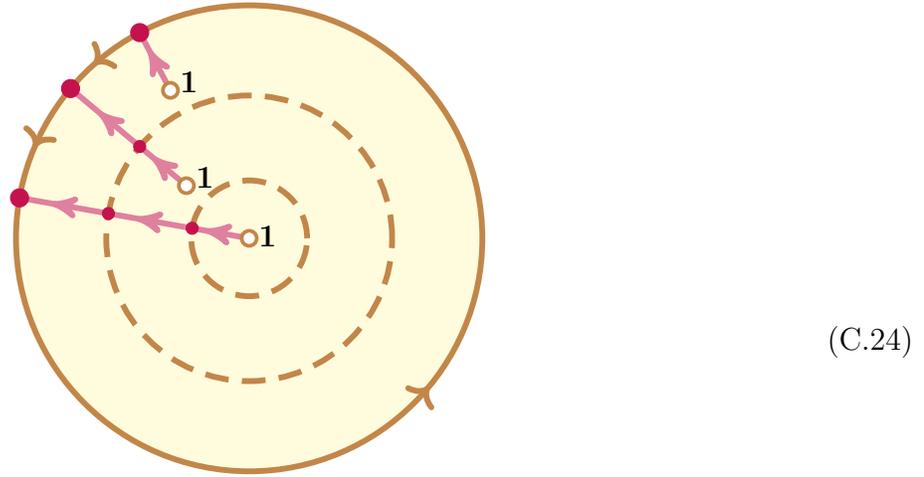
{\begin{tikzpicture}
  \scopeArrow{0.366}{>} \scopeArrow{0.436}{>} \scopeArrow{0.89}{>}
  \filldraw[\colorCircle,fill=\colorSigma,line width=\locpc,postaction={decorate}] (0,0) circle (\locpa);
  \end{scope} \end{scope} \end{scope}
  \draw[\colorCircle,line width=\locpc,dash pattern=on 9pt off 5pt]
             (0,0) circle (\locpb) (0,0) circle (\locpB);
  \scopeArrow{0.69}{\arrowDefect}
  \draw[line width=\locpc,color=\colorTransp,postaction={decorate}] (\locpu:\locpr) -- (\locpu:\locpa);
  \end{scope}
  \scopeArrow{0.29}{\arrowDefect} \scopeArrow{0.75}{\arrowDefect}
  \draw[line width=\locpc,color=\colorTransp,postaction={decorate}] (\locpv:\locps) -- (\locpv:\locpa);
  \end{scope} \end{scope}
  \scopeArrow{0.19}{\arrowDefect} \scopeArrow{0.49}{\arrowDefect} \scopeArrow{0.86}{\arrowDefect}
  \draw[line width=\locpc,color=\colorTransp,postaction={decorate}] (0,0) -- (\locpw:\locpa);
  \end{scope} \end{scope} \end{scope}
  \transV{\locpu:\locpr} {[right=-1pt,yshift=3pt,color=black] {$\one$}};
  \transV{\locpv:\locps} {[right=-1pt,yshift=3pt,color=black] {$\one$}};
  \transV{0,0} {[right=-1pt,yshift=1pt,color=black] {$\one$}};
  \filldraw[color=\colorDefect] (\locpu:\locpa) circle (\locpd)
             (\locpv:\locpa) circle (\locpd) (\locpw:\locpa) circle (\locpd)
             (\locpv:\locpB) circle (\locpD) (\locpw:\locpb) circle (\locpD)
             (\locpw:\locpB) circle (\locpD);
  \end{tikzpicture}}
  \ee
  ~\\[0.1em]
(in the situation shown we have
$n\,{=}\,3$ and $\epsilon_1\,{=}\,\epsilon_2\,{=}\,\epsilon_3\,{=}\,1$).
\end{proof}

It follows that explicit expressions for the \mute\ object $\om_{\DD}$ for the
outer boundary $\partial \DD$ of a \filldisk\ $\DD$ of type $\XX$ can be obtained
by the following procedure: For each defect line $\delta_i$ in $\XX$ take a pair of
variables $(m_{i}, \cc{m_{i}})$ in the categories $\calm_i$ and $\CC\calm_i$ labeling 
the two defect points on $\partial \DD$ at the ends of the defect,
together with the relevant \mute\ objects $\one_{\cala_j}$ and $\cc{D_{\cala_j}}$,
respectively, for the \tadpolecircle s in $\DD$. Build the Deligne product of these
objects and take the coend over the variables $m_i$. Finally apply for every \twocell\
of $\DD$, having $n$ gluing segments on $\partial \DD$ with framings $\{\kappa_{j}\}$,
the corresponding comonads $\ZZ_{[\kappa_{j}]}$ for $n{-}1$ of the
gluing segments (up to canonical isomorphism it does not matter which one of the $n$ 
gluing segments is omitted). The following example illustrates this procedure:

\begin{Example} \label{ex:filldisk}
Consider the \filldisk\ 
    \def\locpa  {2.5}    
    \def\locpA  {1.5}    
    \def\locpc  {2.2pt}  
    \def\locpC  {1.8pt}  
    \def\locpd  {0.12}   
    \def\locpD  {0.08}   
    \def\locpe  {0.09}   
    \def\locpf  {1.1pt}  
    \def\locps  {0.83}   
    \def\locpt  {0.43}   
    \def\locpw  {225}    
    \def\locpx  {132}    
    \def\locpX  {331}    
  \be
  \raisebox{-7.4em}{\begin{tikzpicture}[scale=1.1]
  \scopeArrow{0.13}{>} \scopeArrow{0.51}{>} \scopeArrow{0.80}{>}
  \filldraw[\colorCircle,fill=\colorSigma,line width=\locpc,postaction={decorate}]
             (0,0) circle (\locpa);
  \end{scope} \end{scope} \end{scope}
  \scopeArrow{0.62}{\arrowDefectB}
  \draw[line width=\locpc,color=\colorTransp,postaction={decorate}]
             (\locpw:\locps) -- node[sloped,above,near end, yshift=-2pt]{$\cala$} (\locpw:\locpa);
  \end{scope}
  \filldraw[\colorCircle,fill=white,line width=\locpc]
             (\locpw:\locps) circle (\locpt) node[left=9pt,color=black,yshift=5pt] {$\scriptstyle 2$};
  \scopeArrow{0.91}{<}
  \draw[\colorCircle,line width=\locpf,postaction={decorate}] (\locpw:\locps) circle (\locpt);
  \end{scope}
  \scopeArrow{0.55}{\arrowDefectB}
  \draw[line width=\locpc,color=\colorDefect,postaction={decorate}]
             (\locpx:\locpa) -- node[sloped,above,near end,yshift=-1pt]{$\calm$} (\locpX:\locpa);
  \end{scope}
  \filldraw[color=\colorDefect] 
             (\locpw:\locpa) circle (\locpd) (\locpw:\locps+\locpt) circle (\locpD)
             (\locpx:\locpa) circle (\locpd) (\locpX:\locpa) circle (\locpd);
  \node at (0.5*\locpX+0.5*\locpw:1.08*\locpa) {$\scriptstyle 0$};
  \node at (38:1.07*\locpa) {$\scriptstyle 1$};
  \node at (174:1.06*\locpa) {$\scriptstyle 1$};
  \node at (-\locpa-1.4,0) {$\DD~=$};
  \end{tikzpicture}}
  \label{eq:exa:DDtadLL-a}
  \ee
Using Equation \eqref{eq:end-coend-D} and Lemma \ref{lemma:balancing-psil(F)} we obtain
  \be
  \begin{array}{ll}
  \om_\DD \!\! & \dsty
  = \int^{m \in \calm}\!\! \ZZ_{[0]}(\cc{m} \boti \cc{D_\Cala}) \boti m
  \,= \int^{m \in \calm} \!\!\!\! \int_{\!a \in \cala} \cc{a^{\vee}. m} \boti \cc{D.\, {}^{[3]}a}
  \boti m
  \Nxl5 & \dsty
  \cong \int^{m\in \calm} \!\!\!\! \int^{a\in\cala}\! \cc{a^{\vee}.m} \boti \cc{a} \boti m
  \, \cong \int^{m\in \calm} \!\!\!\! \int^{a\in\cala}\! \cc{m} \boti \cc{a} \boti a.m \,.
  \end{array}
  \label{eq:mute-expl}
  \ee
We have actually already encountered the \defectonemanifold\ that constitutes the
outer boundary of the disk $\DD$: it is the defect circle $\LII_\kappa^{\sss\swarrow}(\calm)$
in \eqref{eq:circle.IIn-M-M} with framing index $\kappa \,{=}\, 0$.
Similarly we obtain the following list of \defectonemanifold s and \mute\ objects 
for all other transparent disks with outer boundaries given by one of the circles
\eqref{eq:circle.IIn-M-M}:
  \be
  \hspace*{-1.5em}\bearll
  \LII_\kappa^{\sss\nearrow}(\calm) :~ & \dsty \om=\,
  \int^{m\in\calm}\!\!\!\int_{a\in\cala} \cc m \boti m\,.\,a^{[\kappa-1]} \boti a
  \,= \int^{m\in\calm}\!\!\!\int_{a\in\cala} \cc m \boti m\,.\,a \boti {}^{[\kappa-1]\!}a \,,
  \Nxl8
  \LII_{-\kappa}^{\sss\swarrow}(\calm) :~ & \dsty \om=\,
  \int^{m\in\calm}\!\!\!\int^{a\in\cala}\!\! \cc m \boti \cc a \boti a^{[\kappa]}.\,m \,,
  \Nxl8
  \LII_\kappa^{\sss\nwarrow}(\calm) :~ & \dsty \om=\,
  \int^{m\in\calm}\!\!\!\int_{a\in\cala} \cc m \boti a \boti {}^{[\kappa-1]\!}a\,.\,m \,,
  \Nxl8
  \LII_{-\kappa}^{\sss\searrow}(\calm) :~ & \dsty  \om=\,
  \int^{m\in\calm}\!\!\!\int^{a\in\cala}\!\! \cc a \boti \cc m \boti m\,.{}^{[\kappa]\!}a \,.
  \eear
  \ee
\end{Example}

Next we show

\begin{Lemma} \label{Lemma:rho}
The functors 
$G \bs i \,{\equiv}\, G_{\XX} \bs i \colon \ZZ(\LL_{\XX}^{}) \,{\to}\, \ZZ(\LL_{\XX} \bs i) $
for fixed type $\XX$ and different values of $i$ commute up to canonical natural
isomorphism, i.e.\ for any pair $i,j$ with $i \,{\ne}\, j$ there is a canonical isomorphism
$\gamma\bs{i,j}\colon G\bs{i,j} \,{\circ}\, G \bs i \,{\cong}\, G \bs{j,i} \,{\circ}\, G\bs j$ 
$($with obvious notation$)$. 
\\
Moreover, these isomorphisms are compatible with the \mute\ objects in the following sense:
For every $i$ there is a canonical isomorphism $\rho\bs i\colon 
G_{\XX}\bs i(\om(\LL_{\XX}^{})) \,{\tocong}\, \om\bs i \,{:=}\, \om(\LL_{\XX} \bs i)$,
and analogous isomorphisms relating the \mute\ objects for the gluing boundaries
$\LL_{\XX}\bs i$ and $\LL_{\XX}\bs{i,j}$ etc., such that the diagram
  \be
  \begin{tikzcd}[column sep=4.7em,row sep=1.2em]
  G\bs{i,j} {\circ}\, G\bs i(\om(\LL_{\XX}^{}))
  \ar{r} {G\bs{i,j}(\rho\bs i)} \ar{dd}[swap]{\gamma\bs{i,j}}
  & G\bs{i,j}(\om\bs i) \ar{rd}[xshift=-4pt]{\rho\bs{i,j}} & ~
  \\
  ~ & ~ & \om\bs{i,j}
  \\
  G\bs{j,i} {\circ}\, G\bs j(\om(\LL_{\XX}^{})) \ar{r} {G\bs{j,i}(\rho\bs j)}
  & G\bs{j,i}(\om\bs j) \ar{ur}[swap,xshift=-10pt]{\rho\bs{j,i}} & ~
  \end{tikzcd}
  \label{eq:coherence-rho}
  \ee
commutes.
\end{Lemma}

\begin{proof}
That the functors $\widetilde G_{n,\kappa}\bs i$ respect the \mute\ objects 
is easily seen graphically. In the transparent case, the relevant situation is
    \def\locpa  {2.6}    
    \def\locpA  {2.0}    
    \def\locpb  {0.85}   
    \def\locpc  {2.2pt}  
    \def\locpC  {1.8pt}  
    \def\locpd  {0.12}   
    \def\locpD  {0.08}   
    \def\locpe  {0.105}  
    \def\locpf  {1.1pt}  
    \def\locpr  {1.51}   
    \def\locpt  {0.23}   
    \def\locpu  {125}    
    \def\locpv  {170}    
  \begin{eqnarray}
  \begin{tikzpicture}
  \scopeArrow{0.42}{>} \scopeArrow{0.77}{>}
  \filldraw[\colorCircle,fill=\colorSigma,line width=\locpc,postaction={decorate}] (0,0) circle (\locpa);
  \end{scope} \end{scope}
  \draw[\colorCircle,line width=\locpc,dash pattern=on 9pt off 5pt] (0,0) circle (\locpb);
  \scopeArrow{0.29}{\arrowDefect} \scopeArrow{0.77}{\arrowDefect}
  \draw[line width=\locpc,color=\colorTransp,postaction={decorate}]
             (\locpv:-0.2) -- (\locpv:\locpa);
  \end{scope} \end{scope}
  \scopeArrow{0.59}{\arrowDefect}
  \draw[line width=\locpc,color=\colorTransp,postaction={decorate}]
             (\locpu:\locpa) -- (\locpu:\locpr);
  \end{scope}
  \filldraw[\colorCircle,fill=white,line width=\locpc]
             (\locpu:\locpr) circle (\locpt) node[right=3pt,color=black] {$\one$}
             (\locpv:-0.2) circle (\locpt) node[above=3pt,xshift=-4pt,color=black] {$\one$};
  \scopeArrow{0.74}{<}
  \draw[\colorCircle,line width=\locpf,postaction={decorate}] (\locpu:\locpr) circle (\locpt);
  \draw[\colorCircle,line width=\locpf,postaction={decorate}] (\locpv:-0.2) circle (\locpt);
  \end{scope}
  \filldraw[color=\colorDefect] (\locpv:\locpa) circle (\locpd)
             (\locpv:\locpt-0.2) circle (\locpD) (\locpv:\locpb) circle (\locpe)
             (\locpu:\locpa) circle (\locpd) (\locpu:\locpr+\locpt) circle (\locpD);
  \node at (-\locpa-0.71,0) {$\ZZ\Big($};
  \node at (\locpa+1.11,0) {$\Big) \qquad =$};
  \begin{scope}[shift={(\locpa+\locpA+3.98,0)}]
  \scopeArrow{0.42}{>} \scopeArrow{0.77}{>}
  \filldraw[\colorCircle,fill=\colorSigma,line width=\locpc,postaction={decorate}] (0,0) circle (\locpA);
  \end{scope} \end{scope}
  \scopeArrow{0.59}{\arrowDefect}
  \draw[line width=\locpc,color=\colorTransp,postaction={decorate}] (\locpv:0.3) -- (\locpv:\locpA);
  \end{scope}
  \scopeArrow{0.59}{\arrowDefect}
  \draw[line width=\locpc,color=\colorTransp,postaction={decorate}]
             (\locpu:\locpA) -- (\locpu:0.7*\locpr);
  \end{scope}
  \filldraw[\colorCircle,fill=white,line width=\locpc]
             (\locpv:0.3) circle (\locpt) node[right=3pt,color=black] {$\one$}
             (\locpu:0.7*\locpr) circle (\locpt) node[right=3pt,color=black] {$\one$};
  \scopeArrow{0.74}{<}
  \draw[\colorCircle,line width=\locpf,postaction={decorate}] (\locpu:0.7*\locpr) circle (\locpt);
  \draw[\colorCircle,line width=\locpf,postaction={decorate}] (\locpv:0.3) circle (\locpt);
  \end{scope}
  \filldraw[color=\colorDefect] (\locpv:\locpA) circle (\locpd) (\locpv:\locpt+0.3) circle (\locpD)
             (\locpu:\locpA) circle (\locpd) (\locpu:0.7*\locpr+\locpt) circle (\locpD);
  \node at (-\locpA-0.78,0) {$\ZZ\Big($};
  \node at (\locpA+0.51,0) {$\Big)$};
  \end{scope}
  \end{tikzpicture}
  \nonumber \\[-2.6em]~ \\[-0.8em] \nonumber~
  \end{eqnarray}
The case of generic type $\XX$ is analogous.
 \\[.2em]
The isomorphism $\gamma\bs{i,j}$
and the commutativity of \eqref{eq:coherence-rho} are seen graphically, by comparing the
following two situations, for which we clearly have $\ZZ(\Sigma\bs{ij})\,{=}\,\ZZ(\Sigma\bs{ji})$:
    \def\locpa  {4.0}    
    \def\locpA  {0.75}   
    \def\locpb  {2.42}   
    \def\locpc  {2.2pt}  
    \def\locpC  {1.8pt}  
    \def\locpd  {0.12}   
    \def\locpD  {0.08}   
    \def\locpe  {0.105}  
    \def\locpf  {1.1pt}  
    \def\locpr  {0.78}   
    \def\locps  {0.86}   
    \def\locpt  {0.26}   
    \def\locpu  {105}    
    \def\locpv  {56}     
    \def\locpw  {151}    
    \def\locpx  {212}    
    \def\locpy  {288}    
    \def\locpz  {251}    
  \be
  \raisebox{-9.5em}{\begin{tikzpicture}
  \scopeArrow{0.26}{>} \scopeArrow{0.52}{>} \scopeArrow{0.72}{>} \scopeArrow{0.85}{>}
  \filldraw[\colorCircle,fill=\colorSigma,line width=\locpc,postaction={decorate}] (0,0) circle (\locpa);
  \end{scope} \end{scope} \end{scope} \end{scope}
  \draw[\colorCircle,line width=\locpc,dash pattern=on 9pt off 5pt] (0,0) circle (\locpb);
  \draw[line width=\locpc,color=white,dashed] (-30:\locpa) arc (-30:35:\locpa);
  \scopeArrow{0.33}{\arrowDefect} \scopeArrow{0.83}{\arrowDefect}
  \draw[line width=\locpc,color=\colorTransp,postaction={decorate}] (\locpw:\locpA) -- (\locpw:\locpa);
  \end{scope} \end{scope}
  \scopeArrow{0.31}{\arrowDefectB} \scopeArrow{0.77}{\arrowDefectB}
  \draw[line width=\locpc,color=\colorTransp,postaction={decorate}] (\locpv:\locpA) -- (\locpv:\locpa);
  \draw[line width=\locpc,color=\colorTransp,postaction={decorate}] (\locpx:\locpA) -- (\locpx:\locpa);
  \draw[line width=\locpc,color=\colorTransp,postaction={decorate}] (\locpy:\locpA) -- (\locpy:\locpa);
  \end{scope} \end{scope}
  \scopeArrow{0.52}{\arrowDefectB}
  \draw[line width=\locpc,color=\colorTransp,postaction={decorate}]
             (\locpu:\locpA) node[right=-4pt,yshift=10pt] {$i$} -- (\locpu:\locpr*\locpb);
  \end{scope}
  \scopeArrow{0.33}{\arrowDefectB} \scopeArrow{0.83}{\arrowDefectB}
  \draw[line width=\locpc,color=\colorTransp,postaction={decorate}]
             (\locpz:\locpA) node[right=-6pt,yshift=-9pt] {$j$} -- (\locpz:\locps*\locpa);
  \end{scope} \end{scope}
  \scopeArrow{0.24}{>} \scopeArrow{0.38}{>} \scopeArrow{0.53}{>} \scopeArrow{0.666}{>}
  \filldraw[\colorCircle,fill=white,line width=\locpc,postaction={decorate}] (0,0) circle (\locpA);
  \end{scope} \end{scope} \end{scope} \end{scope}
  \draw[line width=\locpc,color=white,dashed] (-25:\locpA) arc (-25:30:\locpA);
  \filldraw[\colorCircle,fill=white,line width=\locpC]
	     (\locpu:\locpr*\locpb) circle (\locpt) node[right=3pt,color=black] {$\one$}
	     (\locpz:\locps*\locpa) circle (\locpt) node[right=3pt,color=black] {$\one$};
  \scopeArrow{0.54}{<}
  \draw[\colorCircle,line width=\locpf,postaction={decorate}] (\locpu:\locpr*\locpb) circle (\locpt);
  \draw[\colorCircle,line width=\locpf,postaction={decorate}] (\locpz:\locps*\locpa) circle (\locpt);
  \end{scope}
  \filldraw[color=\colorDefect] (\locpv:\locpa) circle (\locpd) (\locpv:\locpA) circle (\locpd)
             (\locpw:\locpa) circle (\locpd) (\locpw:\locpA) circle (\locpd)
             (\locpx:\locpa) circle (\locpd) (\locpx:\locpA) circle (\locpd)
             (\locpy:\locpa) circle (\locpd) (\locpy:\locpA) circle (\locpd)
             (\locpv:\locpb) circle (\locpe) (\locpw:\locpb) circle (\locpe)
             (\locpx:\locpb) circle (\locpe) (\locpy:\locpb) circle (\locpe)
             (\locpz:\locpb) circle (\locpe)
             (\locpu:\locpA) circle (\locpd) (\locpu:\locpr*\locpb-\locpt) circle (\locpD)
             (\locpz:\locpA) circle (\locpd) (\locpz:\locps*\locpa-\locpt) circle (\locpD);
  \node at (-\locpa-1.8,0) {$\Sigma\bs{ij} ~:=$};
  \end{tikzpicture}}
  \ee
and
  \be
  \raisebox{-9.3em}{\begin{tikzpicture}
  \scopeArrow{0.26}{>} \scopeArrow{0.52}{>} \scopeArrow{0.72}{>} \scopeArrow{0.85}{>}
  \filldraw[\colorCircle,fill=\colorSigma,line width=\locpc,postaction={decorate}] (0,0) circle (\locpa);
  \end{scope} \end{scope} \end{scope} \end{scope}
  \draw[\colorCircle,line width=\locpc,dash pattern=on 9pt off 5pt] (0,0) circle (\locpb);
  \draw[line width=\locpc,color=white,dashed] (-30:\locpa) arc (-30:35:\locpa);
  \scopeArrow{0.33}{\arrowDefect} \scopeArrow{0.83}{\arrowDefect}
  \draw[line width=\locpc,color=\colorTransp,postaction={decorate}] (\locpw:\locpA) -- (\locpw:\locpa);
  \end{scope} \end{scope}
  \scopeArrow{0.31}{\arrowDefectB} \scopeArrow{0.77}{\arrowDefectB}
  \draw[line width=\locpc,color=\colorTransp,postaction={decorate}] (\locpv:\locpA) -- (\locpv:\locpa);
  \draw[line width=\locpc,color=\colorTransp,postaction={decorate}] (\locpx:\locpA) -- (\locpx:\locpa);
  \draw[line width=\locpc,color=\colorTransp,postaction={decorate}] (\locpy:\locpA) -- (\locpy:\locpa);
  \end{scope} \end{scope}
  \scopeArrow{0.52}{\arrowDefectB}
  \draw[line width=\locpc,color=\colorTransp,postaction={decorate}]
             (\locpz:\locpA) node[right=-5pt,yshift=-8pt] {$j$} -- (\locpz:\locpr*\locpb);
  \end{scope}
  \scopeArrow{0.33}{\arrowDefectB} \scopeArrow{0.83}{\arrowDefectB}
  \draw[line width=\locpc,color=\colorTransp,postaction={decorate}]
             (\locpu:\locpA) node[right=-4pt,yshift=10pt] {$i$} -- (\locpu:\locps*\locpa);
  \end{scope} \end{scope}
  \scopeArrow{0.24}{>} \scopeArrow{0.38}{>} \scopeArrow{0.53}{>} \scopeArrow{0.666}{>}
  \filldraw[\colorCircle,fill=white,line width=\locpc,postaction={decorate}] (0,0) circle (\locpA);
  \end{scope} \end{scope} \end{scope} \end{scope}
  \draw[line width=\locpc,color=white,dashed] (-25:\locpA) arc (-25:30:\locpA);
  \filldraw[\colorCircle,fill=white,line width=\locpC]
	     (\locpz:\locpr*\locpb) circle (\locpt) node[right=3pt,color=black] {$\one$}
             (\locpu:\locps*\locpa) circle (\locpt) node[right=3pt,color=black] {$\one$};
  \scopeArrow{0.54}{<}
  \draw[\colorCircle,line width=\locpf,postaction={decorate}] (\locpz:\locpr*\locpb) circle (\locpt);
  \draw[\colorCircle,line width=\locpf,postaction={decorate}] (\locpu:\locps*\locpa) circle (\locpt);
  \end{scope}
  \filldraw[color=\colorDefect] (\locpv:\locpa) circle (\locpd) (\locpv:\locpA) circle (\locpd)
             (\locpw:\locpa) circle (\locpd) (\locpw:\locpA) circle (\locpd)
             (\locpx:\locpa) circle (\locpd) (\locpx:\locpA) circle (\locpd)
             (\locpy:\locpa) circle (\locpd) (\locpy:\locpA) circle (\locpd)
             (\locpv:\locpb) circle (\locpe) (\locpw:\locpb) circle (\locpe)
             (\locpx:\locpb) circle (\locpe) (\locpy:\locpb) circle (\locpe)
             (\locpu:\locpb) circle (\locpe)
             (\locpz:\locpA) circle (\locpd) (\locpz:\locpr*\locpb-\locpt) circle (\locpD)
             (\locpu:\locpA) circle (\locpd) (\locpu:\locps*\locpa-\locpt) circle (\locpD);
  \node at (-\locpa-1.8,0) {$\Sigma\bs{ji} ~:=$};
  \end{tikzpicture}}
  \ee
  ~\\[0.1em]
In these pictures, all but the most relevant labels are omitted, while also the circle 
along which the two bordisms are glued (including its defect points) is indicated as a
dashed circle.
\end{proof}


\subsection{Changing refinements} \label{sec:B2}

We now show that to a change of refinement from $(\Sigma;\Sigma\refi^{})$ to
$(\Sigma;\Sigma\refi')$ there is associated a canonical isomorphism between the respective 
\csfun s $\Zprime(\Sigma;\Sigma\refi^{})$ and $\Zprime(\Sigma;\Sigma\refi')$. This is achieved
in two steps: first we consider refinements of \filldisk s, and afterwards refinements
of arbitrary \defectsurface s. The following terminology will be convenient:

\begin{Definition}\label{def:filldiskrep}
By a \emph{\filldiskrep} $\DR_{\DD,\DD'}$ from $\DD$ to $\DD'$ we mean the operation of
replacing in a \defectsurface\ $\Sigma$ a \filldisk\ $\DD \,{\subset}\, \Sigma $
of some type $\XX$ by a \filldisk\ $\DD'$ of the same type with the same outer boundary.
\end{Definition}

Owing to the factorization result for fine \defectsurface s in  Theorem \ref{thm:T1T2=Ttr},
such a manipulation is completely under control by canonical isomorphisms.

We start by recalling that for each \filldisk\ $\DD \,{=}\, \DD_\XX$ there exists a disk
$\DDtad_{} \,{=}\, \DDtad_\XX$ with the same outer boundary as $\DD$ and with all inner 
boundaries being \tadpolecircle s (compare the example shown in \eqref{eq:exa:DDtadLL}).
We abbreviate by $\LL_\DD$ the gluing part of the outer boundary of $\DD$.

We want to construct an isomorphism 
  \be
  \varphi_{\DD}^{}(\Gamma) :\quad \ZZ(\DD)( \om(\DD))
  \tocong \ZZ(\DDtad)(\om(\DDtad))
  \label{eq:phiDD(Gamma)}
  \ee
of functors from $\TT(\LL_\DD)$ to \vect. The construction given below will
a priori depend on a combinatorial datum $\Gamma$ that is defined as follows: 
Considering the inner boundary circles of $\DD$ as (fattened) vertices, the set of inner
boundary circles together with the defect lines of $\DD$ and their end points on the outer
boundary circle or on generic defects of $\DD$ form a graph $\Gamma_{\!\DD}$. Without loss 
of generality
we assume that this graph is connected (otherwise the arguments below are to be applied to
every connected component separately, and the order in which this is done is irrelevant).
Then we select a subgraph $\Gamma \,{\subset}\, \Gamma_{\!\DD}$ that is a \emph{spanning tree}
in $\Gamma_{\!\DD}$, i.e.\ a rooted tree $\Gamma$ with a minimal number of edges
such that every vertex of $\Gamma_{\!\DD}$ is met by $\Gamma$.
(It is well known that a spanning tree exists for every graph.)
We can also take the root $v_0$ of $\Gamma$ to be lie on the outer boundary $\LL_\DD$.
As an example, the following picture shows such a spanning tree $\Gamma$ for the 
\transdisk\ $\DDt$ shown in \eqref{eq:exampleAtranspdisk}:
    \def\locpa  {2.7}    
    \def\locpA  {0.09}   
    \def\locpc  {2.2pt}  
    \def\locpd  {0.12}   
  \be
  \raisebox{-7.8em}{\begin{tikzpicture}[scale=1.2]
  \coordinate (incirc1) at (-55:0.66*\locpa);
  \coordinate (incirc2) at (3:0.71*\locpa);
  \coordinate (incirc3) at (65:0.55*\locpa);
  \coordinate (incirc4) at (145:0.66*\locpa);
  \coordinate (incirc5) at (225:0.35*\locpa);
  \coordinate (oucirc1) at (40:\locpa);
  \coordinate (oucirc2) at (90:\locpa);
  \coordinate (oucirc3) at (180:\locpa);
  \coordinate (oucirc4) at (225:\locpa);
  \coordinate (oucirc5) at (257:\locpa);
  \scopeArrow{0.19}{>} \scopeArrow{0.38}{>} \scopeArrow{0.57}{>} \scopeArrow{0.68}{>}
  \scopeArrow{0.94}{>}
  \filldraw[\colorCircle,fill=\colorSigma,line width=\locpc,postaction={decorate}] (0,0) circle (\locpa);
  \end{scope} \end{scope} \end{scope} \end{scope} \end{scope}
  \scopeArrow{0.49}{\arrowDefectB}
  \draw[line width=0.8*\locpc,color=\colorTransp,postaction={decorate}] (oucirc5) -- (incirc5);
  \end{scope}
  \scopeArrow{0.56}{\arrowDefect}
  \draw[line width=0.8*\locpc,color=\colorTransp,postaction={decorate}] (oucirc1) -- (incirc3);
  \draw[line width=0.8*\locpc,color=\colorTransp,postaction={decorate}] (oucirc2) -- (incirc3);
  \draw[line width=0.8*\locpc,color=\colorTransp,postaction={decorate}] (oucirc3) -- (incirc5);
  \end{scope}
  \scopeArrow{0.63}{\arrowDefect}
  \draw[line width=0.8*\locpc,color=\colorTransp,postaction={decorate}] (incirc5) -- (incirc1);
  \end{scope}
  \draw[line width=1.3*\locpc,color=\colorGraph] (incirc3) -- (incirc2) -- (incirc1)
            (oucirc4) -- node[above=-2pt,rotate=48] {$e_{01}$}
            (incirc5) -- node[left]{{\boldmath{$\Gamma$}}} (incirc3) -- (incirc4);
  \filldraw[double,\colorGraph,fill=\colorGraph,line width=0.65*\locpc] 
            (oucirc4) circle (\locpA) node[below=1pt,xshift=-4pt]{$v_0$}
            (incirc1) circle (\locpA) (incirc2) circle (\locpA) (incirc3) circle (\locpA)
            (incirc4) circle (\locpA) (incirc5) circle (\locpA);
  \filldraw[color=\colorDefect] (oucirc1) circle (\locpd) (oucirc2) circle (\locpd)
           (oucirc3) circle (\locpd) (oucirc5) circle (\locpd);
  \label{eq:spanningtree}
  \end{tikzpicture}}
  \ee ~\\[0.3em]
Here for clarity the vertices of $\Gamma$ are drawn as encircled points,
and also the remaining defect lines that do not give rise to edges of $\Gamma$ are indicated.

By the \emph{length} of a path in a graph we mean the number of its edges. Then the
\emph{depth} of a vertex $v$ of $\Gamma$ is defined as the length of the 
(unique) path from $v$ to the root of $\Gamma$; in particular the root has depth $0$. 

For any choice of spanning tree $\Gamma \,{\subset}\, \DD$, an isomorphism
\eqref{eq:phiDD(Gamma)} is obtained by the following prescription:
Apply the canonical isomorphism $\Xi$ from Lemma \ref{lem:canisoXi1} for every edge 
$e \,{\in}\, \DD {\setminus} \Gamma$; this results in a transparent disk
with two new tadpole vertices for each defect line not covered by $\Gamma$, which we
denote by $\widetilde\DD(\Gamma)$. Next apply the canonical isomorphisms $\rho\bs i$
from Lemma \ref{Lemma:rho} for every vertex and every edge of this disk
$\widetilde\DD(\Gamma)$ (in arbitrary order), whereby we end up with a \tadpoledisk\ 
$\DDtad(\Gamma)$. As an illustration, in the case of the spanning tree chosen in
\eqref{eq:spanningtree}, the disks $\widetilde\DD(\Gamma)$ and $\DDtad(\Gamma)$ look as follows:
    \def\locpa  {2.7}    
    \def\locpA  {0.26}   
    \def\locpc  {2.2pt}  
    \def\locpC  {1.8pt}  
    \def\locpd  {0.12}   
    \def\locpD  {0.08}   
    \def\locpf  {1.1pt}  
    \def\locpx  {2.2*\locpa,-1.4*\locpa}  
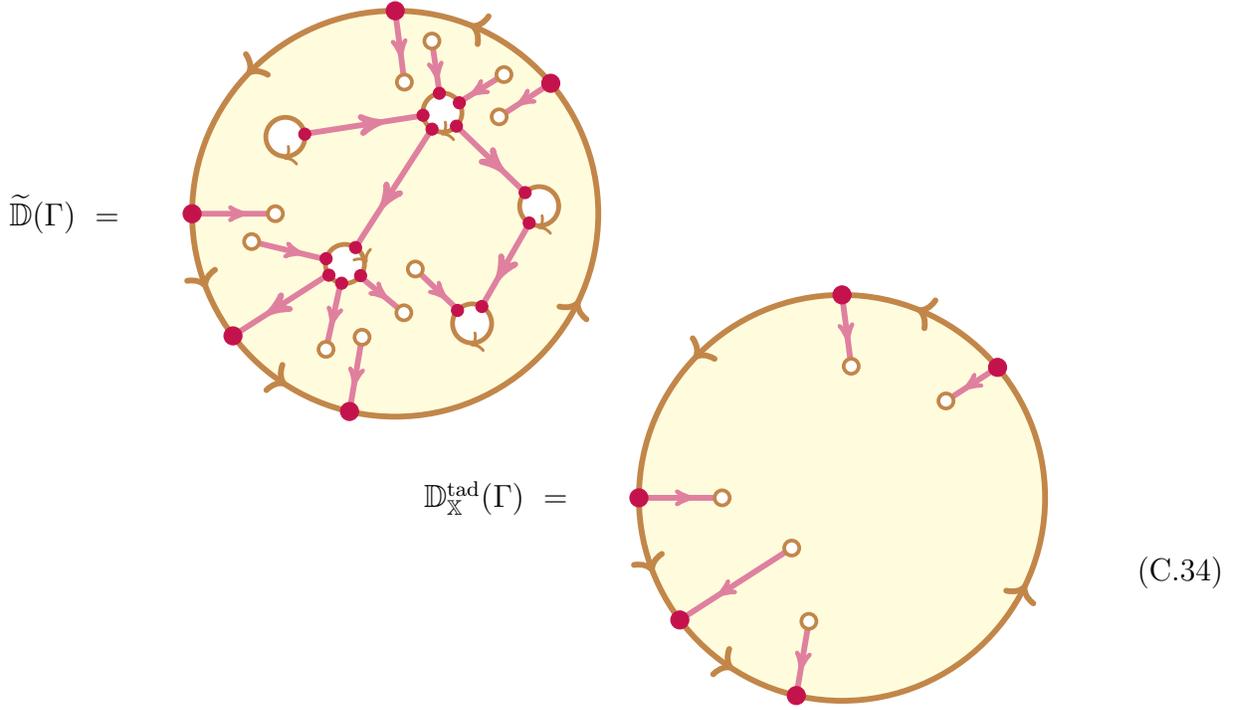
  \begin{eqnarray}
  \begin{tikzpicture}
  \coordinate (incirc1) at (-55:0.66*\locpa);
  \coordinate (Incirc1) at (-70:0.29*\locpa);
  \coordinate (incirc2) at (3:0.71*\locpa);
  \coordinate (incirc3) at (65:0.55*\locpa);
  \coordinate (Incirc3) at (52:0.87*\locpa);
  \coordinate (Incirk3) at (78:0.87*\locpa);
  \coordinate (incirc4) at (145:0.66*\locpa);
  \coordinate (incirc5) at (225:0.35*\locpa);
  \coordinate (Incirc5) at (191:0.72*\locpa);
  \coordinate (Incirk5) at (243:0.75*\locpa);
  \coordinate (Incirq5) at (275:0.49*\locpa);
  \coordinate (oucirc1) at (40:\locpa);
  \coordinate (Oucirc1) at (43:0.7*\locpa);
  \coordinate (oucirc2) at (90:\locpa);
  \coordinate (Oucirc2) at (86:0.65*\locpa);
  \coordinate (oucirc3) at (180:\locpa);
  \coordinate (Oucirc3) at (180:0.59*\locpa);
  \coordinate (oucirc4) at (217:\locpa);
  \coordinate (oucirc5) at (257:\locpa);
  \coordinate (Oucirc5) at (255:0.63*\locpa);
  \node at (-\locpa-1.5,0) {$\widetilde\DD(\Gamma)~=$};
  \scopeArrow{0.19}{>} \scopeArrow{0.38}{>} \scopeArrow{0.56}{>} \scopeArrow{0.66}{>}
  \scopeArrow{0.93}{>}
  \filldraw[\colorCircle,fill=\colorSigma,line width=\locpc,postaction={decorate}] (0,0) circle (\locpa);
  \end{scope} \end{scope} \end{scope} \end{scope} \end{scope}
  \scopeArrow{0.49}{\arrowDefectB}
  \draw[line width=\locpc,color=\colorTransp,postaction={decorate}] (oucirc4) -- (incirc5);
  \end{scope}
  \scopeArrow{0.63}{\arrowDefect}
  \draw[line width=\locpc,color=\colorTransp,postaction={decorate}] (incirc4) -- (incirc3);
  \draw[line width=\locpc,color=\colorTransp,postaction={decorate}] (incirc3) -- (incirc2);
  \draw[line width=\locpc,color=\colorTransp,postaction={decorate}] (incirc2) -- (incirc1);
  \draw[line width=\locpc,color=\colorTransp,postaction={decorate}] (incirc3) -- (incirc5);
  \end{scope}
  \scopeArrow{0.67}{\arrowDefect}
  \draw[line width=0.7*\locpc,color=\colorTransp,postaction={decorate}] (oucirc1) -- (Oucirc1);
  \draw[line width=0.7*\locpc,color=\colorTransp,postaction={decorate}] (oucirc2) -- (Oucirc2);
  \draw[line width=0.7*\locpc,color=\colorTransp,postaction={decorate}] (oucirc3) -- (Oucirc3);
  \draw[line width=0.7*\locpc,color=\colorTransp,postaction={decorate}] (Oucirc5) -- (oucirc5);
  \end{scope}
  \scopeArrow{0.56}{\arrowDefect}
  \draw[line width=0.7*\locpc,color=\colorTransp,postaction={decorate}] (Incirc1) -- (incirc1);
  \draw[line width=0.7*\locpc,color=\colorTransp,postaction={decorate}] (Incirc3) -- (incirc3);
  \draw[line width=0.7*\locpc,color=\colorTransp,postaction={decorate}] (Incirk3) -- (incirc3);
  \draw[line width=0.7*\locpc,color=\colorTransp,postaction={decorate}] (Incirc5) -- (incirc5);
  \end{scope}
  \scopeArrow{0.73}{\arrowDefect}
  \draw[line width=0.7*\locpc,color=\colorTransp,postaction={decorate}] (incirc5) -- (Incirk5);
  \draw[line width=0.7*\locpc,color=\colorTransp,postaction={decorate}] (incirc5) -- (Incirq5);
  \end{scope}
  \draw[line width=\locpc,color=\colorTransp] (oucirc1) -- (Oucirc1)
                     (oucirc2) -- (Oucirc2) (oucirc3) -- (Oucirc3) (oucirc5) -- (Oucirc5)
                     (Incirc1) -- (incirc1) (Incirc3) -- (incirc3) (Incirk3) -- (incirc3)
                     (Incirc5) -- (incirc5) (incirc5) -- (Incirk5) (incirc5) -- (Incirq5);
  \transv{Incirc1}; \transv{Incirc3}; \transv{Incirk3}; \transv{Incirc5}; \transv{Incirk5};
  \transv{Incirq5}; \transv{Oucirc1}; \transv{Oucirc2}; \transv{Oucirc3}; \transv{Oucirc5};
  \filldraw[\colorCircle,fill=white,line width=\locpC]
           (incirc1) circle (\locpA) (incirc2) circle (\locpA) (incirc3) circle (\locpA)
           (incirc4) circle (\locpA) (incirc5) circle (\locpA);
  \scopeArrow{0.82}{<}
  \draw[\colorCircle,line width=\locpf,postaction={decorate}] (incirc1) circle (\locpA);
  \draw[\colorCircle,line width=\locpf,postaction={decorate}] (incirc2) circle (\locpA);
  \draw[\colorCircle,line width=\locpf,postaction={decorate}] (incirc3) circle (\locpA);
  \draw[\colorCircle,line width=\locpf,postaction={decorate}] (incirc4) circle (\locpA);
  \end{scope}
  \scopeArrow{0.09}{<}
  \draw[\colorCircle,line width=\locpf,postaction={decorate}] (incirc5) circle (\locpA);
  \end{scope}
  \filldraw[color=\colorDefect] (oucirc1) circle (\locpd) (oucirc2) circle (\locpd)
           (oucirc3) circle (\locpd) (oucirc4) circle (\locpd) (oucirc5) circle (\locpd)
           (incirc1)+(60:\locpA) circle (\locpD) (incirc1)+(138.2:\locpA) circle (\locpD)
           (incirc2)+(136:\locpA) circle (\locpD) (incirc2)+(239.3:\locpA) circle (\locpD)
           (incirc3)+(30:\locpA) circle (\locpD) (incirc3)+(99:\locpA) circle (\locpD)
           (incirc3)+(188:\locpA) circle (\locpD) (incirc3)+(238:\locpA) circle (\locpD)
           (incirc3)+(316:\locpA) circle (\locpD) (incirc4)+(8:\locpA) circle (\locpD)
           (incirc5)+(57:\locpA) circle (\locpD) (incirc5)+(164.3:\locpA) circle (\locpD)
           (incirc5)+(215:\locpA) circle (\locpD) (incirc5)+(260.3:\locpA) circle (\locpD)
           (incirc5)+(323.2:\locpA) circle (\locpD);
  \begin{scope}[shift={(\locpx+0.5,0.4)}]   
  \coordinate (incirc5) at (225:0.35*\locpa);
  \coordinate (oucirc1) at (40:\locpa);
  \coordinate (Oucirc1) at (43:0.7*\locpa);
  \coordinate (oucirc2) at (90:\locpa);
  \coordinate (Oucirc2) at (86:0.65*\locpa);
  \coordinate (oucirc3) at (180:\locpa);
  \coordinate (Oucirc3) at (180:0.59*\locpa);
  \coordinate (oucirc4) at (217:\locpa);
  \coordinate (oucirc5) at (257:\locpa);
  \coordinate (Oucirc5) at (255:0.63*\locpa);
  \node at (-\locpa-1.6,0) {$\DDtad_{\XX}(\Gamma)~=$}; 
  \scopeArrow{0.19}{>} \scopeArrow{0.38}{>} \scopeArrow{0.56}{>} \scopeArrow{0.66}{>}
  \scopeArrow{0.93}{>}
  \filldraw[\colorCircle,fill=\colorSigma,line width=\locpc,postaction={decorate}] (0,0) circle (\locpa);
  \end{scope} \end{scope} \end{scope} \end{scope} \end{scope}
  \draw[line width=\locpc,color=\colorTransp] (oucirc1) -- (Oucirc1) (oucirc2) -- (Oucirc2)
                       (oucirc3) -- (Oucirc3) (incirc5) -- (oucirc4) (oucirc5) -- (Oucirc5);
  \scopeArrow{0.67}{\arrowDefect}
  \draw[line width=0.7*\locpc,color=\colorTransp,postaction={decorate}] (oucirc1) -- (Oucirc1);
  \draw[line width=0.7*\locpc,color=\colorTransp,postaction={decorate}] (oucirc2) -- (Oucirc2);
  \draw[line width=0.7*\locpc,color=\colorTransp,postaction={decorate}] (oucirc3) -- (Oucirc3);
  \draw[line width=0.7*\locpc,color=\colorTransp,postaction={decorate}] (incirc5) -- (oucirc4);
  \draw[line width=0.7*\locpc,color=\colorTransp,postaction={decorate}] (Oucirc5) -- (oucirc5);
  \end{scope}
  \filldraw[color=\colorDefect] (oucirc1) circle (\locpd) (oucirc2) circle (\locpd)
              (oucirc3) circle (\locpd) (oucirc4) circle (\locpd) (oucirc5) circle (\locpd);
  \transv{Oucirc1}; \transv{Oucirc2}; \transv{Oucirc3}; \transv{incirc5}; \transv{Oucirc5};
  \end{scope}
  \end{tikzpicture}
  \nonumber \\[-2.4em]~ \\[-0.8em] \nonumber~
  \end{eqnarray}

Altogether this defines canonically an isomorphism $\varphi_{\DD}(\Gamma)$ of the form 
\eqref{eq:phiDD(Gamma)}. Next we show that
this isomorphism does in fact not depend on the choices made in its construction:

\begin{Lemma}\label{lem:varphiDD-indep}
\Itemizeiii

\item[{\rm (i)}]
The isomorphism $\varphi_{\DD}(\Gamma) \colon \ZZ(\DD)(- \boti \om(\DD))
\,{\tocong}\, \ZZ(\DDtad)(- \boti \om(\DDtad))$
does not  depend on the order in which the isomorphisms $\rho\bs i$ are applied. 

\item[{\rm (ii)}]
Let $\Gamma$ and $\Gamma'$ be two spanning trees for $\DD$. Then
$\varphi_{\DD}(\Gamma) \,{=}\, \varphi_{\DD}(\Gamma')$.
\end{itemize}
\end{Lemma}

\begin{proof}
(i)\, Obviously, any two isomorphisms $\rho\bs i$ commute if they are applied on two
different vertices. If they are applied on one and the same vertex, the statement follows
directly from Lemma \ref{Lemma:rho}.
\\[.2em]
(ii)\, As above we assume without loss of generality that the graph $\GammaD$ on $\DD$ that
is formed by the defect lines is connected and fix a spanning tree $\Gamma$ for $\GammaD$,
with root vertex $v_0$. Denote by $E_0$ the set of all edges of $\GammaD$ that have one
of their ends on the outer boundary $\LL_\DD$. By construction, exactly one edge
$e_{01} \,{\in}\, E_0$ (as indicated in the picture \eqref{eq:spanningtree}) belongs to
the spanning tree $\Gamma$. Removing $e_{01}$ from $\Gamma$ and replacing it by any other 
edge $e_{0i} \,{\in}\, E_0$ gives another spanning tree, with different root $v_{0i}$, which we 
denote by $\Gamma_{\!0i}$. For instance, the following spanning tree $\Gamma_{\!0i}$ for the
transparent disk \eqref{eq:exampletranspdisk} arises this way from the spanning tree shown in 
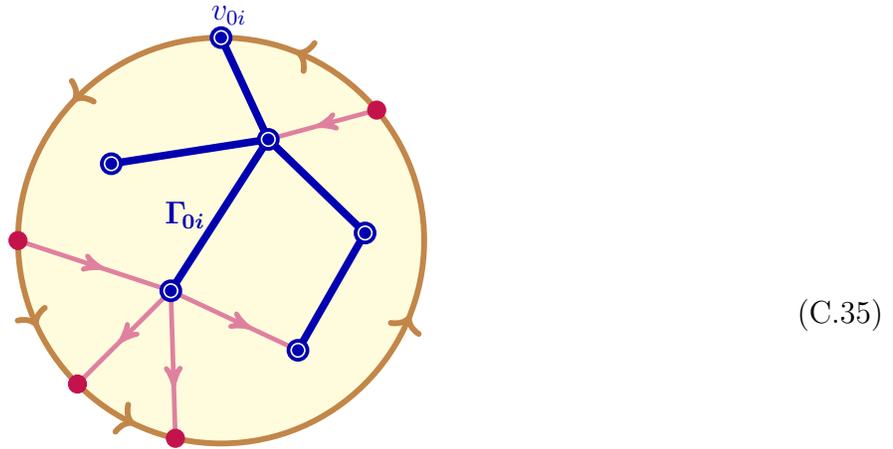
\eqref{eq:spanningtree}:
 \\[-1.5em]~
    \def\locpa  {2.7}    
    \def\locpA  {0.09}   
    \def\locpc  {2.2pt}  
    \def\locpd  {0.12}   
  \be
  \raisebox{-6.1em}{\begin{tikzpicture}[scale=0.9]
  \coordinate (incirc1) at (-55:0.66*\locpa);
  \coordinate (incirc2) at (3:0.71*\locpa);
  \coordinate (incirc3) at (65:0.55*\locpa);
  \coordinate (incirc4) at (145:0.66*\locpa);
  \coordinate (incirc5) at (225:0.35*\locpa);
  \coordinate (oucirc1) at (40:\locpa);
  \coordinate (oucirc2) at (90:\locpa);
  \coordinate (oucirc3) at (180:\locpa);
  \coordinate (oucirc4) at (225:\locpa);
  \coordinate (oucirc5) at (257:\locpa);
  \scopeArrow{0.19}{>} \scopeArrow{0.38}{>} \scopeArrow{0.57}{>} \scopeArrow{0.68}{>}
  \scopeArrow{0.94}{>}
  \filldraw[\colorCircle,fill=\colorSigma,line width=\locpc,postaction={decorate}] (0,0) circle (\locpa);
  \end{scope} \end{scope} \end{scope} \end{scope} \end{scope}
  \scopeArrow{0.49}{\arrowDefectB}
  \draw[line width=0.8*\locpc,color=\colorTransp,postaction={decorate}] (oucirc5) -- (incirc5);
  \end{scope}
  \scopeArrow{0.56}{\arrowDefect}
  \draw[line width=0.8*\locpc,color=\colorTransp,postaction={decorate}] (oucirc1) -- (incirc3);
  \draw[line width=0.8*\locpc,color=\colorTransp,postaction={decorate}] (oucirc3) -- (incirc5);
  \draw[line width=0.8*\locpc,color=\colorTransp,postaction={decorate}] (incirc5) -- (oucirc4);
  \end{scope}
  \scopeArrow{0.63}{\arrowDefect}
  \draw[line width=0.8*\locpc,color=\colorTransp,postaction={decorate}] (incirc5) -- (incirc1);
  \end{scope}
  \draw[line width=1.3*\locpc,color=\colorGraph] (incirc4) -- (incirc3) -- (incirc2) -- (incirc1)
	    (incirc5) -- node[left]{{\boldmath{$\Gamma_{\!0i}$}}} (incirc3) -- (oucirc2);
  \filldraw[double,\colorGraph,fill=\colorGraph,line width=0.65*\locpc] 
	    (oucirc2) circle (\locpA) node[above=0.5pt,xshift=3pt]{$v_{0i}$}
            (incirc1) circle (\locpA) (incirc2) circle (\locpA) (incirc3) circle (\locpA)
            (incirc4) circle (\locpA) (incirc5) circle (\locpA);
  \filldraw[color=\colorDefect] (oucirc1) circle (\locpd) (oucirc3) circle (\locpd)
           (oucirc4) circle (\locpd) (oucirc5) circle (\locpd);
  \end{tikzpicture}}
  \ee ~\\[0.2em]
We are now going to show that $\varphi_{\DD}(\Gamma) \,{=}\, \varphi_{\DD}(\Gamma_{\!0i})$. 
It is enough to assume that $E_0$ has precisely two elements. In this case the statement is
implied by the following result:

\begin{Lemma} \label{lem:Xirho=rhoXi}
The natural isomorphisms between block functors that are indicated in the following picture
commute, for any choice of orientations of the $($suppressed$)$ transparent defect 
lines in the \filldisk\ that is present in the two upper rows of the picture:
    \def\locpa  {1.9}    
    \def\locpA  {0.52}   
    \def\locpb  {0.28}   
    \def\locpc  {2.2pt}  
    \def\locpd  {0.12}   
    \def\locpf  {0.88}   
    \def\locpt  {0.34}   
    \def\locpx  {230}    
    \def\locpy  {344}    
    \def\locpz  {107}    
  \be
  \raisebox{-4.2em}
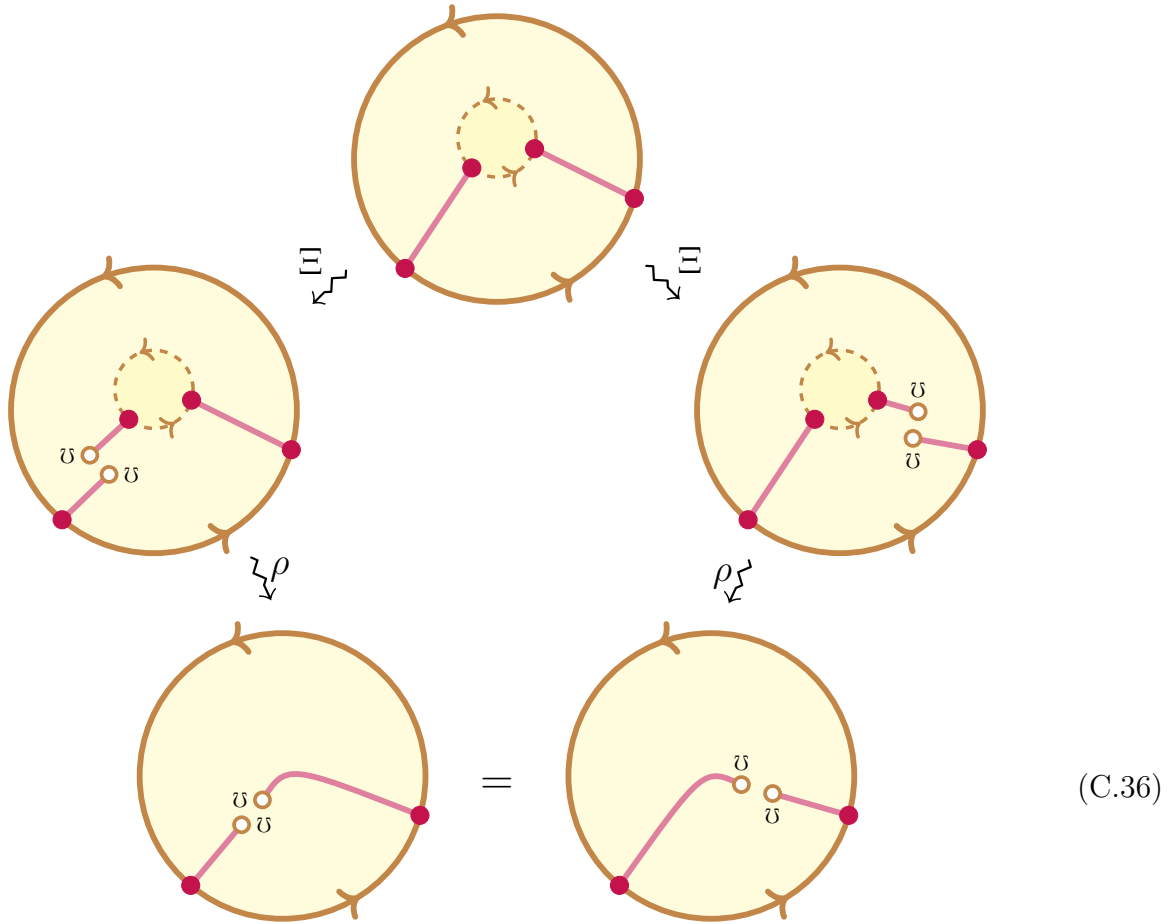
{\begin{tikzpicture}
  \begin{scope}[shift={(0,4*\locpf*\locpa)}]
  \scopeArrow{0.31}{>} \scopeArrow{0.84}{>}
  \filldraw[\colorCircle,fill=\colorSigma,line width=\locpc,postaction={decorate}]
                  (0,0) circle (\locpa);
  \filldraw[dashed,\colorCircle,fill=\colorSigmaDark,line width=0.6*\locpc,postaction={decorate}]
                  (0,\locpb) circle (\locpA);
  \end{scope} \end{scope}
  \draw[line width=\locpc,color=\colorTransp] (\locpx:\locpA)+(0,\locpb) -- (\locpx:\locpa);
  \draw[line width=\locpc,color=\colorTransp] (\locpy:\locpA)+(0,\locpb) -- (\locpy:\locpa);
  \filldraw[color=\colorDefect] 
                  (\locpx:\locpa) circle (\locpd) (\locpx:\locpA)+(0,\locpb) circle (\locpd)
                  (\locpy:\locpa) circle (\locpd) (\locpy:\locpA)+(0,\locpb) circle (\locpd);
  \end{scope} 
  \begin{scope}[shift={(-2.4*\locpa,2*\locpf*\locpa)}]
  \scopeArrow{0.31}{>} \scopeArrow{0.84}{>}
  \filldraw[\colorCircle,fill=\colorSigma,line width=\locpc,postaction={decorate}]
                  (0,0) circle (\locpa);
  \filldraw[dashed,\colorCircle,fill=\colorSigmaDark,line width=0.6*\locpc,postaction={decorate}]
                  (0,\locpb) circle (\locpA);
  \end{scope} \end{scope}
  \draw[line width=\locpc,color=\colorTransp] (\locpx:\locpA)+(0,\locpb) -- (\locpx-15:2*\locpA);
  \draw[line width=\locpc,color=\colorTransp,postaction={decorate}]
                  (\locpx+5:2*\locpA) -- (\locpx:\locpa);
  \draw[line width=\locpc,color=\colorTransp] (\locpy:\locpA)+(0,\locpb) -- (\locpy:\locpa);
  \transv{\locpx-15:2*\locpA} node[left=1pt,color=black] {$\scriptstyle\om$};
  \transv{\locpx+5:2*\locpA} node[right=1pt,color=black] {$\scriptstyle\om$};
  \filldraw[color=\colorDefect] 
                  (\locpx:\locpa) circle (\locpd) (\locpx:\locpA)+(0,\locpb) circle (\locpd)
                  (\locpy:\locpa) circle (\locpd) (\locpy:\locpA)+(0,\locpb) circle (\locpd);
  \end{scope} 
  \begin{scope}[shift={(2.4*\locpa,2*\locpf*\locpa)}]
  \scopeArrow{0.31}{>} \scopeArrow{0.84}{>}
  \filldraw[\colorCircle,fill=\colorSigma,line width=\locpc,postaction={decorate}]
                  (0,0) circle (\locpa);
  \filldraw[dashed,\colorCircle,fill=\colorSigmaDark,line width=0.6*\locpc,postaction={decorate}]
                  (0,\locpb) circle (\locpA);
  \end{scope} \end{scope}
  \draw[line width=\locpc,color=\colorTransp] (\locpy:\locpA)+(0,\locpb) -- (\locpy+15:2*\locpA);
  \draw[line width=\locpc,color=\colorTransp] (\locpy-5:2*\locpA) -- (\locpy:\locpa);
  \draw[line width=\locpc,color=\colorTransp] (\locpx:\locpA)+(0,\locpb) -- (\locpx:\locpa);
  \transv{\locpy+15:2*\locpA} node[above=1pt,color=black] {$\scriptstyle\om$};
  \transv{\locpy-5:2*\locpA} node[below=1pt,color=black] {$\scriptstyle\om$};
  \filldraw[color=\colorDefect] 
                  (\locpx:\locpa) circle (\locpd) (\locpx:\locpA)+(0,\locpb) circle (\locpd)
                  (\locpy:\locpa) circle (\locpd) (\locpy:\locpA)+(0,\locpb) circle (\locpd);
  \end{scope} 
  \begin{scope}[shift={(-1.5*\locpa,-0.8*\locpa)}]
  \scopeArrow{0.31}{>} \scopeArrow{0.84}{>}
  \filldraw[\colorCircle,fill=\colorSigma,line width=\locpc,postaction={decorate}]
                  (0,0) circle (\locpa);
  \end{scope} \end{scope}
  \draw[line width=\locpc,color=\colorTransp] 
                  (\locpy:\locpa) .. controls (0.1*\locpA,0.33*\locpA) .. (\locpx:0.8*\locpA);
  \draw[line width=\locpc,color=\colorTransp] (\locpx:1.63*\locpA) -- (\locpx:\locpa);
  \transv{\locpx:0.8*\locpA} node[left=1pt,color=black] {$\scriptstyle\om$};
  \transv{\locpx:1.63*\locpA} node[right=1pt,color=black] {$\scriptstyle\om$};
  \filldraw[color=\colorDefect] (\locpx:\locpa) circle (\locpd) (\locpy:\locpa) circle (\locpd);
  \end{scope} 
  \begin{scope}[shift={(1.5*\locpa,-0.8*\locpa)}]
  \scopeArrow{0.31}{>} \scopeArrow{0.84}{>}
  \filldraw[\colorCircle,fill=\colorSigma,line width=\locpc,postaction={decorate}]
                  (0,0) circle (\locpa);
  \end{scope} \end{scope}
  \draw[line width=\locpc,color=\colorTransp]
                  (\locpx:\locpa) .. controls (-0.2*\locpA,0.23*\locpA) .. (\locpy:0.8*\locpA);
  \draw[line width=\locpc,color=\colorTransp]
                  (\locpy:1.63*\locpA) -- (\locpy:\locpa);
  \transv{\locpy:0.8*\locpA} node[above=1pt,color=black] {$\scriptstyle\om$};
  \transv{\locpy:1.63*\locpA} node[below=1pt,color=black] {$\scriptstyle\om$};
  \filldraw[color=\colorDefect] (\locpx:\locpa) circle (\locpd) (\locpy:\locpa) circle (\locpd);
  \end{scope} 
  \node at (0,-0.87*\locpa) {\Large $=$};
  \node at (-1.35*\locpf*\locpa,3*\locpf*\locpa)[rotate=-135] {{\LARGE$\rightsquigarrow$}};
  \node at (-0.22-1.35*\locpf*\locpa,0.32+3*\locpf*\locpa) {{\large$\Xi$}};
  \node at (1.3*\locpf*\locpa,3*\locpf*\locpa)[rotate=-45] {{\LARGE$\rightsquigarrow$}};
  \node at (0.39+1.3*\locpf*\locpa,0.33+3*\locpf*\locpa) {{\large$\Xi$}};
  \node at (-1.9*\locpf*\locpa,0.65*\locpf*\locpa)[rotate=-65] {{\LARGE$\rightsquigarrow$}};
  \node at (-1.25-0.98*\locpf*\locpa,-0.45+\locpf*\locpa) {{\large$\rho$}};
  \node at (1.9*\locpf*\locpa,0.65*\locpf*\locpa)[rotate=245] {{\LARGE$\rightsquigarrow$}};
  \node at (1.26+1.04*\locpf*\locpa,-0.54+\locpf*\locpa) {{\large$\rho$}};
  \end{tikzpicture}}
  \label{eq:Xirho=rhoXi}
  \ee
  ~\\[-0.1em]
$($For concreteness, the picture shows the case that the outer boundary $\LL_\DD$ is a gluing
circle, but the statement applies to \filldisk s of arbitrary type $\XX$.$)$

\end{Lemma}

\begin{proof}
The functors $\rho$ and $\Xi$ are both defined using the adjoint equivalence between the
functors $G_\XX$ and $\widetilde G_\XX$ that follow from Lemma \ref{lem:canisoXi1}. If we
factorize the block functors into corresponding composites of $G_\XX$ and $\widetilde G_\XX$,
then the statement reduces to the zigzag identity for this adjoint equivalence. 
\end{proof}
       
\noindent
We continue the proof of Lemma \ref{lem:varphiDD-indep} by induction on the depth of the 
vertices of $\Gamma$. Consider the edges
of $\Gamma$ from the single depth-1 vertex $v_1$ to the depth-$2$ vertices. Pick a
slightly smaller disk $\DD_1 \,{\subset}\, \DD$ that does not contain the vertices
$v_0$ and $v_1$, but contains all other vertices of $\Gamma$ of depth larger than 1,
as indicated in the picture
    \def\locpa  {2.7}    
    \def\locpA  {0.09}   
    \def\locpc  {2.2pt}  
    \def\locpd  {0.12}   
  \be
  \raisebox{-6.0em}
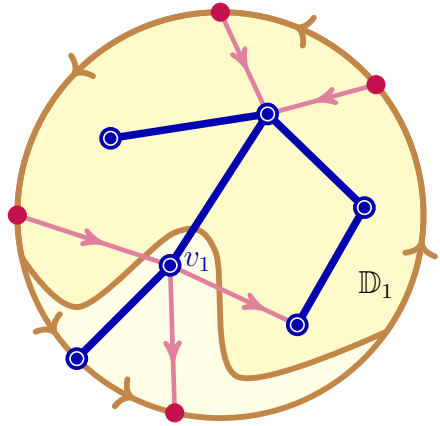
{\begin{tikzpicture}
  \coordinate (incirc1) at (-55:0.66*\locpa);
  \coordinate (incirc2) at (3:0.71*\locpa);
  \coordinate (incirc3) at (65:0.55*\locpa);
  \coordinate (incirc4) at (145:0.66*\locpa);
  \coordinate (incirc5) at (225:0.35*\locpa);
  \coordinate (oucirc1) at (40:\locpa);
  \coordinate (oucirc2) at (90:\locpa);
  \coordinate (oucirc3) at (180:\locpa);
  \coordinate (oucirc4) at (225:\locpa);
  \coordinate (oucirc5) at (257:\locpa);
  \fill[\colorSigmaLight] (0,0) circle (\locpa);
  \filldraw[line width=\locpc,draw=\colorCircle,fill=\colorSigmaDark] 
                  (-35:\locpa) arc (-35:192:\locpa) plot[smooth,tension=0.55]
                  coordinates{(192:\locpa) (213:0.83*\locpa) (201:0.24*\locpa)
                  (262:0.18*\locpa) (276:0.8*\locpa) (325:\locpa)};
  \scopeArrow{0.19}{>} \scopeArrow{0.38}{>} \scopeArrow{0.60}{>} \scopeArrow{0.68}{>}
  \scopeArrow{0.98}{>}
  \draw[\colorCircle,line width=\locpc,postaction={decorate}] (0,0) circle (\locpa);
  \end{scope} \end{scope} \end{scope} \end{scope} \end{scope}
  \scopeArrow{0.49}{\arrowDefectB}
  \draw[line width=0.8*\locpc,color=\colorTransp,postaction={decorate}] (oucirc5) -- (incirc5);
  \end{scope}
  \scopeArrow{0.56}{\arrowDefect}
  \draw[line width=0.8*\locpc,color=\colorTransp,postaction={decorate}] (oucirc1) -- (incirc3);
  \draw[line width=0.8*\locpc,color=\colorTransp,postaction={decorate}] (oucirc2) -- (incirc3);
  \draw[line width=0.8*\locpc,color=\colorTransp,postaction={decorate}] (oucirc3) -- (incirc5);
  \end{scope}
  \scopeArrow{0.74}{\arrowDefect}
  \draw[line width=0.8*\locpc,color=\colorTransp,postaction={decorate}] (incirc5) -- (incirc1);
  \end{scope}
  \draw[line width=1.3*\locpc,color=\colorGraph] (incirc3) -- (incirc2) -- (incirc1)
            (oucirc4) -- (incirc5) -- (incirc3) -- (incirc4);
  \filldraw[double,\colorGraph,fill=\colorGraph,line width=0.65*\locpc] 
            (oucirc4) circle (\locpA)
            (incirc1) circle (\locpA) (incirc2) circle (\locpA) (incirc3) circle (\locpA)
            (incirc4) circle (\locpA) (incirc5) node[right=0.5pt,yshift=1.5pt]{$v_1$} circle (\locpA);
  \filldraw[color=\colorDefect] (oucirc1) circle (\locpd) (oucirc2) circle (\locpd)
           (oucirc3) circle (\locpd) (oucirc5) circle (\locpd);
  \node at (336:0.84*\locpa) {$\DD_1$};
  \end{tikzpicture}}
  \ee ~\\[0.2em]
Define the graph $\Gamma_{\!1}$ as the graph obtained by erasing from 
$\Gamma \,{\cap}\, \DD_1$ the edges $E_{0}$. This graph has, in general, 
several components. In the sequel we assume for simplicity that $\Gamma_{\!1}$
is connected -- if it is not, then each of its (finitely many) components
is to be treated analogously. With this assumption, $\Gamma_{\!1}$
furnishes a spanning tree for the disk $\DD_1$, with root
$v_1'$ at the intersection of $\DD_1$ and the edge of $\Gamma$ that
connects $v_1$ with the (by the assumption just made, unique) depth-$2$ vertex.
Repeating the previous argument we see that the corresponding isomorphism $\varphi_{\DD_1}$ 
remains unchanged if we replace the edge containing $v_1'$ on the spanning tree
$\Gamma_{\!1}$ by a different edge. By iterating this process we can reach any spanning tree
$\Gamma'$. We can thus conclude that $\varphi_{\DD}(\Gamma) \,{=}\, \varphi_{\DD}(\Gamma')$
for all spanning trees $\Gamma$ and $\Gamma'$ for $\DD$.
\end{proof}

In view of this result from now on we just write $\varphi_{\DD}$ for the isomorphism
$\varphi_{\DD}(\Gamma)$, for any choice of spanning tree $\Gamma$. 
Next we observe that our construction is local, in the following sense:

\begin{Proposition} \label{Proposition:main-disk-repl}
Let $\DD$ and $\DD'$ be two \filldisk s of the same type. There is a distinguished family 
  \be
  \varphi_{\DD,\DD'}^{} :\quad
  \ZZ(\DD)(- \boti \om(\DD)) \xrightarrow{~~} \ZZ(\DD')(- \boti \om(\DD))
  \label{eq:nat-iso-phi-transp}
  \ee
of natural isomorphisms, one for each \filldiskrep\ $\DR_{\DD,\DD'}$ with
the following properties:
\Itemizeiii

\item[{\rm 1.}] 
\label{item:coherence} \emph{(Coherence):} 
For any triple $\DD$, $\DD'$, $\DD''$ of \filldisk s of the same type the vertical 
com\-po\-sition 
of the natural transformations $\varphi_{\DD,\DD'}^{}$ and $\varphi_{\DD',\DD''}^{}$ 
is given by
  \be
  \varphi_{\DD',\DD''}^{} \veco \varphi_{\DD,\DD'}^{} = \varphi_{\DD,\DD''}^{} \,.
  \ee

\item[{\rm 2.}] 
\label{item:inside-disk-repl} \emph{(Factorization):}
Given a \filldisk\ that has
the form $\DD \,{=}\, \mathbb Y \,{\circ}\, (\DD_1 \,{\sqcup} \cdots
{\sqcup}\, \DD_n)$, with $\DD_1,\DD_2, ...\,, \DD_n$ non-intersecting \filldisk s in $\DD$
(of types $\XX_i$) and $\mathbb Y$ the \defectsurface\ that results from removing all the
disks $\DD_i$, $i \,{=}\, 1,2,...\,,n$, from $\DD$, we have, for any $n$-tuple of \filldiskrep s
$\big( \DR_{\DD_i^{},\DD_i'} \big)_{i=1,...,n}$ that do not change the outer boundaries
$\partial_\text{outer}\DD_i$, the equality
  \be
  \varphi_{\DD,\DD'}^{} = (\varphi_{\DD_1^{},\DD_1'}^{} \boti \cdots \boti
  \varphi_{\DD_n^{},\DD_n'}^{}) \hoco \ZZ(\mathbb Y)
  \label{eq:fact-varphi}
  \ee
of natural transformations, where $\DD'$ is the \defectsurface\
$\DD' \,{=}\, \mathbb Y \,{\circ}\, (\DD_1'\,{\sqcup}\cdots{\sqcup}\,\DD_n')$
and `$\,\hoco$' is the horizontal composition of natural transformations.
\end{itemize}
\end{Proposition}

\begin{proof}                   
For any pair $\DD$ and $\DD'$ of \filldisk s of the same type $\XX$, we define the isomorphism 
$\varphi_{\DD,\DD'}$ by $\varphi_{\DD,\DD'} \,{:=}\, \varphi_{\DD'}^{-1}\,{\circ}\,\varphi_{\DD}^{}$,
with $ \varphi_{\DD} \colon \ZZ(\DD)(- \boti \om(\DD)) \,{\to}\, \ZZ(\DDtad)(- \boti \om(\DDtad)$
the isomorphism constructed above. It follows directly from this definition that
$\varphi_{\DD,\DD'}$ satisfies coherence. 
To establish factorization, we observe that there is a spanning tree $\Gamma'$ for 
$\mathbb Y$ in $\DD \,{=}\, \mathbb Y \,{\circ}\, (\DD_1 \,{\sqcup} \cdots {\sqcup}\, \DD_n)$
that has exactly one vertex on each boundary component of $\mathbb Y$. We can complete 
this graph $\Gamma'$ to a spanning tree $\Gamma$ of $\DD$ in such a way that
$\Gamma_{\!i} \,{:=}\, \Gamma \,{\cap}\, \DD_i$ is a spanning
tree for $\DD_i$ for every $i \,{\in}\, \{1,2,...\,,n\}$. Since the order in which we apply
the isomorphisms  $\rho$ in the definition of $\varphi_{\DD}$ is irrelevant, we readily see 
that the equality \eqref{eq:fact-varphi} indeed holds.
\end{proof}

Now recall the notion of a \filldiskrep\ $\DR_{\DD,\DD'}$ inside a \defectsurface\ $\Sigma$ 
(which is e.g.\ implicit in the factorization property \eqref{eq:fact-varphi}).
We denote the resulting \defectsurface\ by $\DR(\Sigma) \,{\equiv}\, \DR_{\DD,\DD'}(\Sigma)$.
Let $\Sigma$ be an arbitrary \defectsurface\ and $(\Sigma;\Sigma_{\rm ref_1})$ and 
$(\Sigma;\Sigma_{\rm ref_2})$ be any two refinements of $\Sigma$. 

\begin{Definition}
\Itemizeiii
\item[{\rm (i)}]
A \emph{refinement replacement} from $(\Sigma;\Sigma_{\rm ref_1})$ to
$(\Sigma;\Sigma_{\rm ref_2})$ is a sequence of
$($possibly intersecting$)$ \filldiskrep s $(\DR_1,\DR_2,...\,,\DR_n)$ such that
  \be
  \DR_n( \cdots \DR_1(\Sigma_{\rm ref_1}) \cdots) = \Sigma_{\rm ref_2} \,.
  \ee

\item[{\rm (ii)}]
We call two refinement replacements $(\DR_1,...\,,\DR_n)$ and $(\DR_1',...\,,\DR_{n'}')$
from $(\Sigma;\Sigma_{\rm ref_1})$ to $(\Sigma;\Sigma_{\rm ref_2})$ \emph{equivalent} iff 
the induced natural isomorphisms agree. 
\end{itemize}
\end{Definition}

According to Proposition \ref{Proposition:main-disk-repl}, any \filldiskrep\ $\DR$
in $\Sigma_{\rm ref_1}$ provides us with an isomorphism $\varphi_{\Sigma,\DR(\Sigma)}\colon 
\ZZ(\Sigma_{\rm ref_1})(-\boti\om_1^{}) \,{\to}\, \ZZ(\DR(\Sigma_{\rm ref_1})))(-\boti\om_1')$, 
with $\om_1^{}$ and $\om_1'$ the \mute\ objects for the respective \filldisk s involved.
Hence a refinement replacement $(\DR_1,...\,,\DR_n)$ from $(\Sigma;\Sigma_{\rm ref_1})$ to 
$(\Sigma;\Sigma_{\rm ref_2})$ gives an isomorphism 
  \be
  \varphi^{}_{ \Sigma_{\rm ref_1}, \DR_n( \cdots\, \DR_1(\Sigma_{\rm ref_1}) \,\cdots) } :\quad
  \ZZ(\Sigma_{\rm ref_1})(-\boti\om_1) \xrightarrow{~~} \ZZ(\Sigma_{\rm ref_2})(-\boti\om_2) \,.
  \label{eq:ref-repl-M}
  \ee
As we will see in Lemma \ref{lem:refrepl} below, a refinement replacement 
exists between any two fine refinements that refine a given \defectsurface. 

To proceed we introduce the notion of \emph{common subrefinement}. Let $(\Sigma;\Sigma_1)$
and $(\Sigma;\Sigma_2)$ be refinements that refine the same \defectsurface\ $\Sigma$. Then the 
common subrefinement $(\Sigma;\Sigma_{1,2})$ of $(\Sigma;\Sigma_1)$ and $(\Sigma;\Sigma_2)$
is constructed by combining all transparent defects from $\Sigma_1$ and from $\Sigma_2$
in the following manner: First take the collection of all transparent defects $\delta_2$
of $\Sigma_2$ that are not part of $\Sigma_1$.
We can use the embedding of the defects $\delta_2$ in $\Sigma_2$ to embed $\delta_2$ in the
surface $\Sigma_1$ in such a way that any resulting intersections of transparent defects are
generic (if necessary, deform the defects slightly to achieve this, see Remark 
\ref{refi-vs-transdisk} (iv)). Denote the so obtained surface with defects by $\mathring\Sigma_{1,2}$.

The following prescription makes $\mathring\Sigma_{1,2}$ into a \defectsurface\ $\Sigma_{1,2}$
endowed with a vector field that (just like the representatives of the
framings $\chi_1$ on $\Sigma_1$ and $\chi_2$ on $\Sigma_2$) is homotopic to the one of $\Sigma$:
Consider a tubular neighborhood $N_2$ of all defects $\delta_2$ in $\Sigma_2$, and take,
for each intersection $v \,{\in}\, \mathring\Sigma_{1,2}$ of (the images of) a defect $\delta_2$ 
of $\Sigma_2$ with a defect $\delta_1$ of $\Sigma_1$, a small circle $S_v$ around $v$
that intersects $\delta_1$ outside the image of $N_2$, as indicated in
    \def\locpa  {2.7}    
    \def\locpc  {2.2pt}  
    \def\locpd  {0.12}   
    \def\locpx  {4.4}    
    \def\locpy  {3.2}    
  \be
  \raisebox{-5.0em}
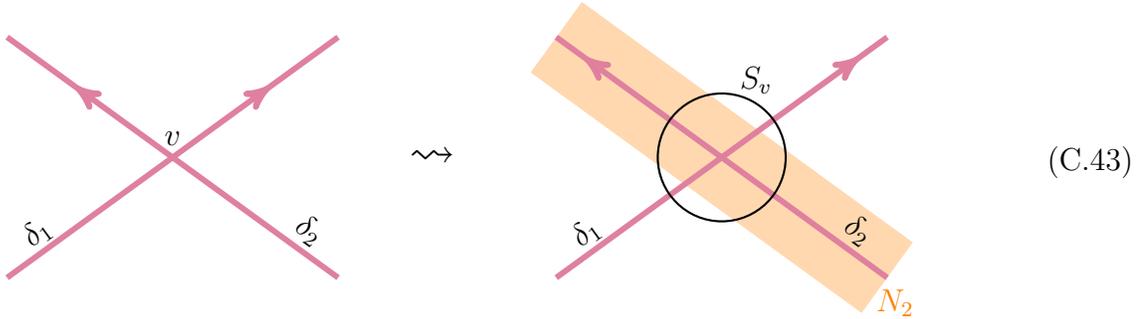
{\begin{tikzpicture}
  \coordinate (defec1a) at (0,0);
  \coordinate (defec1o) at (\locpx,\locpy);
  \coordinate (defec2a) at (0,\locpy);
  \coordinate (defec2o) at (\locpx,0);
  \scopeArrow{0.79}{\arrowDefect}
  \draw[line width=\locpc,color=\colorTransp,postaction={decorate}]
         (defec1a) -- node[sloped,above,very near start,yshift=-3pt,color=black] {$\delta_1$} (defec1o);
  \draw[line width=\locpc,color=\colorTransp,postaction={decorate}]
         (defec2o) -- node[sloped,above,very near start,yshift=-2pt,color=black] {$\delta_2$} (defec2a);
  \node at (0.5*\locpx,0.5*\locpy+0.26) {$v$};
  \end{scope}
  \begin{scope}[shift={(\locpx+2.9,0)}]
  \coordinate (defec1a) at (0,0);
  \coordinate (defec1o) at (\locpx,\locpy);
  \coordinate (defec2a) at (0,\locpy);
  \coordinate (defec2o) at (\locpx,0);
  \draw[line width=15*\locpc,color=orange!31]
         (defec2a) -- (defec2o) node[below=-16pt,xshift=3pt,color=orange] {$N_2$};
  \scopeArrow{0.91}{\arrowDefect}
  \draw[line width=\locpc,color=\colorTransp,postaction={decorate}]
         (defec1a) -- node[sloped,above,very near start,yshift=-3pt,color=black] {$\delta_1$} (defec1o);
  \draw[line width=\locpc,color=\colorTransp,postaction={decorate}]
         (defec2o) -- node[sloped,above,very near start,yshift=-2pt,color=black] {$\delta_2$} (defec2a);
  \end{scope}
  \draw[thick] (0.5*\locpx,0.5*\locpy) circle (11*\locpc) node[above=20pt,xshift=13pt] {$S_v$};
  \node at (-1.65,0.5*\locpy) {{\Large$\rightsquigarrow$}};
  \end{scope}
  \end{tikzpicture}}
  \ee
Since $\Sigma_1$ and $\Sigma_2$ both refine $\Sigma$, there is a homotopy 
$h_t\colon T_p\vSigma \,{\to}\, T_p\vSigma$ for $t \,{\in}\, [0,1]$ and all $p \,{\in}\, \vSigma$
satisfying $h_0 \,{=}\, \id$ and $h_1(\chi_2) \,{=}\, \chi_1$.
Now let $N_2' \,{\subset}\, N_2^{}$ be a smaller tubular neighborhood of the defects $\delta_2$ and 
$b\colon \vSigma \,{\to}\, [0,1]$ a smooth monotonous function that is $0$ on $N_2'$ and $1$
on $\vSigma{\setminus}N_2^{}$. Then by setting $\chi_{1,2}(p) \,{:=}\, h_{b(p)}(\chi_2(p))$ 
for $p \,{\in}\, \vSigma$ we obtain a vector field $\chi_{1,2}$ on $\vSigma$ that
looks like the framing of $\Sigma_i$, for $i\,{\in}\,\{1,2\}$,
around the defects of $\mathring\Sigma_{1,2}$ that correspond to the defects of $\Sigma_i$
and thus defines a framing on $\mathring\Sigma_{1,2}$ of the desired form.

To obtain a proper \defectsurface\ we still have to get rid of the intersections $v$ between 
defect lines. To this end we remove for each such point $v$ the interior of the disk bounded 
by $S_v$ from $\mathring\Sigma_{1,2}$ and replace $S_v$ by a gluing circle, with appropriate
defect points at the intersection of $S_v$ with $\delta_1$ and $\delta_2$. Now notice that the
vector field $\chi_{1,2}$ is such that all the thus obtained gluing circles $\SS_v$ are fillable. This
means that after forgetting all transparent defects the framing of $\Sigma_{1,2}$ is by 
construction homotopic to the framing of $\Sigma$; hence we have indeed constructed
a refinement $(\Sigma;\Sigma_{1,2})$ of $\Sigma$.
Also, if $\Sigma_1$ and $\Sigma_2$ are fine, then so is $\Sigma_{1,2}$.

Moreover, each of the resulting `four-valent' fillable gluing circles $\SS_v$ in the common 
subrefinement $(\Sigma;\Sigma_{1,2})$ can be `resolved' to a pair of two three-valent fillable
gluing circles. This can be done in two specific ways, as indicated in 
    \def\locpa  {1.6}    
    \def\locpc  {2.2pt}  
    \def\locpd  {0.12}   
    \def\locpu  {-10}    
    \def\locpv  {80}     
    \def\locpw  {170}    
    \def\locpx  {260}    
  \be
  \raisebox{-7.9em}
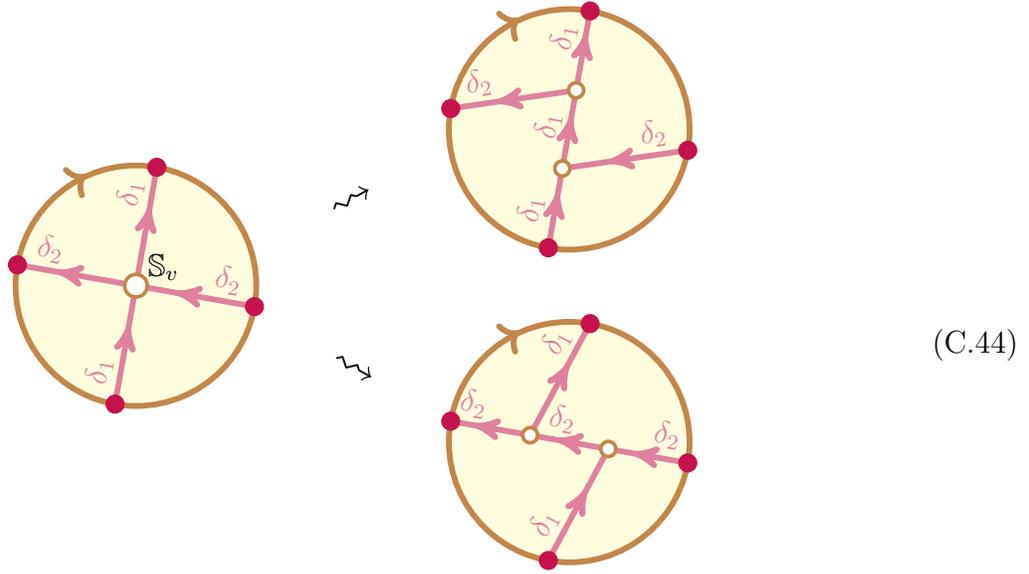
{\begin{tikzpicture}
  \node[rotate=29] at (0,0.7*\locpa) {{\Large$\rightsquigarrow$}};
  \node[rotate=-29] at (0,-0.7*\locpa) {{\Large$\rightsquigarrow$}};
  \begin{scope}[shift={(-1.8*\locpa,0)}]
  \scopeArrow{0.34}{<} 
  \filldraw[\colorCircle,fill=\colorSigma,line width=\locpc,postaction={decorate}]
                  (0,0) circle (\locpa);
  \end{scope} 
  \scopeArrow{0.33}{\arrowDefect} \scopeArrow{0.83}{\arrowDefect}
  \draw[line width=\locpc,color=\colorTransp,postaction={decorate}] 
		  (\locpu:\locpa) -- node[sloped,above,very near end,yshift=-2pt] {$\delta_2$}
		  node[sloped,above,very near start,yshift=-2pt] {$\delta_2$} (\locpw:\locpa);
  \draw[line width=\locpc,color=\colorTransp,postaction={decorate}] 
                  (\locpx:\locpa) -- node[sloped,above,very near end,yshift=-2pt] {$\delta_1$} 
	          node[sloped,above,very near start,yshift=-2pt] {$\delta_1$} (\locpv:\locpa);
  \end{scope} \end{scope}
  \filldraw[color=\colorDefect] (\locpu:\locpa) circle (\locpd) (\locpv:\locpa) circle (\locpd)
                                (\locpw:\locpa) circle (\locpd) (\locpx:\locpa) circle (\locpd);
  \def\transvR {0.15} %
  \transv {0,0} node[right,yshift=7pt,color=black] {$\SS_v$};
  \def\transvR {0.10} 
  \end{scope}
  \begin{scope}[shift={(1.8*\locpa,1.3*\locpa)}]
  \scopeArrow{0.34}{<} 
  \filldraw[\colorCircle,fill=\colorSigma,line width=\locpc,postaction={decorate}] (0,0) circle (\locpa);
  \end{scope} 
  \scopeArrow{0.23}{\arrowDefect} \scopeArrow{0.55}{\arrowDefect} \scopeArrow{0.90}{\arrowDefect}
  \draw[line width=\locpc,color=\colorTransp,postaction={decorate}] 
                  (\locpx:\locpa) -- node[sloped,above,very near end,yshift=-2pt] {$\delta_1$} 
	          node[sloped,above,midway,yshift=-2pt] {$\delta_1$}
	          node[sloped,above,very near start,yshift=-1pt] {$~\delta_1$} (\locpv:\locpa);
  \end{scope} \end{scope} \end{scope}
  \scopeArrow{0.63}{\arrowDefect}
  \draw[line width=\locpc,color=\colorTransp,postaction={decorate}] (\locpu:\locpa)
                  -- node[sloped,above,near start,yshift=-2pt] {$\delta_2$} (\locpx:0.33*\locpa);
  \draw[line width=\locpc,color=\colorTransp,postaction={decorate}] (\locpv:0.33*\locpa)
                  -- node[sloped,above,near end,yshift=-2pt] {$\delta_2$} (\locpw:\locpa);
  \end{scope}
  \filldraw[color=\colorDefect] (\locpu:\locpa) circle (\locpd) (\locpv:\locpa) circle (\locpd)
                                (\locpw:\locpa) circle (\locpd) (\locpx:\locpa) circle (\locpd);
  \transv {\locpv:0.33*\locpa}; \transv {\locpx:0.33*\locpa};
  \end{scope}
  \begin{scope}[shift={(1.8*\locpa,-1.3*\locpa)}]
  \scopeArrow{0.34}{<} 
  \filldraw[\colorCircle,fill=\colorSigma,line width=\locpc,postaction={decorate}] (0,0) circle (\locpa);
  \end{scope} 
  \scopeArrow{0.23}{\arrowDefect} \scopeArrow{0.55}{\arrowDefect} \scopeArrow{0.90}{\arrowDefect}
  \draw[line width=\locpc,color=\colorTransp,postaction={decorate}] 
                  (\locpu:\locpa) -- node[sloped,above,very near end,yshift=-2pt] {$\delta_2~~$}
                  node[sloped,above,midway,yshift=-2pt] {$\delta_2~~$}
                  node[sloped,above,very near start,yshift=-1pt] {$~\delta_2$} (\locpw:\locpa);
  \end{scope} \end{scope} \end{scope}
  \scopeArrow{0.63}{\arrowDefect}
  \draw[line width=\locpc,color=\colorTransp,postaction={decorate}] (\locpx:\locpa) 
                  -- node[sloped,above,near start,yshift=-2pt] {$\delta_1$} (\locpu:0.33*\locpa);
  \draw[line width=\locpc,color=\colorTransp,postaction={decorate}] (\locpw:0.33*\locpa)
                  -- node[sloped,above,near end,yshift=-2pt] {$\delta_1$} (\locpv:\locpa);
  \end{scope}
  \filldraw[color=\colorDefect] (\locpu:\locpa) circle (\locpd) (\locpv:\locpa) circle (\locpd)
                                (\locpw:\locpa) circle (\locpd) (\locpx:\locpa) circle (\locpd);
  \transv {\locpu:0.33*\locpa}; \transv {\locpw:0.33*\locpa};
  \end{scope}
  \end{tikzpicture}}
  \label{eq:resolve}
  \ee ~\\[0.2em]

It will be convenient to have separate terminology for specific manipulations of defect networks:

\begin{Definition}\label{def:resolve}
Let $(\Sigma;\Sigma_1)$ and $(\Sigma;\Sigma_2)$ be refinements refining the same
\defectsurface\ $\Sigma$, and $(\Sigma;\Sigma_{1,2})$ a common subrefinement.
\Itemizeiii

\item[{\rm (i)}]
We call the change of defect mesh shown in \eqref{eq:resolve} the \emph{resolvement} of
the four-valent gluing circle $\SS_v$ to $\Sigma_1$ and to $\Sigma_2$, respectively. 

\item[{\rm (ii)}]
We call a refinement replacement $\DR_{\DD,\widetilde\DD}$ of \emph{creation type},
respectively of \emph{annihilation type}, iff the defects on $\widetilde\DD$ are
obtained by adding defects to, respectively deleting transparent defects from, the disk $\DD$.
\end{itemize}
\end{Definition}

Given any two refinements $(\Sigma;\Sigma_1)$ and $(\Sigma;\Sigma_2)$ and a common subrefinement 
$(\Sigma;\Sigma_{1,2})$, we obtain a specific refinement replacement from $(\Sigma;\Sigma_1)$
to $(\Sigma;\Sigma_2)$ by the following two steps: First, perform resolvements, in the sense of
Definition \ref{def:resolve}, of each four-valent gluing circle that arises in the construction
of $\Sigma_{1,2}$ to $\Sigma_1$. Next perform local creation-type \filldiskrep s by adding
single $\Sigma_2$-defects to $\Sigma_1$, as indicated in the following picture
which shows disks that arise from a tubular neighborhood of the defect line in $\Sigma_2$:
    \def\locpa  {1.6}    
    \def\locpA  {0.6}    
    \def\locpc  {2.2pt}  
    \def\locpd  {0.12}   
    \def\locpu  {-20}    
    \def\locpv  {80}     
    \def\locpw  {160}    
    \def\locpx  {260}    
  \be
  \raisebox{-4.8em}
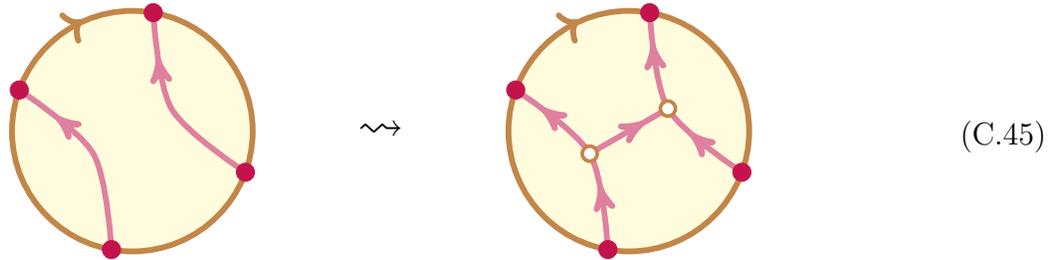
{\begin{tikzpicture}[scale=1.2]
  \node at (0,0) {{\Large$\rightsquigarrow$}};
  \begin{scope}[shift={(-1.8*\locpa-.1,0)}]
  \scopeArrow{0.34}{<} 
  \filldraw[\colorCircle,fill=\colorSigma,line width=\locpc,postaction={decorate}] (0,0) circle (\locpa);
  \end{scope} 
  \scopeArrow{0.77}{\arrowDefect}
  \draw[line width=\locpc,color=\colorTransp,postaction={decorate}] plot[smooth,tension=0.55]
                  coordinates{(\locpu:\locpa) (0.5*\locpu+0.5*\locpv:\locpA) (\locpv:\locpa)};
  \draw[line width=\locpc,color=\colorTransp,postaction={decorate}] plot[smooth,tension=0.55]
                  coordinates{(\locpx:\locpa) (0.5*\locpw+0.5*\locpx:\locpA) (\locpw:\locpa)};
  \end{scope}
  \filldraw[color=\colorDefect] (\locpu:\locpa) circle (\locpd) (\locpv:\locpa) circle (\locpd)
                                (\locpw:\locpa) circle (\locpd) (\locpx:\locpa) circle (\locpd);
  \end{scope}
  \begin{scope}[shift={(1.8*\locpa+.1,0)}]
  \scopeArrow{0.34}{<} 
  \filldraw[\colorCircle,fill=\colorSigma,line width=\locpc,postaction={decorate}] (0,0) circle (\locpa);
  \end{scope} 
  \scopeArrow{0.33}{\arrowDefect} \scopeArrow{0.84}{\arrowDefect}
  \draw[line width=\locpc,color=\colorTransp,postaction={decorate}] plot[smooth,tension=0.55]
                  coordinates{(\locpu:\locpa) (0.5*\locpu+0.5*\locpv:\locpA) (\locpv:\locpa)};
  \draw[line width=\locpc,color=\colorTransp,postaction={decorate}] plot[smooth,tension=0.55]
                  coordinates{(\locpx:\locpa) (0.5*\locpw+0.5*\locpx:\locpA) (\locpw:\locpa)};
  \end{scope} \end{scope}
  \scopeArrow{0.70}{\arrowDefect}
  \draw[line width=\locpc,color=\colorTransp,postaction={decorate}]
                  (0.5*\locpw+0.5*\locpx:\locpA) -- (0.5*\locpu+0.5*\locpv:\locpA);
  \end{scope}
  \filldraw[color=\colorDefect] (\locpu:\locpa) circle (\locpd) (\locpv:\locpa) circle (\locpd)
                                (\locpw:\locpa) circle (\locpd) (\locpx:\locpa) circle (\locpd);
  \transv {.5*\locpu+0.5*\locpv:\locpA}; 
  \transv {.5*\locpw+0.5*\locpx:\locpA};
  \end{scope}
  \end{tikzpicture}}
  \label{eq:repl:creation}
  \ee
Finally use replacements of annihilation type to resolve back to the four-valent gluing circles
of $\Sigma_{1,2}$. We refer to this procedure as a
\emph{standard refinement replacement} from $(\Sigma;\Sigma_1)$ to $(\Sigma;\Sigma_{1,2})$.

\begin{Lemma}\label{lem:refrepl}
For any two refinements $(\Sigma;\Sigma_1)$ and $(\Sigma;\Sigma_2)$ that refine the same 
\defectsurface\ $\Sigma$ there exists a refinement replacement $(\DR_1,...\,,\DR_n)$ from 
$(\Sigma;\Sigma_1)$ to $(\Sigma;\Sigma_2)$.
\end{Lemma}

\begin{proof}
Choose any common subrefinement $(\Sigma;\Sigma_{1,2})$ of $(\Sigma;\Sigma_1)$ and $(\Sigma;\Sigma_2)$.
Composing the standard refinement replacement from $(\Sigma;\Sigma_1)$ to $(\Sigma;\Sigma_{1,2})$
with the inverse of the standard replacement from $(\Sigma;\Sigma_2)$ to $(\Sigma;\Sigma_{1,2})$
gives a refinement replacement from $(\Sigma;\Sigma_1)$ to $(\Sigma;\Sigma_2)$ that 
factors through $(\Sigma;\Sigma_{1,2})$.
\end{proof}

We are now almost ready to show that any two refinement replacements are equivalent. 
Before giving the proof we just introduce some further convenient terminology. 

\begin{Definition}
Let $\Sigma$ be a \defectsurface.
\Itemizeiii

\item[{\rm (i)}]
Let $\delta$ be a set of defect lines on $\Sigma$. We say that a \filldiskrep\ $\DR_{\DD, \DD'}$
on $\Sigma$ \emph{keeps the defects $\delta$} iff each of the defects in $\delta$ corresponds to a
defect on $\DR_{\DD,\DD'}(\Sigma)$, possibly interrupted by gluing circles that are not 
present in $\Sigma$ $($like e.g.\ in the refinement shown in the picture \eqref{eq:picFeb1}$)$.

\item[{\rm (ii)}] 
Analogously we say that a sequence of \filldiskrep s keeps $\delta$ iff each of its members
keeps $\delta$. 

\item[{\rm (iii)}] 
Let $\PP$ be a \twocell\ $($in the sense of Definition $\ref{def:twocell})$ of $\Sigma$.
A disk $\DD$ on $\Sigma$ is said to be \emph{local with respect to $\PP$}, or 
\emph{$\PP$-local}, for short, iff $\,\DD$ does not meet a gluing circle on $\partial \PP$ 
and there are no defects in  $\DD \,{\setminus}\, \PP$. 
A \filldiskrep\ $\DR_{\DD,\DD'}$ on a $\PP$-local disk $\DD$ is said to be \emph{$\PP$-local}
iff $\DR_{\DD,\DD'}$ keeps the defects on $\partial \PP$ and $\DD'$ is $\PP$-local as well
$($that is, no defects are created in $\DD \,{\setminus}\, \PP)$.

\item[{\rm (iv)}] 
Analogously we say that a sequence of \filldiskrep s on $\Sigma$ is \emph{$\PP$-local} iff 
each of its members is $\PP$-local.
\end{itemize}
\end{Definition}

The following picture gives an example of a $\PP$-local \filldiskrep:
    \def\locpa  {3.7}    
    \def\locpA  {2.4}    
    \def\locpb  {0.35}   
    \def\locpB  {0.22}   
    \def\locpc  {2.2pt}  
    \def\locpd  {0.10}   
    \def\locpD  {0.08}   
  \be
  \begin{tikzpicture}[scale=1.1]
  \node at (-0.2,0.65) {\Large$\rightsquigarrow$};
  \begin{scope}[shift={(-1.6*\locpA,0)}]
  \begin{scope}[decoration={random steps,segment length=3mm}]
  \fill[\colorSigma,decorate] (90:\locpa) -- (210:\locpa) -- (330:\locpa) -- cycle;
  \end{scope}
  \draw[\colorDefect,line width=\locpc,dashed] 
           (90:\locpA) -- (90:\locpa) (210:\locpA) -- (210:\locpa) (330:\locpA) -- (330:\locpa);
  \scopeArrow{0.29}{\arrowDefect} \scopeArrow{0.77}{\arrowDefect}
  \filldraw[\colorDefect,fill=\colorSigma,line width=\locpc,postaction={decorate}]
           (90:\locpA) -- (210:\locpA) -- (330:\locpA);
  \end{scope} \end{scope}
  \node at (0.03*\locpA,0.25*\locpA) {\large$\PP$};
  \scopeArrow{0.50}{\arrowDefectB}
  \filldraw[\colorDefect,fill=\colorSigma,line width=\locpc,postaction={decorate}]
           (330:\locpA) -- (90:\locpA);
  \end{scope}
  \filldraw[\colorCircle,fill=white,line width=\locpc] (90:\locpA) circle (\locpb)
           (210:\locpA) circle (\locpb) (330:\locpA) circle (\locpb);
  \scopeArrow{0.53}{>}
  \draw[\colorCircle,line width=0.5*\locpc,postaction={decorate}] (90:\locpA) circle (\locpb);
  \end{scope}
  \scopeArrow{0.82}{>}
  \draw[\colorCircle,line width=0.5*\locpc,postaction={decorate}] (210:\locpA) circle (\locpb);
  \end{scope}
  \scopeArrow{0.76}{>}
  \draw[\colorCircle,line width=0.5*\locpc,postaction={decorate}] (330:\locpA) circle (\locpb);
  \end{scope}
  \filldraw[color=\colorDefect]
           (90:\locpA+\locpb) circle (\locpd) (210:\locpA+\locpb) circle (\locpd)
           (330:\locpA+\locpb) circle (\locpd)
           (90:\locpA)+(240:\locpb) circle (\locpd) (90:\locpA)+(300:\locpb) circle (\locpd)
           (210:\locpA)+(0:\locpb) circle (\locpd) (210:\locpA)+(60:\locpb) circle (\locpd) 
           (330:\locpA)+(120:\locpb) circle (\locpd) (330:\locpA)+(180:\locpb) circle (\locpd);
  \end{scope}
  \begin{scope}[shift={(1.6*\locpA,0)}]
  \begin{scope}[decoration={random steps,segment length=3mm}]
	  \fill[\colorSigma,decorate] (90:\locpa) -- (150:0.68*\locpa) -- (210:\locpa)
          -- (270:0.68*\locpa) -- (330:\locpa) -- cycle;
  \end{scope}
  \fill[\colorDefect,fill=\colorSigma] (90:\locpA) -- (210:\locpA) -- (330:\locpA);
  \filldraw[line width=0.5*\locpc,color=\colorCircle,dotted,fill=\colorSigmaDark,rotate=-60]
          (270:0.25*\locpA) circle (0.8*\locpA cm and 0.4*\locpA cm);
  \node at (-0.11*\locpA,0.1*\locpA) {\large$\DD$};
  \draw[\colorDefect,line width=\locpc,dashed] 
           (90:\locpA) -- (90:\locpa) (210:\locpA) -- (210:\locpa) (330:\locpA) -- (330:\locpa);
  \scopeArrow{0.13}{\arrowDefect} \scopeArrow{0.43}{\arrowDefect} \scopeArrow{0.62}{\arrowDefect}
  \scopeArrow{0.92}{\arrowDefect}
  \draw[\colorDefect,line width=\locpc,postaction={decorate}]
           (90:\locpA) -- (210:\locpA) -- (330:\locpA);
  \end{scope} \end{scope} \end{scope} \end{scope}
  \scopeArrow{0.50}{\arrowDefectB}
  \filldraw[\colorDefect,fill=\colorSigma,line width=\locpc,postaction={decorate}]
           (330:\locpA) -- (90:\locpA);
  \end{scope}
  \scopeArrow{0.65}{\arrowDefect}
  \filldraw[\colorTransp,fill=\colorSigma,line width=\locpc,postaction={decorate}]
	   (150:0.5*\locpA) -- (270:0.5*\locpA);
  \end{scope}
  \filldraw[\colorCircle,fill=white,line width=\locpc] (90:\locpA) circle (\locpb)
           (210:\locpA) circle (\locpb) (330:\locpA) circle (\locpb)
	   (150:0.5*\locpA) circle (\locpB) (270:0.5*\locpA) circle (\locpB);
  \scopeArrow{0.53}{>}
  \draw[\colorCircle,line width=0.5*\locpc,postaction={decorate}] (90:\locpA) circle (\locpb);
  \draw[\colorCircle,line width=0.5*\locpc,postaction={decorate}] (150:0.5*\locpA) circle (\locpB);
  \end{scope}
  \scopeArrow{0.82}{>}
  \draw[\colorCircle,line width=0.5*\locpc,postaction={decorate}] (210:\locpA) circle (\locpb);
  \draw[\colorCircle,line width=0.5*\locpc,postaction={decorate}] (270:0.5*\locpA) circle (\locpB);
  \end{scope}
  \scopeArrow{0.76}{>}
  \draw[\colorCircle,line width=0.5*\locpc,postaction={decorate}] (330:\locpA) circle (\locpb);
  \end{scope}
  \filldraw[color=\colorDefect]
           (90:\locpA+\locpb) circle (\locpd) (210:\locpA+\locpb) circle (\locpd)
           (330:\locpA+\locpb) circle (\locpd)
           (90:\locpA)+(240:\locpb) circle (\locpd) (90:\locpA)+(300:\locpb) circle (\locpd)
           (210:\locpA)+(0:\locpb) circle (\locpd) (210:\locpA)+(60:\locpb) circle (\locpd) 
           (330:\locpA)+(120:\locpb) circle (\locpd) (330:\locpA)+(180:\locpb) circle (\locpd)
           (150:0.5*\locpA)+(60:\locpB) circle (\locpD) (150:0.5*\locpA)+(240:\locpB) circle (\locpD)
           (150:0.5*\locpA)+(300:\locpB) circle (\locpD)
           (270:0.5*\locpA)+(0:\locpB) circle (\locpD) (270:0.5*\locpA)+(120:\locpB) circle (\locpD)
           (270:0.5*\locpA)+(180:\locpB) circle (\locpD);
  \end{scope}
  \end{tikzpicture} 
  \label{eq:exa-filldiskrep}
  \ee ~\\[-0.5em]

\begin{Lemma} \label{Lemma:basic-equivalence}
Let $\Sigma$ be a fine \defectsurface\ and $(\Sigma;\Sigma\refi)$ be a refinement that
refines $\Sigma$. 
\Itemizeiii

\item[{\rm (i)}]
Any sequence $(\DR_{1}, ...\,, \DR_n)$ of \filldiskrep s from $\Sigma$ to
$\Sigma_{\rm ref}$ that keeps the defects of $\Sigma$ is equivalent to a sequence that is
local with respect to all \twocells\ of $\Sigma$.

\item[{\rm (ii)}]
Any two sequences $(\DR_{1}, ...\,, \DR_n)$ and $(\DR_1', ...\,, \DR_{n'}')$ of \filldiskrep s
from $\Sigma$ to $\Sigma\refi$ keeping the defects of $\Sigma$ are equivalent.
\end{itemize}
\end{Lemma}

\begin{proof}
(i)\, We ``localize'' $(\DR_{1}, ...\,, \DR_n)$ as follows with respect to the \twocells\ 
of $\Sigma$. Consider any of the \filldiskrep s $\DR_j$. Since, by assumption, $\DR_j$ keeps
the defects on $\Sigma$, there is a sequence $\{\DR_{j,s}\}$ of \filldiskrep s that are local
with respect to the \twocells\ of $\Sigma$, such that $\{\DR_{j,s}\}$ is equivalent to
$\DR_j$ by the factorization property of Proposition \ref{Proposition:main-disk-repl}. 
An illustration of this localization procedure is given in the following picture, in which 
the disks $\DD_1$ and $\DD_2$ contain a pair of neighboring defect circles that result from
the resolvement of a circle with four defect points, with $\DD_1$ being $\PP_1$-local and
$\DD_2$ being $\PP_2$-local, while $\DD$ is neither $\PP_1$- nor $\PP_2$-local:
    \def\locpa  {4.1}    
    \def\locpA  {3.1}    
    \def\locpb  {0.35}   
    \def\locpB  {0.22}   
    \def\locpc  {2.2pt}  
    \def\locpd  {0.10}   
    \def\locpD  {0.08}   
  \be
  \raisebox{-9.8em}
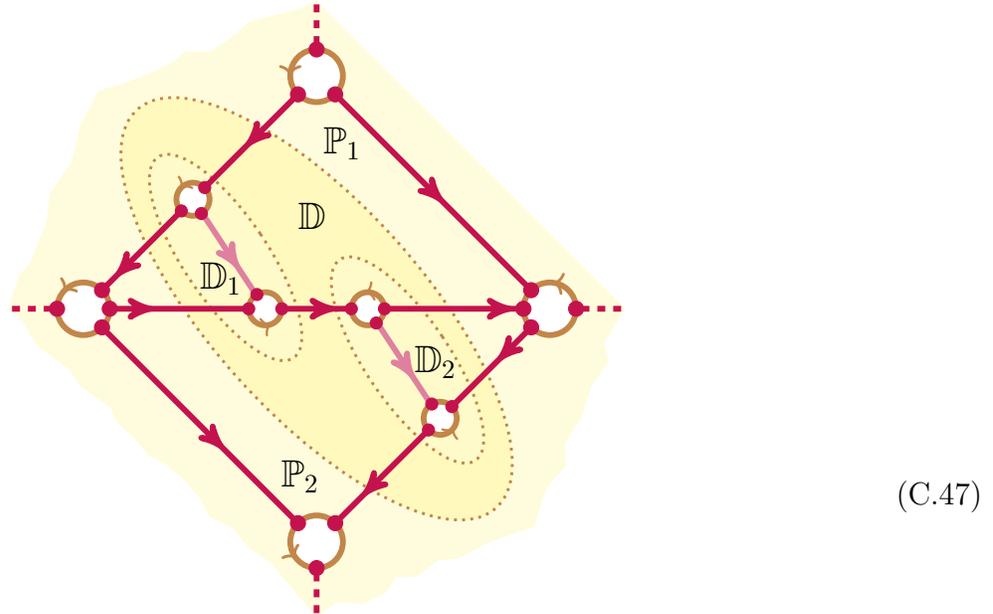
{\begin{tikzpicture}
  \coordinate (circNW) at (-0.53*\locpA,0.47*\locpA);
  \coordinate (circW) at (-0.22*\locpA,0);
  \coordinate (elliW) at (-0.39*\locpA,0.22*\locpA); 
  \coordinate (circSE) at (0.53*\locpA,-0.47*\locpA);
  \coordinate (circE) at (0.22*\locpA,0);
  \coordinate (elliE) at (0.39*\locpA,-0.22*\locpA);
  \begin{scope}[decoration={random steps,segment length=3mm}]
  \fill[\colorSigma,decorate] (0,\locpa) -- (135:\locpa) -- (-\locpa,0)
           -- (0,-\locpa) -- (315:\locpa) -- (\locpa,0);
  \end{scope}
  \filldraw[line width=0.5*\locpc,color=\colorCircle,dotted,fill=\colorSigmaDARK,rotate=42]
           (0,0) circle (0.45*\locpA cm and 1.15*\locpA cm);
  \filldraw[line width=0.5*\locpc,color=\colorCircle,dotted,fill=\colorSigmaDark,rotate=34]
           (elliW) circle (0.18*\locpA cm and 0.52*\locpA cm)
           (elliE) circle (0.18*\locpA cm and 0.52*\locpA cm);
  \draw[\colorDefect,line width=\locpc,dashed] (0,\locpA) -- (0,\locpa) 
           (0,-\locpA) -- (0,-\locpa) (\locpA,0) -- (\locpa,0) (-\locpA,0) -- (-\locpa,0);
  \scopeArrow{0.15}{\arrowDefect} \scopeArrow{0.43}{\arrowDefect} \scopeArrow{0.80}{\arrowDefect}
  \node at (0.11*\locpA,0.71*\locpA) {\large$\PP_1$};
  \node at (-0.07*\locpA,-0.72*\locpA) {\large$\PP_2$};
  \node at (-0.02*\locpA,0.4*\locpA) {\large$\DD$};
  \node at (-0.41*\locpA,0.13*\locpA) {\large$\DD_1$};
  \node at (0.51*\locpA,-0.23*\locpA) {\large$\DD_2$};
  \draw[\colorDefect,line width=\locpc,postaction={decorate}]
           (0,\locpA) -- (-\locpA,0) -- (0,-\locpA);
  \end{scope} \end{scope} \end{scope}
  \scopeArrow{0.27}{\arrowDefect} \scopeArrow{0.61}{\arrowDefect} \scopeArrow{0.90}{\arrowDefect}
  \draw[\colorDefect,line width=\locpc,postaction={decorate}]
           (0,\locpA) -- (\locpA,0) -- (0,-\locpA);
  \end{scope} \end{scope} \end{scope}
  \scopeArrow{0.15}{\arrowDefect} \scopeArrow{0.54}{\arrowDefect} \scopeArrow{0.92}{\arrowDefect}
  \draw[\colorDefect,line width=\locpc,postaction={decorate}] (-\locpA,0) -- (\locpA,0);
  \end{scope} \end{scope} \end{scope}
  \scopeArrow{0.62}{\arrowDefect}
  \draw[\colorTransp,line width=\locpc,postaction={decorate}] (circE) -- (circSE);
  \draw[\colorTransp,line width=\locpc,postaction={decorate}] (circNW) -- (circW);
  \end{scope}
  \filldraw[\colorCircle,fill=white,line width=\locpc]
           (0:\locpA) circle (\locpb) (90:\locpA) circle (\locpb)
           (180:\locpA) circle (\locpb) (270:\locpA) circle (\locpb)
           (circW) circle (\locpB) (circNW) circle (\locpB)
           (circE) circle (\locpB) (circSE) circle (\locpB);
  \scopeArrow{0.23}{>}
  \draw[\colorCircle,line width=0.5*\locpc,postaction={decorate}] (0:\locpA) circle (\locpb);
  \end{scope}
  \scopeArrow{0.28}{>}
  \draw[\colorCircle,line width=0.4*\locpc,postaction={decorate}] (circE) circle (\locpB);
  \end{scope}
  \scopeArrow{0.50}{>}
  \draw[\colorCircle,line width=0.5*\locpc,postaction={decorate}] (90:\locpA) circle (\locpb);
  \end{scope}
  \scopeArrow{0.39}{>}
  \draw[\colorCircle,line width=0.5*\locpc,postaction={decorate}] (180:\locpA) circle (\locpb);
  \draw[\colorCircle,line width=0.4*\locpc,postaction={decorate}] (circNW) circle (\locpB);
  \end{scope}
  \scopeArrow{0.61}{>}
  \draw[\colorCircle,line width=0.5*\locpc,postaction={decorate}] (270:\locpA) circle (\locpb);
  \end{scope}
  \scopeArrow{0.81}{>}
  \draw[\colorCircle,line width=0.4*\locpc,postaction={decorate}] (circW) circle (\locpB);
  \end{scope}
  \scopeArrow{0.91}{>}
  \draw[\colorCircle,line width=0.4*\locpc,postaction={decorate}] (circSE) circle (\locpB);
  \end{scope}
  \filldraw[color=\colorDefect]
           (0:\locpA+\locpb) circle (\locpd) (90:\locpA+\locpb) circle (\locpd)
           (180:\locpA+\locpb) circle (\locpd) (270:\locpA+\locpb) circle (\locpd)
           (0:\locpA)+(135:\locpb) circle (\locpd) (0:\locpA)+(180:\locpb) circle (\locpd)
           (0:\locpA)+(225:\locpb) circle (\locpd)
           (90:\locpA)+(225:\locpb) circle (\locpd) (90:\locpA)+(315:\locpb) circle (\locpd)
           (180:\locpA)+(0:\locpb) circle (\locpd) (180:\locpA)+(45:\locpb) circle (\locpd)
           (180:\locpA)+(315:\locpb) circle (\locpd)
           (270:\locpA)+(45:\locpb) circle (\locpd) (270:\locpA)+(135:\locpb) circle (\locpd)
           (circE)+(0:\locpB) circle (\locpD) (circE)+(180:\locpB) circle (\locpD)
           (circSE)+(45:\locpB) circle (\locpD) (circSE)+(225:\locpB) circle (\locpD)
           (circW)+(0:\locpB) circle (\locpD) (circW)+(180:\locpB) circle (\locpD)
           (circNW)+(45:\locpB) circle (\locpD) (circNW)+(225:\locpB) circle (\locpD)
           (circE)+(300.5:\locpB) circle (\locpD) (circSE)+(120.5:\locpB) circle (\locpD) 
           (circW)+(120.5:\locpB) circle (\locpD) (circNW)+(300.5:\locpB) circle (\locpD);
  \end{tikzpicture}}
  \ee
(ii)\, Owing to (i) we can without loss of generality assume that both sequences are local 
with respect to the defects in $\Sigma$. It is then enough to consider a single \twocell\ $\PP$ of
$\Sigma$. Let $(\DR_1, ...\,, \DR_p)$ and $(\DR_1', ...\,, \DR_{p'}')$ be two $\PP$-local 
sequences of disk replacements. Since by assumption $\Sigma$ is fine, $\PP$ is a disk which,
in turn, implies that there is a disk $\DD$ on $\Sigma$ such that both sequences lie entirely 
in $\DD$. Thus the two
sequences are equivalent by Proposition \ref{Proposition:main-disk-repl}.
Again we give an illustrative example:
    \def\locpa  {4.1}    
    \def\locpA  {3.0}    
    \def\locpb  {0.35}   
    \def\locpB  {0.22}   
    \def\locpc  {2.2pt}  
    \def\locpd  {0.10}   
    \def\locpD  {0.08}   
  \be
  \raisebox{-8.2em}
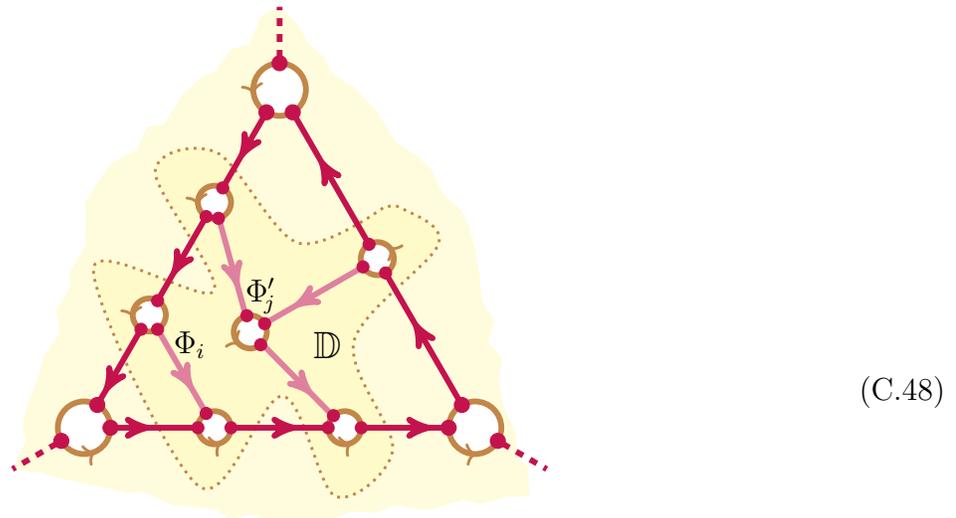
{\begin{tikzpicture}
  \coordinate (circ0) at (210:0.15*\locpA);
  \coordinate (circ1) at (30:0.5*\locpA);
  \coordinate (circ2) at (120:{\locpA/sqrt(3)});
  \coordinate (circ3) at (180:{\locpA/sqrt(3)});
  \coordinate (circ4) at (240:{\locpA/sqrt(3)});
  \coordinate (circ5) at (300:{\locpA/sqrt(3)});
  \coordinate (Circ1a) at (41:0.72*\locpA);
  \coordinate (Circ1b) at (46:0.66*\locpA);
  \coordinate (Circ1x) at (24:0.77*\locpA);
  \coordinate (Dirc1) at (76:0.33*\locpA);
  \coordinate (Circ2a) at (111:0.77*\locpA);
  \coordinate (Circ2b) at (131:0.84*\locpA);
  \coordinate (Dirc2) at (160:0.31*\locpA);
  \coordinate (Circ3) at (166:0.85*\locpA);
  \coordinate (Circ4) at (245:0.84*\locpA);
  \coordinate (Dirc4) at (272:0.36*\locpA);
  \coordinate (Circ5a) at (287:0.84*\locpA);
  \coordinate (Circ5b) at (310:0.78*\locpA);
  \coordinate (Dirc5) at (340:0.35*\locpA);
  \begin{scope}[decoration={random steps,segment length=3mm}]
	  \fill[\colorSigma,decorate] (90:\locpa) -- (144:0.73*\locpa) -- (210:\locpa)
          -- (270:0.68*\locpa) -- (324:\locpa) -- (30:0.68*\locpa) -- cycle;
  \end{scope}
  \filldraw[line width=0.5*\locpc,color=\colorCircle,dotted,fill=\colorSigmaDark]
	  plot[smooth,tension=0.55] coordinates{
          (Circ1b) (Dirc1) (Circ2a) (Circ2b) (Dirc2) (Circ3) (Circ4) (Dirc4)
          (Circ5a) (Circ5b) (Dirc5) (Circ1x) (Circ1a) (Circ1b) };
  \node at (328:0.25*\locpA) {\large$\DD$};
  \draw[\colorDefect,line width=\locpc,dashed] 
           (90:\locpA) -- (90:\locpa) (210:\locpA) -- (210:\locpa) (330:\locpA) -- (330:\locpa);
  \scopeArrow{0.067}{\arrowDefect} \scopeArrow{0.182}{\arrowDefect} \scopeArrow{0.297}{\arrowDefect}
  \scopeArrow{0.39}{\arrowDefect} \scopeArrow{0.515}{\arrowDefect} \scopeArrow{0.63}{\arrowDefect}
  \scopeArrow{0.77}{\arrowDefect} \scopeArrow{0.93}{\arrowDefect}
  \draw[\colorDefect,line width=\locpc,postaction={decorate}]
           (90:\locpA) -- (210:\locpA) -- (330:\locpA) -- cycle;
  \end{scope} \end{scope} \end{scope} \end{scope} \end{scope} \end{scope} \end{scope} \end{scope} 
  \scopeArrow{0.65}{\arrowDefect}
  \draw[\colorTransp,line width=\locpc,postaction={decorate}] (circ0) -- (circ5);
  \draw[\colorTransp,line width=\locpc,postaction={decorate}] (circ1) -- (circ0);
  \draw[\colorTransp,line width=\locpc,postaction={decorate}] (circ2) -- (circ0);
  \draw[\colorTransp,line width=\locpc,postaction={decorate}]
           (circ3) -- node[above,xshift=3pt,color=black]{$\Phi_i$} (circ4);
  \end{scope}
  \filldraw[\colorCircle,fill=white,line width=\locpc] (90:\locpA) circle (\locpb)
           (210:\locpA) circle (\locpb) (330:\locpA) circle (\locpb)
	   (circ0) circle (\locpB) node[above=2pt,xshift=4pt,color=black]{$\Phi_{\!j}'$}
           (circ1) circle (\locpB) (circ2) circle (\locpB)
	   (circ3) circle (\locpB) (circ4) circle (\locpB) (circ5) circle (\locpB);
  \scopeArrow{0.13}{>}
  \draw[\colorCircle,line width=0.5*\locpc,postaction={decorate}] (circ1) circle (\locpB);
  \end{scope}
  \scopeArrow{0.53}{>}
  \draw[\colorCircle,line width=0.5*\locpc,postaction={decorate}] (90:\locpA) circle (\locpb);
  \draw[\colorCircle,line width=0.5*\locpc,postaction={decorate}] (circ2) circle (\locpB);
  \draw[\colorCircle,line width=0.5*\locpc,postaction={decorate}] (circ3) circle (\locpB);
  \end{scope}
  \scopeArrow{0.64}{>}
  \draw[\colorCircle,line width=0.5*\locpc,postaction={decorate}] (circ0) circle (\locpB);
  \end{scope}
  \scopeArrow{0.76}{>}
  \draw[\colorCircle,line width=0.5*\locpc,postaction={decorate}] (330:\locpA) circle (\locpb);
  \end{scope}
  \scopeArrow{0.82}{>}
  \draw[\colorCircle,line width=0.5*\locpc,postaction={decorate}] (210:\locpA) circle (\locpb);
  \draw[\colorCircle,line width=0.5*\locpc,postaction={decorate}] (circ4) circle (\locpB);
  \draw[\colorCircle,line width=0.5*\locpc,postaction={decorate}] (circ5) circle (\locpB);
  \end{scope}
  \filldraw[color=\colorDefect]
           (90:\locpA+\locpb) circle (\locpd) (210:\locpA+\locpb) circle (\locpd)
           (330:\locpA+\locpb) circle (\locpd)
           (90:\locpA)+(240:\locpb) circle (\locpd) (90:\locpA)+(300:\locpb) circle (\locpd)
           (210:\locpA)+(0:\locpb) circle (\locpd) (210:\locpA)+(60:\locpb) circle (\locpd) 
           (330:\locpA)+(120:\locpb) circle (\locpd) (330:\locpA)+(180:\locpb) circle (\locpd)
           (circ0)+(30:\locpB) circle (\locpD) (circ0)+(102:\locpB) circle (\locpD)  
           (circ0)+(309:\locpB) circle (\locpD) 
           (circ1)+(120:\locpB) circle (\locpD) (circ1)+(210:\locpB) circle (\locpD)
           (circ1)+(300:\locpB) circle (\locpD)
           (circ2)+(60:\locpB) circle (\locpD) (circ2)+(240:\locpB) circle (\locpD)
           (circ2)+(284:\locpB) circle (\locpD) 
           (circ3)+(60:\locpB) circle (\locpD) (circ3)+(240:\locpB) circle (\locpD)
           (circ3)+(300:\locpB) circle (\locpD)
           (circ4)+(0:\locpB) circle (\locpD) (circ4)+(120:\locpB) circle (\locpD)
           (circ4)+(180:\locpB) circle (\locpD)
           (circ5)+(0:\locpB) circle (\locpD) (circ5)+(133:\locpB) circle (\locpD) 
           (circ5)+(180:\locpB) circle (\locpD);
  \end{tikzpicture}}
  \ee
This picture shows the \twocell\ $\PP$ and indicates the disk $\DD$ that encloses
both $\PP$-local \filldiskrep s $\DR_i^{}$ and $\DR_j'$,
which consist of one and three transparently labeled defect lines, respectively.
\end{proof}

\begin{Lemma} \label{Lemma:equiv-comm-sub}
Let $\Sigma$ be a fine \defectsurface. Any sequence $(\DR_1,\DR_2,...\,, \DR_n)$ of
fine \filldiskrep s on disks $\{\DD_i\}$ in $\Sigma$ is equivalent to a sequence 
$(\DR_1', ...\,, \DR_{n'}')$ of \filldiskrep s on $\{\DD_i\}$ such that, for some 
$1 \,{\le}\, p \,{<}\, q \,{\le}\, n'$ the \filldiskrep s $\DR_1',...\,, \DR_p'$ 
are replacements of creation type, $\DR_q',...\,, \DR_{n'}'$ are
of annihilation type, and $\DR_{p+1}',...\,, \DR_{q-1}'$ are resolvements of vertices. 
\end{Lemma}

\begin{proof}
We consider iteratively common subrefinements. We then need to show commutativity of a
diagram of the following form, in which the bottom row consists of the original sequence
$(\DR_1, ...\,, \DR_n)$ (depicted for the case $n\,{=}\,4$):
  \begin{eqnarray}
  \begin{tikzcd}[column sep= 1.1em]
  & & & \ar{dr} \ZZ(\Sigma_{1,n}) & & &
  \\
  & & \ZZ(\Sigma_{1,3}) \ar{dr} \ar{ur} & & \ZZ(\Sigma_{n-2,n}) \ar{dr} & &
  \\
  & \ZZ(\Sigma_{1,2}) \ar{dr} \ar{ur} & & \ZZ(\Sigma_{2,3}) \ar{dr}\ar{ur} & &
  \ZZ(\Sigma_{n-1,n}) \ar{dr} &
  \\
  \ZZ(\Sigma_{1}) \ar{rr}[swap]{\varphi^{}_{\Sigma_{1},\DR_{1}(\Sigma_{1})}} \ar{ur}
  & & \ZZ(\Sigma_{2}) \ar{rr}[swap]{\ldots} \ar{ur} & & \ZZ(\Sigma_{n-1}) 
  \ar{rr}[swap]{\varphi^{}_{\Sigma_{n-1},\DR_{n-1}(\Sigma_{n-1})}}\ar{ur} & & \ZZ(\Sigma_{n})
  \end{tikzcd}
  \nonumber\\[-7em]~ \label{eq:ref-diag} \\[4.4em]~\nonumber
  \end{eqnarray}
(In this diagram and in the rest of the proof, to save space we abuse notation
by just writing $\ZZ(\Sigma)$ in place of $\ZZ(\Sigma)(-\boti\om)$.)
We construct the diagram by proceeding from bottom to top. First, the triangle above the 
arrow labeled by $\varphi_{\Sigma_i,\DR_i(\Sigma_i)}$ is obtained by standard refinement
replacements on $\DR_i$: The arrow from $\ZZ(\Sigma_i)$ to $\ZZ(\Sigma_{i,i+1})$ is a 
replacement of creation type to the common subrefinement $\Sigma_{i,i+1}$ of $\Sigma_i$ 
and $\Sigma_{i+1}$, while the arrow from $\ZZ(\Sigma_{i,i+1})$ to $\ZZ(\Sigma_{i+1})$
is an analogous replacement of annihilation type. All three arrows in the
so obtained triangle are replacements inside one and the same disk, and hence
the triangle commutes by Proposition \ref{Proposition:main-disk-repl}.
\\[.2em]
Next consider a square above two consecutive triangles. It involves, besides
$\Sigma_i$ and the subrefinements $\Sigma_{i-1,i}$ and $\Sigma_{i,i+1}$, the common  
standard subrefinement $\Sigma_{i-1,i+1}$ of $\Sigma_{i-1,i}$ and $\Sigma_{i,i+1}$.
The arrow from $\ZZ(\Sigma_{i})$ to $\ZZ(\Sigma_{i,i+1})$ keeps the defects from $\Sigma_{i}$,
and likewise the composite of the other three arrows (with the first of them to be inverted) is a
sequence of \filldiskrep s from $\Sigma_{i}$ to $\Sigma_{i,i+1}$ that keeps the defects from $\Sigma_{i}$. 
Since $\Sigma_{i}$ is by assumption fine, it follows from Lemma \ref{Lemma:basic-equivalence}
that the square commutes. For any of the squares `higher up'' in the diagram, we can
likewise use the defects of the fine surface at the bottom of the square to invoke
Lemma \ref{Lemma:basic-equivalence}.
\\[.2em]
We have thus shown that all triangles and all squares in the diagram \eqref{eq:ref-diag}
commute, and hence the whole diagram commutes. Moreover, by construction the diagram
is of the required type. 
\end{proof}

We are now finally in a position to state

\begin{Proposition} \label{Proposition:main-generic-ref}
Let $(\Sigma;\Sigma_{\rm ref_1})$ and $(\Sigma;\Sigma_{\rm ref_2})$ be two refinements that
refine the same \defectsurface\ $\Sigma$. Any two refinement replacements
$(\DR_1',...\,,\DR_{n'}')$ and $(\DR_1'',...\,,\DR_{n''}'')$
from $\Sigma_{\rm ref_1}$ to $\Sigma_{\rm ref_2}$ are equivalent, i.e.\ they satisfy
  \be
  \varphi^{}_{ \Sigma_{\rm ref_1}, \DR_{n''}''(\cdots\, \DR_1''(\Sigma_{\rm ref_1}) \,\cdots)}
  = \varphi^{}_{ \Sigma_{\rm ref_1}, \DR_{n'}'(\cdots\, \DR_1'(\Sigma_{\rm ref_1}) \,\cdots)} \,.
  \ee
\end{Proposition}

\begin{proof}
We show that any sequence of refinement replacements $(\DR_1',...\,,\DR_{n'}')$ is equivalent
to the standard refinement replacement $(\DR_1,...\,,\DR_n)$, see Lemma  \ref{lem:refrepl}.
By Lemma \ref{Lemma:equiv-comm-sub} we can assume that $(\DR_1',...\,,\DR_{n'}')$ consists
first of replacements  $(\DR_1',...\,,\DR_k')$, for $k\,{\le}\,n'$, of creation type to a
refinement $(\Sigma;\Sigma_{\rm ref_{12}}')$, then of annihilation type refinements 
$(\DR_{l+1}',...\,,\DR_{n'}')$ $k\,{\le}\,l\,{\le}\,n'$, and in between
of resolvements of vertices. Likewise, the standard refinement replacement consists first
of creation type replacements $(\DR_1,...\,,\DR_p)$, for $p \,{\le}\, n$, to the common
subrefinement $\Sigma_{\rm ref_{12}}$ of $\Sigma_{\rm ref_1}$ and $\Sigma_{\rm ref_2}$
and then of annihilation type replacements $(\DR_{p+1},...\,,\DR_n)$ to $\Sigma_{\rm ref_2}$. 
The refinement $\Sigma_{\rm ref_{12}}'$ is necessarily a subrefinement of
$\Sigma_{\rm ref_{12}}$, thus there exists a sequence of refinement replacements
$(\DR_{q_1}, ...\,, \DR_{q_s})$ from $\Sigma_{\rm ref_{12}}$ to $\Sigma_{\rm ref_{12}}'$
keeping the defects from $\Sigma_{\rm ref_{12}}$. Consider then the sequence 
  \be
  (\DR_1, ...\,, \DR_p, \DR_{q_1}, ...\,, \DR_{q_s}, \DR_{q_s}^{-1}, ...\,,
  \DR_{q_1}^{-1}, \DR_{p+1}, ...\,, \DR_n)
  \ee
of replacements from $\Sigma_{\rm ref_{1}}$ to $\Sigma_{\rm ref_{2}}$. This sequence is
clearly equivalent to the standard refinement replacement, and
$(\DR_1, ...\,, \DR_p, \DR_{q_1}, ...\, \DR_{q_s})$ is a sequence of \filldiskrep s
from $\Sigma_{\rm ref_{1}}$ to $\Sigma_{\rm ref_{12}}'$ that keeps the defects of
$\Sigma_{\rm ref_{1}}$, just like $(\DR_{1}', ...\,, \DR_{k}')$. By Lemma
\ref{Lemma:basic-equivalence} they are thus equivalent. In the same way, the sequences
$(\DR_{q_s}^{-1}, ...\,, \DR_{q_1}^{-1}, \DR_{p+1}, ...\,, \DR_n)$ and
$(\DR_{k+1}', ...\,, \DR_{n'}')$ both keep the defects from $\Sigma_{\rm ref_2}$ and are 
thus equivalent as well (apply Lemma \ref{Lemma:basic-equivalence} to the inverses of the
sequences). Thus the statement follows. 
\end{proof}


\end{document}